\documentclass[sts,preprint]{imsart}
\usepackage[show]{simple_Notes}
\usepackage{lipsum}
\usepackage{graphicx} %
\usepackage{placeins} %
\usepackage{natbib}
\usepackage[greek,english]{babel}

\usepackage{comment}
\usepackage{enumerate}
\usepackage{enumitem}
\usepackage{dblfloatfix}
\usepackage{cancel}

\usepackage{fixltx2e}                       
\usepackage[usenames,dvipsnames]{xcolor}
\usepackage{fancyhdr}
\usepackage{amsmath,amsfonts,amsbsy,amsgen,amscd,mathrsfs,amssymb,amsthm}
\usepackage{url}
\usepackage{tikz}
\usepackage{microtype}
\usepackage{enumitem}

\definecolor{dark-gray}{gray}{0.3}
\definecolor{dkgray}{rgb}{.4,.4,.4}
\definecolor{dkblue}{rgb}{0,0,.5}
\definecolor{medblue}{rgb}{0,0,.75}
\definecolor{rust}{rgb}{0.5,0.1,0.1}

\usepackage[colorlinks=true]{hyperref}
\hypersetup{urlcolor=rust}
\hypersetup{citecolor=dkblue}
\hypersetup{linkcolor=dkblue}

\usepackage{booktabs,longtable,tabu} 
\usepackage{multirow} 
\usepackage{float}

\usepackage[full]{textcomp}

\usepackage{bm} 
\usepackage{subcaption}
\usepackage{overpic}

\usepackage{chngcntr}
\usepackage[symbol]{footmisc}
\usepackage{mathtools}
\counterwithout{equation}{section}

\usepackage[hmargin=1in,vmargin=1in]{geometry}

\usepackage{upgreek}
\usepackage{mathrsfs}
\theoremstyle{plain}
\newtheorem{theorem}{Theorem}

\newtheorem{proposition}{Proposition}
\newtheorem{assumption}{Assumption}

\newtheorem{lemma}{Lemma}
\numberwithin{lemma}{section}

\theoremstyle{definition}

\theoremstyle{remark}

\usepackage{soul}

\newcommand{\R}{\mathbb{R}}
\newcommand{\Q}{\mathbb{Q}}

\renewcommand{\hat}[1]{\widehat{#1}}

\newcommand{\I}[1]{\mathbbm{1}_{#1}}

\newcommand{\A}{\mathcal{A}}

\renewcommand{\Q}{\mathcal{Q}}

\renewcommand{\P}{\mathbb{P}}
\newcommand{\E}{\mathbb{E}}

\renewcommand{\I}{\mathcal{I}}

\newcommand{\e}{\epsilon}

\newcommand{\tr}{\operatorname{tr}}

\newcommand{\var}{\operatorname{var}}

\newcommand{\ttop}{^{\top}}

\newcommand{\ts}{\textstyle}

\newcommand{\op}{\textup{op}}

\newcommand{\blambda}{{\boldsymbol \lambda}}

\newcommand{\bdiag}{{\boldsymbol d}}
\newcommand{\mylessthan}{{\  \preceq t}}
\newcommand{\Proof}{\textbf{Proof}. }

\newcommand{\mnorm}[1]{\left\vert\kern-1.5pt\left\vert\kern-1.5pt\left\vert #1\right\vert\kern-1.5pt\right\vert\kern-1.5pt\right\vert}

\makeatletter
\newcommand*{\rom}[1]{\expandafter\@slowromancap\romannumeral #1@}
\makeatother

\makeatletter
\def\namedlabel#1#2{\begingroup
    #2%
    \def\@currentlabel{#2}%
    \phantomsection\label{#1}\endgroup
}
\makeatother

\makeatletter
\newcommand{\vast}{\bBigg@{3}}
\newcommand{\Vast}{\bBigg@{5}}
\makeatother

\begin{document}

\begin{frontmatter}
\centering
\textbf{ \textsf{\Huge  Rates of Bootstrap Approximation for\\[0.2cm] $\hspace{0.00cm}$ Eigenvalues in High-Dimensional PCA}}
~\\
~\\
~\\
\textbf{\textsf{\large Junwen Yao and Miles E. Lopes\footnote{Supported in part by NSF grant DMS 1915786.}}}\\[0.2cm]
\emph{University of California, Davis}

\begin{abstract}
In the context of principal components analysis (PCA), the bootstrap is commonly applied to solve a variety of inference problems, such as constructing confidence intervals for the eigenvalues of the population covariance matrix $\Sigma$. However, when the data are high-dimensional, there are relatively few theoretical guarantees that quantify the performance of the bootstrap. Our aim in this paper is to analyze how well the bootstrap can approximate the joint distribution of the leading eigenvalues of the sample covariance matrix $\hat\Sigma$, and we establish non-asymptotic rates of approximation with respect to the multivariate Kolmogorov metric. Under certain assumptions, we show that the bootstrap can achieve the dimension-free rate of ${\tt{r}}(\Sigma)/\sqrt n$ up to logarithmic factors, where ${\tt{r}}(\Sigma)$ is the effective rank of $\Sigma$, and $n$ is the sample size. 
From a methodological standpoint, our work also illustrates that applying a transformation to the eigenvalues of $\hat\Sigma$ before bootstrapping is an important consideration in high-dimensional settings.

\end{abstract}


\end{frontmatter}

\section{Introduction}
	
The applications of the bootstrap in principal components analysis (PCA) go back almost as far as the advent of the bootstrap itself~\citep{Diaconis:1983}, and over the years such applications have become part of standard practice in  multivariate analysis~\citep{Davison,Jolliffe,Olive}.
With regard to theory, there is also a well-established set of asymptotic results showing that the bootstrap generally works in the context of PCA with low-dimensional data~\citep{BeranS85,EatonT91}.
Furthermore, in aberrant situations where the bootstrap is known to encounter difficulty in low dimensions, such as in the case of tied population eigenvalues, various remedies have been proposed and analyzed~\citep{BeranS85,Dumbgen:1993,Hall:2009}.

 However, in the context of PCA with high-dimensional data, the relationship between theory and practice is quite different.  On one hand, bootstrap methods are popular among practitioners for solving inference problems related to high-dimensional PCA~\citep[e.g.][]{Wagner:2015,Fisher:2016,Webb:2017,Terry:2018,Ralph:2019,Holmes:2019,Stewart:2019}.
Yet, on the other hand, the theory for describing these methods is relatively incomplete.

As a way to develop a more precise understanding of the bootstrap in this context, we focus on the fundamental problem of approximating the joint distribution of the leading eigenvalues $\lambda_1(\hat\Sigma),\dots,\lambda_k(\hat\Sigma)$ of a sample covariance matrix $\hat\Sigma\in\R^{p\times p}$, where $k<p$. (Precise defnitions will be given later.) Because the fluctuations of these eigenvalues are relevant to many inference tasks, this problem plays a central role in multivariate analysis, and is also of broad interest in other areas,
such as signal processing~\citep{Debbah:Couillet:2011}, and finance~\citep{Ruppert:2015}.
For concreteness, we summarize below some examples of inference tasks involving sample eigenvalues. In addition, these tasks are illustrated with real-data examples based on stock market returns in Section~\ref{sec:stocks} of the supplementary material.

\begin{itemize}
\item  \textit{Selecting principal components.} A key step that occurs in any implementation of PCA is to choose the number of principal components, and many established techniques for making this choice are informed by the distributions of eigenvalue-based statistics. Examples of these statistics include eigengaps $\lambda_j(\hat\Sigma)-\lambda_{j+1}(\hat\Sigma)$, and the proportions of explained variance  $(\lambda_1(\hat\Sigma)+\cdots+\lambda_k(\hat\Sigma))/\tr(\hat\Sigma)$, as well as the componentwise proportions $\lambda_j(\hat\Sigma)/\tr(\hat\Sigma)$ for $j=1,\dots,k$. In addition, other selection rules are based on confidence intervals for the eigenvalues $\lambda_1(\Sigma),\dots,\lambda_k(\Sigma)$ of the population covariance matrix $\Sigma\in\R^{p\times p}$, and the construction of such intervals is directly linked to the distribution of the eigenvalues of $\hat\Sigma$.
For a general overview of selection rules, we refer to~\cite{Jolliffe}.

\item \textit{Quantifying uncertainty.} 
The eigenvalues of a population covariance matrix arise as unknown parameters of interest in many situations beyond the selection of principal components.
For instance, these parameters govern the performance of statistical methods for covariance estimation, regression, and classification~\citep{Ledoit:2012,Hsu:2014,Dobriban:2018}. Also, these parameters have domain-specific meaning in applications ranging from portfolio selection to ecology~\citep{Fabozzi,Ecology}.
Consequently, it becomes necessary to quantify the uncertainty associated with the population eigenvalues, such as in constructing confidence intervals for them---and again, this leads to the use of distributional approximation results for the sample eigenvalues. 

\end{itemize}

Although there is an extensive literature on distributional approximations for sample eigenvalues, this body of work primarily work focuses on asymptotic results involving analytical formulas. Roughly speaking, the bulk of the literature can be divided into two parts, dealing either with classical asymptotics where $p$ is held fixed as $n\to\infty$~\citep{Anderson}, or high-dimensional asymptotics where $p/n$ converges to a positive constant as $p$ and $n$ diverge simultaneously~\citep{BaiSilverstein}. In either case, an essential limitation is that asymptotic results do not usually quantify how close the limiting distribution is to the finite-sample distribution. In more practical terms, this means that it is often hard to know if tests statistics and confidence intervals are well calibrated (i.e.~if their actual levels and coverage probabilities are close to the nominal values). A second limitation is that approximations based on analytical formulas are often tied to specific model assumptions, which can make it difficult to adapt such formulas outside of a given model. 

With regard to the second limitation, bootstrap methods have an advantage insofar as they do not rely on formulas, and hence can be applied in a more flexible manner.  Nevertheless, the existing work on bootstrap methods for PCA still tends to suffer from the first limitation above, since the results are generally asymptotic~\citep{BeranS85,EatonT91,ElKarouiP19}. From this standpoint, a key motivation for our work is to provide results that explicitly quantify the accuracy of bootstrap approximation in terms of the sample size $n$ and the effective rank of $\Sigma$. (For example, our results can be used to quantify how close the coverage probabilities of bootstrap confidence intervals are to the nominal values.)  Another motivation is based on the fact that, until quite recently, most of the literature on bootstrap methods for PCA has been limited to low-dimensional settings. Consequently, it is of general interest to establish a more complete theoretical description of bootstrap methods for high-dimensional PCA---a point that was highlighted in a recent survey on this topic~\citep[][\textsection X.C]{JohnstoneP18}.

\subsection{Contributions}\label{sec:contrib}
Let $X_1,\dots,X_n\in\R^p$ be centered i.i.d.~observations with population covariance matrix $\Sigma=\E[X_1X_1\ttop]$. Also, let $\hat\Sigma=\ts\sum_{i=1}^nX_iX_i\ttop / n$ denote the associated sample covariance matrix, and let $\hat\Sigma^{\star}=\ts\sum_{i=1}^nX_i^{\star}(X_i^{\star})\ttop / n$ be its bootstrap version, formed from random vectors $X_1^{\star},\dots,X_n^{\star}$ that are sampled with replacement from the observations. In addition, let the eigenvalues of a symmetric matrix $A\in\R^{p\times p}$ be denoted as $\lambda_1(A)\geq\cdots\geq\lambda_p(A)$, and let $\boldsymbol \lambda_k(A)=(\lambda_1(A),\dots,\lambda_k(A))$ for a fixed integer $k<p$. 

In this notation, our goal is to establish non-asymptotic bounds on the multivariate Kolmogorov distance
\small
$$ \Delta_n= \sup_{t \in \R^k } \Bigg| \P \Big( \sqrt{n} \big( \blambda_k (\hat{\Sigma}) - \blambda_k (\Sigma) \big) \mylessthan \Big) - \P \Big( \sqrt{n} \big( \blambda_k (\hat{\Sigma}^\star) - \blambda_k (\hat\Sigma) \big) \mylessthan \, \Big| \, X \Big) \Bigg|,$$
\normalsize
where the relation $v\preceq w$ between two vectors $v,w\in\R^k$ means $v_j\leq w_j$ for all $j=1,\dots,k$, and  $\P(\cdot\,|X)$ refers to probability that is conditional on $X_1,\dots,X_n$. Under certain conditions, our central result (Theorem~\ref{THM:ALL}) shows that the dimension-free bound
\begin{equation}\label{eqn:informal}
\Delta_n \ \leq \ \frac{C_n{\,\tt{r}}(\Sigma)}{\sqrt n}
\end{equation}
holds with high probability, where $C_n>0$ is a polylogarithmic function of $n$, and the quantity ${\tt{r}}(\Sigma)$ is the effective rank of $\Sigma$, defined by ${\tt{r}}(\Sigma)=\tr(\Sigma) / \lambda_1(\Sigma)$.

There are several aspects of the bound~\eqref{eqn:informal} and the parameter ${\tt{r}}(\Sigma)$ that are worth noting. First, the effective rank satisfies \smash{$1\leq {\tt{r}}(\Sigma)\leq p$} whenever $\Sigma$ is nonzero, and can be interpreted as a proxy for the number of ``dominant'' principal components of $\Sigma$. Hence, even in very high-dimensional settings where $n\ll p$, the bound~\eqref{eqn:informal} shows that the bootstrap can perform well if the number of dominant components is not too large, which is precisely the situation where high-dimensional PCA is of greatest interest.
Meanwhile, even in situations where ${\tt{r}}(\Sigma)$ is moderately large, e.g.~${\tt{r}}(\Sigma)\to\infty$ with ${\tt{r}}(\Sigma)=o(\sqrt n)$, the bound~\eqref{eqn:informal} is still able to quantify the accuracy of the bootstrap. Indeed, both of these points are borne out by our numerical experiments in Section~\ref{sec:experiments}, which confirm that the performance of the bootstrap is governed more by ${\tt{r}}(\Sigma)$ than $p$, and that the bootstrap can still be accurate when ${\tt{r}}(\Sigma)$ is moderately large. More generally, it should also be mentioned that the theoretical role of effective rank in many other aspects of high-dimensional PCA has attracted considerable attention in recent years~\citep[e.g.][]{Lounici14, BuneaX15, KoltchinskiiL17b, JungA18, NaumovSU19, KoltchinskiiLN20}.

As an alternative to approximating the distribution of $\sqrt{n} \big( \blambda_k (\hat{\Sigma}) - \blambda_k (\Sigma) \big)$ by bootstrapping in a direct manner, it can be advantageous to use a \emph{transformation} prior to bootstrapping, which is a fundamental topic in the bootstrap literature~\citep[e.g.][]{Diciccio:1984, Tibshirani:1988,Konishi:1991,Diciccio:1996,Davison,Chernick}. To be more specific, let $h$ be a univariate scalar function, referred to as a transformation,
and for any symmetric matrix $A\in\R^{p\times p}$, let $\boldsymbol h(\boldsymbol\lambda_k(A))=(h(\lambda_1(A)),\dots,h(\lambda_k(A)))$. Then, the conditional distribution of $\sqrt n(\boldsymbol h(\boldsymbol\lambda_k(\hat\Sigma^{\star}))-\boldsymbol h(\boldsymbol\lambda_k(\hat\Sigma)))$ given the observations can be used to approximate the distribution of $\sqrt n(\boldsymbol h(\boldsymbol\lambda_k(\hat\Sigma))-\boldsymbol h(\boldsymbol\lambda_k(\Sigma)))$.
(Additional discussion is provided in Sections~\ref{sec:main} and~\ref{sec:experiments}.) For instance, a classical choice of transformation is $h(x)=\log(x)$, because it is known to be variance-stabilizing under certain conditions when $n\to\infty$ with $p$ held fixed~\citep{BeranS85}. 
With this in mind, a second contribution our analysis is an extended version of the bound~\eqref{eqn:informal} that can accommodate the use of certain transformations (Theorem~\ref{THM:TRANS_GENERAL}). 

From a more methodological standpoint, our numerical experiments also shed new light on the role of transformations in bootstrap methods for high-dimensional PCA. Although we confirm that the classical logarithm transformation can be beneficial in low dimensions, we show that it is less effective when ${\tt{r}}(\Sigma)$ is moderately large. Consequently, we explore some alternative transformations and provide numerical results demonstrating that there are opportunities to improve upon $h(x)=\log(x)$ in high dimensions. To put such empirical findings into perspective, we are not aware any prior work investigating how transformations can be used to enhance  bootstrap methods in this context.

\subsection{Related work} 
Quite recently, there has been an acceleration in the pace of
research on bootstrap methods for high-dimensional sample covariance matrices, as evidenced in the papers~\cite{HanXZ18,JohnstoneP18,ElKarouiP19,LopesBA19,LopesEM19,NaumovSU19}. Among these, the most relevant to our work is \cite{ElKarouiP19}, which examines both the successes and failures of the bootstrap in doing inference with the leading eigenvalues of $\hat\Sigma$. In the negative direction, that paper focuses  on a specialized model with $\lambda_1(\Sigma)>1$ and $\lambda_2(\Sigma)=\cdots=\lambda_p(\Sigma)=1$, which corresponds to a very large effective rank ${\tt{r}}(\Sigma)\asymp p$ that makes dimension reduction via PCA inherently difficult. In the positive direction, that paper deals with a different situation where $\Sigma$ is assumed to have a near low-rank structure of the form 
\begin{equation}\label{eqn:blockassume} 
\Sigma=\begin{pmatrix} A & B \\ B\ttop & C(\eta) \end{pmatrix},
\end{equation}
\normalsize
where $A$ is of size $k\times k$ with $k\asymp 1$, and the diagonal blocks satisfy $\lambda_1(A)\asymp 1$, and  $\lambda_1(C(\eta))\lesssim n^{-\eta}$ for a fixed parameter $\eta>1/2$. Working under an elliptical model, the paper~\citep{ElKarouiP19} shows that the bootstrap consistently approximates the distribution of $\sqrt{n} \big( \blambda_k (\hat{\Sigma}) - \blambda_k (\Sigma) \big)$ in an asymptotic framework where $p/n\lesssim 1$. In relation to our work, the most crucial distinction is that our results quantify the accuracy of the bootstrap with non-asymptotic \emph{rates of approximation}. To illustrate the significance of this, note that our bound~\eqref{eqn:informal} provides an explicit link between the size of ${\tt{r}}(\Sigma)$ and the accuracy bootstrap, whereas in an asymptotic setup, the effect of ${\tt{r}}(\Sigma)$ is hidden---because it ``washes out in the limit''. Our numerical experiments will also confirm that different sizes of ${\tt{r}}(\Sigma)$ can have an appreciable effect on the finite-sample accuracy of the bootstrap. In this way, our work indicates that the quantity ${\tt{r}}(\Sigma)/\sqrt n$ serves as a type of conceptual diagnostic for assessing the reliability of the bootstrap in high-dimensional PCA.

Beyond these points of contrast with~\cite{ElKarouiP19}, there are several distinctions with regard to model assumptions. First, we work in a dimension-free setting where there are no restrictions on the size of $p$ with respect to $n$. Second, the model based on~\eqref{eqn:blockassume} implicitly requires that $\lambda_j(\Sigma)\lesssim n^{-\eta}$ for all $j\geq k+1$, whereas this constraint on $\Sigma$ is not used here. 
Third, it is straightforward to check that in the model based on~\eqref{eqn:blockassume}, the condition $\eta>1/2$ \emph{implies} ${\tt{r}}(\Sigma)=o(\sqrt n)$, which means that our bound~\eqref{eqn:informal} ensures bootstrap consistency in models that subsume the one based on~\eqref{eqn:blockassume}.
(As an example, if $p\asymp e^{m(n)}$ for some sequence of integers satisfying $m(n)=o(\sqrt n)$  and if $\lambda_j(\Sigma) \asymp j^{-1}$, then the bound~\eqref{eqn:informal} implies bootstrap consistency, whereas this is not guaranteed by the previous result even when $p\asymp n$.)

Other works on bootstrap methods related to high-dimensional sample covariance matrices have dealt with models or statistics that are qualitatively different from those considered here. The papers~\cite{HanXZ18,LopesEM19} look at bootstrapping the operator norm error $\sqrt n\|\hat\Sigma-\Sigma\|_\op$, as well as variants of this statistic, such as $\sup_{u\in\mathcal{U}}\sqrt n |u\ttop(\hat\Sigma-\Sigma) u|/u\ttop \Sigma u$, where $\mathcal U$ is a set of sparse vectors in the unit sphere of $\R^p$.
In a different direction, the paper~\citep{LopesBA19} focuses on ``linear spectral statistics'' of the form $\sum_{j=1}^p f(\lambda_j(\hat\Sigma)) / p$, where $f:[0,\infty)\to \R$ is a smooth function. In that paper, it is shown that a type of parametric bootstrap procedure consistently approximates the distributions of such statistics when $p/n$ converges to a positive limit. Lastly, the paper~\cite{NaumovSU19} deals with bootstrapping statistics related the eigenvectors of $\hat\Sigma$. \\[-0.2cm]

\noindent\textbf{Notation.}
For a random variable $X$ and an integer $q\in\{1,2\}$, define the $\psi_q$-Orlicz norm as $\| X \|_{\psi_q} = \inf \: \{ t > 0 \ | \ \E [\exp( |X|^q / t^q)] \leq 2 \}$. The random variable $X$ is said to be sub-exponential if $\|X\|_{\psi_1}$ is finite, and sub-Gaussian if $\|X\|_{\psi_2}$ is finite. In addition, for any $q\geq 1$, the $L_q$ norm of $X$ is defined as $\| X \|_q = ( \E [ |X|^q ] )^{1/q}$.
For any vectors $u,v\in\R^p$, their inner product is $\langle u,v\rangle = \sum_{j=1}^p u_jv_j$.
For any real numbers $a$ and $b$, the expression $a\ll b$ is used in an informal sense to mean that $b$ is much larger than $a$. Also, we use the notation $a\vee b=\max\{a,b\}$ and $a\wedge b=\min\{a,b\}$. If $\{a_n\}$ and $\{b_n\}$ are two sequences on non-negative numbers, then the relation $a_n \lesssim b_n$ means that there is a positive constant $c$ not depending on $n$ such that $a_n \leq c \, b_n$ holds for all large $n$. When both of the conditions $a_n\lesssim b_n$ and $b_n\lesssim a_n$ hold, we write $a_n\asymp b_n$. 

\vspace{0.1cm}

\section{Main results}\label{sec:main}
We consider a sequence of models indexed by $n$, in which all parameters may depend on $n$ except when stated otherwise.  In particular, the dimension $p = p(n)$ is allowed to have arbitrary dependence on $n$. Likewise, if a parameter does not depend on $n$, then it is understood not to depend on $p$ either. One of the few parameters that will be treated as fixed with respect to $n$ is the positive integer $k<p$.

\begin{assumption}[Data-generating model]\label{model_assumptions}  \noindent \\[-0.4cm]
\begin{enumerate} 
\item[\namedlabel{A1}{(a)}.] 
There is a non-zero positive semidefinite matrix $\Sigma\in\R^{p\times p}$, such that the $i$th observation is generated as  $X_i = \Sigma^{1/2} Z_i$ for all $i = 1, \dots, n,$ where $Z_1,\dots,Z_n\in\R^p$ are i.i.d.~random vectors with $\E[Z_1]=0$, and $\E[Z_1Z_1\ttop]=I_p$. \\[-0.2cm]

\item[\namedlabel{A2}{(b)}.] The eigenvalues of $\Sigma$ satisfy
$\displaystyle\min_{1\leq j\leq k} \ts\big(\lambda_j(\Sigma) - \lambda_{j+1}(\Sigma)\big) \gtrsim  \lambda_1(\Sigma) $. \\[-0.2cm]

\item[\namedlabel{A3}{(c)}.] 
Let $u_j\in\R^p$ denote the $j$th eigenvector of $\Sigma$, and let $\Gamma \in\R^{k \times k}$ have entries given by $\Gamma_{jj'}=\E [(\langle u_j, Z_1 \rangle^2 - 1)(\langle u_{j'}, Z_1 \rangle^2 - 1)]$ for all $1\leq j,j'\leq k$. Then, the matrix $\Gamma$ satisfies $\lambda_k(\Gamma)\gtrsim 1$.
\end{enumerate}
\end{assumption}
\noindent In connection with the model described by  Assumption~\ref{model_assumptions}, our results will make reference to a moment parameter defined as
$\beta_q = \max_{1 \leq j \leq p} \| \langle u_j, Z_1 \rangle^2 \|_q$
for any $q\geq 1$. \\

\noindent \textbf{Remarks.} Regarding Assumption~\ref{model_assumptions}.\ref{A2}, it ensures that there is some degree of separation between the leading eigenvalues of $\Sigma$.  In less compact notation, the assumption states that there is a fixed constant $c>0$ such that the inequality \smash{$\lambda_j(\Sigma)-\lambda_{j+1}(\Sigma)\geq c \lambda_1(\Sigma)$} holds for all $j=1,\dots,k$, and all large $n$. (There is no restriction on the size of $c$.) In general, a separation condition on the leading eigenvalues is unavoidable, because it is known both theoretically and empirically that the bootstrap can fail to approximate the distribution of $\sqrt{n} \big( \blambda_k (\hat{\Sigma}) - \blambda_k (\Sigma) \big)$ if the leading population eigenvalues are not distinct~\citep{Beran:correction,Hall:2009}. In more technical terms, the source of this issue can be explained briefly as follows: If $\mathcal{S}^{p\times p}$ denotes the space of  real symmetric $p\times p$ matrices, and if $\lambda_j(\cdot)$ is viewed as a functional from $\mathcal{S}^{p\times p}$ to $\R$, then $\lambda_j(\cdot)$ becomes non-differentiable at $\Sigma$ in the case when  $\lambda_j(\Sigma)$ is a repeated eigenvalue (i.e.~with multiplicity larger than 1). In turn, this lack of smoothness makes it difficult for the bootstrap to approximate the distribution of $\sqrt n(\lambda_j(\hat\Sigma)-\lambda_j(\Sigma))$.

To interpret Assumption~\ref{model_assumptions}.\ref{A3}, the matrix $\Gamma$ serves a technical role as a surrogate for the correlation matrix of $\sqrt{n} \big( \blambda_k (\hat{\Sigma}) - \blambda_k (\Sigma) \big)$. Hence, the lower bound $\lambda_k(\Gamma)\gtrsim 1$ can be viewed as a type of non-degeneracy condition for the distribution of interest. The  proposition below gives examples of well-established models in which Assumption \ref{model_assumptions}.\ref{A3} holds. Namely, parts $(i)$ and $(ii)$ below respectively correspond to \emph{Mar\v{c}enko-Pastur models} and \emph{elliptical models}. The latter case also illustrates that the entries of the vector $Z_1$ are not required to be independent.

\begin{proposition} \label{PROP1} 
\noindent \namedlabel{PROP1_i}{(i)} (Mar\v{c}enko-Pastur case). Suppose that Assumption~\ref{model_assumptions}.\ref{A1} holds. In addition, suppose that the entries of $Z_1$ are independent, and there is a constant $\kappa>1$ not depending on $n$ such that $\min_{1\leq j\leq p }\E[Z_{1j}^4]\geq \kappa$. Then, Assumption \ref{model_assumptions}.\ref{A3} holds.

\noindent \namedlabel{PROP1_ii}{(ii)}  (Elliptical case). Let $V$ \!be a random vector that is uniformly distributed on the unit sphere of $\R^p$,
and let $\xi $ be a non-negative scalar random variable independent of $V$ that satisfies $\E[\xi ^2]=p$ and $\E[\xi^4]<\infty$. Under these conditions, if $Z_1$ has the same distribution as $\xi V$, then
Assumption \ref{model_assumptions}.\ref{A3} holds.
\end{proposition}

\noindent The proof of Proposition~\ref{PROP1} is given in Section~\ref{appendix_G} of the supplementary material. \\

\noindent \textbf{Bootstrap approximation.} 
The following theorem is the central result of the paper, and quantifies the accuracy of the bootstrap when it is used to approximate the distribution of $\sqrt{n}( \blambda_k (\hat{\Sigma})- \blambda_k (\Sigma) )$.

\begin{theorem} \label{THM:ALL}
Suppose that Assumption \ref{model_assumptions} holds and let $q = 5\log(kn)$. Then, there is a constant $c > 0$ not depending on $n$ such that the event
\scriptsize
\begin{equation}\label{eqn:thm:all} 
 \sup_{t \in \R^k } \Bigg| \P \Big( \sqrt{n} \big( \blambda_k (\hat{\Sigma}) - \blambda_k (\Sigma) \big) \mylessthan \Big) - \P \Big( \sqrt{n} \big( \blambda_k (\hat{\Sigma}^\star) - \blambda_k (\hat\Sigma) \big) \mylessthan \, \Big| \, X \Big) \Bigg| \
%
\leq \ \frac{c \, \log(n) \,  \beta_{3q}^3\,{\tt{r}}(\Sigma) }{\sqrt n} 
\end{equation}
\normalsize
holds with probability at least $1 - c/n$.
\end{theorem}

\noindent \textbf{Remarks.} The proof of Theorem~\ref{THM:ALL} is given in Section~\ref{appendix_A} of the supplemenatary material. It is possible to provide a more concrete understanding of the bound~\eqref{eqn:thm:all} by looking at how the factors ${\tt{r}}(\Sigma)$ and $\beta_{3q}$ behave in some well-known situations. For instance, consider the class of matrices $\Sigma$ whose eigenvalues have a polynomial decay profile of the form $\lambda_j(\Sigma)\asymp j^{-\gamma}$, for some fixed constant $\gamma>0$. This class offers a convenient point of reference, because it interpolates between models that have low-dimensional structure and those that do not. Specifically, the effective rank can be related to $\gamma$ as 
$$ {\tt{r}}(\Sigma) \asymp \begin{cases} & \!\!1 \text{ \ \ \  \ \ \ \  \ \   if } \gamma >1  \\
&\!\!\log(p) \  \text{ \ \  if } \gamma = 1  \\ 
& \!\!  p^{1-\gamma} \  \text{ \ \  \ \   if } \gamma<1.
\end{cases}
$$

With regard to the parameter $\beta_{3q}$, its dependence on $q$ is simple to describe in some commonly considered cases.
If the entries of $Z_1$ are i.i.d.~and sub-Gaussian, then $\beta_{3q}$ grows at most linearly in $q$, with $\beta_{3q}\lesssim q \|Z_{11}\|_{\psi_2}^2$. Alternatively, if the entries of $Z_1$ are i.i.d.~and sub-exponential, then $\beta_{3q}$ grows at most quadratically in $q$, with $\beta_{3q}\lesssim q^2 \|Z_{11}\|_{\psi_1}^2$. (See Chapter 2 of~\cite{Vershynin:2018} for further details.) Hence, a direct consequence of Theorem~\ref{THM:ALL} in such cases is that bootstrap consistency holds when $\gamma>1/2$, $p\asymp n$ and $\|Z_{11}\|_{\psi_1}\lesssim 1$. Likewise, when $\gamma>1$, the bound in Theorem~\ref{THM:ALL} nearly achieves the \emph{parametric rate} $n^{-1/2}$ and is not influenced by the size of $p$ at all. This conclusion also conforms with the numerical results that we present in Section~\ref{sec:experiments}. 

From a more practical standpoint, it is possible to gauge the size of ${\tt{r}}(\Sigma)$ in an empirical way, by either estimating ${\tt{r}}(\Sigma)$ directly, or estimating upper bounds on it. Some examples of upper bounds on ${\tt{r}}(\Sigma)$ for which straightforward estimation methods are known to be effective in high dimensions include $\tr(\Sigma)/\max_{1\leq j\leq p}\Sigma_{jj}$ and $\tr(\Sigma)^2/\|\Sigma\|_F^2$. (Although guarantees can be established for direct estimates of ${\tt{r}}(\Sigma)$ in high-dimensions, such results can involve a more complex set of considerations than the upper bounds just mentioned.) \\

\noindent \textbf{Transformations.} To briefly review the idea of transformations, they are often used to solve inference problems involving a parameter $\theta$ and an estimator $\hat\theta$ for which the distribution of $(\hat\theta-\theta)$ is difficult to approximate. In certain situations, this difficulty can be alleviated if there is a monotone function $h$ for which the distribution of $(h(\hat\theta)-h(\theta))$ is easier to approximate. In turn, this allows for more accurate inference on the ``transformed parameter'' $h(\theta)$, and then the results can be inverted to do inference on $\theta$. 
In light of this, our next result shows that the rates of bootstrap approximation established in Theorem~\ref{THM:ALL} remain essentially unchanged when using the class of fractional power transformations from $[0,\infty)$ to $[0,\infty)$. This class will be denoted by $\mathcal{H}$, so that if $h\in\mathcal{H}$, then $h(x)=x^{a}$ for some $a\in (0,1]$.

Beyond the class of transformations just mentioned, the bootstrap can be combined with another type of transformation known as \emph{partial standardization}~\citep{LopesLM18}. Letting $h\in\mathcal{H}$ be a given function, and letting  $\varsigma_j^2 = \var \big( h( \lambda_j (\hat\Sigma)) \big)$ for each $j=1,\dots,p$, this technique is well suited to bootstrapping ``max statistics'' of the form
\begin{equation}\label{eqn:maxstat}
M = \max_{1\leq j\leq k}  \frac{ h(\lambda_j (\hat{\Sigma})) - h(\lambda_j (\Sigma)) }{\varsigma_j^\tau},
\end{equation}
where $\tau\in[0,1]$ is a parameter that can be viewed as a degree of standardization. The ability to approximate the distribution of $M$ is relevant to the construction of simultaneous confidence intervals for $\lambda_1(\Sigma),\dots,\lambda_k(\Sigma)$. It also turns out that the choice of $\tau$ encodes a trade-off between the coverage accuracy and the width of such intervals, and that choosing an intermediate value $\tau\in (0,1)$ can offer benefits in relation to $\tau=0$ and $\tau=1$. This will be discussed in greater detail later in Section~\ref{sec:experiments}.

In order to state our extension of Theorem~\ref{THM:ALL} in a way that handles both partial standardization and transformations $h\in\mathcal{H}$ in a unified way, we need to introduce a bit more notation. First, when considering the bootstrap counterpart of a partially standardized statistic such as~\eqref{eqn:maxstat}, the vector $\boldsymbol{\varsigma}_k^{\tau} = ( \varsigma_1^{\tau}, \dots, \varsigma_k^{\tau})$ can be replaced with the estimate $\boldsymbol{\hat\varsigma}_k^{\tau} = ( \hat\varsigma_1^{\tau}, \dots, \hat\varsigma_k^{\tau})$, whose entries are defined by $\hat\varsigma_j^2 = \var \big( h(\lambda_j (\hat\Sigma^\star) )\big| X \big)$ for all $j=1,\dots,p$. Second,
the expression $v/u$ involving vectors $v$ and $u$ denotes the vector obtained by entrywise division,  $(v/u)_j=v_j/u_j$. (To handle the possibility zero denominators, events of the form $\{V/\hat{\boldsymbol \varsigma}_k^{\tau}\preceq t\}$ are understood as $\{V\preceq t\odot \hat{\boldsymbol \varsigma}_k^{\tau}\}$, where $V\in\R^k$ is random, $t\in\R^k$ is fixed, and $\odot$ is entrywise multiplication. Lemma~\ref{lemma:trans_var_lower_bound} in the supplementary material also shows that such cases occur with negligible probability.) Lastly, recall that we write $\boldsymbol h (v) = (h(v_1), \dots, h(v_k))$ for a $k$-dimensional vector $v$ and \smash{transformation $h$.}

\begin{theorem} \label{THM:TRANS_GENERAL}
Suppose that Assumption \ref{model_assumptions} holds. Fix a transformation $h\in\mathcal{H}$ and a constant $\tau\in[0,1]$ with respect to $n$, and let $q = 5\log(kn)$. Then, there is a constant $c > 0$ not depending on $n$, such that the event
\scriptsize
\begin{align*}
\sup_{t \in \R^k } \bigg| \P \Bigg( \frac{\boldsymbol h (\blambda_k(\hat\Sigma)) -\boldsymbol h( \blambda_k(\Sigma))}{\boldsymbol\varsigma_k^{\tau}} \mylessthan \bigg) %
- \P \bigg( \frac{\boldsymbol h (\blambda_k(\hat\Sigma^{\star})) -\boldsymbol h( \blambda_k(\hat\Sigma))}{\boldsymbol{\hat{\varsigma}}_k^{\tau}} \mylessthan\,\bigg|X \bigg) \Bigg|  
& \ \leq \
\frac{c\, \log(n)\, \beta_{3q}^{5}\, {\tt{r}}(\Sigma) }{\sqrt n} 
\end{align*}
\normalsize
holds with probability at least $1 - c/n$.
\end{theorem}
\noindent \textbf{Remarks.}  The proof of Theorem~\ref{THM:TRANS_GENERAL} is given in Section \ref{appendix_D} of the supplementary material. To comment on the technical relationship between Theorems~\ref{THM:ALL} and \ref{THM:TRANS_GENERAL}, it is important to call attention to the differences between asymptotic and non-asymptotic analysis. When using asymptotics, the process of showing that bootstrap consistency for $\sqrt{n} \big( \blambda_k (\hat{\Sigma}) - \blambda_k (\Sigma) \big)$ implies the same for $(\boldsymbol h (\blambda_k(\hat\Sigma)) -\boldsymbol h( \blambda_k(\Sigma)))/\boldsymbol\varsigma_k^{\tau} $ can typically be handled with a brief argument, based on the delta method and the consistency of the estimate $\boldsymbol{\hat{\varsigma}}_k^{\tau}$. However, when taking a non-asymptotic approach, this process is much more involved. For instance, it is necessary to establish fine-grained error bounds for $\boldsymbol{\hat{\varsigma}}_k^{2}$, such as in showing that the uniform relative error $\|(\boldsymbol{\hat{\varsigma}}_k^{2}-\boldsymbol{\varsigma}_k^{2})/\boldsymbol\varsigma_k^2\|_{\infty}$ is likely to be at most of order $n^{-1/2}\beta_{2q}^3\lambda_1(\Sigma)^2{\tt{r}}(\Sigma)$ up to logarithmic factors.

\vspace{0.2cm}
	
\section{Numerical Experiments}\label{sec:experiments}
	
In this section, we focus on the application of constructing simultaneous confidence intervals for $\lambda_1(\Sigma),\dots,\lambda_k(\Sigma)$. This will be done in a variety of settings, corresponding to different values of $n$ and $p$, as well as different values of effective rank, and different choices of transformations. In a nutshell, there are two overarching conclusions to take away from the experiments: {\textbf{(1)}} In situations where $n\ll p$ and ${\tt{r}}(\Sigma)\asymp 1$, the bootstrap generally produces intervals with accurate coverage, which provides a confirmation of our theoretical results. {\textbf{(2)}} The classical log transformation mostly works well in low dimensions, but it can lead to coverage that is substantially below the nominal level when ${\tt{r}}(\Sigma)$ is moderately large. Nevertheless, we show that it is possible to find transformations that offer more reliable coverage in this challenging case.
More generally, this indicates that alternative transformations are worth exploring in high-dimensional settings. \\[-0.3cm]

\subsection{Simulation settings}
The eigenvalues of the population covariance matrix $\Sigma$ were chosen to have two different decay profiles: 
\begin{enumerate}
\item[(a)] A polynomial decay profile $\lambda_j(\Sigma) = j^{-\gamma}$ for all $j = 1, \dots, p$, with $\gamma \in \{ 0.7, 1.0, 1.3 \}$. \\[-0.2cm]
\item[(b)] An exponential decay profile $\lambda_j(\Sigma) = \delta^j$ for all $j = 1, \dots, p$, with $\delta \in \{ 0.7, 0.8, 0.9 \}$. 
\end{enumerate}
As a clarification, it is important to note that the effective rank of $\Sigma$ increases for larger values of $\delta$, but decreases for larger values of $\gamma$. For the purposes of simulations, the choices (a) and (b) have the valuable property that the eigenvalues are parameterized in the same way for every choice of $p$, which facilitates the comparison of results across different dimensions. The matrix of eigenvectors for $\Sigma$ was drawn uniformly from the set of $p\times p$ orthogonal matrices.
The dimension $p$ was taken from $\{10, 50, 100, 200\}$, and the sample size $n$ ranged from 50 to 500.  For each triple $(n, p, \gamma)$ or $(n,p,\delta)$, the data $X_1,\dots,X_n$ were generated in an i.i.d.~manner with the following choices for the distribution of $X_1$: 
\begin{enumerate} 
\item[(i)] The vector $X_1=\Sigma^{1/2}\xi V$ was generated with $V$ being uniformly distributed on the unit sphere of $\R^p$, and $\xi^2$ being an exponential random variable independent of $V$ with $\E[\xi^2]=p$. \\[-0.2cm]
\item[(ii)] The vector $X_1$ was generated from the Gaussian distribution $N(0,\Sigma)$. 
\end{enumerate}
For each parameter setting, we generated 1000 realizations of the dataset $X_1,\dots,X_n$, and for each such realization, we generated $B:=1000$ sets of bootstrap samples of size $n$.
When constructing simultaneous confidence intervals for $\lambda_1(\Sigma),\dots,\lambda_k(\Sigma)$, the value of $k$ was set to $5$.

\subsection{Bootstrap confidence intervals}\label{sec:bootci}

For any $\alpha\in (0,1)$, we aim to construct approximate versions of ideal random intervals $\mathcal{I}_1,\dots,\mathcal{I}_k$ that satisfy
\begin{equation}\label{eqn:desiredCI}
\P \Bigg( \bigcap_{j = 1}^k \big\{ \lambda_j (\Sigma) \in \I_j \big\} \Bigg) \ \geq \ 1 - \alpha.
\end{equation}

\noindent
To this end, consider the following max and min statistics, based on any choice of partial standardization parameter $\tau\in[0,1]$ and transformation $h$,
\begin{align*}
M & \ =\ \max_{1 \leq j \leq k}  \frac{ h(\lambda_j (\hat{\Sigma})) - h(\lambda_j (\Sigma)) }{\varsigma_j^\tau}\\
L & \ =\ \min_{1\leq j\leq k} \frac{h(\lambda_j (\hat{\Sigma})) - h(\lambda_j (\Sigma)) }{\varsigma_j^\tau}.
\end{align*}
Letting $q_{M}(\alpha)$ and $q_L(\alpha)$ denote the respective $\alpha$-quantiles of $M$ and $L$ for any $\alpha\in (0,1)$, it follows that the desired condition~\eqref{eqn:desiredCI} holds if each interval  $\mathcal{I}_j$ is defined as
\begin{equation}\label{eqn:idealintervals}
\mathcal{I}_j =
h^{-1}\bigg(\bigg[ h(\lambda_j(\hat\Sigma))-\varsigma_j^{\tau}q_M(1-\ts\frac{\alpha}{2}) \ , \ h(\lambda_j(\hat\Sigma))-\varsigma_j^{\tau}q_L(\ts\frac{\alpha}{2}) \bigg]\bigg),
\end{equation}
with $h^{-1}([a,b])$ being understood as the preimage of $[a,b]$ under $h$.

To construct bootstrap intervals $\hat{\mathcal{I}}_1,\dots,\hat{\mathcal{I}}_k$ based on~\eqref{eqn:idealintervals}, it is only necessary to replace $q_M(1-\frac{\alpha}{2})$, $q_{L}(\frac{\alpha}{2})$, and $\varsigma_1,\dots,\varsigma_k$ with estimates. 
In detail, each $\varsigma_j$ is first estimated using the sample variance of  $B$ bootstrap replicates of the form $h(\lambda_j(\hat\Sigma^{\star}))$. Next, the empirical $1-\frac{\alpha}{2}$ quantile of $B$ bootstrap replicates of the form $M^{\star}=\max_{1 \leq j \leq k}  [h(\lambda_j (\hat{\Sigma}^{\star})) - h(\lambda_j (\hat\Sigma))]/\hat\varsigma_j^{\tau} $ is taken as an estimate of $q_M(1-\frac{\alpha}{2})$, and similarly for $q_L(\frac{\alpha}{2})$.

Regarding the use of transformations, the following three options were included in the experiments:
\begin{itemize}
\item \emph{log transformation}: $h(x)=\log(x)$ with $\tau=0$.\\[-0.2cm]
\item \emph{standardization}: $h(x)=x$ with $\tau=1$.\\[-0.2cm]
\item \emph{square-root transformation}: $h(x)=x^{1/2}$ with $\tau\in[0,1]$ chosen data-adaptively. 
\end{itemize}
In the case of the log transformation, the choice of $\tau=0$ corresponds to the way that this transformation has been used in the classical literature~\citep{BeranS85}, while in the case of standardization, the choice of $\tau=1$ is definitional. For the square-root transformation, the use of a data-adaptive selection rule for $\tau\in[0,1]$ is more nuanced, and can be informally explained in terms of the following ideas developed previously in~\citep{LopesLM18, LinLM2020}.

In essence, this choice can be understood in terms of a trade-off between two competing effects that occur in the extreme cases of $\tau=1$ and $\tau =0$.  
When using $\tau=1$, the random variables $[\lambda_j(\hat{\Sigma})^{1/2} - \lambda_j (\Sigma)^{1/2}]/\varsigma_j$ with $\varsigma_j^2=\var(\lambda_j(\hat\Sigma)^{1/2})$ and $j=1,\dots,k$ are on approximately ``equal footing'', which makes the behavior of the statistic $M$ sensitive to their joint distribution (and likewise for $L$).
By contrast, when $\tau=0$ is used, the variables $[\lambda_j(\hat{\Sigma})^{1/2} - \lambda_j (\Sigma)^{1/2}]$ will tend to be on different scales, and the variable on the largest scale, say $j^{\prime}$, will be the maximizer for $M$ relatively often. In this situation, the statistic $M$ is governed more strongly by the marginal distribution of $[\lambda_{j^{\prime}}( \hat{\Sigma})^{1/2} - \lambda_{j^{\prime}} (\Sigma)^{1/2}]$. So, from this heuristic point of view, the choice of $\tau=0$ can simplify the behavior of $M$ relative to the case of $\tau=1$, making the distribution of $M$ easier to approximate.
However, the choice of $\tau=0$ also has the drawback that it can lead to simultaneous confidence intervals that are excessively wide, because the widths are no longer adapted to the different values $\varsigma_1,\dots,\varsigma_k$ (since $\varsigma_1^0=\cdots=\varsigma_k^0=1$).

To strike a balance between these competing effects, we used the following simple rule to select $\tau$ in the case of the square-root transformation. For a candidate value of $\tau$, let $\hat{\mathcal{I}}_1(\tau),\dots,\hat{\mathcal{I}}_k(\tau)$ denote the associated bootstrap intervals defined beneath equation~\eqref{eqn:idealintervals} (so that the dependence on $\tau$ is explicit), and let $|\hat\I_1(\tau)|,\dots,|\hat\I_k(\tau)|$ denote their widths. Also define $\hat\mu(\tau) =  \sum_{j=1}^k |\hat\I_j(\tau)| / k$ and $\hat\sigma(\tau)^2 = \sum_{i = 1}^k (|\hat\I_j(\tau)| - \hat\mu(\tau) )^2 / k$. In this notation, we selected the value of $\tau$ that minimized $\hat\mu(\tau)+\hat\sigma(\tau)$ over the set of candidates $\{0.0, 0.1,\dots,0.9, 1.0\}$. Different variants of this type of criterion minimization rule have also been observed to be effective in other contexts~\citep{LopesLM18,LinLM2020}. \\[-0.3cm]

\begin{figure}[H]	
\vspace{0.5cm}
	\quad\quad\quad 
	\begin{overpic}[width=0.29\textwidth]{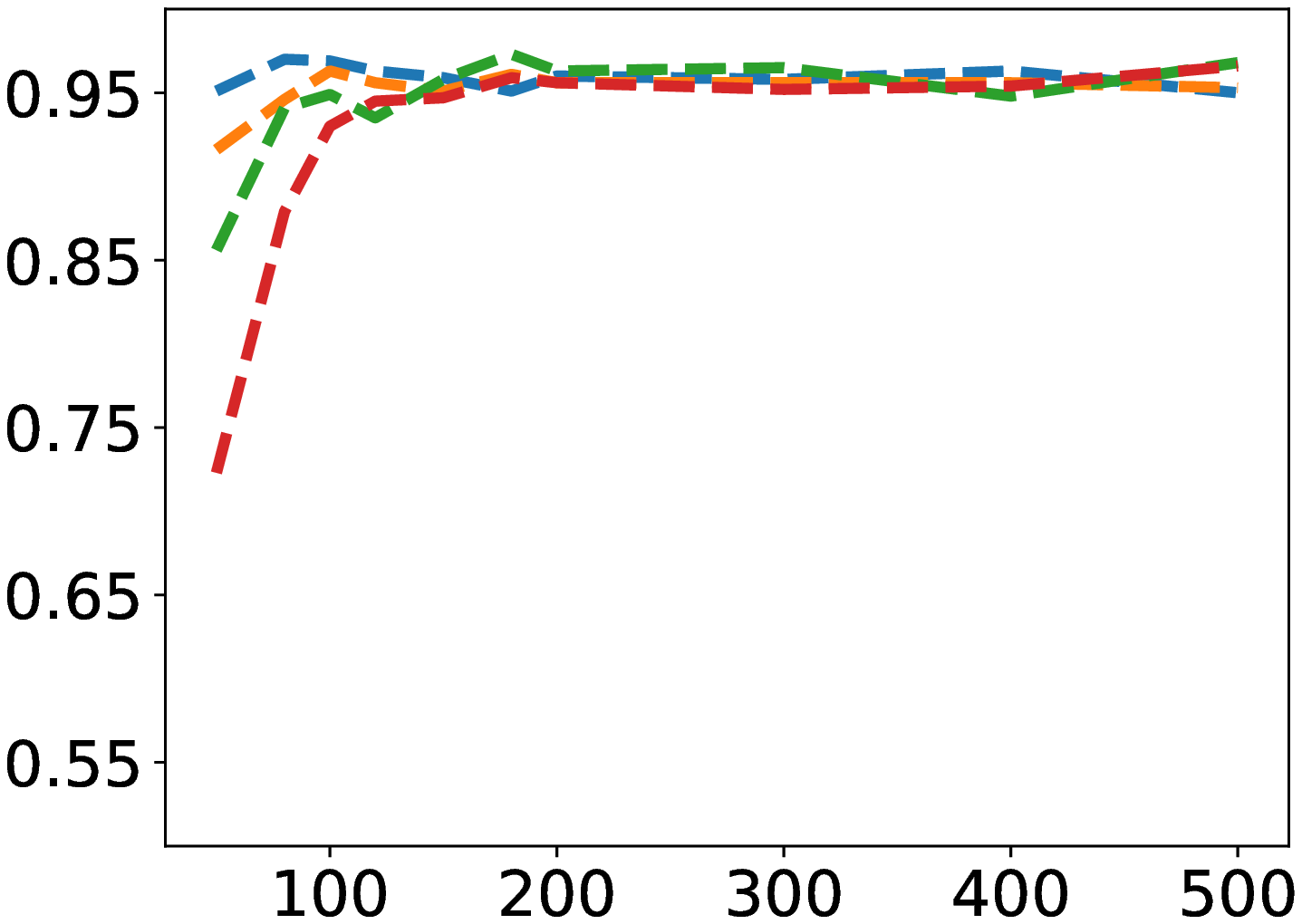} 
\put(25,80){ \ul{\ \ \  \ $\gamma=0.7$ \ \ \ \    }}
		\put(-20,-5){\rotatebox{90}{ {\small \ \ \ log transformation  \ \ }}}
\end{overpic}
	~
	\DeclareGraphicsExtensions{.png}
	\begin{overpic}[width=0.29\textwidth]{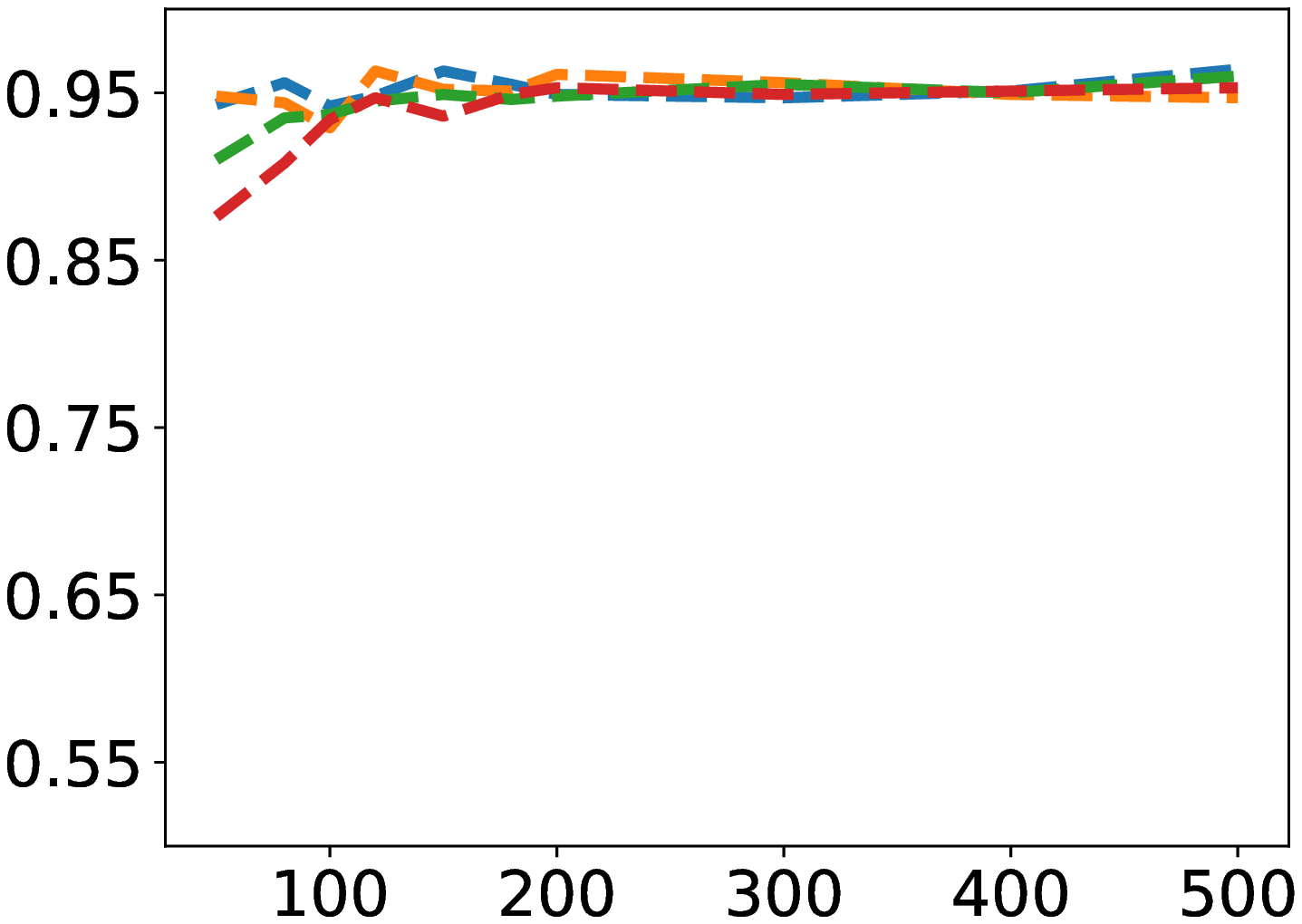} 
		\put(25,80){ \ul{\ \ \  \ $\gamma=1.0$ \ \ \ \    }}
	\end{overpic}
	~	
	\begin{overpic}[width=0.29\textwidth]{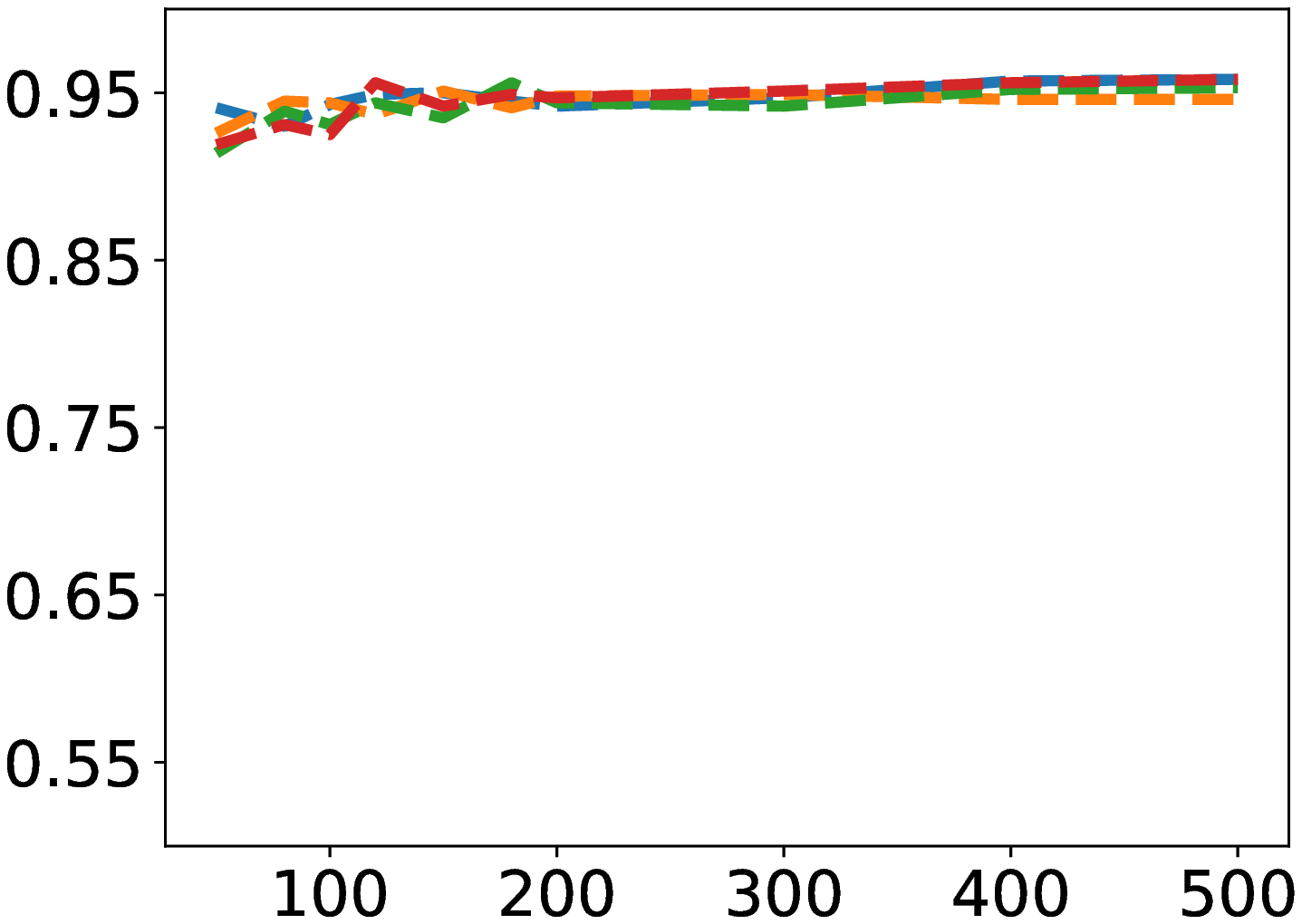} 
		\put(25,80){ \ul{\ \ \  \ $\gamma=1.3$ \ \ \ \    }}
	\end{overpic}	
\end{figure}

\vspace{-0.7cm}

\begin{figure}[H]	
	\quad\quad\quad 
	\begin{overpic}[width=0.29\textwidth]{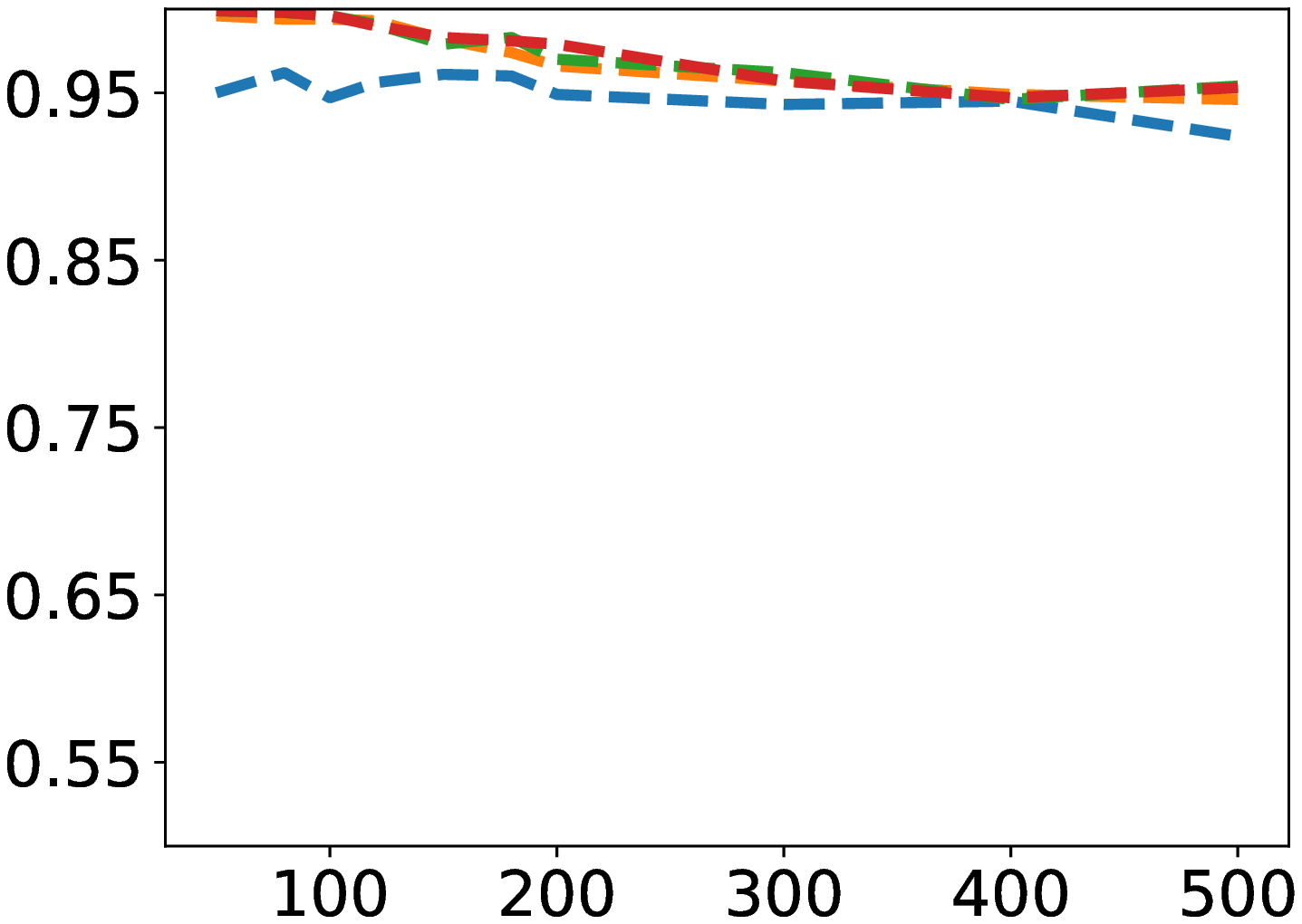} 
	\put(-20,-1){\rotatebox{90}{ {\small \ \ \ standardization \  \ \ }}}
	\end{overpic}
	~
	\DeclareGraphicsExtensions{.png}
	\begin{overpic}[width=0.29\textwidth]{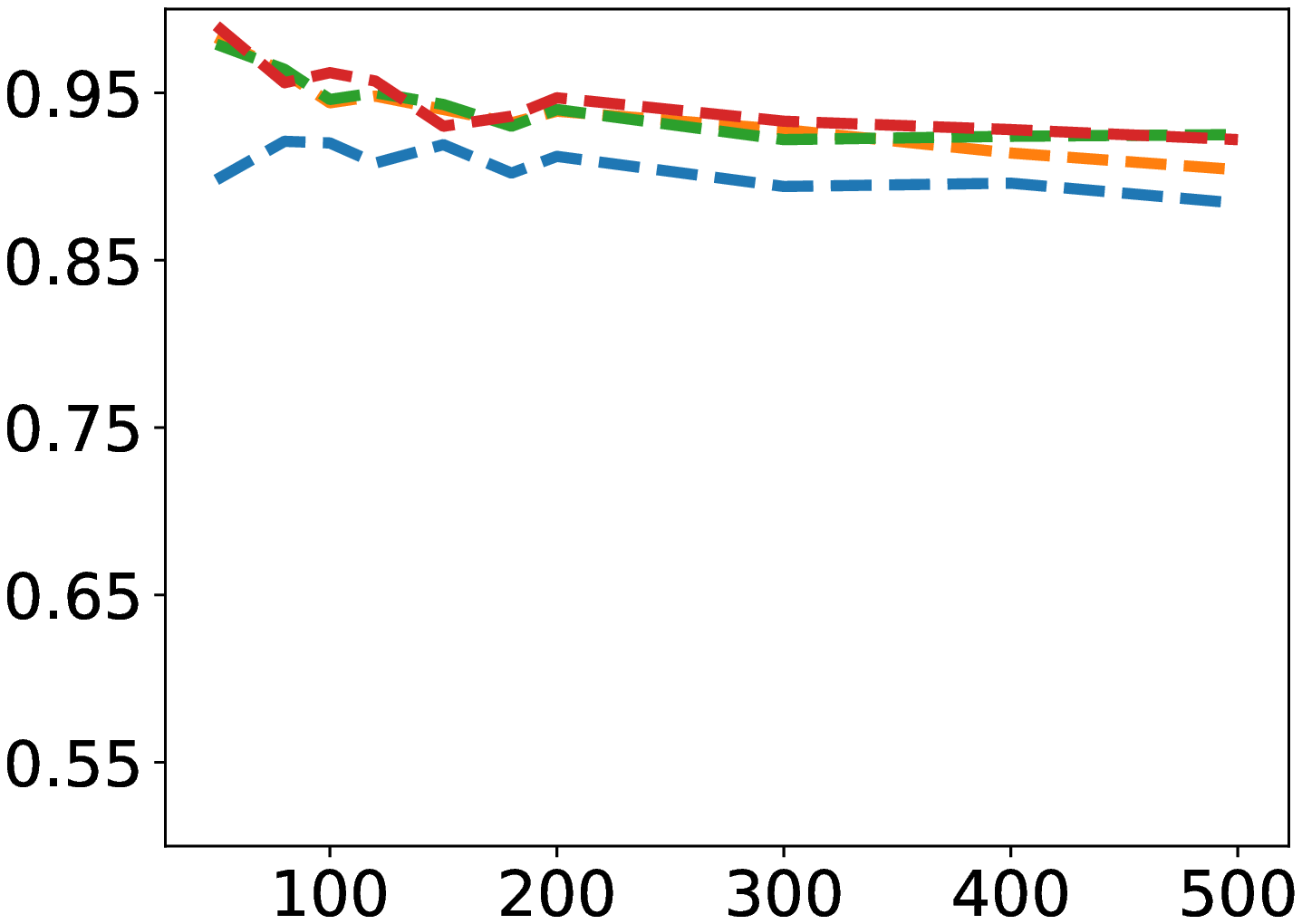} 
	\end{overpic}
	~	
	\begin{overpic}[width=0.29\textwidth]{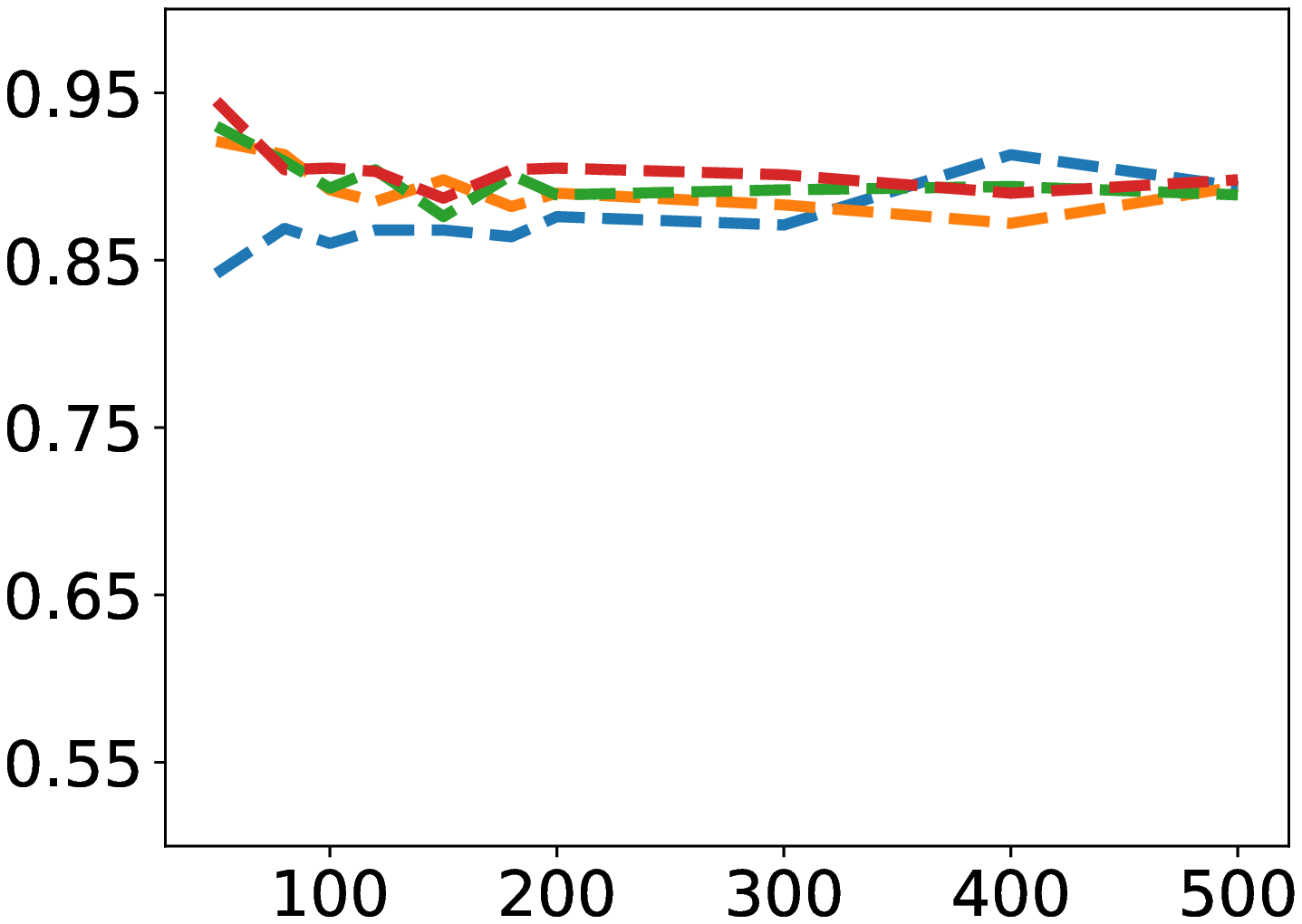} 
	\end{overpic}	
\end{figure}

\vspace{-0.7cm}

\begin{figure}[H]	
	\quad\quad\quad 
	\begin{overpic}[width=0.29\textwidth]{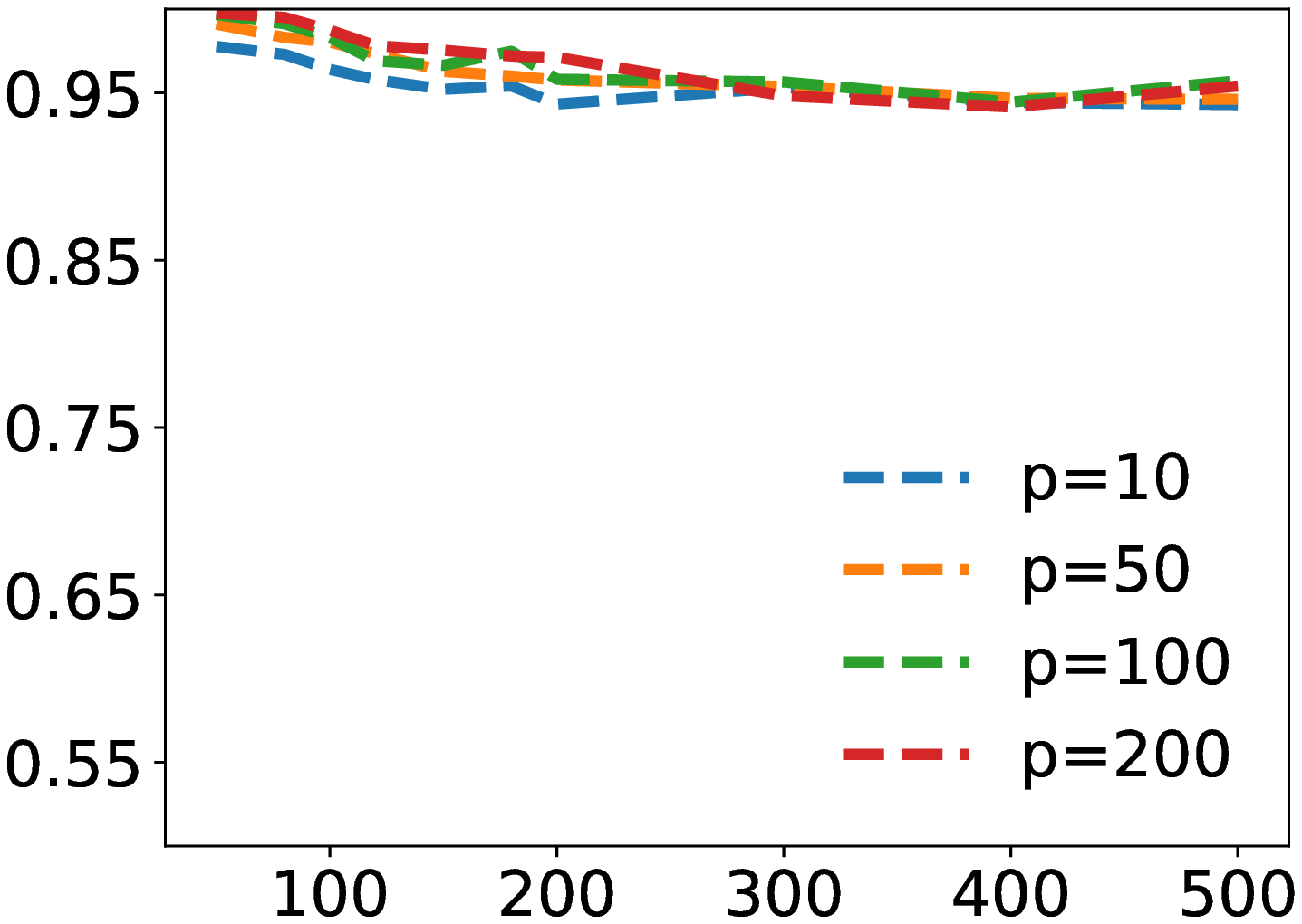} 
		\put(-21,1){\rotatebox{90}{\ $\sqrt{ \ \ }$}}
	    \put(-20,-3){\rotatebox{90}{  { \ \ \ \ \ \ \small transformation \ \ } }}

	\end{overpic}
	~
	\DeclareGraphicsExtensions{.png}
	\begin{overpic}[width=0.29\textwidth]{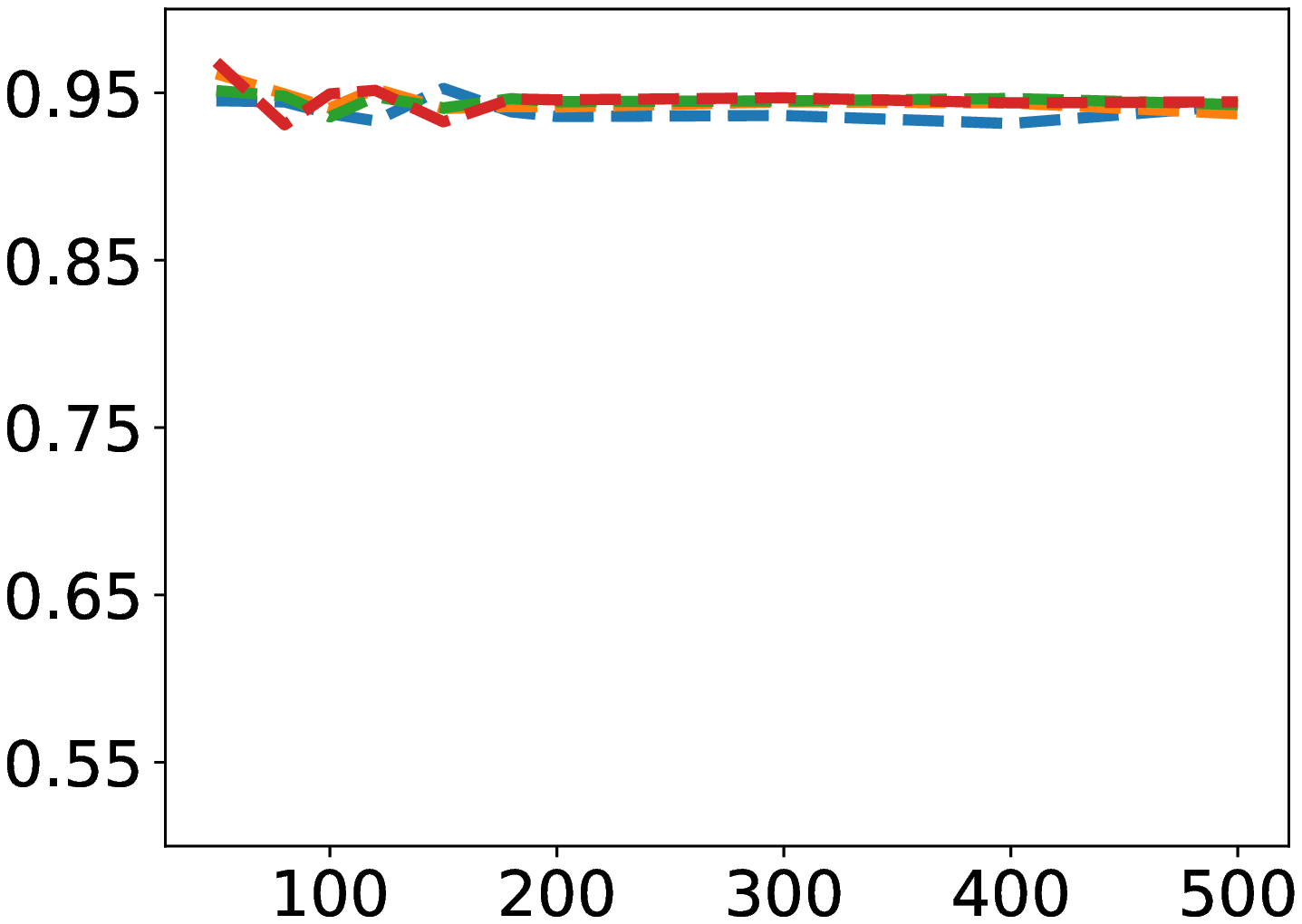} 
	\end{overpic}
	~	
	\begin{overpic}[width=0.29\textwidth]{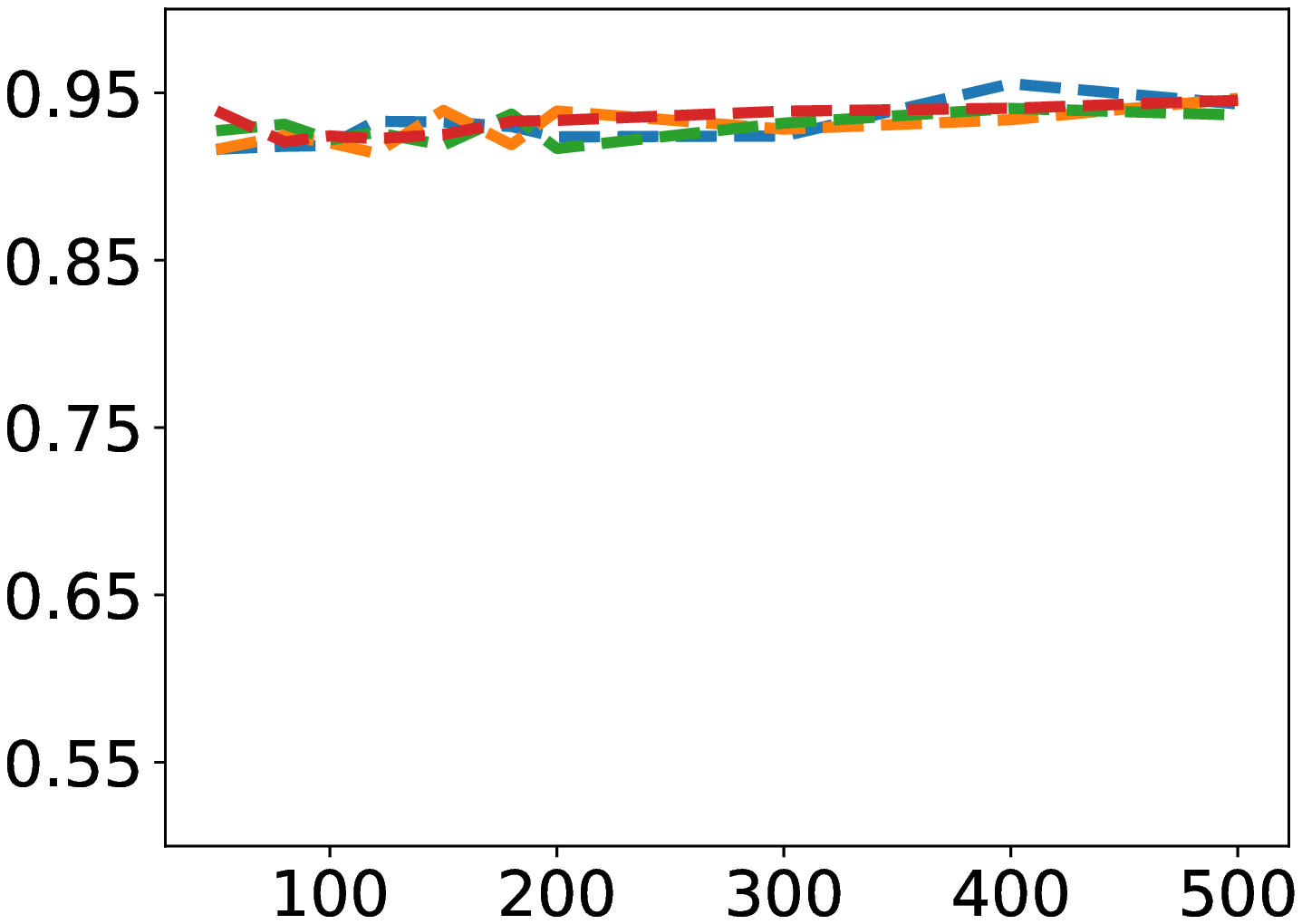} 
	\end{overpic}	
	\vspace{+.2cm}	
	\caption{(Simultaneous coverage probability versus $n$ in simulation model (i) with a polynomial decay profile). In each panel, the $y$-axis measures~$\P(\cap_{j=1}^5 \{\lambda_j(\Sigma)\in\hat{\mathcal{I}}_j\})$ based on a nominal value of 95\%, and the $x$-axis measures $n$. The colored curves correspond to the different values of $p$, indicated in the legend. }\label{fig1} 

\end{figure}

\subsection{Discussion of coverage}

Figure \ref{fig1} contains nine panels displaying the results for the simultaneous coverage probability $\P(\cap_{j=1}^5 \{\lambda_j(\Sigma)\in\hat{\mathcal{I}}_j\})$, based on a nominal value of 95\% (i.e.~$\alpha=0.05$) in the case of the simulation model (i) with a polynomial decay profile for the population eigenvalues. The figure summarizes a large amount of information, because it shows how the coverage depends on $n$, $p$, the eigenvalue decay parameter $\gamma$, and the three transformations described above. For each panel, the $x$-axis measures $n$, and the $y$-axis measures $\P(\cap_{j=1}^5 \{\lambda_j(\Sigma)\in\hat{\mathcal{I}}_j\})$. Results corresponding to the dimensions $p=10, 50, 100, 200$ are plotted with colored curves that are labeled in the legend. The three rows of panels from top to bottom correspond to the log transformation, ordinary standardization, and the square-root transformation. The three columns of panels from left to right correspond to the eigenvalue decay parameters $\gamma=0.7, 1.0, 1.3$.  In addition, Figure~\ref{fig3} displays analogous results for exponentially decaying population eigenvalues in model (i). Lastly, results for model (ii), as well as for a nominal value of 90\% (instead of 95\%), are provided in Section~\ref{appendix_last} of the supplementary material.

%



\begin{figure}[H]	
\vspace{0.5cm}
	\quad\quad\quad 
	\begin{overpic}[width=0.29\textwidth]{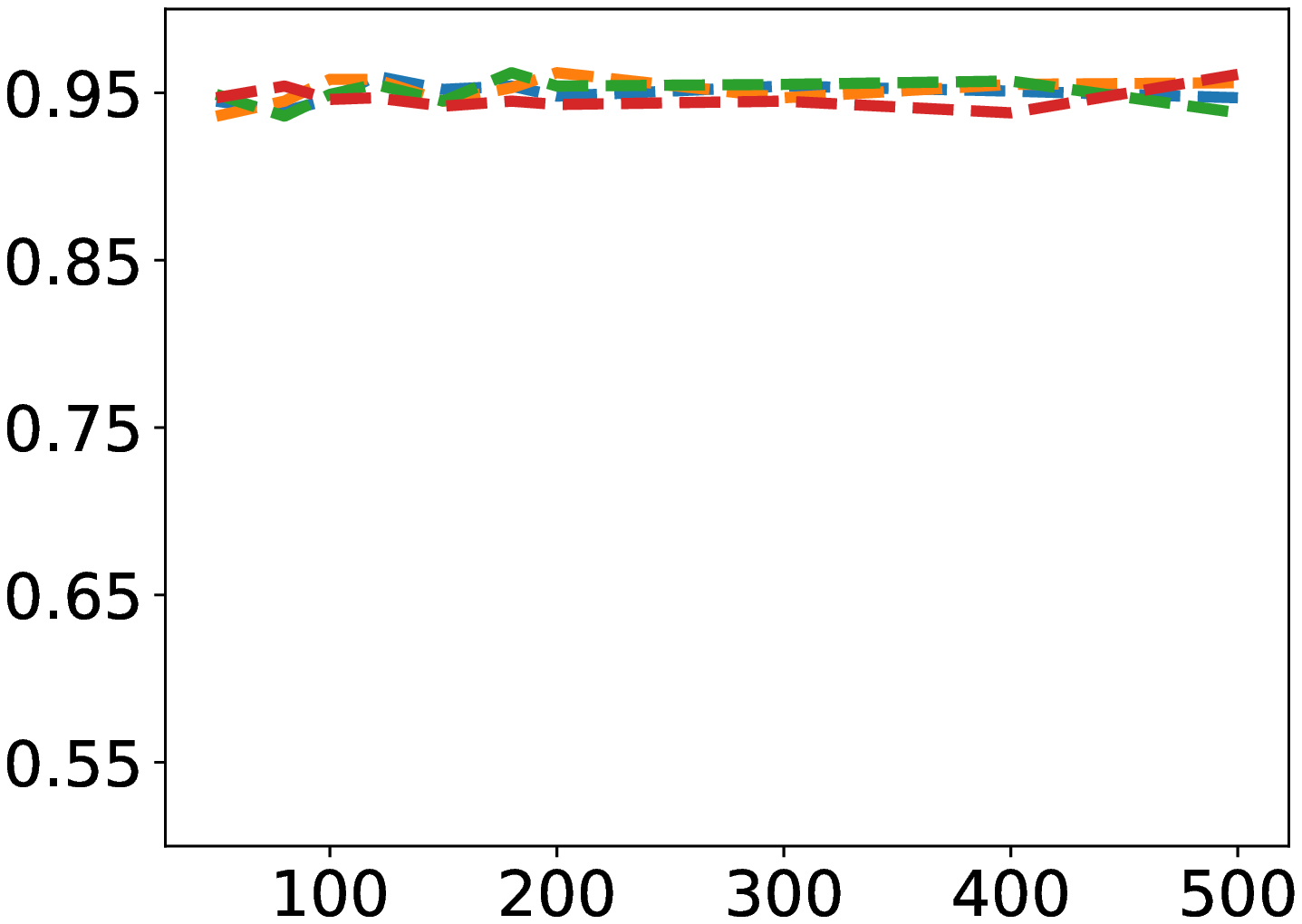} 
\put(25,80){ \ul{\ \ \  \ $\delta=0.7$ \ \ \ \    }}
		\put(-20,-5){\rotatebox{90}{ {\small \ \ \ log transformation  \ \ }}}
\end{overpic}
	~
	\DeclareGraphicsExtensions{.png}
	\begin{overpic}[width=0.29\textwidth]{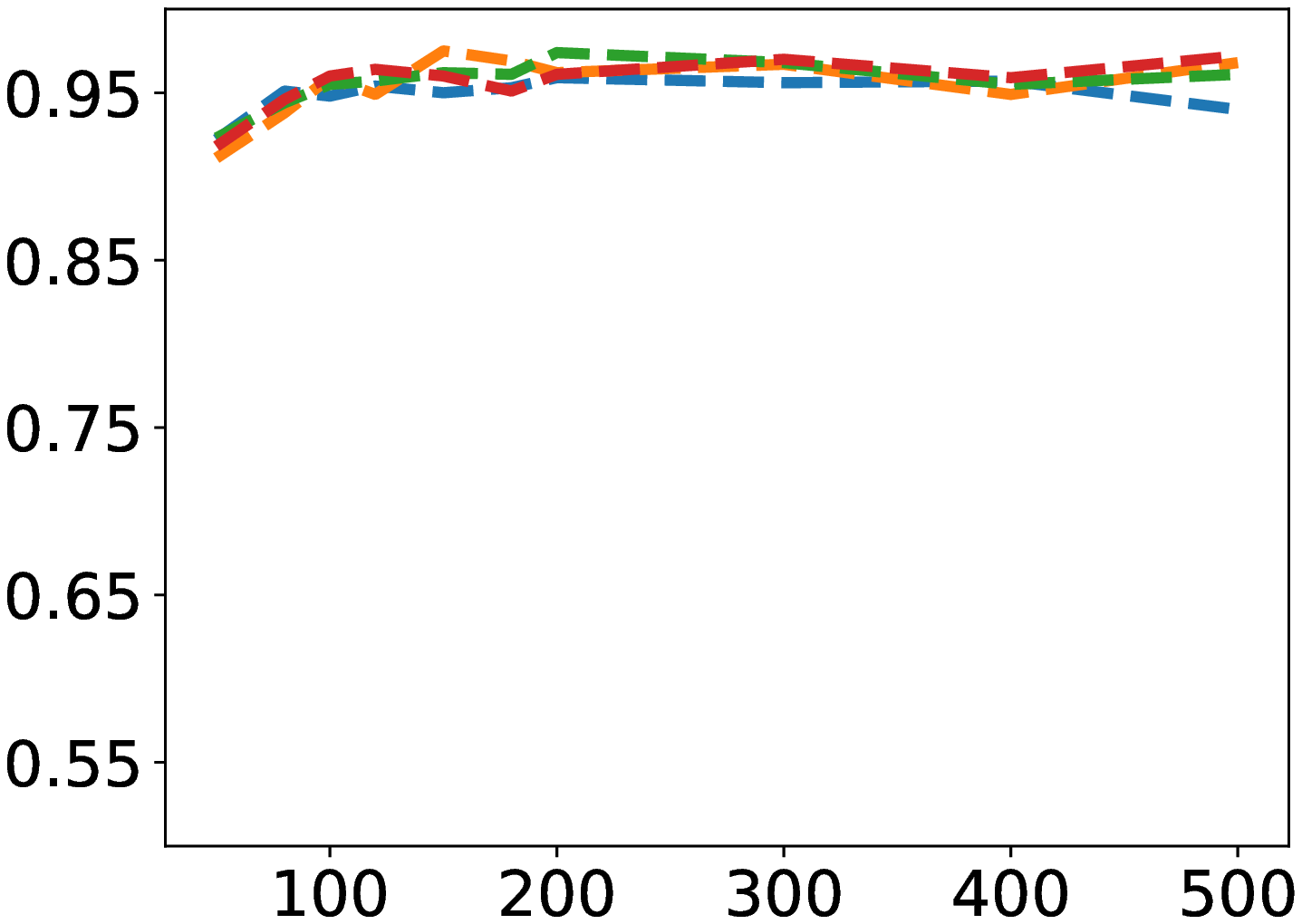} 
		\put(25,80){ \ul{\ \ \  \ $\delta=0.8$ \ \ \ \    }}
	\end{overpic}
	~	
	\begin{overpic}[width=0.29\textwidth]{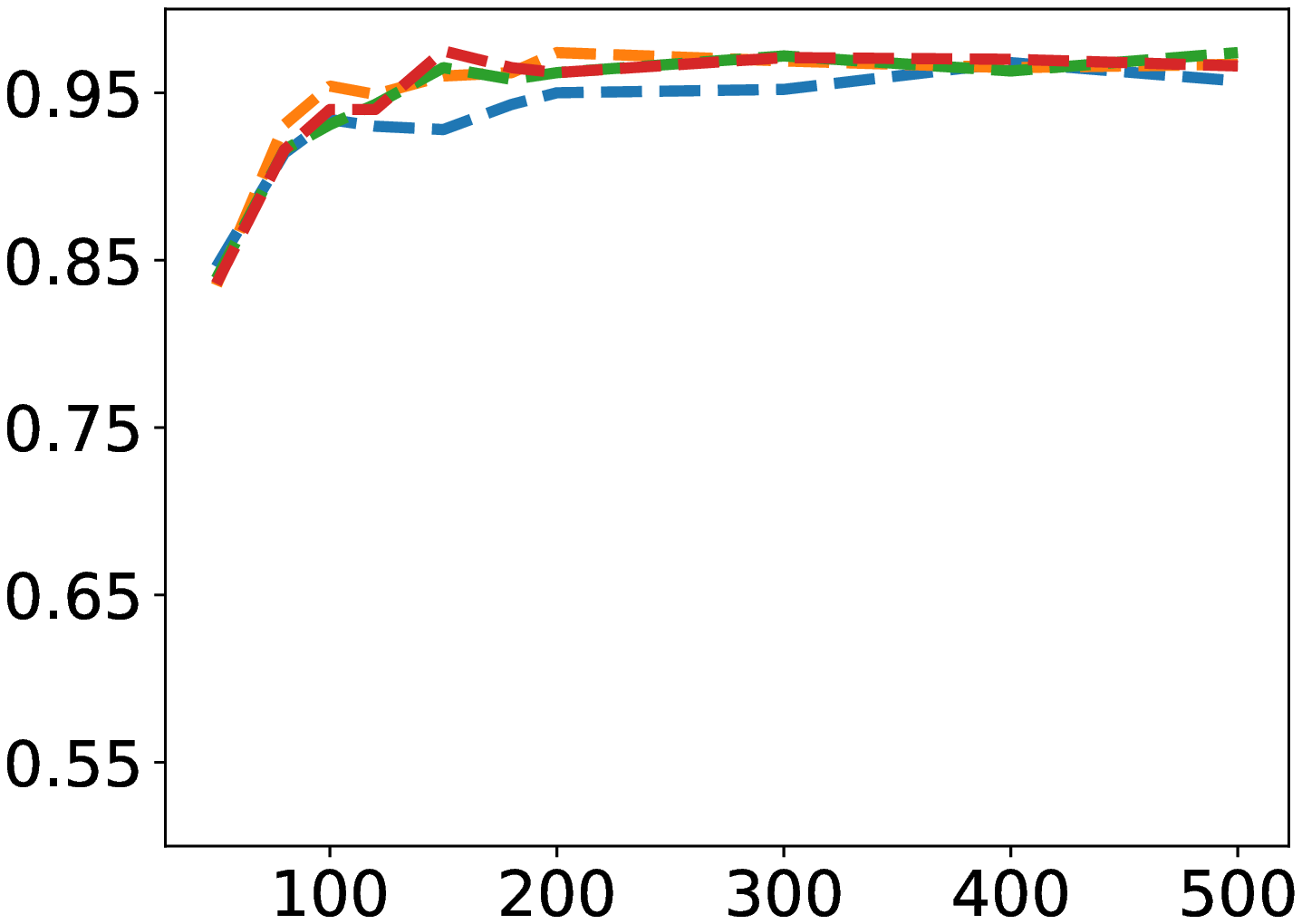} 

		\put(25,80){ \ul{\ \ \  \ $\delta=0.9$ \ \ \ \    }}
 		 				
	\end{overpic}	
\end{figure}
\vspace{-0.7cm}

\begin{figure}[H]	
	\quad\quad\quad 
	\begin{overpic}[width=0.29\textwidth]{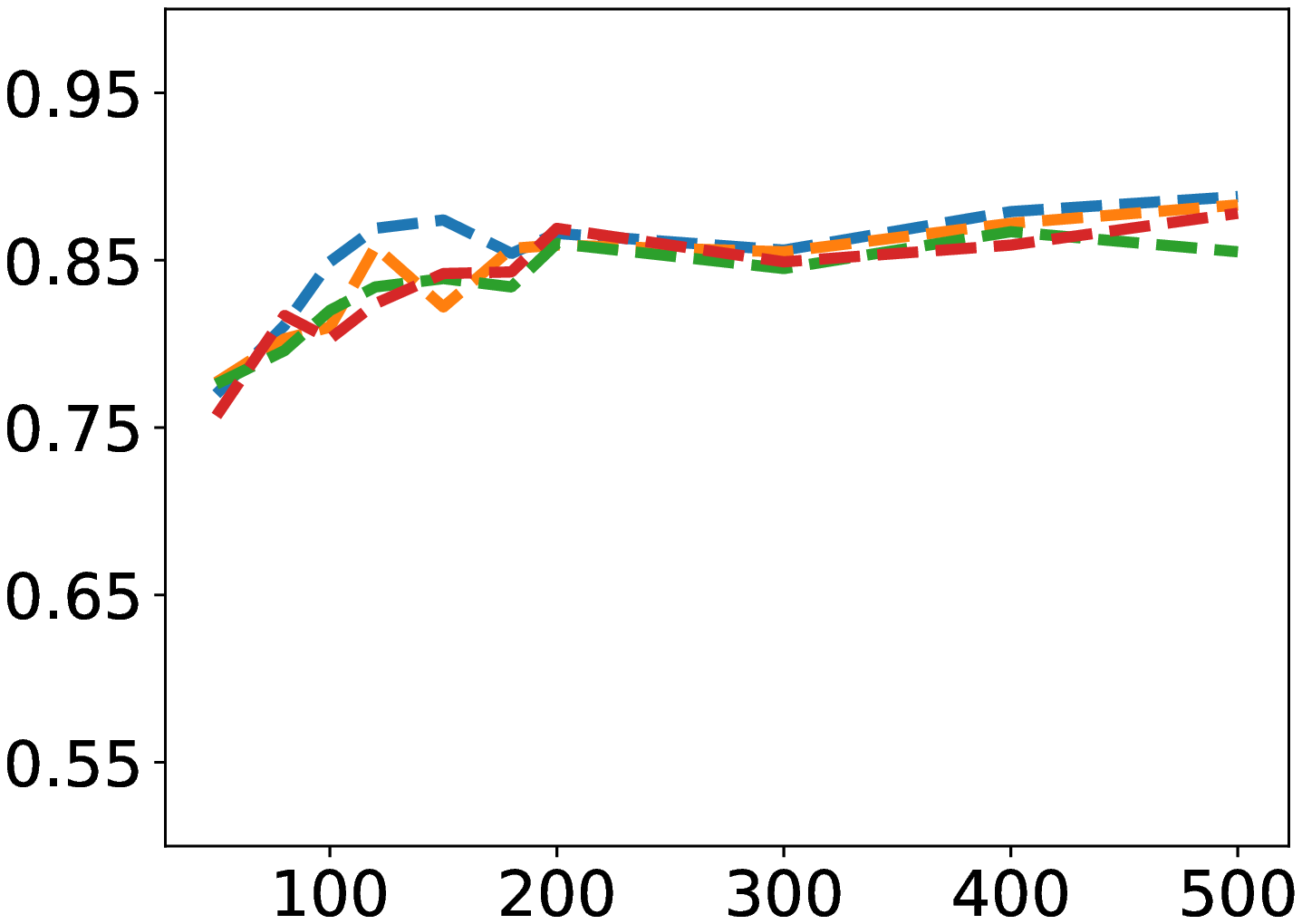} 
	\put(-20,-1){\rotatebox{90}{ {\small \ \ \ standardization \  \ \ }}}
	\end{overpic}
	~
	\DeclareGraphicsExtensions{.png}
	\begin{overpic}[width=0.29\textwidth]{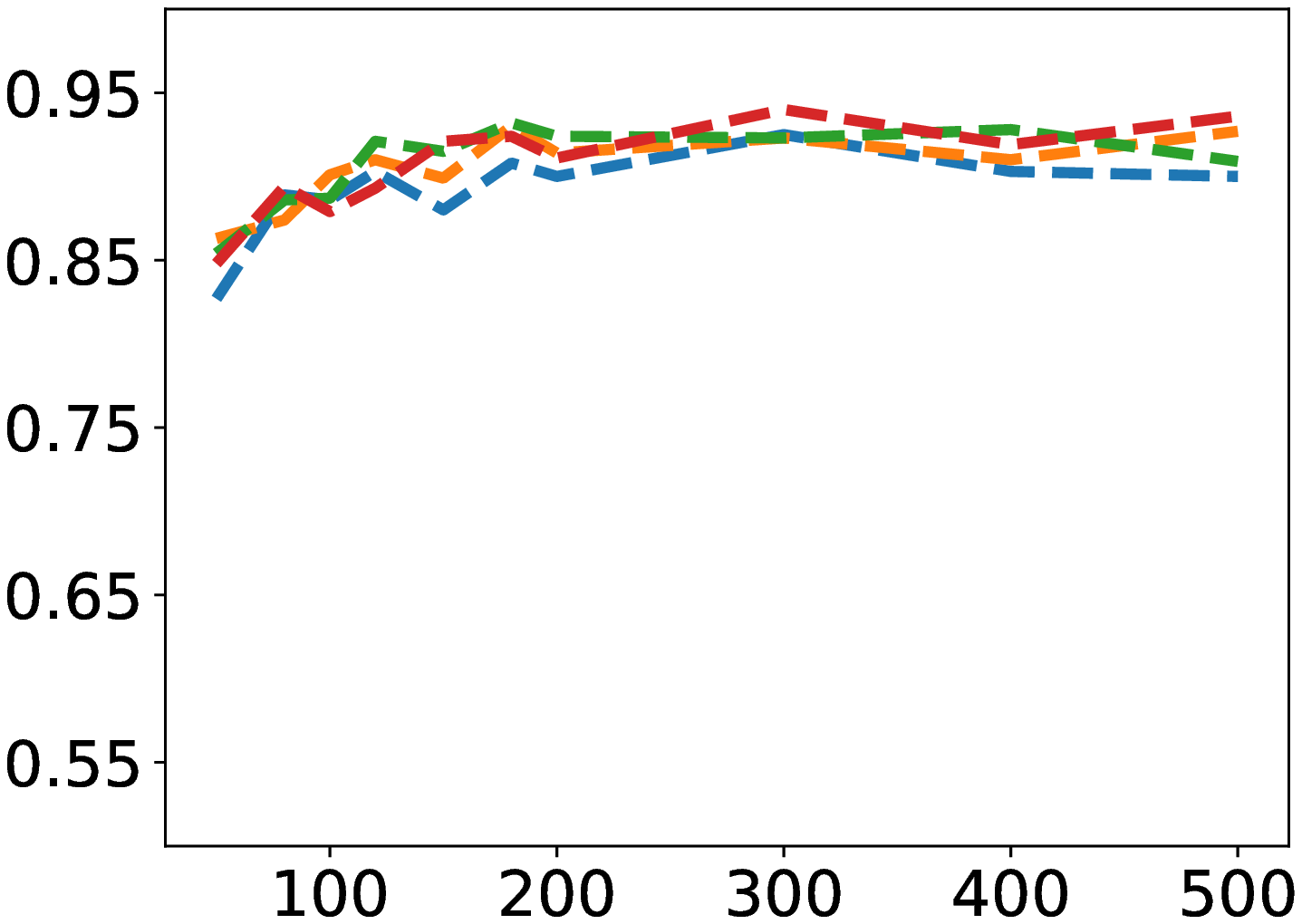} 
	\end{overpic}
	~	
	\begin{overpic}[width=0.29\textwidth]{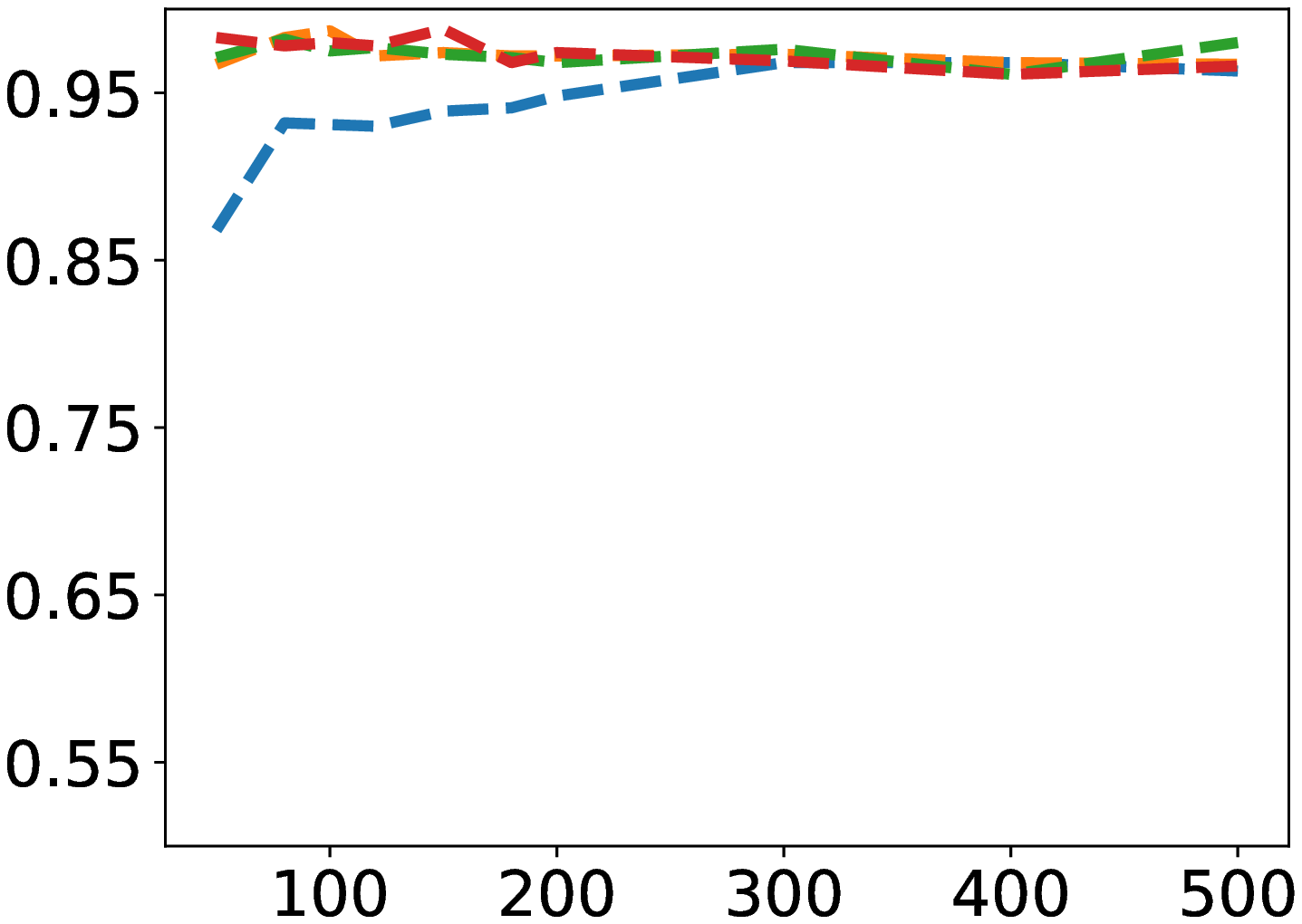} 
		 				
	\end{overpic}	
\end{figure}

\vspace{-0.7cm}

\begin{figure}[H]	
	\quad\quad\quad 
	\begin{overpic}[width=0.29\textwidth]{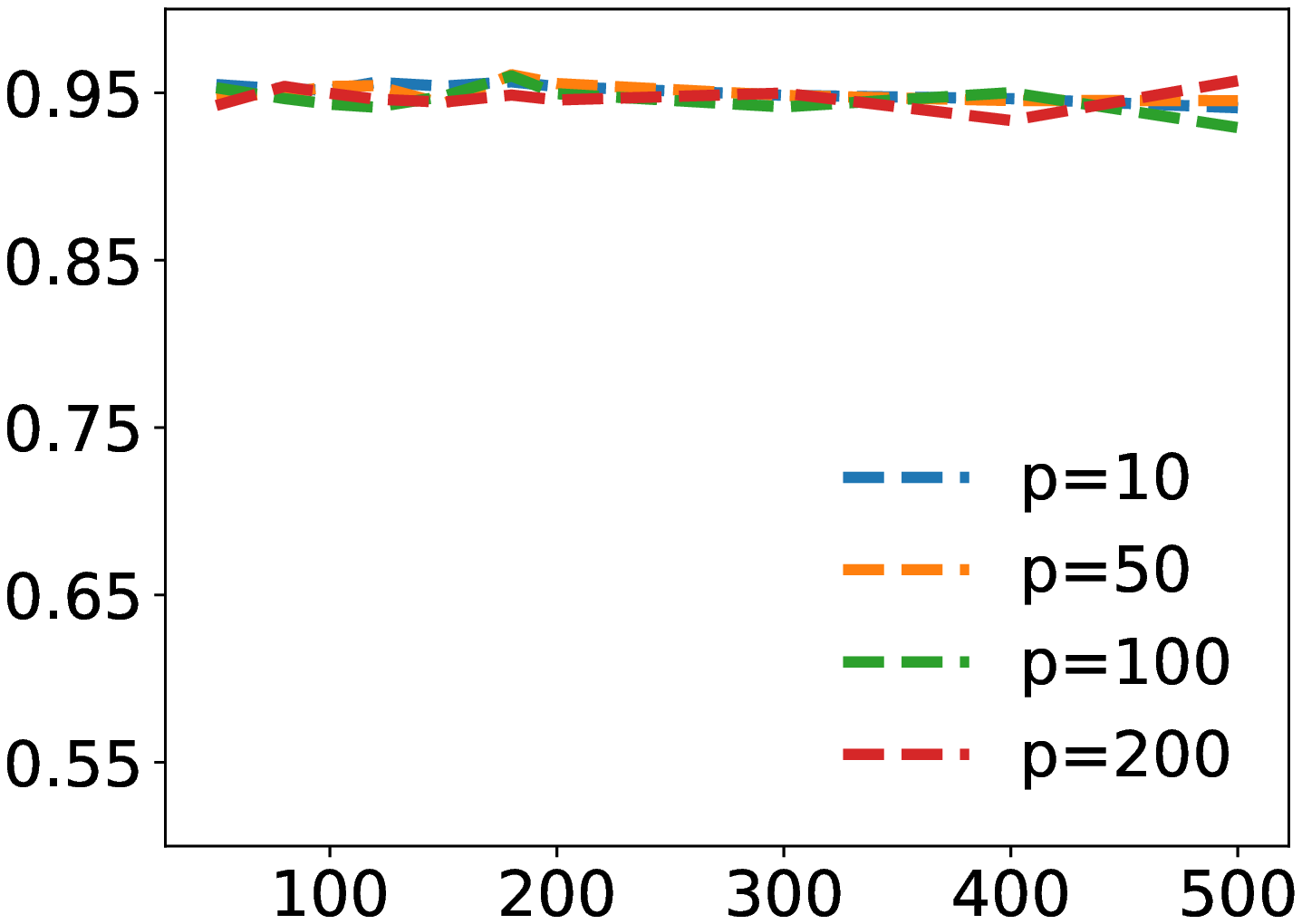} 
		\put(-21,1){\rotatebox{90}{\ $\sqrt{ \ \ }$}}
	    \put(-20,-3){\rotatebox{90}{  { \ \ \ \ \ \ \small transformation \ \ } }}

	\end{overpic}
	~
	\DeclareGraphicsExtensions{.png}
	\begin{overpic}[width=0.29\textwidth]{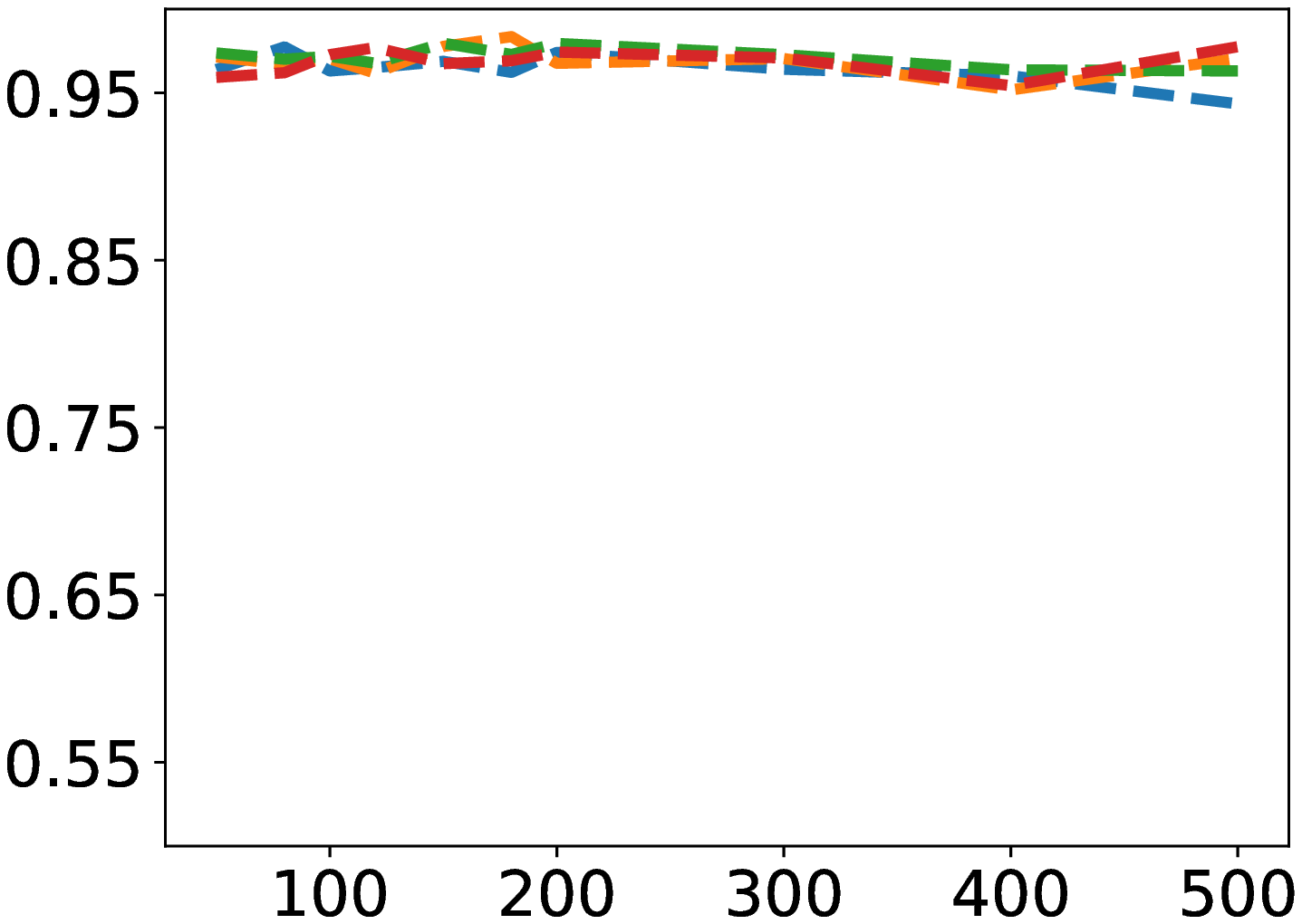} 
	\end{overpic}
	~	
	\begin{overpic}[width=0.29\textwidth]{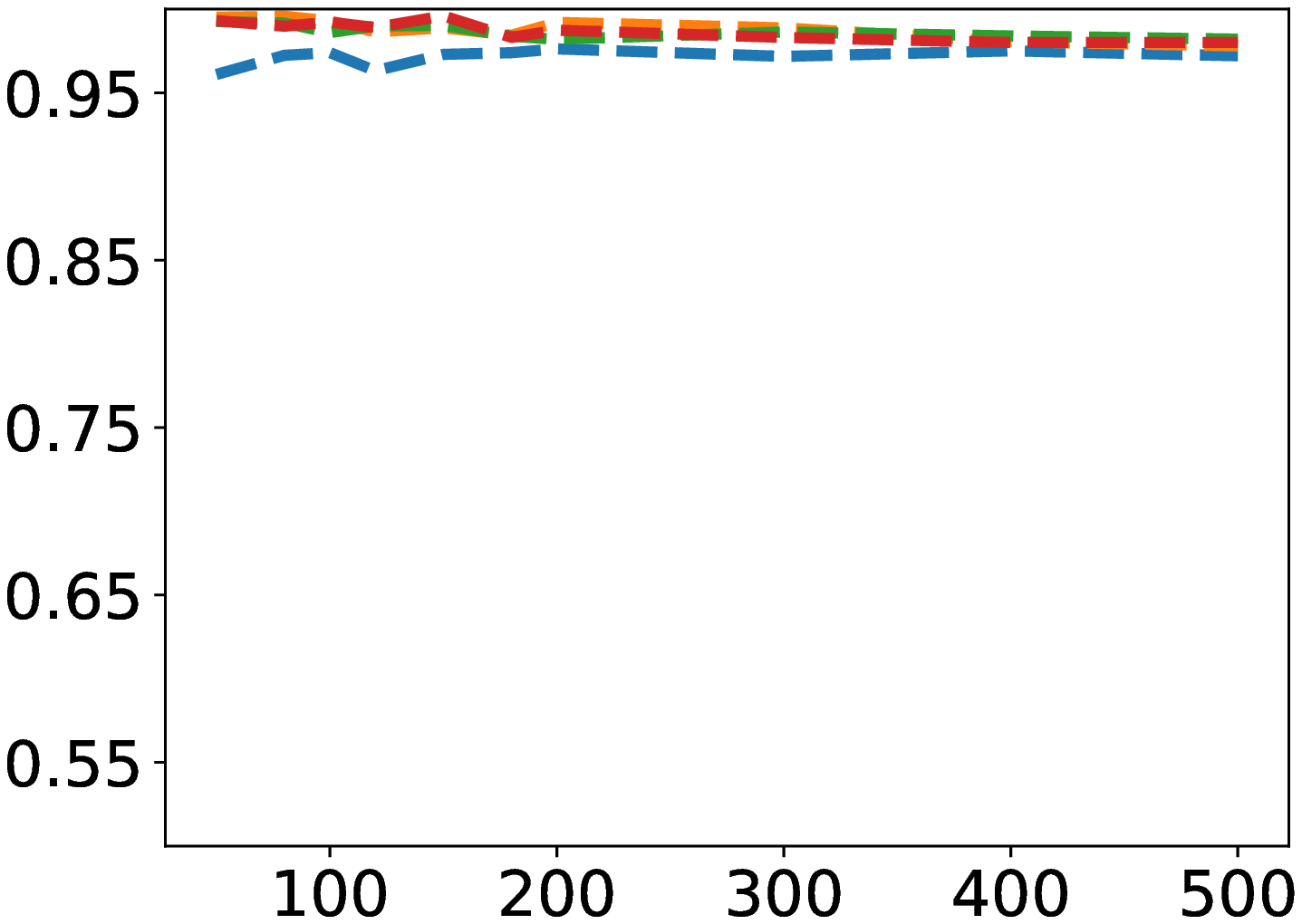} 
		 				
	\end{overpic}	
	\vspace{+.2cm}
	\caption{(Simultaneous coverage probability versus $n$ in simulation model (i) with an exponential decay profile). The plotting scheme is the same as described in the caption of Figure~\ref{fig1}, except that the three columns correspond to values of the eigenvalue decay parameter $\delta$.}
	\label{fig3}
\end{figure}

There are several notable patterns in Figures~\ref{fig1} to discuss. The first is that faster rates of decay tend to lead to better coverage accuracy---as anticipated by our theoretical results. In particular, when the eigenvalue decay parameter is set to $\gamma=1.3$, the coverage is rather accurate even when $n\ll p$. Furthermore, the accuracy is essentially unaffected by the dimension $p$ in this situation, as indicated by the overlap of the four colored curves. On the other hand, as the decay parameter becomes smaller, 
the three transformations perform in different ways. For instance, when $\gamma=0.7$, $p=200$, and $n<200$, the log transformation yields coverage that clearly falls short of the nominal level. By contrast, the standardization and square-root transformations tend to err more safely in the conservative direction when $\gamma=0.7$. To give some indication of the difficulty of $\gamma=0.7$, it should be noted that if $\gamma$ were decreased slightly to 0.5 with $p\gtrsim n$, this would imply $\ts {\tt{r}}(\Sigma) / \sqrt{n} \asymp \sqrt{p/n}\gtrsim 1$, in which case bootstrap consistency would not be guaranteed. When considering all three cases $\gamma=0.7,0.8,0.9$ collectively, the square-root transformation seems to yield the best overall coverage results if conservative errors are viewed as preferable to anti-conservative ones.

Turning to the coverage results for exponential spectrum decay, the log and square-root transformations continue to follow the pattern that faster decay improves coverage accuracy. Also, the log transformation maintains its tendency to err in the anti-conservative direction, while the square-root transformation maintains its tendency to err in the conservative direction. Meanwhile, ordinary standardization yields larger errors in the anti-conservative direction than it did in the previous context.

\begin{figure}[H]
\vspace{0.5cm}
	\quad\quad\quad 
	\begin{overpic}[width=0.29\textwidth]{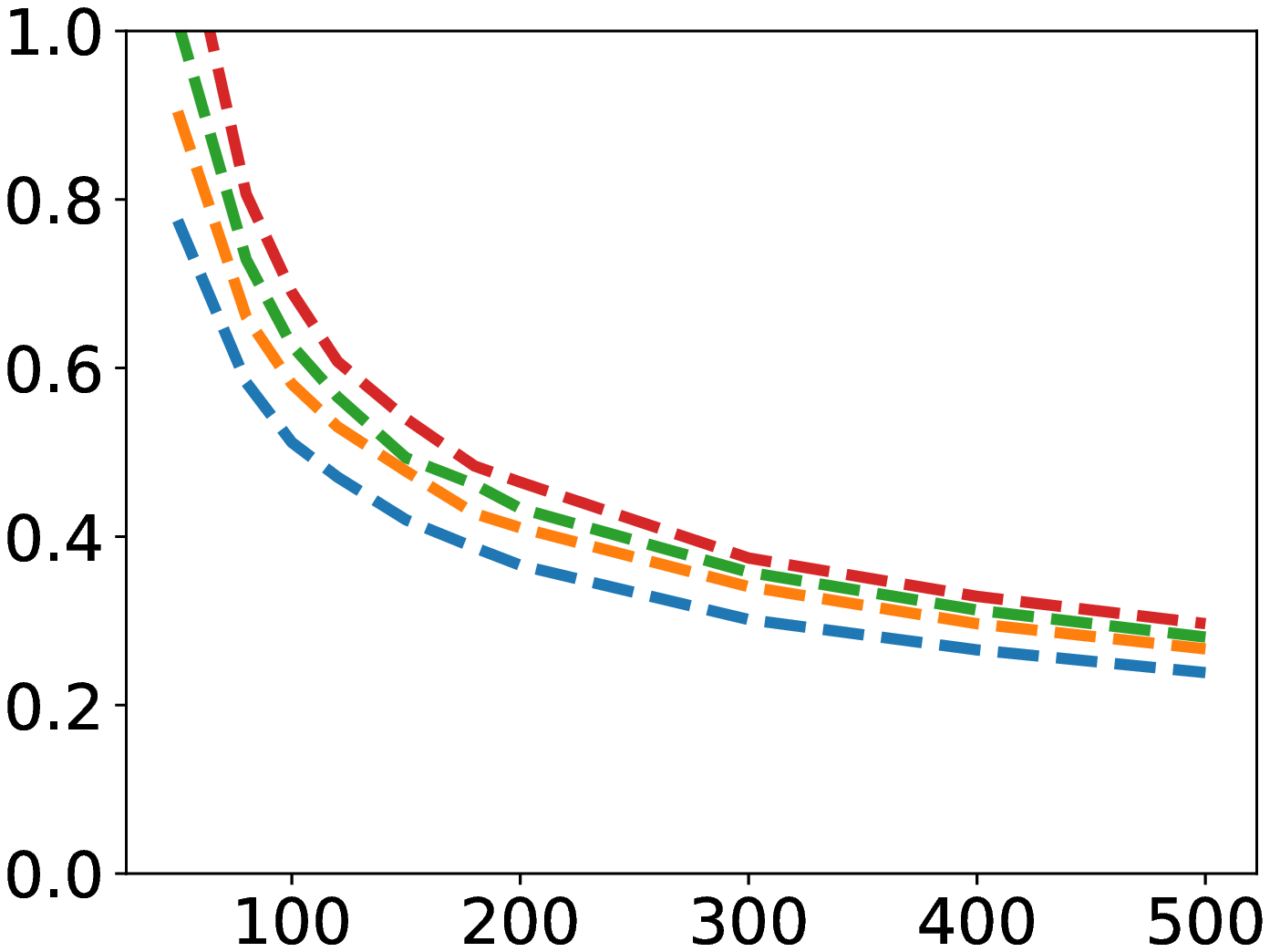} 
		\put(25,80){ \ul{\ \ \  \ $\gamma=0.7$ \ \ \ \    }}
		\put(-20,-5){\rotatebox{90}{ {\small \ \ \ log transformation  \ \ }}}
	\end{overpic}
	~
	\DeclareGraphicsExtensions{.png}
	\begin{overpic}[width=0.29\textwidth]{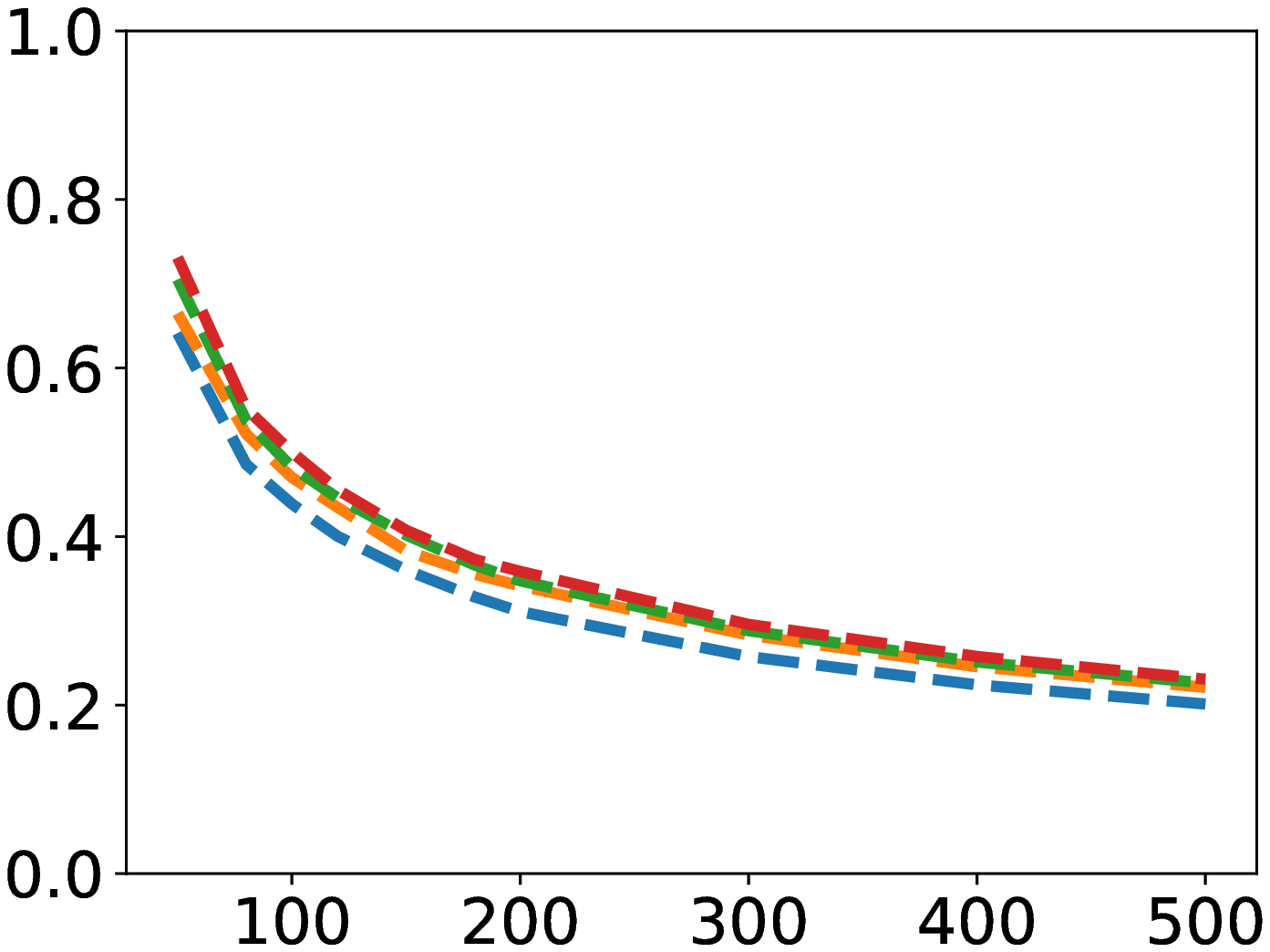} 
		\put(25,80){ \ul{\ \ \  \ $\gamma=1.0$ \ \ \ \    }}
	\end{overpic}
	~	
	\begin{overpic}[width=0.29\textwidth]{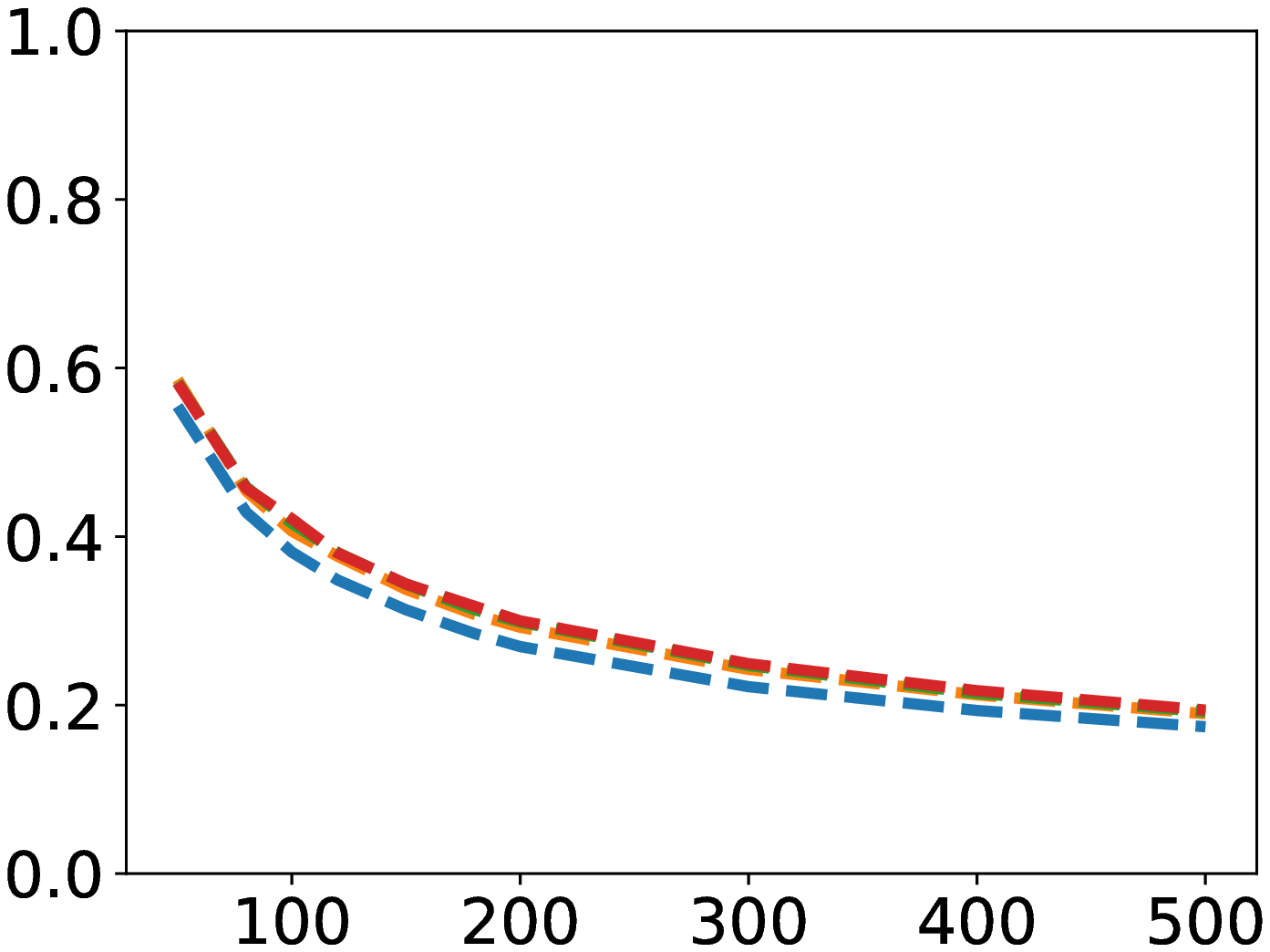} 
		\put(25,80){ \ul{\ \ \  \ $\gamma=1.3$ \ \ \ \    }}
	\end{overpic}	
\end{figure}

\vspace{-0.7cm}

\begin{figure}[H]	
	\quad\quad\quad 
	\begin{overpic}[width=0.29\textwidth]{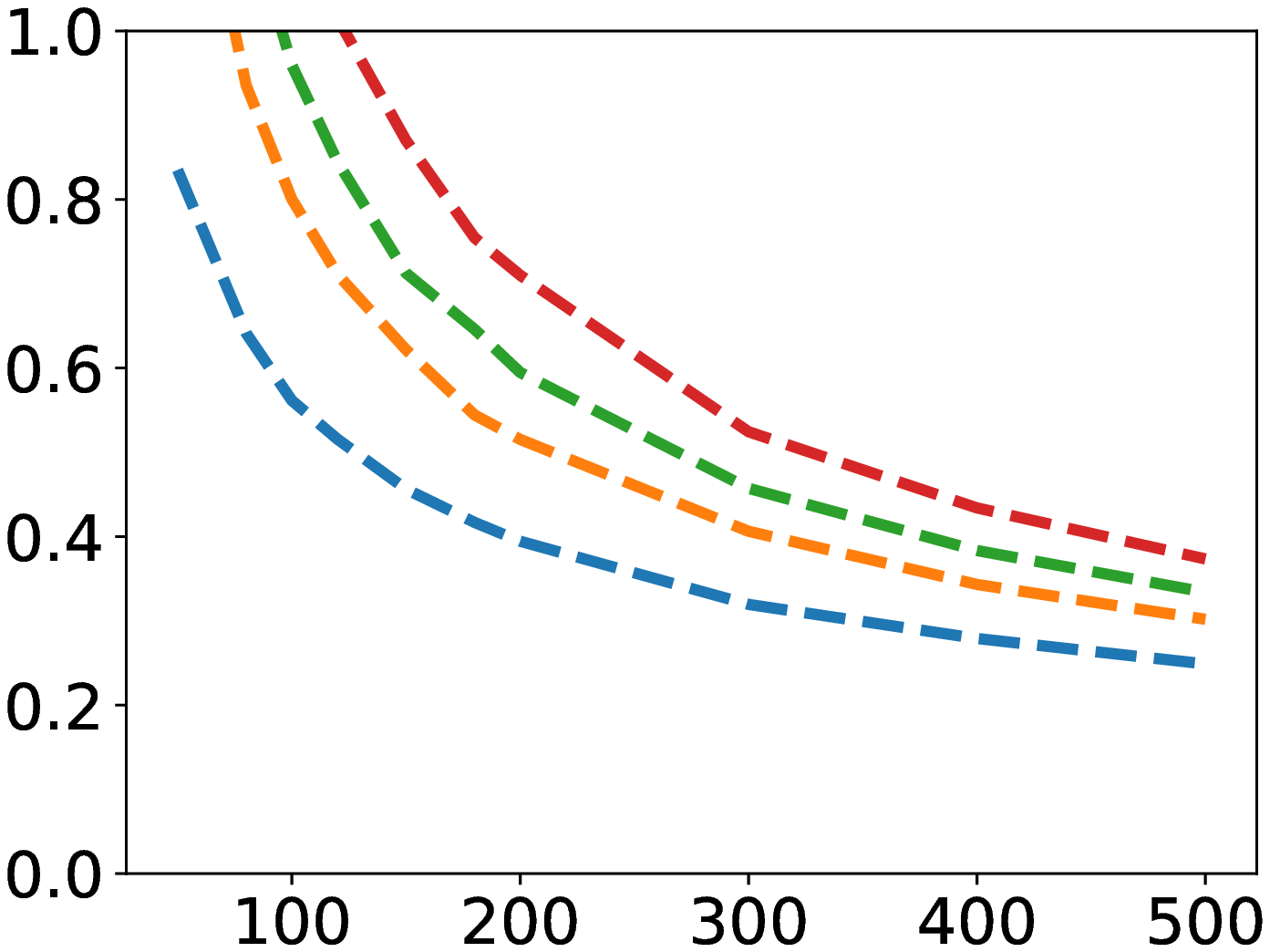} 
	    \put(-20,-1){\rotatebox{90}{ {\small \ \ \ standardization \  \ \ }}}
	\end{overpic}
	~
	\DeclareGraphicsExtensions{.png}
	\begin{overpic}[width=0.29\textwidth]{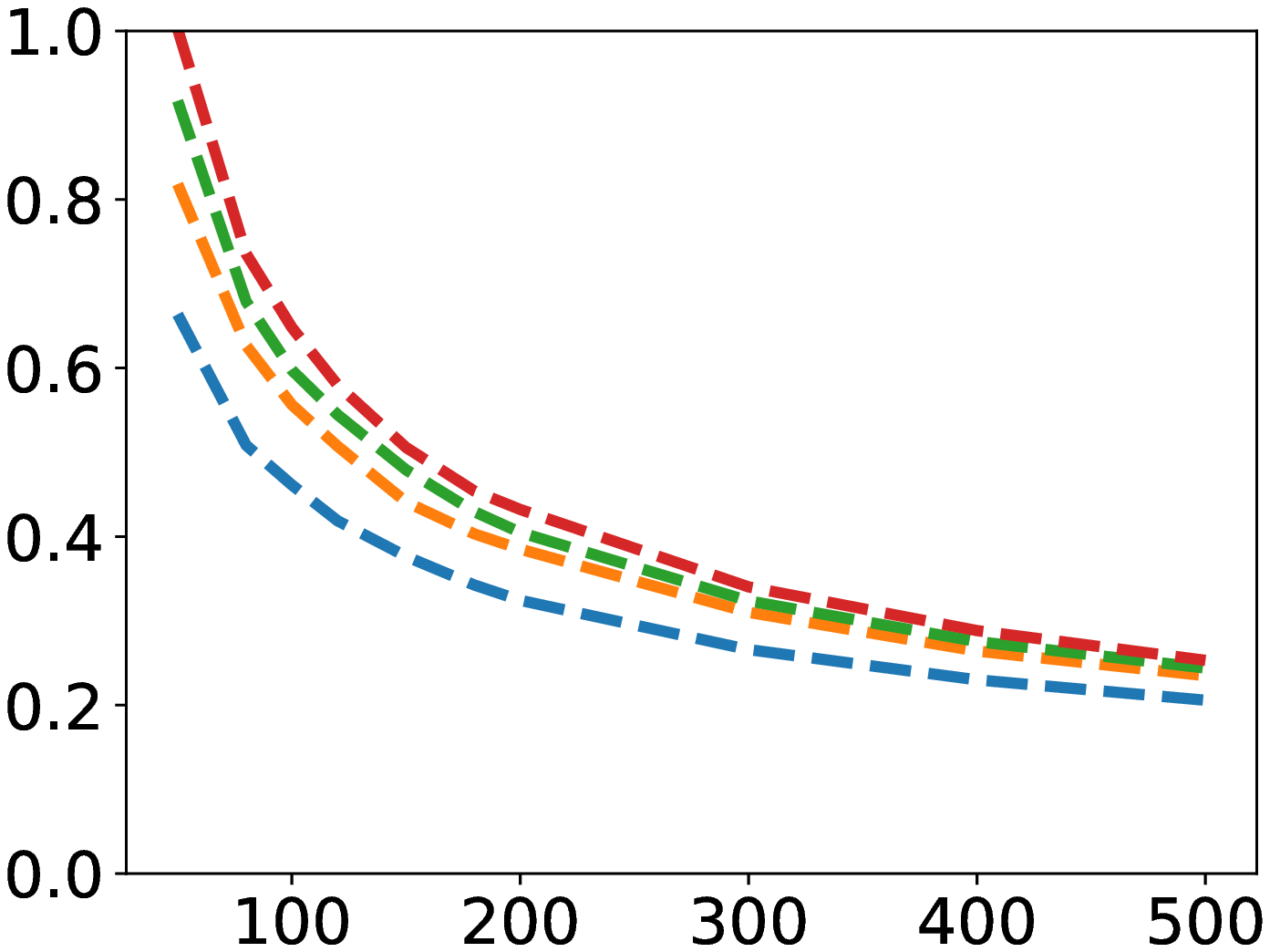} 
	\end{overpic}
	~	
	\begin{overpic}[width=0.29\textwidth]{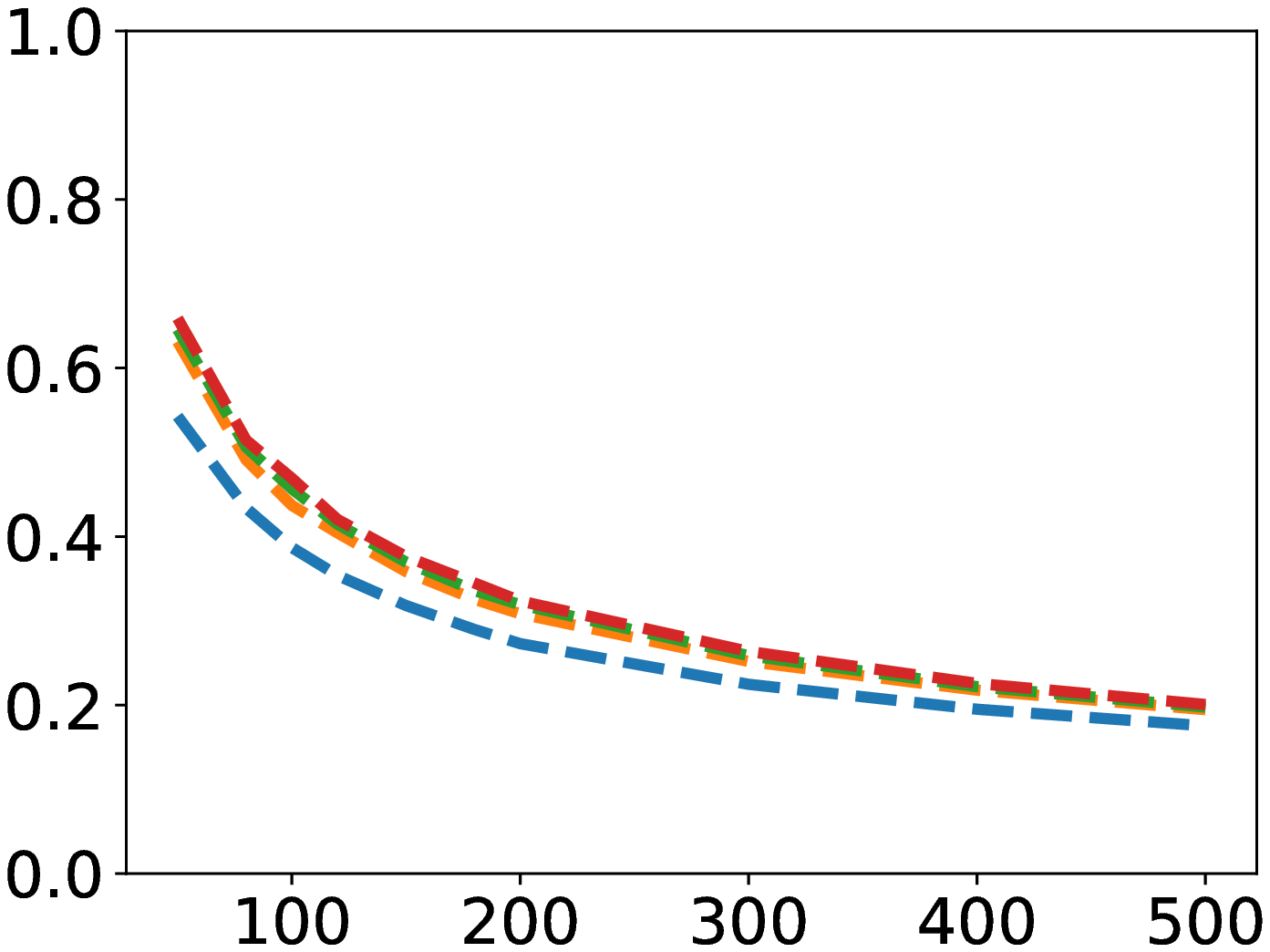} 
	\end{overpic}	
	%
\end{figure}

\vspace{-0.7cm}

\begin{figure}[H]	
	\quad\quad\quad 
	\begin{overpic}[width=0.29\textwidth]{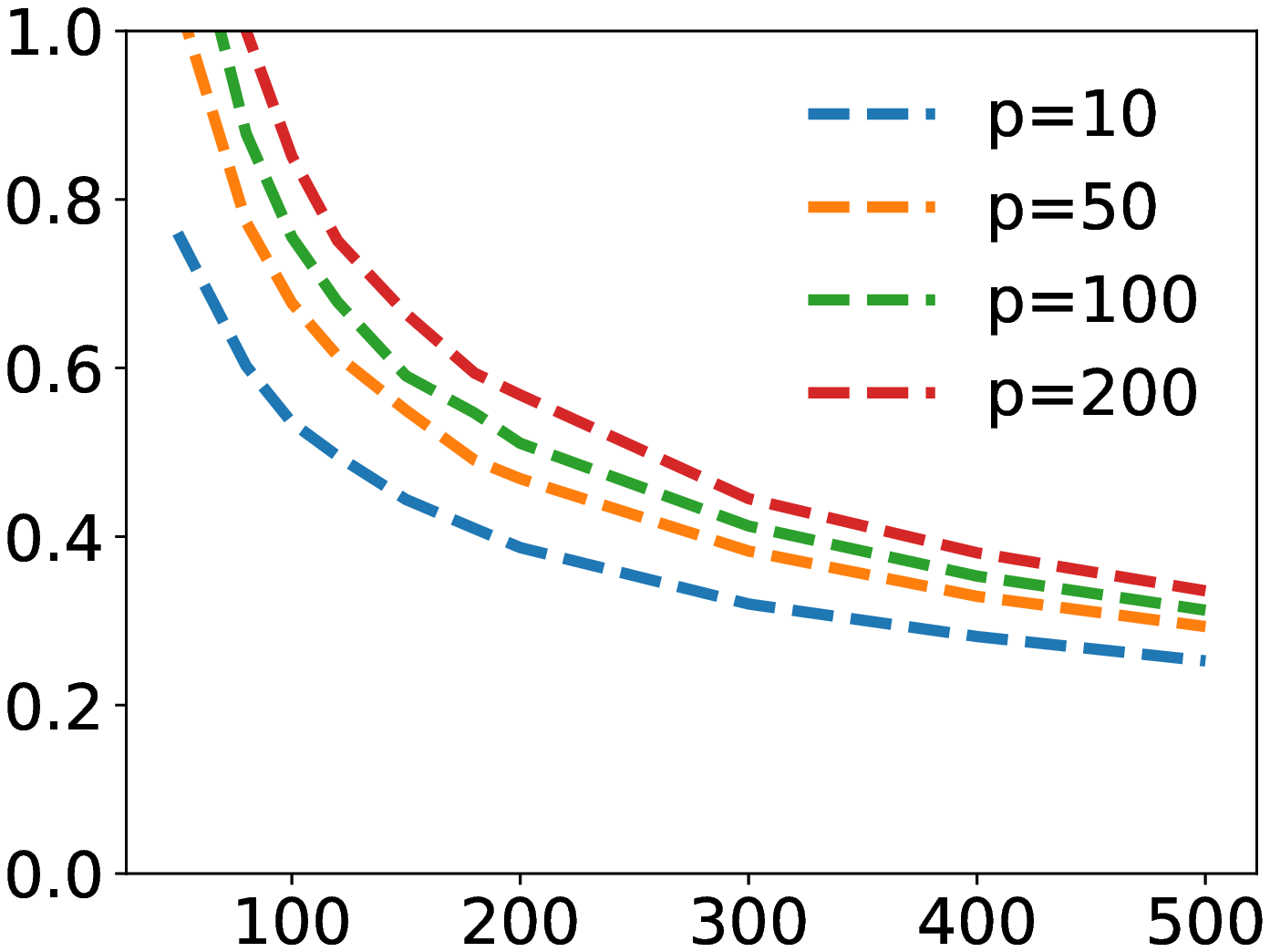} 
	    \put(-21,1){\rotatebox{90}{\ $\sqrt{ \ \ }$}}
	    \put(-20,-3){\rotatebox{90}{  { \ \ \ \ \ \ \small transformation \ \ } }}

	\end{overpic}
	~
	\DeclareGraphicsExtensions{.png}
	\begin{overpic}[width=0.29\textwidth]{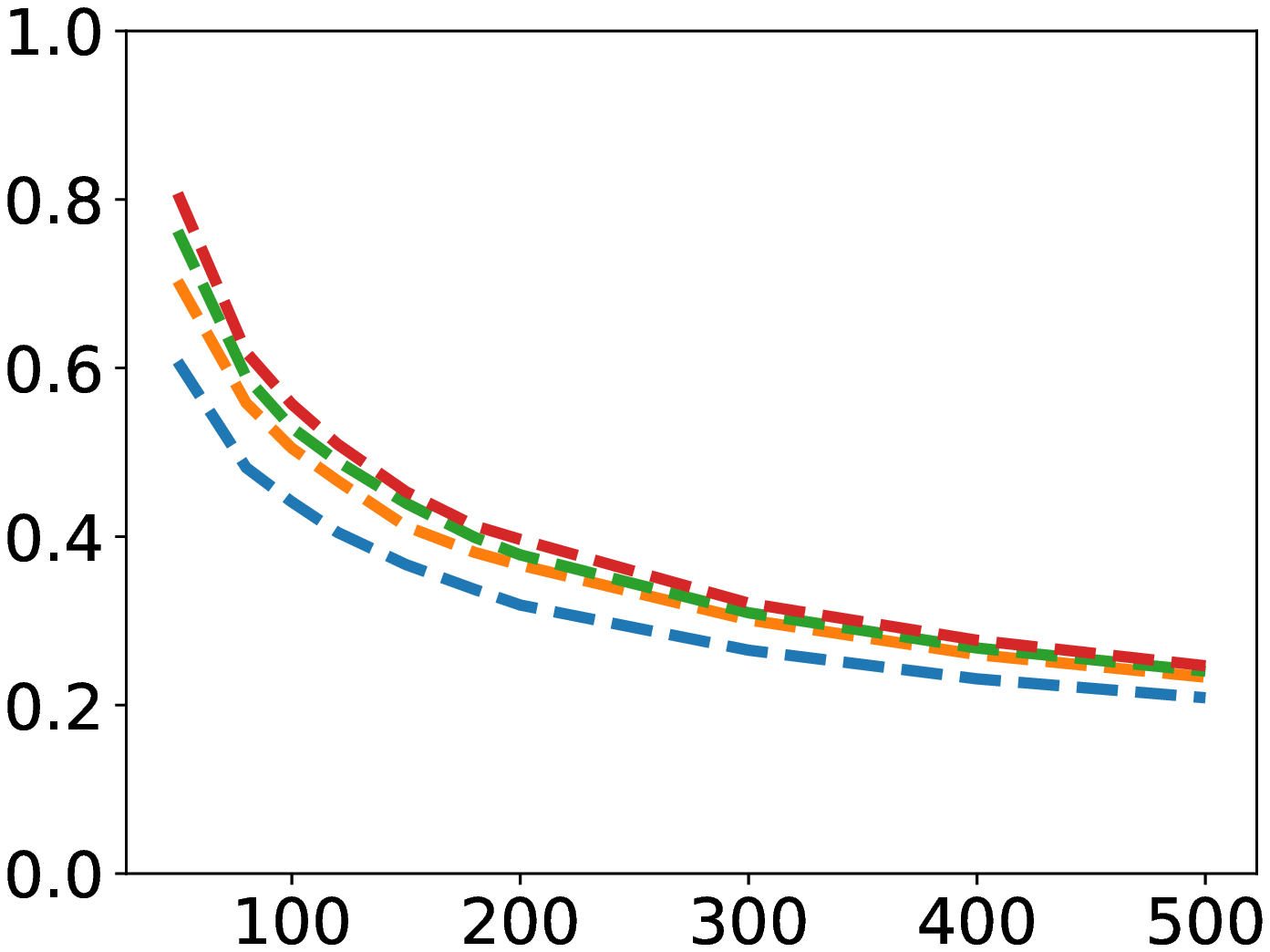} 
	\end{overpic}
	~	
	\begin{overpic}[width=0.29\textwidth]{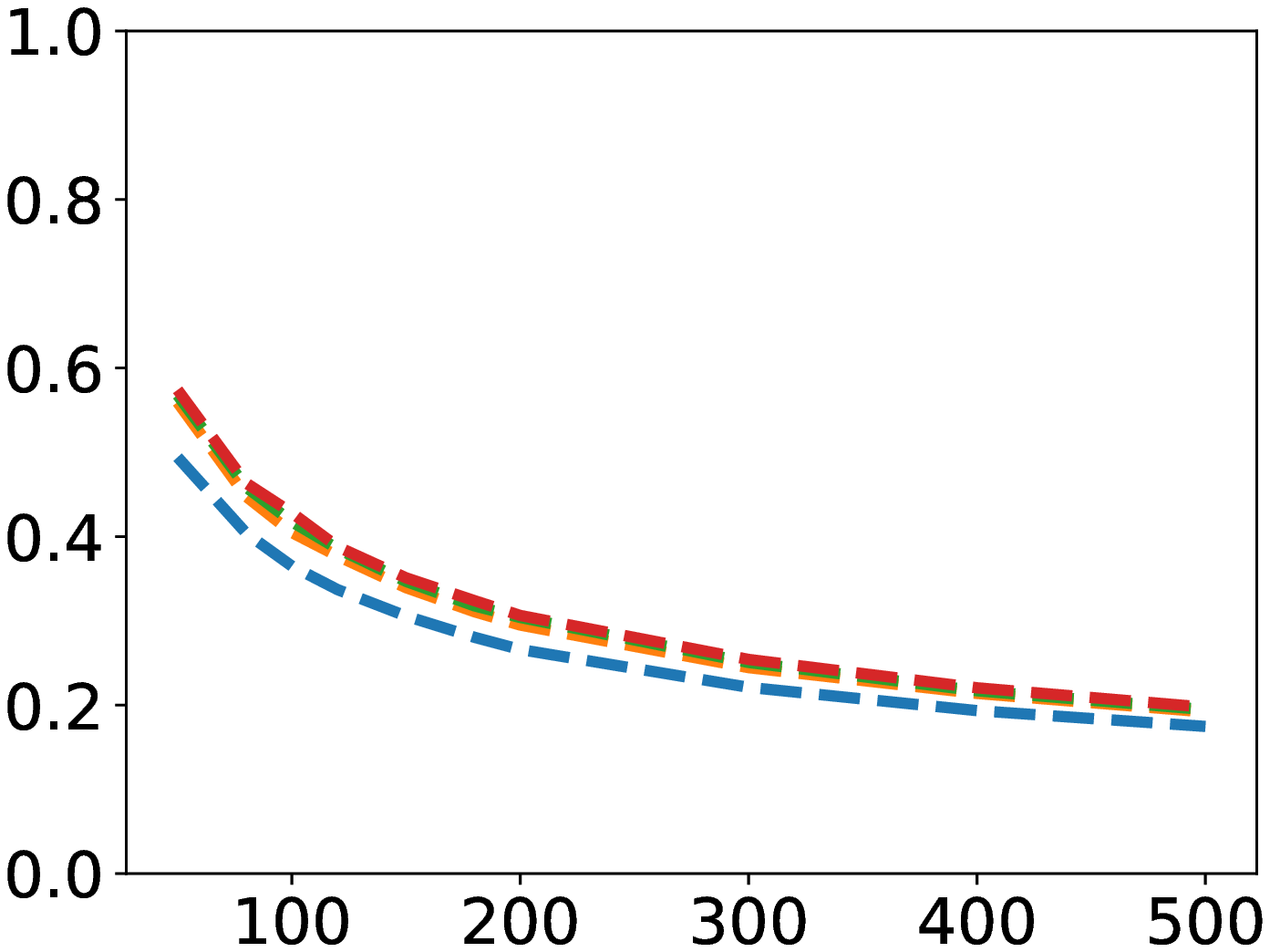} 
	\end{overpic}	
	\vspace{+.1cm}		
	\caption{(Average width versus $n$ in simulation model (i) with a polynomial decay profile). In each of the nine panels, the $y$-axis measures the average width $\E[|\hat{\I}_1|+\cdots+|\hat{\I}_{5}|] / 5$, and the $x$-axis measures $n$. The colored curves correspond to the different values of $p=10$, $50$, $100$, $200$, indicated in the legend. The three rows and three columns correspond to labeled choices of transformations and values of the eigenvalue decay parameter $\gamma$.} \label{fig5}
\end{figure}

\subsection{Discussion of width}
Beyond coverage probability, interval width is another important factor to consider when appraising confidence intervals. In Figures~\ref{fig5}-\ref{fig_width_delta_ii}, the average width~$\E[|\hat{\mathcal{I}}_1|+\cdots+|\hat{\mathcal{I}}_k|] / k$ is plotted on the $y$-axis as a function of the sample size $n$ on the $x$-axis, with the underlying parameter settings being organized in the same manner as in Figures~\ref{fig1}-\ref{fig4}. (Corresponding results for settings based on model (ii) and a nominal value of 90\% are presented in  Section~\ref{appendix_last} of the supplementary material.) With regard to the three transformations, they produce intervals that have roughly similar widths across most parameter settings. However, at a more fine-grained level, the results in the case of polynomial spectrum decay show that the log transformation tends to yield slightly shorter widths than the square-root transformation, which in turn, tends to yield slightly shorter widths than ordinary standardization. In the case of exponential spectrum decay with $\delta=0.9$, the same pattern is also apparent, while for smaller values of $\delta$, there is not much difference among the transformations.

Aside from the transformations, there are two other general trends to notice. Within each of the 18 panels of Figures~\ref{fig5}-\ref{fig_width_delta_ii}, there is a monotone relationship between width and the dimension $p$, with the width generally increasing as the dimension increases. Similarly, the width generally also increases as the effective rank ${\tt{r}}(\Sigma)$ increases. \\[-0.2cm]

\begin{figure}[H]	
\vspace{0.5cm}
	\quad\quad\quad 
	\begin{overpic}[width=0.29\textwidth]{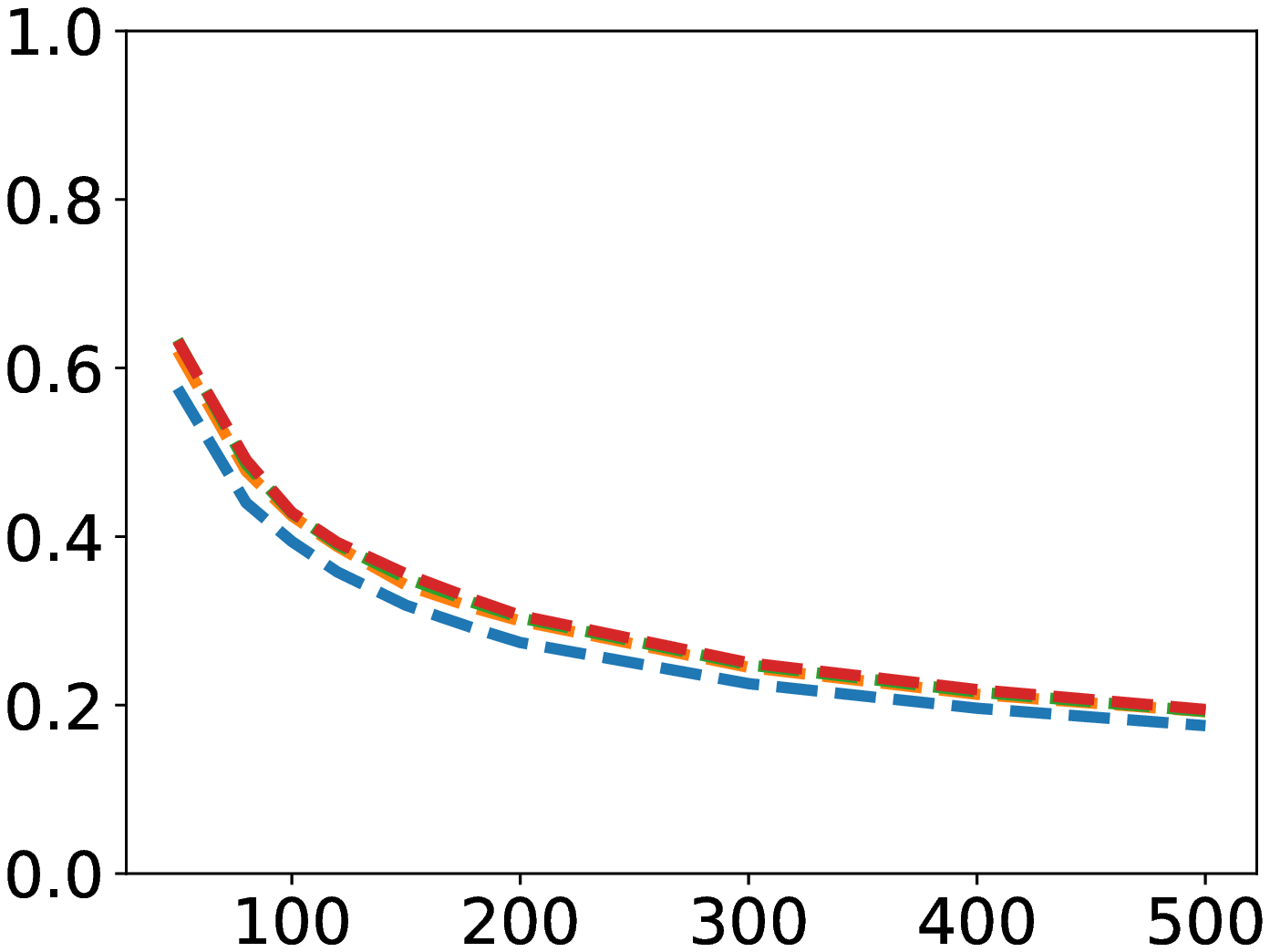} 
\put(25,80){ \ul{\ \ \  \ $\delta=0.7$ \ \ \ \    }}
		\put(-20,-5){\rotatebox{90}{ {\small \ \ \ log transformation  \ \ }}}
\end{overpic}
	~
	\DeclareGraphicsExtensions{.png}
	\begin{overpic}[width=0.29\textwidth]{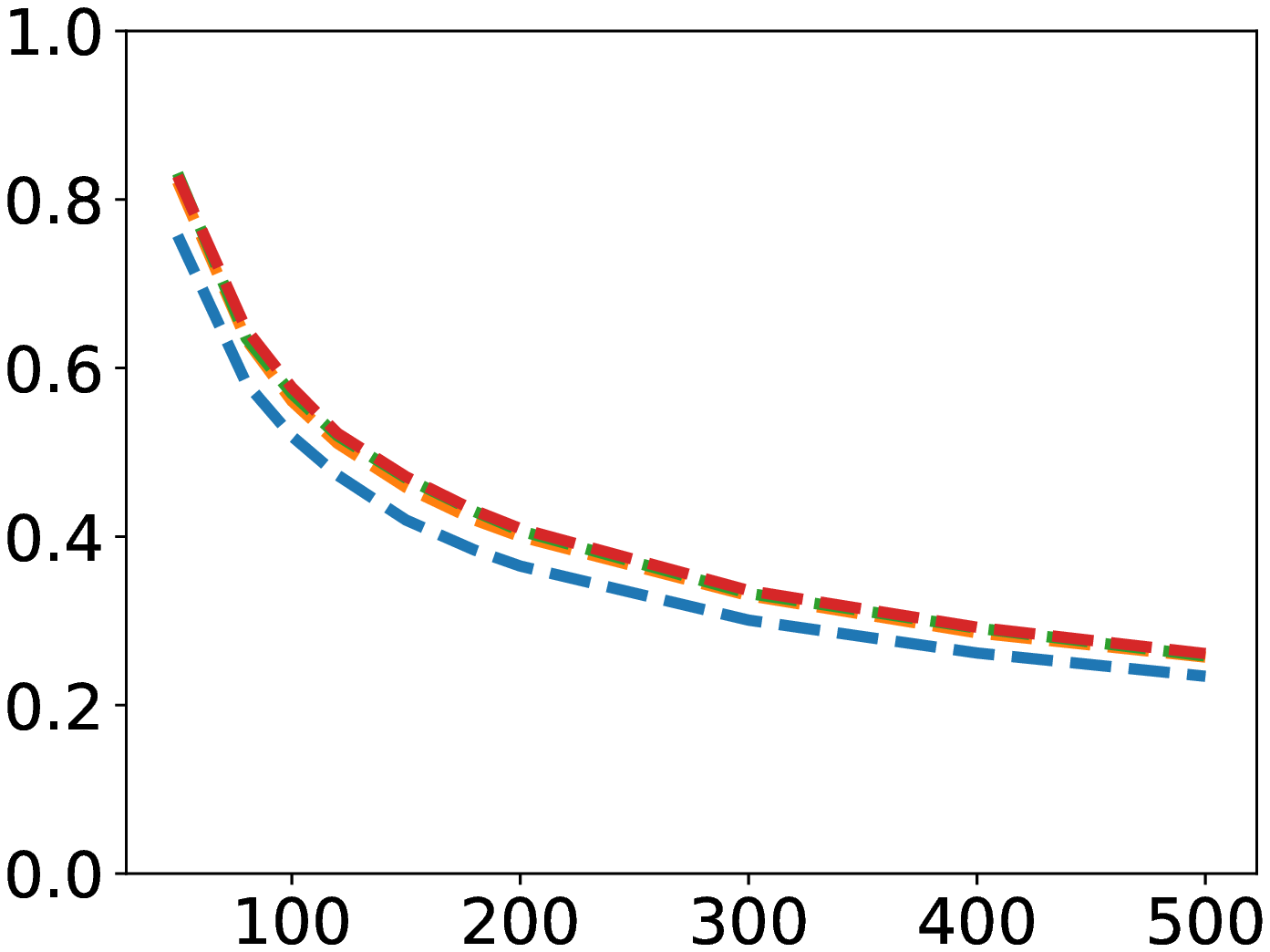} 
		\put(15,80){ \ul{\ \ \  \ $\delta=0.8$ \ \ \ \    }}
	\end{overpic}
	~	
	\begin{overpic}[width=0.29\textwidth]{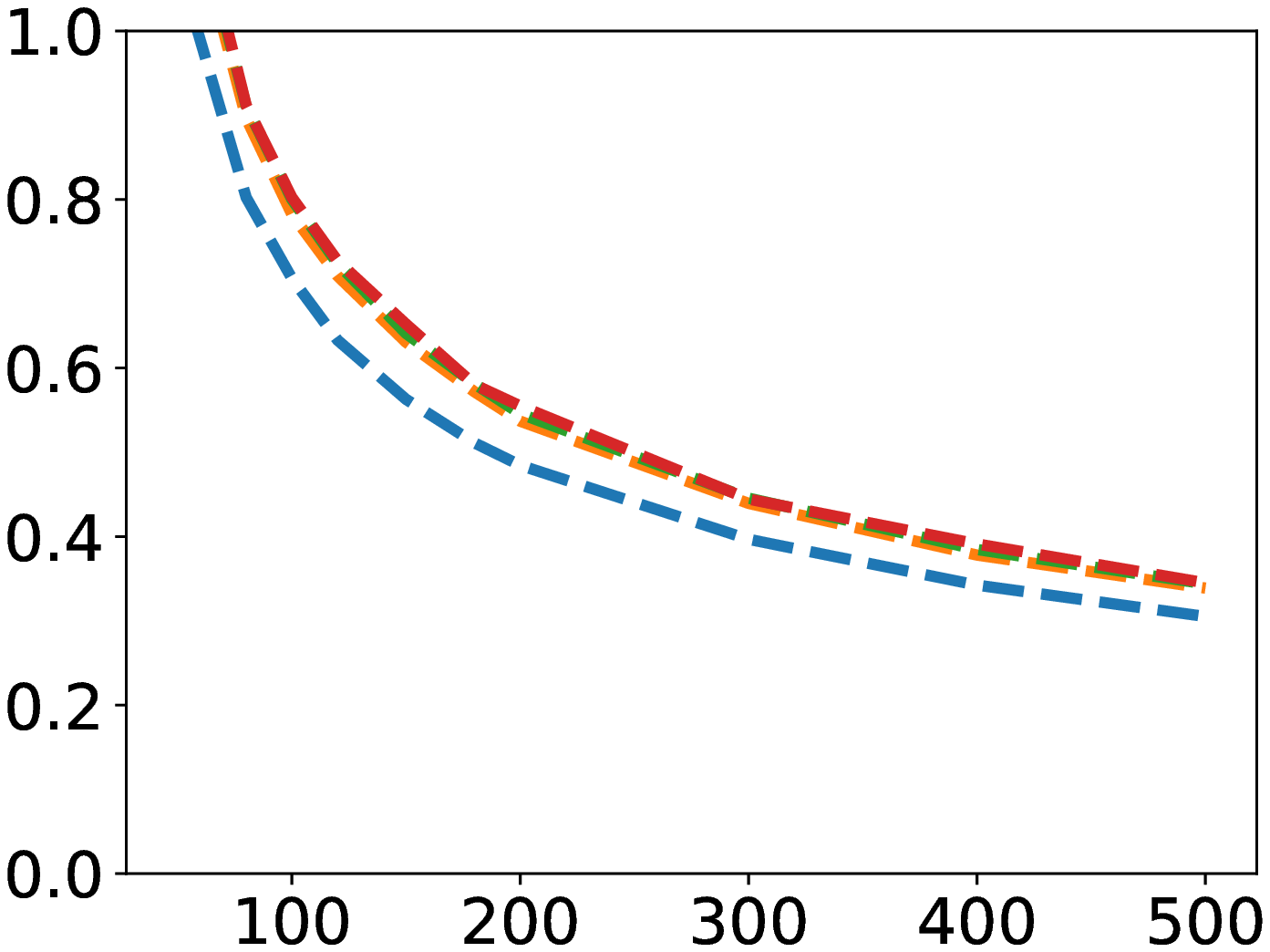} 
		\put(15,80){ \ul{\ \ \  \ $\delta=0.9$ \ \ \ \    }}
	\end{overpic}	
\end{figure}

\vspace{-0.5cm}

\begin{figure}[H]	
	\quad\quad\quad 
	\begin{overpic}[width=0.29\textwidth]{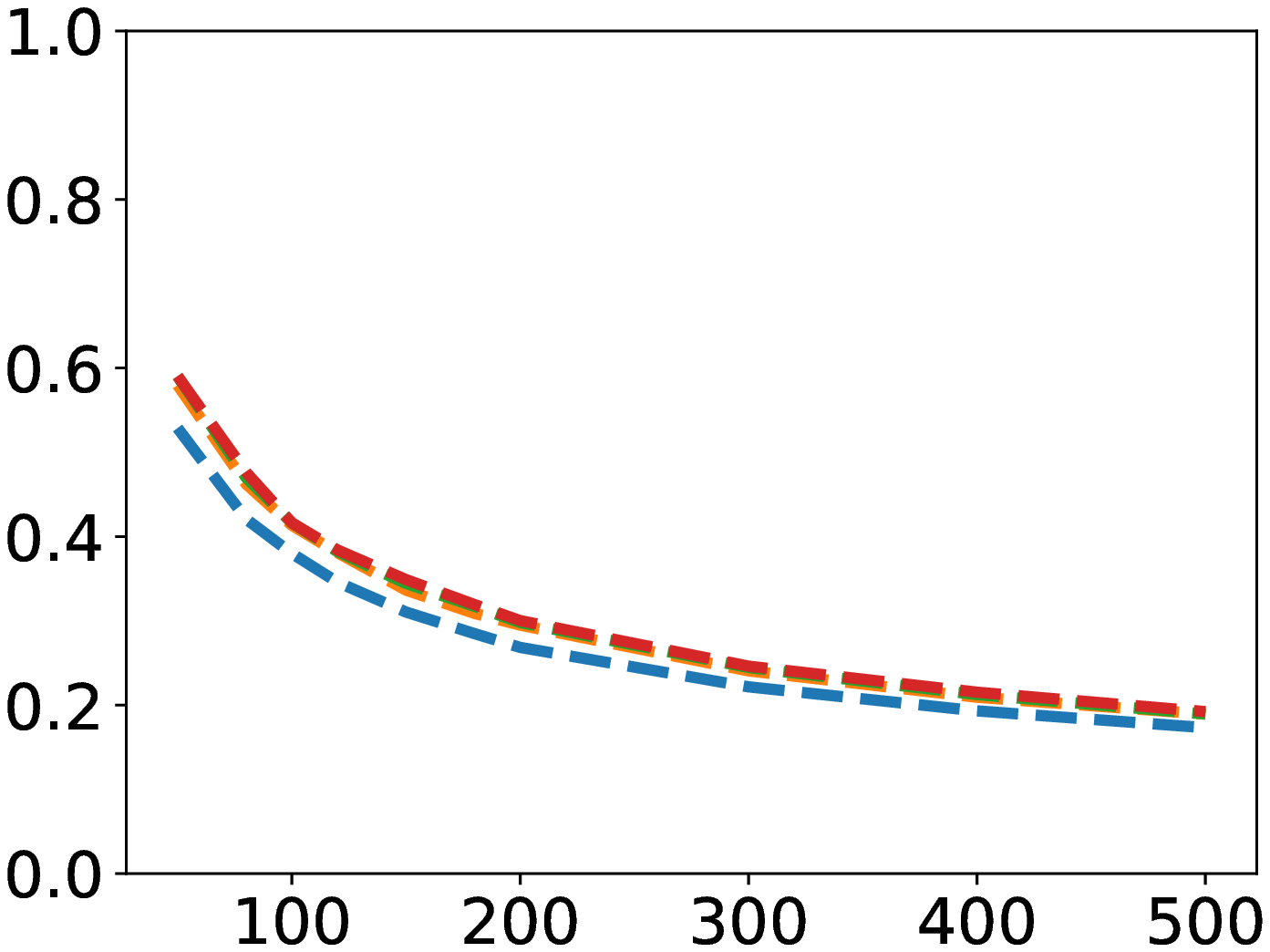} 
	\put(-20,-1){\rotatebox{90}{ {\small \ \ \ standardization \  \ \ }}}
	\end{overpic}
	~
	\DeclareGraphicsExtensions{.png}
	\begin{overpic}[width=0.29\textwidth]{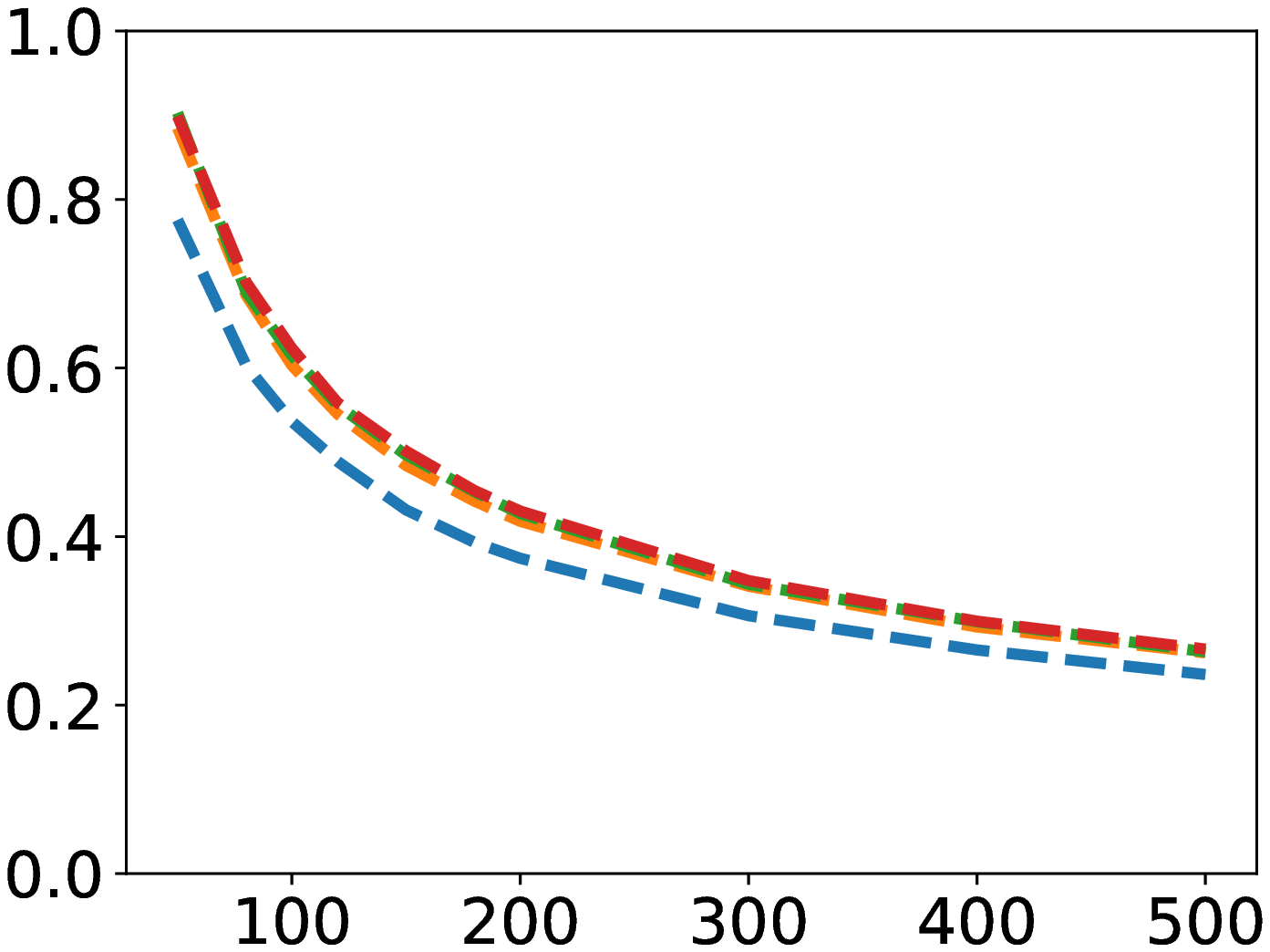} 
	\end{overpic}
	~	
	\begin{overpic}[width=0.29\textwidth]{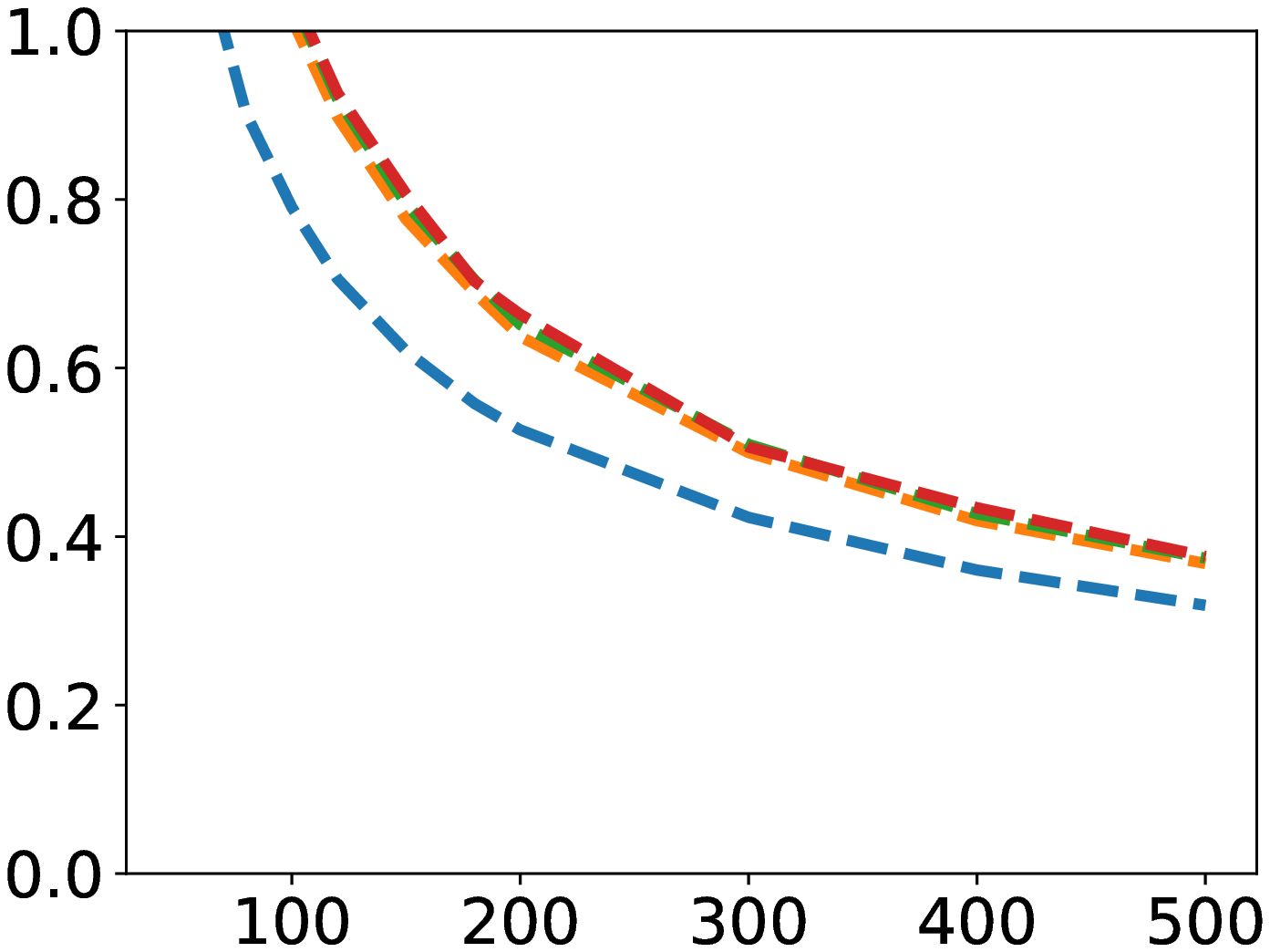} 
		 				
	\end{overpic}	
\end{figure}

\vspace{-0.5cm}

\begin{figure}[H]	
	\quad\quad\quad 
	\begin{overpic}[width=0.29\textwidth]{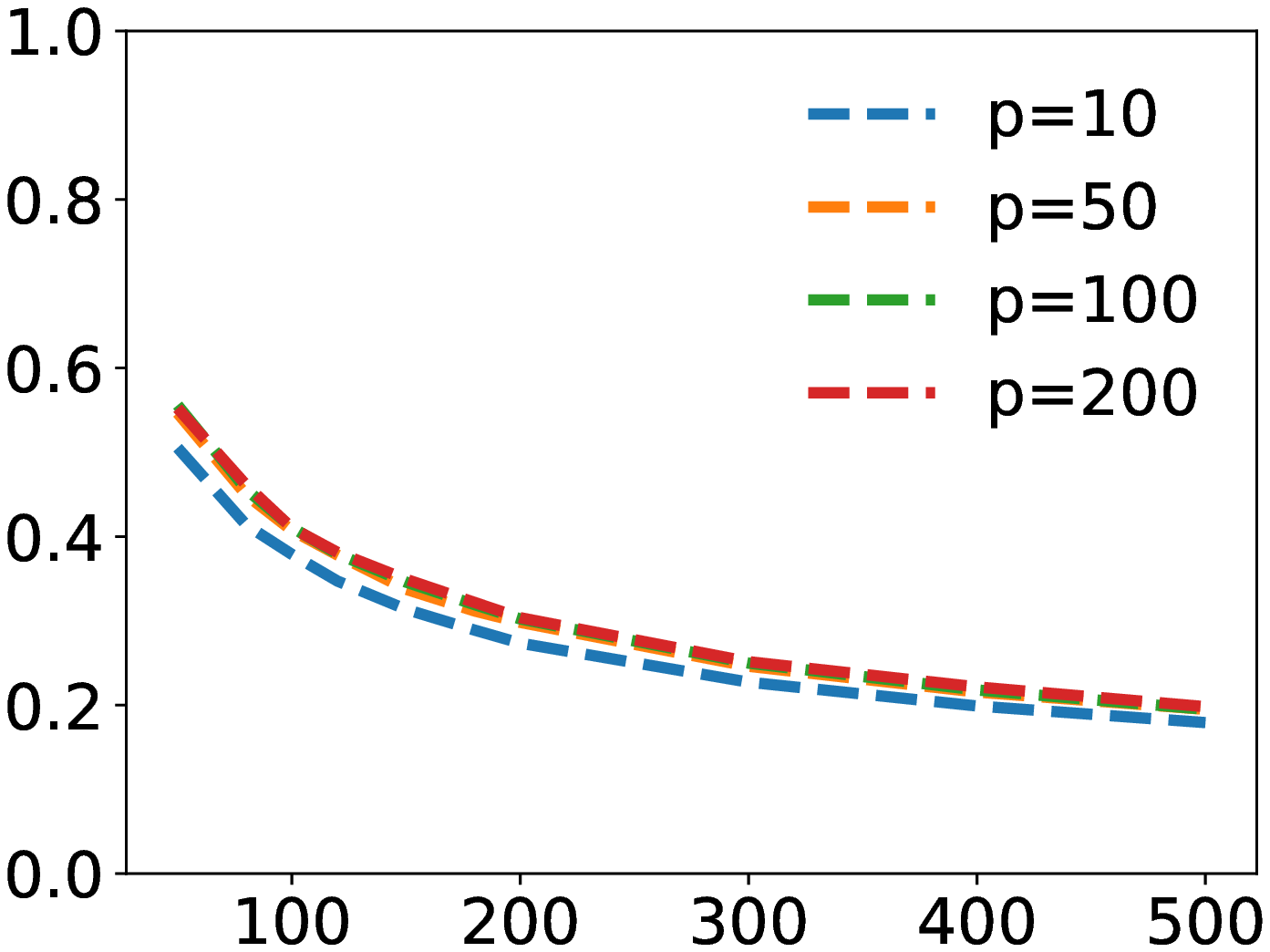} 
		\put(-21,1){\rotatebox{90}{\ $\sqrt{ \ \ }$}}
	    \put(-20,-3){\rotatebox{90}{  { \ \ \ \ \ \ \small transformation \ \ } }}

	\end{overpic}
	~
	\DeclareGraphicsExtensions{.png}
	\begin{overpic}[width=0.29\textwidth]{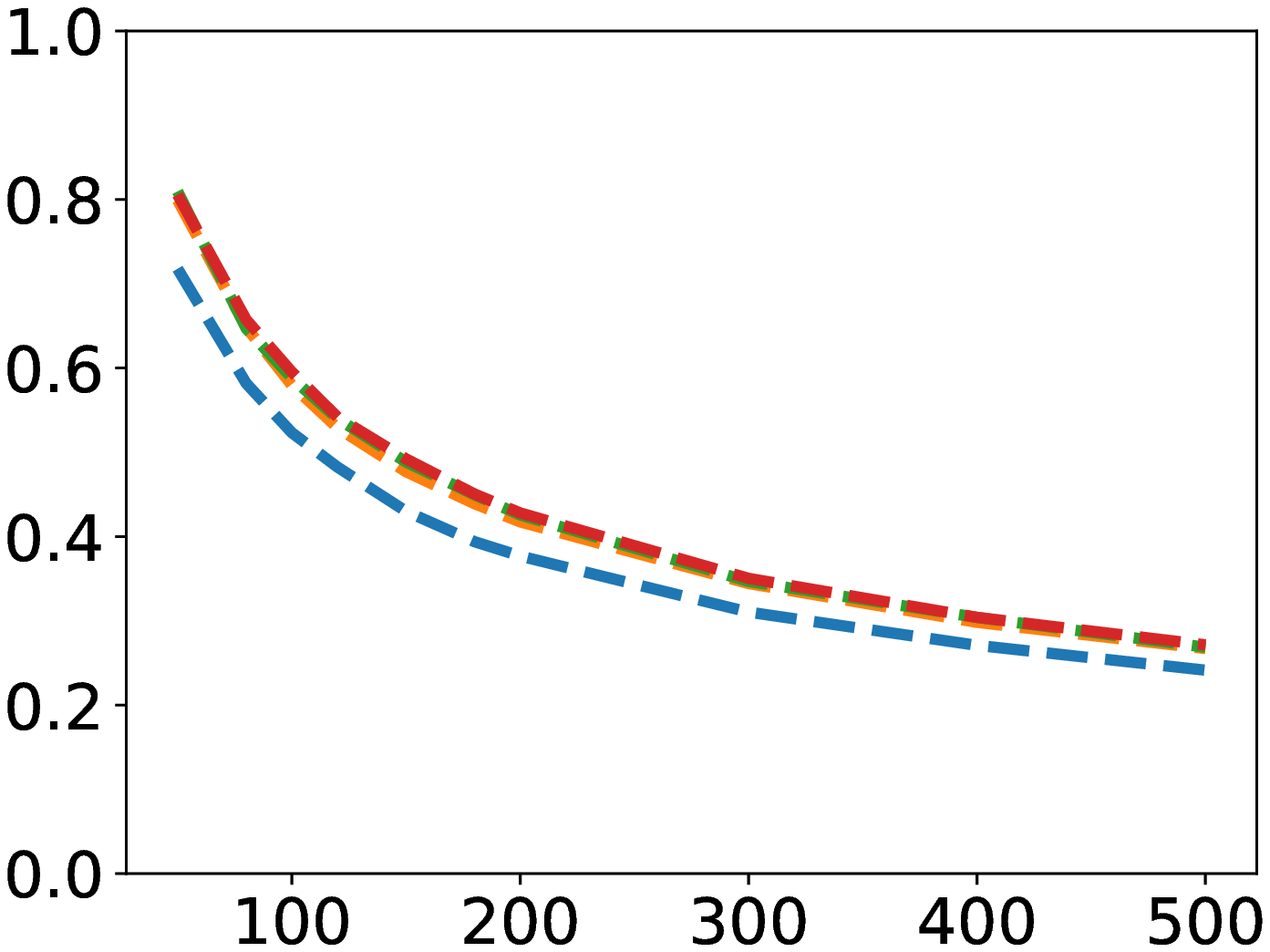} 
	\end{overpic}
	~	
	\begin{overpic}[width=0.29\textwidth]{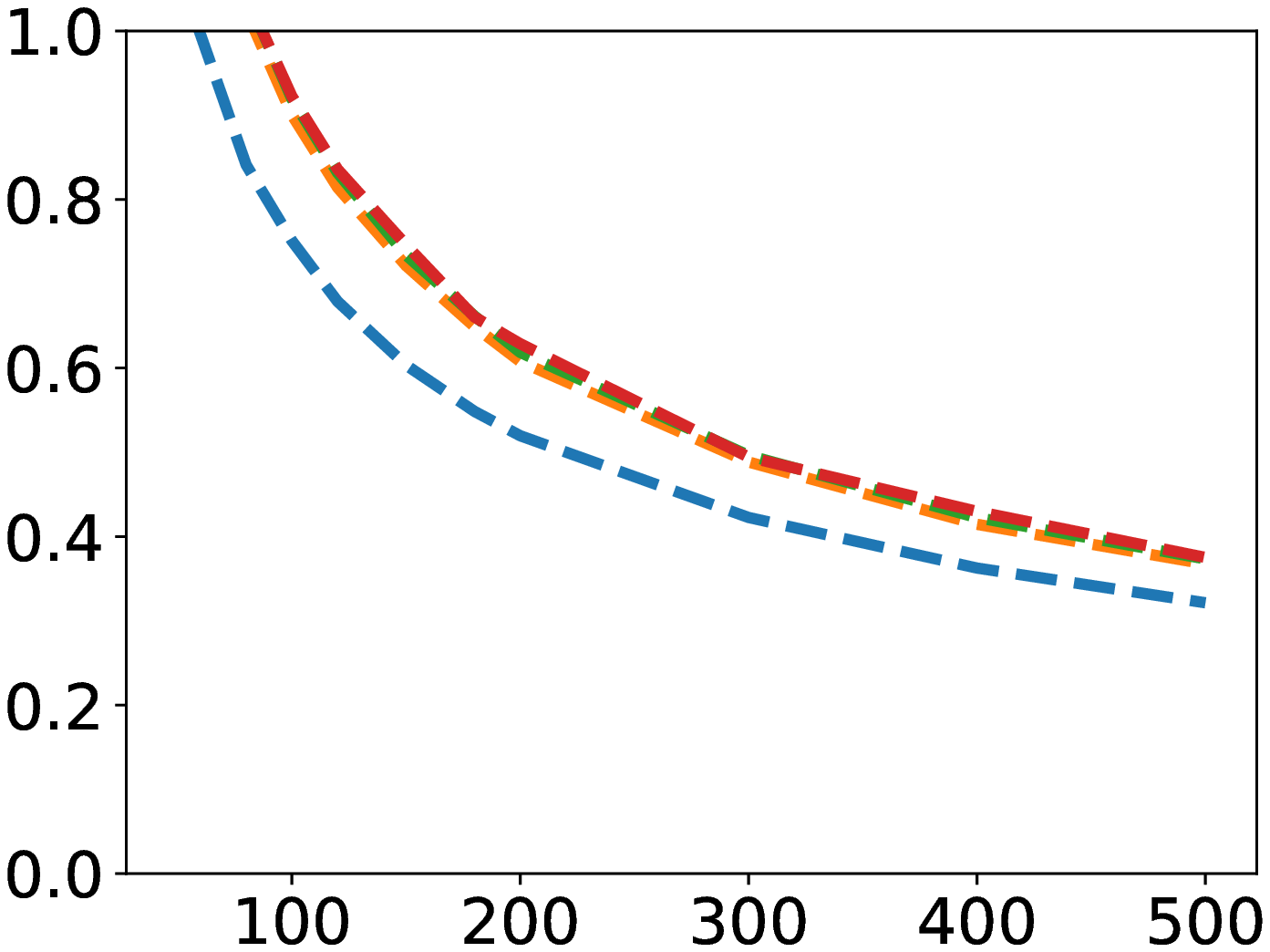} 
		 				
	\end{overpic}	
	\vspace{+.2cm}
	\caption{(Average width versus $n$ in simulation model (i) with an exponential decay profile). The plotting scheme is the same as described in the caption of Figure~\ref{fig5}, except that the three columns correspond to values of the eigenvalue decay parameter $\delta$.} 
	\label{fig7}
\end{figure}

\bibliography{reference}

\clearpage


\section*{Supplementary Material}
\noindent\textbf{Organization.} In Appendices \ref{appendix_A} and \ref{appendix_D}, we prove Theorems \ref{THM:ALL} and \ref{THM:TRANS_GENERAL} respectively. The building blocks for these results are Theorems~\ref{THM:GAUSS} and~\ref{THM:BOOT}, which are presented in Appendices~\ref{appendix_B} and~\ref{appendix_C} respectively. Technical lemmas are given in Appendix \ref{appendix_H}. Proposition~\ref{PROP1} from the main text is proved in Appendix \ref{appendix_G}. Background results are stated in Appendix~\ref{sec:background}. Appendix \ref{appendix_last} provides additional plots of simulation results. Appendix~\ref{sec:stocks} presents real-data examples based on stock market returns. Appendix~\ref{app:computation} describes the computational cost of implementing the bootstrap. Lastly, Appendix~\ref{app:assumption} provides sensitivity analysis with regard to Assumption 1(b).

\noindent\textbf{Conventions.} Throughout the proofs of Theorems~\ref{THM:ALL} and~\ref{THM:TRANS_GENERAL}, we may assume without loss of generality that $n\geq 3$, and that for any constant $\e\in (0,1)$ fixed with respect to $n$, the following inequality holds
\begin{equation} \label{n_large}
\frac{\log(n) \, \beta_{3q}^3 \, {\tt{r}} (\Sigma) }{ n^{1/2} } \ \leq \ \e,
\end{equation}
where $q=5 \log(kn)$. 
(If $n<3$ or if \eqref{n_large} does not hold, then the constant $c$ in the statements of Theorems~\ref{THM:ALL} and~\ref{THM:TRANS_GENERAL} may be taken as $c=3\vee \frac{1}{\e}$, which makes the results trivially true.) In addition, we will frequently re-use the symbol $c$ to denote a constant that does not depend on $n$, and we will allow its value to vary with each appearance.

\noindent \textbf{Notation.} The spectral decomposition for $\Sigma$ will be written as
\begin{equation}
    \Sigma = U\Lambda U\ttop,
\end{equation}
where $\Lambda=\text{diag}(\lambda_1(\Sigma),\dots,\lambda_p(\Sigma))\in\R^{p\times p}$, and the $j$th column of $U\in\R^{p\times p}$ is the $j$th eigenvector of $\Sigma$.  For a vector $v\in \R^p$, the $\ell_{\infty}$ and $\ell_2$-norms are denoted as $\| v \|_\infty =\max_{1\leq j\leq p}|v_j|$ and $\| v \|_2=(\sum_{j=1}^p v_j^2)^{1/2}$. For a $p_1 \times p_2$ matrix $M$, the following three norms will be used: $\| M\|_\op = \sup_{\| v\|_2 = 1} \| Mx \|_2$,  $\| M\|_1 = \max_{1 \leq j \leq p_2} \sum_{i = 1}^{p_1} |M_{ij}|$, and $\| M\|_\infty = \max_{1 \leq i \leq p_1} \sum_{j = 1}^{p_2} |M_{ij}|$. If $A$ is a symmetric $p\times p$ matrix and $j\in\{1,\dots,p\}$, then $\Lambda_j(A)$ denotes the $j\times j$ diagonal matrix formed by the largest $j$ eigenvalues of $A$,
$$\Lambda_j (A) = \text{diag} \big( \lambda_1(A), \dots, \lambda_j(A) \big).$$
Also, for each $j = 1, \dots, p$, let $d_j (A)=A_{jj}$ be the $j$th diagonal element of $A$, and let $\bdiag_j (A)$ be the vector of the first $j$ diagonal entries
$$ \bdiag_j (A) = (d_1 (A), \dots, d_j (A)) .$$
If $B$ is another symmetric $p\times p$ matrix, then the relation $B\succcurlyeq A$ means that $B-A$ is positive semidefinite.
The symbol $1_j$ denotes the $j$-dimensional all-ones vector. For a univariate scalar function $h$ with first derivative $h'$, define $\boldsymbol h'(v) = (h'(v_1), \dots, h'(v_p))$. Lastly, for two vectors $u$ and $v$, the symbol $u \odot v$ denotes the vector obtained from entrywise multiplication, $(u \odot v)_j = u_j v_j$.

\setcounter{section}{0}
\setcounter{equation}{0}
\def\theequation{S\arabic{section}.\arabic{equation}}
\def\thesection{S\arabic{section}}

\fontsize{12}{14pt plus.8pt minus .6pt}\selectfont

\section{Proof of Proposition~\ref{PROP1}} \label{appendix_G}

\noindent
\textbf{Proof of Proposition~\ref{PROP1}}(\emph{i})\textbf{.} For each $l=1,\dots, p$, let $\kappa_l=\E[Z_{1l}^4]$. Also recall that in part (\emph{i}) of the proposition, the entries of $Z_1$ are assumed to be independent. For any pair of indices $(j,j')$ satisfying $1 \leq j, j' \leq k$, it follows from the equation (9.8.6) in \cite{BaiS10} that the corresponding entry of $\Gamma\in\R^{k\times k}$ is 
\begin{align*}
\Gamma_{jj'}
& \ = \ \E \big[ \big( Z_1^\top (u_j u_j^\top) Z_1 - \tr(u_j u_j^\top) \big) \big(Z_1^\top (u_{j'} u_{j'}^\top) Z_1 - \tr(u_{j'} u_{j'}^\top) \big) \big] \\[0.2cm]
& \ = \ 2 \cdot 1 \{ j = j' \} +  \sum_{l = 1}^p (\kappa_l - 3) u_{j l }^2 u_{j'l}^2.
\end{align*}
If we let $H$ denote the $p \times k$ matrix whose $j$th column is equal to $(u_{j1}^2,\dots,u_{jp}^2)$, then the previous entrywise expression for $\Gamma$ can be written in matrix form as 
\[ \Gamma \ = \ 2 I_k + H^\top (D - 3 I_p) H ,\]
where $D = \text{diag} (\kappa_1, \dots, \kappa_p)$. Now recall the assumption that there is a constant $\kappa>1$ not depending on $n$ such that $\min_{1\leq l\leq p}\kappa_l\geq \kappa$. If we consider the case when $\kappa\geq 3$, then it is clear that $D-3I_p\succcurlyeq 0$, which implies
$\Gamma \ \succcurlyeq \  2I_k,$
and hence $\lambda_k(\Gamma)\geq 2$. \\

On the other hand, in the case when $\kappa<3$, we have
\begin{equation}\label{eqn:almostdone}
\begin{split}
\Gamma & \ \succcurlyeq \  2I_k-(3-\kappa)H\ttop H\\[0.2cm]
& \ \succcurlyeq \big(2-(3-\kappa)\|H\|_\op^2\big) I_k.
\end{split}
\end{equation}
With regard to the quantity $\|H\|_\op^2$, observe that
\begin{align*}
\| H \|_\op^2 
& \ \leq \ \| H \|_1 \| H \|_\infty \\[0.2cm]
& \ = \ \max_{1 \leq j \leq k} \| u_j \|_2^2 \cdot \max_{1 \leq i \leq p} \sum_{j = 1}^k u_{ji}^2 \\
& \ \leq \ 1.
\end{align*}
Combining this with \eqref{eqn:almostdone} gives
$\Gamma  \ \succcurlyeq \  (2-(3-\kappa))I_k,$
and hence $\lambda_k(\Gamma)\geq \kappa-1$, which completes the proof. \qed
~\\

\noindent
\textbf{Proof of Proposition~\ref{PROP1}}(\emph{ii})\textbf{.}
For $1 \leq j, j' \leq k$, it follows from Lemma A.1 in \cite{HuLLZ19} that 
\begin{align*}
\Gamma_{jj'}
& \ =\ \text{cov}\Big(Z_1\ttop(u_ju_j\ttop) Z_1\, , \, Z_1\ttop(u_{j'}u_{j'}\ttop)Z_1\Big)\\[0.2cm]
& \ =\ \frac{\E[\xi ^4]}{p(p+2)} \big( 1 + 2 \tr (u_j u_j^\top u_{j'} u_{j'}^\top) \big) - 1 \\[0.2cm]
& \ =\ 2 \cdot 1 \{j = j'\} \cdot \frac{\E[\xi ^4]}{p(p+2)} + \frac{\E[\xi ^4]}{p(p+2)} - 1.
\end{align*}
When written in matrix form, this is equivalent to
\[ \Gamma \ = a\,I_k + b\,1_k1_k\ttop \]
where  $a=\frac{2 \: \E[\xi ^4]}{p(p+2)}$  and $b= \frac{\E[\xi ^4]}{p(p+2)} - 1$.
Consequently, one eigenvalue of $\Gamma$ is equal to $a+kb$, and the rest are equal to $a$. Furthermore, due to the basic inequality $\E[\xi ^4]\geq (\E[\xi ^2])^2=p^2$,  and the fact that $p\geq 2$, we have the lower bounds
$$a \ \geq \ \ts\frac{2}{1+2/p}\geq \ 1$$
and
\begin{align*}
a+kb 
& \ \geq \ \ts \ \frac{2}{1+2/p}+ k\Big(\frac{1}{1+2/p}-1\Big)\\[0.2cm]
& \ = \ \ts\frac{2(1-k/p)}{1+2/p}\\[0.2cm]
& \ \geq \ \ts 1-\ts\frac{k}{p}\\
& \ \gtrsim \ 1,
\end{align*}
 which completes the proof.\qed

\section{Gaussian Approximation} \label{appendix_B}

\begin{theorem}[Gaussian approximation] \label{THM:GAUSS} 
Suppose that the conditions of Theorem~\ref{THM:ALL} hold and let $\zeta \sim N(0, \Lambda_k\Gamma\Lambda_k)$. Then, 
\begin{align}\label{eqn:gaussapproxdisp}
\sup_{t \in \R^k } \Big| \P \Big( \sqrt{n} \big( \blambda_k (\hat{\Sigma}) - \blambda_k (\Sigma) \big) \mylessthan \Big) - \P \Big( \zeta  \mylessthan \Big) \Big| \ \lesssim \ \frac{\log(n) \, \beta_{2q}^3 \,{\tt{r}}(\Sigma) }{n^{1/2}}.
\end{align}
\end{theorem}

\noindent
\Proof For each $i=1,\dots,n$, let $\tilde X_i=\Lambda^{1/2}U\ttop Z_i$, and let the associated sample covariance matrix be denoted as
\[ \tilde{\Sigma} \ = \ \frac{1}{n} \sum_{i = 1}^n \tilde{X}_i \tilde{X}_i^\top.\]
This implies $\tilde\Sigma=U\ttop \hat\Sigma U$, and so the eigenvalues of $\hat\Sigma$ may be equivalently written as $ \lambda_j (\hat{\Sigma}) = \lambda_j (\tilde{\Sigma})$ for every $j=1,\dots,p$. 
For future reference, it will also be helpful to note that $\E[\tilde\Sigma]=\Lambda$. \\

To partition $\tilde{\Sigma}$ into suitable blocks, we will write
\begin{align}
\tilde{\Sigma} \ = \ \Bigg(\!\!\!\begin{matrix} \tilde{\Sigma}[1,1] & \! \! \tilde{\Sigma}[1,2] \\ \ \ \tilde{\Sigma}[1,2]^\top & \!\! \tilde{\Sigma}[2,2] \end{matrix}\Bigg), \label{gauss_partition}
\end{align}
where the matrix $\tilde{\Sigma}[1,1]$ is of size $k\times k$, and the matrix $\tilde{\Sigma}[2,2]$ is of size $(p-k)\times (p-k)$. \\

Based on the notation above, the desired result can be broken down as follows:
\[ \sup_{t \in \R^k } \Big| \P \Big( \sqrt{n} \big( \blambda_k (\hat{\Sigma}) - \blambda_k (\Sigma) \big) \mylessthan \Big) - \P \Big( \zeta  \mylessthan \Big) \Big| \ \leq \  \textup{I}+\textup{II} \]
where the two terms on the right are defined as  
\begin{align}
\textup{I} 
\ & = \ \sup_{t \in \R^k } \Big| \P \Big( \sqrt{n} \big( \blambda_k (\tilde{\Sigma}) - \blambda_k (\Sigma) \big) \mylessthan \Big) -  \P \Big( \sqrt{n} \big( \bdiag_k (\tilde{\Sigma}[1,1]) - \blambda_k (\Sigma) \big) \mylessthan \Big) \Big|,\label{eqn:Idef} \\[0.2cm]
\textup{II}
\ & = \ \sup_{t \in \R^k } \Big| \P \Big( \sqrt{n} \big( \bdiag_k (\tilde{\Sigma}[1,1]) - \blambda_k (\Sigma) \big) \mylessthan \Big) -  \P \Big( \zeta  \mylessthan \Big)\Big|.\label{eqn:IIdef}
\end{align}

For each of these terms, Lemmas \ref{lemma:gauss4}, 
and \ref{lemma:gauss3} respectively give the following bounds
\begin{align*}
 \textup{I} & \ \lesssim \ \frac{\log(n) \, \beta_{2q}^3 \, {\tt{r}} (\Sigma) }{ n^{1/2} }\\[0.2cm]
\textup{II} & \ \lesssim \ \frac{\beta_4^3}{n^{1/2}},
\end{align*}
completing the proof. \qed

\subsection{Lemmas for Gaussian approximation}

\begin{lemma} \label{lemma:gauss1}
Suppose that the conditions of Theorem \ref{THM:GAUSS} hold. Then, there is a constant $c>0$ not depending on $n$ such that the event
\[ \Big\| \sqrt n \Big(\blambda_k (\tilde{\Sigma}) - \blambda_k ( \tilde{\Sigma}[1,1])\Big) \Big\|_{\infty}
\ \leq \ \frac{c \, \log (n) \, \beta_{2q} \, \lambda_1 (\Sigma) \, {\tt{r}}(\Sigma) }{n^{1/2}} \]
holds with probability at least $1 - \frac{c}{n^4}$.
\end{lemma}

\noindent
\Proof Recall the partition \eqref{gauss_partition} of the matrix $\tilde\Sigma$ in the proof of Theorem \ref{THM:GAUSS}. By Wielandt's inequality (Lemma \ref{background0}), the following inequality holds when the event $\{\lambda_k (\tilde{\Sigma}[1,1]) > \lambda_1 (\tilde{\Sigma}[2,2])\}$ occurs,
\begin{equation}
 \max_{1 \leq j \leq k} \big| \lambda_j (\tilde{\Sigma}) - \lambda_j ( \tilde{\Sigma}[1,1]) \big| 
 \ \leq \ \frac{\| \tilde{\Sigma}[1,2] \|_\op^2}{\lambda_k (\tilde{\Sigma}[1,1]) - \lambda_1 (\tilde{\Sigma}[2,2])} .\label{eqn:wielandt1}
\end{equation}
We will derive a high-probability upper bound for the right side of~\eqref{eqn:wielandt1} by separately handling the numerator and denominator.\\

To control the numerator in the bound~\eqref{eqn:wielandt1}, we may apply Lemma \ref{tl1} to conclude that 
\[  \big\| \| \tilde{\Sigma}[1,2] \|_\emph{\op}^2 \big\|_q 
 \ \lesssim \ \frac{q \, \beta_{2q} \, \lambda_1(\Sigma) \, \tr (\Sigma)}{ n ^{1 - 3/(2q)}} .\]
Then, for any $t > 0$, Chebyshev's inequality yields the tail bound
\[ \P \big( \| \tilde{\Sigma}[1,2] \|_\op^2 \geq e\,t \big) \ \leq \ \frac{ e^{-q} \ \big\| \| \tilde{\Sigma}[1,2] \|_\op^2 \big\|_q^q}{ t^q} .\]
Recalling the choice $q = 5 \log (kn)$ and taking $t = \frac{c}{n}\log(n) \beta_{2q} \lambda_1(\Sigma) \tr (\Sigma) $ for a sufficiently large constant $c$, it follows that the event
\begin{equation} \label{gauss1:bound1}
\| \tilde{\Sigma}[1,2] \|_\op^2  \ \leq \  \frac{c \, \log (n) \, \beta_{2q} \, \lambda_1(\Sigma)^2 \, {\tt{r}}(\Sigma) }{n} 
\end{equation}
holds with probability at least $1 - \frac{1}{n^4}$. \\

To handle the denominator in the bound~\eqref{eqn:wielandt1}, let $\Pi\in\R^{k\times p}$ denote the matrix whose $i$th row is the $i$th standard basis vector in $\R^p$. Then, Weyl's inequality implies
\begin{align*}
\max_{1 \leq j \leq k} |\lambda_j (\tilde{\Sigma}[1,1]) - \lambda_j (\Sigma)| 
& \ \leq \ \big\| \Pi (\tilde{\Sigma} - \Lambda) \Pi^\top \big\|_\op \\[0.2cm]
& \ \leq \ \| \tilde{\Sigma} - \Lambda \|_\op.
\end{align*}
Similarly, we have
\[ \max_{k+1 \leq j \leq p} |\lambda_{j-k} (\tilde{\Sigma}[2,2]) - \lambda_j(\Sigma)| \ \leq \ \| \tilde{\Sigma} - \Lambda \|_\op .\]
Combining the last two steps yields the following bound for some positive constant $c_1$ not depending on $n$, 
\begin{align} \label{eig_gap_lower_bound}
\lambda_k (\tilde{\Sigma}[1,1]) - \lambda_1 (\tilde{\Sigma}[2,2])
& \ = \ \big( \lambda_k (\tilde{\Sigma}[1,1]) - \lambda_k(\Sigma) \big) + \big( \lambda_k(\Sigma)  - \lambda_{k+1}(\Sigma) \big) - \big( \lambda_1 (\tilde{\Sigma}[2,2]) - \lambda_{k+1}(\Sigma) \big) \notag  \\[0.2cm]
& \ \geq \ c_1 \lambda_1(\Sigma) - 2 \| \tilde{\Sigma} - \Lambda \|_\op,
\end{align}
where Assumption \ref{model_assumptions}.\ref{A2} has been used in the last line.\\

To complete the proof, it suffices to show that $\|\tilde \Sigma-\Lambda\|_\op\leq \frac{c_1}{4}\lambda_1(\Sigma)$ holds with high probability. This may be accomplished using Lemma \ref{background2}, which gives the following bound,

\begin{equation}\label{eqn:complicated}
\big(\E \| \tilde{\Sigma} - \Lambda \|_\op^q \big)^{1/q} \ \lesssim \ \Bigg( \frac{ \sqrt{q} \, \big(\E \| \tilde{X}_1 \|_2^{2q} \big)^{\frac{1}{2q} } \, \big\| \E [\tilde{X}_1 \tilde{X}_1^\top] \big\|_\op^{ 1/2 } }{ n^{1/2 - 3/(2q)} }\Bigg) \, \bigvee \,\bigg( \frac{q \, \big(\E \| \tilde{X}_1 \|_2^{2q} \big)^{1/q}}{n^{1-3/q}} \bigg).
\end{equation}
This bound can be simplified by noting that $\E [\tilde{X}_1 \tilde{X}_1^\top]  = \Lambda$ and
\begin{equation}
\begin{split}
 \big(\E \| \tilde{X}_1 \|_2^{2q}\big)^{1/q}  & \ = \  \big\| \|\tilde X_1\|_2^2\big\|_q  \\[0.2cm]
& \ = \  \bigg\| \sum_{j = 1}^p \lambda_j \langle u_j, Z_1 \rangle^2 \bigg\|_q\\[0.2cm]
& \   \leq \ \tr (\Sigma) \, \beta_q .
\end{split}
\end{equation}
Therefore, the bound~\eqref{eqn:complicated} reduces to
\begin{align} \label{gauss1:bound_interm}
\big(\E \| \tilde{\Sigma} - \Lambda \|_\op^q \big)^{1/q} 
& \ \lesssim \ 
\sqrt{ \frac{q \, \beta_q \, \lambda_1(\Sigma)^2 \, {\tt{r}} (\Sigma)}{n^{1 - 3/q} } } \bigvee \frac{q \, \beta_q \, \lambda_1(\Sigma) \, {\tt{r}}(\Sigma)}{n^{1-3/q}} \notag \\[0.2cm]
& \ \lesssim \ 
\sqrt{ \frac{q \, \beta_q \, \lambda_1(\Sigma)^2 \, {\tt{r}} (\Sigma)}{n^{1 - 3/q} } },
\end{align}
where the second line has used the condition~\eqref{n_large}. \\

To apply the previous bound, Chebyshev's inequality implies that for any $t>0$,
\[ \P \big( \| \tilde{\Sigma} - \Lambda \|_\op \geq e \ t \big)
\ \leq \
\frac{e^{-q} \ \E \| \tilde{\Sigma} - \Lambda \|_\op^q}{t^q} .\]
Recalling the choice $q = 5 \log (kn)$ and taking $t = \frac{c\lambda_1 (\Sigma)}{\sqrt n} \sqrt{\log (n)\beta_q {\tt{r}} (\Sigma)}$ for a sufficiently large constant $c$ not depending on $n$, we conclude that the event
\begin{equation} \label{gauss1:bound2}
\| \tilde{\Sigma} - \Lambda \|_\op 
\ \leq \ 
c \lambda_1 (\Sigma)\sqrt{ \frac{\log(n)\, \beta_q \, {\tt{r}} (\Sigma)}{n} } 
\end{equation}
holds with probability at least $1 - \frac{1}{n^4}$. Furthermore, if we apply this bound to~\eqref{eig_gap_lower_bound} and use the condition~\eqref{n_large}, then  the bound
\begin{align} \label{gauss1:bound3}
\lambda_k (\tilde{\Sigma}[1,1]) - \lambda_1 (\tilde\Sigma[2,2])
& \ \geq \ \frac{c_1}{2} \lambda_1 (\Sigma)
\end{align}
holds with probability at least $1 - \frac{1}{n^4}$. \\

The stated result is obtained by combining the bounds~\eqref{eqn:wielandt1}, \eqref{gauss1:bound1}, and \eqref{gauss1:bound3}. \qed

~\\
The next result shows that the random vector $\sqrt{n} \big( \blambda_k (\tilde{\Sigma}[1,1]) - \blambda_k(\Sigma) \big)$ is well approximated by $\sqrt{n} \big( \bdiag_k (\tilde{\Sigma}[1,1]) - \blambda_k(\Sigma) \big)$ in an entrywise sense.
~\\

\begin{lemma} \label{lemma:gauss2}
Suppose that the conditions of Theorem \ref{THM:GAUSS} hold. Then, there is a constant $c > 0$ not depending on $n$ such that the event
\begin{align*} \label{eqn:lemma5}
\Big\| \sqrt{n} \Big( \blambda_k (\tilde{\Sigma}[1,1]) - \bdiag_k (\tilde{\Sigma}[1,1]) \Big) \Big\|_\infty
\ \leq \ \frac{c \, \log(n) \, \beta_{2q} \, \lambda_1 (\Sigma) \, {\tt{r}} (\Sigma) }{n^{1/2}}
\end{align*}
holds with probability at least $1 - \frac{c}{n^4}$.
\end{lemma}

\noindent
\Proof We will use an iterative argument. As an initial step, first partition the matrix $\tilde{\Sigma}[1,1]$ as
\begin{align*}
\tilde{\Sigma}[1,1]
= \begin{pmatrix}
d_1 (\tilde{\Sigma}[1,1]) & V_1 \\
V_1^\top & D_1
\end{pmatrix},
\end{align*}
where $d_1 (\tilde{\Sigma}[1,1])$ is a scalar, and $D_1$ is a $(k-1) \times (k-1)$ matrix. To lighten the notational burden of subscripts when handling matrices of different sizes, we will write $\blambda(A)=(\lambda_1(A),\dots,\lambda_r(A))$, as well as $\bdiag(A)=(A_{11},\dots,A_{rr})$ for any symmetric matrix $A\in\R^{r\times r}$ and integer $r\geq 1$. By the triangle inequality, we have
\begin{equation}\label{eqn:deltadeltaprime}
\small
\begin{split}
\Big\| \blambda(\tilde{\Sigma}[1,1]) - \bdiag(\tilde{\Sigma}[1,1])  \Big\|_\infty
& \ \leq \
\bigg\| \blambda(\tilde{\Sigma}[1,1]) -  \begin{bmatrix} d_1(\tilde{\Sigma}[1,1])\\ \blambda(D_1)  \end{bmatrix} \bigg\|_\infty 
+ \
\bigg\| \begin{bmatrix} d_1(\tilde{\Sigma}[1,1])\\ \blambda(D_1)  \end{bmatrix} - \bdiag (\tilde{\Sigma}[1,1]) \bigg\|_\infty ,\\[0.2cm]
& \ = \ \textup{T} + \textup{T}' 
\end{split}
\end{equation}
where we have defined the random variables $\textup{T}$ and $\textup{T}'$ through the last line.
Later on, we will show that the event
\begin{equation} \label{eqn:gauss:lemma1}
\textup{T} \ \leq \ \frac{c \,\log (n) \, \beta_{2q} \,\tr (\Sigma)  }{\, n }
\end{equation}
holds with probability at least $1 - \frac{2}{kn^4}$. \\

Next, to handle $\textup{T}'$, we partition the matrix $D_1$ as
\begin{align*}
D_1
= \begin{pmatrix}
d_1 (D_1) & V_2 \\
V_2^\top & D_2
\end{pmatrix},
\end{align*}
where $d_1 (D_1) = d_2 (\tilde{\Sigma}[1,1])$ is a scalar, and $D_2$ is a $(k-2) \times (k-2)$ matrix.  
Proceeding in a similar manner to~\eqref{eqn:deltadeltaprime}, and noting that $d_1(\tilde\Sigma[1,1])$ is the same as the first entry of $\boldsymbol d(\tilde \Sigma[1,1])$, we have
%
\begin{equation}
\small
\begin{split}
\textup{T}' & \ \ \leq \ \
\bigg\| \blambda(D_1)  - \begin{bmatrix} d_1(D_1)\\ \blambda(D_2)  \end{bmatrix}\bigg\|_\infty
 \ \ + \ \ \
\bigg\|\begin{bmatrix} d_1(D_1)\\ \blambda(D_2)  \end{bmatrix} - \bdiag(D_1) 
\bigg\|_\infty,\\[0.2cm]
& \ \ = \ \ \textup{T}'' \ + \ \textup{T}''',
\end{split}
\end{equation}
where the random variables $\textup{T}''$ and $\textup{T}'''$ have been defined through the last line.
The argument that will be used to prove~\eqref{eqn:gauss:lemma1} can also be used to show that the event
\begin{equation*} 
\textup{T}''
\ \ \leq \ \ \frac{ c \,\log (n) \, \beta_{2q} \,  \tr (\Sigma) }{ n }
\end{equation*}
holds with probability at least $1 - \frac{c}{kn^4}$. Likewise, we can combine $k-1$ iterations of this process by summing the bounds on $\textup{T}$, $\textup{T}'',\textup{T}'''',\dots$ appearing at each iteration. \\

To complete the proof, it remains to validate the claim \eqref{eqn:gauss:lemma1}. Let $D_0 = \tilde{\Sigma}[1,1]$, and observe that
\begin{equation}\label{eqn:deltaeq}
\textup{T} \ = \ |\lambda_1 (D_{0}) - d_1 (D_{0})| \ \vee \max_{2 \leq j \leq k} | \lambda_j (D_0) -\lambda_{j-1} (D_1) |.
\end{equation}
Wielandt's inequality (Lemma~\ref{background0}) gives the following bounds when the event $\{d_1(D_{0})>\lambda_1(D_1)\}$ occurs,
$$
|\lambda_1 (D_{0}) - d_1 (D_{0})| \ \leq \ \frac{\| V_1\|_\op^2}{d_1 (D_{0}) - \lambda_1 (D_1)} 
$$
and 
$$
|\lambda_{j} (D_0) - \lambda_{j-1} (D_1) | \ \leq \ \frac{\| V_1 \|_\op^2}{d_1 (D_{0}) - \lambda_{1} (D_1)}, 
$$
for all $j=2,\dots,k$. So, in light of the formula for $\textup{T}$ given in~\eqref{eqn:deltaeq}, we conclude that the bound
\begin{equation}\label{eqn:V1bound}
\textup{T} \ \leq \ \frac{\| V_1 \|_\op^2}{d_1 (D_{0}) - \lambda_{1} (D_1)}.
\end{equation}
holds whenever the event $\{d_1 (D_{0}) > \lambda_{1} (D_1)\}$ holds. By using the argument at \eqref{gauss1:bound3} from the proof of Lemma~\ref{lemma:gauss1}, it can be shown that the denominator in the previous bound satisfies
$$d_1 (D_{0}) - \lambda_{1} (D_1) \geq c\lambda_1(\Sigma)$$
with probability at least $1-\frac{1}{kn^4}$ for some constant $c>0$ not depending on $n$. Next, in order to derive an upper bound on $\|V_1\|_\op$, note that for any value of $k$, the matrix  $V_1\ttop$ is contained in the submatrix of $\tilde \Sigma$ indexed by $\{2,\dots,p\}\times \{1\}$. Hence, $\|V_1\|_\op$ is upper bounded by the operator norm of that submatrix, which is the same as $\tilde\Sigma[1,2]\ttop$ in the particular case when $k=1$. Due to this observation, it follows from Lemma~\ref{tl1} that there is a constant $c>0$ not depending on $n$, such that for any choice of $k$, the bound
\begin{equation}\label{eqn:V1nextbound}
\|V_1\|_\op^2 \leq \frac{c \,q \, \beta_{2q} \, \lambda_1 (\Sigma) \tr (\Sigma) }{n}
\end{equation}
holds with probability at least $1-\frac{1}{kn^4}$, where we continue to use $q=5 \log(kn)$. Thus, combining the last few steps establishes the claim~\eqref{eqn:gauss:lemma1}. \\

To comment on how this argument can be applied iteratively to $\textup{T}''$ and its successors, the previous reasoning involving Weilandt's inequality shows that the bound
$$\textup{T}'' \ \leq \ \frac{\|V_2\|_\op^2}{d_1(D_1)-\lambda_1(D_2)}$$
holds whenever the event $\{d_1(D_1)>\lambda_1(D_2)\}$ holds.
In turn, a lower bound on the denominator $d_1(D_1)-\lambda_1(D_2)$ that is proportional to $\lambda_1(\Sigma)$ can be established in the same manner as for $d_1(D_0)-\lambda_1(D_1)$. Meanwhile, the numerator can be handled by noting that for any value of $k$, the matrix $V_2\ttop$ is contained in the submatrix of $\tilde\Sigma$ indexed by $\{3,\dots,p\}\times\{1,2\}$. The latter matrix is the same as $\tilde\Sigma[1,2]\ttop$ in the particular case of $k=2$, and consequently, Lemma~\ref{tl1} can be used to establish a bound on $\|V_2\|_\op^2$ that is of the same form as~\eqref{eqn:V1nextbound}.
This completes the proof.\qed

~\\
The next lemma provides a Gaussian approximation result for $\sqrt{n} \big( \bdiag_k (\tilde{\Sigma}[1,1]) - \blambda_k(\Sigma) \big)$. 
~\\

\begin{lemma}\label{lemma:gauss3}
Suppose that the conditions of Theorem \ref{THM:GAUSS} hold, and let $\textup{II}$ be as defined in~\eqref{eqn:IIdef}. Then, 
\[ 
\textup{II} \ \lesssim \ \frac{ \beta_4^3}{n^{1/2} } .\]
\end{lemma}

\noindent
\Proof For each $i=1,\dots,n$, let the random vector $W_i\in\R^k$ have its $j$th entry defined as $W_{ij} = \langle u_j, Z_i\rangle^2-1$, where $u_j$ is the $j$th eigenvector of $\Sigma$. Letting $Y_i = \Lambda_k W_i$, we have $\sqrt{n} \big( \bdiag_k (\tilde{\Sigma}[1,1]) - \blambda_k(\Sigma) \big) \ = \ \frac{1}{\sqrt{n}} \sum_{i = 1}^n Y_i$. Also, note that the covariance matrix of $Y_i$ is given by $\E[Y_iY_i\ttop] = \Lambda_k \Gamma \Lambda_k$, with $\Gamma$ as defined in Assumption~\ref{model_assumptions}.\ref{A3}. \\

Applying Bentkus' multivariate Berry-Esseen theorem (Lemma \ref{background5}) yields
\begin{equation}\label{eqn:berry}
 \sup_{t\in\R^k} \bigg| \P \Big( \ts\frac{1}{\sqrt{n}} \ts\sum_{i = 1}^n Y_i \preceq t\Big) - \P \big( \zeta  \preceq t \big) \bigg| \ \lesssim \displaystyle \ \frac{  \E \big\| (\Lambda_k \Gamma \Lambda_k)^{-1/2} Y_1 \big\|_2^3}{n^{1/2} }.
 \end{equation}
The proof is complete once we derive a bound on $\E \big\| (\Lambda_k \Gamma \Lambda_k)^{-1/2} Y_1 \big\|_2^3$. 
Due to Assumption \ref{model_assumptions}.\ref{A3}, we have
\[ \big\| (\Lambda_k \Gamma \Lambda_k)^{-1/2} Y_1 \big\|_2^2 \ = \ W_1^\top \Gamma^{-1} W_1 \ \leq \ c \| W_1 \|_2^2 \ ,\]
almost surely for some constant $c>0$ not depending on $n$.
In turn, Lyapunov's inequality implies
\begin{align*}
\E \big\| (\Lambda_k \Gamma \Lambda_k)^{-1/2} Y_1 \big\|_2^3
\ & \ \lesssim \   \E \big[(W_1^\top W_1)^2 \big]^{3/4} \\[0.2cm]
& \ = \ \bigg\|\sum_{j = 1}^k \big( \langle u_j, Z_1 \rangle^2 - 1 \big)^2 \bigg\|_2^{3/2} \\[0.2cm]
& \ \leq \ \bigg(\sum_{j = 1}^k \big\| \langle u_j, Z_1 \rangle^2 - 1 \big)^2 \big\|_2\bigg)^{3/2} \\[0.2cm]
& \ \lesssim \ \beta_4^3. 
\end{align*}
Substituting this bound into~\eqref{eqn:berry} completes the proof. \qed

\begin{lemma}\label{lemma:gauss4} 
Suppose that the conditions of Theorem \ref{THM:GAUSS} hold, and let $\textup{I}$ be as defined in~\eqref{eqn:Idef}. Then, 
\begin{equation}\label{eqn:lambdatod}
\textup{I} \ \lesssim \ 
\frac{\log(n) \, \beta_{2q}^3 \, {\tt{r}} (\Sigma)}{n^{1/2}} .
\end{equation}
\end{lemma}
\noindent
\Proof
For any $\epsilon > 0$ and  $t \in \R^k$, define the two events
\begin{align*}
\A (t) & = \Big\{ \sqrt{n} \big( \bdiag_k (\tilde{\Sigma}[1,1]) - \blambda_k(\Sigma) \big) \mylessthan  \Big\}, \\[0.3cm]
\mathcal{B} (\epsilon) & = \Big\{ \sqrt{n} \big\| \blambda_k(\tilde{\Sigma}) - \bdiag_k (\tilde{\Sigma}[1,1]) \|_{\infty} \geq \epsilon  \Big\}.
\end{align*}

It follows from Lemma~\ref{background7} that 
\begin{equation}\label{eqn:Bcomp}
\begin{split}
\textup{I} \ \leq \ \P \big( \mathcal{B} (\epsilon) \big)  \ + \  \sup_{t \in \R^k} \Big| \P \big( \A (t + \epsilon \, 1_k) \big) - \P \big( \A (t - \epsilon \, 1_k) \big) \Big|.
\end{split}
\end{equation}
The first term $\P(\mathcal{B}(\epsilon))$ can be handled by Lemmas~\ref{lemma:gauss1} and~\ref{lemma:gauss2}, which show that there is a constant $c>0$ not depending on $n$ such that  if $\epsilon = \frac{c}{\sqrt n} \, \log(n) \, \beta_{2q} \, \lambda_1 (\Sigma) \, {\tt{r}} (\Sigma)$, then 
$$ \P(\mathcal{B}(\epsilon)) \ \lesssim  \ \ts \frac{1}{n}. $$
Next, the anti-concentration term in~\eqref{eqn:Bcomp} can be bounded through an approximation involving the Gaussian vector $\zeta\sim N(0,\Lambda_k\Gamma\Lambda_k)$,
\begin{align*}
 \ \quad \, \sup_{t \in \R^k} \bigg| \P \big( \A (t + \epsilon 1_k) \big) - \P \big( \A (t - \epsilon 1_k) \big) \bigg|
 & \ \leq \
 2 \, \textup{II} \ + \ \sup_{t \in \R^k} \bigg| \P \Big( \zeta  \mylessthan + \epsilon 1_k \Big) - \P \big(\zeta  \mylessthan - \epsilon 1_k \big) \bigg| \\[0.3cm]
& \ = \ 2 \, \textup{II}\ + \ \textup{J}(\epsilon),
\end{align*}
where the quantity  $\textup{J}(\epsilon)$ is defined by the last line, and $\textup{II}\lesssim \beta_4^3/n^{1/2}$ holds by Lemma~\ref{lemma:gauss3}.
To handle $\textup{J}(\epsilon)$, we need the following lower bound
\begin{align*}
\min_{1 \leq j \leq k} \var(\zeta_j) \ 
 = \ \min_{1 \leq j \leq k} \lambda_j(\Sigma)^2 \, \Gamma_{jj} \
 \gtrsim \ \lambda_k(\Sigma)^2,
\end{align*}
where we have used Assumption~\ref{model_assumptions}.\ref{A3}. Based on this lower bound, Nazarov's inequality (Lemma \ref{background4}) yields
\begin{equation}\label{eqn:NazarovJ2}
\begin{split}
 \textup{J}(\epsilon) & \ \lesssim \ \frac{\epsilon}{\lambda_k(\Sigma)}\\[0.2cm]
& \lesssim \ \frac{ \log(n) \, \beta_{2q} \, {\tt{r}} (\Sigma) \, \lambda_1(\Sigma) /\lambda_k (\Sigma) }{n^{1/2} } .
\end{split}
\end{equation}
Combining the previous bounds and noting that Assumption~\ref{model_assumptions}.\ref{A2} implies $\lambda_1(\Sigma)/\lambda_k(\Sigma)\lesssim 1$, we obtain the stated result. \qed
~\\

\section{Bootstrap Approximation} \label{appendix_C}

To introduce some notation, first note that the bootstrapped vectors $X_i^\star$ can be represented theoretically as $X_i^{\star}=\Sigma^{1/2}Z_i^{\star}$ for each $i=1,\dots,n$, where $Z_1^{\star},\dots,Z_n^{\star}$ are sampled with replacement from $Z_1,\dots,Z_n$. Also, let $\hat\Sigma^{\star} = \frac{1}{n} \sum_{i = 1}^n X_i^\star (X_i^\star)^\top$, and let the diagonal matrix of the largest $k$ sample eigenvalues be denoted as 
\[ \hat\Lambda_k = \text{diag} \big( \lambda_1 (\hat\Sigma), \dots, \lambda_k (\hat\Sigma) \big). \] 
For each $i=1,\dots,n$, recall the vector $W_i\in\R^k$ whose $i$th entry is defined as $W_{ij} = \langle u_j, Z_i\rangle^2-1$, where $u_j$ is the $j$th eigenvector of $\Sigma$. Also let $\bar W=\frac{1}{n}\sum_{i=1}^n W_i$. In light of the fact that $\Gamma=\E[W_1W_1\ttop]$ with $\E[W_1]=0$, the empirical counterpart of $\Gamma$ is defined as
\begin{equation}\label{eqn:hatgammadef}
\hat\Gamma \ = \ \frac{1}{n}\sum_{i=1}^n (W_i-\bar W)(W_i-\bar W)\ttop.
\end{equation}
Lastly, recall that $\P(\cdot|X)$ and $\E[\cdot|X]$ refer to probability and expectation that are conditional on $X_1,\dots,X_n$. 
\begin{theorem}[Bootstrap approximation] \label{THM:BOOT}\!\!\!Suppose that Assumption \ref{model_assumptions} holds, and let \smash{$q = 5 \log (kn)$}. Also, let $\xi\in\R^k $ be a random vector that is conditionally distributed as $N \big(0, \Lambda_k \hat\Gamma \Lambda_k \big)$, given the observations $X_1,\dots,X_n$. Then, there is a constant $c > 0$ not depending on $n$ such that the event
\begin{align}\label{eqn:bootapproxdisp}
\sup_{t \in \R^k } \Big| \P \Big( \sqrt{n} \big( \blambda_k (\hat{\Sigma}^\star) - \blambda_k (\hat\Sigma) \big) \mylessthan \, \Big| \, X \Big) - \P \Big( \xi  \mylessthan \, \Big| \, X \Big) \Big| \ \leq \ \frac{ c\, \log(n) \,  \beta_{3q}^3  \,{\tt{r}}(\Sigma)  }{n^{1/2}} 
\end{align}
holds with probability at least $1 - \frac{c}{n}$. 
\end{theorem}

\noindent
\Proof For each $i=1,\dots,n$, define 
$$\tilde{X}_i^\star = \Lambda^{1/2} U^\top Z_i^\star,$$
as well as the matrix
\[ \tilde{\Sigma}^\star  = \frac{1}{n} \sum_{i = 1}^n (\Lambda^{1/2} U^\top Z_i^\star) (\Lambda^{1/2} U^\top Z_i^\star)^\top .\]
Based on this definition, we have $\tilde \Sigma^{\star}=U\ttop \hat\Sigma^{\star}U$, and hence
\[  \lambda_j ( \tilde{\Sigma}^\star)=\lambda_j ( \hat\Sigma^{\star} )\]
for all $j = 1,\dots,p$.
By analogy with the proof of Theorem \ref{THM:GAUSS}, we partition $\tilde\Sigma^{\star}$ as
\begin{align} \label{boot_partition}
\tilde{\Sigma}^\star = \begin{pmatrix}\!\! \tilde{\Sigma}^\star[1,1] & \tilde{\Sigma}^\star[1,2] \\ \ \tilde{\Sigma}^\star[1,2]^\top & \tilde{\Sigma}^\star[2,2] \end{pmatrix},
\end{align}
where the matrix $\tilde{\Sigma}^\star[1,1]$ is of size $k \times k$, and the matrix $\tilde{\Sigma}^\star[2,2]$ is of size $(p-k) \times (p-k)$. \\

To proceed, consider the bound
\[ \sup_{t \in \R^k } \Big| \P \Big( \sqrt{n} \big( \blambda_k (\hat{\Sigma}^\star) - \blambda_k (\hat\Sigma) \big) \mylessthan \, \Big| \, X \Big) - \P \big( \xi  \mylessthan \, \big| \, X \big) \Big| \   \leq \  \hat{\textup{I}} \ + \ \hat{\textup{II}}  \]
where the terms on the right are defined as  
\begin{align}
\hat{\textup{I}} & \ = \ \sup_{t \in \R^k } \Big| \P \Big( \sqrt{n} \big( \blambda_k (\tilde{\Sigma}^\star) - \blambda_k (\hat\Sigma) \big) \mylessthan \, \Big| \, X \Big) -  \P \Big( \sqrt{n} \big( \bdiag_k (\tilde{\Sigma}^\star[1,1]) - \bdiag_k (\tilde{\Sigma}[1,1]) \big) \mylessthan \, \Big| \, X \Big) \Big|,\label{eqn:hatIdef} \\[0.3cm]
\hat{\textup{II}} & \ = \ \sup_{t \in \R^k } \Big| \P \Big( \sqrt{n} \big( \bdiag_k (\tilde{\Sigma}^\star[1,1]) - \bdiag_k (\tilde{\Sigma}[1,1]) \big) \mylessthan \, \Big| \, X \Big) -  \P \Big( \xi  \mylessthan \, \Big| \, X \Big)\Big|.\label{eqn:hatIIdef}
\end{align}
Lemmas \ref{lemma:boot5} and \ref{lemma:boot3} ensure that there is a constant $c>0$ not depending on $n$ such that the following events
\begin{align}
\hat{\textup{I}} & \ \leq \ \frac{c \, \log(n) \, \beta_{3q}^3 \,{\tt{r}}(\Sigma) }{n^{1/2}} \\[0.3cm]
\hat{\textup{II}} & \ \leq \ \frac{c\,\beta_{3 q}^3}{n^{1/2}}. 
\end{align}
each hold with probability $1 - \frac{c}{n}$. Combining these bounds gives the stated result. \qed

\subsection{Lemmas for bootstrap approximation}

\begin{lemma} \label{lemma:boot1}
Suppose that the conditions of Theorem \ref{THM:BOOT} hold. Then, there is a constant $c>0$ not depending on $n$ such that the event
\begin{equation}\label{eqn:condcoupling1}
\P\bigg(\Big\| \sqrt n\Big(\boldsymbol\lambda_k (\tilde{\Sigma}^\star) - \boldsymbol\lambda_k (\tilde{\Sigma}^\star[1,1]) \Big)\Big\|_{\infty}
\, \geq  \, \frac{c \, \log (n) \beta_{2q} \, \lambda_1 (\Sigma) \, {\tt{r}} (\Sigma) }{n^{1/2}} \, \bigg| X\bigg) \ \leq \ \frac{c}{n^4}
\end{equation}
holds with probability at least $1 - \frac{c}{n}$.
\end{lemma}

\noindent
\Proof Let $\pi(X)$ denote the conditional probability on the left side of~\eqref{eqn:condcoupling1}. By Markov's inequality, we have
$$\P(\pi(X)\geq \ts\frac{1}{n^4}) \ \leq  \  n^4 \E[\pi(X)],$$
and so it is sufficient to show that there is a constant $c>0$ not depending on $n$ such that $\E[\pi(X)]\leq c/n^5$. 
In other words, it is enough to show that the event
$$\Big\| \sqrt n\Big(\boldsymbol\lambda_k (\tilde{\Sigma}^\star) - \boldsymbol\lambda_k (\tilde{\Sigma}^\star[1,1]) \Big)\Big\|_{\infty}
\, \leq \, \frac{c \, \log (n) \beta_{2q} \, \lambda_1 (\Sigma) \, {\tt{r}} (\Sigma) }{n^{1/2}} $$
holds with probability at least $1-c/n^5$. \\

Whenever the event $\big\{\lambda_k (\tilde{\Sigma}^\star[1,1]) > \lambda_1 (\tilde{\Sigma}^\star[2,2]) \big\}$ holds, Wielandt's inequality (Lemma \ref{background0})
implies
\begin{equation} \label{boot1:bound1}
\max_{1 \leq j \leq k} | \lambda_j (\tilde{\Sigma}^\star) - \lambda_j (\tilde{\Sigma}^\star[1,1]) \big|
\ \leq \
\frac{ \| \tilde{\Sigma}^\star[1,2] \|^2_\op }{\lambda_k (\tilde{\Sigma}^\star[1,1]) - \lambda_1 (\tilde{\Sigma}^\star[2,2])} . 
\end{equation}
The denominator in~\eqref{boot1:bound1} can be controlled by analogy with the proof of Lemma \ref{lemma:gauss1}: It follows from Weyl's inequality and Assumption~\ref{model_assumptions}.\ref{A2} that
\begin{equation} \label{boot_eig_lower_bound}
\begin{split}
\lambda_k (\tilde{\Sigma}^\star[1,1]) - \lambda_1 (\tilde{\Sigma}^\star[2,2])
& \ \geq \ \big( \lambda_k (\Sigma) - \lambda_{k+1} (\Sigma) \big) - 2 \| \tilde\Sigma^\star - \Lambda \|_\op\\[0.2cm]
& \ \geq c_1\lambda_1(\Sigma) - 2 \| \tilde\Sigma^\star - \Lambda \|_\op,
\end{split}
\end{equation}
for some constant $c_1>0$ not depending on $n$.
So, controlling the denominator in~\eqref{boot1:bound1} amounts to showing that the random variable $\| \tilde{\Sigma}^\star - \Lambda \|_\op$ is small with high probability.  We will proceed by upper bounding  $\|\tilde\Sigma^{\star}-\tilde\Sigma\|_\op+\|\tilde\Sigma-\Lambda\|_\op$ with high probability.
Applying Lemma \ref{tl4} with the choice $q=5\log(kn)$, it follows from a simple marginalization argument that the event 
\begin{equation} \label{boot1:bound2}
\| \tilde{\Sigma}^\star - \tilde{\Sigma} \|_\op \ \leq \ c\lambda_1 (\Sigma) \sqrt{ \frac{\log(n)\, \beta_{q} \,  {\tt{r}} (\Sigma) }{n} }
\end{equation}
holds with probability at least $1 - \frac{c}{n^5}$. Next, recalling the bound \eqref{gauss1:bound2} in Lemma \ref{lemma:gauss1}, there is a constant $c>0$ not depending on $n$, such that the event
\begin{equation} \label{boot1:prev}
\| \tilde{\Sigma} - \Lambda \|_\op \ \leq \ c\lambda_1 (\Sigma) \sqrt{ \frac{\log(n) \, \beta_{q} \, {\tt{r}} (\Sigma) }{n} }
\end{equation}
holds with probability at least $1 - \frac{c}{n^5}$. Combining the last two bounds with~\eqref{boot_eig_lower_bound} and the condition~\eqref{n_large} implies that the event $\{\lambda_k(\tilde\Sigma^{\star}[1,1])-\lambda_1(\tilde\Sigma^*[2,2])>\frac{c_1}{2}\lambda_1(\Sigma)\}$ holds with probability at least $1-c/n^5$, where $q=5\log(kn)$.

Lastly, to address the numerator in~\eqref{boot1:bound1}, it follows from Lemma \ref{tl3} and a simple marginalization argument that the event
\begin{equation}\label{eqn:bootoffdiag} 
 \|\tilde{\Sigma}^\star[1,2] \|_\op^{2}  \ \leq \ \frac{c \, \log(n) \, \beta_{2q} \, \lambda_1 (\Sigma) \, \tr (\Sigma) }{n} 
\end{equation}
holds with probability at least $1 - \frac{c}{n^5}$. This completes the proof.
\qed

~\\

\begin{lemma} \label{lemma:boot2} 
Suppose that the conditions of Theorem \ref{THM:BOOT} hold. Then, there is a constant $c>0$ not depending on $n$ such that the event 
\begin{equation}
\P\bigg(\Big\| \sqrt{n} \Big( \blambda_k (\tilde\Sigma^\star[1,1]) - \bdiag_k (\tilde{\Sigma}^\star[1,1]) \Big) \Big\|_\infty
\, \geq \, \frac{c \, \log(n)\, \beta_{2q} \, \lambda_1 (\Sigma) \, {\tt{r}} (\Sigma) }{n^{1/2}} \, \bigg| X\bigg) \ \leq  \ \frac{c}{n^4}
\end{equation}
holds with probability at least $1-\frac{c}{n}$.
\end{lemma}
\noindent
\Proof 
As in the proof of Lemma~\ref{lemma:gauss2}, we will lighten the use of subscripts by writing $\blambda(A)=(\lambda_1(A),\dots,\lambda_r(A))$ and $\bdiag(A)=(A_{11},\dots,A_{rr})$ for any symmetric matrix $A\in\R^{r\times r}$ and integer $r\geq 1$.  By the same reasoning used at the beginning of the proof of Lemma~\ref{lemma:boot1}, it is sufficient to show that
the event
\begin{align*}
\Big\| \sqrt{n} \Big( \blambda (\tilde\Sigma^\star[1,1]) - \bdiag(\tilde{\Sigma}^\star[1,1]) \Big) \Big\|_\infty
\ \leq \ \frac{c \, \log(n)\, \beta_{2q} \, \lambda_1 (\Sigma) \, {\tt{r}} (\Sigma) }{n^{1/2}} 
\end{align*}
holds with probability at least $1 -  \frac{c}{n^5}$. Overall, the current proof is similar to that of Lemma \ref{lemma:gauss2}. Let $D_0^\star = \tilde{\Sigma}^\star[1,1]$, and for each $r = 1, \dots, k-1$, partition the matrix $\tilde\Sigma^{\star} [1,1]$ recursively as 
\begin{align*}
D_{r-1}^\star
\ = \ \begin{pmatrix}
d_1 (D_{r-1}^\star) & V_r^\star \\
(V_r^\star)^\top & D_r^\star
\end{pmatrix},
\end{align*}
where $d_1 (D_{r-1}^\star)$ is a scalar, and $D_r^\star$ is of size $(k-r) \times (k-r)$. Based on the proof of Lemma~\ref{lemma:gauss2}, it suffices to show that there is a constant $c > 0$ not depending on $n$, such that for any $r=1,\dots,k$, the event
\begin{align*}
\bigg\|\blambda (D_{r-1}^\star) - \begin{bmatrix} d_1(D_{r-1}^\star)\\ \blambda(D_r^\star)  \end{bmatrix}\bigg\|_\infty \ \leq \ 
\frac{c \, \log(n) \, \beta_{2q} \, \lambda_1 (\Sigma) \, {\tt{r}} (\Sigma) }{n}
\end{align*}
holds with probability at least $1 - \frac{c}{kn^5}$. \\

Using the reasoning that led to~\eqref{eqn:V1bound} in the proof of Lemma~\ref{lemma:gauss2}, Wielandt's inequality (Lemma \ref{background0}) implies that if the event
$\{d_1(D_{r-1}^\star) >\lambda_1(D_r^{\star})\}$ holds, then the following event also holds
\begin{equation}\label{eqn:weilandtagain}
 \bigg\| \blambda (D_{r-1}^\star) - \begin{bmatrix} d_1(D_{r-1}^\star)\\ \blambda(D_r^\star)  \end{bmatrix} \bigg\|_\infty \ \leq \ \frac{\| V_r^\star \|_\op^2}{d_1 (D_{r-1}^\star) - \lambda_1 (D_r^\star)}.
 \end{equation}

The numerator in this bound~\eqref{eqn:weilandtagain} can be controlled with Lemma \ref{tl3} and a simple marginalization argument, which imply that under the choice $q=5\log(kn)$, there is a constant $c>0$ not depending on $n$ such that the event 
\[ \| V_r^\star \|_\op^2 \ \leq \ \frac{c \, \log (n) \, \beta_{2q} \, \lambda_1(\Sigma) \, \tr (\Sigma) }{n}  \]
holds with probability at least $1 - \frac{1}{k n^5}$. \\

To control the denominator in the bound~\eqref{eqn:weilandtagain}, Weyl's inequality and the reasoning in the proof of Lemma~\ref{lemma:gauss2} based on Assumption~\ref{model_assumptions}.\ref{A2} imply
\begin{align*}
d_1 (D_{r-1}^\star) - \lambda_1 (D_r^\star) 
& \ \geq \ c_2 \lambda_1 (\Sigma) - 2 \| \tilde\Sigma^\star - \Lambda \|_\op,
\end{align*}
for some constant $c_2>0$ not depending on $n$.
Next, by combining the bound~\eqref{gauss1:bound2} and Lemma \ref{tl4}, it follows that the event
\begin{align} \label{lemma:boot2:bound}
\| \tilde\Sigma^\star - \Lambda \|_\op \ \leq \ c\lambda_1 (\Sigma) \sqrt{ \frac{\log(n) \, \beta_q \, {\tt{r}} (\Sigma)}{n} } 
\end{align}
holds with probability at least $1 - \frac{c}{kn^5}$. Therefore, condition~\eqref{n_large} implies that the lower bound
$$d_1 (D_{r-1}^\star) - \lambda_1 (D_r^\star)
 \ \geq \ \ts\frac{c_2}{2} \lambda_1 (\Sigma)
 $$
holds with probability at least $1 - \frac{1}{kn^5}$. This completes the proof.

\qed

~\\[-1cm]

\begin{lemma} \label{lemma:boot3}
Suppose that the conditions of Theorem \ref{THM:BOOT} hold, and let $\hat{\textup{II}}$ be as defined in~\eqref{eqn:hatIIdef}. Then, there is a constant $c > 0$ not depending on  $n$ such that the event
\[ 
\hat{\textup{II}} \ \leq \   \frac{c  \, \beta_{3 q}^3}{n^{1/2}} \]
holds with probability at least $1 - \frac{c}{n}$.
\end{lemma}

\noindent
\Proof Let $W_1,\dots,W_n$ and $\bar W$ be as defined at the beginning of Appedix~\ref{appendix_C}. Also, define the vector $W_i^{\star}\in\R^k$ with $j$th entry $W_{ij}^{\star}=\langle u_j, Z_i^{\star}\rangle^2-1$, and define $Y_i^\star = \Lambda_k (W_i^\star-\bar W)$. These definitions give the relation
\[ \sqrt{n} \big( \bdiag_k (\tilde{\Sigma}^{\star}[1,1]) - \bdiag_k (\tilde\Sigma[1,1]) \big) \ = \ \frac{1}{ \sqrt{n}} \sum_{i = 1}^n Y_i^\star .\]
Also note that $\E[Y_i^{\star}|X]=0$, and
$$\E\big[Y_i^{\star}(Y_i^{\star})\ttop\big|X\big] \ = \ \Lambda_k\hat\Gamma\Lambda_k.$$
Due to Bentkus' multivariate Berry-Esseen theorem (Lemma \ref{background5}), we have
\begin{align} \label{eqn:boot3:1} 
\sup_{t\in\R^k} \bigg| \P \bigg( \frac{1}{\sqrt{ n }} \sum_{i = 1}^{ n } Y_i^\star \preceq t \; \bigg| \; X \bigg) - \P \big( \xi  \preceq t \; \big| \; X \big) \bigg| \leq \frac{c\, \, \E \big[ \| (\Lambda_k \hat\Gamma \Lambda_k)^{-1/2} Y_1^\star \|_2^3 \, \big| \, X \big]  }{{ n }^{1/2} } .
\end{align}
Applying Lemma \ref{tl6} with $q = 5 \log (kn)$ and using Chebyshev's inequality, there is a constant $c>0$ not depending on $n$ such that the event
\[  \E \big[ \| (\Lambda_k \hat\Gamma \Lambda_k)^{-1/2} Y_1^\star \|_2^3 \, \big| \, X \big] \ \leq \   c  \, \beta_{3 q}^3 \]
holds with probability at least $1 - \frac{c}{n}$. \qed

\begin{lemma} \label{lemma:boot5}
Suppose that the conditions of Theorem \ref{THM:BOOT} hold, and let $\hat{\textup{I}}$ be as defined in~\eqref{eqn:hatIdef}. Then, there is a constant $c > 0$ not depending on $n$ such that the event
\begin{align}\label{eqn:boot6}
\hat{\textup{I}} \ \leq \ \frac{c \, \log(n) \, \beta_{3q}^3 \, {\tt{r}} (\Sigma)}{n^{1/2}} 
\end{align}
\normalsize
holds with probability at least $1 - \frac{c}{n}$.
\end{lemma}

\noindent
\Proof In a similar manner to the proof of Lemma~\ref{lemma:gauss4}, the left side of~\eqref{eqn:boot6} can be bounded by the sum $ 2\hat{\textup{II}}+\hat{\textup{G}}_1(\epsilon) +2\hat{\textup{G}}_2(\e)+2\hat{\textup{G}}_3$, with the last three terms defined for any $\e>0$ according to
\begin{align*}
\hat{\textup{G}}_1(\epsilon) & \ =\  \P \bigg( \Big\|\sqrt{n} \Big( \blambda_k (\tilde\Sigma^\star)  - \blambda_k(\hat\Sigma) \Big) \ - \ \sqrt{n} \Big( \bdiag_k (\tilde{\Sigma}^\star[1,1]) - \bdiag_k(\tilde{\Sigma}[1,1]) \Big) \Big\|_\infty \geq \epsilon \, \bigg| \, X \bigg), \\[0.3cm]
\textup{G}_2(\epsilon) & \ =\ \sup_{t \in \R^k} \Big| \P \Big( \zeta  \mylessthan + \epsilon 1_k \Big) - \P \Big( \zeta  \mylessthan - \epsilon 1_k \Big) \Big|, \\[0.3cm]
\hat{\textup{G}}_3 & \ =\ \sup_{t \in \R^k} \Big| \P \big(  \zeta  \mylessthan  \big)- \P \big( \xi  \mylessthan \, \big| \, X \big) \Big|.
\end{align*}
First, recall from Lemma~\ref{lemma:boot3} that there is a constant $c>0$ not depending on $n$ such that $\hat{\textup{II}}\leq c\beta_{3 q}^3/n^{1/2}$ holds with probability at least $1-c/n$. \\

Next, when handling the terms $\hat{\textup{G}}_1(\e)$, $\textup{G}_2(\e)$, and $\hat{\textup{G}}_3$ below, we will take $\e$ to be of the form $\epsilon = \frac{c}{\sqrt n}\log(n) \beta_{2q} \lambda_1 (\Sigma) {\tt{r}}(\Sigma)$.
Using the choice $q = 5\log(kn)$, Lemmas~\ref{lemma:gauss1},~\ref{lemma:gauss2},~\ref{lemma:boot1} and~\ref{lemma:boot2} imply that there is a constant $c>0$ not depending on $n$ such that the event
\[ \hat{\textup{G}}_1(\e)
\ \leq \ \frac{c}{n} \]
holds with probability at least $1 - \frac{c}{n}$.
With regard to $\textup{G}_2(\epsilon)$, observe that it is deterministic and equal to $\textup{J}(\epsilon)$ in the proof of Lemma \ref{lemma:gauss4}. Therefore, the bound~\eqref{eqn:NazarovJ2} gives
\[ \textup{G}_2(\epsilon)
\ \lesssim
\frac{ \log(n) \, \beta_{2q}  \, {\tt{r}}(\Sigma)}{n^{1/2}}.\]
Lastly, the bound~\eqref{eqn:gaussiancompdesired} (to be established in Appendix~\ref{appendix_A}) implies that the event
\begin{align*}
\hat{\textup{G}}_3 & \leq \ \frac{c \, \log(n) \, \beta_{2q}^2 }{n^{1/2}}
\end{align*}
holds with probability at least $1 - \frac{c}{n}$. Combining the last several bounds yields the stated result.
\qed

\section{Proof of Theorem \ref{THM:ALL} } \label{appendix_A}

By comparing the bounds~\eqref{eqn:gaussapproxdisp} and~\eqref{eqn:bootapproxdisp} in Theorems \ref{THM:GAUSS} and \ref{THM:BOOT}, it is enough to show that there is a constant $c>0$ not depending on $n$ such that the event
\begin{equation}\label{eqn:gaussiancompdesired}
\sup_{t \in\R^k} \Big| \P \big( \xi \preceq t \, \big| \, X \big) - \P \big( \zeta  \preceq t \big) \Big| 
\ \leq \ \frac{c \, \log (n)\, \beta_{2q}^2}{n^{1/2}}
\end{equation}
holds with probability at least $1 - \frac{c}{n}$.
To this end, define the three matrices
\begin{align*}
C_k & \ = \  \Lambda_k \Gamma \Lambda_k,\\[0.2cm]
\widehat{C}_k & \ = \  \Lambda_k \widehat\Gamma\Lambda_k,\\[0.2cm]
\hat B_k & \ = \  C_k^{-1/2} \hat{C}_k C_k^{-1/2} - I_k.
\end{align*}
By Lemma~\ref{lem:Zhilova}, there is a constant $c>0$ not depending on $n$ such that the bound
\begin{align*}
\sup_{t \in \R^k} \Big| \P \big( \xi \preceq t \, \big| \, X \big) - \P \big( \zeta \preceq t \big) \Big| 
& \ \leq \ c \, \| \hat B_k \|_\op
\end{align*}
holds almost surely.
Furthermore, due to Assumptions \ref{model_assumptions}.\ref{A2} and \ref{model_assumptions}.\ref{A3}, the following bounds also hold almost surely,
\begin{align*}
\| \hat B_k \|_\op 
& \ \leq \ \big\| C_k^{-1/2} \big\|_\op^2 \| \hat{C}_k - C_k \|_\op \\[0.2cm]
& \ \leq \ c \big(\ts\frac{\lambda_1(\Sigma)}{\lambda_k(\Sigma)}\big)^2 \big\| \widehat\Gamma - \Gamma \big\|_\op\\[0.2cm]
& \ \leq \ c\, \big\| \widehat\Gamma - \Gamma \big\|_\op.
\end{align*}
Next, Lemma~\ref{lem:hatgamma} implies that the event
\[  \| \hat\Gamma - \Gamma \|_\op \ \leq \ \frac{c \, q \, \beta_{2q}^2 }{n^{1/2} } \]
holds with probability at least $1 - e^{-q}$, and in light of the stated choice of $q=5\log(kn)$, the proof is complete. \qed

\section{Proof of Theorem~\ref{THM:TRANS_GENERAL}} \label{appendix_D}
A simple rescaling argument can be used to show that the left side of the bound in Theorem~\ref{THM:TRANS_GENERAL} is the same as
\begin{align}
\footnotesize
 \sup_{t \in \R^k } \bigg| \P \bigg( \boldsymbol h (\blambda_k(\hat\Sigma)/\lambda_1(\Sigma)) -\boldsymbol h( \blambda_k(\Sigma)/\lambda_1(\Sigma)) \mylessthan \bigg) 
- \P \bigg( \frac{\boldsymbol h (\blambda_k(\hat\Sigma^{\star})/\lambda_1(\Sigma)) -\boldsymbol h( \blambda_k(\hat\Sigma)/\lambda_1(\Sigma))}{(\boldsymbol{\hat{\varsigma}}_k/\boldsymbol{\varsigma}_k)^{\tau}} \mylessthan \, \bigg| \, X \bigg) \bigg|,
\end{align}
where we note that the quantity $(\boldsymbol{\hat{\varsigma}}_k/\boldsymbol{\varsigma}_k)^{\tau}$ is invariant to rescaling of the observations. Hence, without loss of generality, we may assume that $\lambda_1 (\Sigma) = 1$ in the remainder of this appendix. Consequently, Assumption~\ref{model_assumptions}.\ref{A2} implies there is a positive constant $c_\circ$ not depending on $n$ such that $\lambda_j (\Sigma) \in [c_\circ, 1]$ holds for all $j=1,\dots,k$.  (These points will sometimes be used without being explicitly mentioned in this appendix.) The proof is completed by combining Lemmas~\ref{lemma:trans:gauss} and~\ref{lemma:transcompare} given below.\qed

\subsection{Lemmas for bootstrap with transformations}

\begin{lemma} \label{lemma:trans:gauss}
Suppose that the conditions of Theorem \ref{THM:TRANS_GENERAL} hold. 
\begin{itemize}
\item[(i)] Let $\zeta\in\R^k$ be a random vector that is distributed as $N (0, \Lambda_k \Gamma \Lambda_k)$. Then,
\begin{equation*}
\small
\sup_{t \in \R^k} \Big| \P \Big( \sqrt{n} \big( \boldsymbol h (\blambda_k(\hat\Sigma)) - \boldsymbol h (\blambda_k(\Sigma)) \big) \mylessthan \Big) - \P \Big( \boldsymbol h' (\blambda_k(\Sigma)) \odot \zeta \mylessthan \Big) \Big| 
\ \lesssim \
\frac{\log (n) \, \beta_{2q}^3 \, {\tt{r}} (\Sigma)}{n^{1/2} } .
\end{equation*}
%
%
%
%
%
\item[(ii)] Let $\xi\in\R^k$ be a random vector that is conditionally distributed as $N \big(0, \Lambda_k \hat\Gamma \Lambda_k \big)$ given the observations $X_1,\dots,X_n$. Then, there is a constant $c > 0$ not depending on $n$ such that the event
\small
\[ 
\sup_{t \in \R^k} \bigg| \P \bigg( \frac{\sqrt{n} \big( \boldsymbol h (\blambda_k(\hat\Sigma^\star)) - \boldsymbol h (\blambda_k(\hat\Sigma)) \big)}{(\boldsymbol{\hat{\varsigma}}_k / \boldsymbol{\varsigma}_k )^\tau} \mylessthan \, \Big| \, X \bigg) - \P \bigg( \frac{\boldsymbol h' (\blambda_k(\hat\Sigma)) \, \odot \, \xi}{(\boldsymbol{\hat{\varsigma}}_k / \boldsymbol{\varsigma}_k )^\tau} \mylessthan \, \bigg| \, X \bigg) \bigg| 
\ \leq \ 
\frac{c\,\log (n) \, \beta_{3q}^3 \, {\tt{r}} (\Sigma)}{n^{1/2} } \]
\normalsize
holds with probability at least $1 - \frac{c}{n}$. 
\end{itemize}
\end{lemma}

\noindent 
\Proof Part (\emph{i}). Applying a Taylor expansion for each $j = 1, \dots, k$, we obtain\\[-0.2cm]
\begin{align*}
 \sqrt{n} \big( h \big( \lambda_j (\hat\Sigma) \big) - h \big( \lambda_j (\Sigma) \big) \big)  & = \sqrt{n} \, h' \big( \lambda_j (\Sigma) \big) \, \big( \lambda_j (\hat\Sigma) - \lambda_j (\Sigma) \big) \\[0.2cm]
& \ \ \ \ \ + \sqrt{n} \big( \lambda_j (\hat\Sigma) - \lambda_j(\Sigma) \big)^2 \int_0^1 (1 - t) \, h'' \big( \lambda_j(\Sigma) + t (\lambda_j(\hat\Sigma) - \lambda_j(\Sigma)) \big) \textup{d}t.
\end{align*}
Let $S$, $T$, $V$, and $R$ be random vectors in $\R^k$ whose entries are defined by
\begin{align*}
S_j & \ = \ \sqrt{n} \big( h \big( \lambda_j (\hat\Sigma) \big) - h \big( \lambda_j (\Sigma) \big) \big), \\[0.2cm]
T_j & \ = \ h ' \big( \lambda_j (\Sigma) \big) \zeta_j, \\[0.2cm]
V_j & \ = \ \sqrt{n} \, h ' \big( \lambda_j (\Sigma) \big) \, \big( \lambda_j (\hat\Sigma) - \lambda_j(\Sigma) \big) , \\[0.1cm]
R_j & \ = \sqrt{n} \big( \lambda_j (\hat\Sigma) - \lambda_j(\Sigma) \big)^2 \int_0^1 (1 - t) \, h '' \big( \lambda_j(\Sigma) + t (\lambda_j(\hat\Sigma) - \lambda_j(\Sigma)) \big) \textup{d}t.
\end{align*}
Lemma \ref{background6} implies that for any $r > 0$,
\begin{align} \label{thm:transform:gauss:bound0}
\small
\sup_{t \in \R^k} \Big| \P \big( S \mylessthan \big) - \P \big( T \mylessthan \big) \Big| 
\ \lesssim \ 
 \sup_{t \in \R^k} \Big| \P \big( V \mylessthan \big) - \P \big( T \mylessthan \big) \Big| \  +  \frac{r }{\min_{1 \leq j \leq k}  \sqrt{ \var ( T_j )} } \ + \ \P ( \| R \|_\infty \geq r ).
\end{align}
%
%
%
%
Theorem \ref{THM:GAUSS} gives a bound on the first term,
\begin{align} \label{thm:transform:gauss:bound1}
\sup_{t \in \R^k}  \Big| \P \big( V \mylessthan \big) - \P \big( T \mylessthan \big) \Big| 
\ \lesssim \ \frac{\log(n) \, \beta_{2q}^3 \, {\tt{r}} (\Sigma) }{n^{1/2}}.
\end{align}
Next, recall from the discussion at the beginning of this appendix that we may assume there is a constant $c_{\circ}>0$ not depending on $n$ such that $\lambda_j(\Sigma)\in[c_{\circ},1]$ holds for all $j=1,\dots,k$. In turn, combining this with Assumption~\ref{model_assumptions}.\ref{A3} implies 
\begin{align} \label{thm:transform:gauss:bound2}
\var (T_j) 
& \ = \  \lambda_j(\Sigma)^2 \cdot \Gamma_{jj} \cdot h ' \big( \lambda_j (\Sigma) \big)^2 \notag \\[0.2cm]
& \ \gtrsim \ 1 . 
\end{align}

It remains to find a bound for the last term in \eqref{thm:transform:gauss:bound0} and to choose a suitable value for $r$. For each $j=1,\dots,k$, define the event 
$$\A_j = \Big\{ \big| \lambda_j(\hat\Sigma) - \lambda_j(\Sigma) \big| \leq \frac{c_{\circ}}{2} \Big\}.$$ 
Then,
\begin{align*}
\P \big( |R_j|  \geq r \big)
& \ \leq \ \P  \big( \{|R_j| \geq r\} \cap \A_j \big) + \P \big(  \A_j^c \big).
\end{align*}
Using the choice $q=5\log(kn)$ and Weyl's inequality, the bound \eqref{gauss1:bound_interm} implies that
\begin{align} \label{sample_eig_bound}
\max_{1 \leq j \leq k} \big\| \lambda_j (\hat\Sigma) - \lambda_j(\Sigma) \big\|_q
& \ \leq \ \big\| \| \tilde\Sigma - \Lambda \|_\op \big\|_q  \notag \\[0.2cm]\notag
& \ \lesssim \ \frac{\lambda_1(\Sigma) \sqrt{ \log(n) \beta_q{\tt{r}} (\Sigma) }  }{n^{1/2}}\notag\\[0.2cm]
& \ \leq \ \frac{c_{\circ}}{2e},
\end{align}
where the condition~\eqref{n_large} has been used in the last step, and the rationale for the constant $\frac{c_{\circ}}{2e}$ will be seen in the next step. By Chebyshev's inequality, we have
\begin{equation} \label{thm:trans:gauss:tail_prob}
\max_{1\leq j\leq k} \P (\A_j^c) \ \leq \ \max_{1\leq j\leq k} \, \frac{\big\| \lambda_j (\hat\Sigma) - \lambda_j(\Sigma) \big\|_q^q}{(c_{\circ}/2)^q} \ \leq \ \frac{(c_{\circ}/(2e))^q}{(c_{\circ}/2)^q} \ =  \  e^{-q} \ \leq \ \frac{1}{kn}.
\end{equation}
In order to handle $\P  \big( \{|R_j| \geq r\} \cap \A_j \big)$, first let $C$ denote the following supremum
\begin{equation}\label{eqn:Cdef}
C=\sup\Big\{ |h''(x)| \, \Big| \, c_{\circ}/2 \leq x\leq   1+ c_{\circ}/2\Big\},
\end{equation}
which is finite and does not depend on $n$.
Taking $r = \frac{c}{\sqrt n}\lambda_1(\Sigma)^2 \log (kn) \beta_q {\tt{r}} (\Sigma)$ for a sufficiently large constant $c$, we have
\begin{align} \label{intersection_bound}
\P  \big( \{ |R_j| \geq r \} \cap \A_j \big) 
& \ \leq \ \P \Big( \sqrt{n} \big( \lambda_j (\hat\Sigma ) - \lambda_j (\Sigma ) \big)^2\cdot \ts\frac{C}{2}  \geq  r \Big) \notag \\[0.2cm]
& \ \leq \ \P \bigg( \big| \lambda_j (\hat\Sigma ) - \lambda_j (\Sigma ) \big| \geq \frac{\lambda_1(\Sigma)}{n^{1/2}} \sqrt{(2c/C)\log(kn) \beta_q {\tt{r}} (\Sigma)} \bigg) \notag \\[0.2cm]
& \ \leq \ \ts \frac{1}{kn}, 
\end{align}
where the previous step uses~\eqref{sample_eig_bound}.
Combining the last several steps shows that $\P(|R_j|\geq r) \ \leq \ \frac{2}{kn}$ holds for all $j=1,\dots,k$
and so a union bound gives
\begin{align} \label{thm:transform:gauss:bound3}
\P (\| R \|_\infty \geq r ) \ \lesssim \ \ts \frac{1}{n}.
\end{align}
Substituting the bounds \eqref{thm:transform:gauss:bound1}, \eqref{thm:transform:gauss:bound2}, and \eqref{thm:transform:gauss:bound3} into \eqref{thm:transform:gauss:bound0} proves the statement (\emph{i}). \\

\noindent Part (\emph{ii}). First note that by rescaling, it is enough to show that the event
\[ 
\small
\sup_{t \in \R^k} \bigg| \P \bigg( \sqrt{n} \big( \boldsymbol h (\blambda_k(\hat\Sigma^\star)) - \boldsymbol h (\blambda_k(\hat\Sigma)) \big) \mylessthan \, \Big| \, X \bigg) - \P \bigg( \boldsymbol h' (\blambda_k(\hat\Sigma)) \, \odot \, \xi \mylessthan \, \bigg| \, X \bigg) \bigg| 
\ \leq \ 
\frac{c\,\log (n) \, \beta_{3q}^3 \, {\tt{r}} (\Sigma)}{n^{1/2} } \]
holds with probability at least $1-c/n$.
Similar to \eqref{thm:transform:gauss:bound0}, the left side can be decomposed as
\small{
\begin{align} \label{thm:transform:boot:bound0}
\sup_{t \in \R^k} \Big| \P \big( S^\star \mylessthan \, | \, X \big) & - \P \big( T^\star \mylessthan \, |\, X \big) \Big| \notag \\[0.2cm]
& \ \lesssim \ \sup_{t \in \R^k} \Big| \P \big( V^\star \mylessthan \,|\, X \big) - \P \big( T^\star \mylessthan \,|\, X \big) \Big| +  \frac{r }{\min_{1 \leq j \leq k} \sqrt{\smash[b]{ \var ( T_j^\star \,|\, X )} } }\ + \ \P ( \| R^\star \|_\infty \geq r \,|\,X), 
\end{align}
}
%
\large
where $S^\star$, $T^\star$, $R^\star$, $V^\star$ are vectors in $\R^k$ defined as
\begin{align*}
S_j^\star & \ = \ \sqrt{n} \big( h \big( \lambda_j (\hat\Sigma^\star) \big) - h \big( \lambda_j (\hat\Sigma) \big) \big), \\[0.2cm]
T_j^\star & \ = \ h ' \big( \lambda_j (\hat\Sigma) \big) \, \xi_j, \\[0.2cm]
V_j^\star & \ = \ \sqrt{n} \, h ' \big( \lambda_j (\hat\Sigma) \big) \, \big( \lambda_j (\hat\Sigma^\star) - \lambda_j(\hat\Sigma) \big) , \\[0.1cm]
R_j^\star & \ = \sqrt{n} \big( \lambda_j (\hat\Sigma^\star) - \lambda_j(\hat\Sigma) \big)^2 \int_0^1 (1 - t) \, h '' \big( \lambda_j(\hat\Sigma) + t (\lambda_j(\hat\Sigma^\star) - \lambda_j(\hat\Sigma)) \big) \textup{d}t.
\end{align*}
The first term of the bound \eqref{thm:transform:boot:bound0} is handled by Theorem~\ref{THM:BOOT}. For the middle term, first note that $\var(T_j^{\star}\,|\,X) = \lambda_j(\Sigma)^2 \cdot \hat\Gamma_{jj} \cdot h'(\lambda_j(\hat\Sigma))^2$. 
Using the condition~\eqref{n_large}, it follows from \eqref{gauss1:bound2} and Lemma~\ref{lem:hatgamma} that for each $j=1,\dots,k$, the events 
\[\lambda_j(\hat\Sigma) \geq c \, \lambda_j(\Sigma),\\[0.cm] \]
\[ \hat\Gamma_{jj} \geq c \, \Gamma_{jj},\\[0.cm] \]
\[ h'(\lambda_j(\hat\Sigma)) \geq c \, h'(\lambda_j(\Sigma)) \]
each hold with probability at least $1 - \frac{c}{kn}$. Hence, Assumptions~\ref{model_assumptions}.\ref{A2} and~\ref{A3} imply that the event
$$\min_{1\leq j\leq k} \var(T_j^{\star}\,|\,X) \geq c$$
holds with probability at least $1-\frac{c}{n}$. Next, define the event 
\[ \A_j^\star = \Big\{ \big| \lambda_j(\hat\Sigma^\star) - \lambda_j(\hat\Sigma) \big| \leq \frac{c_{\circ}}{2} \Big\} \]
with the constant $c_\circ$ having the same definition as in Part (\emph{i}). 
For the last term on the right side of \eqref{thm:transform:boot:bound0}, we claim that the event 
\[ \P(\|R^\star\|_\infty \geq r \,|\, X) \ \leq \ \frac{c}{n} \]
holds with probability at least $1 - \frac{c}{n}$ when $r$ is appropriately chosen. This can be established with a union bound
\begin{equation}\label{eqn:Astarbound}
 \P(\|R^\star\|_\infty \geq r \,|\, X) \ \leq \
\ \sum_{j=1}^k \P \big( \big\{ |R_j^\star| \geq r \big\} \cap \mathcal{A}_j^\star \, \big| \, X \big) + \P \big( (\mathcal{A}_j^\star)^c \, | \, X \big).
\end{equation}
Analogously to \eqref{thm:trans:gauss:tail_prob}, we can apply Lemma~\ref{tl4} with $q = 5 \log(kn)$ and the condition~\eqref{n_large} to conclude that for each $j=1,\dots,k$ the event
\[ \max_{1 \leq j \leq k} \P \big( (\mathcal{A}_j^\star)^c \, | \, X \big) \ \leq \ \frac{c}{kn}\]
holds with probability at least $1 - \frac{c}{kn}$. For the first term on the right side of~\eqref{eqn:Astarbound}, we may use an argument similar to the one leading up to \eqref{intersection_bound} with $r = \frac{c}{\sqrt{n}}\lambda_1(\Sigma)^2 \log(kn) \beta_q {\tt{r}}(\Sigma)$  to show that for each $j=1,\dots,k$ the event
\[ \P\big( \{|R_j^\star| \geq r\} \cap \mathcal{A}_j^\star \, \big| \, X\big) \ \leq \ \frac{c}{kn} \]
holds with probability at least $1 - \frac{c}{kn}$. Combining the last few bounds completes the proof.
\qed

\begin{lemma}\label{lemma:transcompare}
Suppose the conditions of Theorem~\ref{THM:TRANS_GENERAL} hold. Then, there is a constant $c>0$ not depending on $n$ such that the event
\begin{equation} \sup_{t \in \R^k} \bigg| \P \Big( \boldsymbol h'(\blambda_k(\Sigma)) \odot \zeta  \mylessthan \Big) - \P \bigg( \, \frac{\boldsymbol h' (\blambda_k(\hat\Sigma))\odot \xi }{(\boldsymbol{\hat{\varsigma}}_k / \boldsymbol{\varsigma}_k )^\tau}  \mylessthan \, \bigg| \, X \bigg) \bigg|  \ \leq \ \frac{c \, \log (n) \, \beta_{2q}^5 \, {\tt{r}} (\Sigma) }{n^{1/2}}    
\end{equation}
holds with probability at least $1-\frac{c}{n}$.
\end{lemma}
\noindent \Proof
Define the matrices
\begin{align*}
\mathsf{C}_k
& \ =\ \textup{diag} \big( \boldsymbol h'(\blambda_k(\Sigma)) \big) \Big( \Lambda_k \Gamma \Lambda_k \Big) \textup{diag} \big( \boldsymbol h'(\blambda_k(\Sigma)) \big) \\[0.2cm]
\hat{\mathsf{C}}_k 
& \ =\ \textup{diag} \bigg(  \frac{\boldsymbol h' (\blambda_k(\hat\Sigma))}{(\boldsymbol{\hat{\varsigma}}_k / \boldsymbol{\varsigma}_k )^\tau} \bigg) \Big( \Lambda_k \widehat\Gamma\Lambda_k \Big)\textup{diag} \bigg(  \frac{\boldsymbol h' (\blambda_k(\hat\Sigma))}{(\boldsymbol{\hat{\varsigma}}_k / \boldsymbol{\varsigma}_k )^\tau} \bigg) \\[0.2cm]
\hat{\mathsf{B}}_k
& \ =\ \mathsf{C}_k^{-1/2} \hat{\mathsf{C}}_k \mathsf{C}_k^{-1/2} - I_k.
\end{align*}
By Lemma~\ref{lem:Zhilova}, the following bound holds almost surely
\begin{equation}\label{eqn:Pinsker}
\begin{split}
\sup_{t \in \R^k} \bigg| \P \bigg( \, \frac{\boldsymbol h' (\blambda_k(\hat\Sigma))\odot \xi}{(\hat{\boldsymbol{\varsigma}}_k / \boldsymbol{\varsigma}_k )^\tau}   \mylessthan \, \bigg| \, X \bigg) - \P \Big( \boldsymbol h'(\blambda_k(\Sigma)) \odot \zeta  \mylessthan \Big) \bigg| 
& \ \leq \ c \, \| \hat{\mathsf{B}}_k \|_\op
\end{split}
\end{equation}
for some constant $c > 0$ not depending on $n$. 
 Using several applications of the triangle inequality, 
 the following bound holds almost surely,
\footnotesize
\begin{align}  \label{trans:all:bound1}
\| \hat{\mathsf{B}}_k \|_\op 
& \ \leq \ \big\| \mathsf{C}_k^{-1/2} \big\|_\op^2 \| \mathsf{C}_k - \hat{\mathsf{C}}_k\|_\op \notag \\[0.2cm]
& \ \leq \ \big\| \mathsf{C}_k^{-1/2} \big\|_\op^2 \cdot  \Big\|\boldsymbol h'(\blambda_k(\Sigma)) \Big\|_\infty\cdot \Big\|\Lambda_k\Gamma\Lambda_k\Big\|_{\op}\cdot \Big\|\boldsymbol h'(\blambda_k(\Sigma)) - \frac{\boldsymbol h' (\blambda_k(\hat\Sigma))}{(\boldsymbol{\hat{\varsigma}}_k / \boldsymbol{\varsigma}_k )^\tau}\Big\|_{\infty} \\[0.2cm]
& \ \ \ \ \ + \ \ \  \  \big\| \mathsf{C}_k^{-1/2} \big\|_\op^2 \cdot \Big\|\boldsymbol h'(\blambda_k(\Sigma)) \Big\|_\infty\cdot\Big\|\Lambda_k(\Gamma-\hat\Gamma)\Lambda_k\Big\|_{\op}\cdot \Big\|\frac{\boldsymbol h' (\blambda_k(\hat\Sigma))}{(\boldsymbol{\hat{\varsigma}}_k / \boldsymbol{\varsigma}_k )^\tau}\Big\|_{\infty}\notag\\[0.2cm]
& \ \ \ \ \ + \ \ \  \  \big\| \mathsf{C}_k^{-1/2}\big\|_{\op}^2 \cdot \Big\|\boldsymbol h'(\blambda_k(\Sigma)) - \frac{\boldsymbol h' (\blambda_k(\hat\Sigma))}{(\boldsymbol{\hat{\varsigma}}_k / \boldsymbol{\varsigma}_k )^\tau}\Big\|_{\infty} \cdot \Big\|\Lambda_k\hat\Gamma\Lambda_k\Big\|_{\op} \cdot \Big\| \frac{\boldsymbol h' (\blambda_k(\hat\Sigma))}{(\boldsymbol{\hat{\varsigma}}_k / \boldsymbol{\varsigma}_k )^\tau} \Big\|_\infty.\notag 
\end{align}
\large
The leading factor satisfies
$$\big\| \mathsf{C}_k^{-1/2} \big\|_\op^2 \ \lesssim \ 1$$
because we may assume that there is a constant $c_{\circ}>0$ not depending on $n$  such that $\lambda_j(\Sigma)\in[c_{\circ},1]$ holds for all $j=1,\dots,k$ (as discussed at the beginning of this appendix), and also because $\lambda_k(\Gamma)\gtrsim 1$ by Assumption~\ref{model_assumptions}.\ref{A3}.  Also, the quantity $\| \hat\Gamma - \Gamma \|_\op$ is handled by Lemma~\ref{lem:hatgamma}, which shows there is a constant $c>0$ not depending on  $n$ such that the event
\begin{equation} \label{trans:all:bound2}
\| \hat\Gamma - \Gamma \|_\op
\ \leq \ \frac{c \, \log (n) \, \beta_{2q}^2 }{n^{1/2} }
\end{equation}
holds with probability at least $1-\frac{c}{n}$. Combining this with the condition~\eqref{n_large} and the bound $\|\Gamma\|_\op\lesssim \beta_2^2$, it follows that
\begin{align} \label{trans:all:bound3}
\|\Lambda_k \hat\Gamma\Lambda_k \|_\op
& \ \leq \ \| \Lambda_k\Gamma\Lambda_k \|_\op + \|\Lambda_k( \hat\Gamma - \Gamma )\Lambda_k\|_\op \notag  \ \leq \ c \, \beta_2^2
\end{align} 
holds with probability at least $1 - \frac{c}{n}$.
To handle the remaining quantities in the bound~\eqref{trans:all:bound1}, note that the triangle inequality yields
\begin{equation} \label{cmp_bound}
\footnotesize
\bigg\|\frac{\boldsymbol h' (\blambda_k(\hat\Sigma))}{(\boldsymbol{\hat{\varsigma}}_k / \boldsymbol{\varsigma}_k )^\tau} - \boldsymbol h'(\blambda_k(\Sigma)) \bigg\|_\infty
\ \leq \ 
\big\| (\boldsymbol{\hat{\varsigma}}_k / \boldsymbol{\varsigma}_k)^{-\tau} \big\|_\infty \, \big\| \boldsymbol h'(\blambda_k(\hat\Sigma)) - \boldsymbol h'(\blambda_k(\Sigma))  \big\|_\infty
+ 
\big\|\boldsymbol h'(\blambda_k(\Sigma)) \big\|_\infty \, \big\| (\boldsymbol{\hat{\varsigma}}_k / \boldsymbol{\varsigma}_k)^{-\tau} - 1_k \big\|_\infty.
\end{equation}
%
%
\noindent \!\!\!\! Using Lemmas~\ref{lemma:eig_var},~\ref{lemma:trans_var}, and~\ref{lemma:trans_var_lower_bound}, as well as some elementary inequalities, the following bounds hold with probability at least $1-c/n$ for any $\tau\in[0,1]$,
\begin{align} \label{trans_step1}
\big\| (\hat{\boldsymbol \varsigma}_k / \boldsymbol \varsigma_k)^{-\tau} - 1_k \big\|_\infty
& \ \leq \ \Big\|1_k - \boldsymbol\varsigma_k / \hat{\boldsymbol \varsigma}_k \Big\|_{\infty}\notag 
\\[0.2cm]
& \ \leq \ \Big\|\ts\frac{\hat{\boldsymbol \varsigma}_k^2-\boldsymbol\varsigma_k^2}{ \hat{\boldsymbol\varsigma}_k \varsigma_k }\Big\|_{\infty}\\[0.2cm]
& \ \leq \ \frac{c \, \log(n) \, \beta_{2q}^3 \, {\tt{r}}(\Sigma)}{n^{1/2}}.\notag
\end{align} 
Furthermore, using the condition~\eqref{n_large}, this implies that the event
\begin{equation} \label{trans_step2}
\big\| (\boldsymbol{\hat{\varsigma}}_k / \boldsymbol{\varsigma}_k)^{-\tau} \big\|_\infty \leq c
\end{equation}
also holds with probability at least $1-c/n$.
Next, we derive upper bounds for $\| \boldsymbol h' (\lambda_j(\hat\Sigma)) \big\|_{\infty}$ and $\|\boldsymbol h' (\lambda_j(\hat\Sigma)) - \boldsymbol h' (\lambda_j(\Sigma)) \big\|_{\infty}$. Using \eqref{gauss1:bound2} in the proof of Lemma~\ref{lemma:gauss1}, it follows that the bounds
\begin{align} \label{trans_step4}
\big\| \boldsymbol h' (\blambda_k(\hat\Sigma)) - \boldsymbol h' (\blambda_k(\Sigma)) \big\|_\infty
& \ \leq \ c \, \| \hat\Sigma - \Sigma\|_\op \notag \\[0.2cm]
& \ \leq \ c \lambda_1(\Sigma) \sqrt{ \frac{\log(n) \, \beta_q \, {\tt{r}}(\Sigma) }{n} }
\end{align}
hold with probability at least $1 - \frac{c}{n}$. 
Substituting bounds \eqref{trans_step1}, \eqref{trans_step2}, and \eqref{trans_step4}, into \eqref{cmp_bound}, it follows that the event
\begin{equation} \label{trans:all:bound4}
\bigg\|\frac{\boldsymbol h' (\blambda_k(\hat\Sigma))}{(\boldsymbol{\hat{\varsigma}}_k / \boldsymbol{\varsigma}_k )^\tau} - \boldsymbol h'(\blambda_k(\Sigma)) \bigg\|_\infty
\ \leq \ \frac{c \, \log(n) \, \beta_{2q}^3 \, {\tt{r}} (\Sigma)}{n^{1/2} }
\end{equation}
holds with probability at least $1 - \frac{c}{n}$. Lastly, observe that similar reasoning implies that the event
\begin{align} \label{trans_step5}
\big\| \boldsymbol h' (\blambda_k(\hat\Sigma)) \big\|_\infty
& \ \leq \ \big\| \boldsymbol h' (\blambda_k(\Sigma)) \big\|_\infty + \big\| \boldsymbol h' (\blambda_k(\hat\Sigma)) - \boldsymbol h' (\blambda_k(\Sigma)) \big\|_\infty  \ \leq \ c 
\end{align}
holds with probability at least $1 - \frac{c}{n}$. 
By combining the last several bounds with~\eqref{trans:all:bound1}, the proof is complete. \qed

\begin{lemma} \label{lemma:eig_var}
Suppose that the conditions of Theorem~\ref{THM:TRANS_GENERAL} hold. Then,
\begin{equation}\label{eqn:expecbound}
 \max_{1 \leq j \leq k} \big| \E [\lambda_j (\hat\Sigma)] - \lambda_j (\Sigma) \big|
 \ \lesssim \ \frac{ \log (n) \, \beta_{2q} \tr (\Sigma) }{n}\tag{\emph{i}}
\end{equation}
and
\begin{equation}\label{eqn:varbound}
\max_{1 \leq j \leq k} \big| \var \big( \lambda_j(\hat\Sigma) \big) - \lambda_j (\Sigma)^2 \, \Gamma_{jj} / n \big|
\ \lesssim \ \frac{ \log (n) \, \beta_{2q}^2 \,\lambda_1(\Sigma) \tr(\Sigma) }{n^{3/2}}.\tag{\emph{ii}}
\end{equation}
Also, there is a constant $c>0$ not depending on $n$ such that the events
\begin{equation}\label{eqn:starexpecbound}
\max_{1 \leq j \leq k} \big| \E [\lambda_j (\hat\Sigma^\star) | X] - \lambda_j (\hat\Sigma) \big|
 \ \leq \ \frac{c \, \log (n) \, \beta_{2q} \tr (\Sigma) }{n}\tag{\emph{iii}}
\end{equation}
and
\begin{equation}\label{eqn:starvarbound}
\max_{1 \leq j \leq k} \big| \var \big( \lambda_j(\hat\Sigma^\star) \big| X \big) - \lambda_j (\Sigma)^2 \, \hat\Gamma_{jj} / n \big|
\ \leq \ \frac{c \,  \log (n) \, \beta_{2q}^2 \, \lambda_1(\Sigma) \tr(\Sigma) }{n^{3/2}}\tag{\emph{iv}}
\end{equation}
each hold with probability at least $1 - \frac{c}{n}$. 
\end{lemma}

\noindent
\Proof Part~\eqref{eqn:expecbound}: Recalling the choice $q=5\log(kn)$ from the statement of Theorem~\ref{THM:TRANS_GENERAL}, it follows from Lemmas~\ref{lemma:gauss1} and \ref{lemma:gauss2} that there is a constant $c$ not depending on $n$ such that the event
\begin{equation} \label{lemmaB12}
\max_{1 \leq j \leq k} \big| \lambda_j (\hat\Sigma)- d_j (\tilde\Sigma[1,1]) \big|
\ \leq \ \frac{c \, \log (n) \, \beta_{2q} \tr (\Sigma) }{n}
\end{equation}
holds with probability at least $1 - \frac{c}{n^4}$. To simplify presentation, let $\delta$ be a number of the form $\delta =  \frac{c}{n}\log (n)\beta_{2q} \tr (\Sigma)$ and define the event
$$\mathcal{E}_j = \Big\{ | \lambda_j (\hat\Sigma) - d_j (\tilde\Sigma[1,1]) | \leq \delta \Big\}$$
for each $j=1,\dots,k$.
Also, as a temporary short hand, let $U_j$ and $V_j$ denote the random variables
$$U_j = \lambda_j(\hat\Sigma) - \lambda_j(\Sigma) \text{ \ \ \ \ \ \  and \ \ \ \ \ \  } V_j = d_j(\tilde\Sigma[1,1])- \lambda_j(\Sigma).$$ 
Noting that $\E[d_j(\tilde\Sigma[1,1])] = \lambda_j (\Sigma)$, we have
\begin{equation} \label{eig_mean_bound1}
\begin{split}
\big| \E [\lambda_j (\hat\Sigma)] - \lambda_j (\Sigma) \big|
& \ = \ |\E[U_j]-\E[V_j]|\\[0.2cm]
& \ \leq \ \| U_j - V_j \|_2 \\[0.2cm]
& \ \leq \ \big\| |U_j - V_j| \cdot 1\{ \mathcal{E}_j\} \big\|_2 + \big\| |U_j-V_j| \cdot 1\{ \mathcal{E}_j^c\} \big\|_2 \\[0.2cm]
& \ \lesssim \ \delta + \|\lambda_1(\tilde \Sigma)\|_4 \, \big( \P(\mathcal{E}_j^c) \big)^{1/4},
\end{split}
\end{equation}
where the fact $d_j(\tilde\Sigma[1,1])=d_j(\tilde\Sigma)\leq\lambda_1(\tilde\Sigma)=\lambda_1(\hat\Sigma)$ has been used in the last step. Due to \eqref{lemmaB12}, we have $\P (\mathcal{E}_j^c) \leq c / n^4$, and so it is adequate to derive a simple upper bound on $\|\lambda_1(\tilde \Sigma)\|_4$ using \eqref{gauss1:bound_interm} and condition~\eqref{n_large} as follows
\begin{align} 
\|\lambda_1(\tilde\Sigma)\|_4 
& \ \leq \ \lambda_1(\Sigma) + \big\| \| \tilde\Sigma - \Lambda \|_\op \big\|_q, \notag \\[0.2cm]
& \ \leq \ c \, \lambda_1(\Sigma) \label{lambda1_4th}.
\end{align}
Hence, using the last bound in \eqref{eig_mean_bound1} completes the proof of Part~\eqref{eqn:expecbound}.

~\\

\noindent Part~\eqref{eqn:varbound}: First note that
\begin{equation} \label{var_upper_bound}
\begin{split}
\big| \var \big( \lambda_j(\hat\Sigma) \big) - \lambda_j (\Sigma)^2 \, \Gamma_{jj} / n \big|
&  \ = \  |\var(U_j) - \var(V_j)|\\[0.2cm]
& \ \leq \ 2 \,( \| U_j\|_2 + \| V_j \|_2) \, \| U_j - V_j \|_2.
\end{split}
\end{equation}
Also note that the intermediate steps in Part~\eqref{eqn:expecbound} give
$$ \| U_j - V_j \|_2 \ \lesssim \ \frac{\log (n) \, \beta_{2q} \tr (\Sigma)}{n}.$$
From the proof of Lemma~\ref{lemma:gauss3}, we have the identity  
$$V_j  \ = \  d_j(\tilde\Sigma[1,1]) - \lambda_j(\Sigma) \ = \  \frac{1}{n} \sum_{i = 1}^n (\Lambda_k W_i)_j,$$
which leads to
\begin{align*}
\| V_j \|_2^2
& \ = \ \Big\| \frac{1}{n} \sum_{i = 1}^n (\Lambda_k W_i)_j\Big\|_2^2 \ \lesssim \ \frac{\beta_2^2 \, \lambda_1(\Sigma)^2}{n}.
\end{align*}
Also, we can bound $\|U_j\|_2^2$ as follows
\begin{align} \label{bound_Uj}
\|U_j\|_2^2
& \ \leq \ 2\|V_j\|_2^2+2\|U_j-V_j\|_2^2 \notag \\[0.2cm]
& \ \lesssim \ \frac{\beta_2^2\, \lambda_1(\Sigma)^2}{n}+\bigg(\frac{\log(n) \beta_{2q} \tr(\Sigma)}{\sqrt n}\bigg)^2 \ \frac{1}{n} \notag \\[0.2cm]
& \ \lesssim \ \frac{\beta_2^2\,\lambda_1(\Sigma)^2}{n},
\end{align}
where the third step uses the condition~\eqref{n_large}. Combining the last several bounds completes the proof of Part~\eqref{eqn:varbound}.
~\\

\noindent Part~\eqref{eqn:starexpecbound}:
The proof is similar to that of Part~\eqref{eqn:expecbound}. Lemmas~\ref{lemma:boot1} and~\ref{lemma:boot2} imply that the event
\begin{equation} \label{lemmaC12}
\P\bigg(\Big\| \sqrt n\Big(\blambda_k (\tilde{\Sigma}^\star) - \bdiag_k (\tilde{\Sigma}^\star[1,1]) \Big)\Big\|_{\infty}
\, \geq  \, \frac{c \, \log (n) \beta_{2q} \tr (\Sigma) }{n^{1/2}} \, \bigg| \, X\bigg) \ \leq \ \frac{c}{n^4}
\end{equation}
holds with probability at least $1 - \frac{c}{n}$. Letting $\delta$ have the same form as in Part~\eqref{eqn:expecbound}, define the event
\[ \mathcal{E}_j^\star = \Big\{ | \lambda_j (\hat\Sigma^\star) - d_j (\tilde\Sigma^\star[1,1]) | \leq \delta \Big\} \]
for each $j = 1, \dots, k$. Also, define 
\[ U_j^\star = \lambda_j(\hat\Sigma^\star) - d_j(\tilde\Sigma[1,1]) \quad \text{and} \quad V_j^\star = d_j(\tilde\Sigma^\star[1,1]) - d_j(\tilde\Sigma[1,1]) \]
as the bootstrap counterparts of $U_j$ and $V_j$. To proceed, note that $\E[d_j(\tilde\Sigma^{\star}[1,1])|X]=d_j(\tilde\Sigma[1,1])$, and so
%
\begin{align}\label{eqn:bootdiffexpec}
\big| \E [\lambda_j (\hat\Sigma^\star) | X] - \lambda_j (\hat\Sigma) \big|    
&  \ \leq \ |\E[U_j^{\star}|X]-\E[V_j^{\star}|X]| + |d_j(\tilde\Sigma[1,1])-\lambda_j(\hat\Sigma)|.
\end{align}
The first term on the right side can be handled similarly to \eqref{eig_mean_bound1},
\begin{align*}
|\E[U_j^{\star}|X]-\E[V_j^{\star}|X]|
& \ \leq \ \E \big[ ( U_j^\star - V_j^\star)^2 1\{\mathcal{E}_j^{\star}\}\big| X\big]^{1/2} \ + \ \E \big[ ( U_j^\star - V_j^\star)^2 1\{(\mathcal{E}_j^{\star})^c\}\big| X\big]^{1/2} \\[0.2cm]
& \ \leq \ c \Big( \delta + \big( \E \big[ \lambda_1^4(\tilde\Sigma^\star) \, \big| \, X \big] \big)^{1/4} \, \big( \P (\cup_{j=1}^k(\mathcal{E}_j^{\star})^c | X) \big)^{1/4} \, \Big),
\end{align*}
where the fact $d_j(\tilde\Sigma^{\star}[1,1])=d_j(\tilde\Sigma^{\star})\leq\lambda_1(\tilde\Sigma^{\star})$ has been used in the last step. 
Due to \eqref{lemmaC12}, we have $\P(\cup_{j=1}^k(\mathcal{E}_j^{\star})^c|X) \leq c /n^4$ with probability at least $1-c/n$. Furthermore, it follows from Lemma~\ref{tl4} and~\eqref{gauss1:bound2} that the event
\begin{equation}\label{eqn:lambda4norm}
\begin{split}
\big( \E \big[ \lambda_1(\tilde\Sigma^\star)^4 \, \big| \, X \big] \big)^{1/4}
& \ \leq \ \lambda_1(\Sigma) + \|\tilde\Sigma-\Sigma\|_\op+ \big( \E \big[ \| \tilde\Sigma^\star - \tilde\Sigma \|_\op^q \, \big| \, X\big] \big)^{1/q} \\[0.2cm]
& \ \leq \  \lambda_1(\Sigma) + c \, \lambda_1 (\Sigma) \sqrt{ \frac{\log(n) \beta_q {\tt{r}}(\Sigma)}{n} } \\[0.2cm]
& \ \leq \ c \, \lambda_1 (\Sigma)
\end{split}
\end{equation}
holds with probability at least $1 - \frac{c}{n}$, where the last line has used condition~\eqref{n_large}. Thus, the event
$$\max_{1\leq j\leq k}|\E[U_j^{\star}|X]-\E[V_j^{\star}|X]| \ \leq \ \frac{c \, \log (n) \, \beta_{2q} \tr (\Sigma) }{n}$$
holds with probability at least $1-c/n$. For the second term on the right side of~\eqref{eqn:bootdiffexpec}, note that~\eqref{lemmaB12} implies that the event
\begin{equation*}
    \max_{1\leq j\leq k}|d_j(\tilde\Sigma[1,1])-\lambda_j(\hat\Sigma)| \ \leq \ \frac{c \, \log (n) \, \beta_{2q} \tr (\Sigma) }{n}
\end{equation*}
holds with probability at least $1-c/n$. Applying the last several steps  into~\eqref{eqn:bootdiffexpec} completes the proof of Part~\eqref{eqn:starexpecbound}.\\

\noindent Part~\eqref{eqn:starvarbound}: 
The proof is essentially analogous to that of Part~\eqref{eqn:varbound}, and so the details are omitted.\qed

~\\

\begin{lemma} \label{lemma:trans_var}
Suppose that the conditions of Theorem~\ref{THM:TRANS_GENERAL} hold. Then,
\begin{equation} \label{eqn:var_diff_bound}
\max_{1 \leq j \leq k} \big| \var \big(h ( \lambda_j (\hat\Sigma) ) \big) - h' ( \lambda_j (\Sigma) )^2 \var \big( \lambda_j (\hat\Sigma) \big) \big| 
\ \lesssim \ \frac{ \log(n) \, \beta_{2q}^2 \, \lambda_1(\Sigma) \, \tr(\Sigma)}{n^{3/2}}. \tag{\emph{i}}
\end{equation}
In addition, there is a constant $c > 0$ not depending on $n$ such that the events
\begin{equation} \label{eqn:starvar_diff_bound}
\max_{1 \leq j \leq k} \big|  \var \big(h ( \lambda_j (\hat\Sigma^\star) ) \big| X \big) - h' ( \lambda_j (\hat\Sigma) )^2 \var \big( \lambda_j (\hat\Sigma^\star) \big| X \big) \big| 
\ \leq \ \frac{c \,\log(n) \, \beta_{2q}^2 \, \lambda_1(\Sigma) \, \tr(\Sigma)}{n^{3/2}} \tag{\emph{ii}}
\end{equation}
and
\begin{equation} \label{eqn:bootvar_diff_bound}
\max_{1 \leq j \leq k} \big|  \var \big(h ( \lambda_j (\hat\Sigma^\star) ) \big| X \big) - \var \big(h ( \lambda_j (\hat\Sigma) ) \big)  \big| 
\ \leq \ \frac{c \, \log(n) \, \beta_{2q}^{3} \, \lambda_1(\Sigma) \, \tr(\Sigma)}{n^{3/2}} \tag{\emph{iii}}
\end{equation}
each hold with probability at least $1 - \frac{c}{n}$.
\end{lemma}

\noindent
\Proof Part~\eqref{eqn:var_diff_bound}: For each $j=1,\dots,k$ define the random variables
\begin{equation}\label{eqn:SjTjdef}
 S_j = h(\lambda_j(\hat\Sigma)) - h(\lambda_j(\Sigma)) \quad \text{and} \quad T_j = h'(\lambda_j(\Sigma)) \, \big( \lambda_j(\hat\Sigma) - \lambda_j(\Sigma) \big).
 \end{equation}
In this notation, we have
\begin{align} \label{eqn:var_bound}
\max_{1 \leq j \leq k} \big| \var \big(h ( \lambda_j (\hat\Sigma) ) \big) - h' ( \lambda_j (\Sigma) )^2 \var \big( \lambda_j (\hat\Sigma) \big) \big| 
& \ = \ \max_{1 \leq j \leq k} \big| \var(S_j) - \var(T_j) \big| \notag \notag \\[0.2cm]
& \ \leq \ \max_{1 \leq j \leq k} 2 (\| S_j \|_2 + \|T_j \|_2) \| S_j - T_j \|_2.
\end{align}
As a way of handling $\|S_j-T_j\|_2$, first note that the argument used to establish~\eqref{thm:transform:gauss:bound3} can also be used to show that the event 
\begin{align} \label{eqn:trans_general:bound}
\max_{1\leq j\leq k}|S_j-T_j|
& \ \leq \  \frac{c \,\log(n) \, \beta_{2q}\tr(\Sigma)}{n} 
\end{align}
holds with probability at least $1 - \frac{c}{n^4}$.  Using the bound \eqref{eqn:trans_general:bound} and an argument analogous to the one leading up to \eqref{eig_mean_bound1}, we obtain 
\begin{equation} \label{eqn:Sj_Tj_diff_bound}
\| S_j - T_j \|_2 \ \lesssim \ \frac{\log(n) \beta_{2q} \tr(\Sigma)}{n} + \frac{\| S_j \|_4 + \| T_j \|_4}{n} .
\end{equation}
With regard to $\|T_j\|_4$ we use~\eqref{lambda1_4th} to obtain the following conservative but adequate bound,
\begin{align} \label{eqn:Tj_bound}
\|T_j \|_4 
& \ \lesssim \ \| \lambda_j(\hat\Sigma) - \lambda_j(\Sigma) \|_4 \notag \\[0.2cm]
& \ \leq \ \| \lambda_j(\hat\Sigma) \|_4 + \lambda_j(\Sigma) \notag \\[0.2cm]
& \ \lesssim \ \lambda_1(\Sigma).
\end{align}
To handle $\| S_j\|_4$, first note that the concavity of $h$ implies that $S_j\leq T_j$ almost surely, and so
$$0 \ \leq \ h(\lambda_j(\hat\Sigma)) \ \leq \ T_j+h(\lambda_j(\Sigma)).$$
In turn, this yields
\begin{align}\label{eqn:Sj_bound}
\| S_j \|_4 
& = \|h(\lambda_j(\hat\Sigma))-h(\lambda_j(\Sigma))\|_4\notag\\[0.2cm]
& \ \leq \ \| h (\lambda_j(\hat\Sigma)) \|_4  \ + \ h (\lambda_j(\Sigma)) \notag \\[0.2cm]
& \ \leq \|T_j\|_4 \ + \ 2h(\lambda_j(\Sigma))\notag\\[0.2cm]
& \ \lesssim \ \lambda_1(\Sigma).
\end{align}
So, substituting \eqref{eqn:Tj_bound} and \eqref{eqn:Sj_bound} into \eqref{eqn:Sj_Tj_diff_bound} gives
\begin{equation} \label{eqn:Sj_Tj_diff_bound_result}
\| S_j -T_j \|_2 \ \lesssim \ \frac{\log (n) \, \beta_{2q} \tr(\Sigma)}{n}. 
\end{equation}

Now we turn to bounding $\|T_j\|_2$ and $\|S_j\|_2$ in~\eqref{eqn:var_bound}. Using  Lemma~\ref{lemma:eig_var} and the condition~\eqref{n_large} we have
\begin{align*}
\| T_j \|_2 
& \ \lesssim \ \sqrt{\var\big( \lambda_j(\hat\Sigma) \big)} + \big| \E[\lambda_j(\hat\Sigma)] - \lambda_j(\Sigma) \big| \\[0.2cm]
& \ \lesssim \ \frac{\lambda_j(\Sigma) \, \beta_2}{\sqrt{n}} + \Big( \frac{\log(n) \, \beta_{2q} \, {\tt{r}}(\Sigma)}{\sqrt{n}} \Big) \frac{\lambda_1(\Sigma)}{\sqrt{n}} \\[0.2cm]
& \ \lesssim \ \frac{\lambda_1(\Sigma) \, \beta_2}{\sqrt{n}}
\end{align*}
Likewise, we may bound $\|S_j\|_2$ as
\begin{align*}
\| S_j \|_2 
& \ \leq \ \|T_j\|_2 + \| S_j - T_j \|_2 \\[0.2cm]
& \ \lesssim \ \frac{\lambda_1(\Sigma) \, \beta_2}{\sqrt{n}} + \Big( \frac{\log(n) \, \beta_{2q} \, {\tt{r}}(\Sigma)}{\sqrt{n}} \Big) \frac{\lambda_1(\Sigma)}{\sqrt{n}} \\[0.2cm]
& \ \lesssim \ \frac{\lambda_1(\Sigma) \, \beta_2}{\sqrt{n}},
\end{align*}
which completes the proof of Part~\eqref{eqn:var_diff_bound}.\\

\noindent
Part~\eqref{eqn:starvar_diff_bound}: Define random vectors $S^\star, T^\star\in\R^k$ with entries given by
\[ S_j^\star = h(\lambda_j(\hat\Sigma^\star)) - h(\lambda_j(\hat\Sigma)) \quad \text{and} \quad T_j^\star = h'(\lambda_j(\hat\Sigma)) \big( \lambda_j(\hat\Sigma^\star) - \lambda_j(\hat\Sigma) \big) .\]
As an initial step, it can be shown that there is a constant $c>0$ not depending on $n$ such that the event
\begin{equation} \label{eqn:Sj_Tj_prob_bound}
\P \bigg( \| S^\star - T^\star \|_\infty \ \geq \ \frac{c \log(n) \beta_{2q} \tr(\Sigma) }{n} \, \bigg| \, X \bigg) \ \leq \ \frac{c}{n^4}
\end{equation}
holds with probability at least $1 - \frac{c}{n}$. Verifying this is similar to the proof of Lemma~\ref{lemma:boot1}, and can be handled by showing that 
\[ \| S^\star - T^\star \|_\infty \ \leq \ \frac{c \, \log(n) \, \beta_{2q} \, \tr(\Sigma) }{n} \]
holds with probability at least $1 - \frac{c}{n^5}$. In turn, this can be shown using the condition~\eqref{n_large} and the entrywise Taylor expansion
\begin{align*}
S_j^\star - T_j^\star = \big( \lambda_j(\hat\Sigma^\star) - \lambda_j(\hat\Sigma) \big)^2 \, \int_0^1 (1-t) \, h'' \big( \lambda_j (\hat\Sigma) + t (\lambda_j(\hat\Sigma^\star) - \lambda_j(\hat\Sigma)) \big) \textup{d}t.
\end{align*}
along with \eqref{boot1:bound2} and \eqref{boot1:prev}. \\

To proceed with the rest of the proof, we can bound the left side of \eqref{eqn:starvar_diff_bound} in a manner that is similar to~\eqref{eqn:var_bound},
\begin{equation} \label{eqn:starvar_bound}
\small
 \max_{1 \leq j \leq k} \big| \var(S_j^\star|X) - \var(T_j^\star|X) \big|
 \leq \ \max_{1 \leq j \leq k} 2 \big( (\E[(S_j^\star)^2 |X])^{1/2} + (\E[(T_j^\star)^2 |X])^{1/2} \big) \big(\E\big[ (S_j^\star - T_j^\star)^2 \, \big| \, X \big]\big)^{1/2}.
\end{equation}
It follows from an argument analogous to \eqref{eig_mean_bound1} and an application of \eqref{eqn:Sj_Tj_prob_bound} that the event
\[ 
\small
\max_{1\leq j \leq k} \big(\E\big[ (S_j^\star - T_j^\star)^2 \, \big| \, X \big]\big)^{1/2} \ \leq \ \frac{c}{n} \, \bigg(\log(n) \, \beta_{2q} \tr(\Sigma) + \max_{1\leq j\leq k}\Big((\E[(S_j^\star)^4 |X])^{1/4} + (\E[(T_j^\star)^4 |X])^{1/4}\Big) \bigg) \]
holds with probability at least $1-c/n$.
Bounds on the conditional fourth moments can also be derived similarly to the way that the bounds  \eqref{eqn:Tj_bound} and \eqref{eqn:Sj_bound} were. Namely, by using~\eqref{gauss1:bound2} and~\eqref{eqn:lambda4norm}, it can be shown that the events
\[ \max_{1\leq j\leq k}(\E[(T_j^\star)^4 |X])^{1/4} \ \leq \ c \, \lambda_1 (\Sigma) \]
and
\[ \max_{1\leq j\leq k}(\E[(S_j^\star)^4 |X])^{1/4} \ \leq \ c \, \lambda_1 (\Sigma) \]
each hold with probability at least $1 - \frac{c}{n}$. Hence, the event 
\[ \max_{1\leq j\leq k}\big(\E\big[ (S_j^\star - T_j^\star)^2 \, \big| \, X \big]\big)^{1/2} \ \leq \ \frac{c \, \log (n) \, \beta_{2q} \tr(\Sigma) }{n} \]
holds with probability at least $1 - \frac{c}{n}$. \\

To address the conditional $L_2$ norms of $T_j^{\star}$ and $S_j^{\star}$ in~\eqref{eqn:starvar_bound}, first note that
\begin{align*}
(\E[(T_j^\star)^2 |X])^{1/2} 
& \ \leq \ h'(\lambda_j(\Sigma)) \big(  \var(\lambda_j(\hat\Sigma^\star) | X)^{1/2} + \big| \E[\lambda_j(\hat\Sigma^\star) | X] - \lambda_j(\hat\Sigma) \big| \big),
\end{align*}
In turn, Lemmas~\ref{lemma:eig_var} and~\ref{lem:hatgamma} as well as the condition~\eqref{n_large} imply that the event
\begin{align*}
\max_{1\leq j\leq k} (\E[(T_j^\star)^2 |X])^{1/2} 
& \ \leq \ \frac{c \, \lambda_j(\Sigma)\, \beta_2}{\sqrt{n}} + \Big( \frac{\log(n) \beta_{2q} {\tt{r}}(\Sigma) }{\sqrt{n}} \Big) \frac{c \, \lambda_1(\Sigma)}{\sqrt{n}}  \\[0.2cm]
& \ \leq \ \frac{c \, \lambda_1(\Sigma) \, \beta_2}{\sqrt{n}}
\end{align*} 
holds with probability at least $1 - \frac{c}{n}$. Furthermore, this implies that the event
\begin{align*}
\max_{1\leq j\leq k} (\E[(S_j^\star)^2 |X])^{1/2}
& \ \leq \ (\E[(T_j^\star)^2 |X])^{1/2} +  \big(\E\big[ (S_j^\star - T_j^\star)^2 \, \big| \, X \big]\big)^{1/2} \\[0.2cm]
& \ \leq \ \frac{c \, \lambda_1(\Sigma) \, \beta_2}{\sqrt{n}} + \Big( \frac{\log(n) \beta_{2q} {\tt{r}}(\Sigma)}{\sqrt{n}}\Big) \frac{c\, \lambda_1(\Sigma)}{\sqrt{n}} \\[0.2cm]
& \ \leq \ \frac{c \, \lambda_1(\Sigma) \, \beta_2}{\sqrt{n}}
\end{align*}
holds with probability at least $1 - \frac{c}{n}$. The proof is completed by combining the last several steps with~\eqref{eqn:starvar_bound}.  \\

\noindent
Part~\eqref{eqn:bootvar_diff_bound}: Using several applications of the triangle inequality, we have
\begin{align*}
\big| \var\big(h(\lambda_j(\hat\Sigma^\star))\big|X\big) - \var\big(h(\lambda_j(\hat\Sigma)\big) \big|
& \ \leq \  \big| \var\big(h(\lambda_j(\hat\Sigma^\star))\big|X\big) - h'(\lambda_j(\hat\Sigma))^2 \var\big(\lambda_j(\hat\Sigma^\star) \big|X\big) \big| \\[0.2cm] 
& \ \ \ \ \ + \  \big| \var\big(h(\lambda_j(\hat\Sigma)\big) - h'(\lambda_j(\Sigma))^2  \var\big(\lambda_j(\hat\Sigma) \big) \big| \\[0.2cm]
& \ \ \ \ \ + \  \big| h'(\lambda_j(\hat\Sigma))^2 \var\big(\lambda_j(\hat\Sigma^\star) |X\big) - h'(\lambda_j(\Sigma))^2  \var\big(\lambda_j(\hat\Sigma)\big) \big|
\end{align*}
The first and second terms on the right side have been handled by Parts~\eqref{eqn:starvar_diff_bound} and~\eqref{eqn:var_diff_bound} respectively. To handle the third term, we may use the triangle inequality to obtain
\begin{align*}
\big| h'(\lambda_j(\hat\Sigma))^2 \var\big(\lambda_j(\hat\Sigma^\star) \big|X\big) - h'(\lambda_j(\Sigma))^2  \var\big(\lambda_j(\hat\Sigma)\big)\big|
\ \leq \ \textup{M}_j + \textup{M}_j',
\end{align*}
where the two terms on the right are defined as 
\begin{align*}
\textup{M}_j & \ = \ h'(\lambda_j(\hat\Sigma))^2 \cdot \big| \var\big(\lambda_j(\hat\Sigma^\star) \big|X\big) - \var\big(\lambda_j(\hat\Sigma)\big) \big| \\[0.2cm]
\textup{M}_j' & \ = \ \var\big(\lambda_j(\hat\Sigma) \big) \cdot \big|h'(\lambda_j(\hat\Sigma))^2 - h'(\lambda_j(\Sigma))^2 \big| 
\end{align*}
Using \eqref{gauss1:bound2} and the mean value theorem, it can be shown that the event
\[\max_{1\leq j\leq k} \big| h'(\lambda_j(\hat\Sigma))^2 - h'(\lambda_j(\Sigma))^2 \big| \ \leq \ c \, \sqrt{ \frac{ \log(n) \, \beta_q \, {\tt{r}}(\Sigma) }{n } } \]
holds with probability at least $1 - \frac{c}{n}$. Also, using Lemma~\ref{lemma:eig_var} and the condition~\eqref{n_large}, it follows that 
\begin{align*}
\max_{1\leq j\leq k} \var \big(\lambda_j(\hat\Sigma) \big)
& \ \lesssim \ \frac{\lambda_1(\Sigma)^2 \, \beta_{2}^2}{n}.
\end{align*}
Hence, the event
\begin{equation}\label{eqn:Mjprimebound}
 \max_{1\leq j\leq k} \textup{M}_j' \ \leq \ \frac{c \, \lambda_1(\Sigma)^2 \, \log(n) \, \beta_{2q}^{5/2} \, {\tt{r}}(\Sigma) }{n^{3/2}}
 \end{equation}
holds with probability at least $1 - \frac{c}{n}$. \\

Regarding the quantity $\textup{M}_j$, observe that the bound
\begin{align*}
\big| \var(\lambda_j(\hat\Sigma^\star) |X) - \var(\lambda_j(\hat\Sigma)) \big|
& \ \leq \ \big| \var(\lambda_j(\hat\Sigma^\star) |X) - \lambda_j(\Sigma)^2 \hat\Gamma_{jj} / n \big| \\[0.2cm]
& \ \ \ \ \ + \ \ \ \big| \var(\lambda_j(\hat\Sigma)) - \lambda_j(\Sigma)^2 \Gamma_{jj} / n \big| \\[0.2cm]
& \ \ \ \ \ + \ \ \ \frac{\lambda_j(\Sigma)^2 \, \| \hat\Gamma - \Gamma \|_\op}{n}
\end{align*}
holds almost surely.
Consequently, it follows from Lemmas~\ref{lemma:eig_var} and~\ref{lem:hatgamma}, that the event
\begin{equation}\label{eqn:Mjbound}
 \max_{1\leq j\leq k}\textup{M}_j \ \leq \ \frac{c \, \log(n) \, \beta_{2q}^2 \, \lambda_1(\Sigma) \tr(\Sigma)}{n^{3/2}}
 \end{equation}
holds with probability at least $1 - \frac{c}{n}$. Combining~\eqref{eqn:Mjprimebound} and \eqref{eqn:Mjbound} completes the proof of Part~\eqref{eqn:bootvar_diff_bound}.\qed

~\\

\begin{lemma} \label{lemma:trans_var_lower_bound}
Suppose that the conditions of Theorem~\ref{THM:TRANS_GENERAL} hold. Then, 
\[ \min_{1\leq j\leq k} \var \big( h (\lambda_j(\hat\Sigma)) \big)
 \ \gtrsim \ \ts\frac{\lambda_1(\Sigma)^2}{n}, \]
and there is a constant $c \geq 1 $ not depending on $n$ such that  the event
\[ \min_{1\leq j\leq k} \var \big( h (\lambda_j(\hat\Sigma^\star)) \big| X \big) 
\ \geq \ \ts\frac{\lambda_1(\Sigma)^2}{cn} \]
holds with probability at least $1 - \frac{c}{n}$. 
\end{lemma}

\noindent
\Proof For any $j = 1, \dots, k$, and any $t>0$, Chebyshev's inequality and the triangle inequality give
\begin{align} \label{var_lower_step1} 
\var \big( h (\lambda_j(\hat\Sigma)) \big)
& \ \geq \ t^2 \, \P \Big( \big| h (\lambda_j(\hat\Sigma)) - \E\big[h (\lambda_j(\hat\Sigma))\big] \big| \geq t \Big) \notag \\[0.2cm]
& \ \geq t^2 \, \P \Big( \sqrt{n} \big| h (\lambda_j(\hat\Sigma)) - h (\lambda_j(\Sigma)) \big| \geq \sqrt{n} \big( t + \big| \E\big[h (\lambda_j(\hat\Sigma))\big] - h (\lambda_j(\Sigma)) \big| \big) \Big).
\end{align}
Next, let $S_j$ and $T_j$ be as defined in~\eqref{eqn:SjTjdef}. Applying Part~\eqref{eqn:expecbound} of Lemma~\ref{lemma:eig_var} along with \eqref{eqn:Sj_Tj_diff_bound_result} shows that the bounds
\begin{align*}
\big| \E\big[h (\lambda_j(\hat\Sigma))\big] - h (\lambda_j(\Sigma))  \big| 
%
& \ \leq \ |\E[T_j]| + |\E[S_j]-\E[T_j]|\\[0.2cm]
&\ \lesssim \ | h' (\lambda_j(\Sigma))| \, \big|\E \big[ \lambda_j (\hat\Sigma) - \lambda_j(\Sigma) \big] \big|+ \|S_j - T_j \|_2 \\[0.2cm]
& \ \lesssim \frac{ \log(n) \, \beta_{2q} \tr (\Sigma)}{n}, 
\end{align*}
hold for every $j = 1, \dots, k$. Hence, by taking $t = \lambda_1(\Sigma)/\sqrt n$ in \eqref{var_lower_step1}, and using Lemma~\ref{lemma:trans:gauss} with the condition~\eqref{n_large}, there is a constant $c>0$ not depending on $n$ such that
\begin{align*}
\var \big( h (\lambda_j(\hat\Sigma)) \big)
& \ \geq \ t^2 \, \P \Big( \sqrt{n} \big| h (\lambda_j(\hat\Sigma)) - h (\lambda_j(\Sigma)) \big| \geq c \Big) \\[0.2cm]
& \ \gtrsim \ \frac{\lambda_1(\Sigma)^2}{n} \bigg\{ \P \Big( h' (\lambda_j(\Sigma)) \, \zeta_j \geq c \Big) - \frac{ \log(n) \, \beta_{2q}^3 {\tt{r}} (\Sigma)}{n^{1/2} } \bigg\} \\[0.2cm]
& \ \gtrsim \ \frac{\lambda_1(\Sigma)^2}{n},
\end{align*}
for every $j = 1, \dots,k$. Lastly, the proof for the corresponding lower bound on~$\var(h(\lambda_j(\hat\Sigma^{\star}))|X)$ is analogous and so the details are omitted. \qed

\section{Proof of Technical Lemmas }\label{appendix_H}
\begin{lemma}\label{tl1}
Suppose that Assumption~\ref{model_assumptions} holds and let $q\geq 5\log(kn)$. Then,
\[ \big\| \| \tilde{\Sigma}[1,2] \|_\emph{\op}^2 \big\|_q 
\ \lesssim \ 
\frac{q \, \beta_{2q} \, \lambda_1 (\Sigma) \tr (\Sigma)}{ n ^{1 - 3/(2q)}}. \]
\end{lemma}

\noindent
\Proof We will need two auxiliary matrices to extract $\tilde{\Sigma}[1,2]$ from $\tilde{\Sigma}$. Define matrices $\Pi_1\in\R^{k\times p}$ and $\Pi_2\in\R^{p\times (p-k)}$ according to
\begin{align*}
\Pi_1 \ =\ \begin{bmatrix} I_k & \bm{0} \end{bmatrix} 
\quad \text{and} \quad 
\Pi_2 \ =\ \begin{bmatrix} \bm{0} \\ I_{p-k} \end{bmatrix},
\end{align*}
which allow us to write $\tilde{\Sigma}[1,2] = \Pi_1 \tilde{\Sigma} \Pi_2$. Next, let $\mathbb{B}^{k}$ and $\mathbb{B}^{p-k}$ denote the unit $\ell_2$-balls in $\R^k$ and $\R^{p-k}$, and let $\mathbb{T}=\mathbb{B}^{k}\times \mathbb{B}^{p-k}$. With this notation in hand, it follows that
\begin{align}
\| \tilde{\Sigma}[1,2] \|_\op 
& \ = \ \sup_{(u,v)\in\mathbb{T}}  u^\top \Pi_1 \tilde{\Sigma} \Pi_2 v \notag \\[0.2cm]
& \ = \ \sup_{(u,v)\in\mathbb{T}}  (U \Lambda^{1/2} \Pi_1^\top u)^\top \bigg( \frac{1}{n} \sum_{i = 1}^n Z_i Z_i^\top \bigg) (U \Lambda^{1/2} \Pi_2 v ) . \label{eqn:tl1:1}
\end{align}
Observe that the vectors $U \Lambda^{1/2} \Pi_1^\top u$ and $U \Lambda^{1/2} \Pi_2 v$ both lie in the ellipsoid $ \mathcal{E}:=U \Lambda^{1/2} (\mathbb{B}^p)$, and are orthogonal. Therefore,
\begin{equation}\label{eqn:nonquadsum}
 \| \tilde{\Sigma}[1,2] \|_\op \ \leq \ \sup_{ \substack{\small (\mathsf{u}, \mathsf{v})\in\mathcal{E}\times \mathcal{E}\\ \langle \mathsf{u},\mathsf{v}\rangle=0} } \frac{1}{n} \sum_{i = 1}^n \langle \mathsf{u}, Z_{i} \rangle \langle \mathsf{v}, Z_{i} \rangle.
 \end{equation}
It will be convenient to write the summands in terms of the matrix  $\mathcal{Q}(\mathsf{u},\mathsf{v}) :=\frac{1}{2} (\mathsf{u}\mathsf{v}^\top + \mathsf{v}\mathsf{u}^\top)$, namely
\[ \langle \mathsf{u}, Z_{ i} \rangle \langle \mathsf{v}, Z_{ i} \rangle = Z_{ i}^\top \mathcal{Q}(\mathsf{u},\mathsf{v})  Z_{ i} .\]
Under this definition, it can be checked that if $\mathsf{u}$ and $\mathsf{v}$ are any pair of orthogonal vectors, then $\E [Z_{ i}^\top \mathcal{Q}(\mathsf{u},\mathsf{v})Z_{i}] = 0$. 
Hence, if we subtract $\E [Z_{ i}^\top \mathcal{Q}(\mathsf{u},\mathsf{v})Z_{i}]$ from the $i$th term in~\eqref{eqn:nonquadsum} for each $i=1,\dots,n$, and subsequently drop the constraint $\langle \mathsf{u},\mathsf{v}\rangle=0$ from the supremum,  then we obtain the bound
\begin{align}\label{eqn:newcenteredbound}
\| \tilde{\Sigma}[1,2] \|_\op 
& \ \leq \  \sup_{ \substack{(\mathsf{u}, \mathsf{v})\in\mathcal{E}\times\mathcal{E} }} \bigg( \frac{1}{n} \sum_{i = 1}^n  Z_{ i}^\top \mathcal{Q}(\mathsf{u},\mathsf{v})  Z_{i} - \E [Z_i^\top \mathcal{Q}(\mathsf{u},\mathsf{v})  Z_i]  \bigg) .
\end{align}
Next, for a given pair of vectors $\mathsf{u}, \mathsf{v}\in\mathcal{E}$, define two associated vectors
\begin{align*}
\mathsf{w}  & \ = \ \ts \mathsf{u}/2 + \ts \mathsf{v}/2, \\
\tilde{\mathsf{w}}  & \ = \ \ts \mathsf{u}/2 - \ts \mathsf{v}/2,
\end{align*}
which satisfy the algebraic relation
\[\mathcal{Q}(\mathsf{u},\mathsf{v})  = \mathsf{w} \mathsf{w}^\top - \tilde{\mathsf{w}} \tilde{\mathsf{w}}^\top. \]
To proceed, define the matrix $A\in\R^{p\times2p}$ as the column concatenation $A = \begin{pmatrix} U \Lambda^{1/2}, U \Lambda^{1/2} \end{pmatrix}$, and note that both $\mathsf{w}$ and $\tilde{\mathsf{w}}$ lie in the ellipsoid $\mathcal{E}':=A (\mathbb{B}^{2p})$. As a result, if we write $\xi_i = A^\top Z_i$, then we have
\begin{align}
\| \tilde{\Sigma}[1,2] \|_\op 
& \ \leq \ \sup_{(\mathsf{w}, \tilde{\mathsf{w}})\in\mathcal{E}'\times\mathcal{E}'} \bigg( \bigg| \frac{1}{n} \sum_{i = 1}^n \langle Z_i, \mathsf{w} \rangle^2 - \E [\langle Z_i, \mathsf{w} \rangle^2]  \bigg| + \bigg| \frac{1}{n} \sum_{i = 1}^n \langle Z_i, \tilde{\mathsf{w}} \rangle^2 - \E [\langle Z_i, \tilde{\mathsf{w}} \rangle^2]  \bigg| \bigg) \notag \\[0.2cm]
& \ \leq \  \sup_{w \in \mathbb{B}^{2p} } \bigg| \frac{2}{n} \sum_{i = 1}^n \langle A^\top Z_i, w \rangle^2 - \E [\langle A^\top Z_i, w \rangle^2] \bigg| \notag \\[0.2cm]
& \ = \ \bigg\| \frac{2}{n} \sum_{i = 1}^n \xi_i \xi_i^\top - \E [\xi_1 \xi_1^\top] \bigg\|_\op. \label{eqn:tl1:2}
\end{align}
Next, Lemma \ref{background2} implies that
\begin{equation} \label{eqn:tl1:3}
\footnotesize
\bigg\| \Big\| \frac{1}{n} \sum_{i = 1}^n \xi_i \xi_i^\top - \E [\xi_i \xi_i^\top] \Big\|_\op \bigg\|_q 
\ \leq \ 
c \, \bigg( \sqrt{ \frac{q}{n^{1-3/q}} } \, \big\| \E [\xi_1  \xi_1^\top] \big\|_\op^{1/2}  \, \big( \E \| \xi_1  \|_2^{2q} \big)^{ \frac{1}{2q} } \bigg) \bigvee  \bigg( \frac{q}{n^{1-3/q}} \, \big( \E \| \xi_1  \|_2^{2q} \big)^{1/q} \bigg)  .
\end{equation}
We can further compute
\begin{align} \label{eqn:tl1:4}
\big\| \E [\xi_1  \xi_1^\top] \big\|_\op 
& \ =\ \big\| A^\top \E [Z_1 Z_1^\top] A \big\|_\op \notag \\[0.2cm]
& \ = \ \| A \|_\op^2 \notag \\[0.2cm]
& \ \lesssim \  \lambda_1(\Sigma),
\end{align}
and
\begin{align} \label{eqn:tl1:5}
\big( \E \| \xi_1 \|_2^{2q} \big)^{1/q}
& \ =\ 2 \big\| Z_1^\top U \Lambda U^\top Z_1 \big\|_q \notag \\[0.2cm]
& \ =\ 2 \bigg\| \sum_{j = 1}^p \lambda_j(\Sigma) \langle u_j ,Z_1 \rangle^2 \bigg\|_q \notag \\[0.2cm]
& \ \lesssim \ \tr (\Sigma) \beta_q .
\end{align}
To finish, we use the relation $\big\|\|\tilde\Sigma_{12}\|_\op^2\big\|_q=\big\|\|\tilde\Sigma_{12}\|_\op\big\|_{2q}^2$. Specifically, the previous two bounds can be substituted into \eqref{eqn:tl1:2} and \eqref{eqn:tl1:3} while replacing $q$ with $2q$ and using the condition~\eqref{n_large}.\qed
~\\

\begin{lemma} \label{tl3} 
Suppose that Assumption~\ref{model_assumptions} holds and let $q\geq 5\log(kn)$. Then, there is a constant $c > 0$ not depending on $n$ such that the event 
\[ \Big( \E \big[ \| \tilde{\Sigma}^\star[1,2] \|_\textup{op}^{2q} \, \big| \, X \big] \Big)^{1/q}
\ \leq \ \frac{c \, q \, \beta_{2q} \, \lambda_1(\Sigma) \, \tr (\Sigma) }{n^{1-3/(2q)}} \]
holds with probability at least $1 - ce^{-q}$.
\end{lemma}

\noindent

\noindent \Proof We will follow the same notation that was used in the proof of Lemma \ref{tl1}. Repeating the argument from that proof up to~\eqref{eqn:newcenteredbound} and using $2q$ in place of $q$ we have
\begin{align}
\Big(\E\big[\| \tilde{\Sigma}^\star[1,2] \|_\op^{2q}\big|X\big]\Big)^{\frac{1}{2q}}
& \ \leq \ \Bigg(\E\bigg[ \bigg(\sup_{ \substack{(\mathsf{u},\mathsf{v})\in\mathcal{E}\times\mathcal{E} }} \frac{1}{n} \sum_{i =1}^n (Z_i^\star)^\top \Q(\mathsf{u},\mathsf{v}) Z_i^\star -\E [ Z_i^\top \Q(\mathsf{u},\mathsf{v}) Z_i ]\bigg)^{2q}\bigg|X\bigg]\Bigg)^{\frac{1}{2q}} \notag
\end{align}
Next, observe that %
$$\E \big[ (Z_1^\star)^\top \Q(\mathsf{u},\mathsf{v}) Z_1^\star \, \big| \, X \big] \ = \  \frac{1}{n} \sum_{i =1}^n Z_i^\top \Q(\mathsf{u},\mathsf{v}) Z_i $$ 
and so the triangle inequality for the conditional $L_{2q}$ norm gives
\small
\begin{align}
\Big(\E\big[\| \tilde{\Sigma}^\star[1,2] \|_\op^{2q}\big|X\big]\Big)^{\frac{1}{2q}}
& \ \leq\  \Bigg(\E\bigg[\bigg(\sup_{\substack{(\mathsf{u},\mathsf{v})\in\mathcal{E}\times\mathcal{E}}}  \frac{1}{n} \sum_{i =1}^n (Z_i^\star)^\top \Q(\mathsf{u},\mathsf{v}) Z_i^\star  - \E \big[ (Z_i^\star)^\top \Q(\mathsf{u},\mathsf{v}) Z_i^\star \, \big| \, X \big]\bigg)^{2q}\bigg|X\bigg]\Bigg)^{\frac{1}{2q}}\label{eqn:lemma:bootstrap1}\\[0.2cm]
&  \ \ \ \ \  \ \ \  + \sup_{ \substack{(\mathsf{u},\mathsf{v})\in\mathcal{E}\times\mathcal{E} } } \bigg( \frac{1}{n} \sum_{i =1}^n Z_i^\top \Q(\mathsf{u},\mathsf{v}) Z_i  - \E [ Z_i^\top \Q(\mathsf{u},\mathsf{v}) Z_i ] \bigg). \notag
\end{align}
\large
With regard to the second term in the last bound, the proof of Lemma \ref{tl1} shows (via Chebyshev's inequality and condition~\eqref{n_large}) that the event
\begin{equation} \label{eqn:lemma:bootstrap2}
 \sup_{\substack{(\mathsf{u},\mathsf{v})\in\mathcal{E}\times\mathcal{E} }} \bigg( \frac{1}{n} \sum_{i =1}^n Z_i^\top\Q(\mathsf{u},\mathsf{v}) Z_i  - \E [ Z_1^\top\Q(\mathsf{u},\mathsf{v}) Z_1 ] \bigg) 
\ \leq \ 
c\sqrt{ \frac{q \beta_q \lambda_1(\Sigma) \tr (\Sigma)}{n^{1 - 3/q}} }.
\end{equation}
holds with probability at least $1-ce^{-q}$, for some constant $c>0$ not depending on $n$.
Therefore, the proof of the current lemma is complete once we derive a similar bound for the first term on the right side of~\eqref{eqn:lemma:bootstrap1}. \\

Let $\xi_i^\star = A^\top Z_i^\star$, with $A$ as defined in the proof of Lemma~\ref{tl1}. Then, by following the argument leading up to~\eqref{eqn:tl1:2} and applying Lemma \ref{background2}, there is a constant $c>0$ not depending on $n$ such that
\footnotesize{
\begin{align}
 \, \bigg( \E\bigg[ \bigg( \sup_{\substack{(\mathsf{u},\mathsf{v})\in\mathcal{E}\times\mathcal{E} }}  &\frac{1}{n} \sum_{i =1}^n (Z_i^\star)^\top \Q(\mathsf{u},\mathsf{v})  Z_i^\star  - \E \big[ (Z_1^\star)^\top \Q(\mathsf{u},\mathsf{v}) Z_1^\star \, \big| \, X \big]  \bigg)^{2q} \, \bigg| \, X \bigg] \bigg)^{\frac{1}{2q}} \notag \\[0.4cm]
& \ \leq \ c \,\bigg( \E \bigg[  \bigg\| \frac{1}{n} \sum_{i = 1}^n \xi_i^\star (\xi_i^\star)^\top - \E \big[ \xi_1^\star (\xi_1^\star)^\top \big| \, X \big] \bigg\|_\op^{2q} \, \bigg| \, X \bigg] \bigg)^{\frac{1}{2q}} \notag \\[0.3cm]
& \ \leq \ c\, \bigg( \sqrt{ \frac{ 2q }{n^{1 - \frac{3}{2q}} }} \cdot \big( \big\| \E [ \xi_1^\star (\xi_1^\star)^\top \, | \, X ] \big\|_\op \big)^{ 1/2 } \cdot \big( \E \big[ \| \xi_1^\star \|_2^{4q} \, \big| \, X \big] \big)^{ \frac{1}{4q} }  \bigg) \bigvee \bigg( \frac{q}{n^{1-\frac{3}{2q}}} \cdot \big( \E \big[ \| \xi_1^\star \|_2^{4q} \, \big| \, X \big] \big)^{ \frac{1}{2q} }  \bigg). \label{eqn:lemma:bootstrap3}
\end{align}}
\large
Next, the triangle inequality implies
\begin{align*}
\big\| \E \big[ \xi_1^\star (\xi_1^\star)^\top \, \big| \, X \big] \big\|_\op 
& \ =\ \bigg\| \frac{1}{n} \sum_{i = 1}^n (A^\top Z_i) (A^\top Z_i)^\top \bigg\|_\op \\[0.2cm]
& \ \leq \ \bigg\| \frac{1}{n} \sum_{i =1}^n (A^\top Z_i) (A^\top Z_i)^\top - A^\top A \bigg\|_\op + 4\lambda_1(\Sigma),
 \end{align*}
where we have used $\| A^\top A\|_\op \leq 4 \lambda_1(\Sigma)$ and $\E [ (A^\top Z_1) (A^\top Z_1)^\top] = A^\top A$. By applying Lemma \ref{background2} to the first term (along with the condition~\eqref{n_large} and Chebyshev's inequality), it follows that the event
\begin{equation} \label{eqn:lemma:bootstrap4}
\big\| \E \big[ \xi_1^\star (\xi_1^\star)^\top \, \big| \, X \big] \big\|_\op \leq c \lambda_1 (\Sigma)
\end{equation}
holds with probability at least $1 - e^{-q}$, for some constant $c > 0$ not depending on $n$. \\

It remains to develop an upper bound for $\big( \E \big[ \| \xi_1^\star \|_2^{4q} \, \big| \, X \big]  \big)^{\frac{1}{2q}}$. According to the definition of $\xi_i^\star$, we have
\begin{align*}
\| \xi_1^\star \|_2^2 
& \ = \ (Z_1^\star)^\top A A^\top Z_1^\star \\[0.2cm]
& \ =\ 2 (Z_1^\star)^\top \Sigma Z_1^\star ,
\end{align*}
and so
\begin{align*}
\big( \E \big[ \| \xi_1^\star \|_2^{4q} \, \big| \, X \big]  \big)^{\frac{1}{2q}} 
& \ =\ 2\bigg( \frac{1}{n} \sum_{i = 1}^n  (Z_i^\top \Sigma Z_i)^{2q} \bigg)^{ \frac{1}{2q} }\\[0.2cm]
& \ = \ 2S,
\end{align*}
where the non-negative random variable $S$ is defined by the last line. For any $t>0$, Chebyshev's inequality implies
\begin{equation}\label{eqn:normstarbound}
\begin{split}
 \P (S \geq e \, t ) & \ \leq \  \frac{e^{-2q} \, \E [S^{2q}]}{t^{2q}}\\[0.2cm]
 & \ = \  \frac{e^{-2q} \, \| Z_1^\top \Sigma Z_1 \|_{2q}^{2q} }{t^{2q}}\\[0.2cm]
 & \ =\ e^{-2q} \, \Bigg(\frac{1}{t}\bigg\|\displaystyle \sum_{j = 1}^p \lambda_j(\Sigma) \langle u_j, Z_1 \rangle^2 \bigg\|_{2q}\Bigg)^{2q}\\[0.2cm]
 & \ \leq \  e^{-2q} \, \Big(\frac{\beta_{2q}  \tr (\Sigma)}{t}\Big)^{2q}.
 \end{split}
 \end{equation}
Combining the last few steps and using the choice $t = \beta_{2q} \tr (\Sigma)$, it follows that there is a constant $c >0 $ not depending on $n$ such that the event 
\begin{equation} \label{eqn:lemma:bootstrap5}
\big( \E \big[ \| \xi_1^\star \|_2^{4q} \, \big| \, X \big]  \big)^{\frac{1}{2q}} \ \leq \ c \, \beta_{2q} \, \tr (\Sigma)
\end{equation}
holds with probability at least $1 - e^{-2q}$. Hence, we may substitute  \eqref{eqn:lemma:bootstrap4} and \eqref{eqn:lemma:bootstrap5}  into \eqref{eqn:lemma:bootstrap3}, and use the condition~\eqref{n_large} to conclude that the event
\begin{equation}\label{eqn:bigstep}
\small
\bigg( \E\bigg[ \bigg( \sup_{ \substack{(\mathsf{u},\mathsf{v})\in\mathcal{E}\times\mathcal{E} }} \frac{1}{n} \sum_{i =1}^n (Z_i^\star)^\top \Q(\mathsf{u},\mathsf{v})  Z_i^\star  - \E \big[ (Z_1^\star)^\top \Q(\mathsf{u},\mathsf{v}) Z_1^\star \, \big| \, X \big]  \bigg)^{2q} \, \bigg| \, X \bigg] \bigg)^{\frac{1}{2q}} 
\ \leq \ 
c\, \sqrt{ \frac{2q \, \beta_{2q} \, \lambda_1(\Sigma) \, \tr (\Sigma) }{n^{ 1- \frac{3}{2q}}} }
\end{equation}
%
%
%
%
holds with probability at least $1 - c e^{-q}$. Finally, the proof is completed by combining~\eqref{eqn:bigstep} and~\eqref{eqn:lemma:bootstrap2} with~\eqref{eqn:lemma:bootstrap1}.
\qed

~\\

\begin{lemma} \label{tl4}
Suppose that Assumption~\ref{model_assumptions} holds and let $q\geq 5\log(kn)$.  Then, there is a constant $c > 0$ not depending on $n$, such that the event 
\[ \big( \E \big[ \| \tilde{\Sigma}^\star - \tilde{\Sigma}  \|_\textup{op}^q \, \big| \, X \big] \big)^{1/q} \ \leq \ c \, \sqrt{  \frac{q \, \beta_q \, \lambda_1(\Sigma) \, \tr (\Sigma) }{n^{1 - 3/q}}  } \]
holds with probability at least $1 - ce^{-q}$.
\end{lemma}

\noindent
\Proof Letting $\xi_1^\star = \Lambda^{1/2} U^\top Z_1^\star$ and using Lemma \ref{background2}, we have
\begin{equation}\label{eqn:firststepnormdiff}
\small
\big( \E \big[ \| \tilde{\Sigma}^\star - \tilde{\Sigma}  \|_\op^q \, \big| \, X \big] \big)^{1/q} 
\ \leq \ c \cdot \bigg( \sqrt{ \frac{q}{n^{1 - 3/q} }} \, \| \tilde{\Sigma} \|_\op^{1/2} \, \big( \E \big[ \| \xi_1^\star \|_2^{2q} \, \big| \, X \big] \big)^{\frac{1}{2q}} \bigg)  \bigvee \bigg( \frac{q}{n^{1-3/q}} \, \big( \E \big[ \| \xi_1^\star \|_2^{2q} \, \big| \, X \big] \big)^{\frac{1}{q}}  \bigg).
\end{equation}
%
Using \eqref{gauss1:bound_interm} in the proof of Lemma \ref{lemma:gauss1}, as well as the condition~\eqref{n_large}, it follows that the bound
\begin{align} \label{eqn:tl4:1}
\| \tilde{\Sigma} \|_\op
& \ \leq \ \| \tilde{\Sigma} - \Lambda \|_\op + \| \Lambda \|_\op \notag \\[0.2cm]
& \ \leq \ c \, \lambda_1 (\Sigma) \sqrt{ \frac{q \, \beta_q \, {\tt{r}} (\Sigma)}{n^{1 - 3/q} } } + \lambda_1(\Sigma) \notag \\[0.2cm]
& \ \leq \ c \, \lambda_1(\Sigma)
\end{align}
holds with probability at least $1 - e^{-q}$. Also, recall that the argument leading up to~\eqref{eqn:lemma:bootstrap5} implies that the event
\begin{equation}\label{eqn:normstaragain}
    \big(\E \big[ \| \xi_1^\star \|_2^{2q} \, \big| \, X \big] \big)^{\frac{1}{2q}} \ \leq \ c\sqrt{\beta_{q}\tr(\Sigma)}
\end{equation}
holds with probability at least $1-e^{-q}$, for some constant $c>0$ not depending on $n$. Combining~\eqref{eqn:tl4:1} and~\eqref{eqn:normstaragain} with~\eqref{eqn:firststepnormdiff} and the condition~\eqref{n_large} completes the proof.
\qed

\begin{lemma}\label{lem:hatgamma}
Suppose that Assumption~\ref{model_assumptions} holds and let $q\geq 5\log(kn)$. Then, there is a constant $c>0$ not depending on $n$ such that the bound
\begin{equation}
    \|\hat\Gamma -\Gamma\|_\op \ \leq \ \frac{c \, q \, \beta_{2q}^2}{n^{1/2}}
\end{equation}
holds with probability at least $1-e^{-q}$.
\end{lemma}

\noindent \Proof
Let $W_1,\dots,W_n$ and $\bar W$ be as defined at the beginning of Appendix~\ref{appendix_C}.
Note that $\Gamma=\E[W_1W_1\ttop]$, and that the definition of $\hat\Gamma$ in~\eqref{eqn:hatgammadef} gives
$$\hat\Gamma = \frac{1}{n}\sum_{i=1}^nW_iW_i\ttop - \bar W\bar W\ttop.$$
We may apply Lemma \ref{background2} to obtain
\footnotesize
\begin{align} \label{eqn:tl6:1}
\Big( \E \big\| \hat\Gamma - \Gamma \big\|_\op^q \Big)^{1/q}
& \ \leq \ \bigg( \E \bigg\| \frac{1}{n} \sum_{i = 1}^n W_i W_i^\top - \E [W_1 W_1^\top] \bigg\|_\op^q \bigg)^{1/q} \ + \ \Big(\E\|\bar W\bar W\ttop\|_\op^q\Big)^{1/q} \notag \\[0.2cm]
& \ \lesssim \ \bigg( \sqrt{ \frac{q}{n^{1-3/q} } } \, \big( \| \E[ W_1 W_1^\top] \|_\op \big)^{1/2} \, \big( \E \| W_1 \|_2^{2q} \big)^{ 1/ (2q) } \bigg) \bigvee \bigg( \frac{q}{n^{1-3/q}} \, \big( \E \| W_1 \|_2^{2q} \big)^{ 1/q }  \bigg) \ + \ \Big(\E\|\bar W\|_2^{2q}\Big)^{1/q}.
\end{align}
\large
It is straightforward to verify 
that 
\[ \| \E[ W_1 W_1^\top] \|_\op \ \leq \ \tr(\Gamma) \ \lesssim \ \beta_2^2.\]
Also, we have
\begin{align*}
\big( \E \| W_1 \|_2^{2q} \big)^{1/q}
& \ = \ \bigg\| \sum_{j = 1}^k \big( \langle u_j, Z_1 \rangle^2 - 1 \big)^2 \bigg\|_q \\[0.2cm]
& \ \leq \ \sum_{j = 1}^k \big\| \langle u_j, Z_1 \rangle^2 - 1  \big\|_{2q}^2 \\[0.2cm]
& \ \lesssim \, \beta_{2q}^2.
\end{align*}

To handle the term involving $\bar W$, observe that the previous bound and Lemma~\ref{background1} lead to
\begin{equation} \label{Wbar}
\begin{split}
    \Big(\E\|\bar W\|_2^{2q}\Big)^{\frac{1}{2q}} & \ \lesssim \ q\Bigg(\big(\E\|\bar W\|_2^2\big)^{1/2}+ \Big(\ts\sum_{i=1}^n \E \|\ts\frac{1}{n}W_i\|_2^{2q}\Big)^{\frac{1}{2q}}\Bigg) \\[0.2cm]
    & \ \lesssim \ q\Bigg(\ts\frac{\big(\E\|W_1\|_2^2\big)^{1/2}}{n^{1/2}}+\ts\frac{\big(\E\|W_1\|_2^{2q}\big)^{\frac{1}{2q}}}{n^{1-\frac{1}{2q}}}\Bigg)\\[0.2cm]
    & \ \lesssim  \ q\Bigg(\ts\frac{\beta_{2}}{n^{1/2}}+\ts\frac{\beta_{2q}}{n^{1-\frac{1}{2q}}}\Bigg)\\[0.2cm]
    & \ \lesssim \ \frac{q \, \beta_{2q}}{n^{1/2}}
\end{split}
\end{equation}
Using condition~\eqref{n_large}, we may substitute the last few bounds into \eqref{eqn:tl6:1} to obtain
\begin{equation} \label{tl6:bound}
\Big( \E \big\| \hat\Gamma - \Gamma \big\|_\op^q \Big)^{1/q} \ \lesssim \ \ \frac{q \, \beta_{2q}^2}{n^{1/2}}
\end{equation}
Combining this result with Chebyshev's inequality completes the proof.

\qed

~\\

For the next lemma, let $Y_1^{\star}$ be as defined at the beginning of the proof of Lemma~\ref{lemma:boot3}.

\begin{lemma} \label{tl6}
Suppose that Assumption~\ref{model_assumptions} holds and let $q\geq 5\log(kn)$. Then, there is a constant $c>0$ not depending on $n$, such that the event
\[ \E \Big[ \big\| (\Lambda_k \hat\Gamma \Lambda_k)^{-1/2} Y_1^\star \big\|_2^3 \, \Big| \, X \Big] \ \leq \ c \, \beta_{3q}^3 .\]
holds with probability at least $1 - c e^{-q}$. 
\end{lemma}

\noindent
\Proof 
Note that
\begin{align*}
\big\| (\Lambda_k \hat\Gamma \Lambda_k)^{-1/2} Y_1^\star \big\|_2^2 
& \ = \ (W_1^\star-\bar{W})^\top \hat\Gamma^{-1} (W_1^\star-\bar{W}) \\[0.2cm]
& \ \leq \ \frac{ (W_1^\star - \bar{W})^\top (W_1^\star - \bar{W})}{\lambda_k (  \hat\Gamma)}\\[0.2cm]
& \ \leq \ \frac{2}{\lambda_k(\hat\Gamma)} \bigg( \sum_{j = 1}^k \big( \langle u_j, Z_1^\star \rangle^2 - 1 \big)^2 + \bar{W}^\top \bar{W} \bigg). 
\end{align*}
This implies
\begin{equation}\label{eqn:lambdakhatgamma}
    \begin{split}
        \E\Big[\big\| (\Lambda_k \hat\Gamma \Lambda_k)^{-1/2} Y_1^\star \big\|_2^3\Big|X\Big] & \ \leq \ \frac{c}{\lambda_k(\hat\Gamma)^{3/2}}\bigg(\frac{1}{n}\sum_{i=1}^n\Big( \sum_{j = 1}^k \big( \langle u_j, Z_i \rangle^2 - 1 \big)^2\Big)^{3/2} + (\bar{W}^\top \bar{W})^{3/2} \bigg)\\[0.2cm]
        & \ = \ \frac{c \, \big(T + (\bar{W}^\top \bar{W})^{3/2} \big) }{\lambda_k(\hat\Gamma)^{3/2}},
    \end{split}
\end{equation}
where the first line has used the convexity of the function $x\mapsto  x^{3/2}$, and the non-negative random variable $T$ is defined in the second line. By the triangle inequality for the $L^q$ norm, we have
\begin{equation}\label{eqn:Tbound}
\begin{split}
    \|T\|_q & \ \leq \Bigg(\bigg\|\sum_{j = 1}^k \big( \langle u_j, Z_1 \rangle^2 - 1 \big)^2\bigg\|_{3q/2}^{3q/2}\Bigg)^{1/q}\\[0.3cm]
    & \ \lesssim \ \max_{1\leq j\leq k} \Big\|\langle u_j,Z_1\rangle^2-1\Big\|_{3q}^{3}\\[0.2cm]
    & \ \lesssim \ \beta_{3q}^3.
    \end{split}
\end{equation}
So, by Chebyshev's inequality, there is a constant $c>0$ not depending on $n$ such that the event
$$T \, \leq \, c \, \beta_{3q}^3$$
holds with probability at least $1-e^{-q}$. The bound \eqref{Wbar} in the proof of Lemma~\ref{lem:hatgamma} and the condition~\eqref{n_large} also imply that there is a constant $c>0$ not depending on $n$ such that the bound
\[ (\bar{W}^\top \bar{W})^{3/2} \ \leq \ \Big(\frac{c \, q \, \beta_{2q}  }{n^{1/2}}\Big)^{3/2} \ \leq \ c \]
holds with probability at least $1 - e^{-q}$. \\

Now we turn to showing that $\lambda_k \big( \hat\Gamma \big)$ is greater than a positive constant with high probability. Due to Weyl's inequality
we have
$\lambda_k(\hat\Gamma) 
 \ \geq \ \lambda_k(\Gamma) - \|\hat\Gamma-\Gamma\|_\op$.
Using Lemma~\ref{lem:hatgamma}, Assumption~\ref{model_assumptions}.\ref{A3}, and the condition~\eqref{n_large}, it follows that there is a constant $c>0$ not depending on $n$ such that
the bound
\begin{equation}\label{eqn:lambdakhatgammafinally}
\lambda_k(\hat\Gamma)\geq c,
\end{equation}
holds with probability at least $1-e^{-q}$.
Combining the bounds~\eqref{eqn:Tbound} and \eqref{eqn:lambdakhatgammafinally} with~\eqref{eqn:lambdakhatgamma} completes the proof.
\qed

\section{Background Results}\label{sec:background}

\begin{lemma}[Theorem 1 in \cite{Talagrand89}] \label{background1}
Let $X_1,\dots,X_n$ be independent centered random elements of a Banach space with norm $\|\cdot\|$. Then, there is an absolute constant $c>0$ such that the following inequality holds for any $q\geq 1$,
\begin{align*}
\bigg( \E \Big\| \sum_{i = 1}^n X_i \Big\|^q \bigg)^{1/q} 
\ \leq \ \frac{cq}{1 + \log(q)}  \bigg(\E \Big\| \sum_{i = 1}^n X_i \Big\|  + \Big( \E  \displaystyle\max_{1 \leq i \leq n} \| X_i \|^q \Big)^{1/q}\bigg).
\end{align*}
\end{lemma}

~\\

\begin{lemma}[Weilandt's inequality \cite{EatonT91,Wielandt67}] \label{background0}
Consider a real symmetric $p\times p$ matrix 
\begin{equation*}
A = \begin{pmatrix} B & C \\ C^\top & D \end{pmatrix} ,    
\end{equation*}
where $B$ is $k \times k$ and $D$ is $(p - k) \times (p- k)$.
If $\lambda_k (B) > \lambda_1 (D)$, then
\[ 0 \leq \lambda_j (A) - \lambda_j (B) \leq \frac{\lambda_1 (C C^\top )}{\lambda_j (B) - \lambda_1 (D) }, \quad j = 1, \dots, k\]
and
\[ 0 \leq \lambda_{p - k - i} (D) - \lambda_{p - i} (A) \leq \frac{\lambda_1 (C C^\top )}{\lambda_k (B) - \lambda_{p - k - i} (D) }, \quad i = 0, \dots, p - k - 1 .\]
\end{lemma}
~\\

\begin{lemma}[Proposition 1 in \cite{LopesEM19}] \label{background2}
Let $\xi , \dots, \xi_n \in \R^p$ be i.i.d. random vectors, let $q\geq 3$, and define the quantity
\[ r(q) = q \cdot \frac{\big( \E \| \xi_1 \|_2^{2q} \big)^{1/q}}{\big\| \E [\xi_1 \xi_1^\top] \big\|_\op} .\]
Then, there is an absolute constant $c > 0$ such that
\[ \bigg( \E \bigg\| \frac{1}{n} \sum_{i = 1}^n \xi_i \xi_i^\top - \E [\xi_1 \xi_1^\top] \bigg\|_\op^q \bigg)^{1/q} \ \leq \ c \cdot \big\| \E [\xi_1 \xi_1^\top] \big\|_\op \cdot \bigg(\sqrt{ \frac{r(q)}{n^{1-3/q}}} \bigvee \frac{r(q)}{n^{1-3/q} } \bigg) \]
\end{lemma}

The following anti-concentration lemma originates from~\cite{Nazarov03}, and was further elucidated in \cite[Theorem 1]{ChernozhukovCK17}.

\begin{lemma} [Nazarov's inequality] \label{background4}
Let $Y = (Y_1, \dots, Y_p)$ be a centered Gaussian random vector in $\R^p$ and suppose that the parameter $\underline{\sigma}^2= \min_{1\leq j\leq p}\E [Y_j^2] $ is positive. Then for every $y \in \R^p$ and $\delta > 0$,
\[ \P (Y \preceq y + \delta 1_k) - \P (Y \preceq y) \ \leq \ \frac{\delta}{\underline{\sigma}} \big( \sqrt{2 \log (p)} + 2 \big) . \]
\end{lemma}

The following lemma is Bentkus' multivariate Berry-Esseen theorem.

\begin{lemma} [Theorem 1.1 in \citep{Bentkus03}] \label{background5}
Let $V_1, \dots, V_n$ be i.i.d. random vectors $\R^d$, with zero mean and identity covariance matrix. Furthermore, let $\zeta$ be a standard Gaussian vector in $\R^d$, and let $\A$ denote the collection of all Borel convex subsets of $\R^d$. Then, there is an absolute constant $c > 0$ such that
\[ \sup_{A \in \A} \Big| \P \Big( \frac{1}{\sqrt{n}} \sum_{i = 1}^n V_i \in A \Big) - \P (\zeta \in A) \Big| 
\ \leq \
\frac{c \cdot d^{1/4} \cdot \E [\| V_1 \|_2^3]}{n^{1/2} } .\]
\end{lemma}

For the statement of Lemma~\ref{background7} below, we need to introduce a bit of notation. For any $r>0$ and set $A\subset\R^p$,  define the outer $r$-neighborhood as 
 \smash{$A^r=\big\{x\in\R^p\, |\, d(x,A)\leq r\},$}
 where $d(x,A)=\inf\{\|x-y\| \, | \, y\in A\}$, and $\|\cdot\|$ is any norm on $\R^p$. The corresponding inner $r$-neighborhood may be defined as
 $A^{-r}=\big\{x\in A\, |\, B(x,r)\subset A\},$
 where $B(x,r)=\{y\in\R^p| \|x-y\|\leq r\}$.

\begin{lemma}[Lemma 7.3 in~\citep{Lopes2020}]\label{background7}
Let $\|\cdot\|$ be any norm on $\R^p$, and let $\zeta,\xi\in\R^p$ be any two random vectors. Then, the following inequality holds for any Borel set $A\subset\R^p$, and any $r>0$,
\[
|\P(\zeta \in A)-\P(\xi\in A )| \ \leq \  \P\big(\xi\in (A^r\setminus A^{-r})\big) \ + \   \P\big(\|\zeta-\xi\|\geq r\big).
\]
\end{lemma}

The following lemma is a consequence of Pinsker's inequality and the proof of Lemma A.7 in the paper \cite{SpokoinyZ15}.
\begin{lemma}[\cite{SpokoinyZ15}]\label{lem:Zhilova}
Let $\zeta$ and $\tilde \zeta$ be centered Gaussian vectors in $\R^k$ with respective covariance matrices $C$ and $\tilde C$. Also, suppose that $C$ is invertible, and let $B=C^{-1/2}\tilde C C^{-1/2}-I_k$. Then, there is an absolute constant $c>0$ such that
\begin{equation*}
    \sup_{t\in\R^k}\Big|\P(\zeta\preceq t) - \P(\tilde\zeta\preceq t)\Big| \ \leq \ c\sqrt{k}\|B\|_{\op}.
\end{equation*}
\end{lemma}

The last background lemma follows from the proof of~\cite[][Lemma D.3]{LopesWL20}, Lemma~\ref{background7}, and Nazarov's inequality~(Lemma~\ref{background4}).
\begin{lemma}\label{background6}
Let $U$, $V$, and $R$ be random vectors in $\R^k$ that satisfy $U = V + R$. Also, let $W= (W_1, \dots, W_k)$ be a centered Gaussian random vector in $\R^k$ and suppose that the parameter $\underline{\sigma}^2=\min_{1\leq j\leq k}\E [W_j^2]$ is positive. Then there is an absolute constant $c > 0$, such that the following bound holds for any $r > 0$, 
\[ \sup_{t \in \R^k} \big| \P (U \mylessthan) - \P(W \mylessthan) \big|
\ \leq \ 
3 \cdot \sup_{t \in \R^k} \big| \P (V \mylessthan) - \P (W \mylessthan) \big| + \frac{c r \sqrt{\log(k)} }{\underline{\sigma}} + \P ( \| R \|_\infty \geq r ).
\]
\end{lemma}

\newpage

\section{Additional Numerical Results} \label{appendix_last}

\paragraph{Nominal value of 95\% in model (ii).}
The following four figures are presented in the same manner as in the main text for a 95\% nominal value, except that they are based on simulation model (ii).

\begin{figure}[H]	
\vspace{0.7cm}
	\quad\quad\quad 
	\begin{overpic}[width=0.29\textwidth]{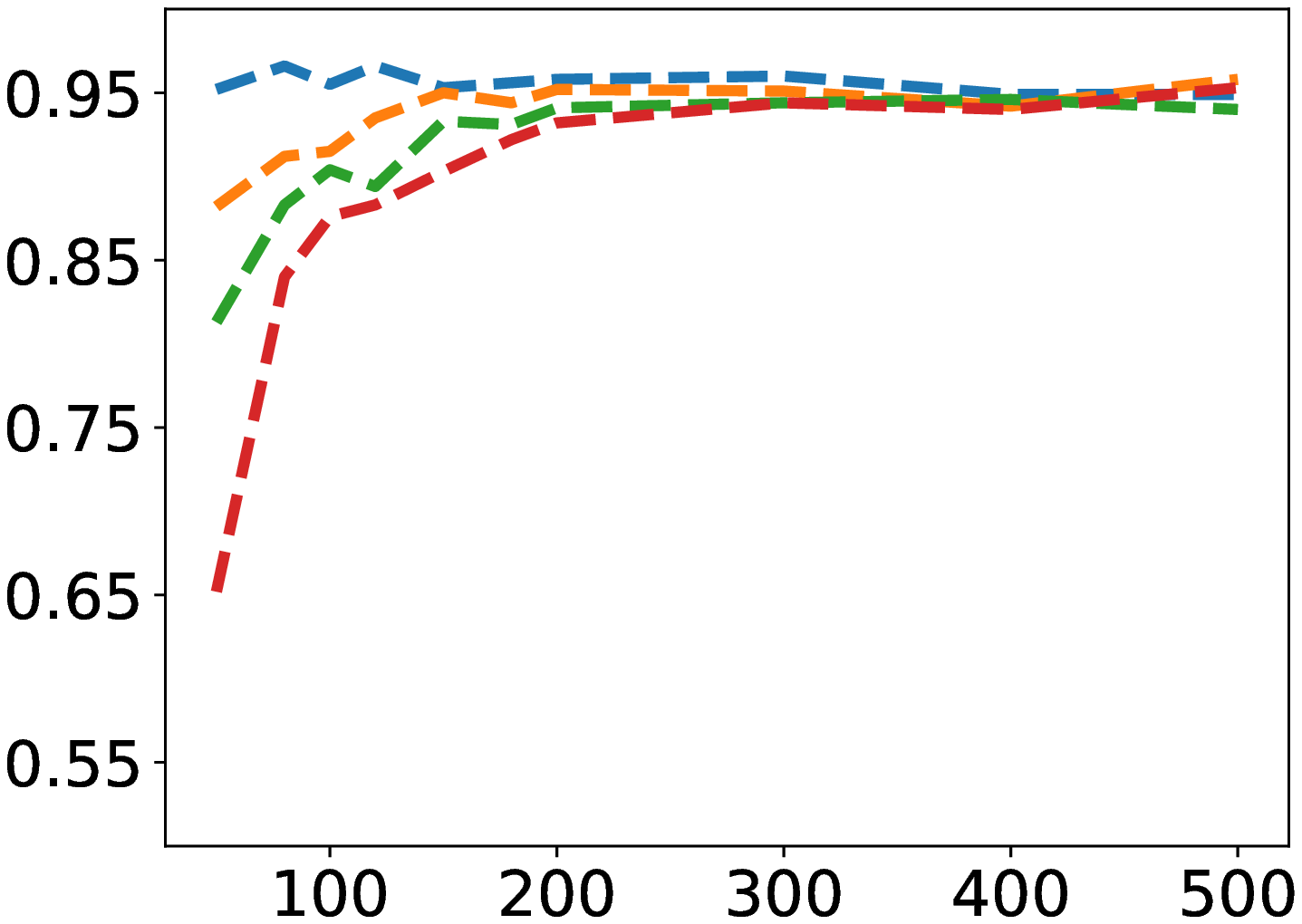} 
		\put(25,80){ \ul{\ \ \  \ $\gamma=0.7$ \ \ \ \    }}
		\put(-20,-5){\rotatebox{90}{ {\small \ \ \ log transformation  \ \ }}}
	\end{overpic}
	~
	\DeclareGraphicsExtensions{.png}
	\begin{overpic}[width=0.29\textwidth]{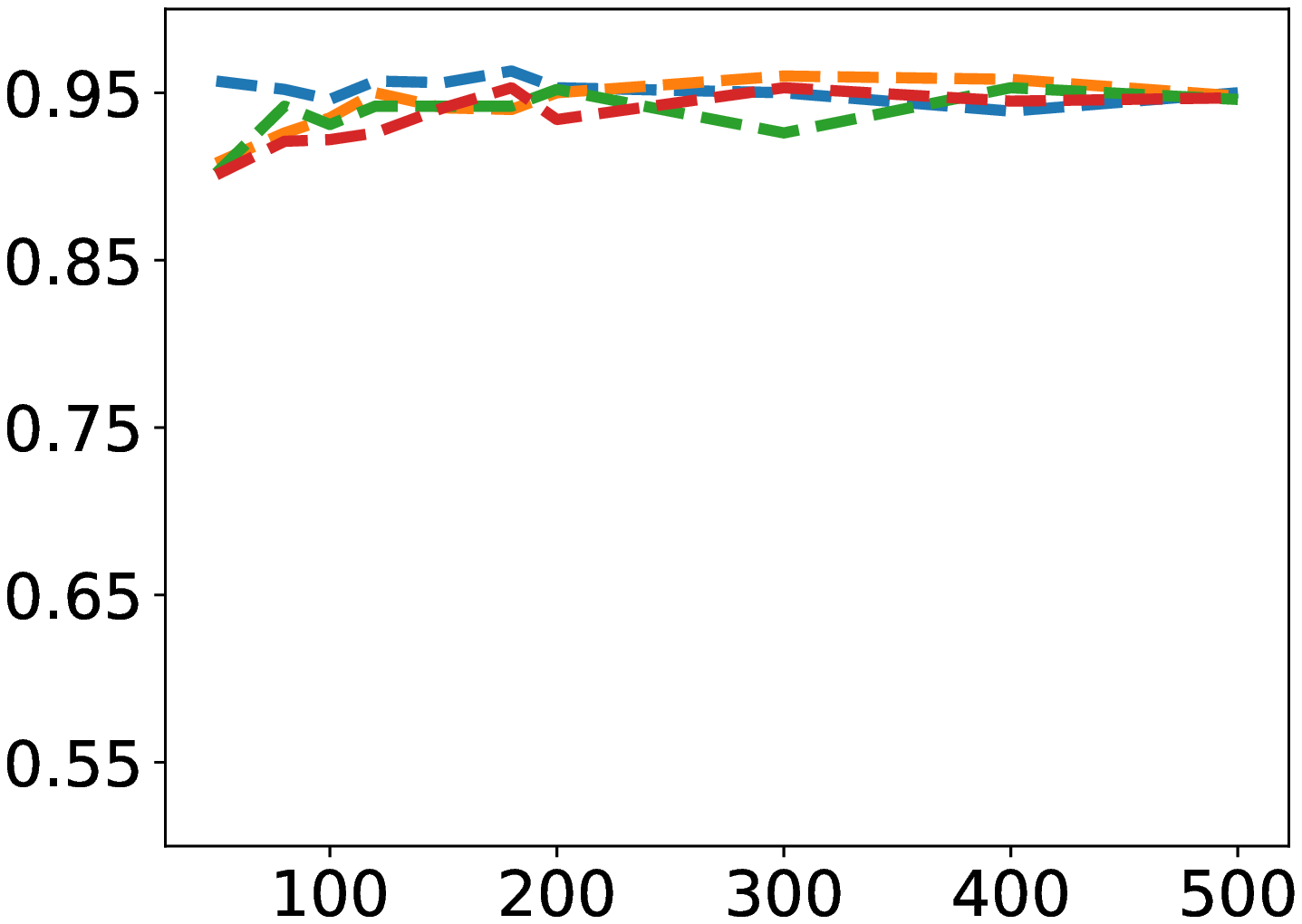} 
		\put(25,80){ \ul{\ \ \  \ $\gamma=1.0$ \ \ \ \    }}
	\end{overpic}
	~	
	\begin{overpic}[width=0.29\textwidth]{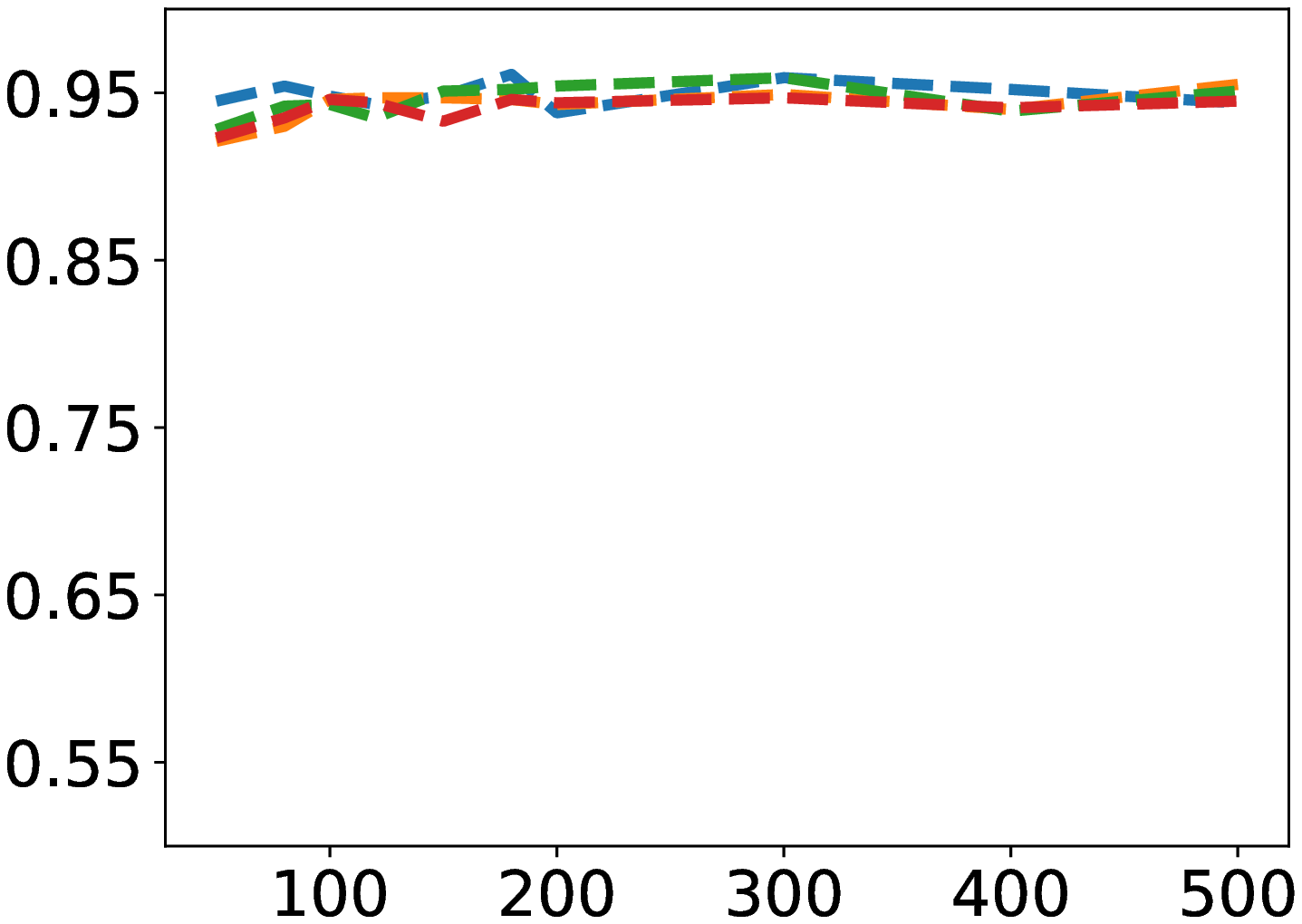} 
			\put(25,80){ \ul{\ \ \  \ $\gamma=1.3$ \ \ \ \    }}
	\end{overpic}	
%
%
\end{figure}

\vspace{-0.5cm}

\begin{figure}[H]	
	\quad\quad\quad 
	\begin{overpic}[width=0.29\textwidth]{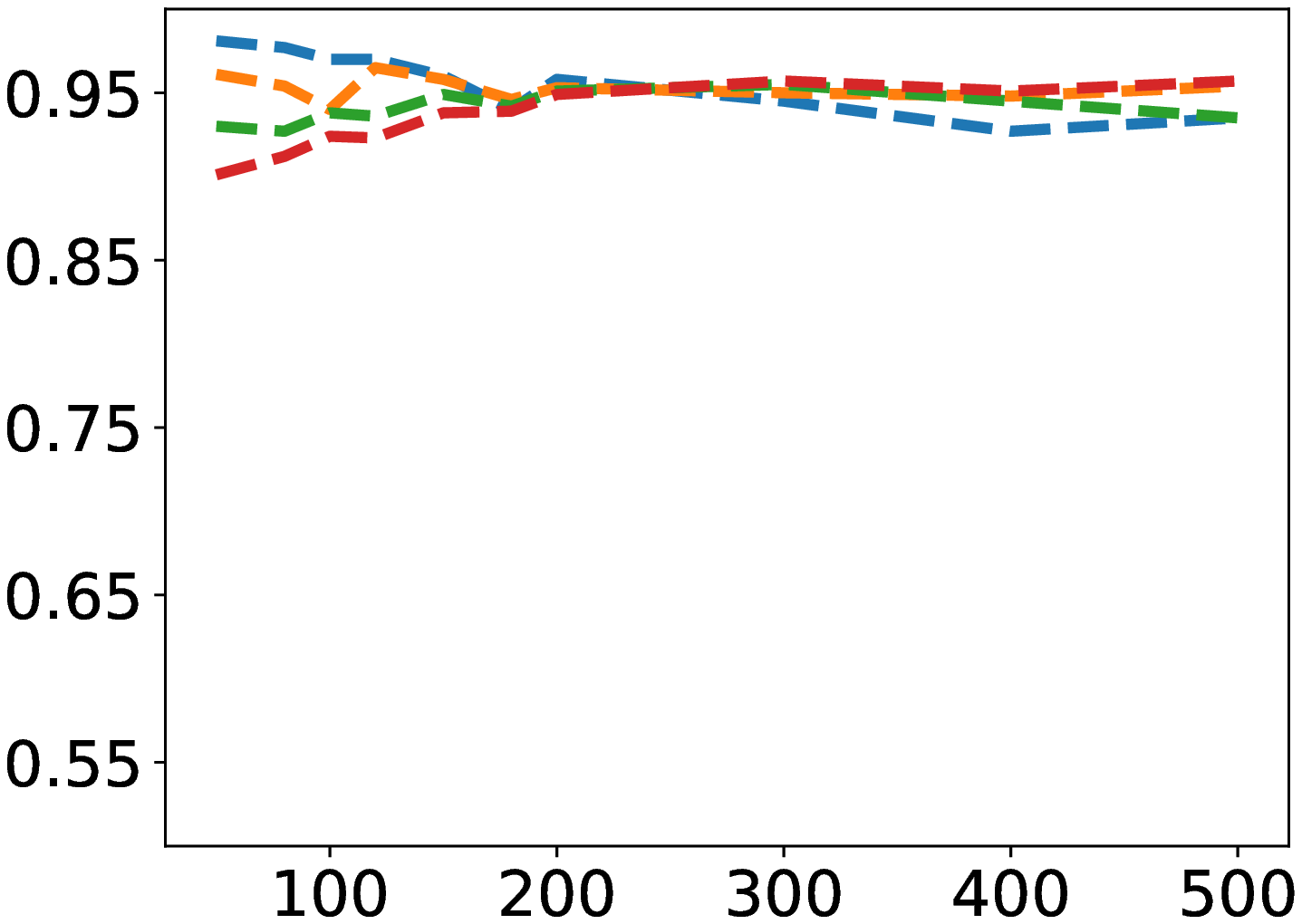} 
	\put(-20,-1){\rotatebox{90}{ {\small \ \ \ standardization \  \ \ }}}
	\end{overpic}
	~
	\DeclareGraphicsExtensions{.png}
	\begin{overpic}[width=0.29\textwidth]{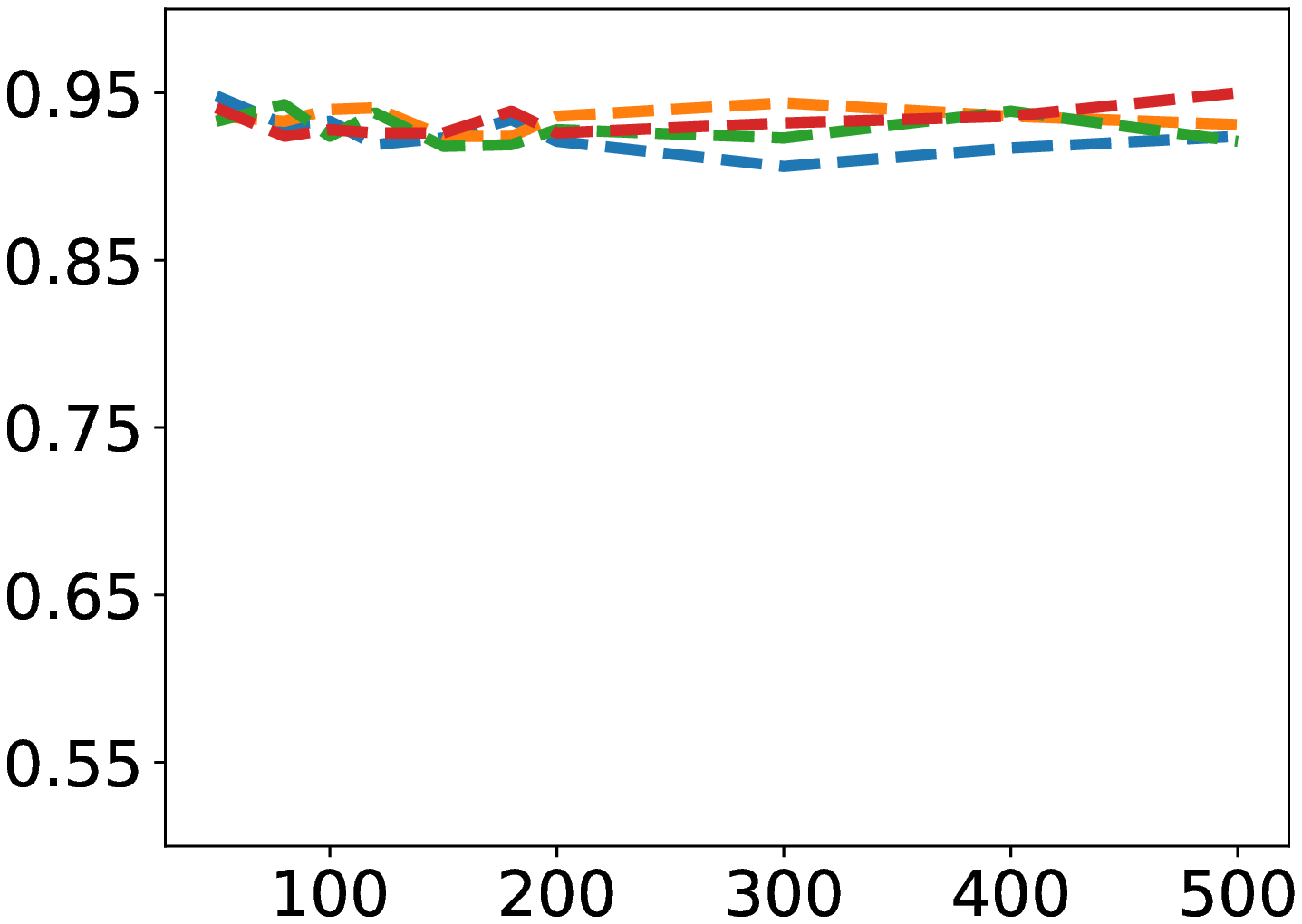} 
	\end{overpic}
	~	
	\begin{overpic}[width=0.29\textwidth]{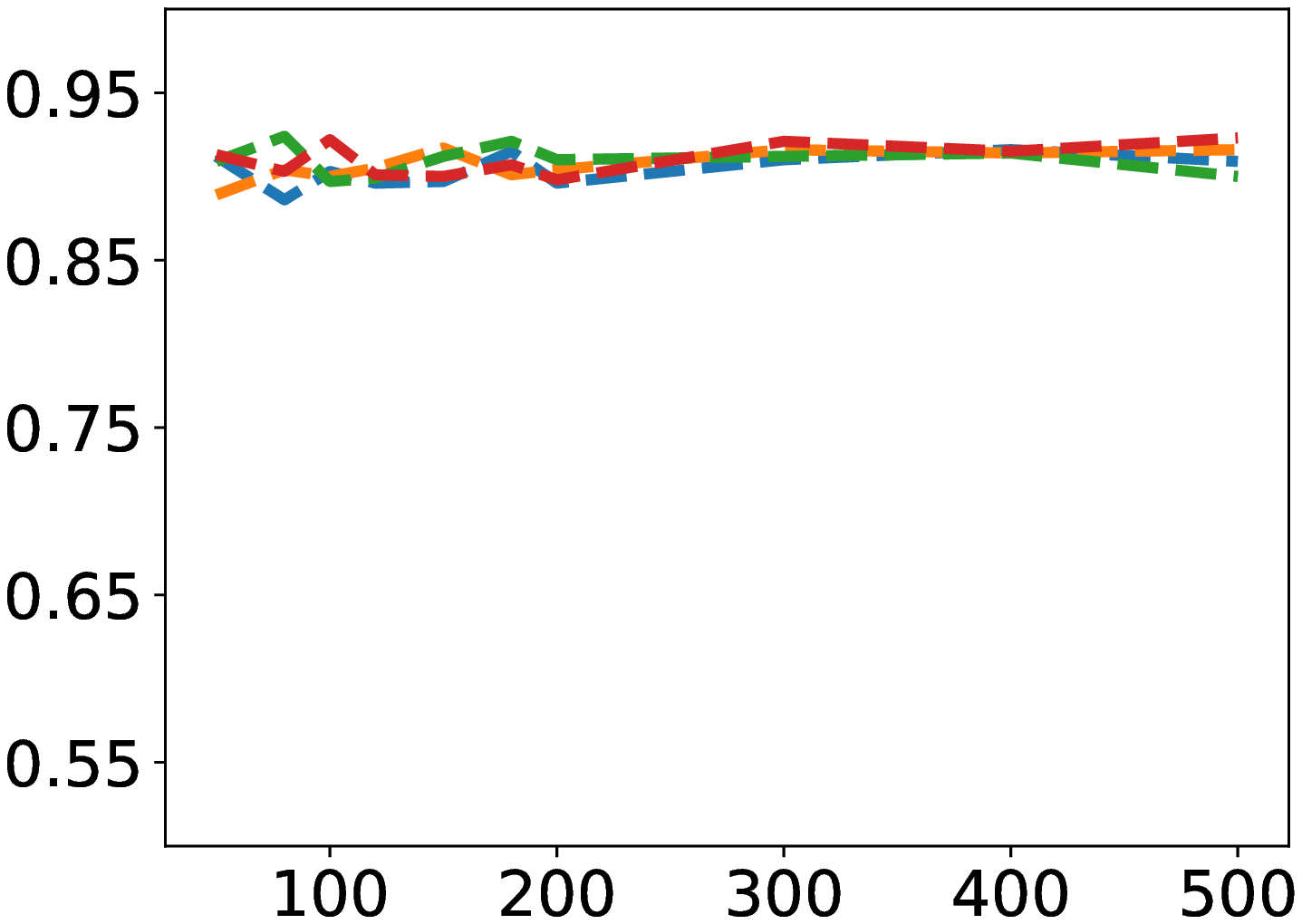} 
	\end{overpic}	
%
%
\end{figure}

\vspace{-0.5cm}

\begin{figure}[H]	
	\quad\quad\quad 
	\begin{overpic}[width=0.29\textwidth]{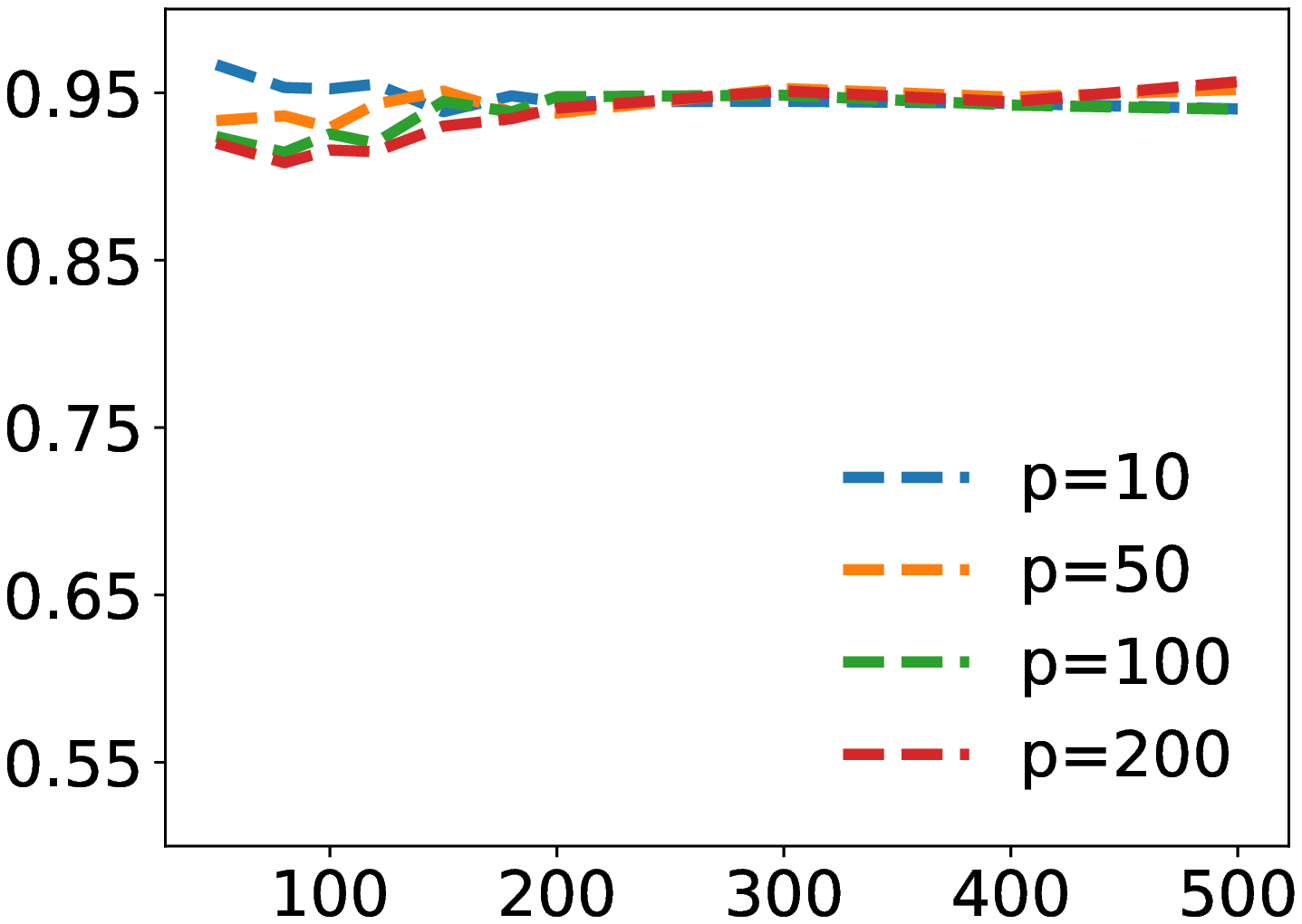} 
		\put(-21,1){\rotatebox{90}{\ $\sqrt{ \ \ }$}}
	    \put(-20,-3){\rotatebox{90}{  { \ \ \ \ \ \ \small transformation \ \ } }}
	\end{overpic}
	~
	\DeclareGraphicsExtensions{.png}
	\begin{overpic}[width=0.29\textwidth]{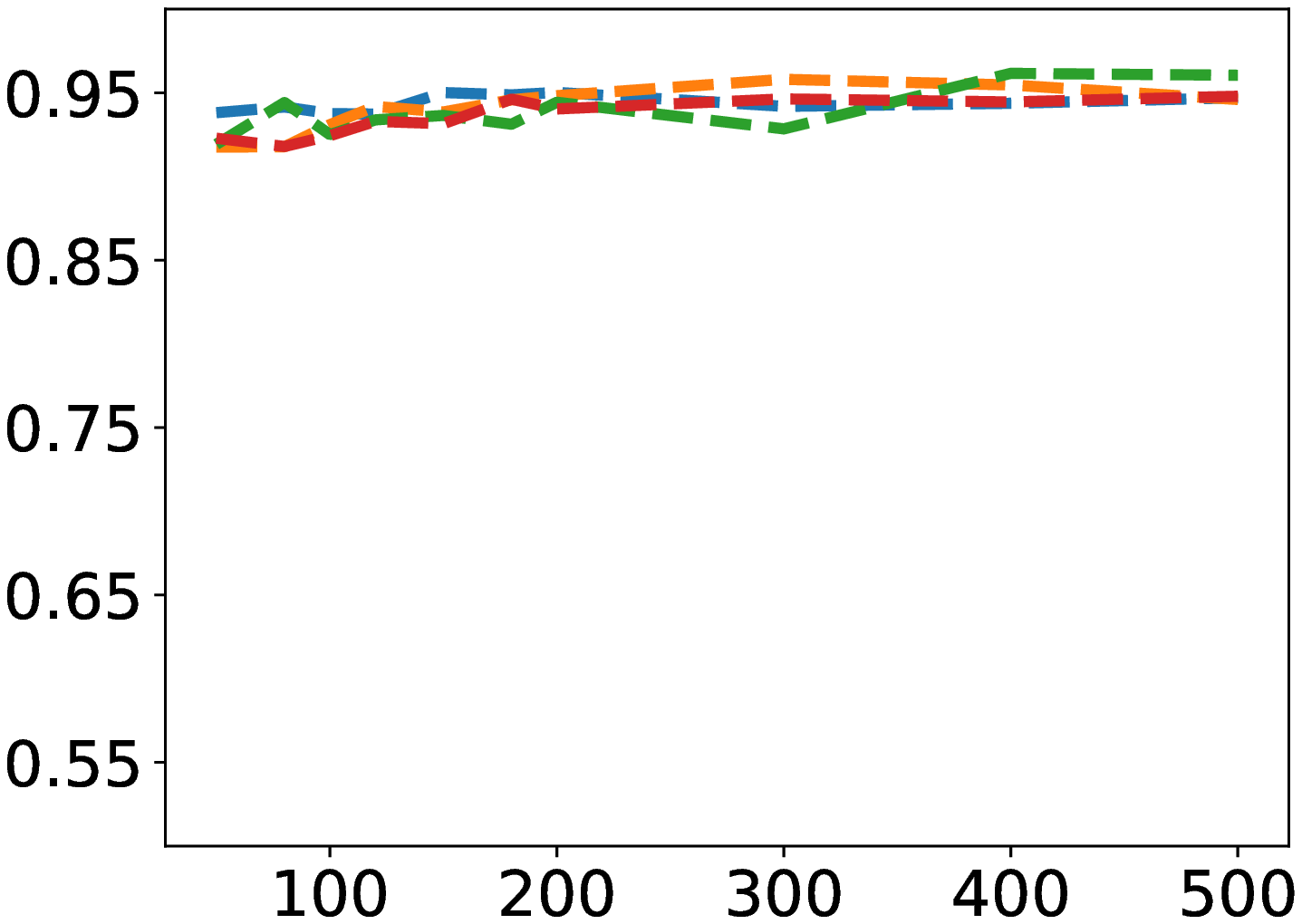} 
	\end{overpic}
	~	
	\begin{overpic}[width=0.29\textwidth]{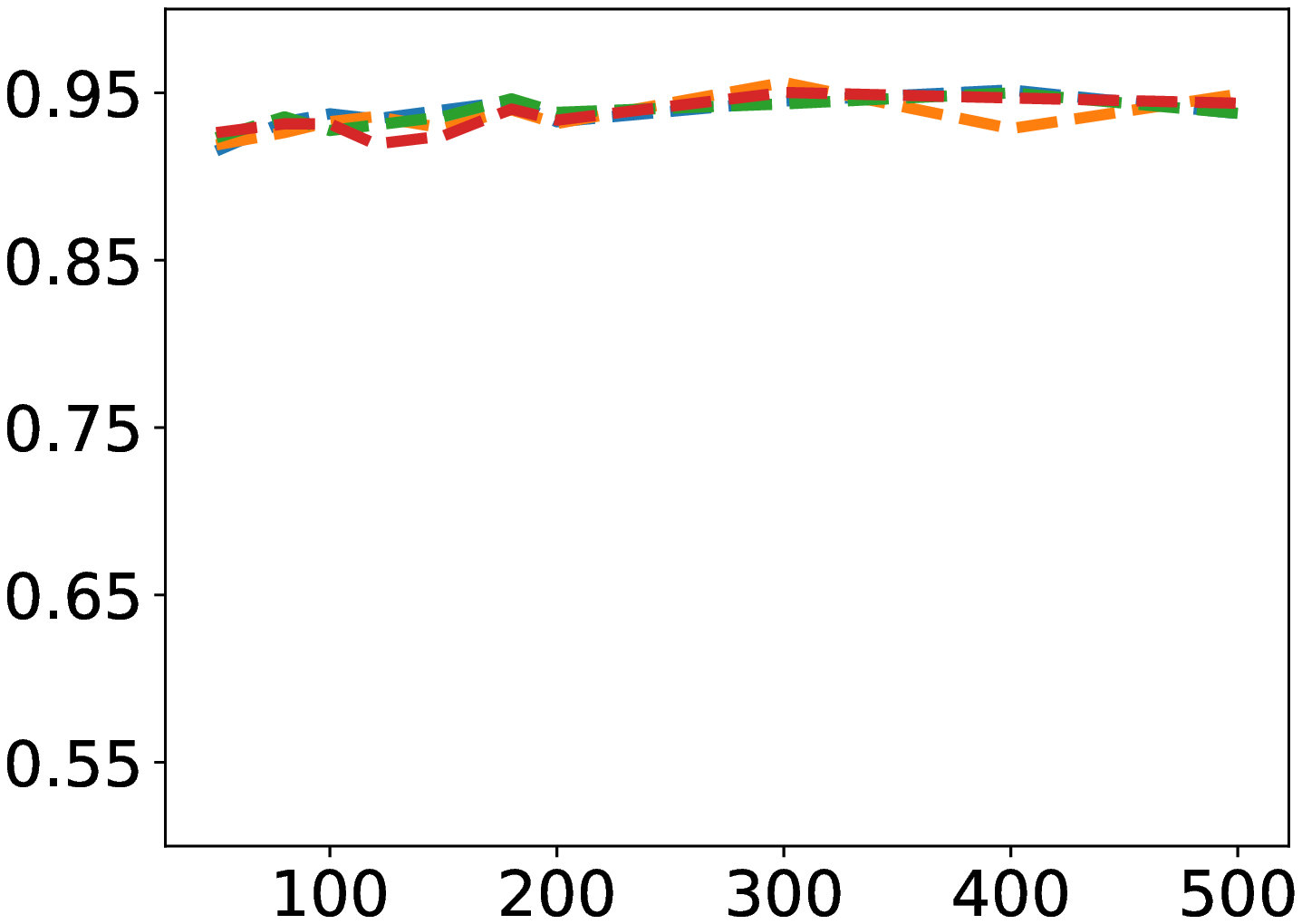} 
	\end{overpic}	
    \vspace{+.2cm}	
	\caption{(Simultaneous coverage probability versus $n$ in simulation model (ii) with a polynomial decay profile). The plotting scheme is the same as described in the caption of Figure~\ref{fig1} in the main text.} 
	\label{fig2}
\end{figure}

\newpage
\begin{figure}[H]	
\vspace{0.5cm}
	\quad\quad\quad 
	\begin{overpic}[width=0.29\textwidth]{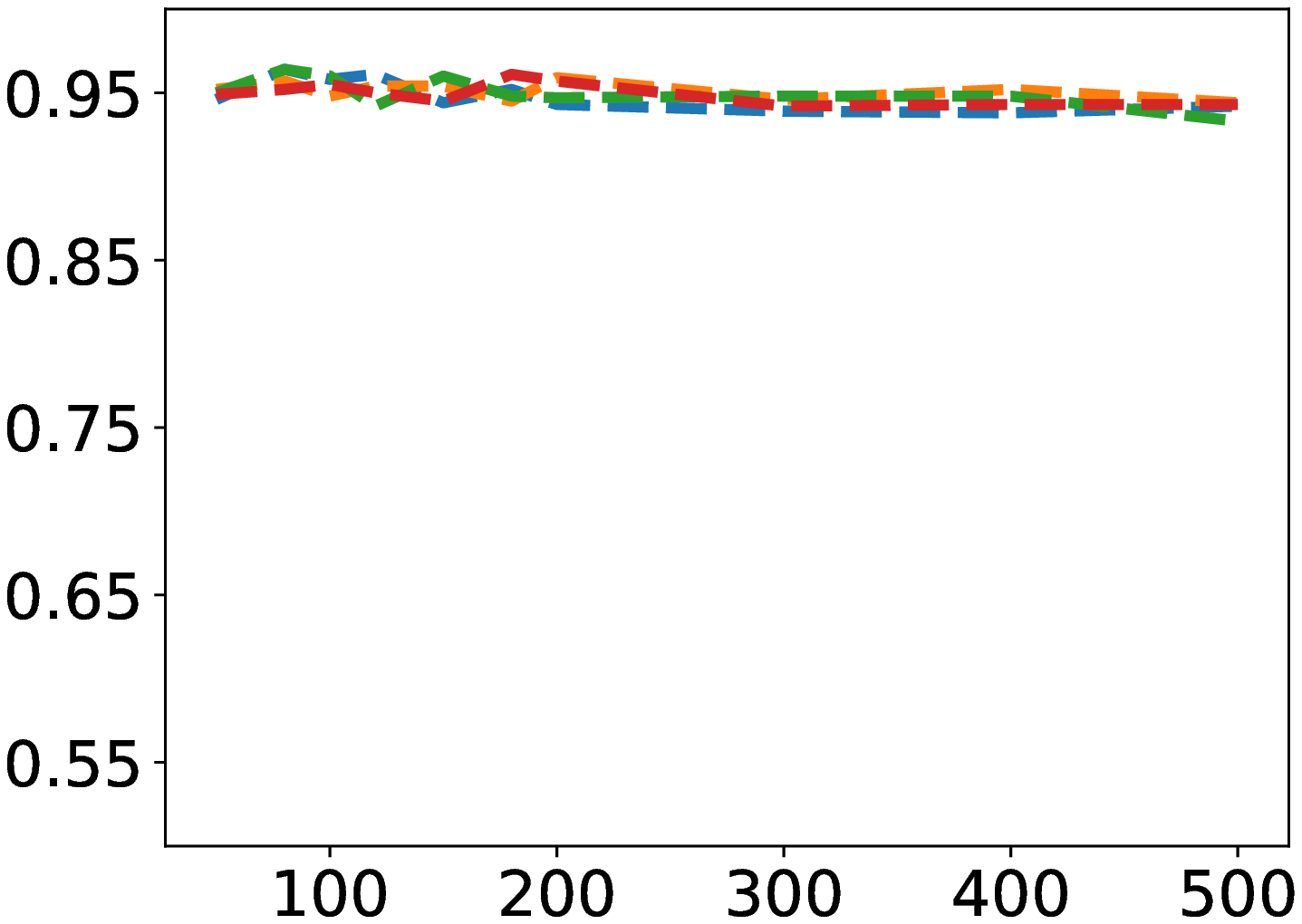} 
\put(25,80){ \ul{\ \ \  \ $\delta=0.7$ \ \ \ \    }}
		\put(-20,-5){\rotatebox{90}{ {\small \ \ \ log transformation  \ \ }}}
\end{overpic}
	~
	\DeclareGraphicsExtensions{.png}
	\begin{overpic}[width=0.29\textwidth]{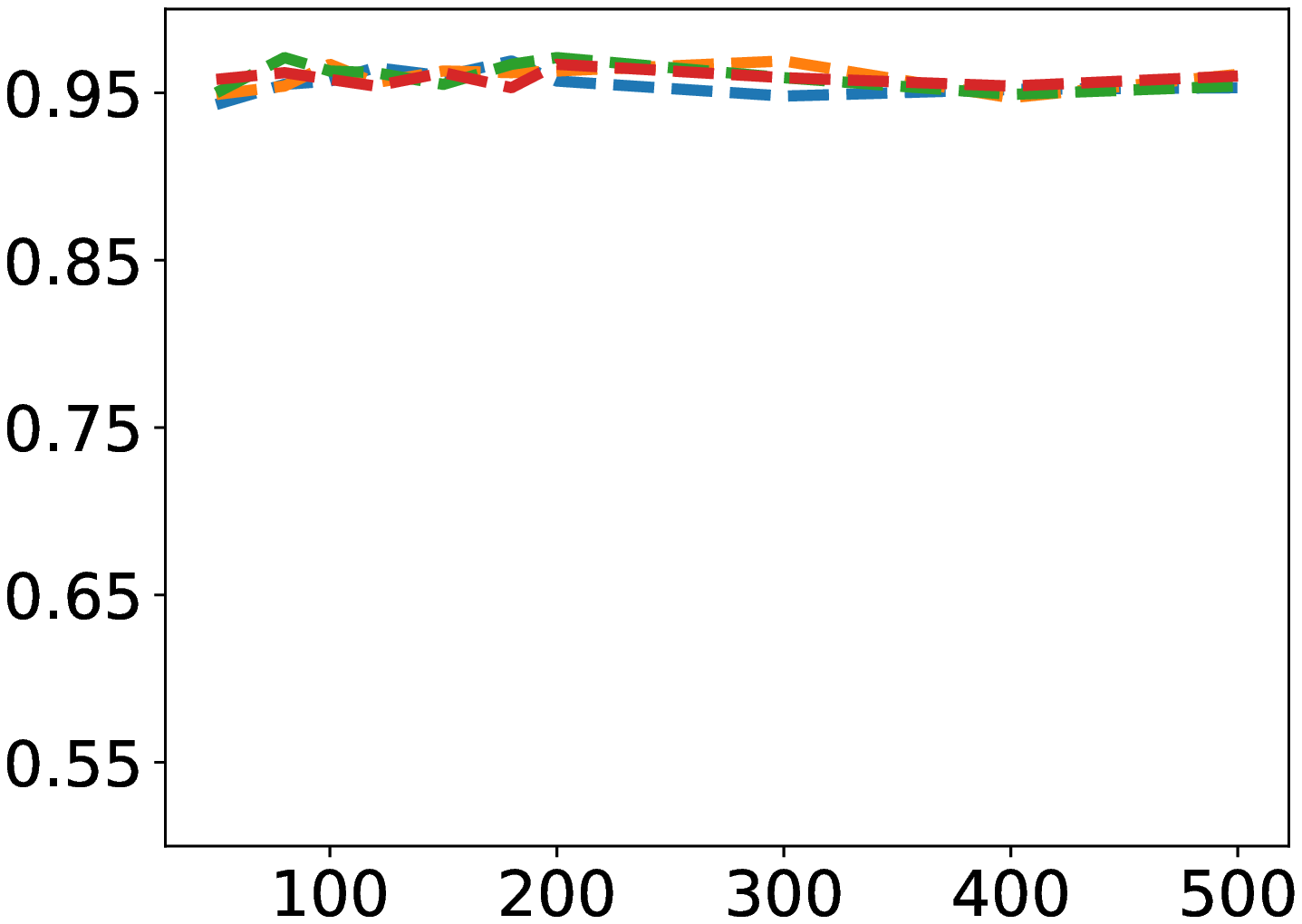} 
		\put(25,80){ \ul{\ \ \  \ $\delta=0.8$ \ \ \ \    }}
	\end{overpic}
	~	
	\begin{overpic}[width=0.29\textwidth]{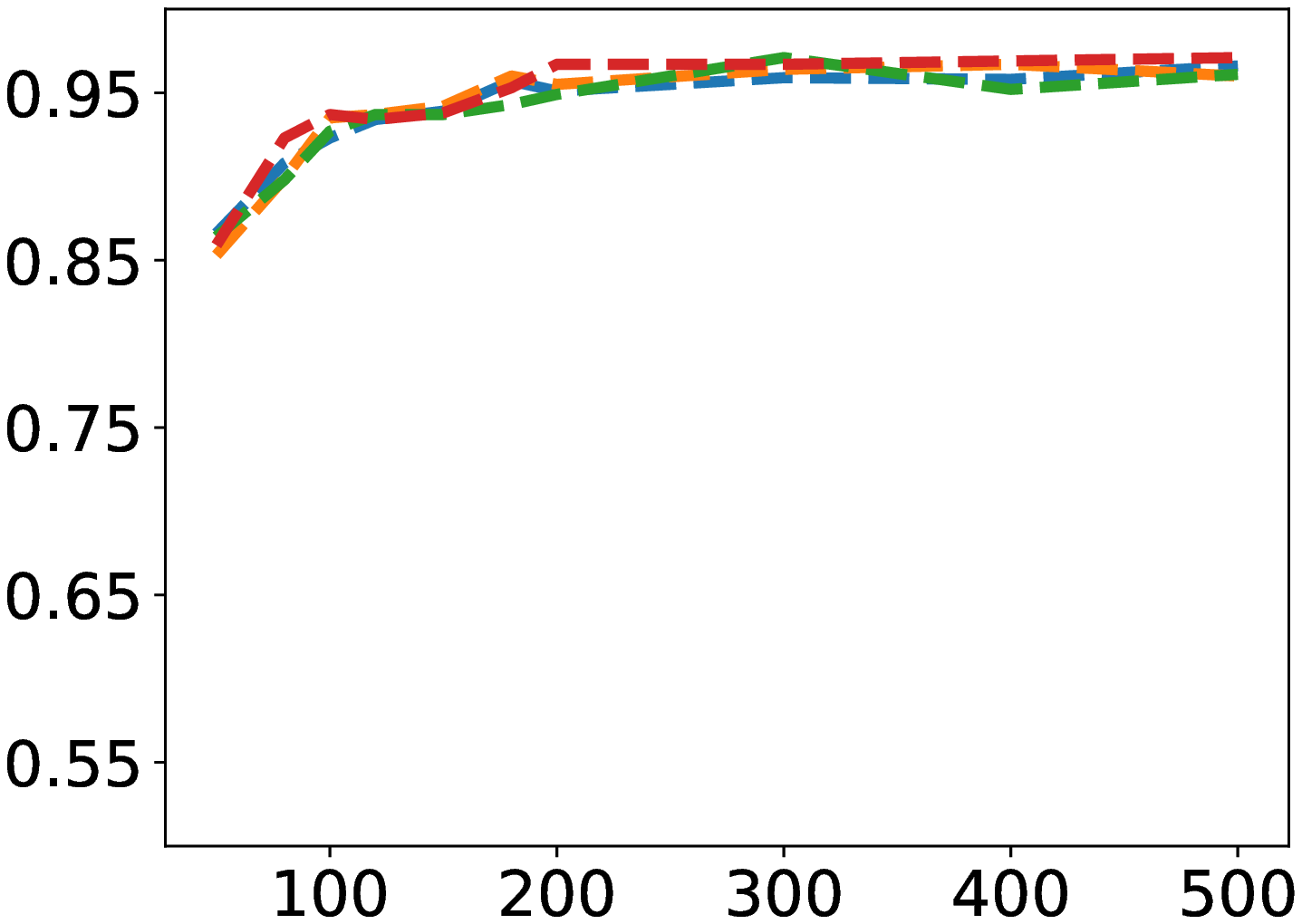} 

		\put(25,80){ \ul{\ \ \  \ $\delta=0.9$ \ \ \ \    }}
 		 				
	\end{overpic}	
\end{figure}

\vspace{-0.5cm}

\begin{figure}[H]	
	\quad\quad\quad 
	\begin{overpic}[width=0.29\textwidth]{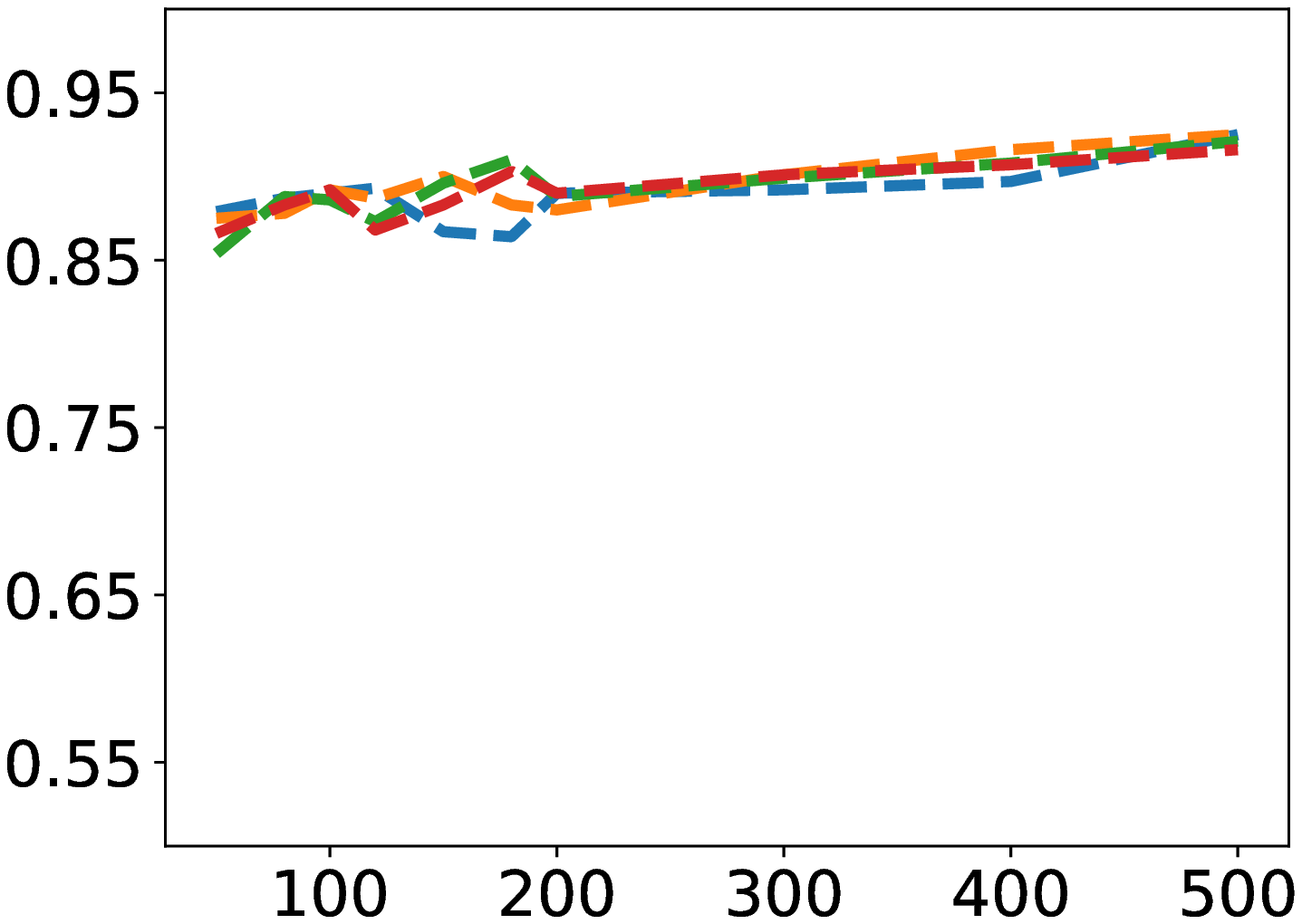} 
	\put(-20,-1){\rotatebox{90}{ {\small \ \ \ standardization \  \ \ }}}
	\end{overpic}
	~
	\DeclareGraphicsExtensions{.png}
	\begin{overpic}[width=0.29\textwidth]{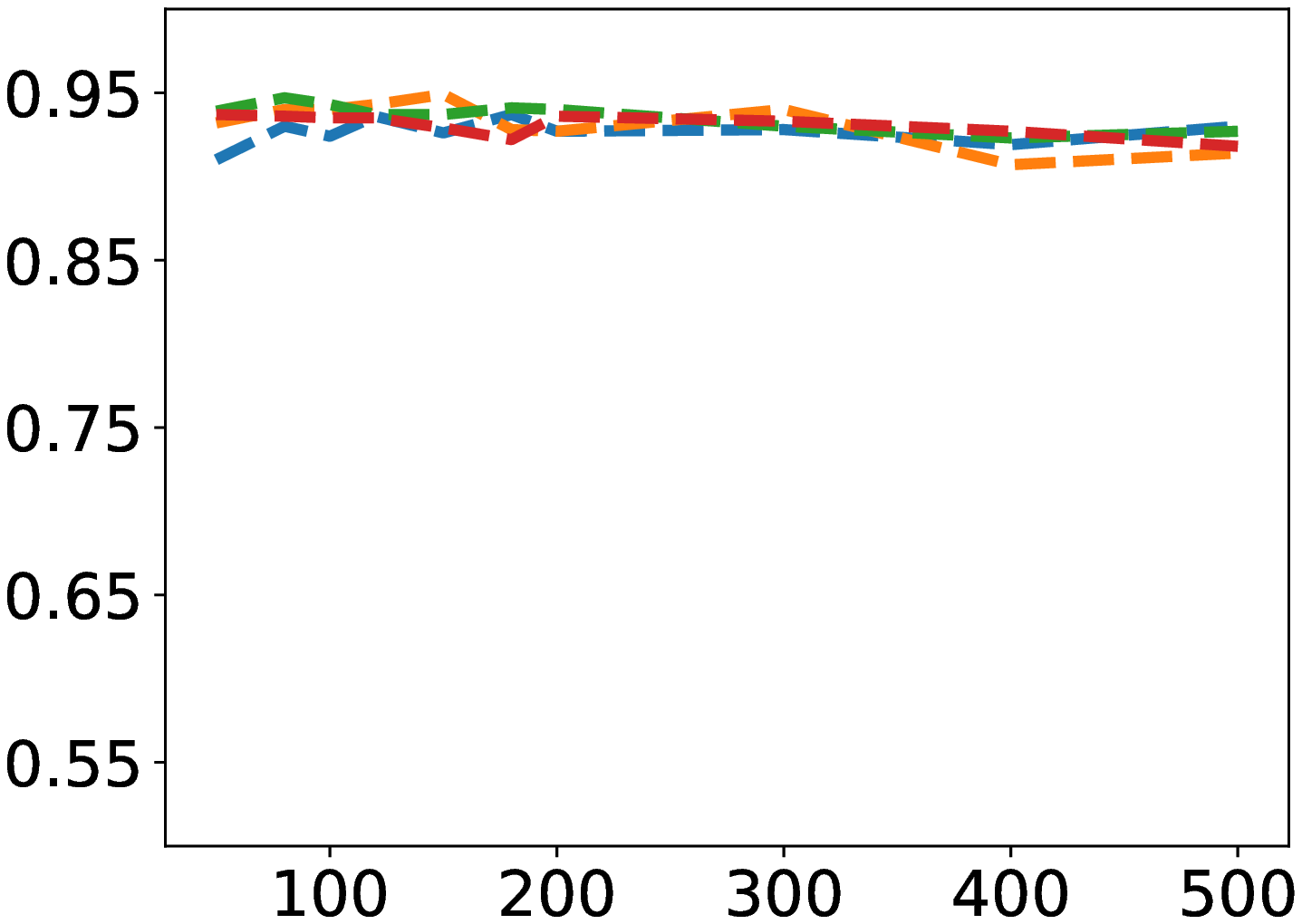} 
	\end{overpic}
	~	
	\begin{overpic}[width=0.29\textwidth]{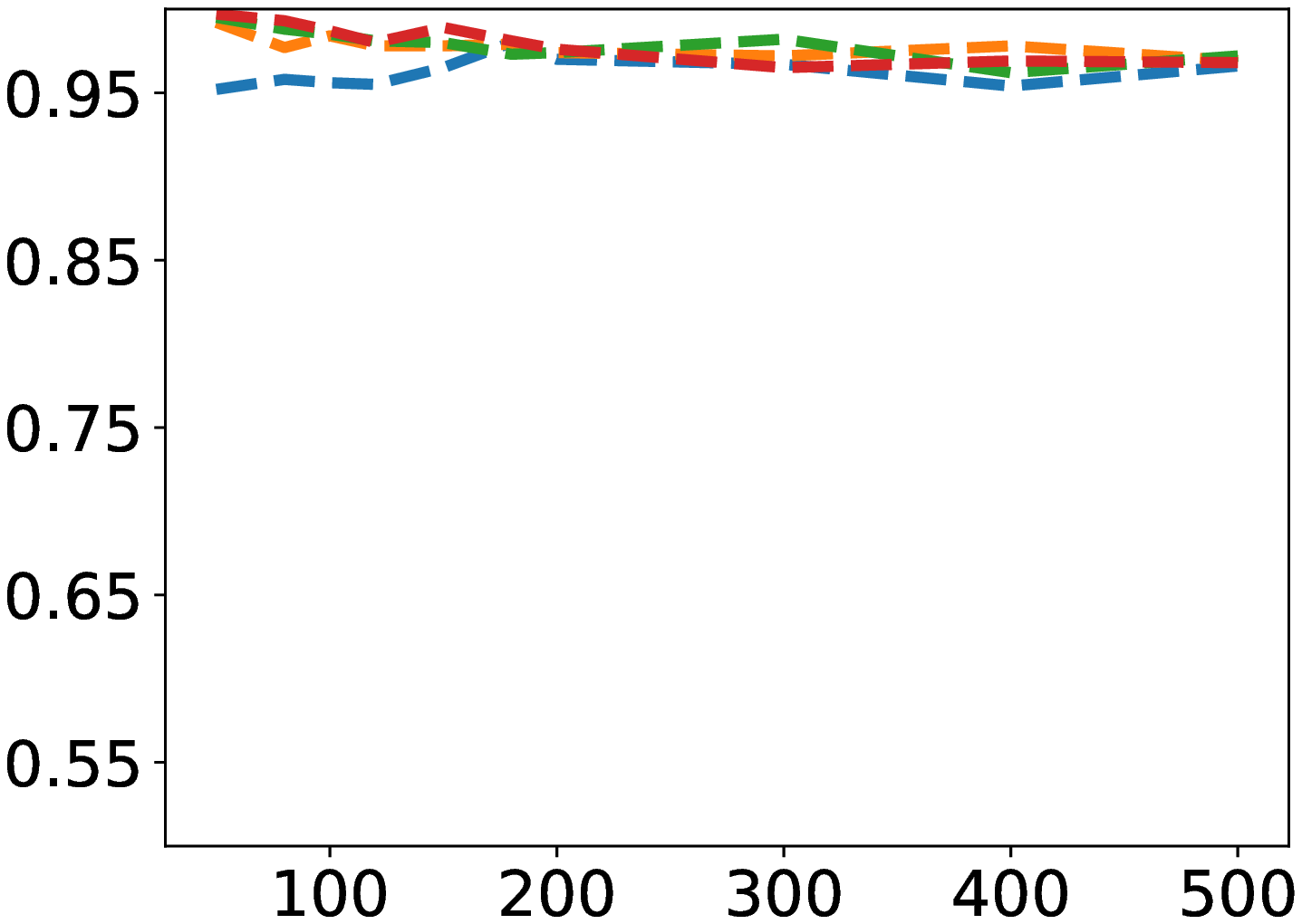} 
	 				
	\end{overpic}	
\end{figure}

\vspace{-0.5cm}

\begin{figure}[H]	
	\quad\quad\quad 
	\begin{overpic}[width=0.29\textwidth]{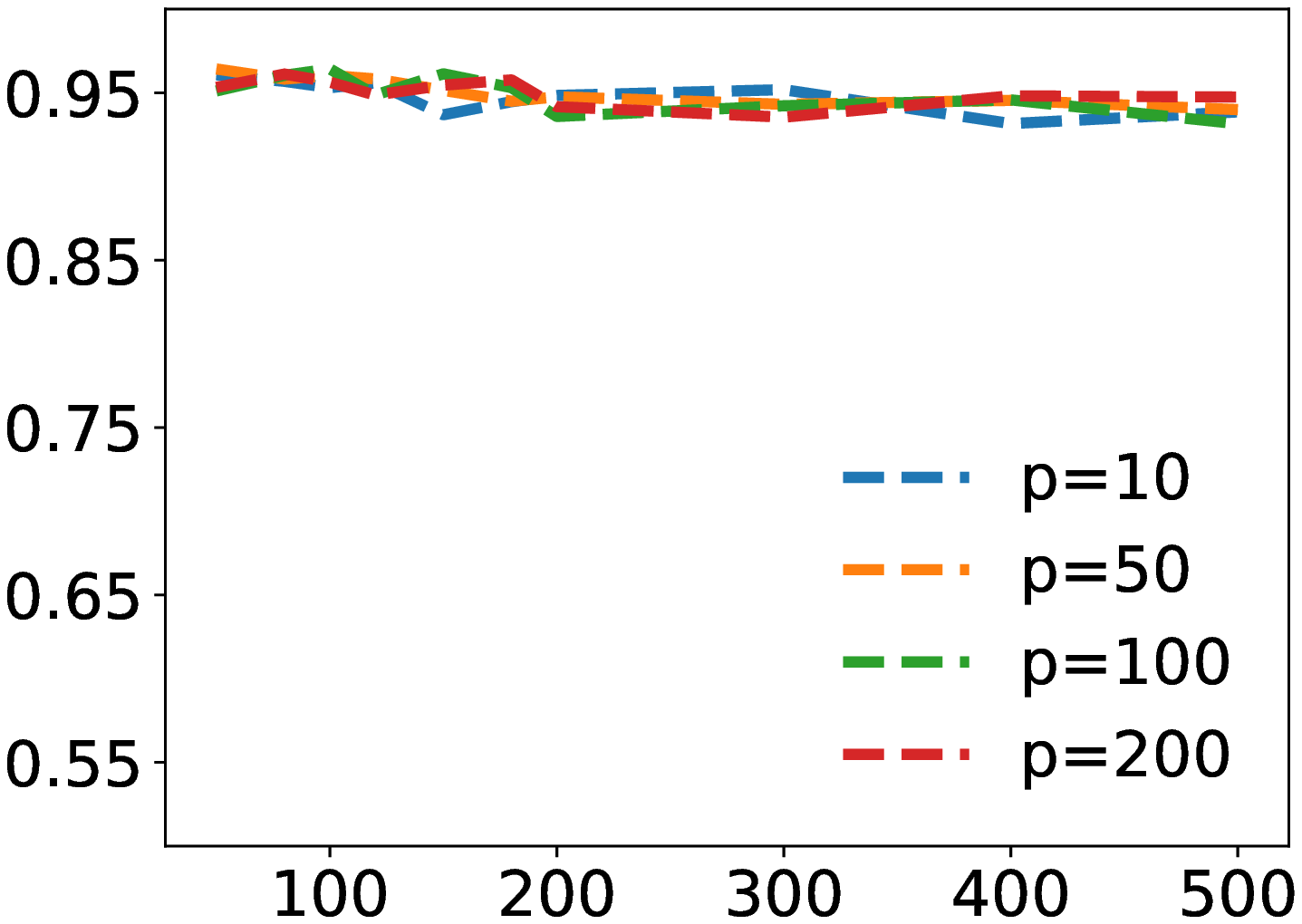} 
		\put(-21,1){\rotatebox{90}{\ $\sqrt{ \ \ }$}}
	    \put(-20,-3){\rotatebox{90}{  { \ \ \ \ \ \ \small transformation \ \ } }}

	\end{overpic}
	~
	\DeclareGraphicsExtensions{.png}
	\begin{overpic}[width=0.29\textwidth]{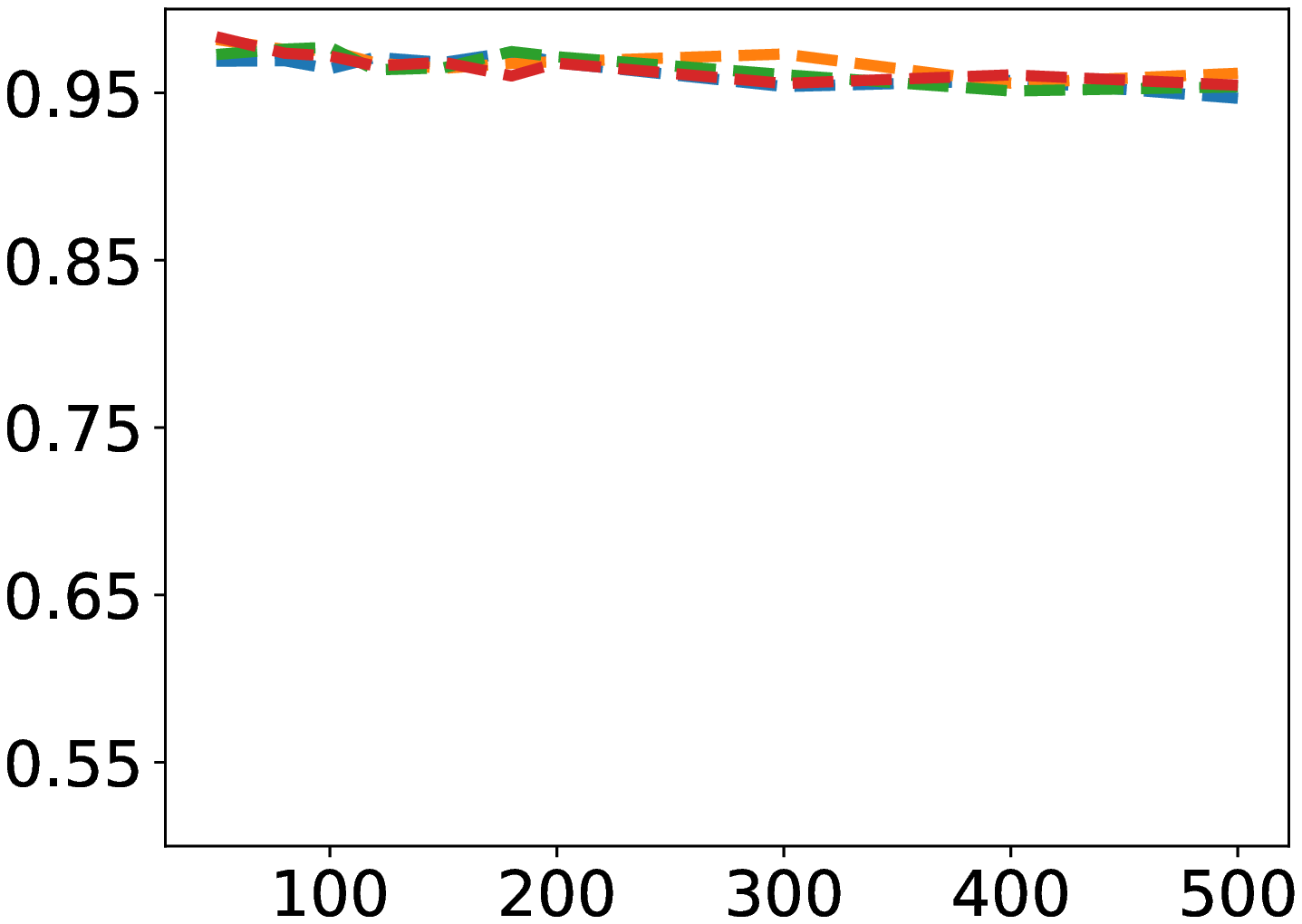} 
	\end{overpic}
	~	
	\begin{overpic}[width=0.29\textwidth]{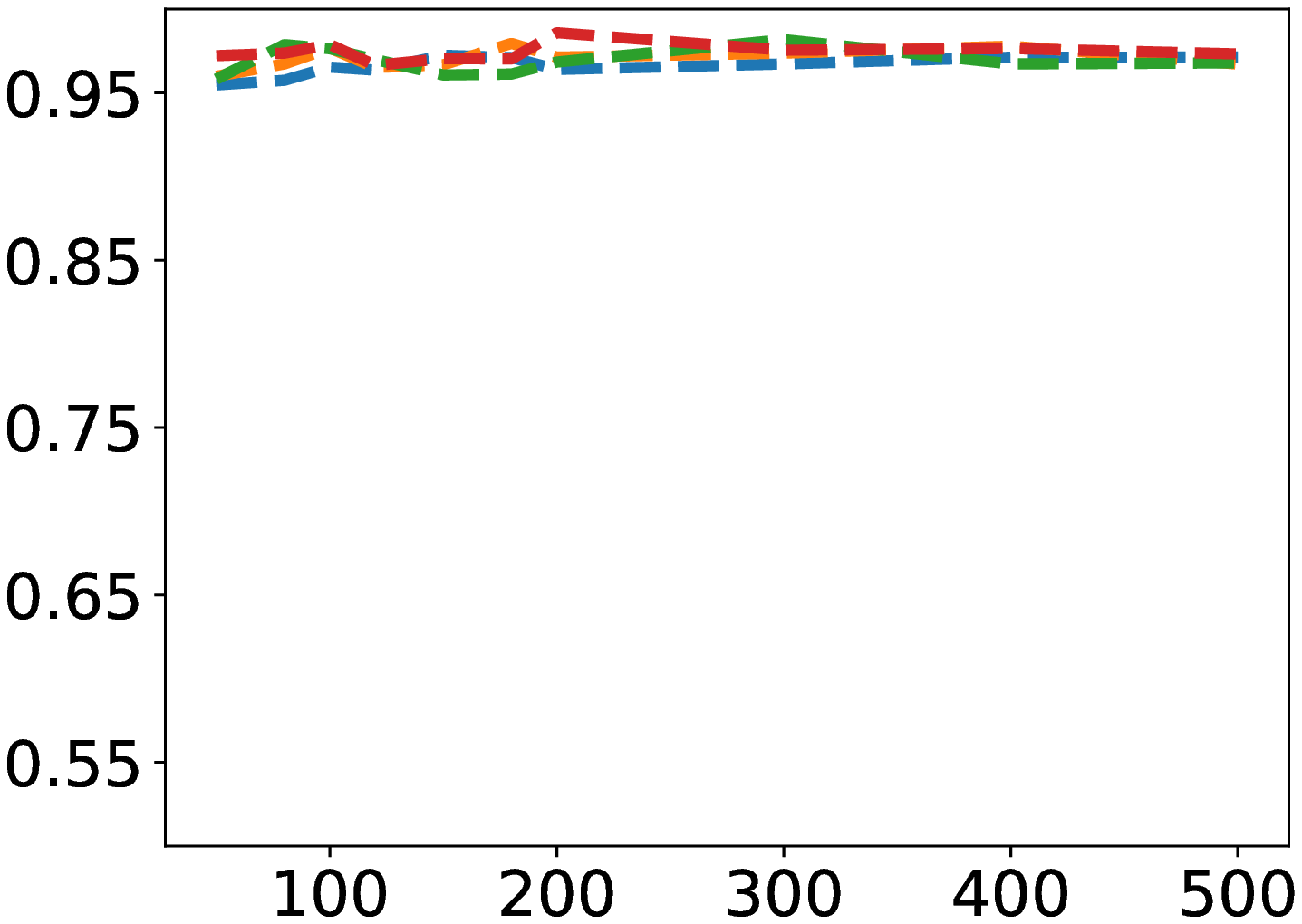} 
		 				
	\end{overpic}	
	\vspace{+0.2cm}	
	\caption{(Simultaneous coverage probability versus $n$ in simulation model (ii) with an exponential decay profile). The plotting scheme is the same as described in the caption of Figure~\ref{fig3} in the main text.}
	\label{fig4}
\end{figure}

\newpage

\begin{figure}[H]	
\vspace{0.5cm}

	\quad\quad\quad 
	\begin{overpic}[width=0.29\textwidth]{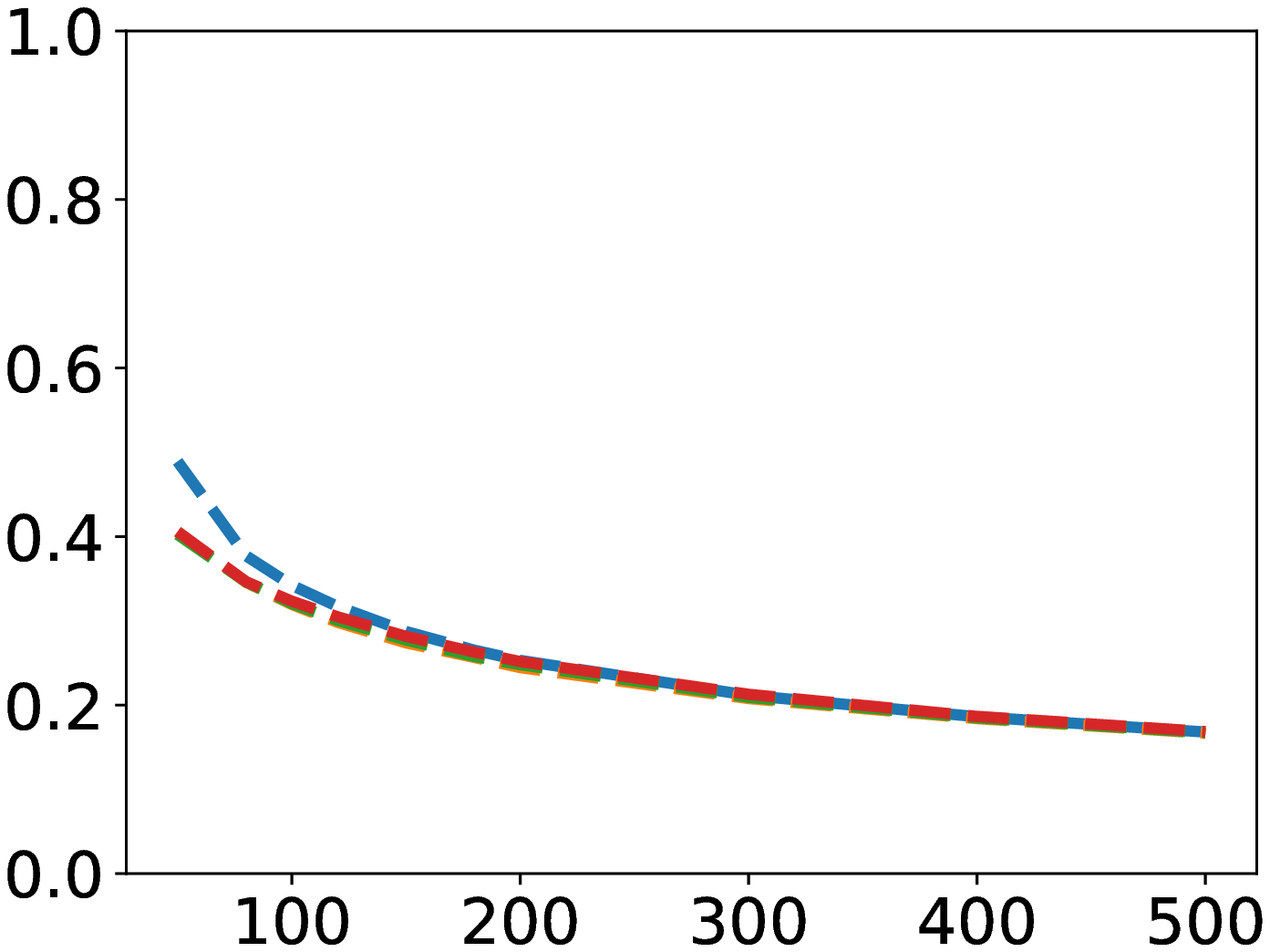} 
		\put(25,80){ \ul{ \ \ \  \ $\gamma=0.7$ \ \ \ \    }}	
		\put(-20,-5){\rotatebox{90}{ {\small \ \ \ log transformation  \ \ }}}
	\end{overpic}
	~
	\DeclareGraphicsExtensions{.png}
	\begin{overpic}[width=0.29\textwidth]{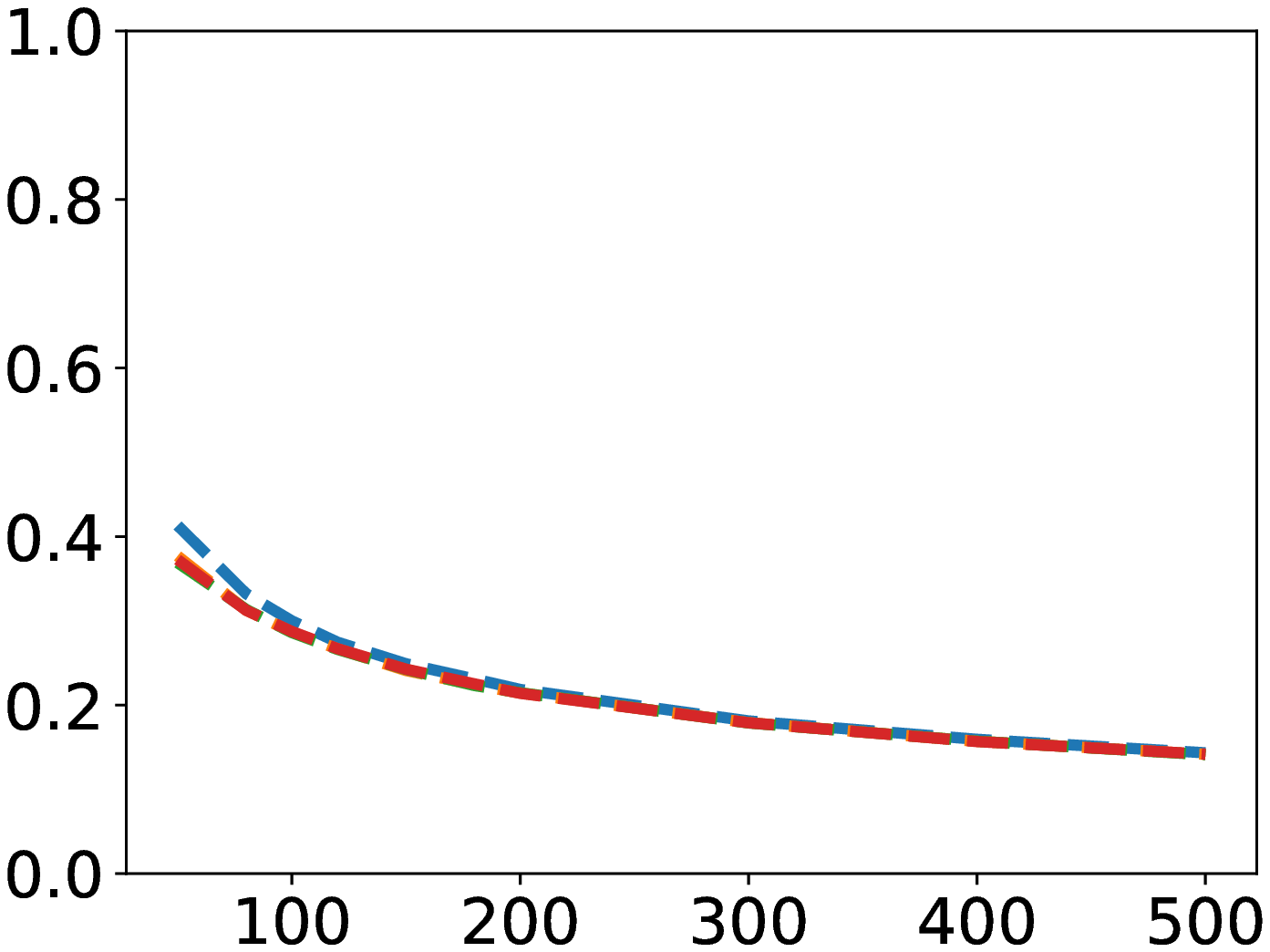} 
		\put(25,80){\ul{\ \ \  \ $\gamma=1.0$ \ \ \ \    }}
	\end{overpic}
	~	
	\begin{overpic}[width=0.29\textwidth]{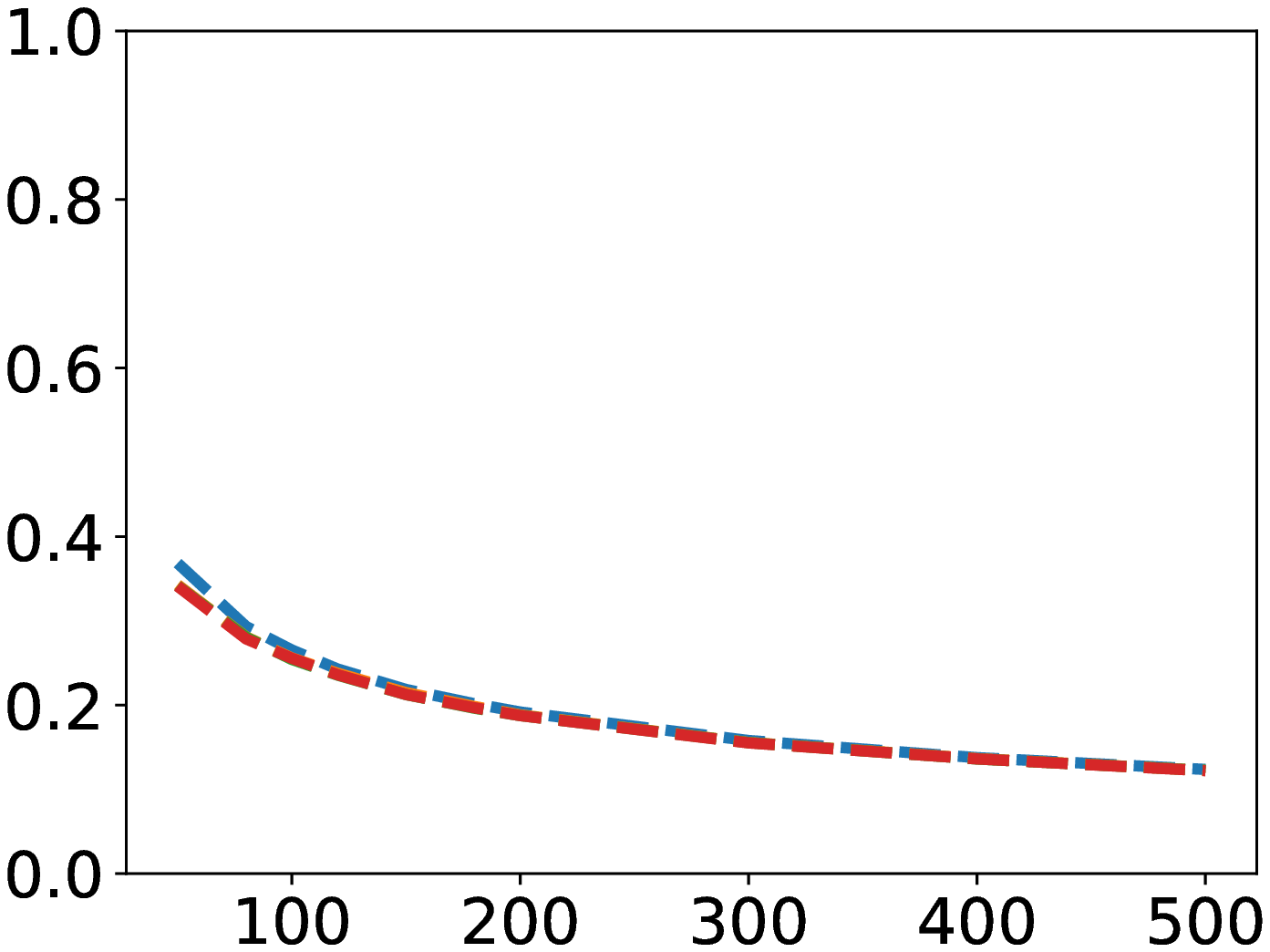} 
		\put(25,80){\ul{\ \ \  \ $\gamma=1.3$ \ \ \ \    }}
	\end{overpic}	
	%
\end{figure}

\vspace{-0.5cm}

\begin{figure}[H]	
	\quad\quad\quad 
	\begin{overpic}[width=0.29\textwidth]{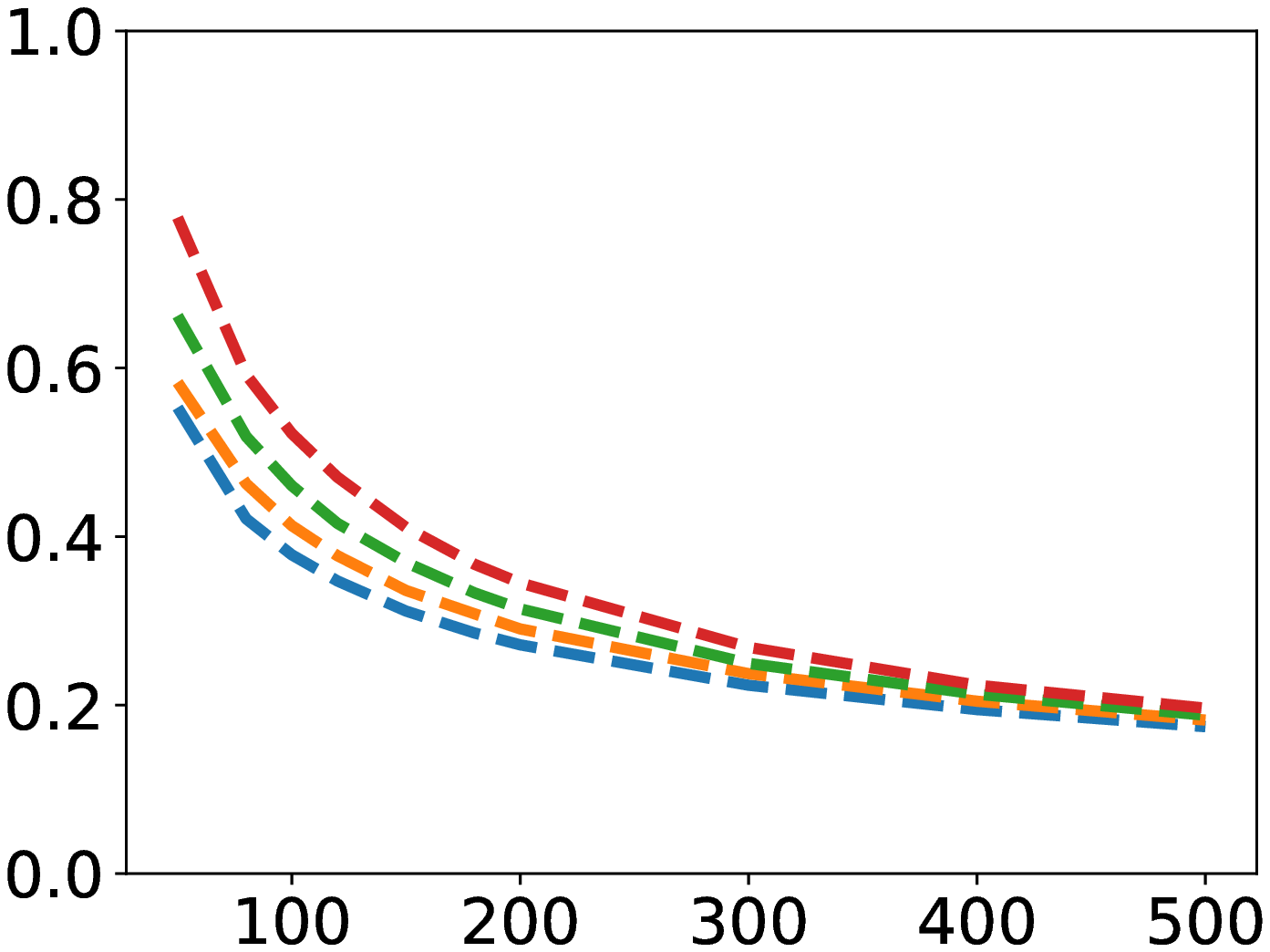} 
	\put(-20,-1){\rotatebox{90}{ {\small \ \ \ standardization \  \ \ }}}
	\end{overpic}
	~
	\DeclareGraphicsExtensions{.png}
	\begin{overpic}[width=0.29\textwidth]{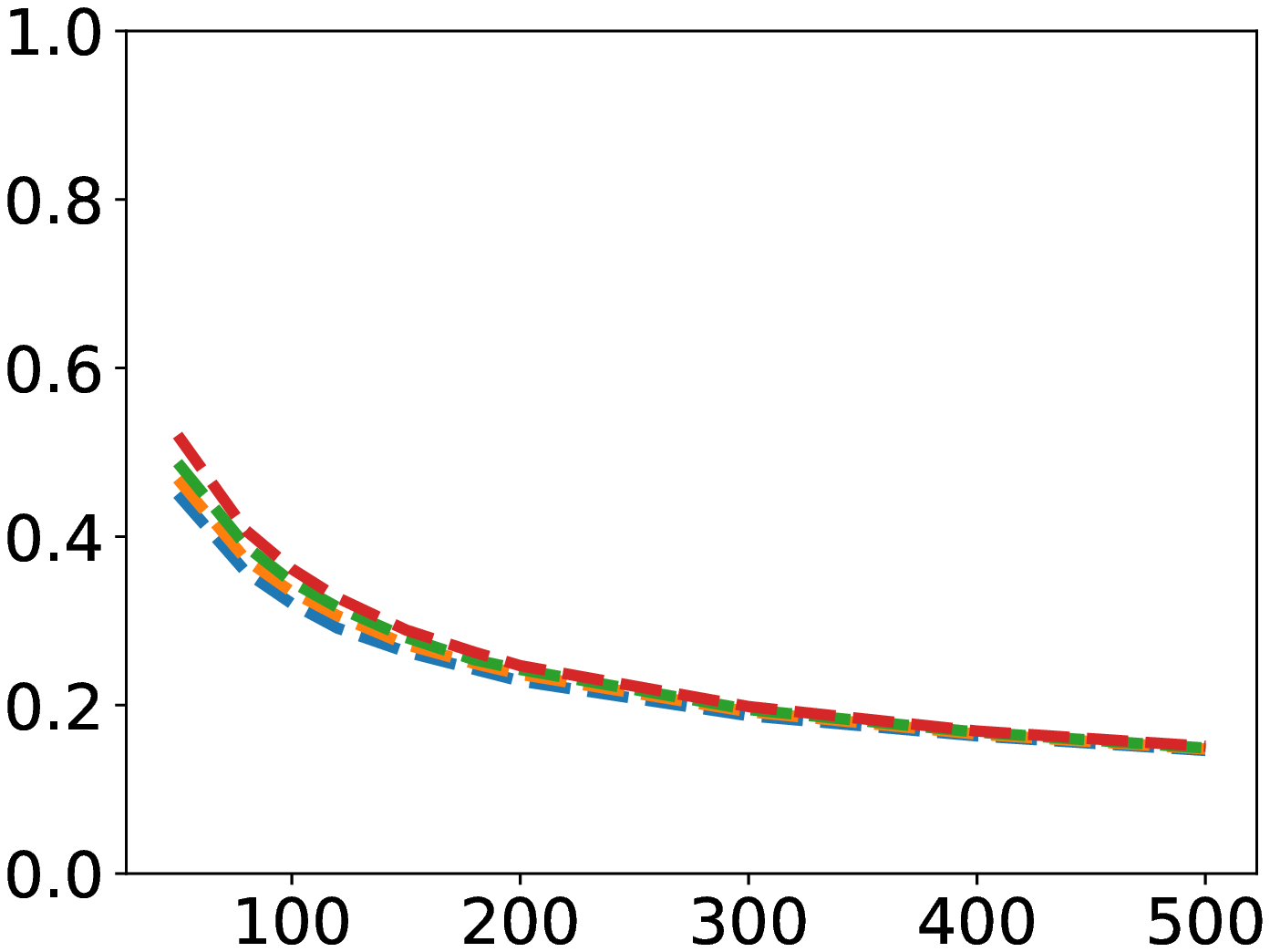} 
	\end{overpic}
	~	
	\begin{overpic}[width=0.29\textwidth]{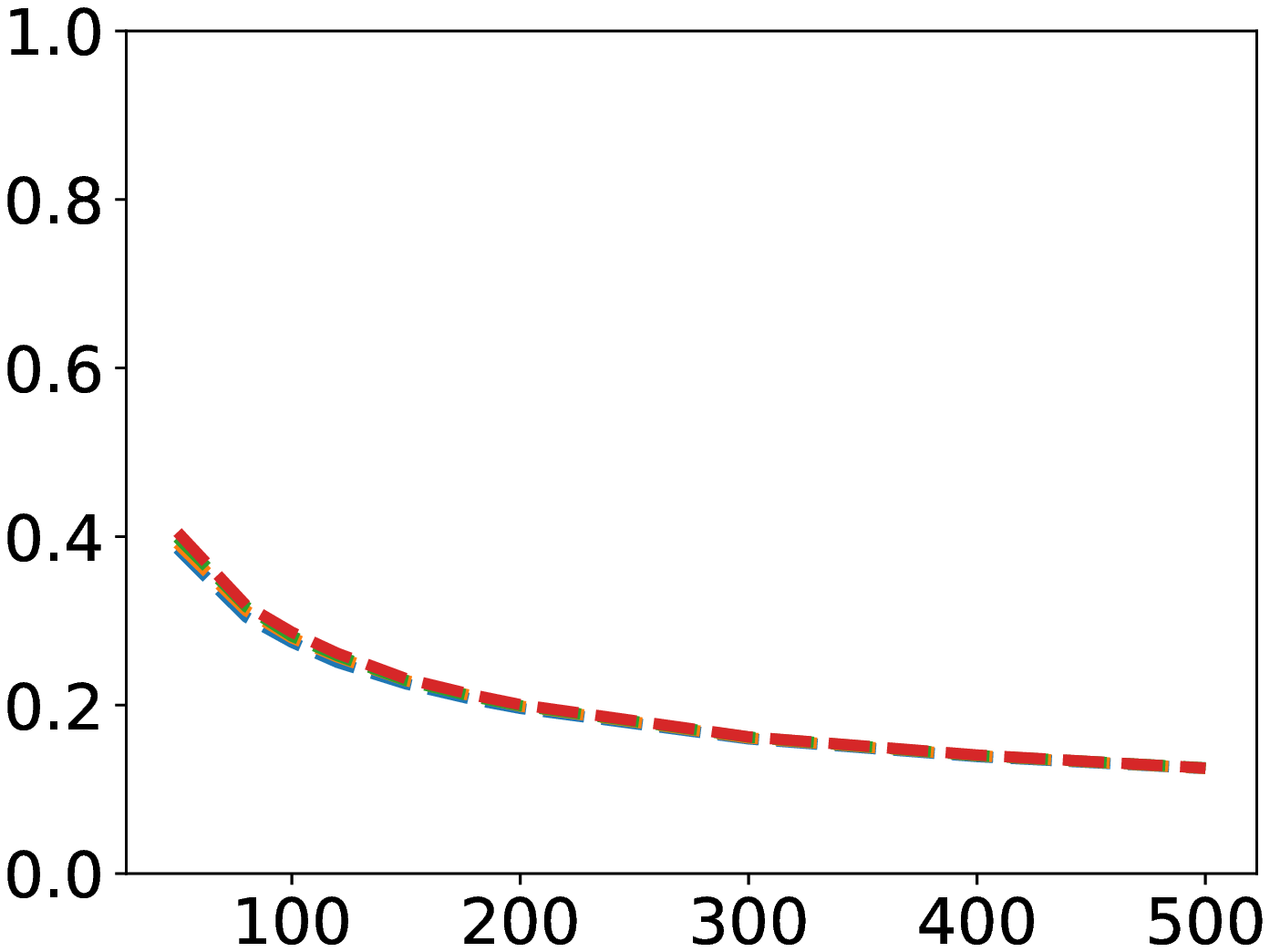} 
	\end{overpic}	

\end{figure}

\vspace{-0.5cm}

\begin{figure}[H]	
	\quad\quad\quad 
	\begin{overpic}[width=0.29\textwidth]{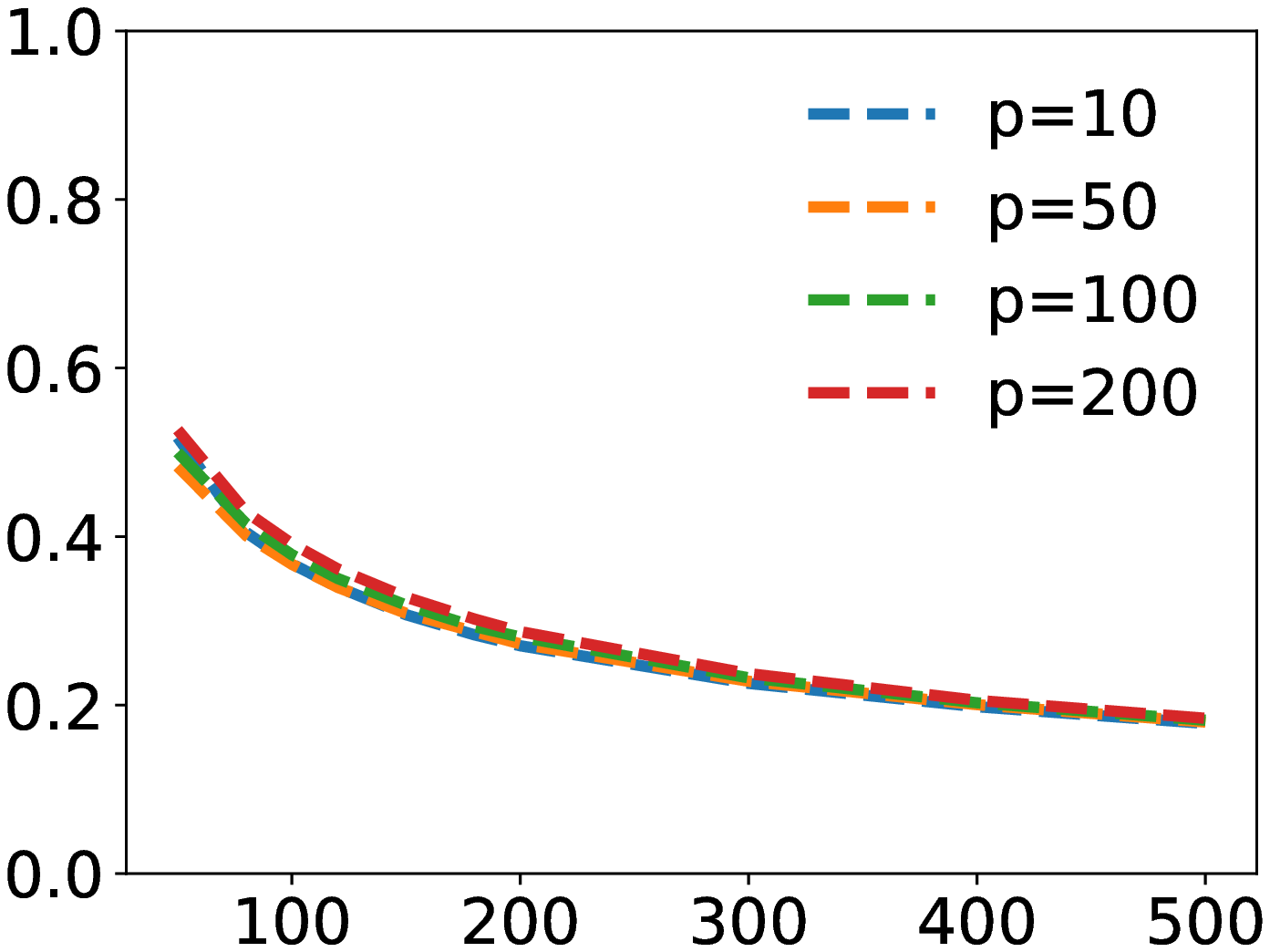} 
		\put(-21,1){\rotatebox{90}{\ $\sqrt{ \ \ }$}}
	    \put(-20,-3){\rotatebox{90}{  { \ \ \ \ \ \ \small transformation \ \ } }}

	\end{overpic}
	~
	\DeclareGraphicsExtensions{.png}
	\begin{overpic}[width=0.29\textwidth]{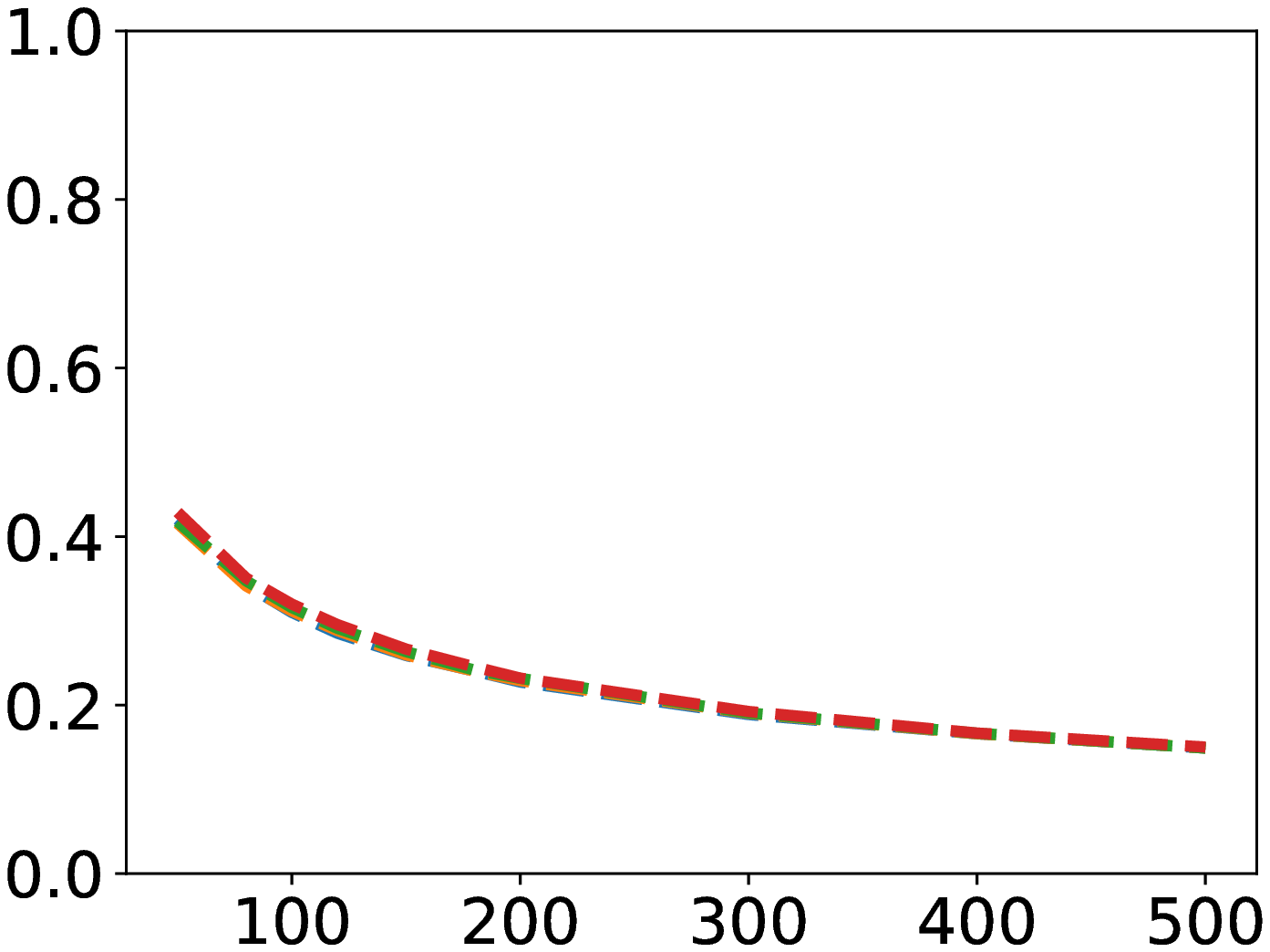} 
	\end{overpic}
	~	
	\begin{overpic}[width=0.29\textwidth]{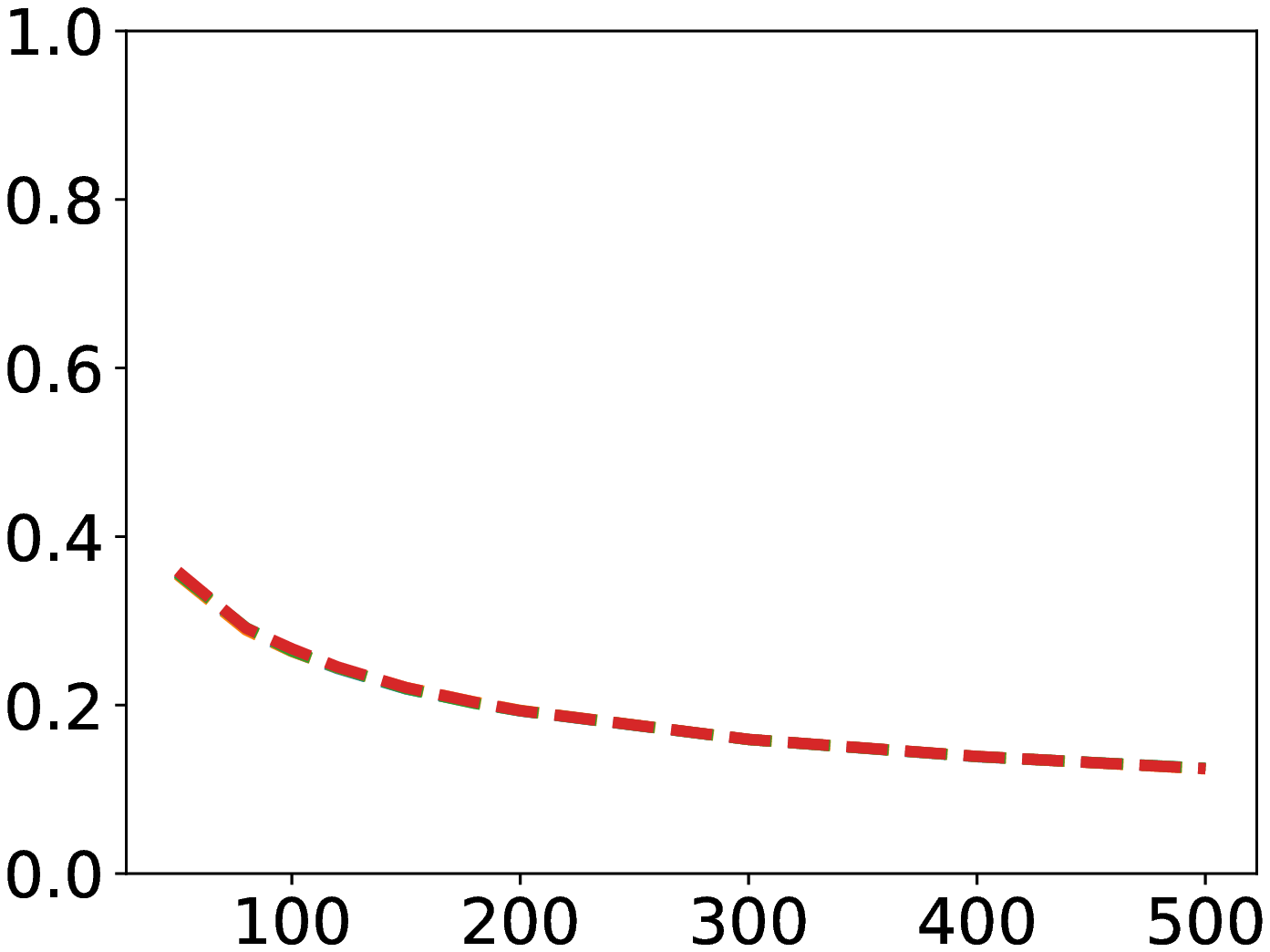} 
	\end{overpic}	
	\vspace{+0.2cm}	
	\caption{(Average width versus $n$ in simulation model (ii) with a polynomial decay profile). The plotting scheme is the same as described in the caption of Figure~\ref{fig5} in the main text.} 
	\label{fig_width_gamma_ii}
\end{figure}

\newpage

%

\begin{figure}[H]	
\vspace{0.5cm}
	\quad\quad\quad 
	\begin{overpic}[width=0.29\textwidth]{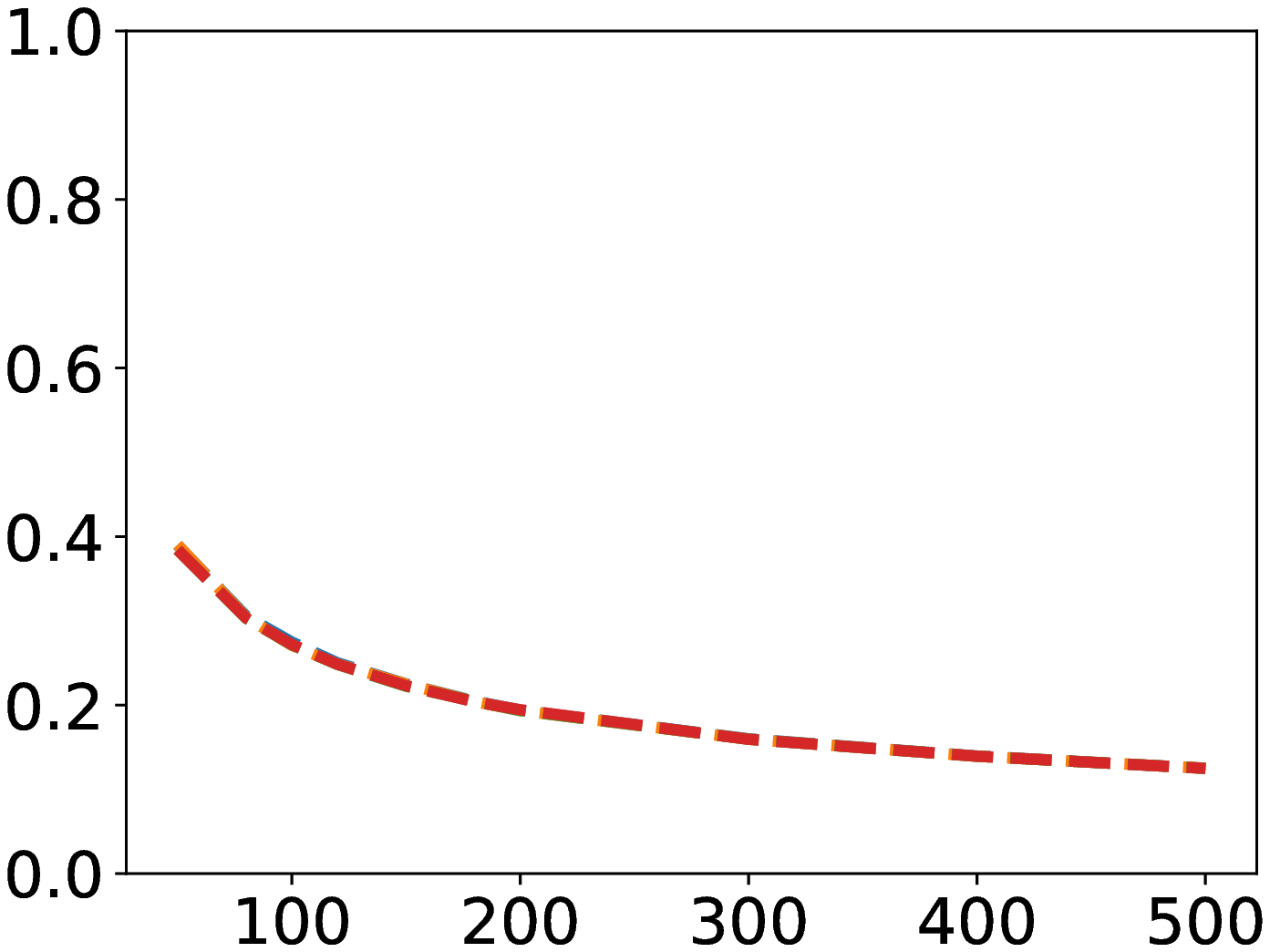} 
\put(25,80){ \ul{\ \ \  \ $\delta=0.7$ \ \ \ \    }}
		\put(-20,-5){\rotatebox{90}{ {\small \ \ \ log transformation  \ \ }}}
\end{overpic}
	~
	\DeclareGraphicsExtensions{.png}
	\begin{overpic}[width=0.29\textwidth]{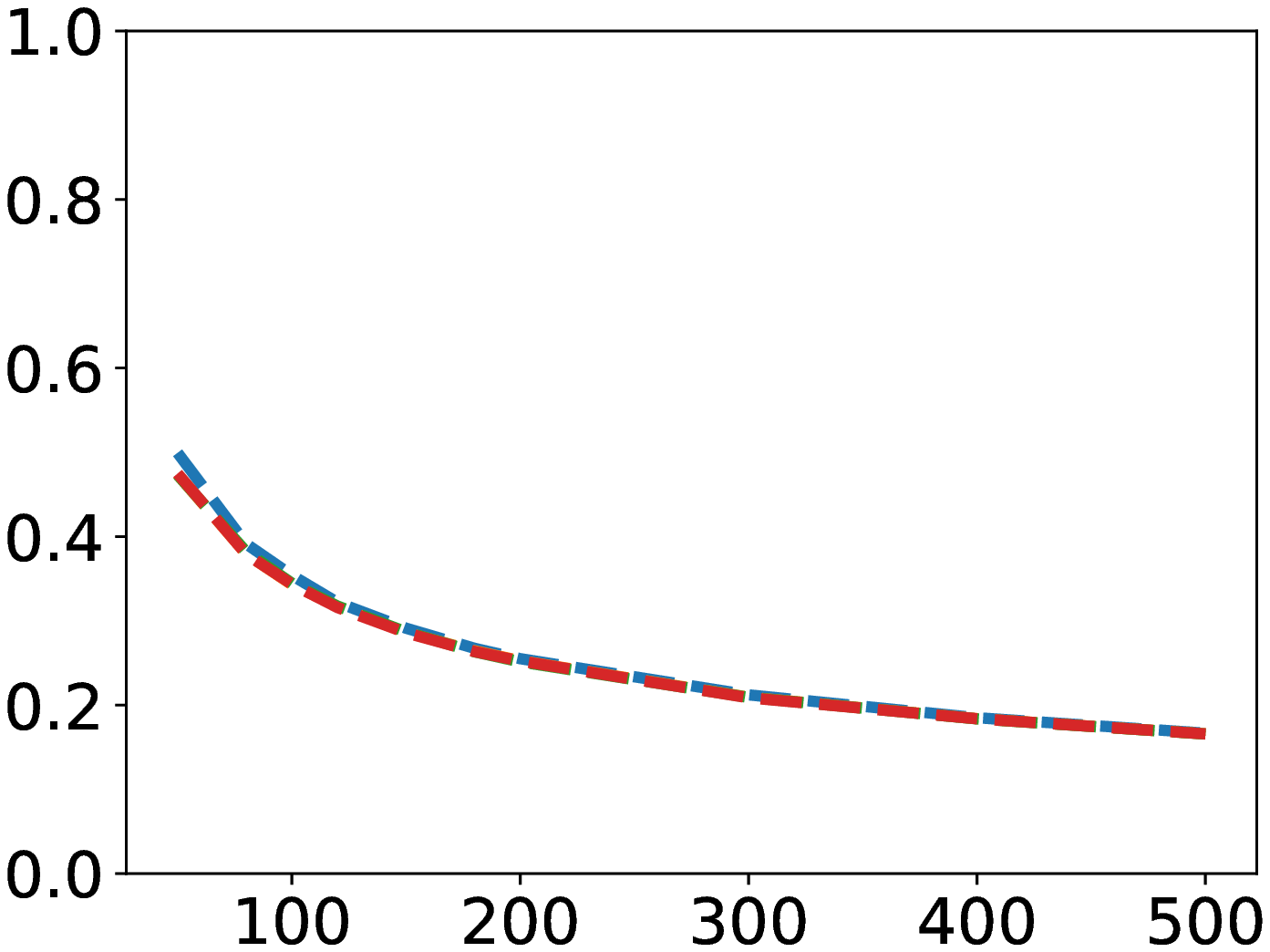} 
		\put(25,80){ \ul{\ \ \  \ $\delta=0.8$ \ \ \ \    }}
	\end{overpic}
	~	
	\begin{overpic}[width=0.29\textwidth]{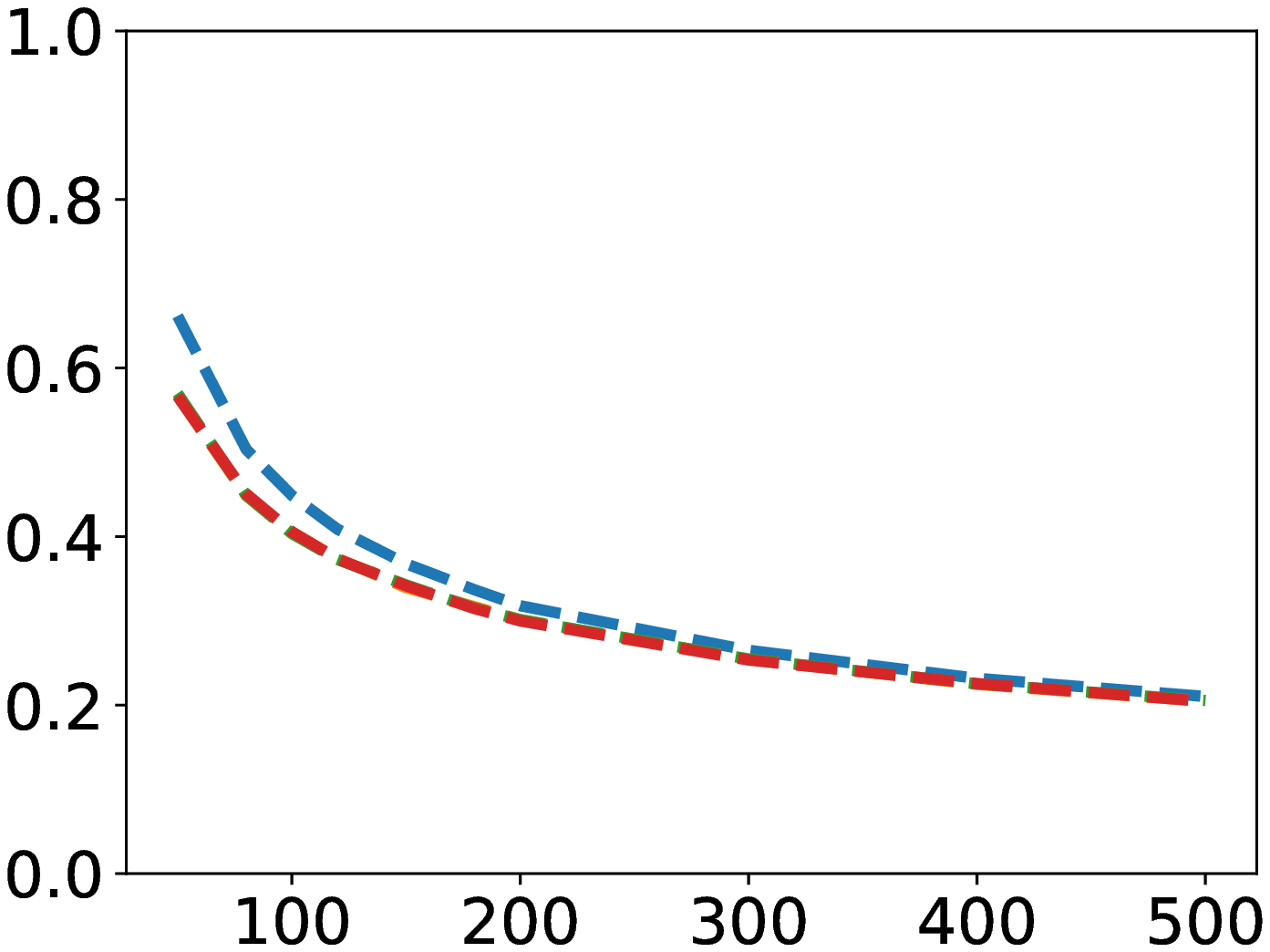} 
		\put(25,80){ \ul{\ \ \  \ $\delta=0.9$ \ \ \ \    }}
	\end{overpic}	
\end{figure}

\vspace{-0.5cm}

\begin{figure}[H]	
	\quad\quad\quad 
	\begin{overpic}[width=0.29\textwidth]{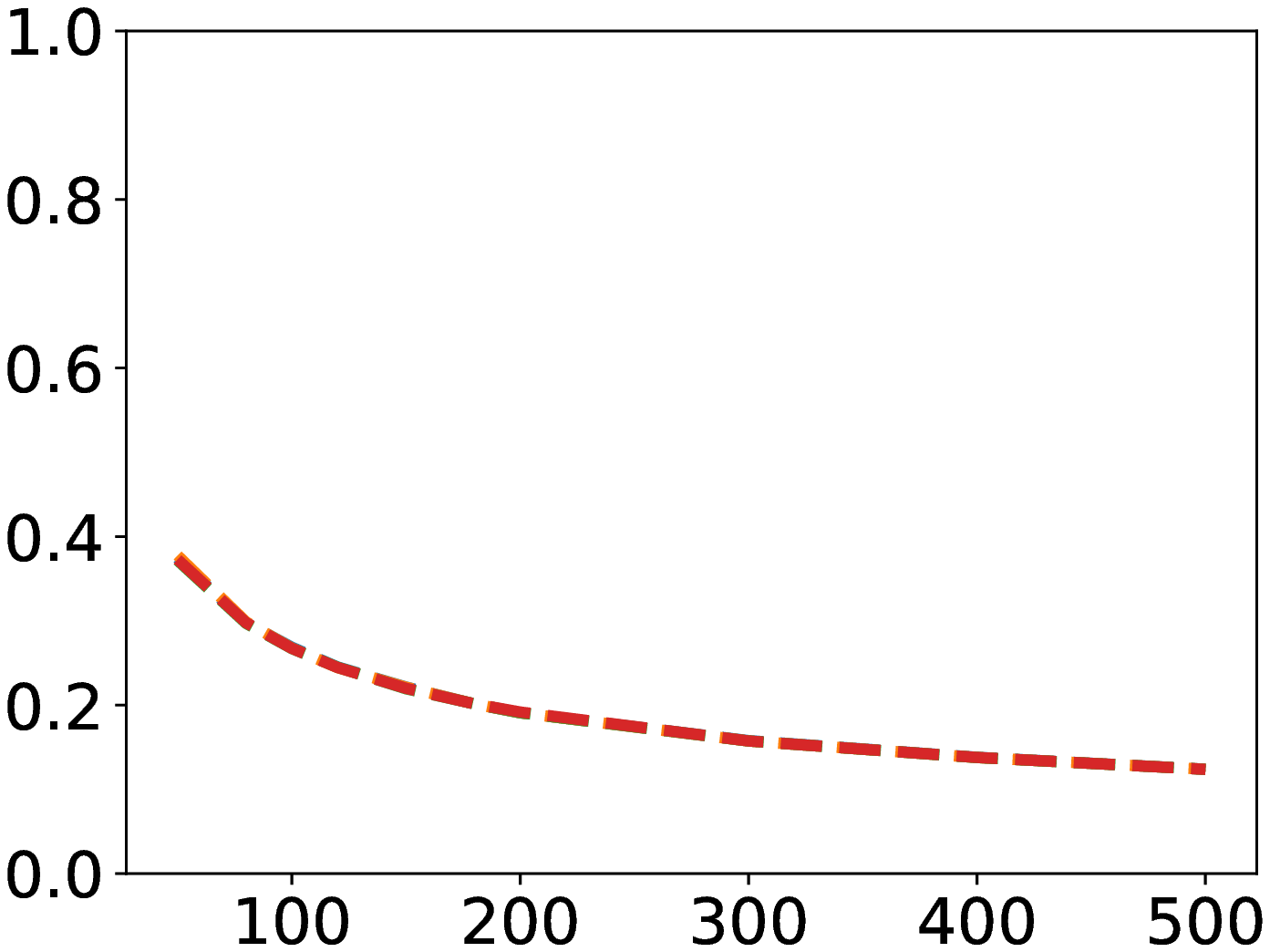} 
	\put(-20,-1){\rotatebox{90}{ {\small \ \ \ standardization \  \ \ }}}
	\end{overpic}
	~
	\DeclareGraphicsExtensions{.png}
	\begin{overpic}[width=0.29\textwidth]{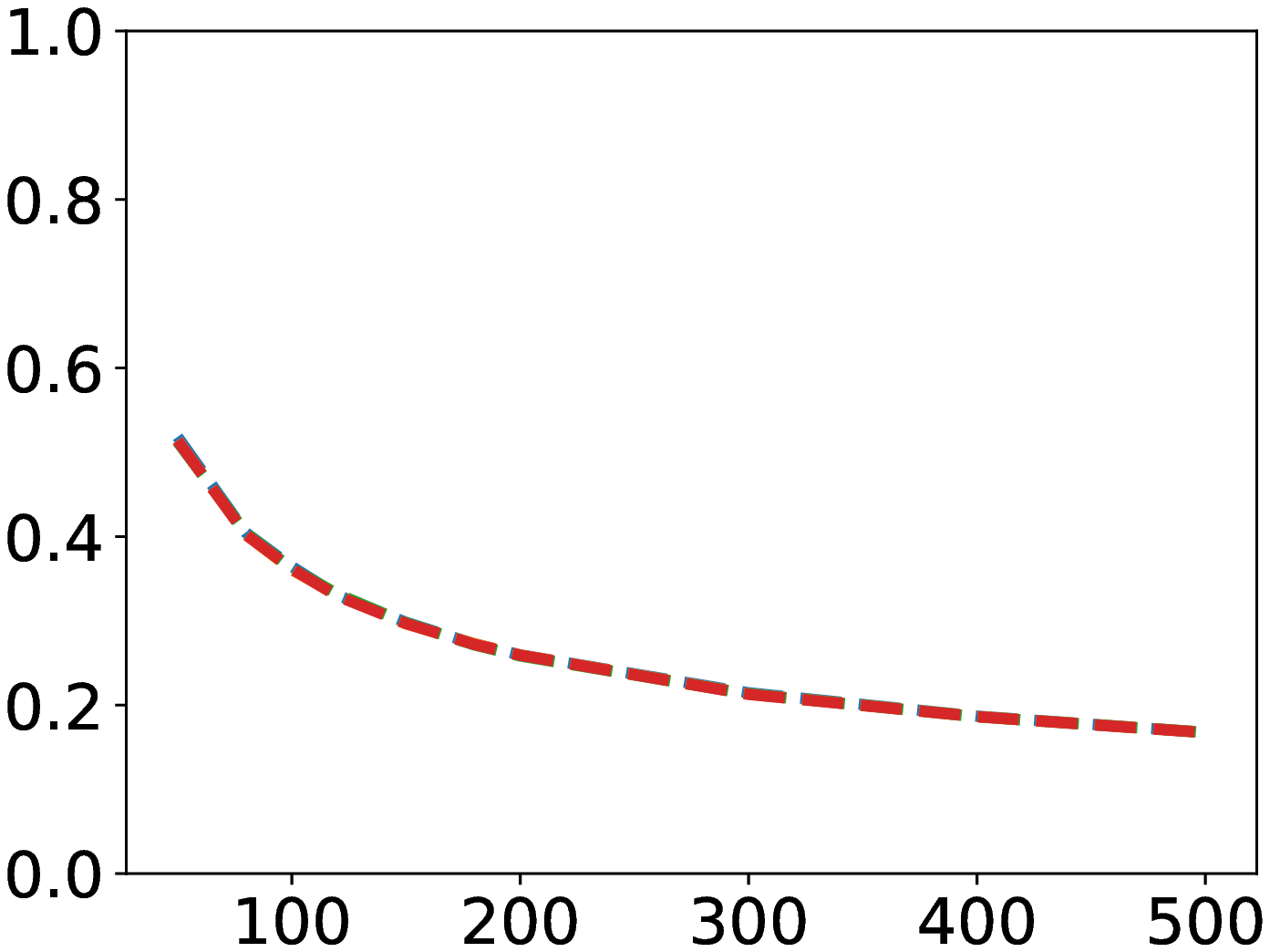} 
	\end{overpic}
	~	
	\begin{overpic}[width=0.29\textwidth]{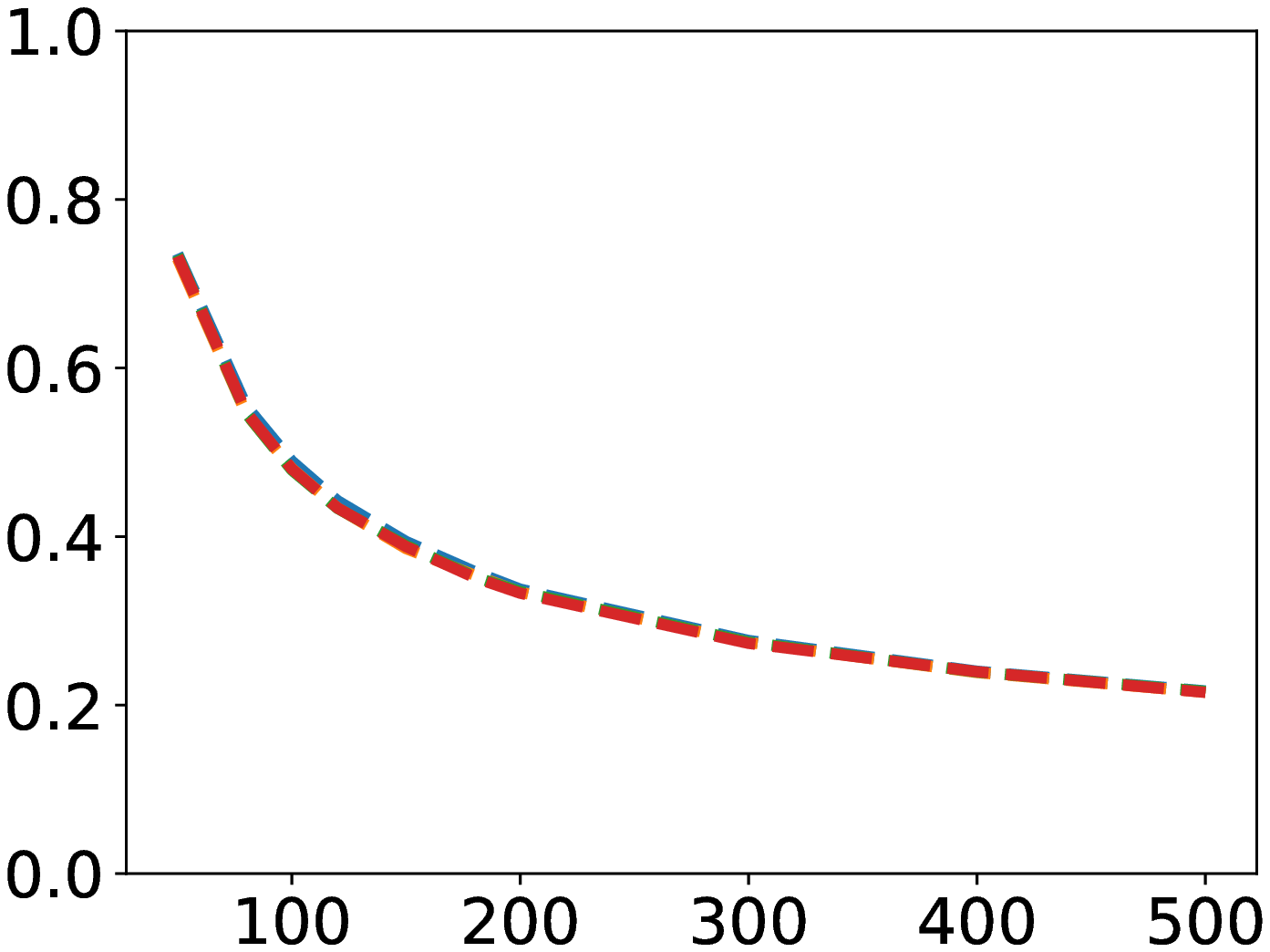} 
	\end{overpic}	
\end{figure}

\vspace{-0.5cm}

\begin{figure}[H]	
	\quad\quad\quad 
	\begin{overpic}[width=0.29\textwidth]{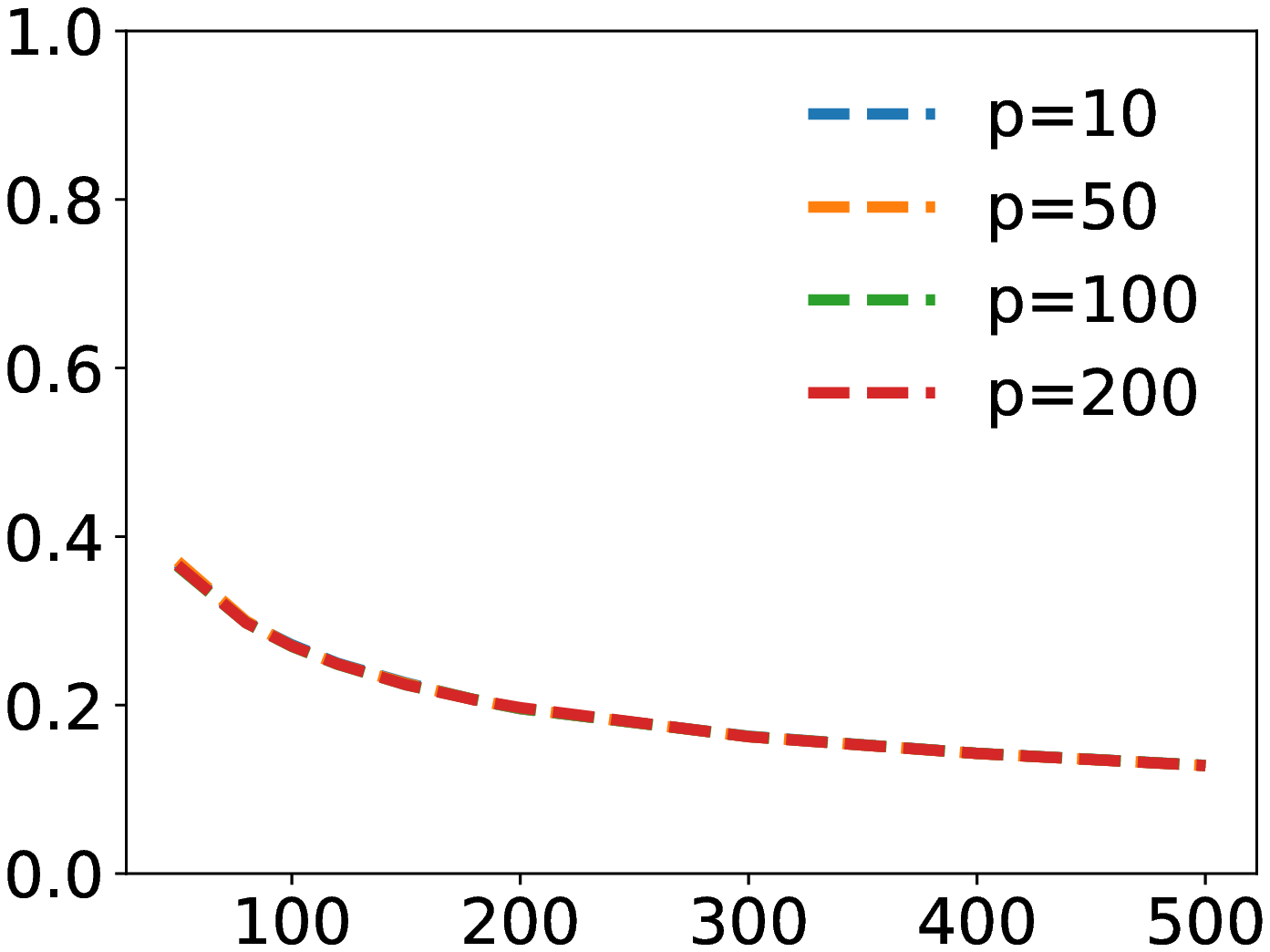} 
		\put(-21,1){\rotatebox{90}{\ $\sqrt{ \ \ }$}}
	    \put(-20,-3){\rotatebox{90}{  { \ \ \ \ \ \ \small transformation \ \ } }}

	\end{overpic}
	~
	\DeclareGraphicsExtensions{.png}
	\begin{overpic}[width=0.29\textwidth]{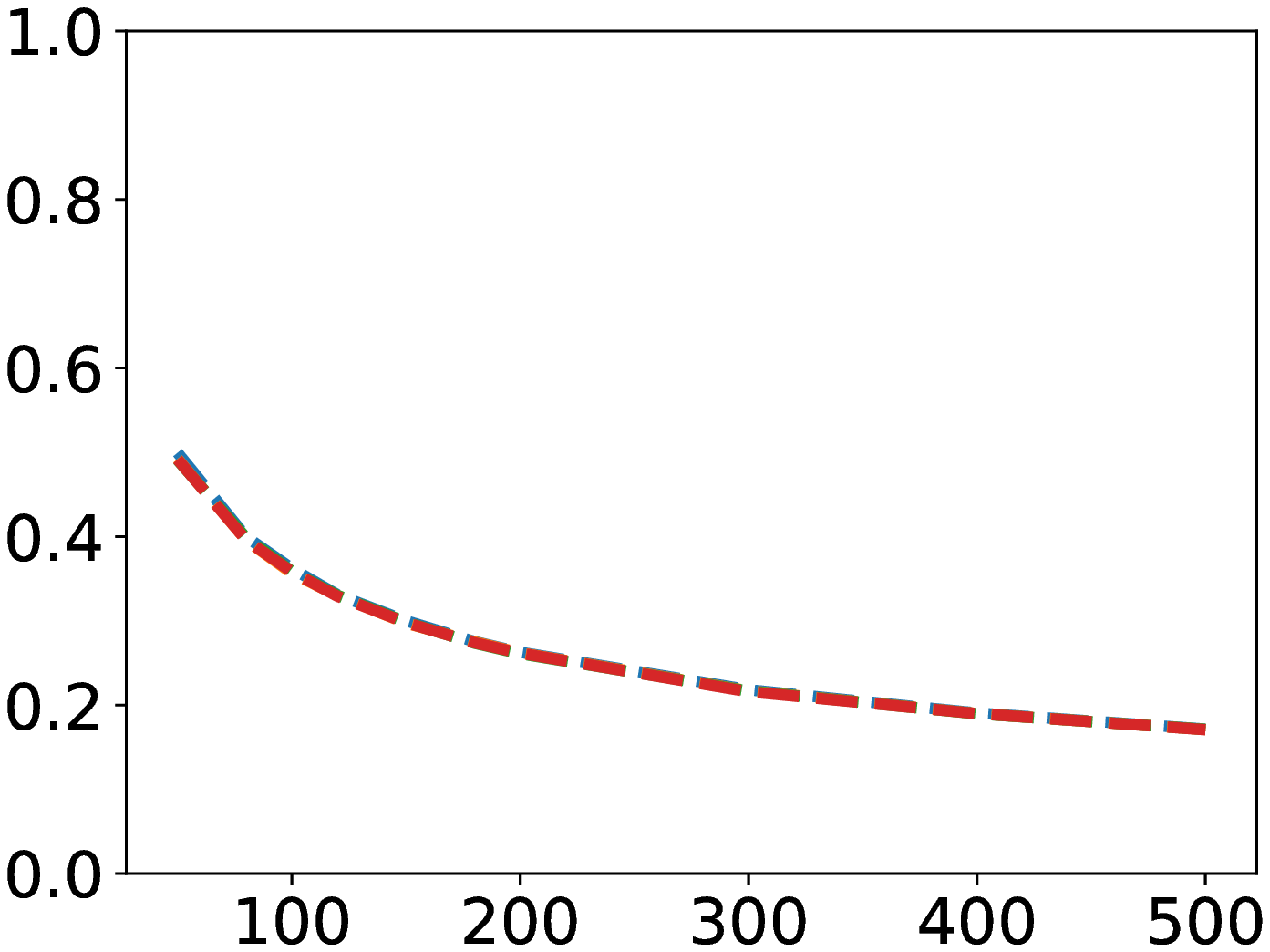} 
	\end{overpic}
	~	
	\begin{overpic}[width=0.29\textwidth]{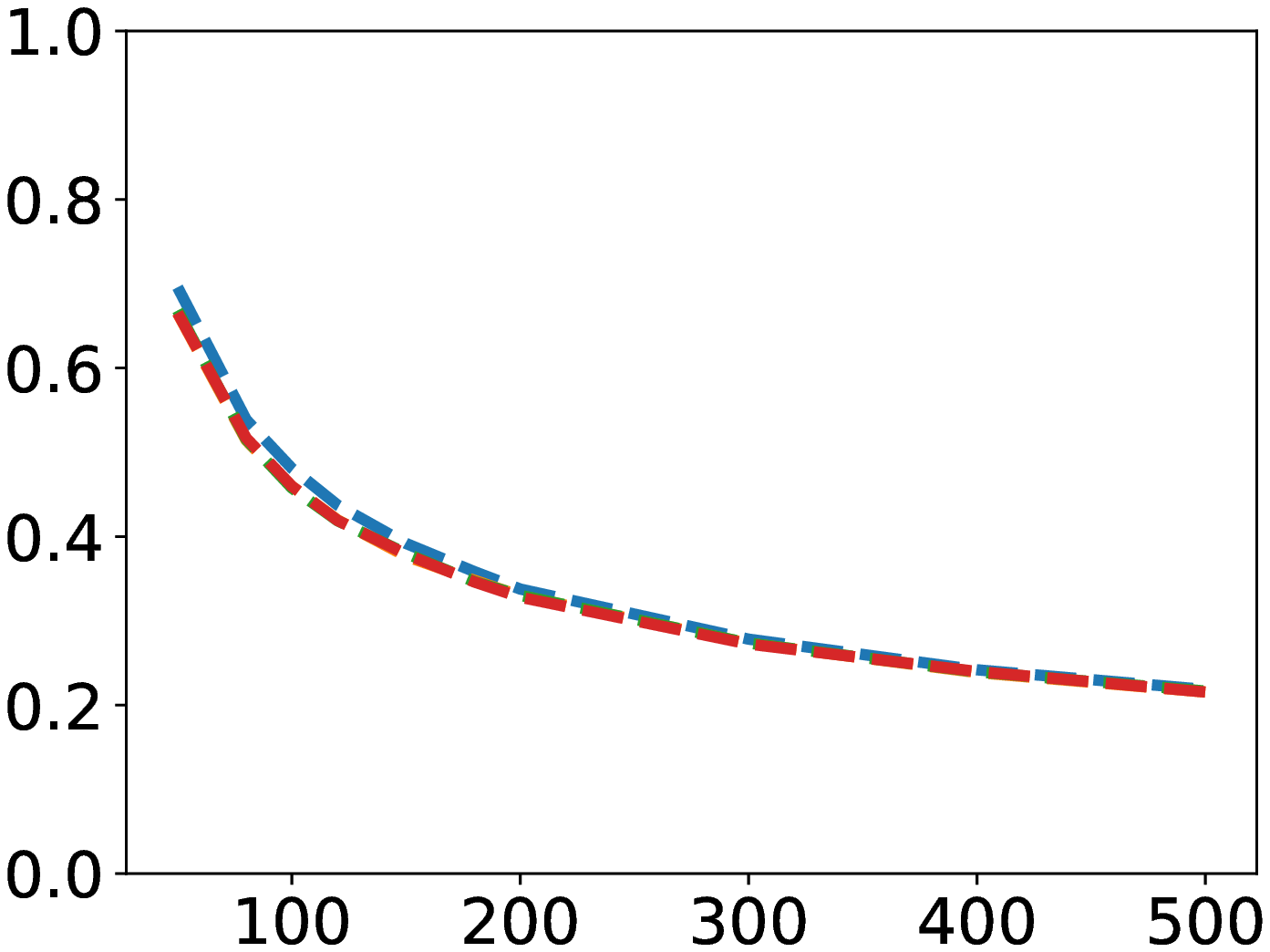} 
		 				
	\end{overpic}	
    \vspace{+0.2cm}	
	\caption{(Average width versus $n$ in simulation model (ii) with an exponential decay profile). The plotting scheme is the same as described in the caption of Figure~\ref{fig7} in the main text.}
	\label{fig_width_delta_ii}
\end{figure}

\newpage

\paragraph{Nominal value of 90\% in models (i) and (ii).}
The following eight figures are presented in the same manner as in the main text, except that they use a nominal value of 90\%, and are based on both models (i) and (ii).

\begin{figure}[H]	
\vspace{0.5cm}
	\quad\quad\quad 
	\begin{overpic}[width=0.29\textwidth]{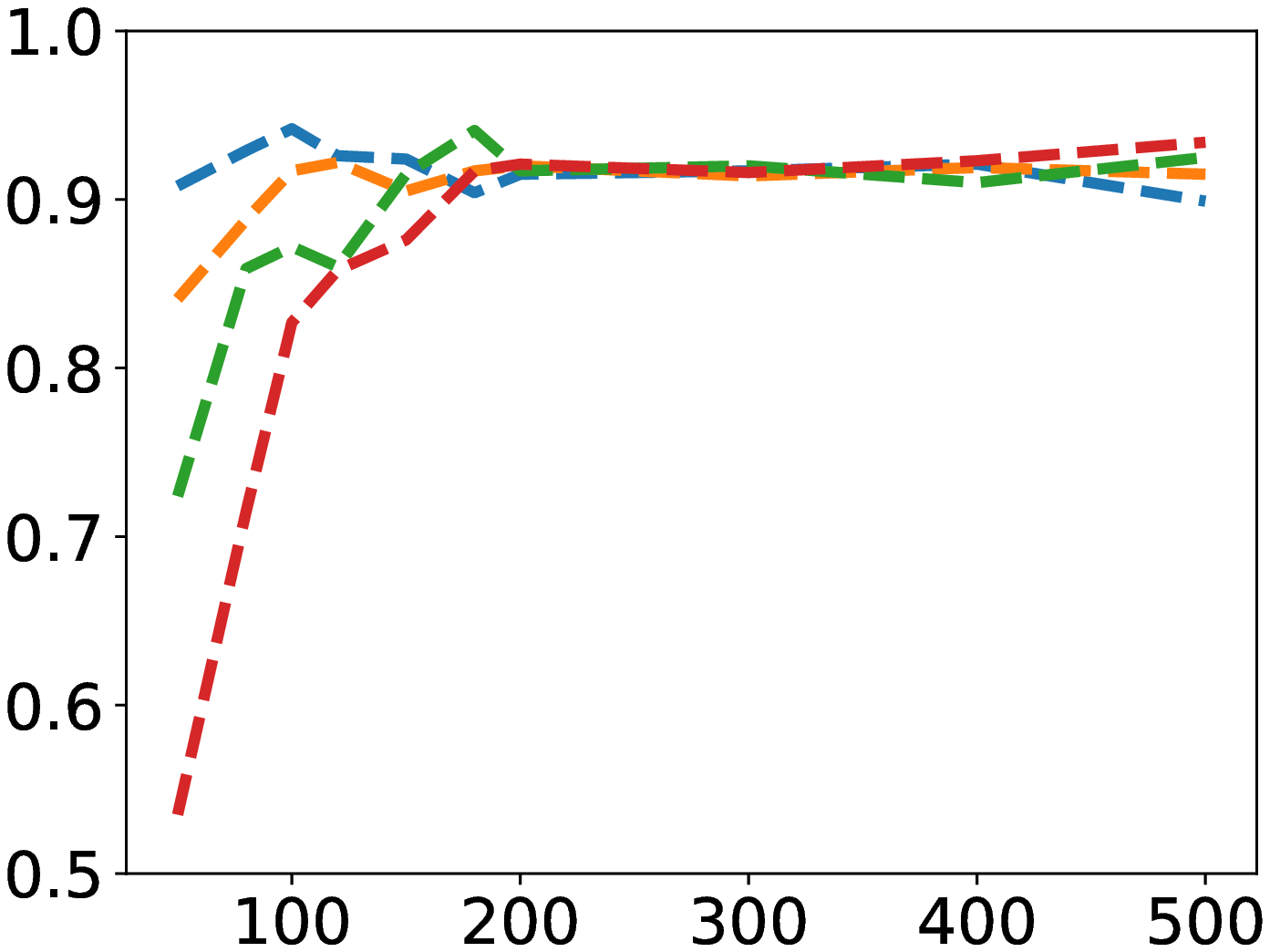} 
    \put(25,80){ \ul{\ \ \  \ $\gamma=0.7$ \ \ \ \    }}
	\put(-20,-5){\rotatebox{90}{ {\small \ \ \ log transformation  \ \ }}}
\end{overpic}
	~
	\DeclareGraphicsExtensions{.png}
	\begin{overpic}[width=0.28\textwidth]{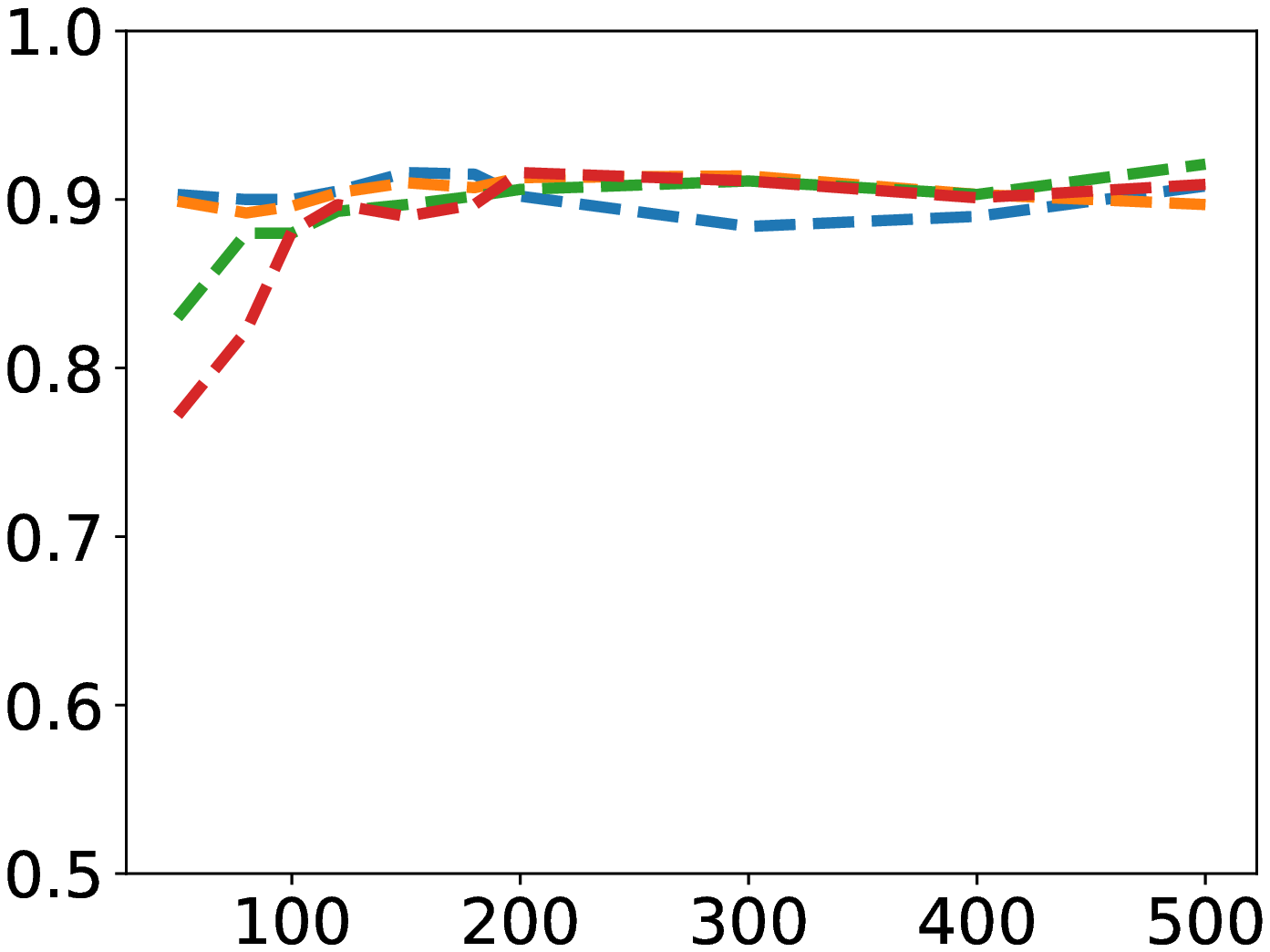} 
	\put(25,80){ \ul{\ \ \  \ $\gamma=1.0$ \ \ \ \    }}
	\end{overpic}
	~	
	\begin{overpic}[width=0.28\textwidth]{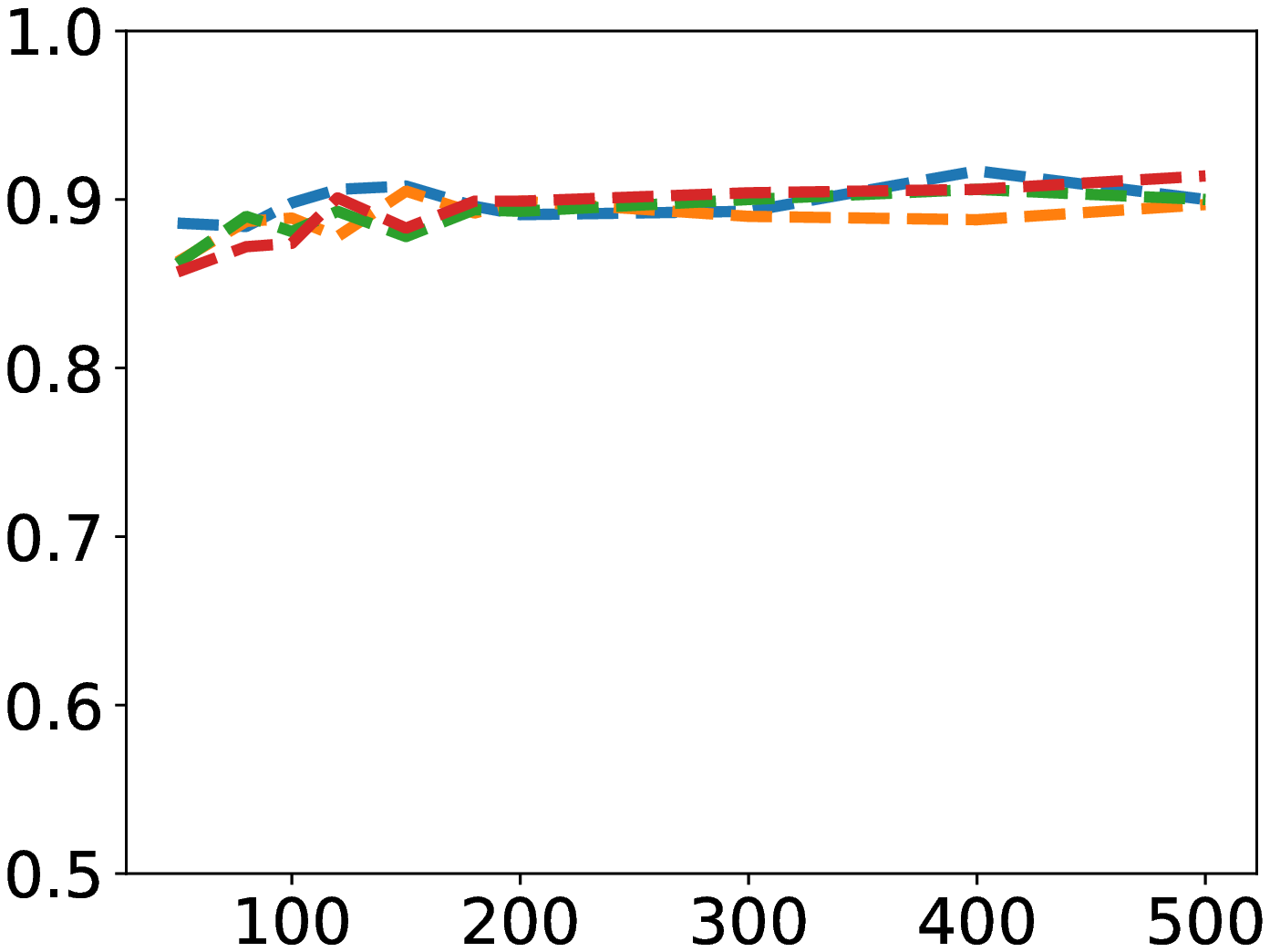} 
	\put(25,80){ \ul{\ \ \  \ $\gamma=1.3$ \ \ \ \    }}
	\end{overpic}	
\end{figure}

\vspace{-0.5cm}

\begin{figure}[H]	
	\quad\quad\quad 
	\begin{overpic}[width=0.29\textwidth]{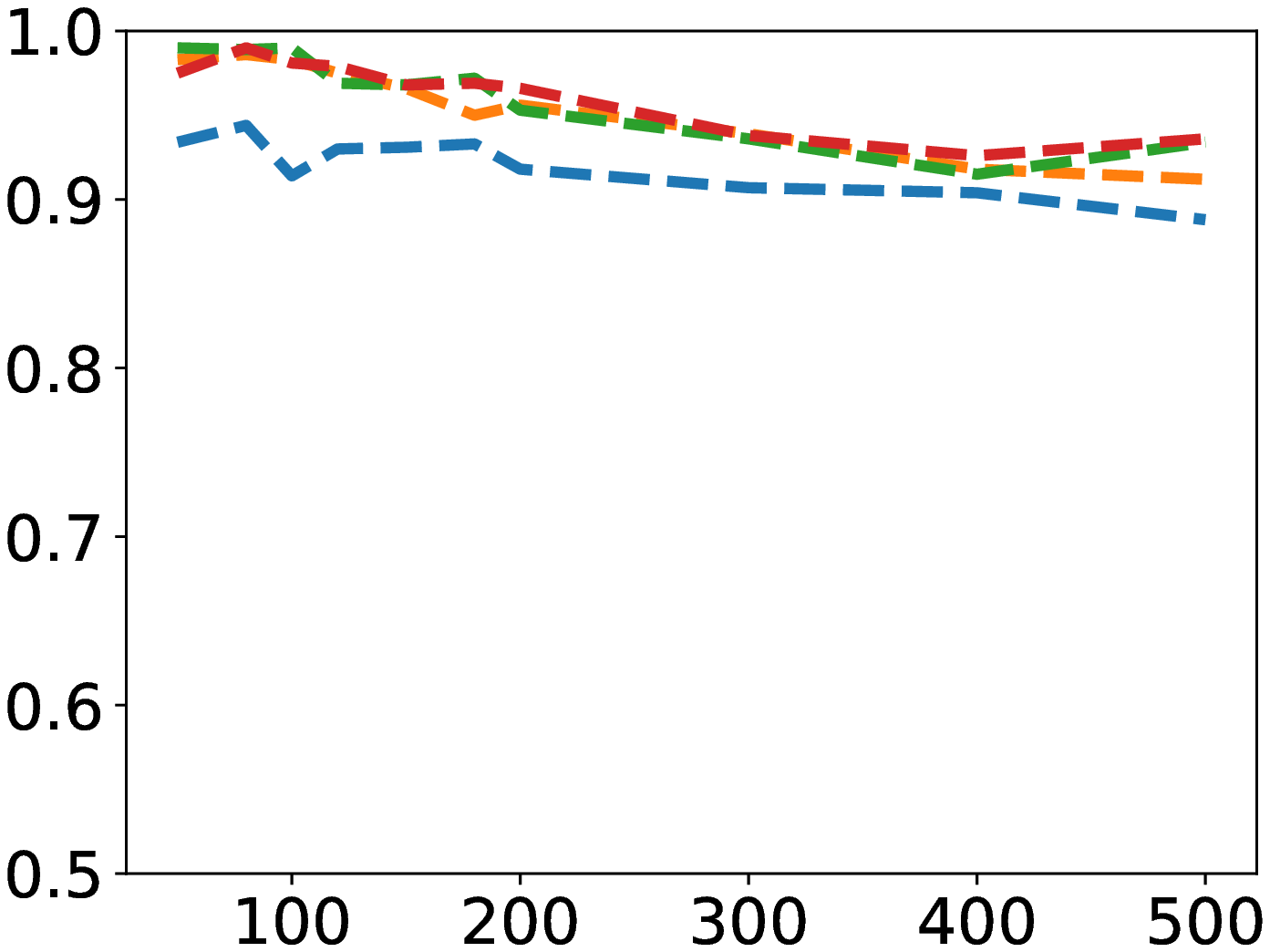} 
	\put(-20,-1){\rotatebox{90}{ {\small \ \ \ standardization \  \ \ }}}
	\end{overpic}
	~
	\DeclareGraphicsExtensions{.png}
	\begin{overpic}[width=0.29\textwidth]{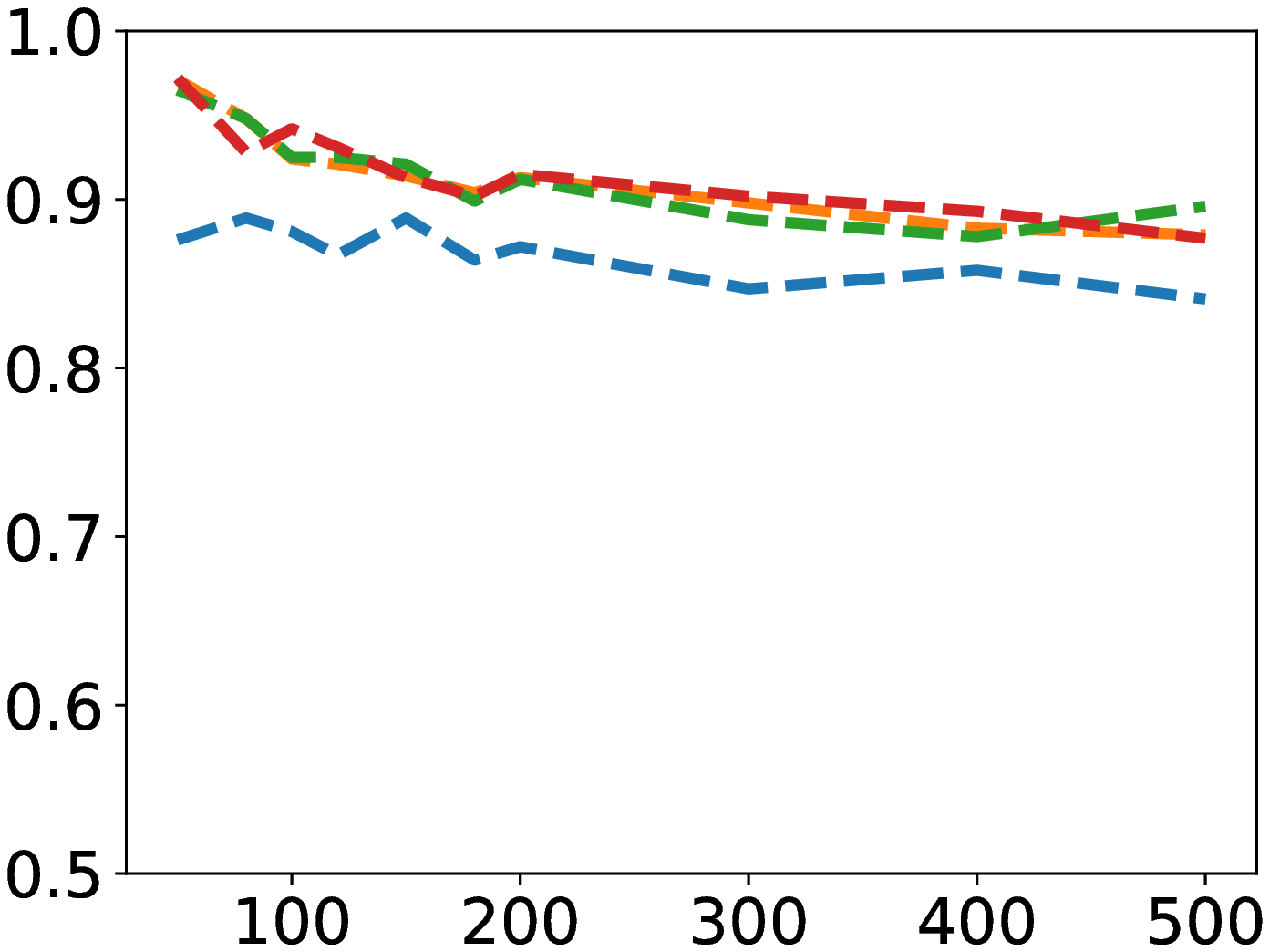} 
	\end{overpic}
	~	
	\begin{overpic}[width=0.29\textwidth]{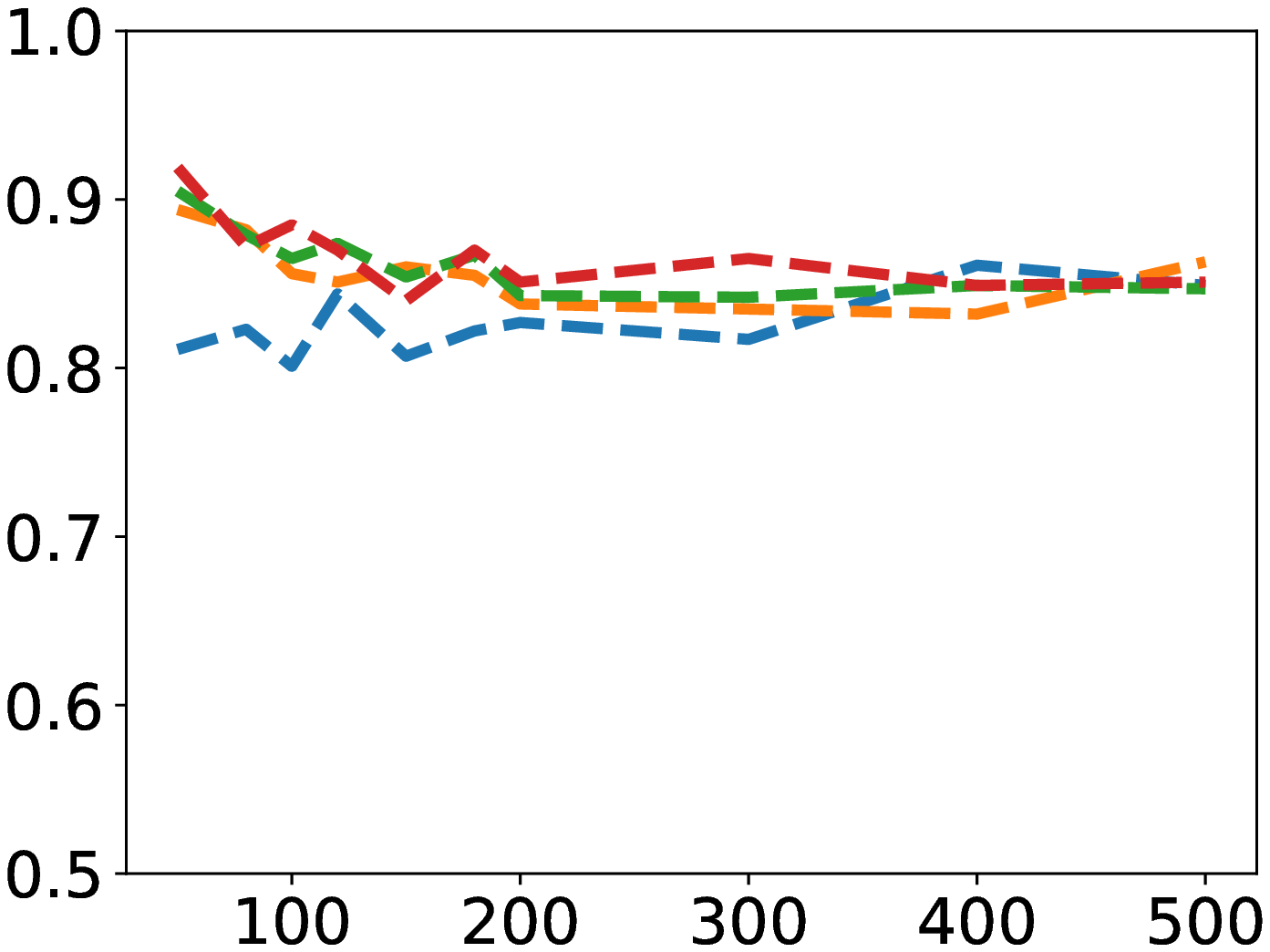} 
	\end{overpic}	
\end{figure}

\vspace{-0.5cm}

\begin{figure}[H]	
	\quad\quad\quad 
	\begin{overpic}[width=0.29\textwidth]{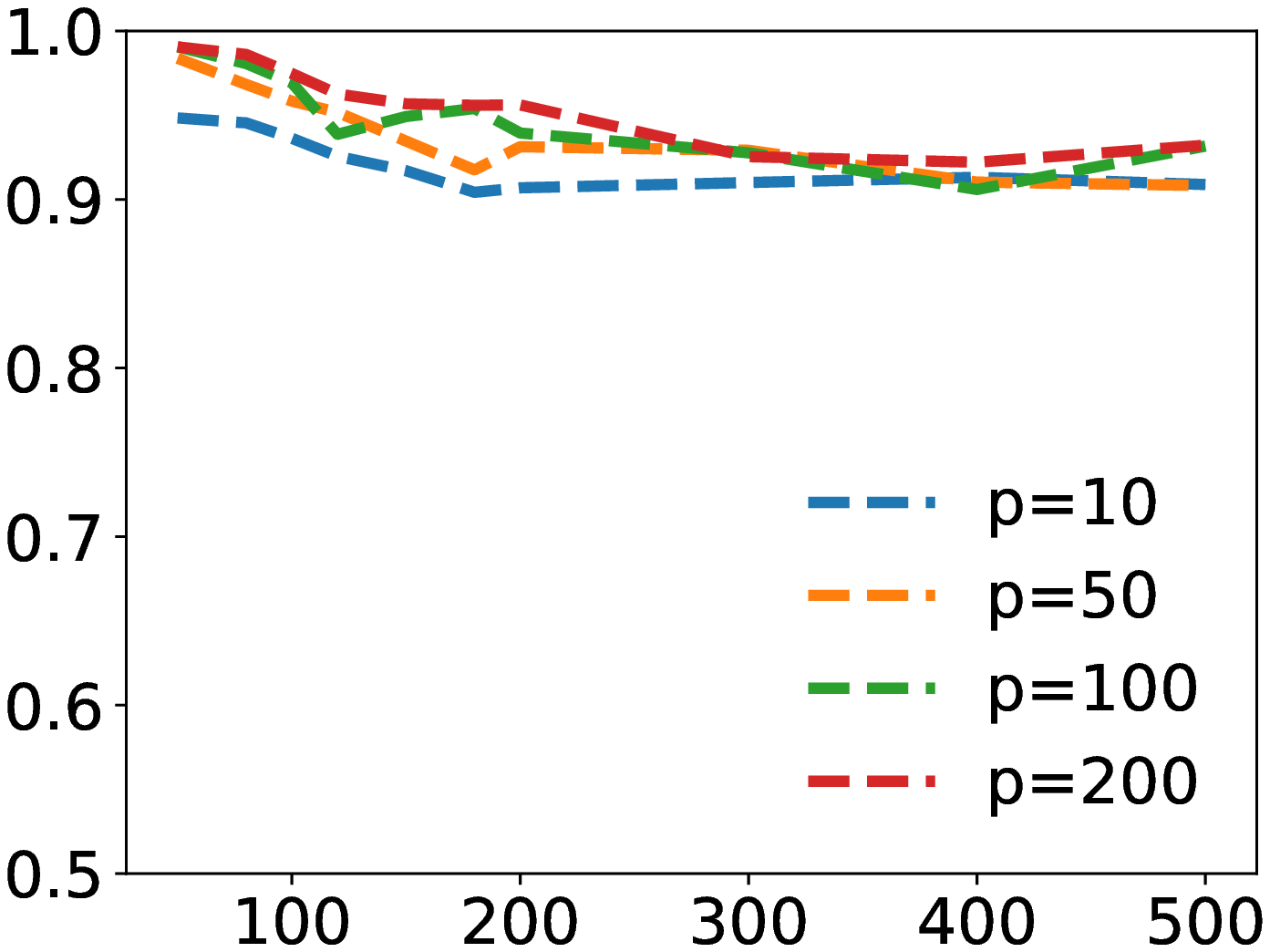} 
	\put(-21,1){\rotatebox{90}{\ $\sqrt{ \ \ }$}}
	\put(-20,-3){\rotatebox{90}{  { \ \ \ \ \ \ \small transformation \ \ } }}

	\end{overpic}
	~
	\DeclareGraphicsExtensions{.png}
	\begin{overpic}[width=0.29\textwidth]{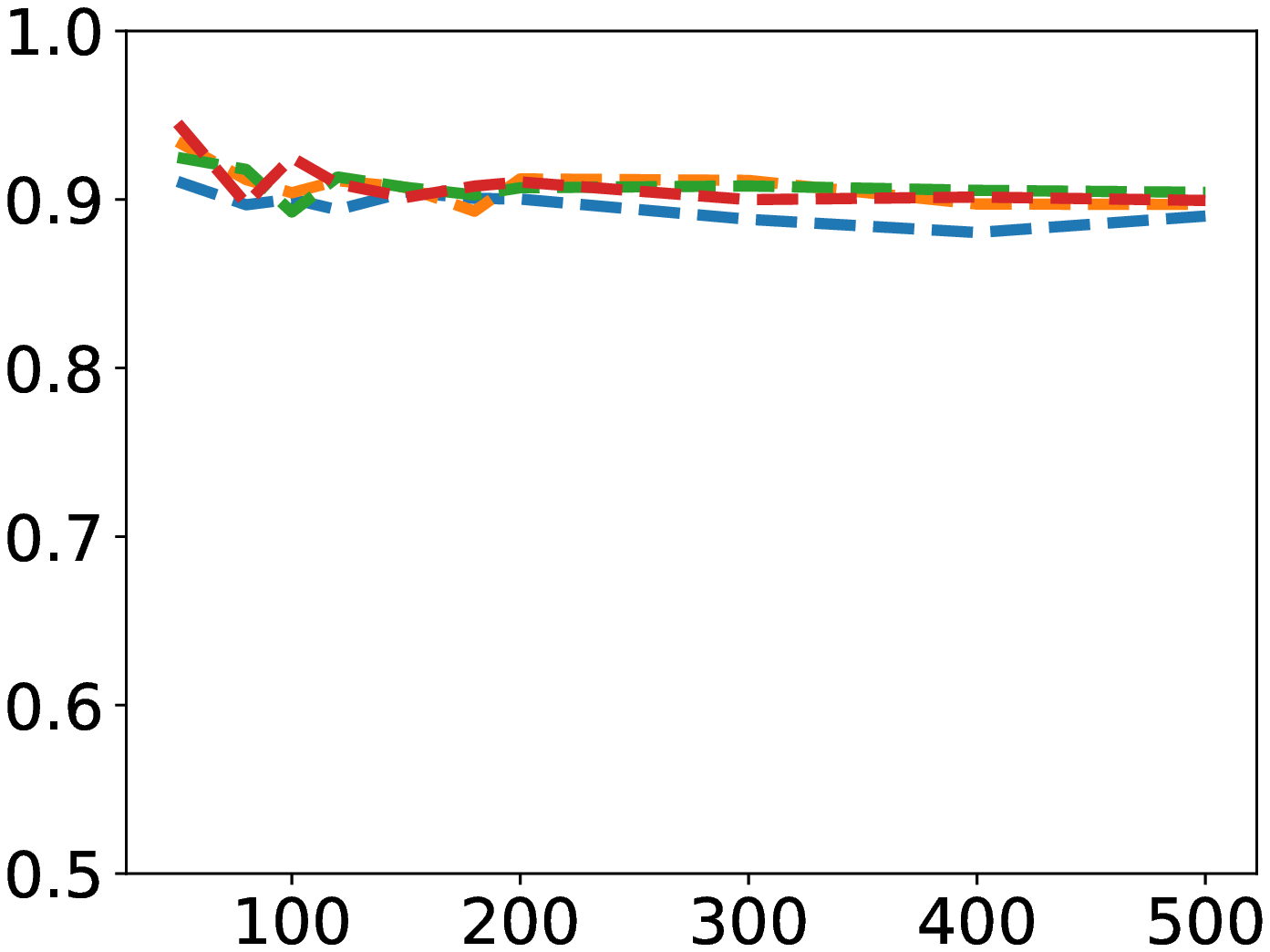} 
	\end{overpic}
	~	
	\begin{overpic}[width=0.29\textwidth]{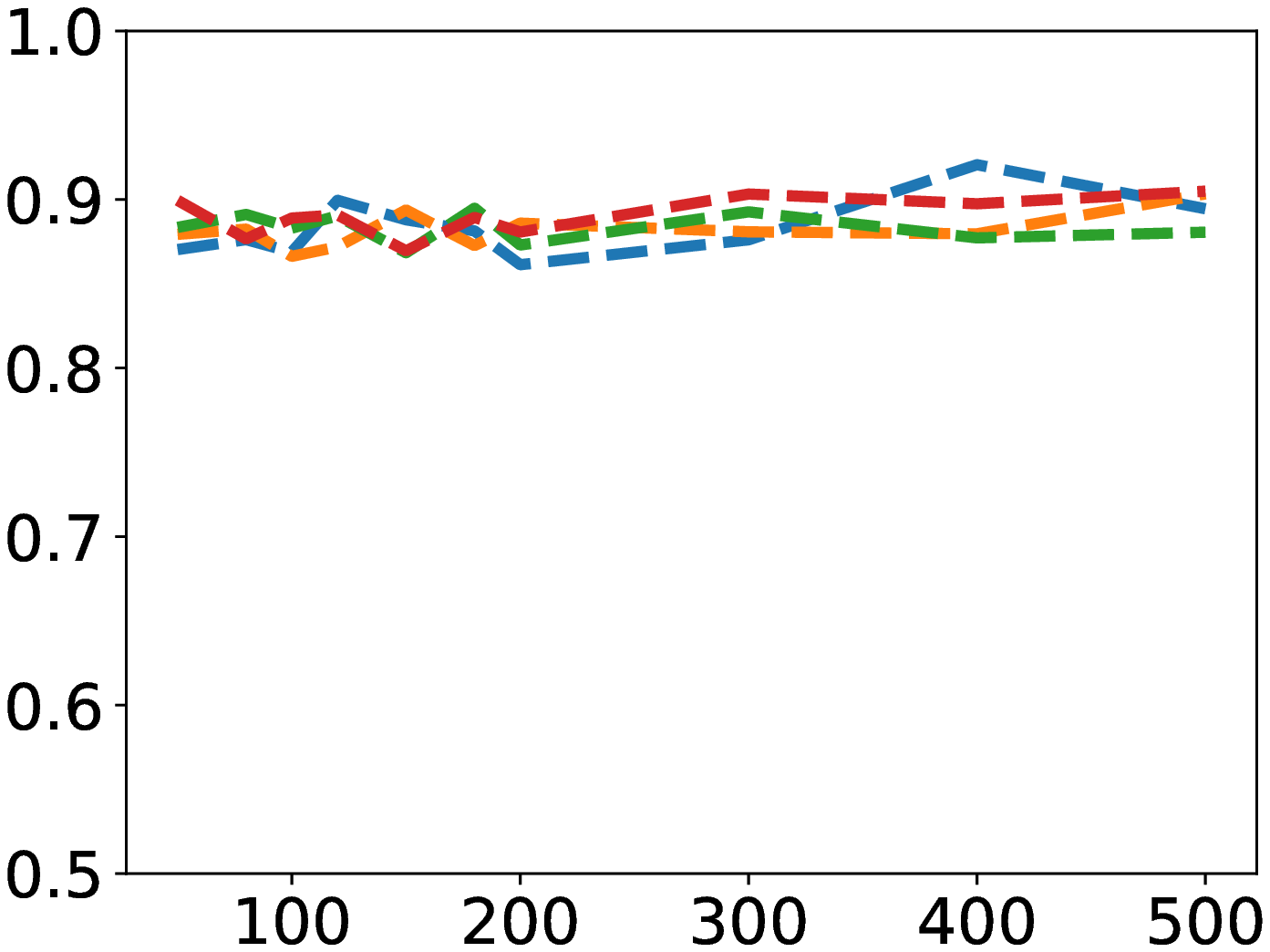} 
	\end{overpic}	
    \vspace{+0.2cm}	
	\caption{(Simultaneous coverage probability versus $n$ in simulation model (i) with a polynomial decay profile). In each panel, the $y$-axis measures~$\P(\cap_{j=1}^5 \{\lambda_j(\Sigma)\in\hat{\mathcal{I}}_j\})$ based on a nominal value of 90\%, and the $x$-axis measures $n$. The colored curves correspond to the different values of $p$, indicated in the legend. The three rows and three columns correspond to labeled choices of transformations and values of the eigenvalue decay parameter $\gamma$. }
	\label{SUPP:fig1} 
\end{figure}

\newpage
\begin{figure}[h]	
\vspace{0.5cm}
	\quad\quad\quad 
	\begin{overpic}[width=0.29\textwidth]{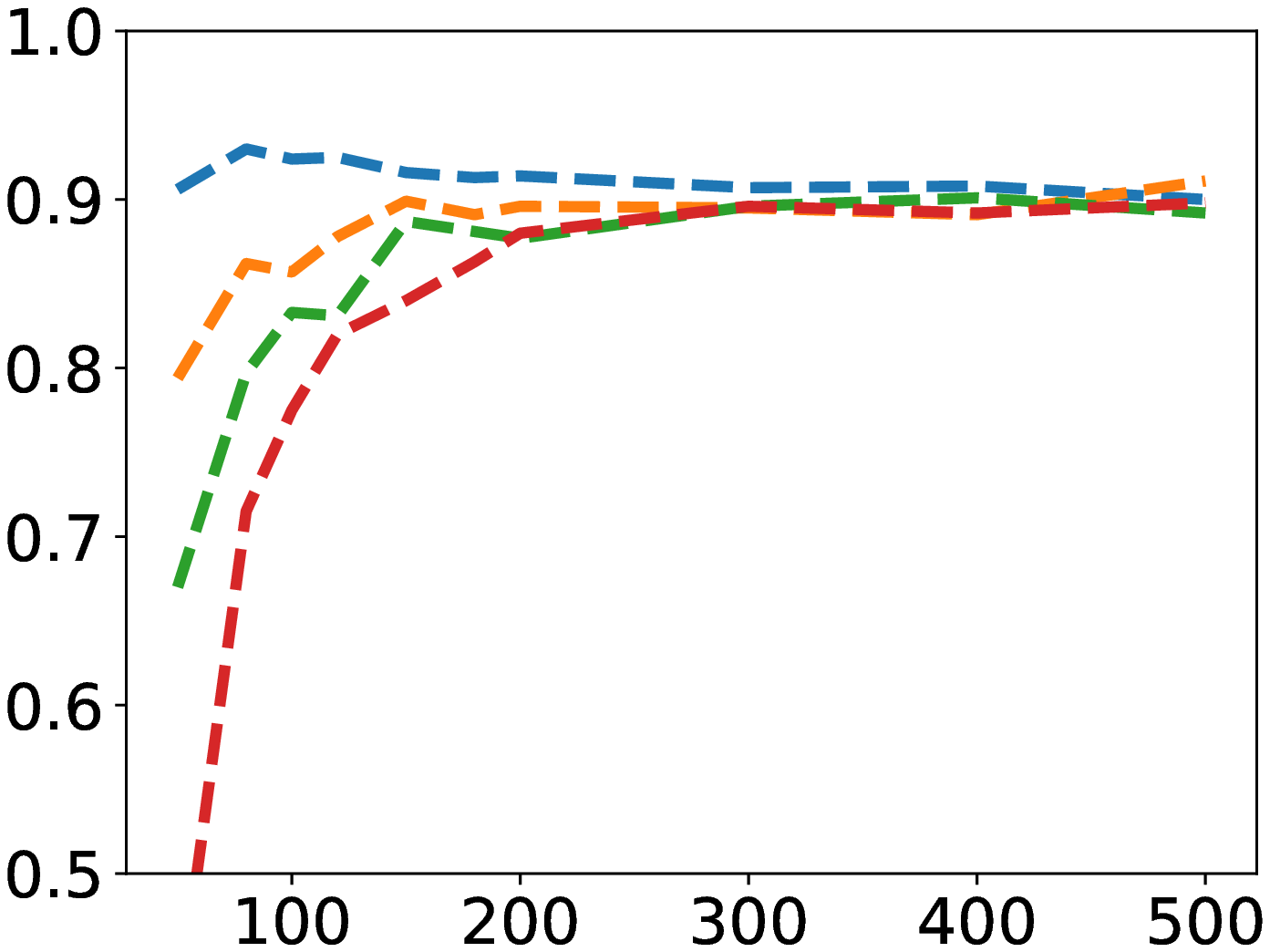} 
	\put(25,80){ \ul{\ \ \  \ $\gamma=0.7$ \ \ \ \    }}
	\put(-20,-5){\rotatebox{90}{ {\small \ \ \ log transformation  \ \ }}}
	\end{overpic}
	~
	\DeclareGraphicsExtensions{.png}
	\begin{overpic}[width=0.29\textwidth]{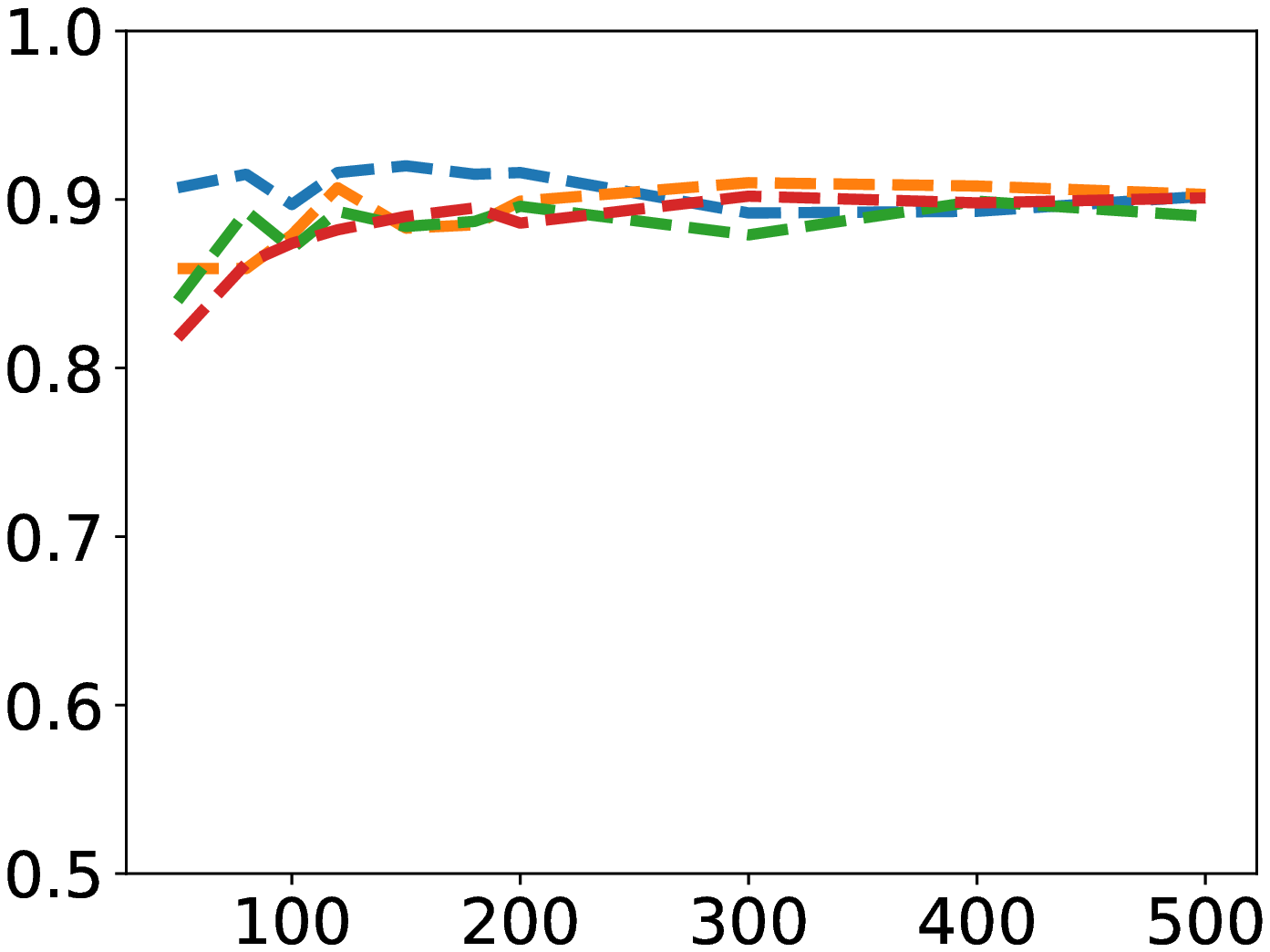} 
	\put(25,80){ \ul{\ \ \  \ $\gamma=1.0$ \ \ \ \    }}
	\end{overpic}
	~	
	\begin{overpic}[width=0.29\textwidth]{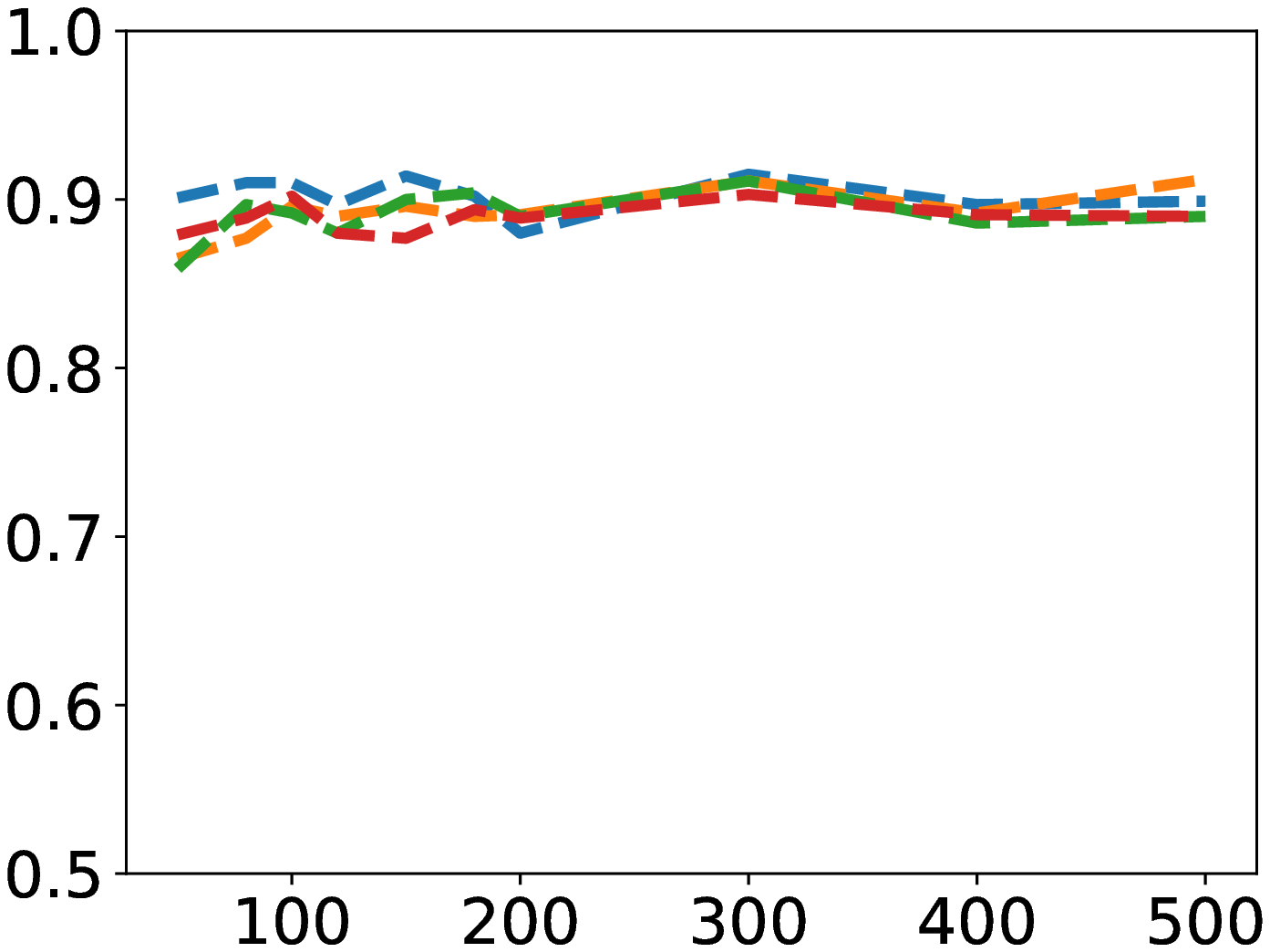} 
	\put(25,80){ \ul{\ \ \  \ $\gamma=1.3$ \ \ \ \    }}
	\end{overpic}	
%
%
\end{figure}

\vspace{-0.5cm}

\begin{figure}[h]	
	\quad\quad\quad 
	\begin{overpic}[width=0.29\textwidth]{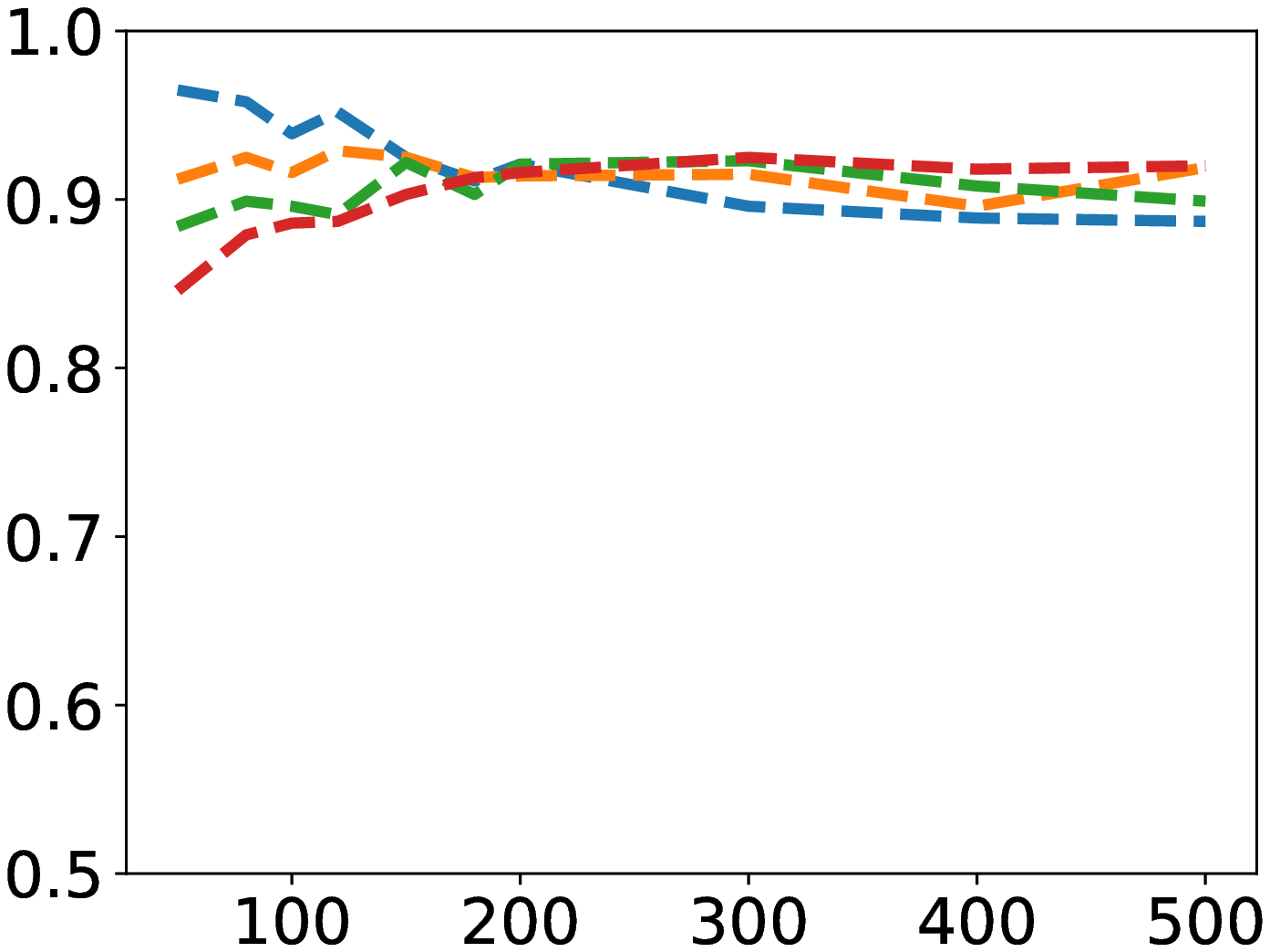} 
	\put(-20,-1){\rotatebox{90}{ {\normalsize \ \ \ standardization \  \ \ }}}
	\end{overpic}
	~
	\DeclareGraphicsExtensions{.png}
	\begin{overpic}[width=0.28\textwidth]{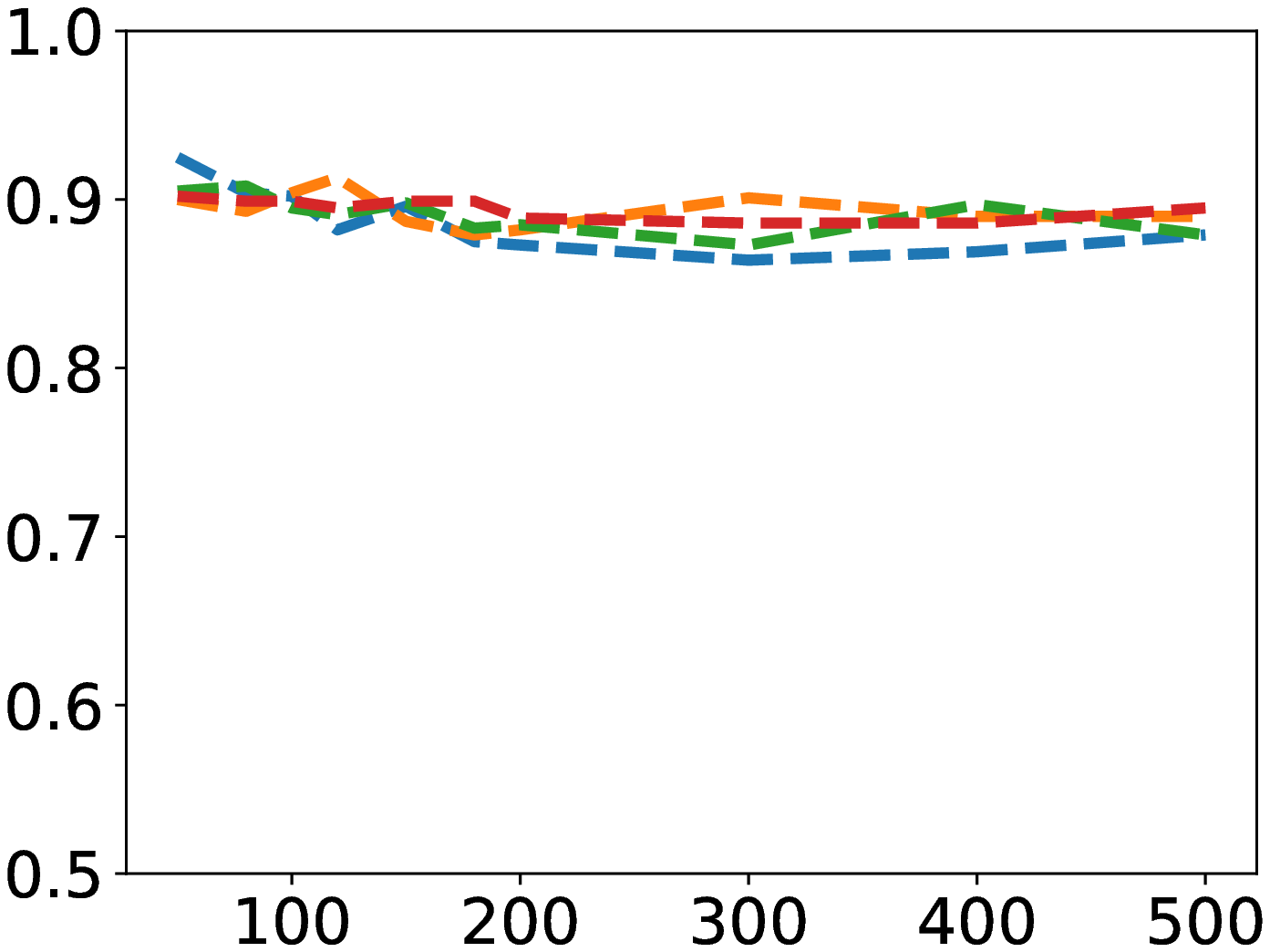} 
	\end{overpic}
	~	
	\begin{overpic}[width=0.28\textwidth]{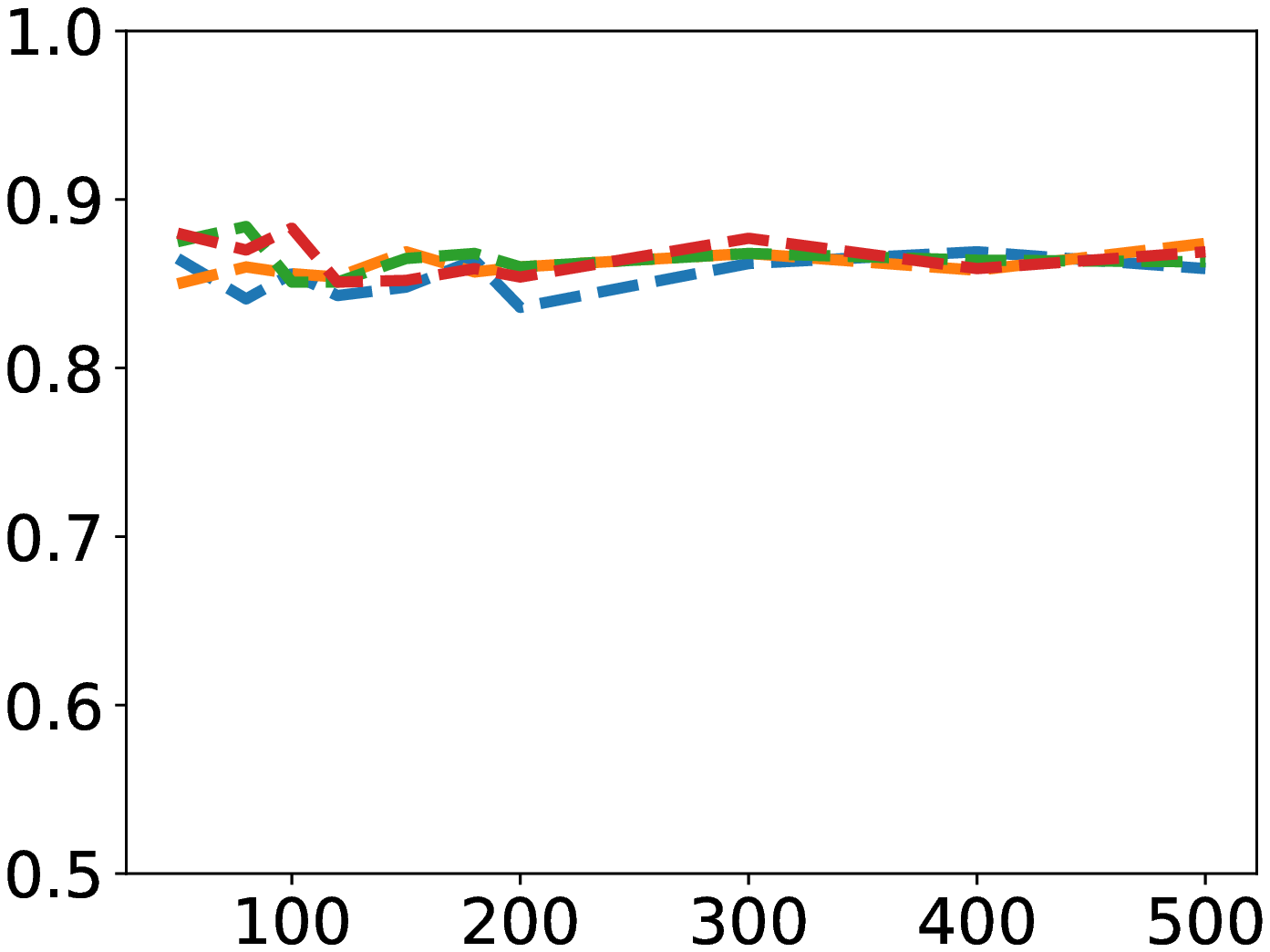} 
	\end{overpic}	
%
%
\end{figure}

\vspace{-0.5cm}

\begin{figure}[H]	
	\quad\quad\quad 
	\begin{overpic}[width=0.29\textwidth]{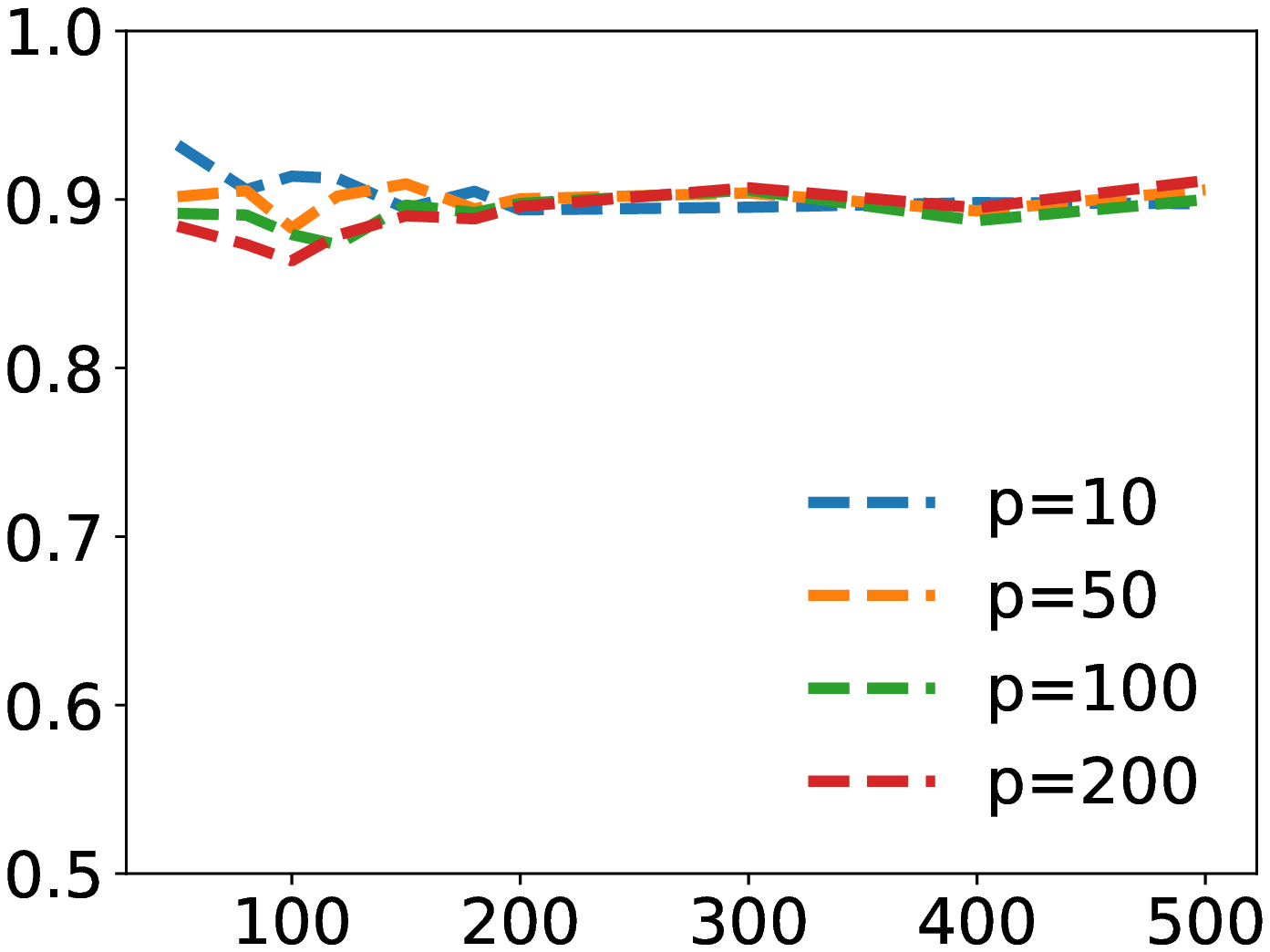} 
	\put(-21,1){\rotatebox{90}{\ $\sqrt{ \ \ }$}}
	\put(-20,-3){\rotatebox{90}{  { \ \ \ \ \ \ \small transformation \ \ } }}
	\end{overpic}
	~
	\DeclareGraphicsExtensions{.png}
	\begin{overpic}[width=0.29\textwidth]{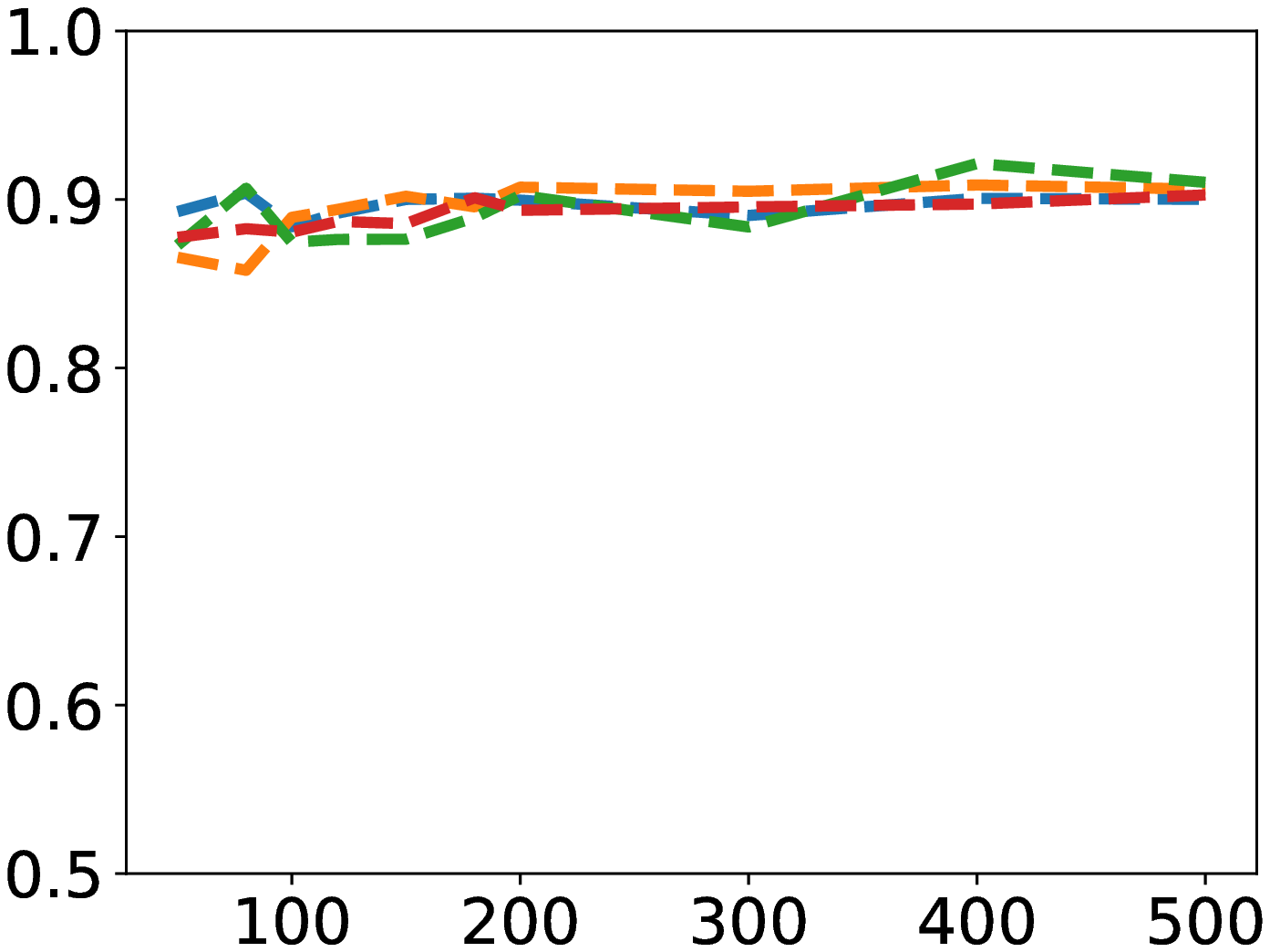} 
	\end{overpic}
	~	
	\begin{overpic}[width=0.29\textwidth]{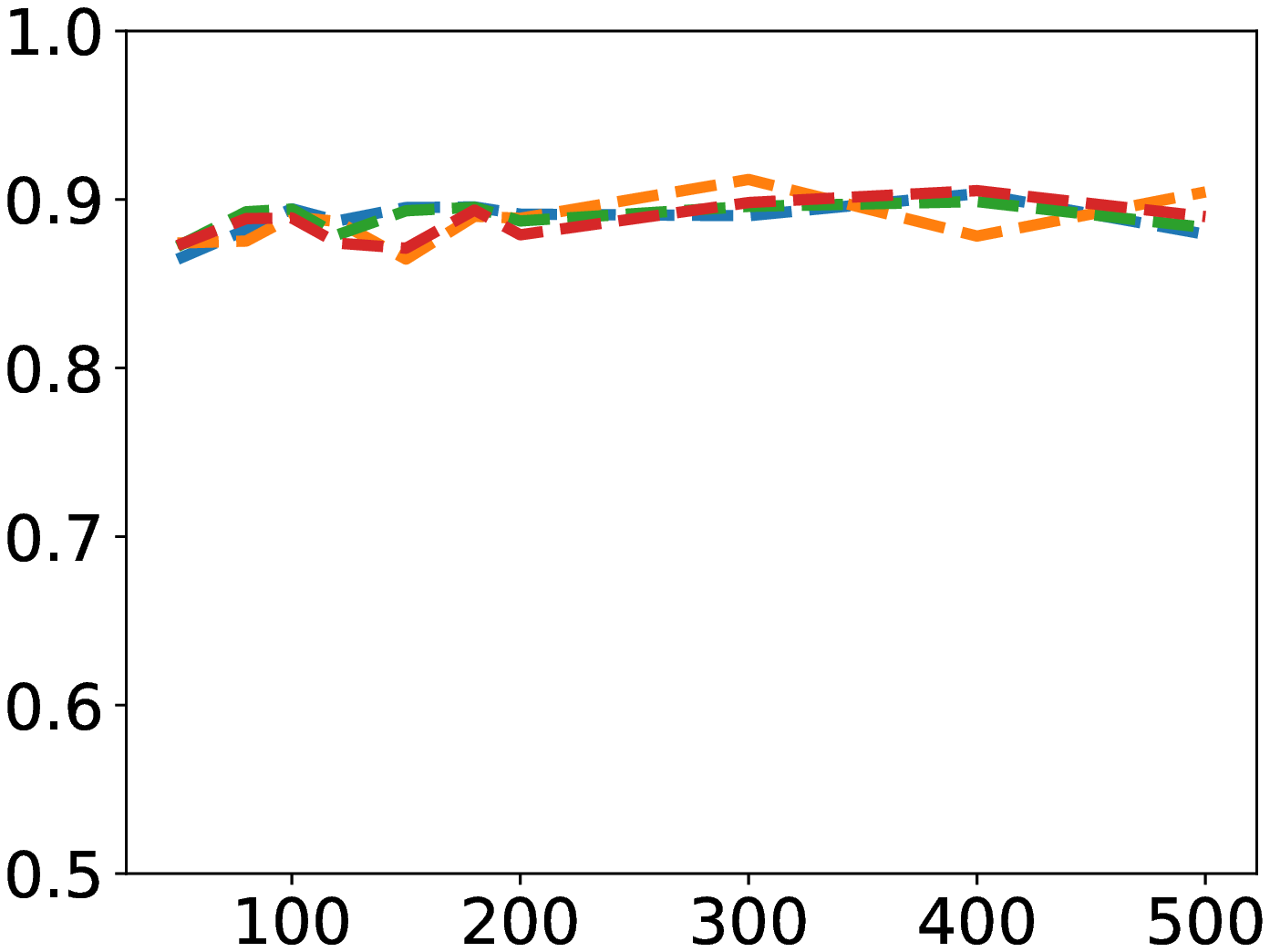} 
	\end{overpic}	
	\vspace{+0.2cm}	
	\caption{(Simultaneous coverage probability versus $n$ in simulation model (ii) with a polynomial decay profile). The plotting scheme is the same as described in the caption of Figure~\ref{SUPP:fig1} above.} 
	\label{SUPP:fig2}
\end{figure}

\newpage

\begin{figure}[H]	
\vspace{0.5cm}
	\quad\quad\quad 
	\begin{overpic}[width=0.29\textwidth]{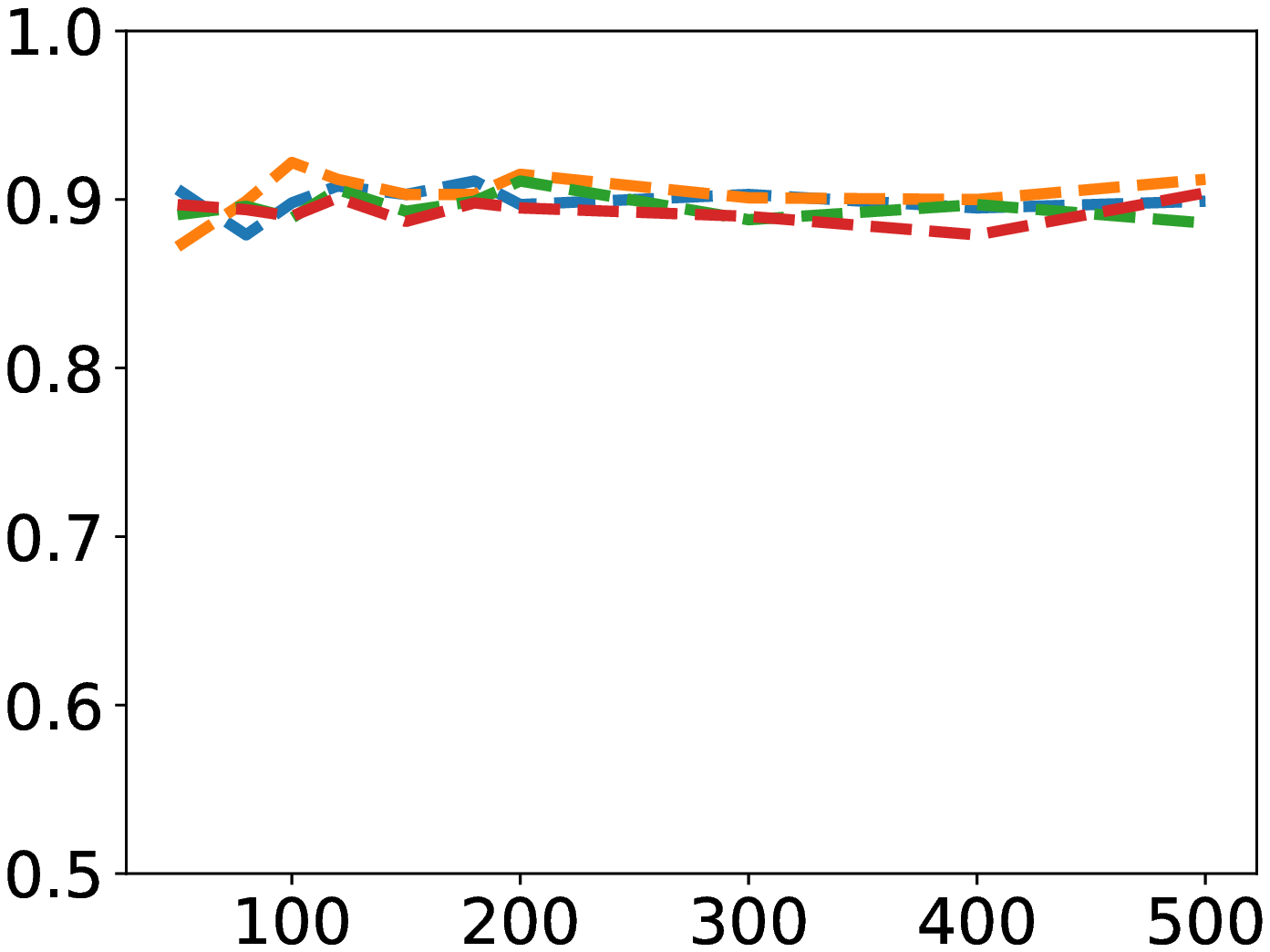} 
    \put(25,80){ \ul{\ \ \  \ $\delta=0.7$ \ \ \ \    }}
	\put(-20,-5){\rotatebox{90}{ {\small \ \ \ log transformation  \ \ }}}
\end{overpic}
	~
	\DeclareGraphicsExtensions{.png}
	\begin{overpic}[width=0.29\textwidth]{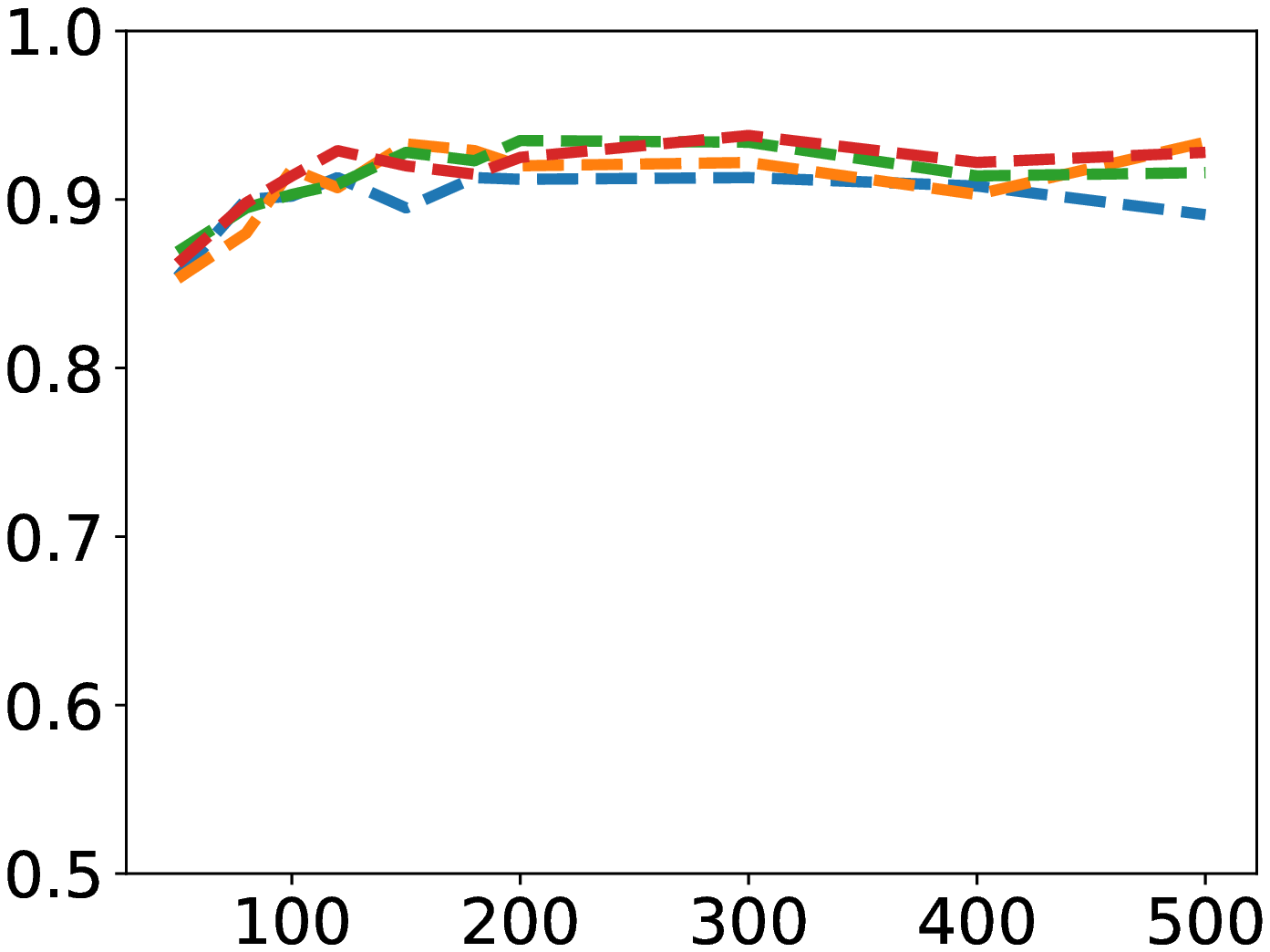} 
	\put(25,80){ \ul{\ \ \  \ $\delta=0.8$ \ \ \ \    }}
	\end{overpic}
	~	
	\begin{overpic}[width=0.29\textwidth]{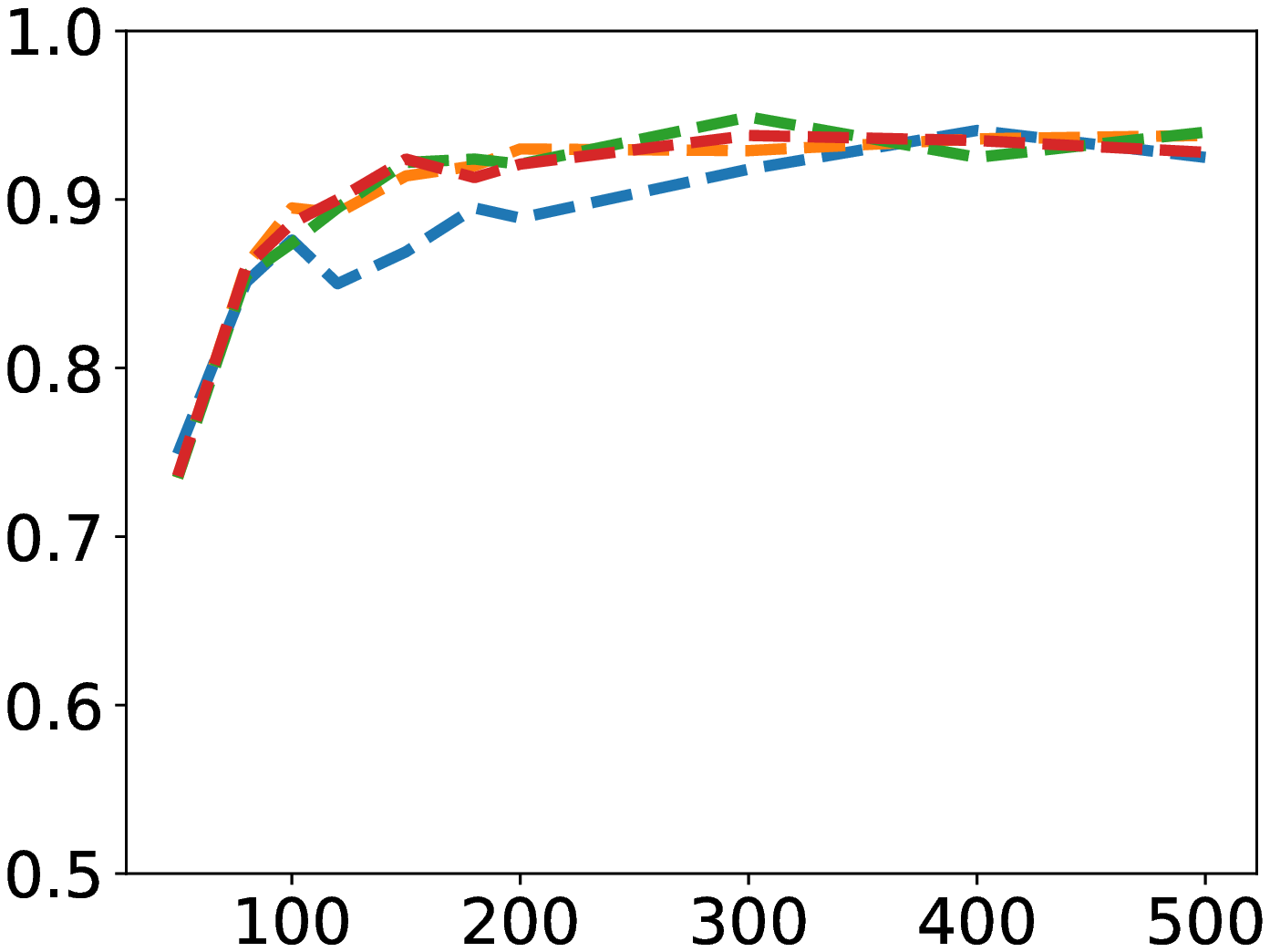} 

	\put(25,80){ \ul{\ \ \  \ $\delta=0.9$ \ \ \ \    }}
 		 				
	\end{overpic}	
\end{figure}

\vspace{-0.5cm}

\begin{figure}[H]	
	\quad\quad\quad 
	\begin{overpic}[width=0.29\textwidth]{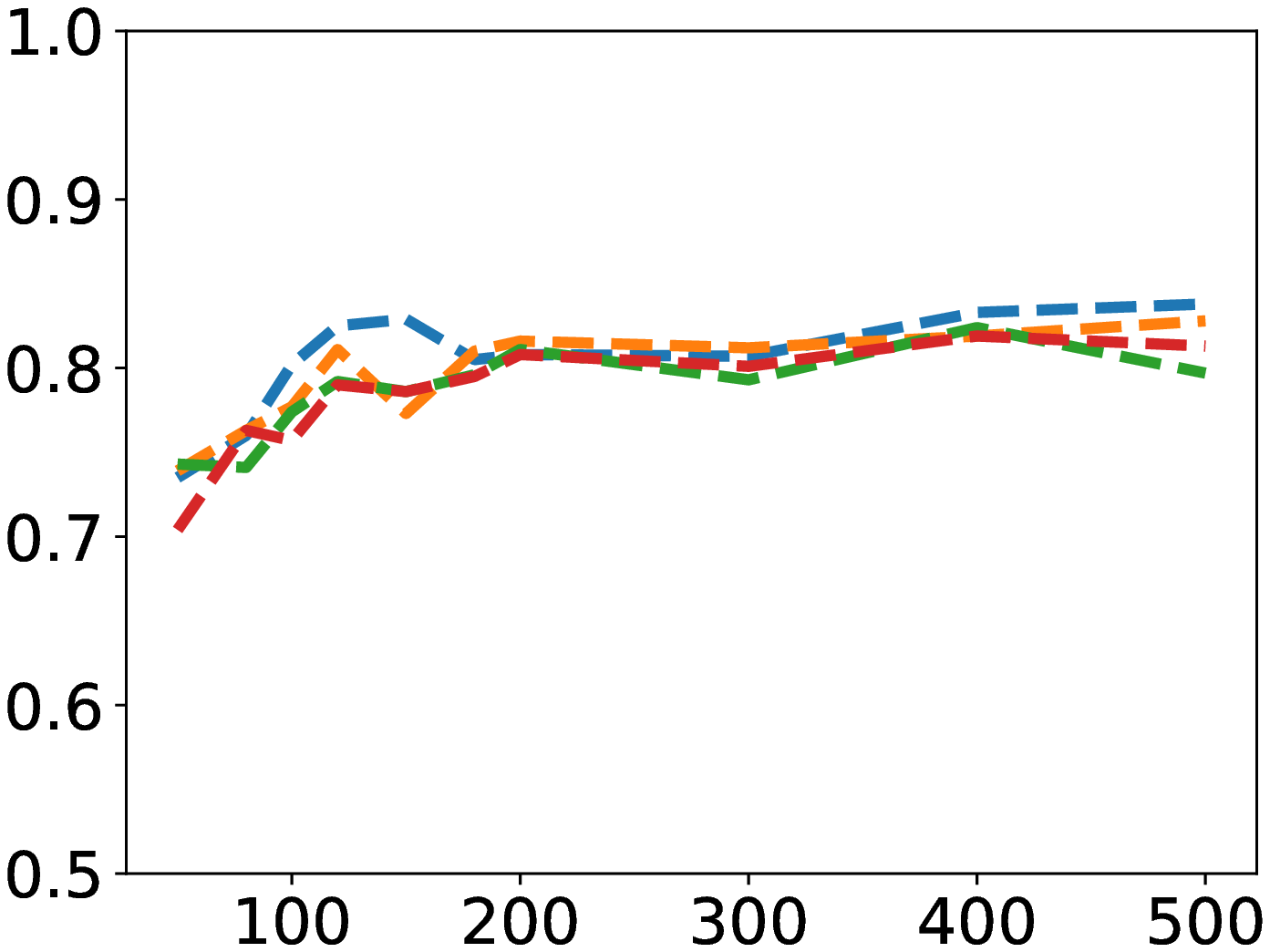} 
	\put(-20,-1){\rotatebox{90}{ {\small \ \ \ standardization \  \ \ }}}
	\end{overpic}
	~
	\DeclareGraphicsExtensions{.png}
	\begin{overpic}[width=0.29\textwidth]{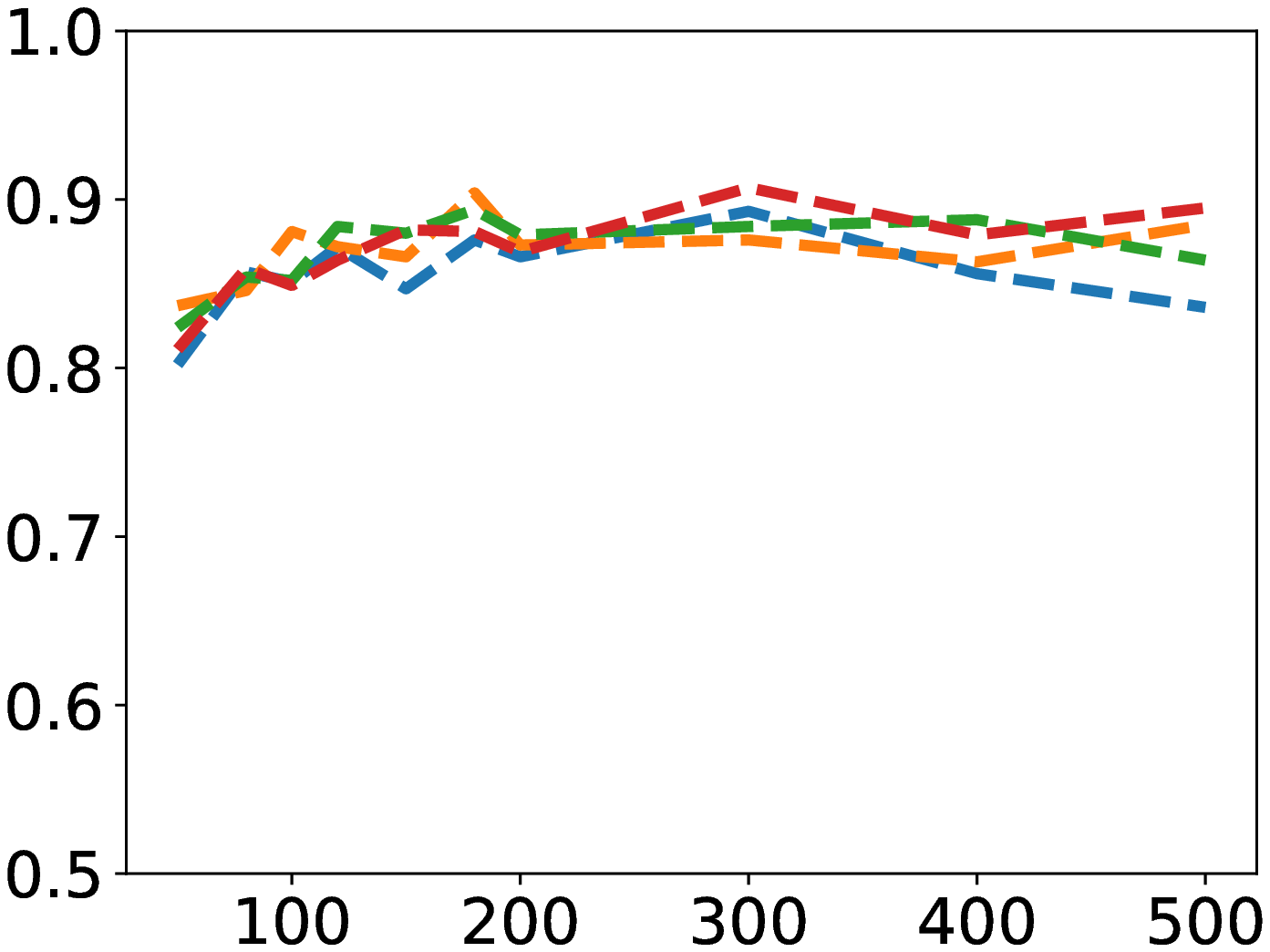} 
	\end{overpic}
	~	
	\begin{overpic}[width=0.29\textwidth]{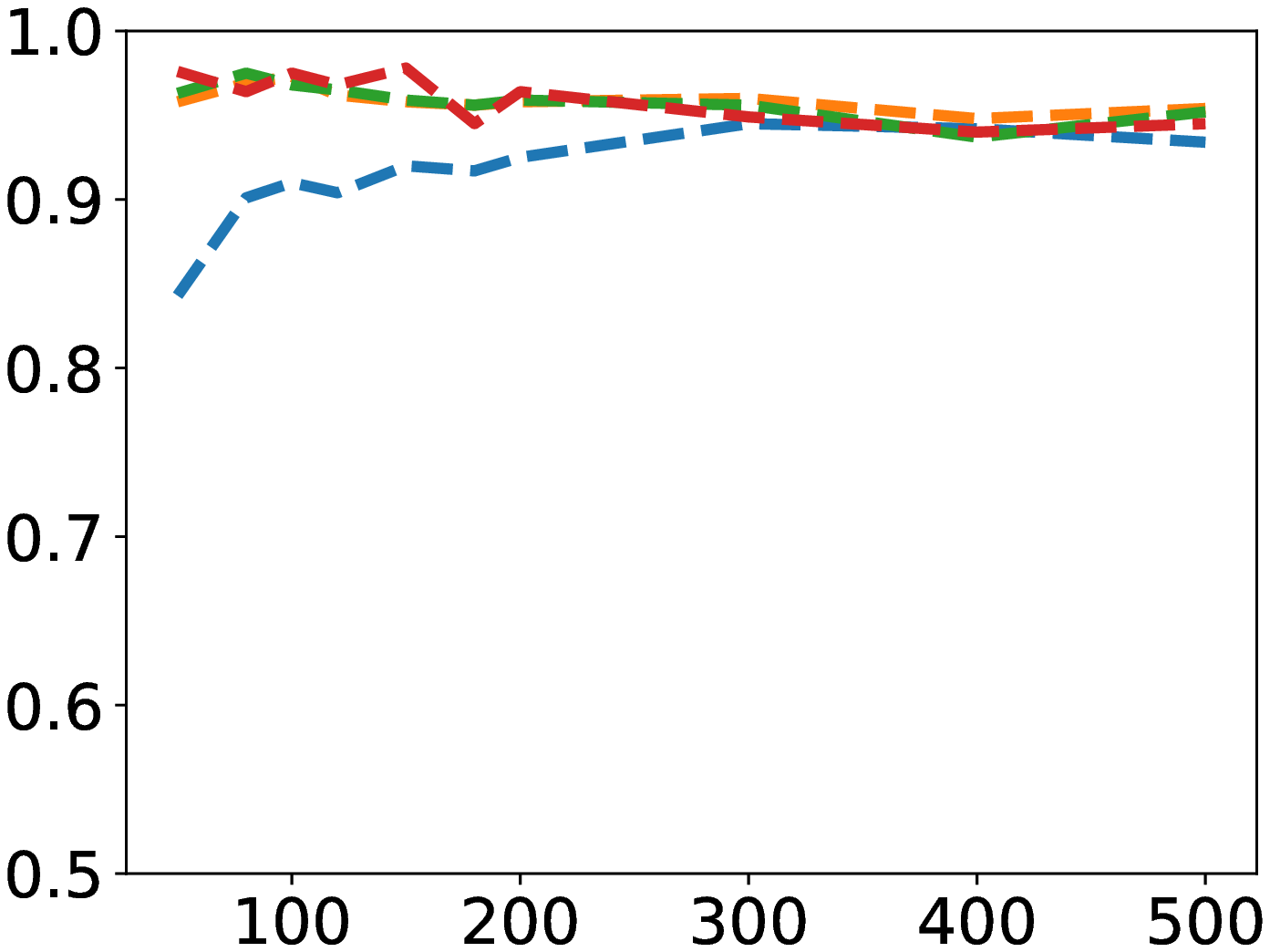} 
		 				
	\end{overpic}	
\end{figure}

\vspace{-0.5cm}

\begin{figure}[H]	
	\quad\quad\quad 
	\begin{overpic}[width=0.29\textwidth]{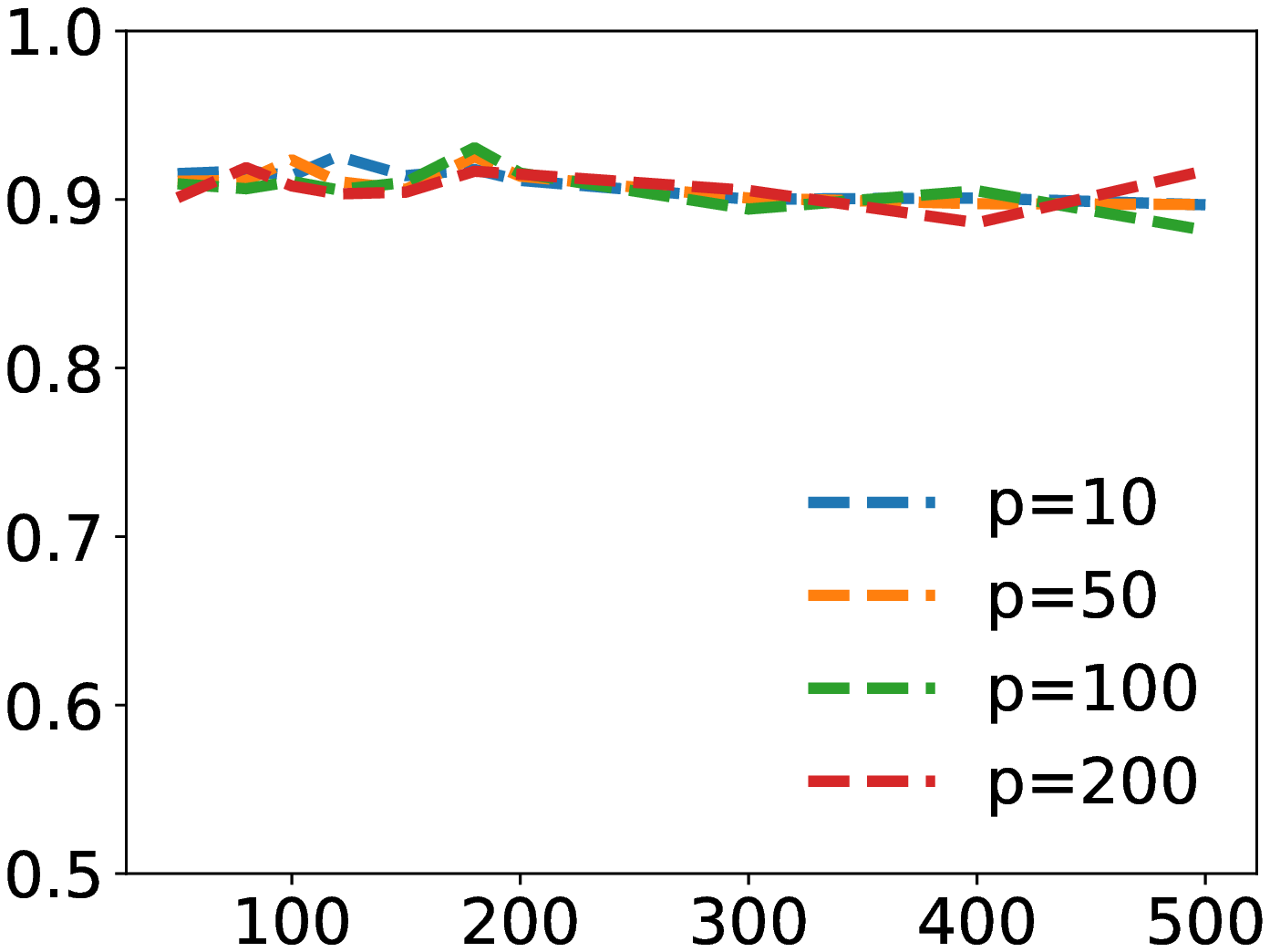} 
	\put(-21,1){\rotatebox{90}{\ $\sqrt{ \ \ }$}}
	\put(-20,-3){\rotatebox{90}{  { \ \ \ \ \ \ \small transformation \ \ } }}

	\end{overpic}
	~
	\DeclareGraphicsExtensions{.png}
	\begin{overpic}[width=0.29\textwidth]{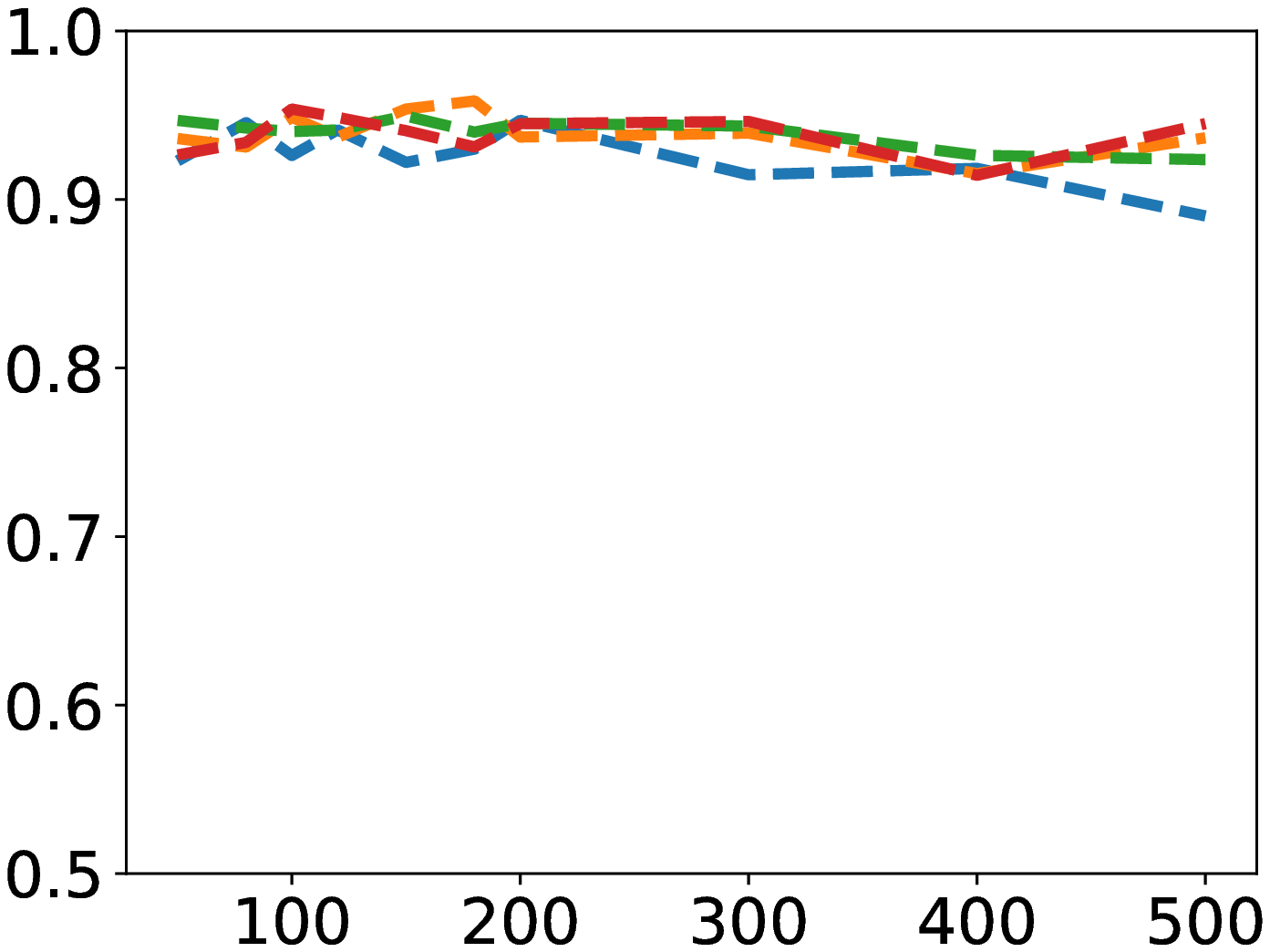} 
	\end{overpic}
	~	
	\begin{overpic}[width=0.29\textwidth]{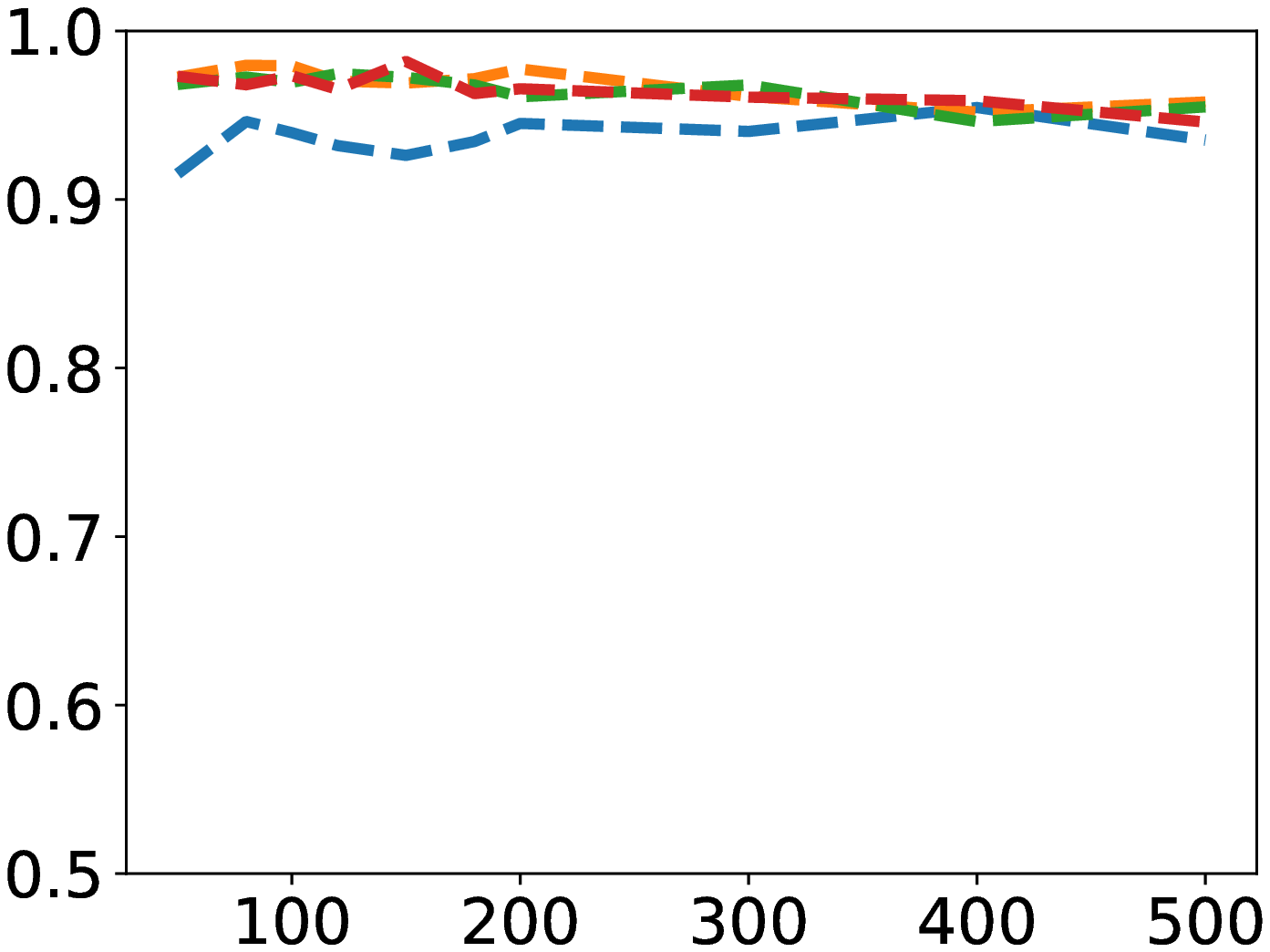} 
		 				
	\end{overpic}	
	\vspace{+0.2cm}	
	\caption{(Simultaneous coverage probability versus $n$ in simulation model (i) with an exponential decay profile). The plotting scheme is the same as described in the caption of Figure~\ref{SUPP:fig1} above, except that the three columns correspond to values of the eigenvalue decay parameter $\delta$.}
	\label{SUPP:fig3}
\end{figure}

\newpage
\begin{figure}[H]	
\vspace{0.5cm}
	\quad\quad\quad 
	\begin{overpic}[width=0.29\textwidth]{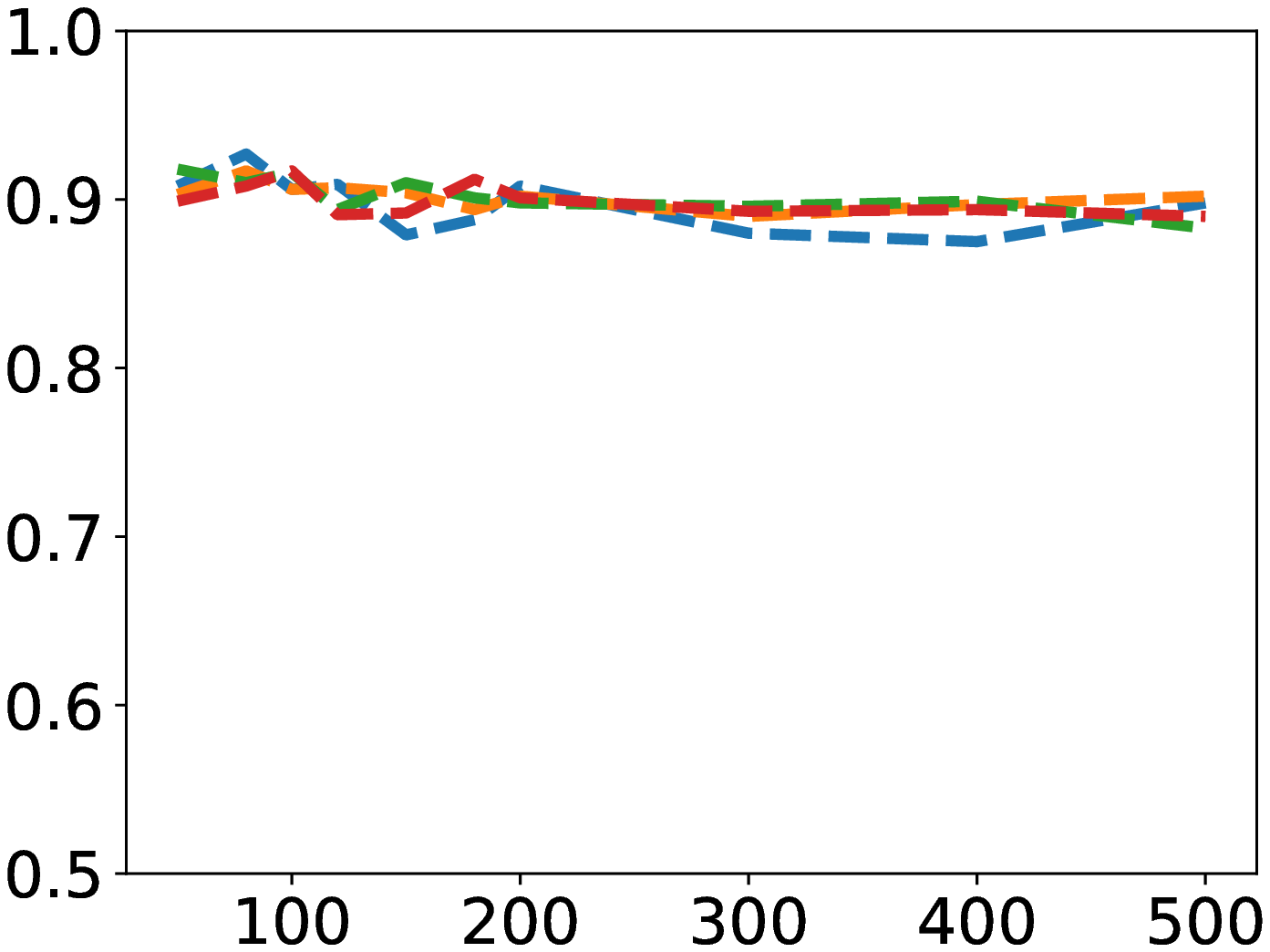} 
    \put(25,80){ \ul{\ \ \  \ $\delta=0.7$ \ \ \ \    }}
	\put(-20,-5){\rotatebox{90}{ {\small \ \ \ log transformation  \ \ }}}
\end{overpic}
	~
	\DeclareGraphicsExtensions{.png}
	\begin{overpic}[width=0.29\textwidth]{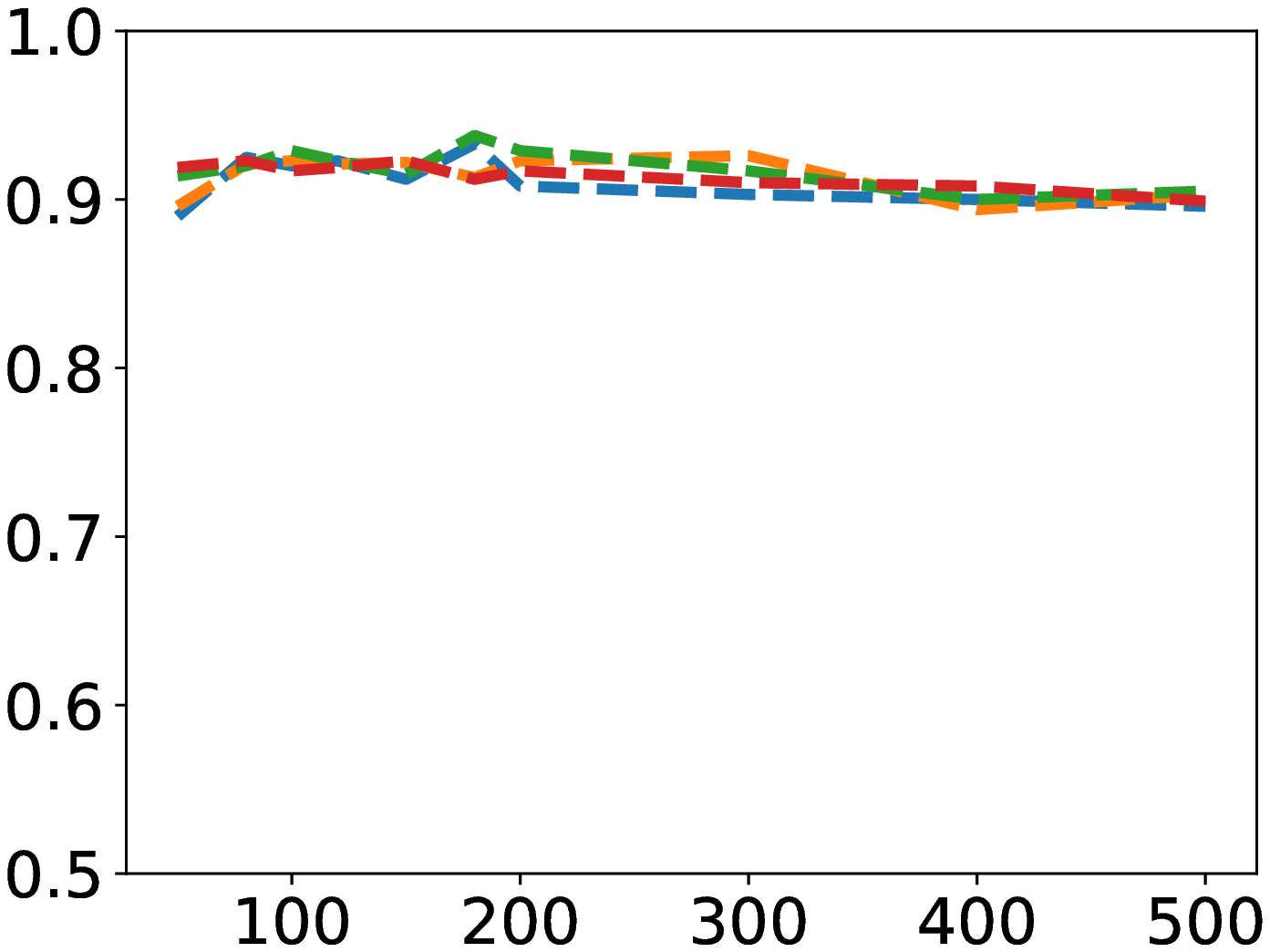} 
	\put(25,80){ \ul{\ \ \  \ $\delta=0.8$ \ \ \ \    }}
	\end{overpic}
	~	
	\begin{overpic}[width=0.29\textwidth]{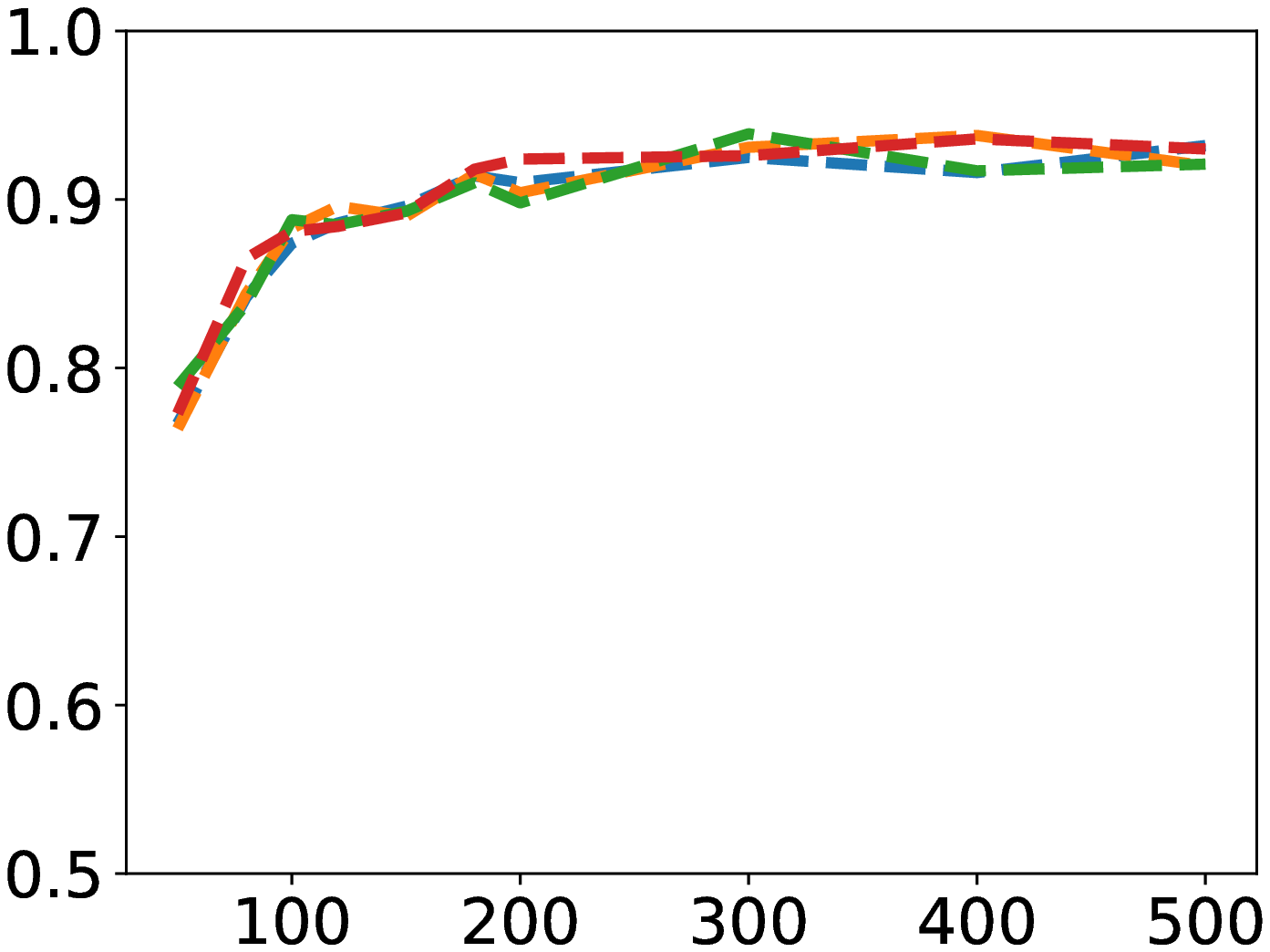} 

	\put(25,80){ \ul{\ \ \  \ $\delta=0.9$ \ \ \ \    }}
 		 				
	\end{overpic}	
\end{figure}

\vspace{-0.5cm}

\begin{figure}[H]	
	\quad\quad\quad 
	\begin{overpic}[width=0.29\textwidth]{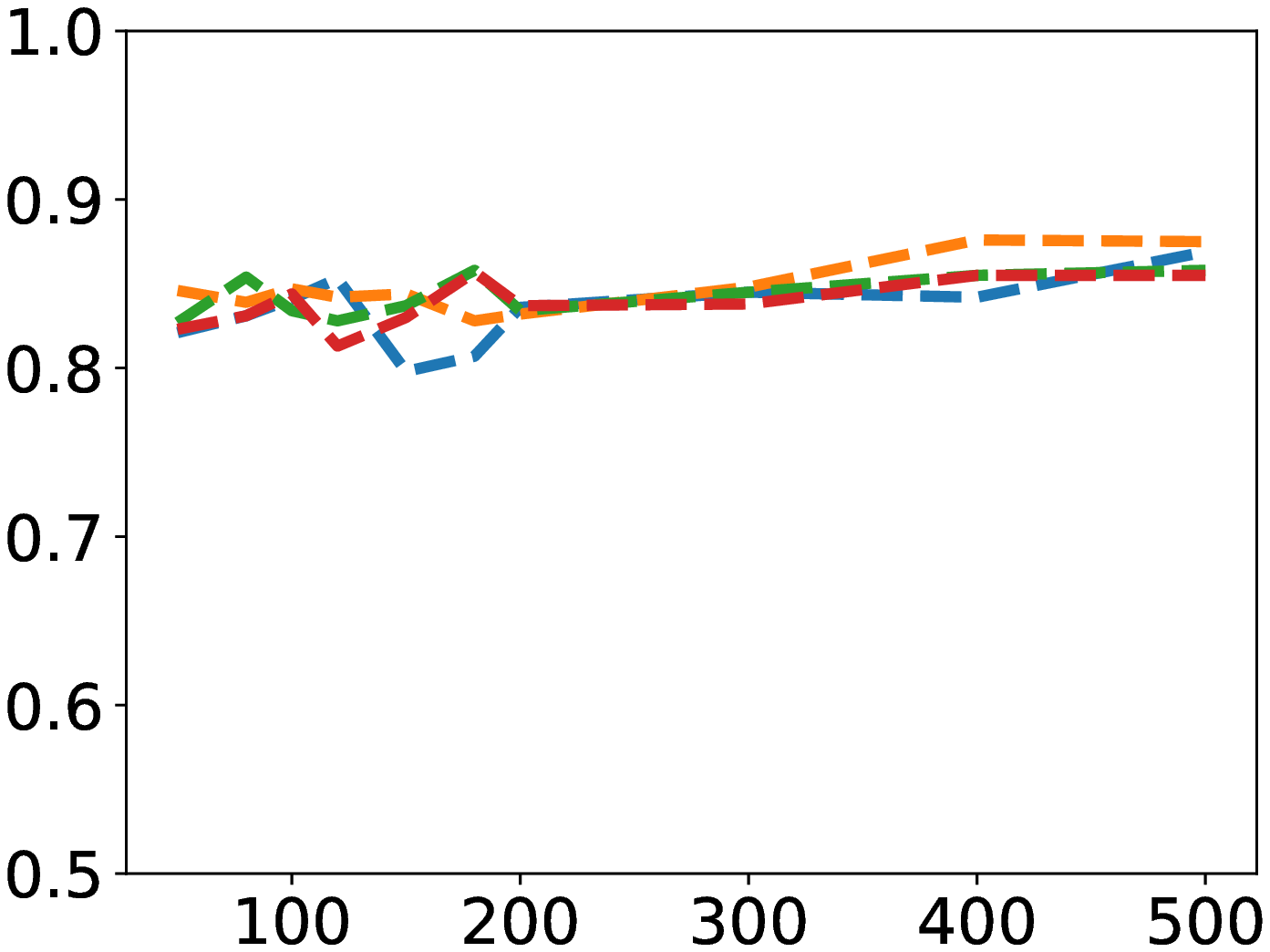} 
	\put(-20,-1){\rotatebox{90}{ {\small \ \ \ standardization \  \ \ }}}
	\end{overpic}
	~
	\DeclareGraphicsExtensions{.png}
	\begin{overpic}[width=0.29\textwidth]{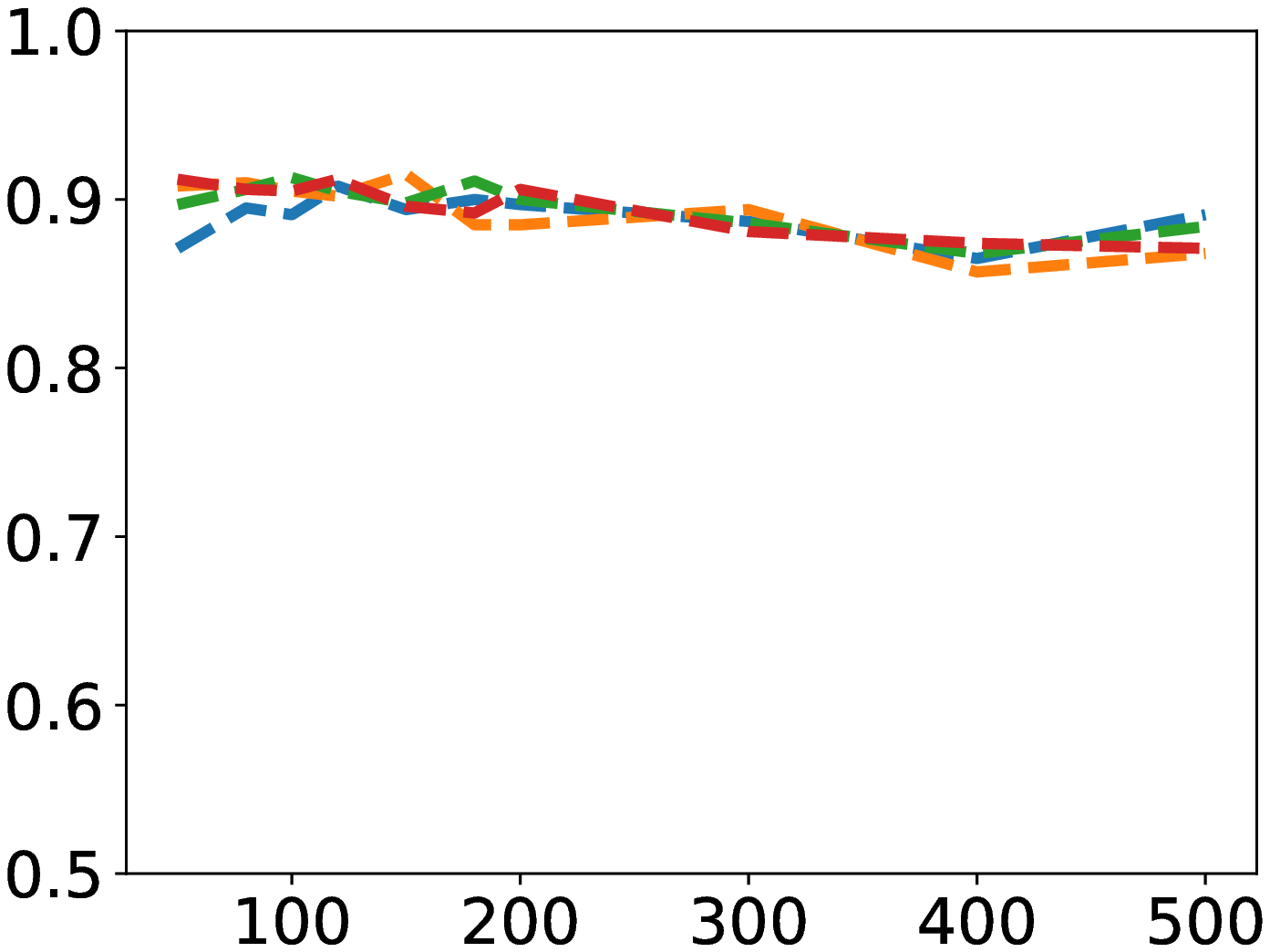} 
	\end{overpic}
	~	
	\begin{overpic}[width=0.29\textwidth]{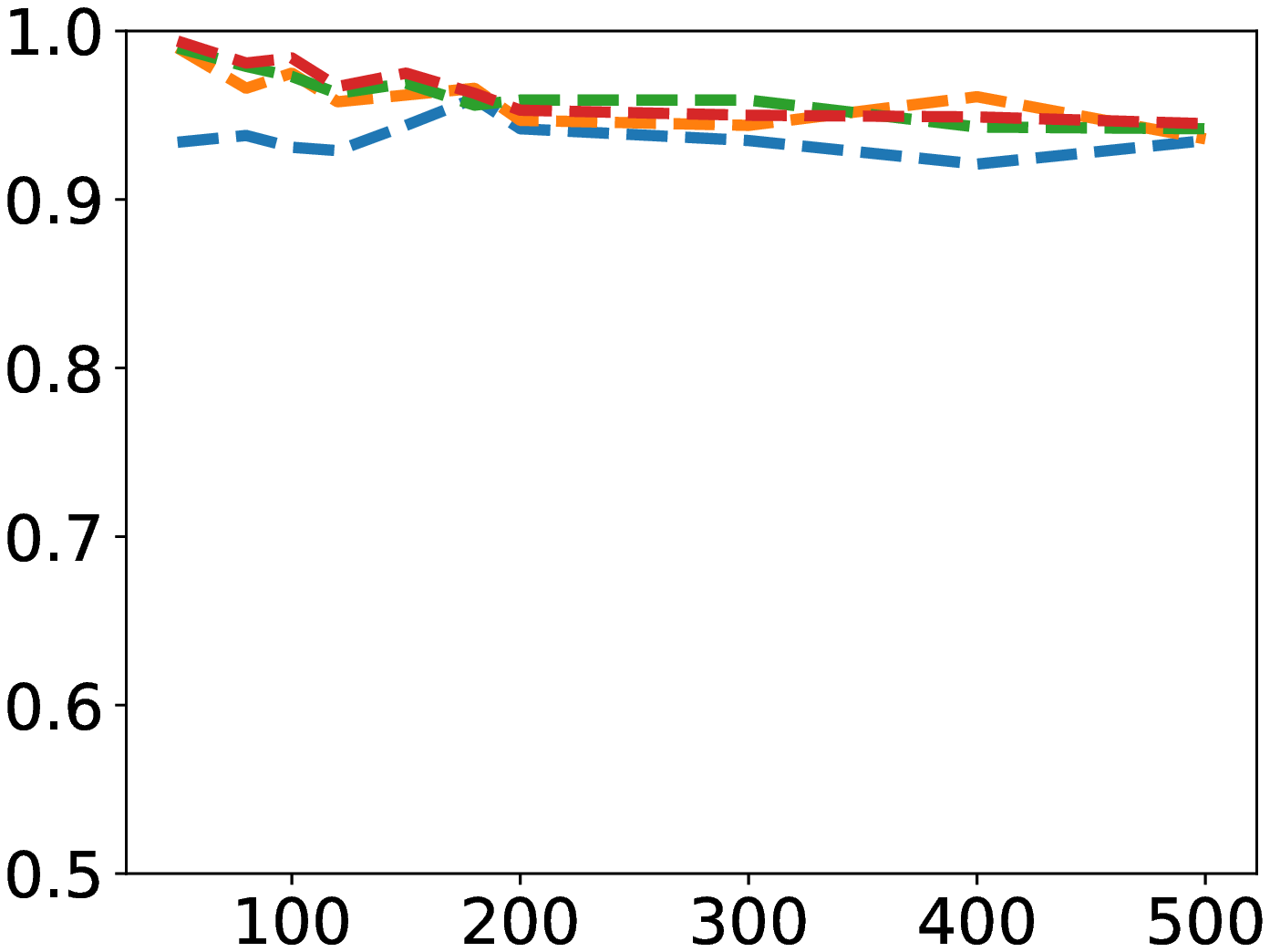} 
	 				
	\end{overpic}	
\end{figure}

\vspace{-0.5cm}

\begin{figure}[H]	
	\quad\quad\quad 
	\begin{overpic}[width=0.29\textwidth]{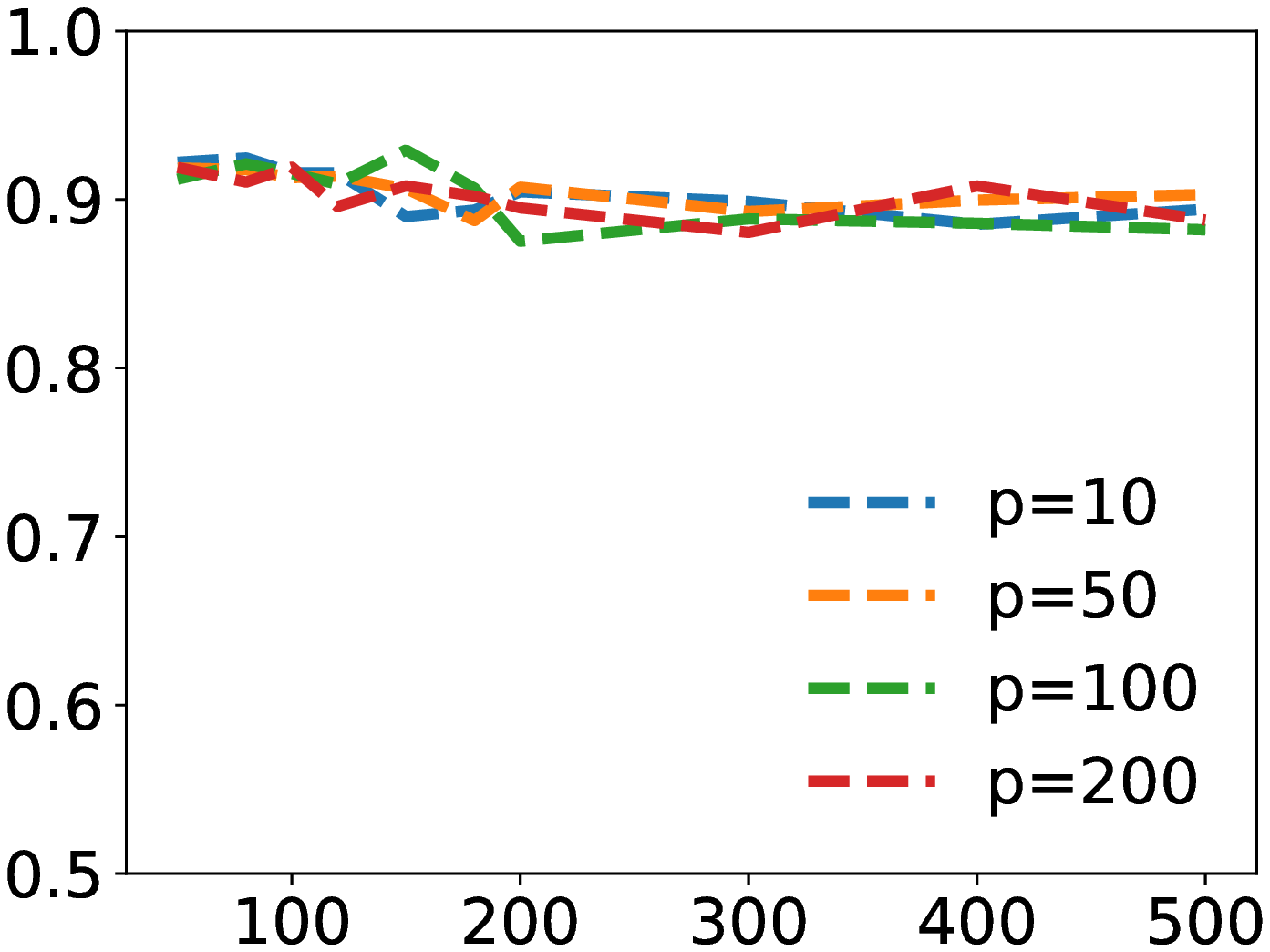} 
		\put(-21,1){\rotatebox{90}{\ $\sqrt{ \ \ }$}}
	    \put(-20,-3){\rotatebox{90}{  { \ \ \ \ \ \ \small transformation \ \ } }}

	\end{overpic}
	~
	\DeclareGraphicsExtensions{.png}
	\begin{overpic}[width=0.29\textwidth]{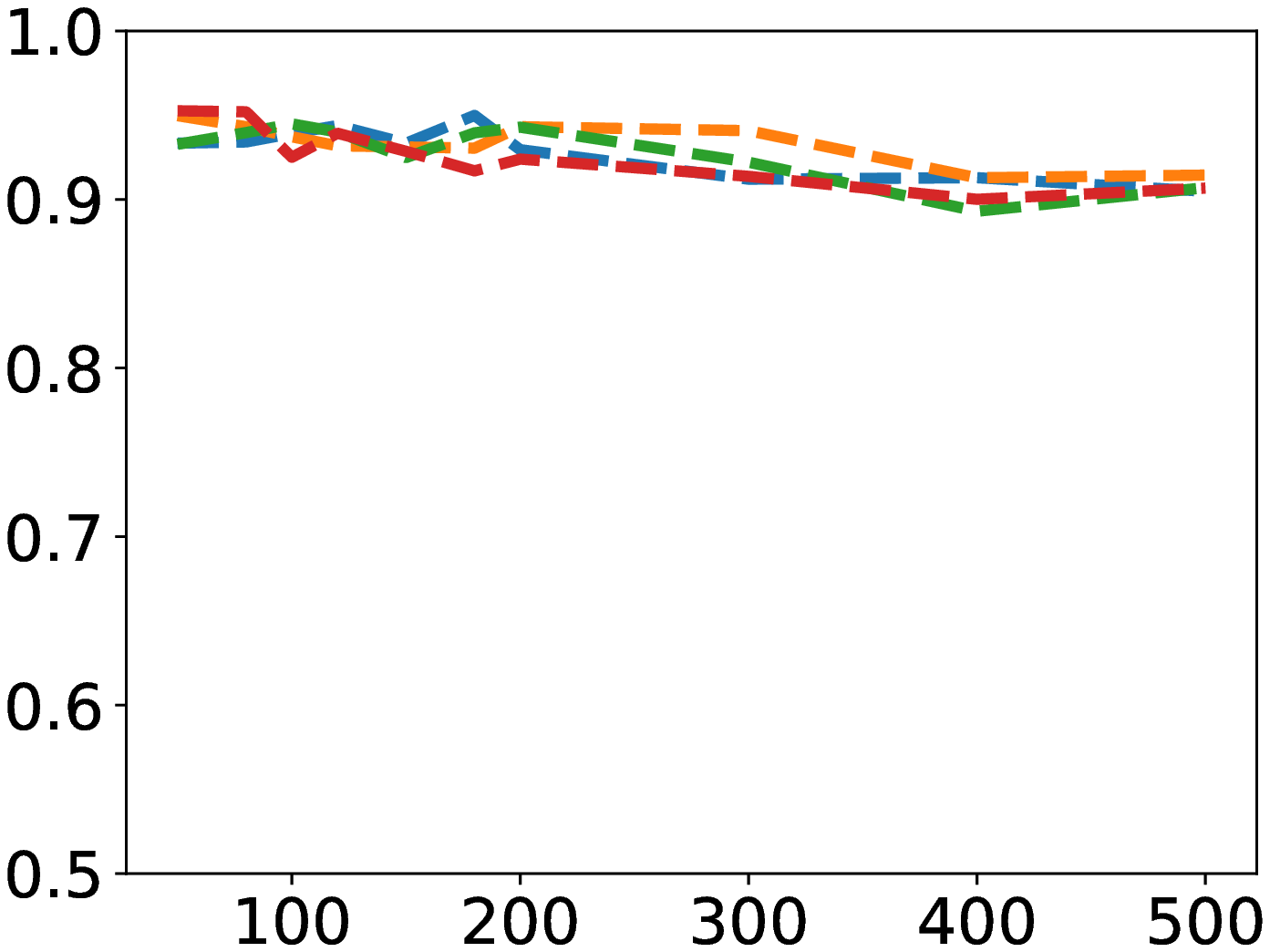} 
	\end{overpic}
	~	
	\begin{overpic}[width=0.29\textwidth]{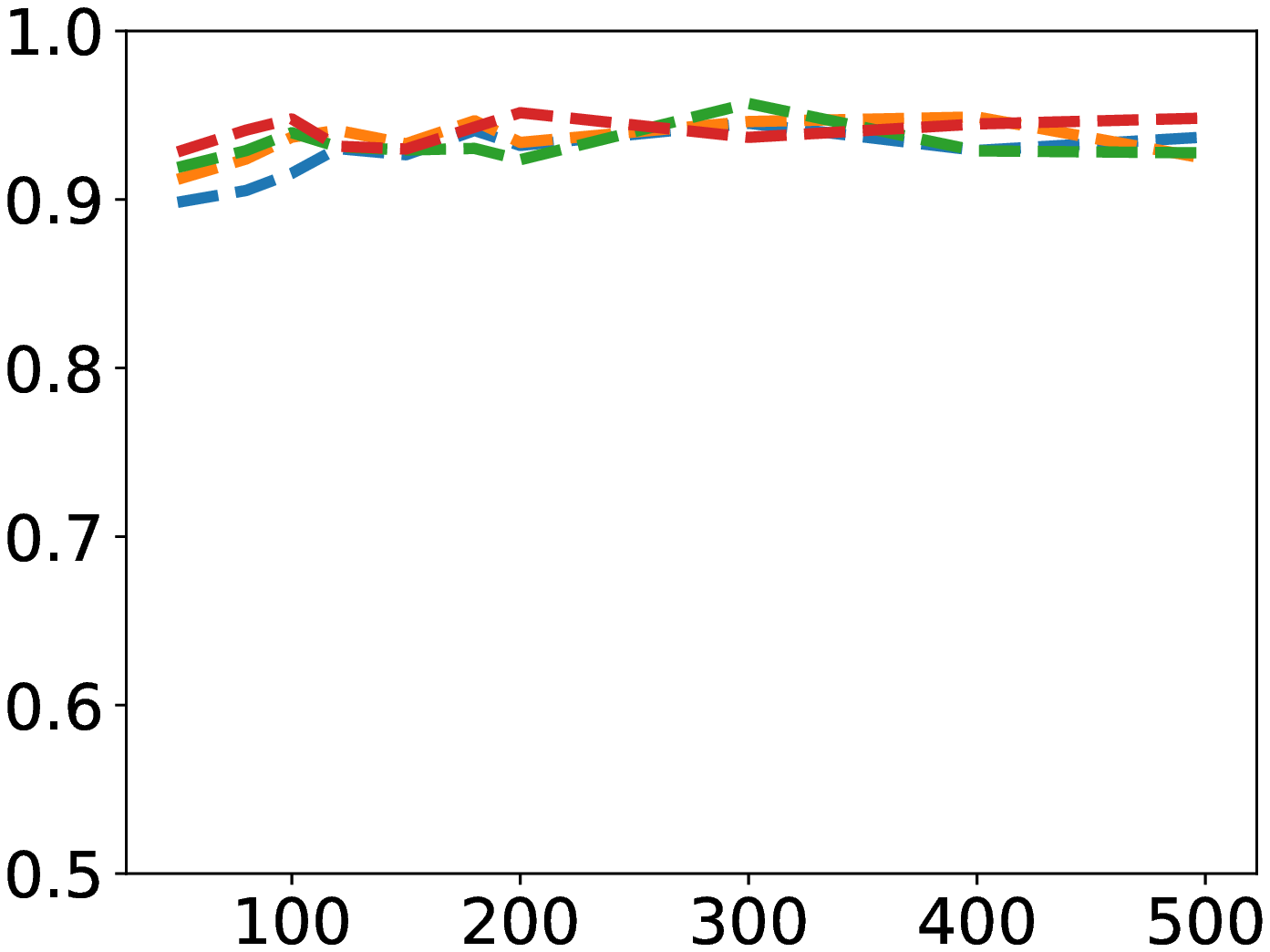} 
		 				
	\end{overpic}	
	\vspace{+0.2cm}	
	\caption{(Simultaneous coverage probability versus $n$ in simulation model (ii) with an exponential decay profile). The plotting scheme is the same as described in the caption of Figure~\ref{SUPP:fig1} above.}
	\label{SUPP:fig4}
\end{figure}

\newpage

\begin{figure}[H]	
\vspace{1cm}
	\quad\quad\quad 
	\begin{overpic}[width=0.29\textwidth]{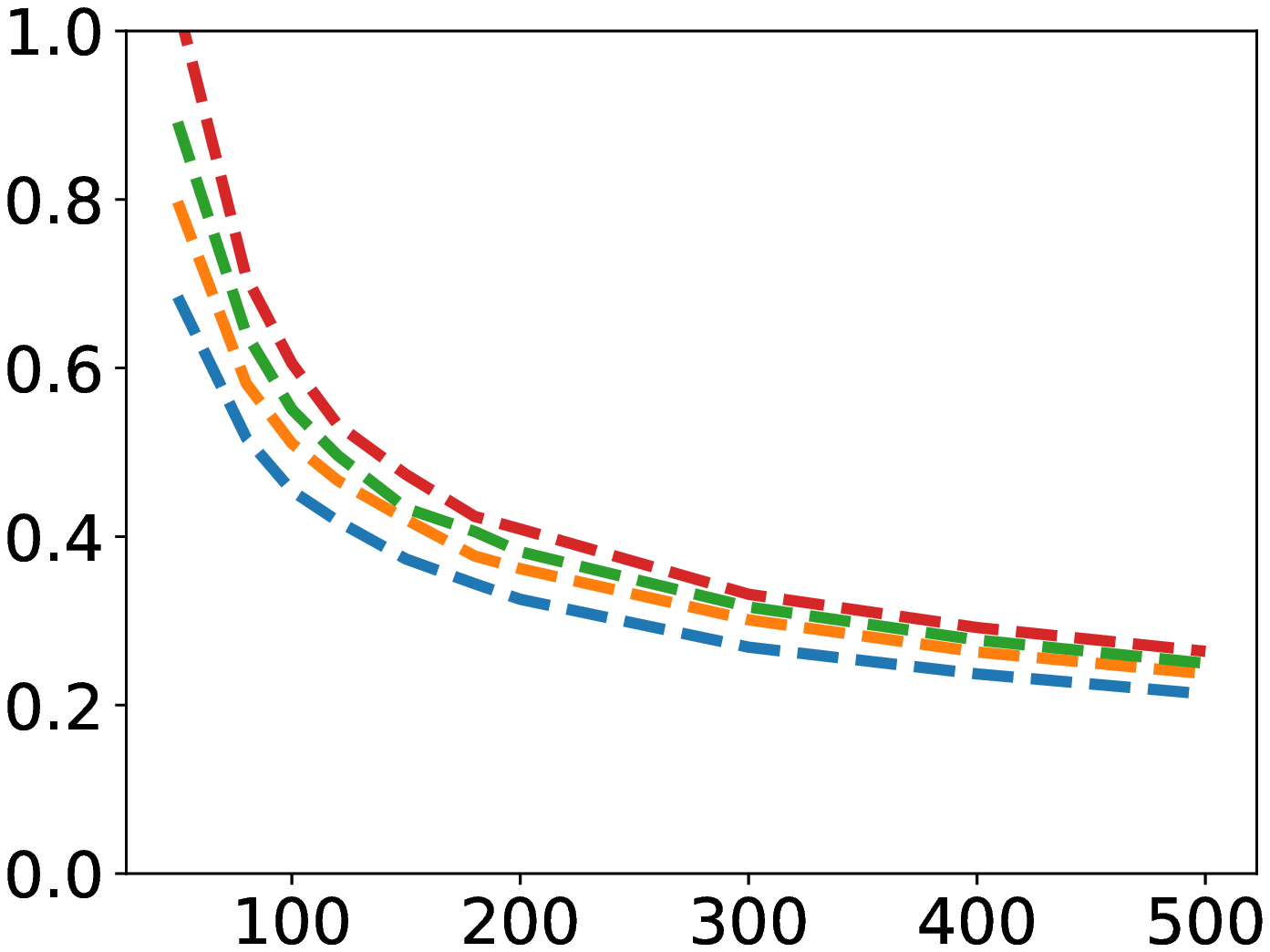} 
	\put(25,80){ \ul{\ \ \  \ $\gamma=0.7$ \ \ \ \    }}
	\put(-20,-5){\rotatebox{90}{ {\small \ \ \ log transformation  \ \ }}}
	\end{overpic}
	~
	\DeclareGraphicsExtensions{.png}
	\begin{overpic}[width=0.29\textwidth]{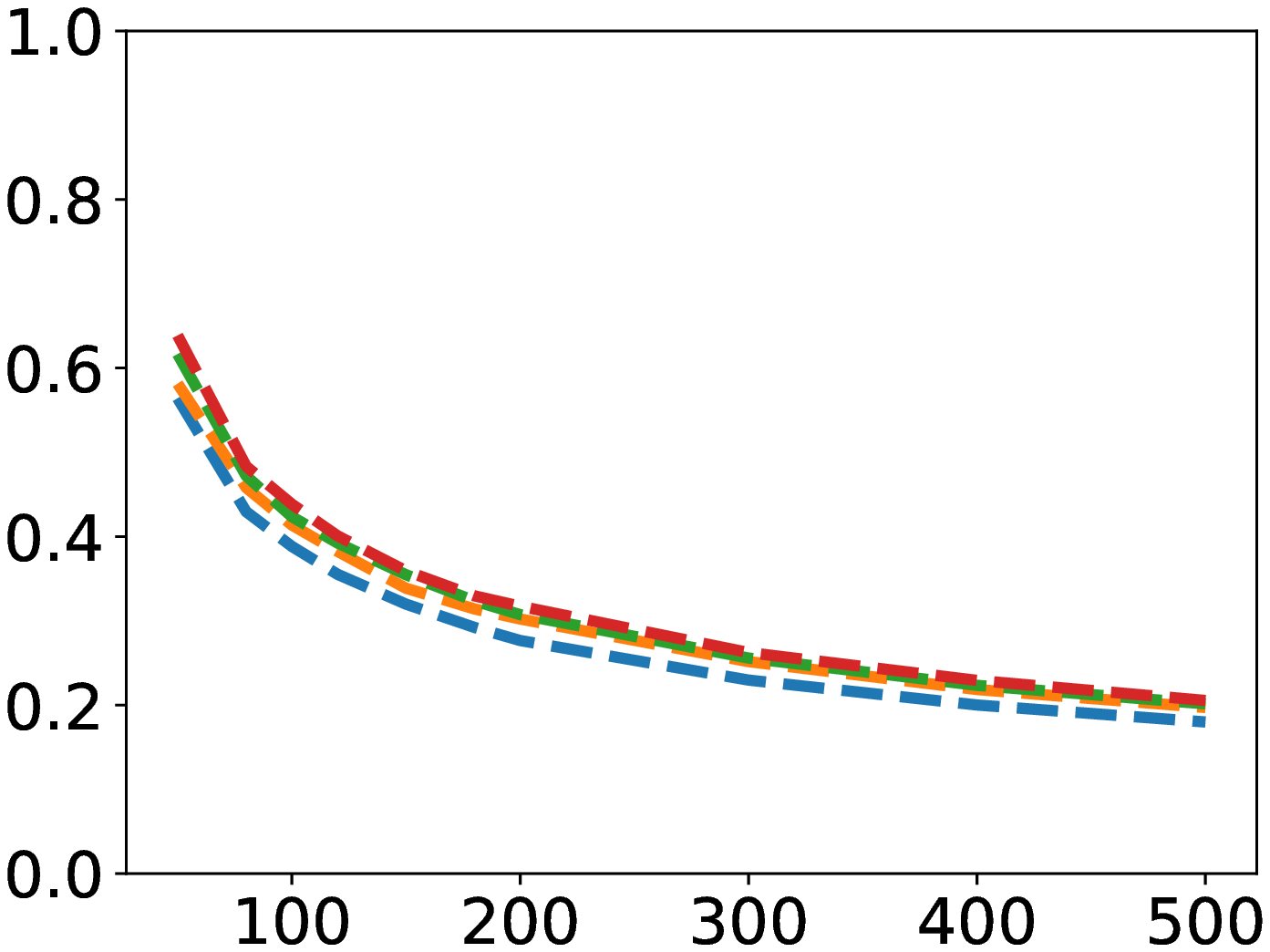} 
	\put(25,80){ \ul{\ \ \  \ $\gamma=1.0$ \ \ \ \    }}
	\end{overpic}
	~	
	\begin{overpic}[width=0.29\textwidth]{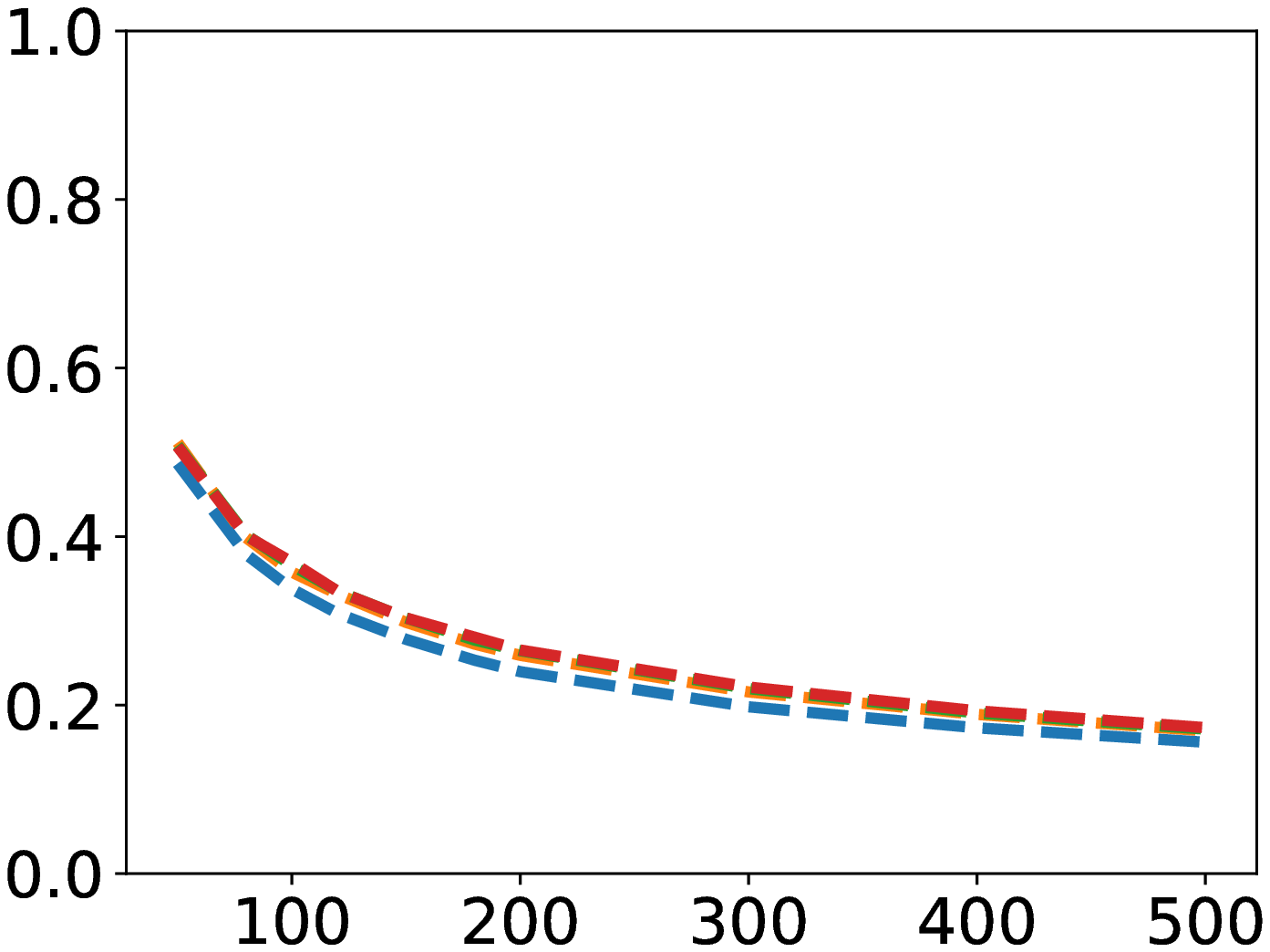} 
	\put(25,80){ \ul{\ \ \  \ $\gamma=1.3$ \ \ \ \    }}
	\end{overpic}	
%
%
\end{figure}

\vspace{-0.5cm}

\begin{figure}[H]	
	\quad\quad\quad 
	\begin{overpic}[width=0.29\textwidth]{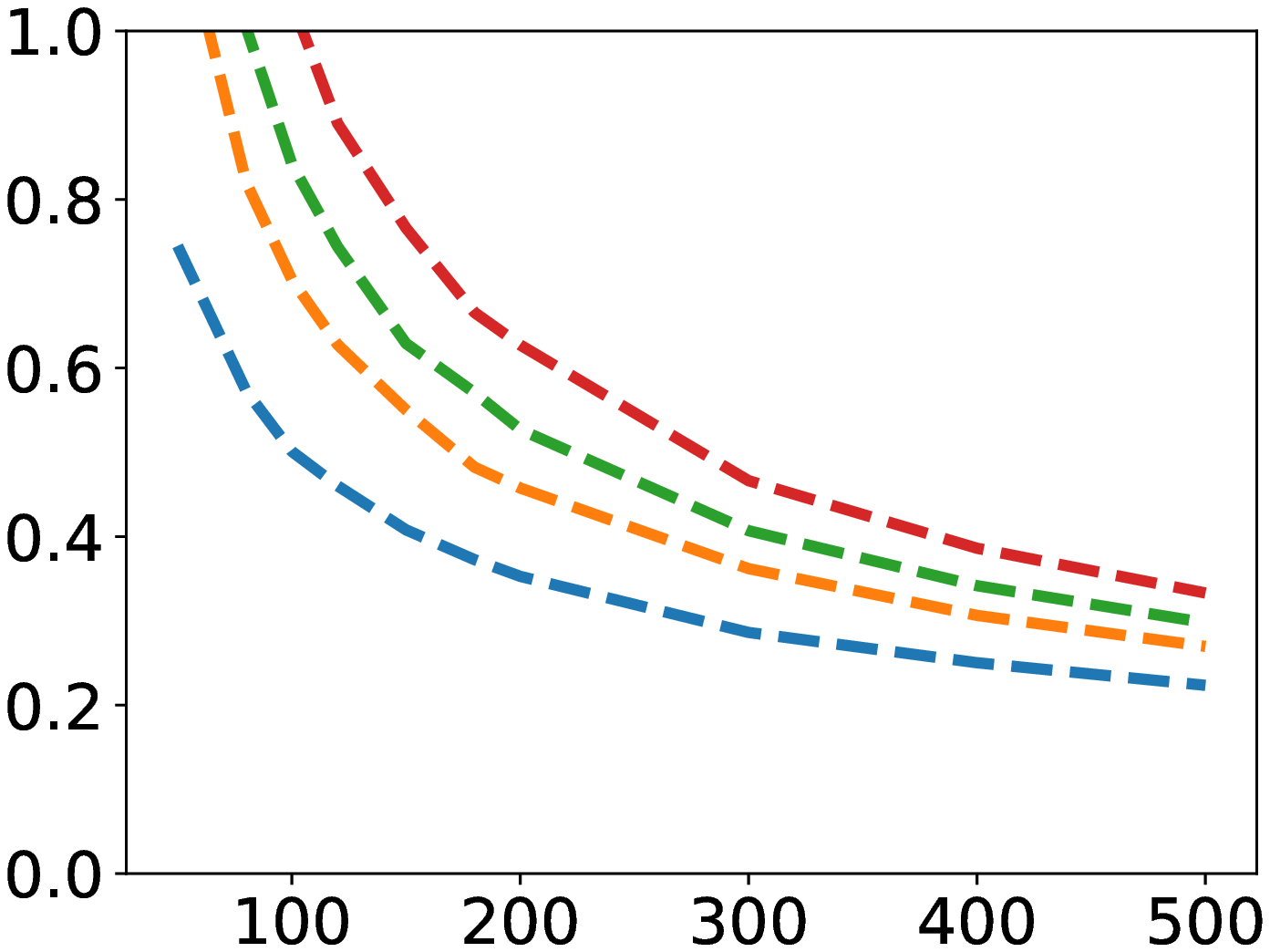} 
	\put(-20,-1){\rotatebox{90}{ {\small \ \ \ standardization \  \ \ }}}
	\end{overpic}
	~
	\DeclareGraphicsExtensions{.png}
	\begin{overpic}[width=0.29\textwidth]{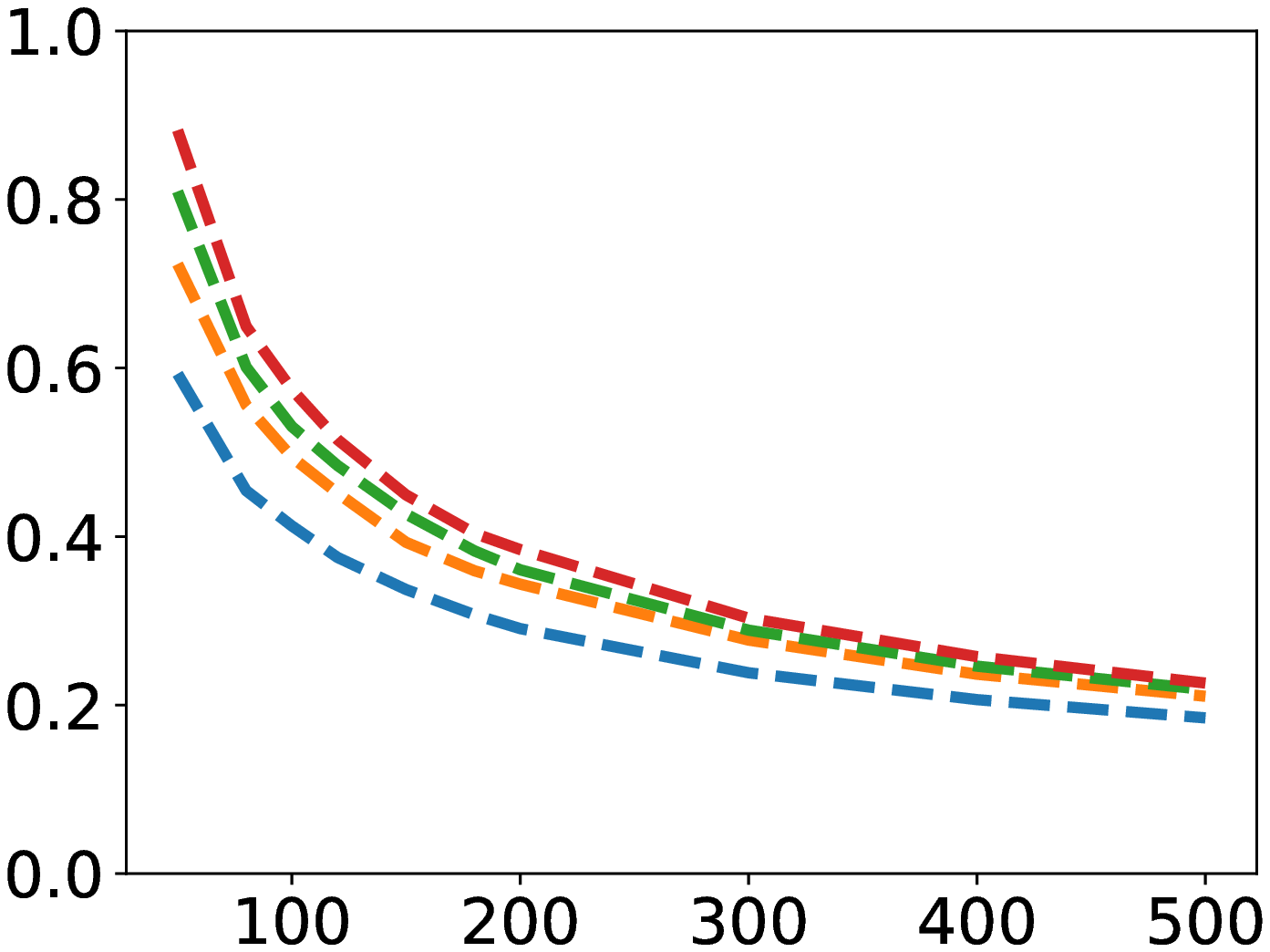} 
	\end{overpic}
	~	
	\begin{overpic}[width=0.29\textwidth]{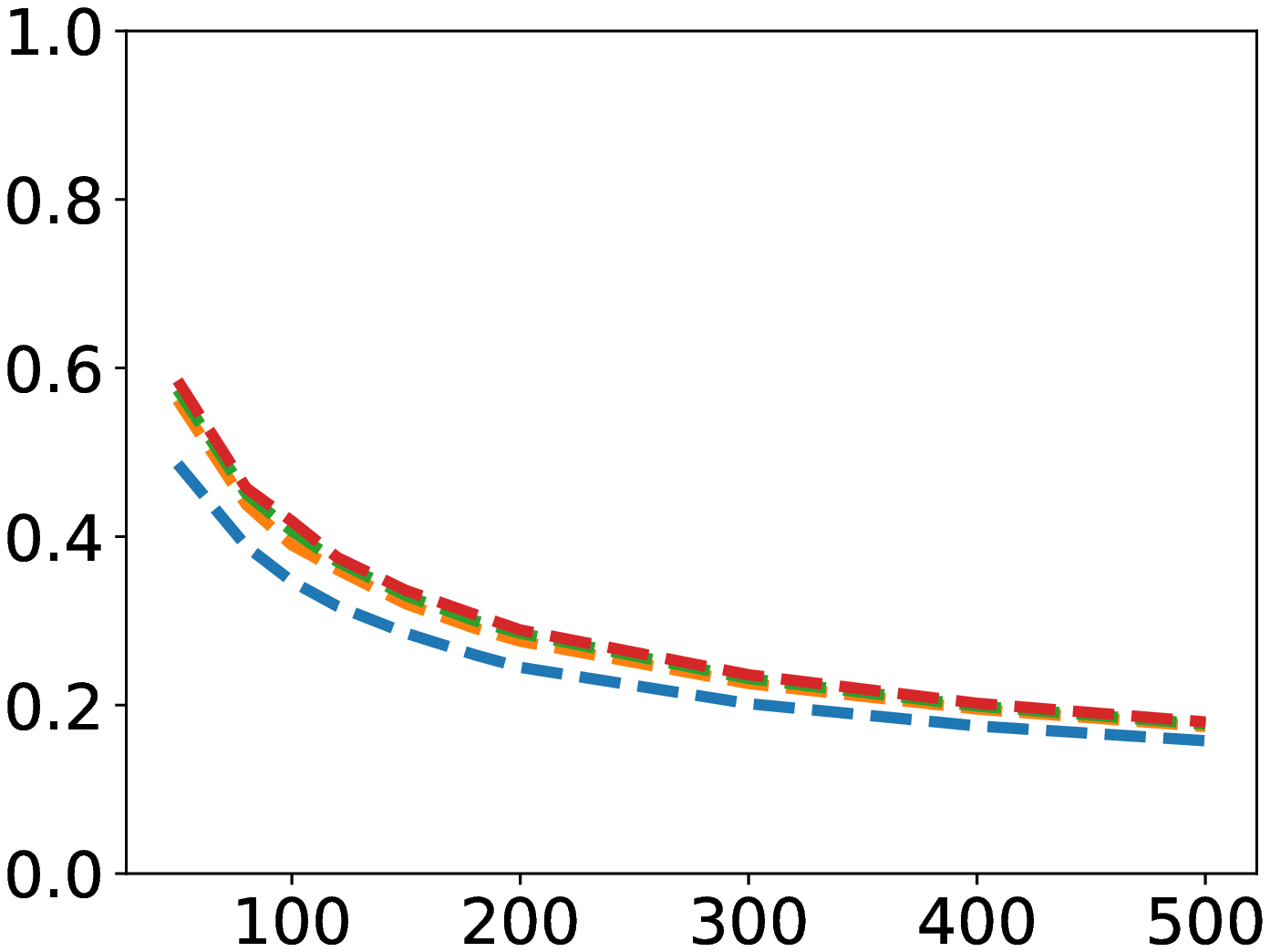} 
	\end{overpic}	
	%
\end{figure}

\vspace{-0.5cm}

\begin{figure}[H]	
	\quad\quad\quad 
	\begin{overpic}[width=0.29\textwidth]{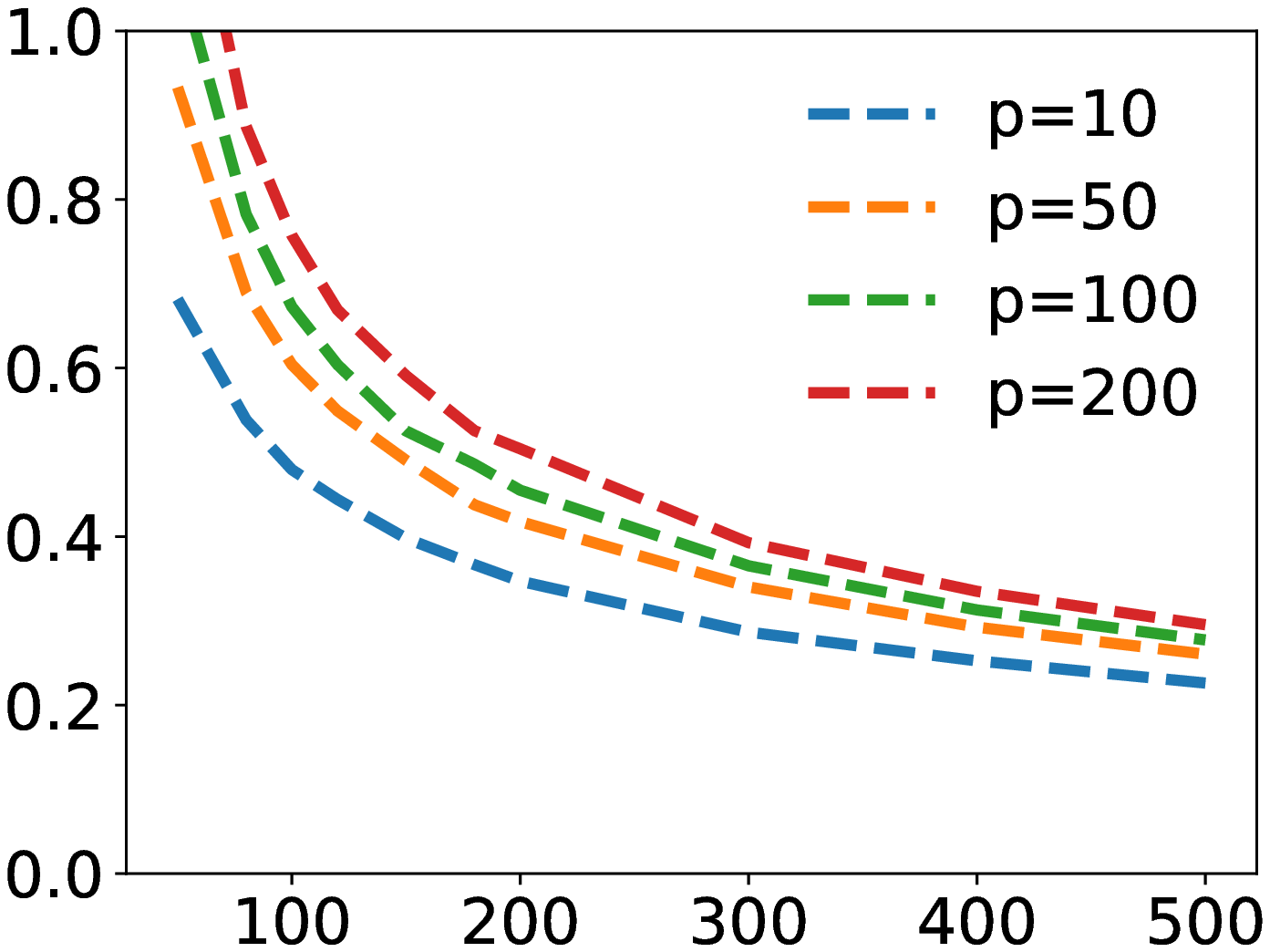} 
	\put(-21,1){\rotatebox{90}{\ $\sqrt{ \ \ }$}}
	\put(-20,-3){\rotatebox{90}{  { \ \ \ \ \ \ \small transformation \ \ } }}

	\end{overpic}
	~
	\DeclareGraphicsExtensions{.png}
	\begin{overpic}[width=0.29\textwidth]{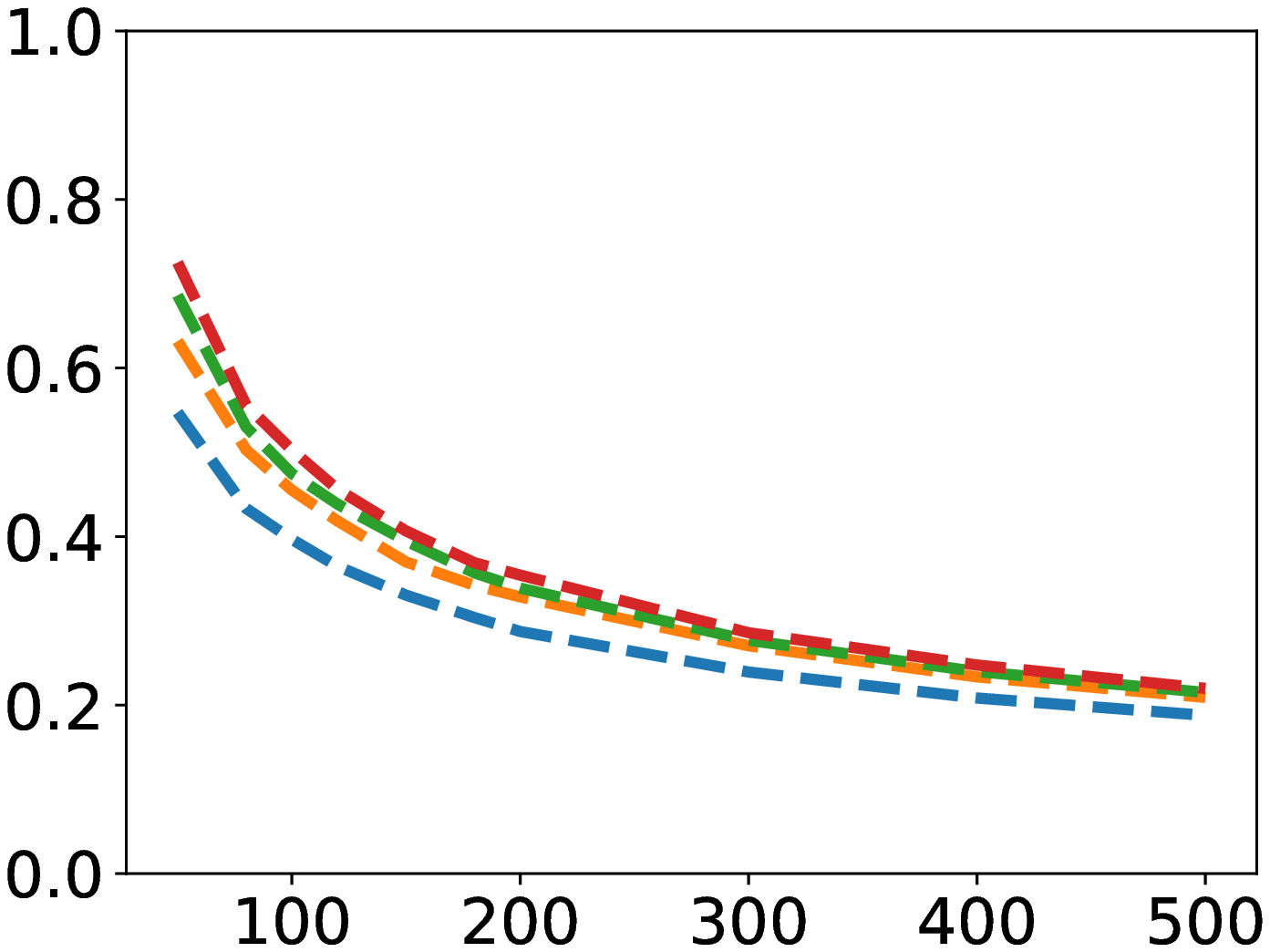} 
	\end{overpic}
	~	
	\begin{overpic}[width=0.29\textwidth]{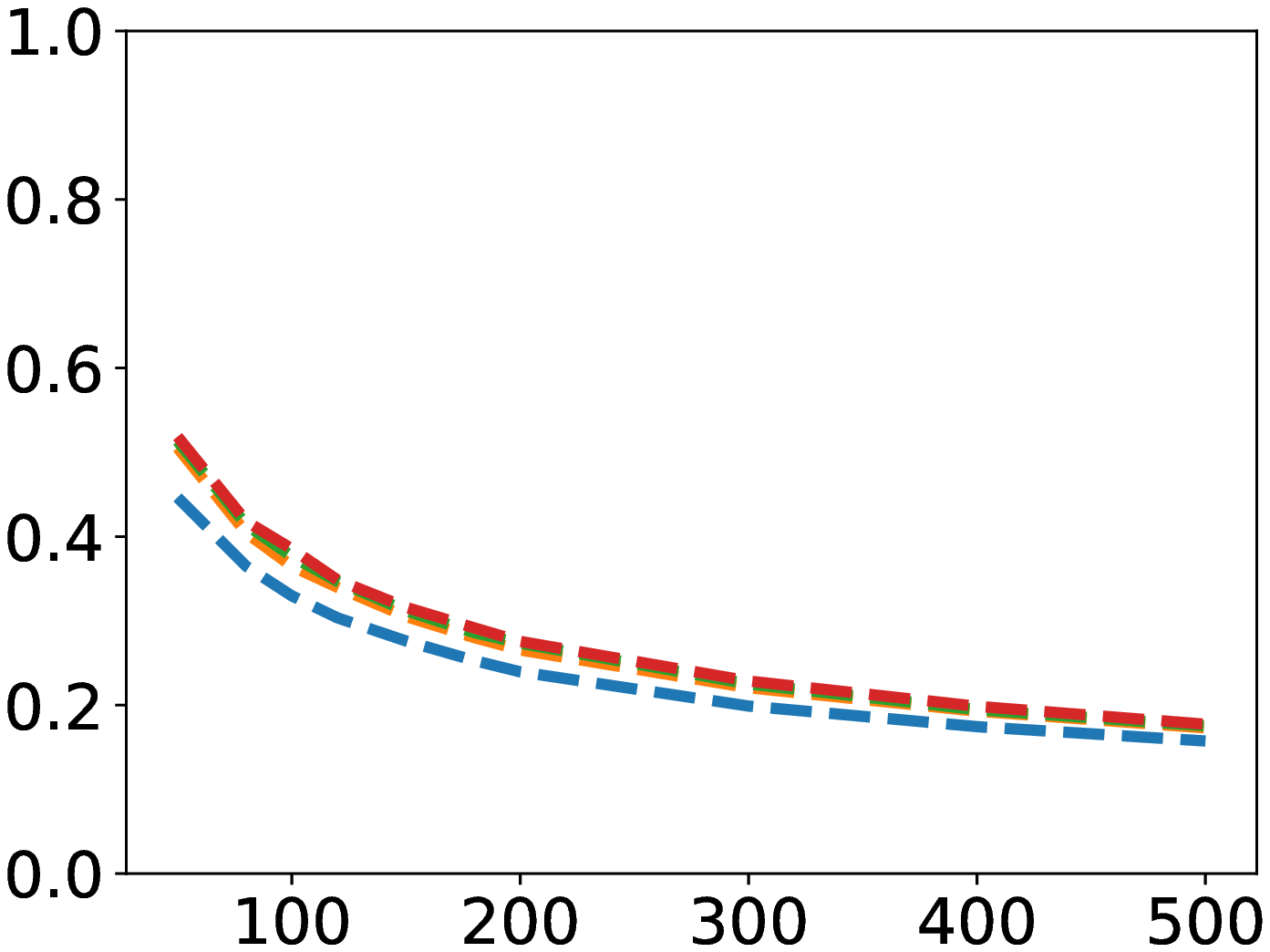} 
	\end{overpic}	
	\vspace{+.2cm}	
	\caption{(Average width versus $n$ in simulation model (i) with a polynomial decay profile). In each of the nine panels, the $y$-axis measures the average width $\E[|\hat{\I}_1|+\cdots+|\hat{\I}_{5}|] / 5$, and the $x$-axis measures $n$. The colored curves correspond to the different values of $p=10$, $50$, $100$, $200$, indicated in the legend. The three rows and three columns correspond to labeled choices of transformations and values of the eigenvalue decay parameter $\gamma$.} 
	\label{SUPP:fig5}
\end{figure}

\newpage

\begin{figure}[H]	
\vspace{0.5cm}

	\quad\quad\quad 
	\begin{overpic}[width=0.29\textwidth]{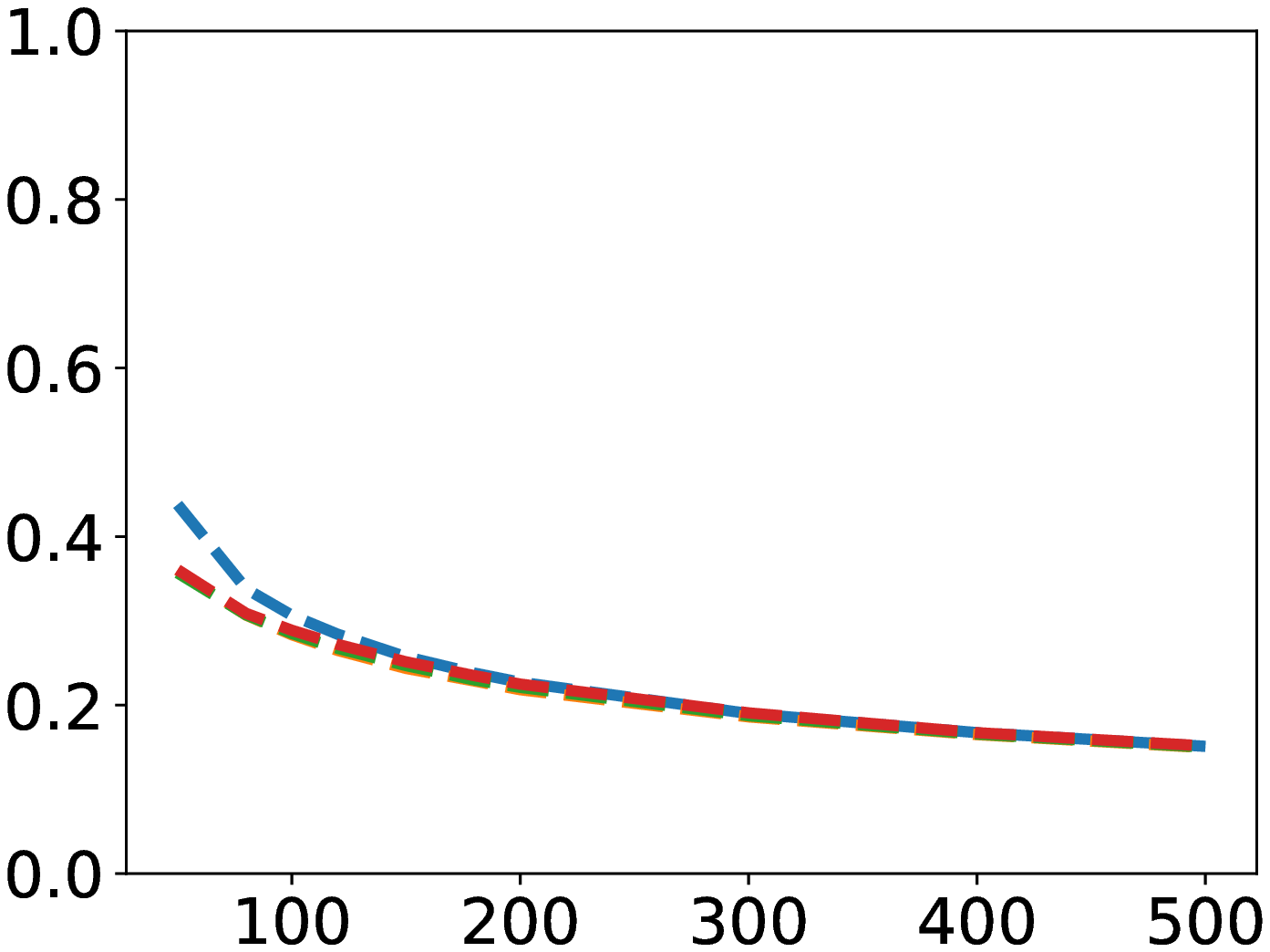} 
	\put(25,80){ \ul{\ \ \  \ $\gamma=0.7$ \ \ \ \    }}	
	\put(-20,-5){\rotatebox{90}{ {\small \ \ \ log transformation  \ \ }}}
	\end{overpic}
	~
	\DeclareGraphicsExtensions{.png}
	\begin{overpic}[width=0.29\textwidth]{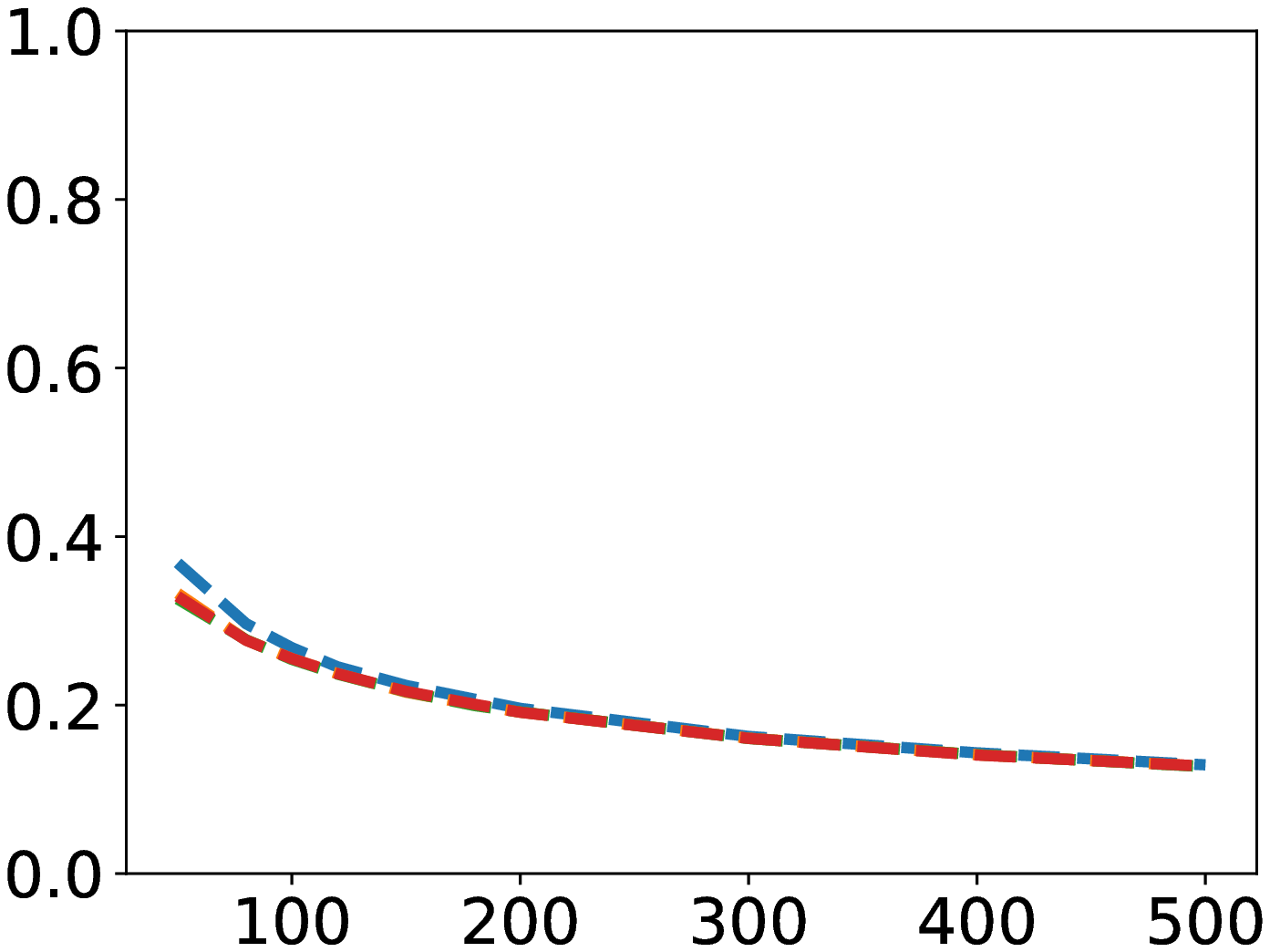} 
	\put(25,80){ \ul{\ \ \  \ $\gamma=1.0$ \ \ \ \    }}
	\end{overpic}
	~	
	\begin{overpic}[width=0.29\textwidth]{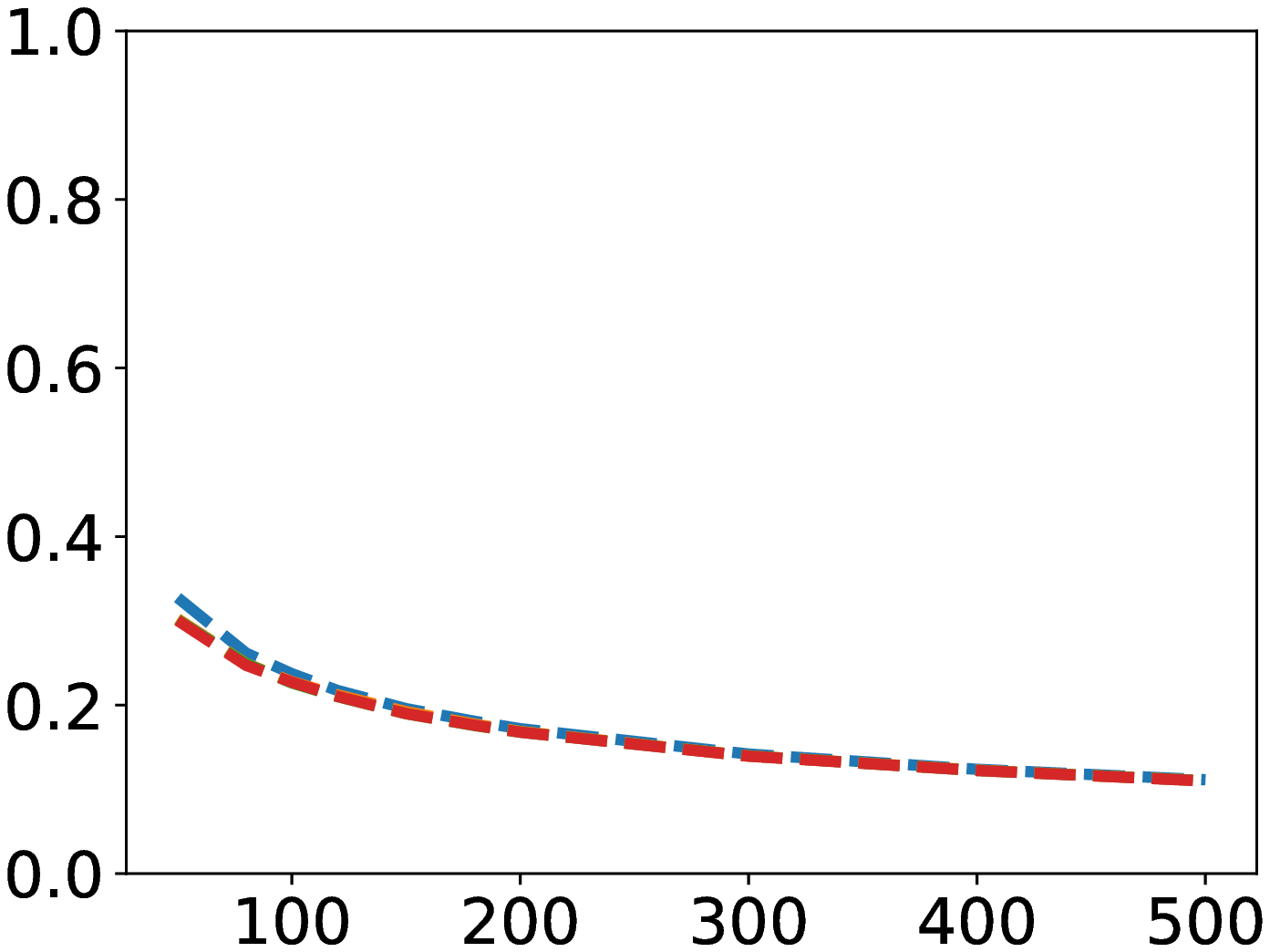} 
	\put(25,80){ \ul{\ \ \  \ $\gamma=1.3$ \ \ \ \    }}
	\end{overpic}	
	%
\end{figure}

\vspace{-0.5cm}

\begin{figure}[H]	
	\quad\quad\quad 
	\begin{overpic}[width=0.29\textwidth]{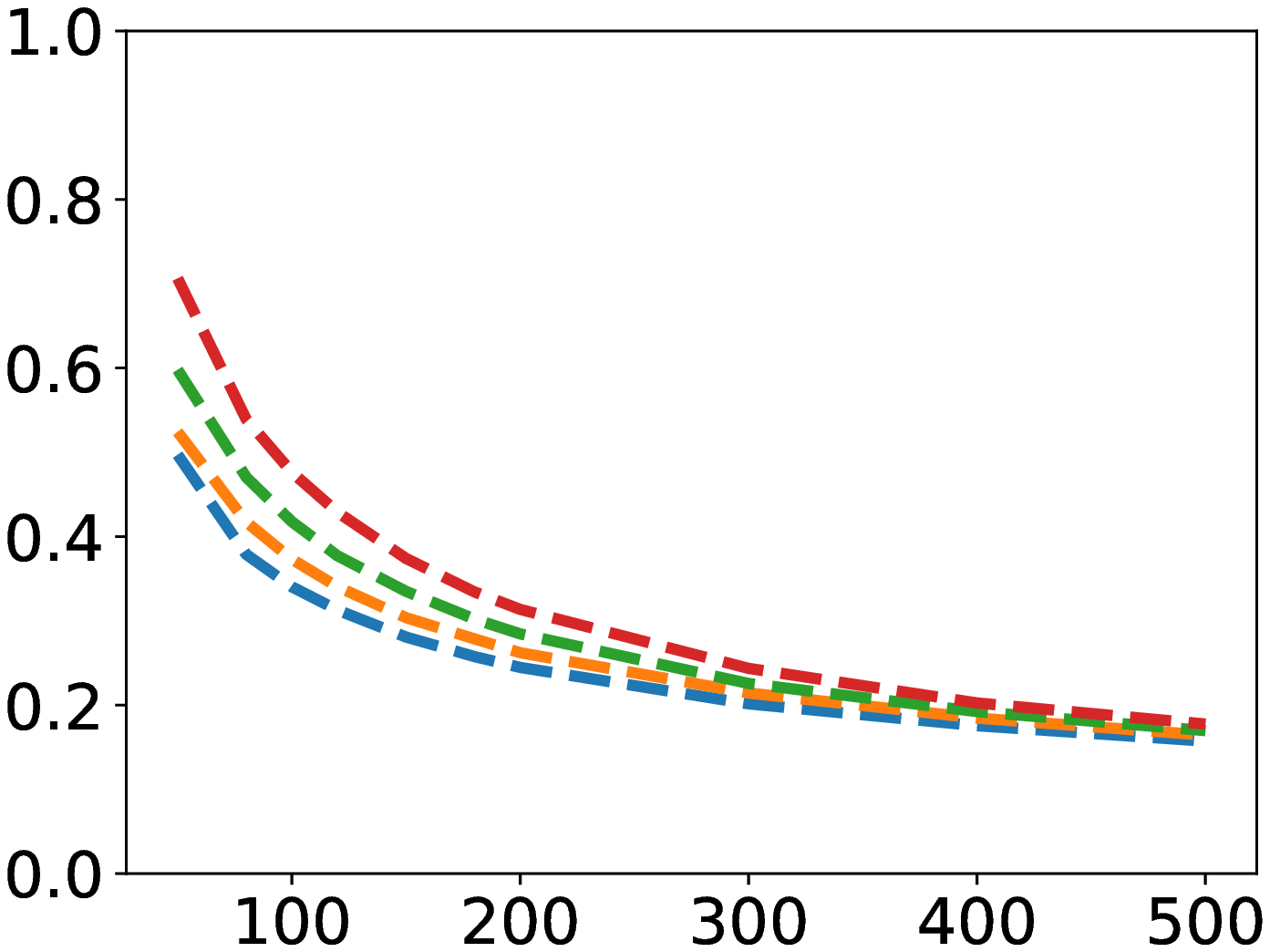} 
	\put(-20,-1){\rotatebox{90}{ {\small \ \ \ standardization \  \ \ }}}
	\end{overpic}
	~
	\DeclareGraphicsExtensions{.png}
	\begin{overpic}[width=0.29\textwidth]{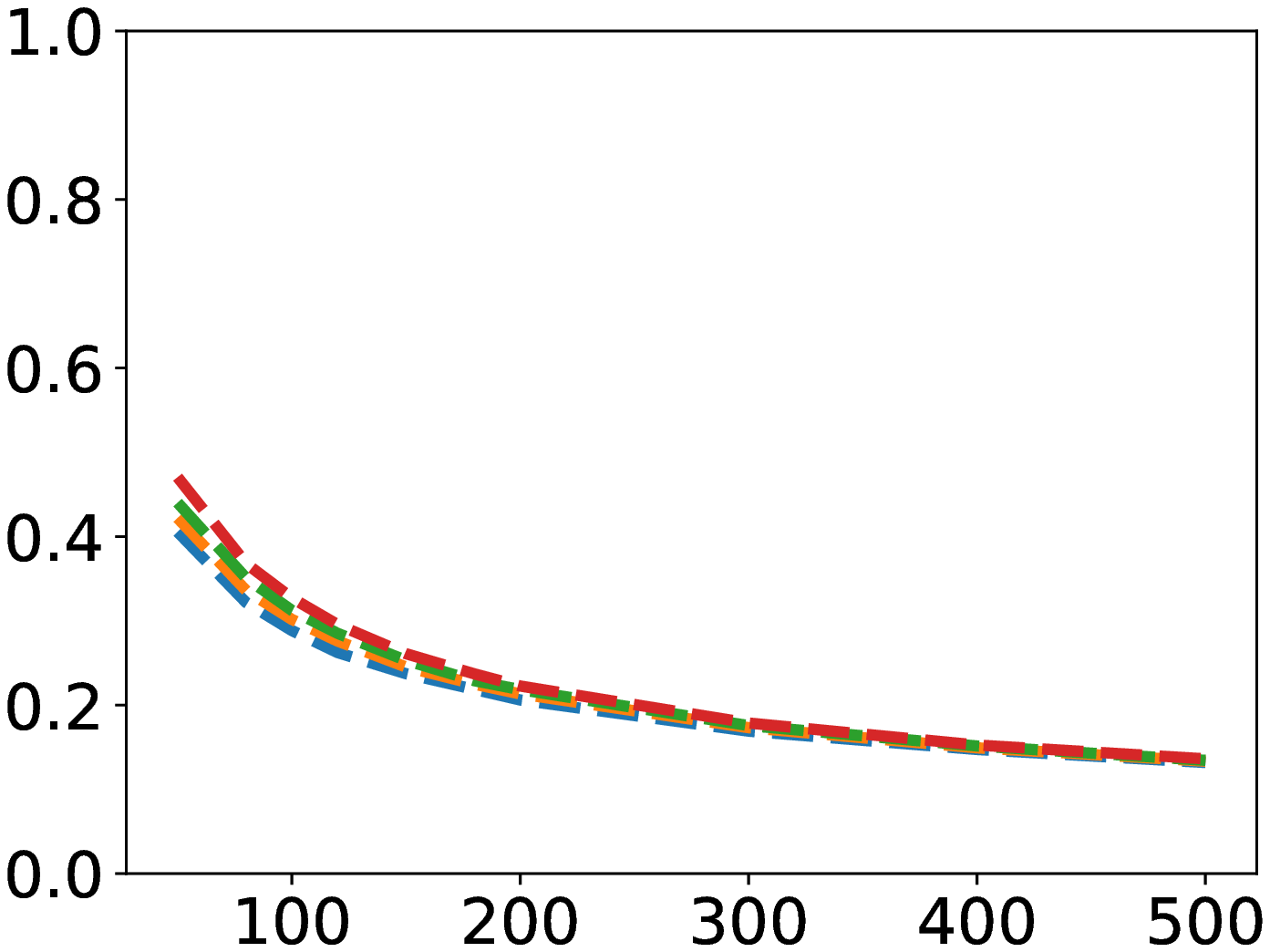} 
	\end{overpic}
	~	
	\begin{overpic}[width=0.29\textwidth]{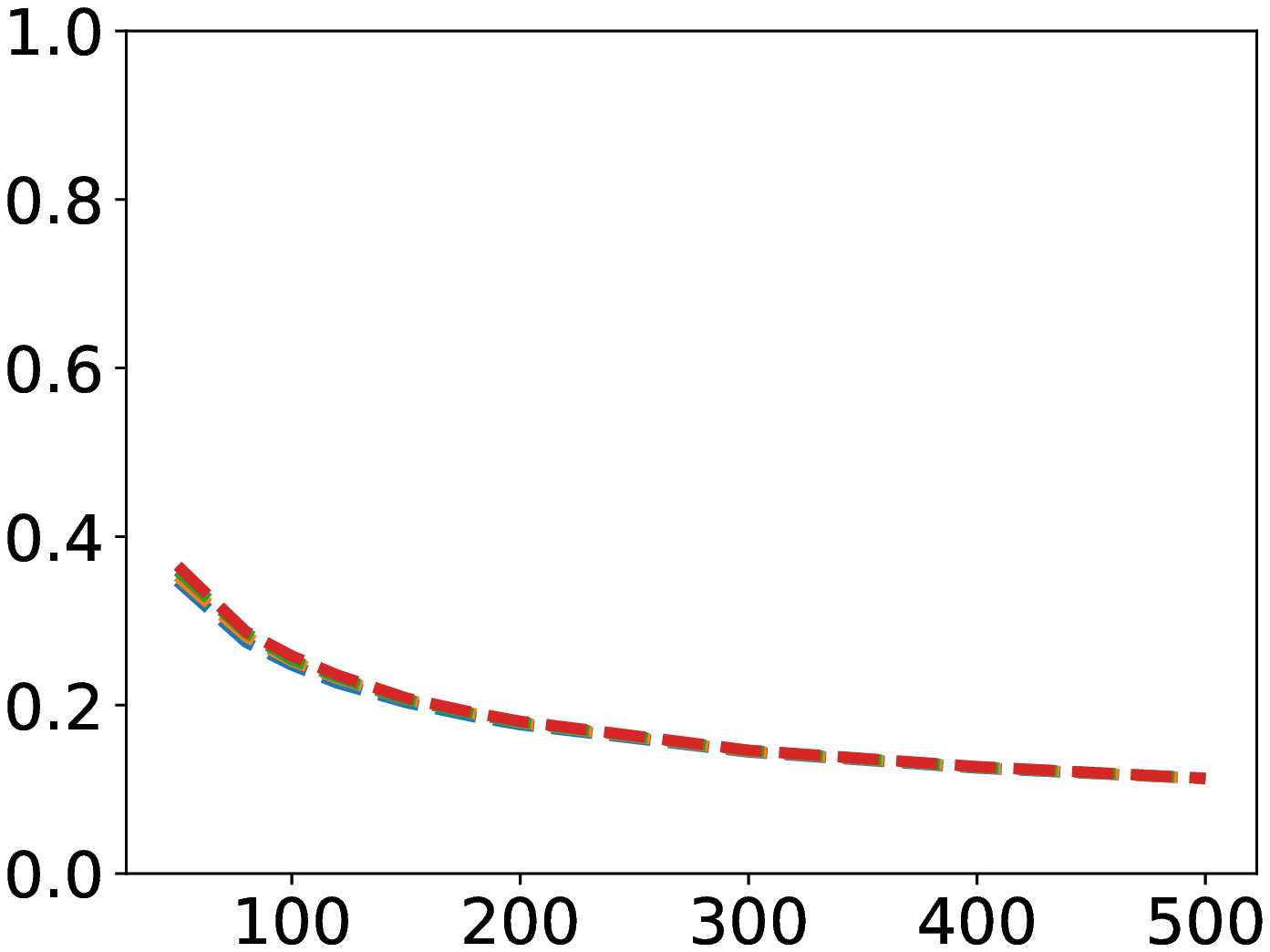} 
	\end{overpic}	

\end{figure}

\vspace{-0.5cm}

\begin{figure}[H]	
	\quad\quad\quad 
	\begin{overpic}[width=0.29\textwidth]{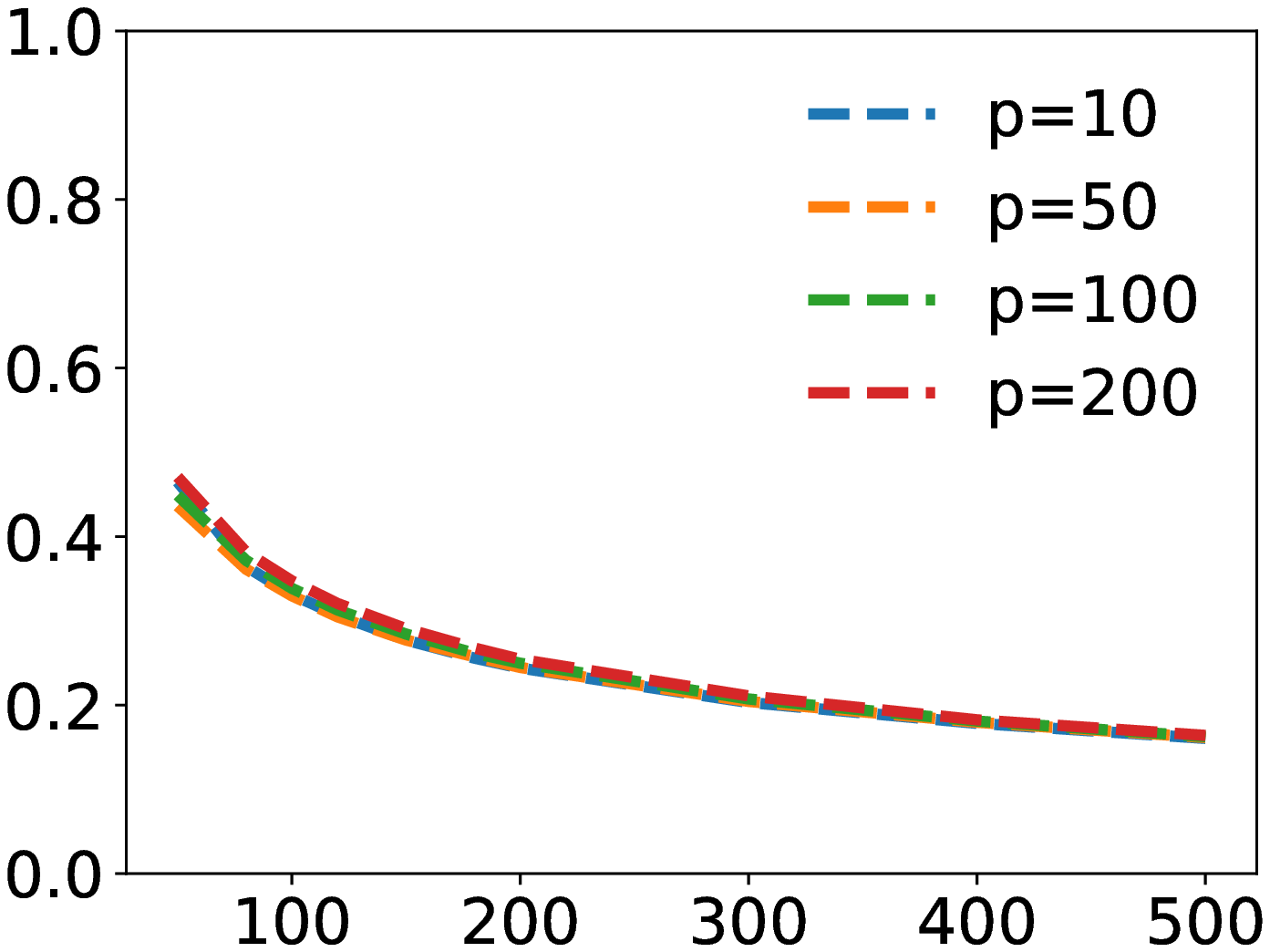} 
    \put(-21,1){\rotatebox{90}{\ $\sqrt{ \ \ }$}}
	\put(-20,-3){\rotatebox{90}{  { \ \ \ \ \ \ \small transformation \ \ } }}

	\end{overpic}
	~
	\DeclareGraphicsExtensions{.png}
	\begin{overpic}[width=0.29\textwidth]{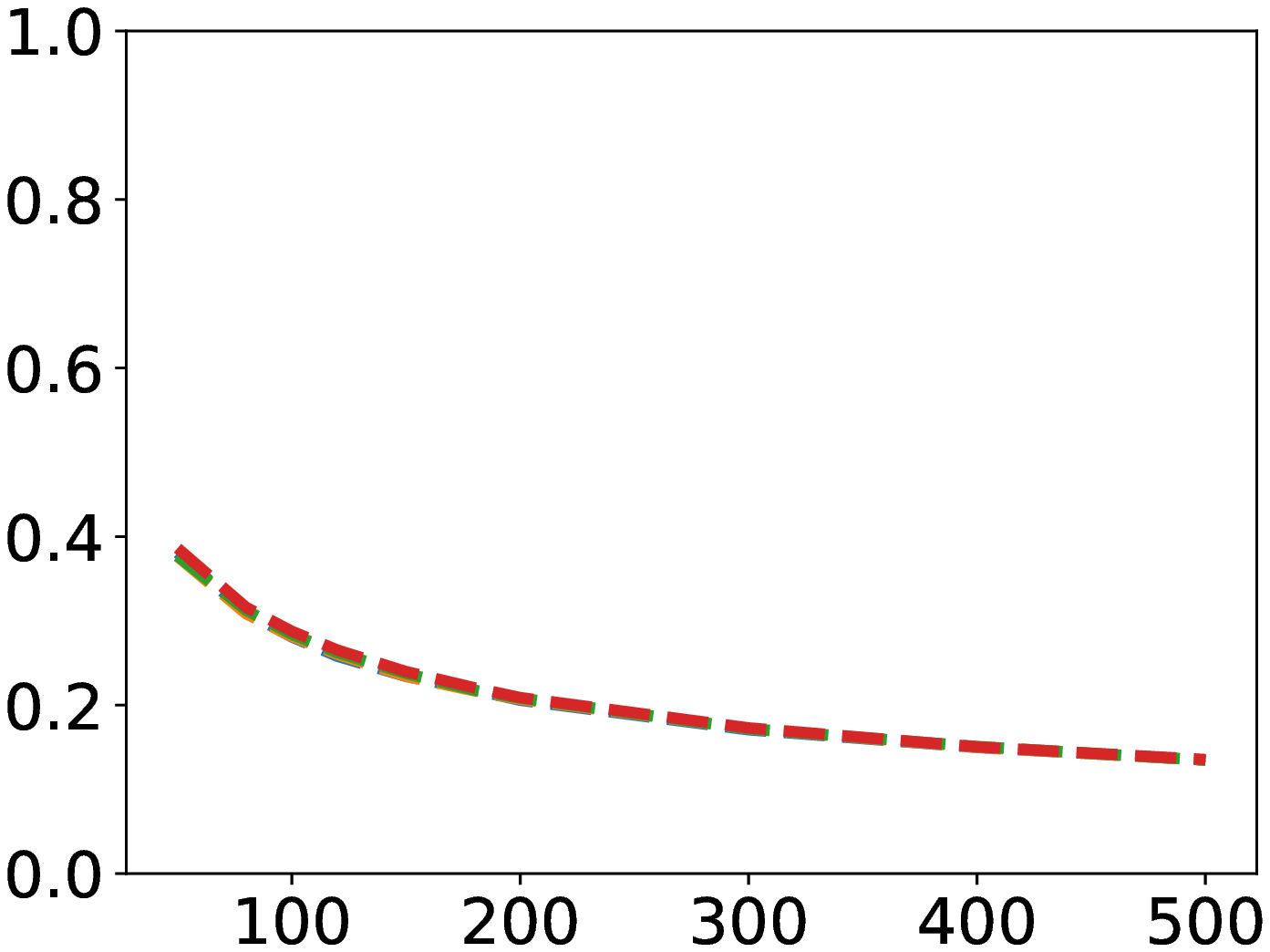} 
	\end{overpic}
	~	
	\begin{overpic}[width=0.29\textwidth]{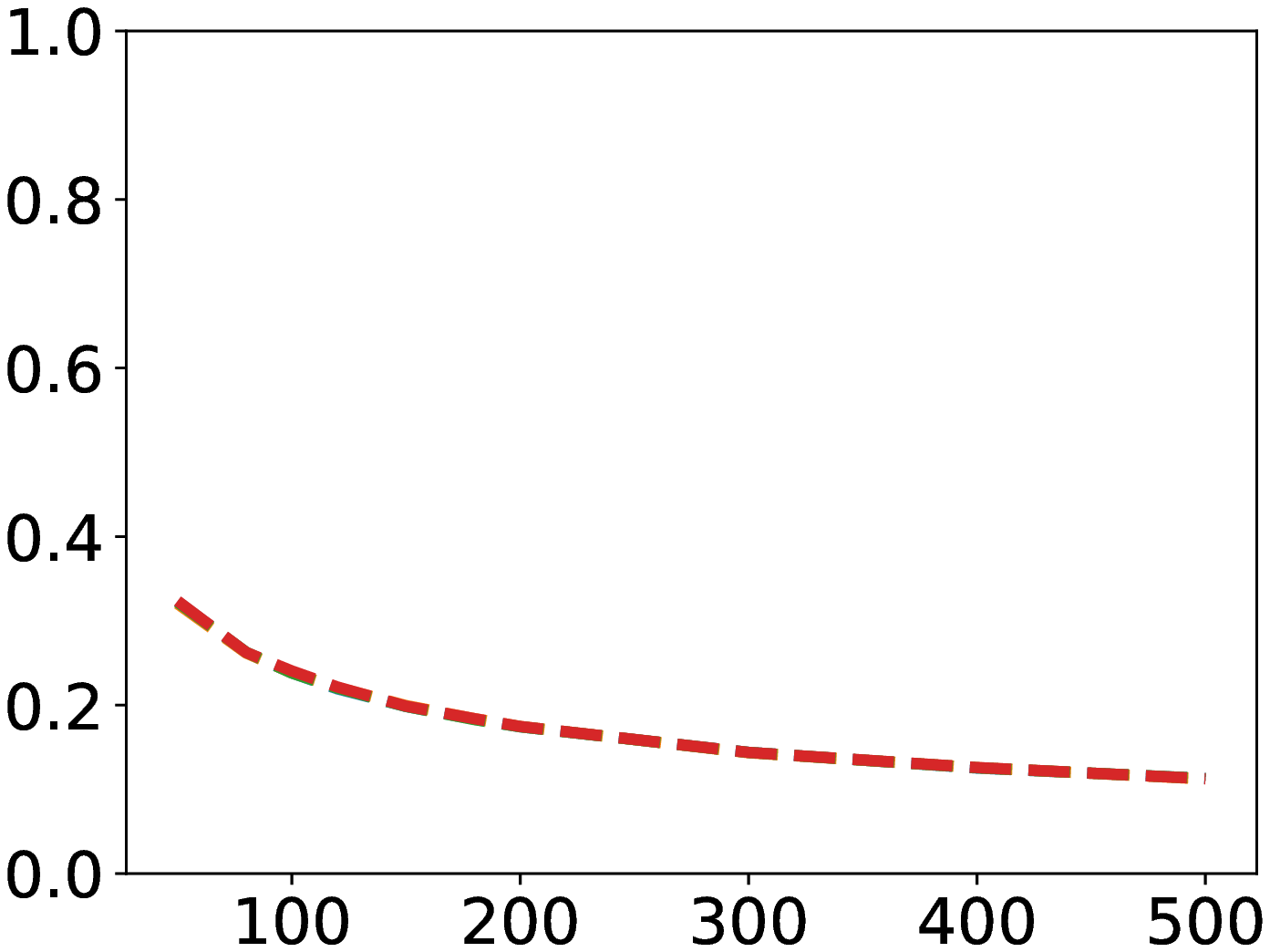} 
	\end{overpic}	
	\vspace{+0.2cm}	
	\caption{(Average width versus $n$ in simulation model (ii) with a polynomial decay profile). The plotting scheme is the same as described in the caption of Figure~\ref{SUPP:fig5} above.} 
	\label{SUPP:fig6}
\end{figure}

\newpage
\begin{figure}[H]	
\vspace{0.5cm}
	\quad\quad\quad 
	\begin{overpic}[width=0.29\textwidth]{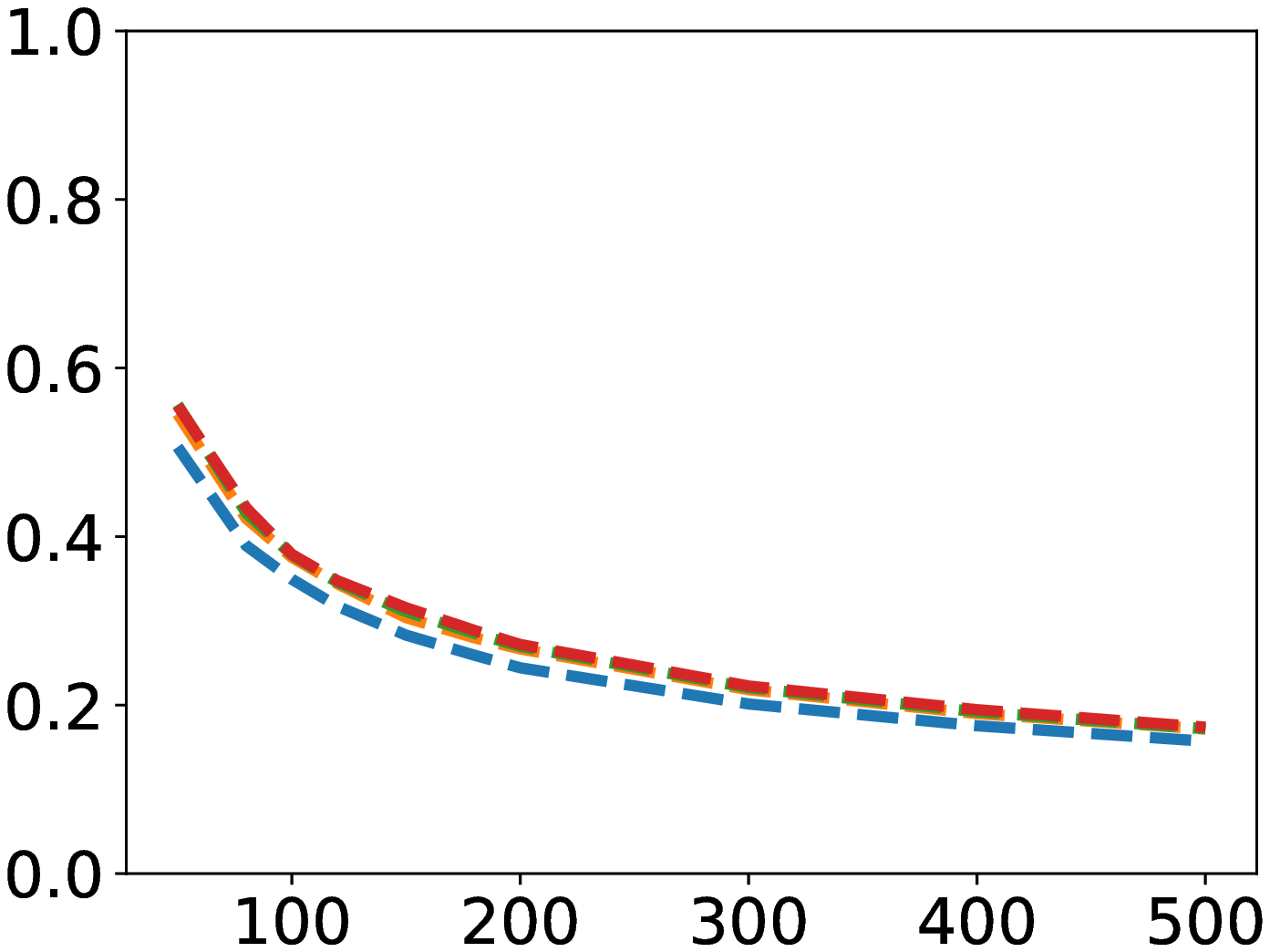} 
    \put(25,80){ \ul{\ \ \  \ $\delta=0.7$ \ \ \ \    }}
	\put(-20,-5){\rotatebox{90}{ {\small \ \ \ log transformation  \ \ }}}
\end{overpic}
	~
	\DeclareGraphicsExtensions{.png}
	\begin{overpic}[width=0.29\textwidth]{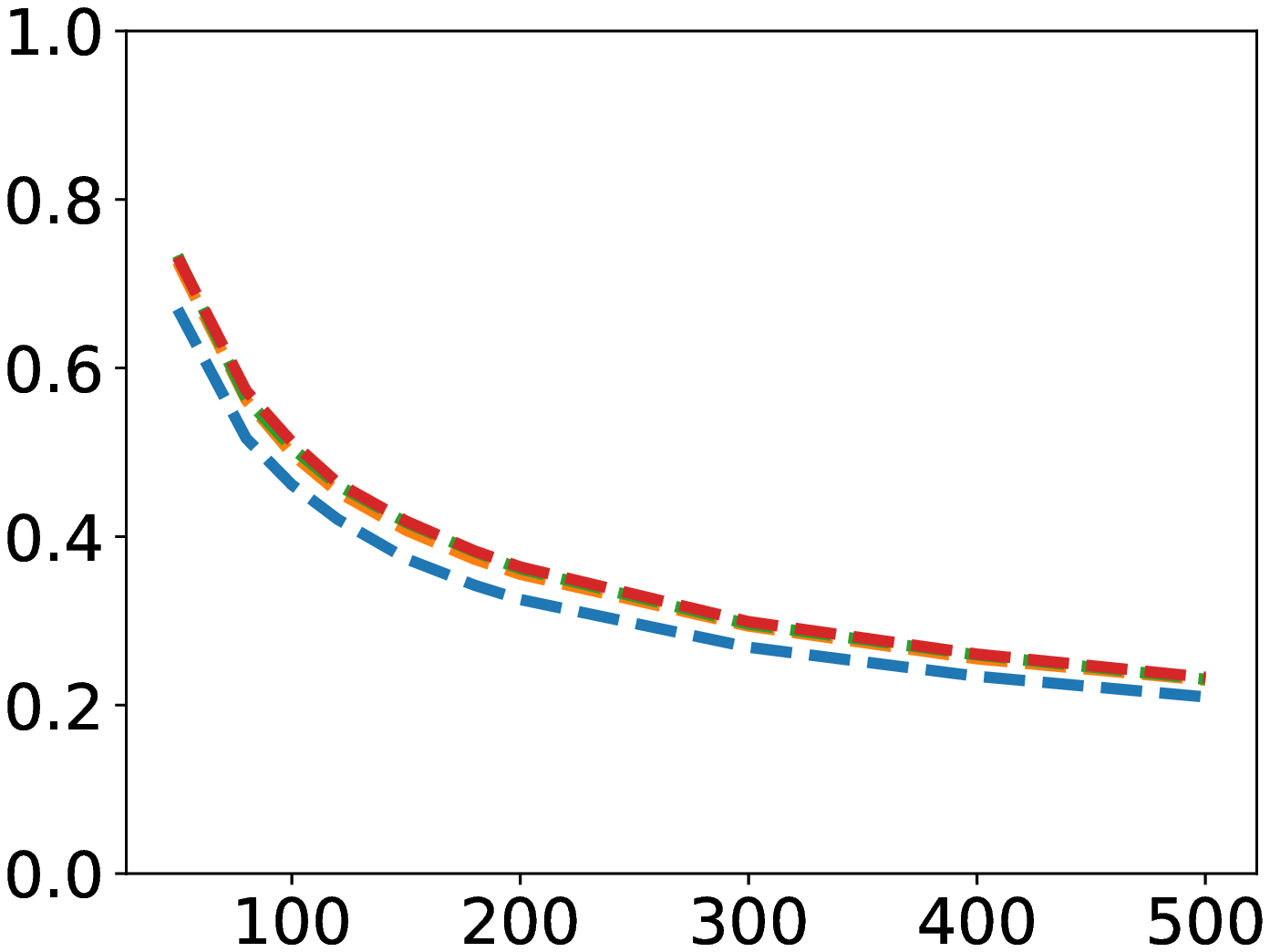} 
    \put(25,80){ \ul{\ \ \  \ $\delta=0.8$ \ \ \ \    }}
	\end{overpic}
	~	
	\begin{overpic}[width=0.29\textwidth]{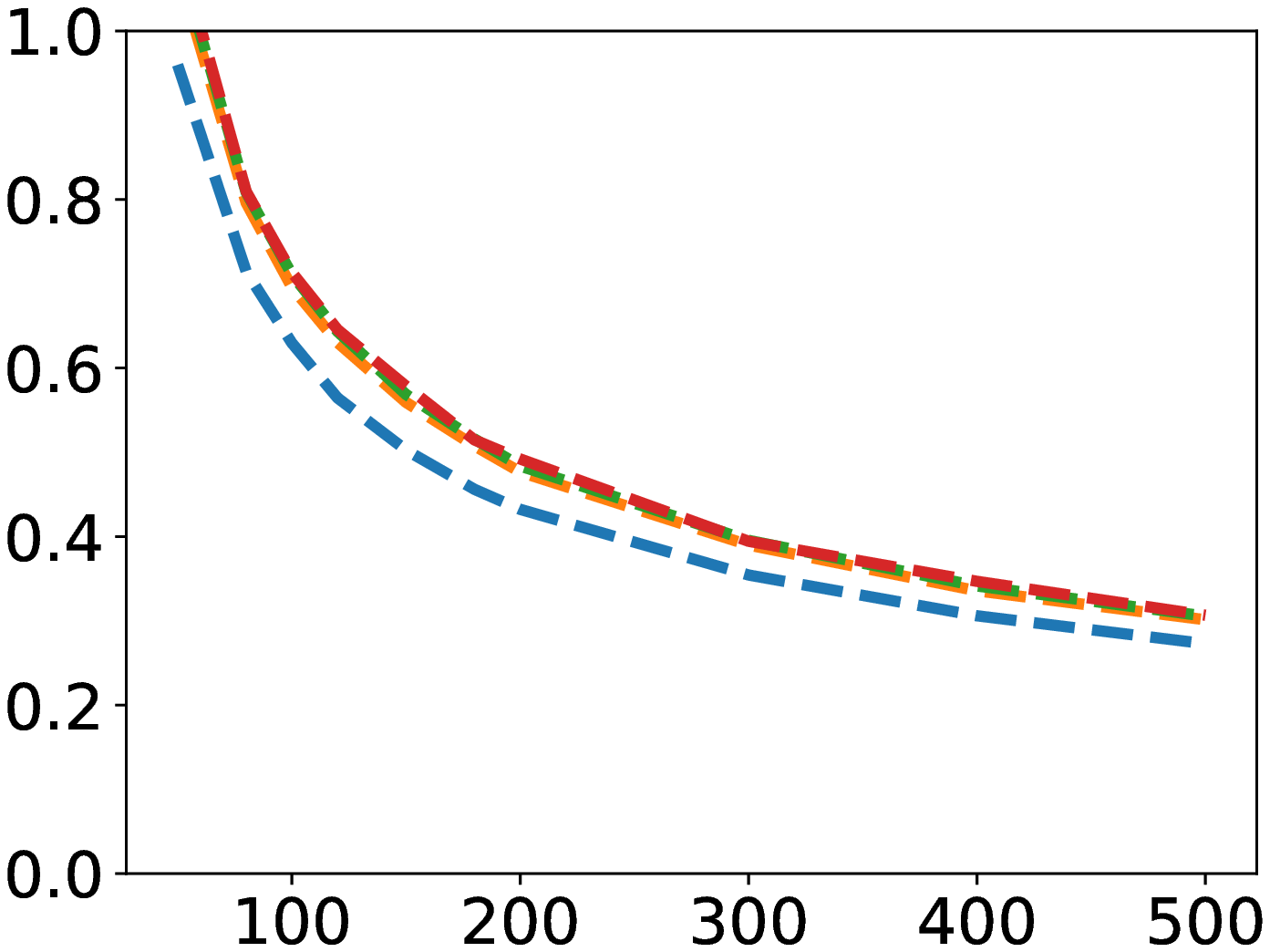} 
	\put(25,80){ \ul{\ \ \  \ $\delta=0.9$ \ \ \ \    }}
	\end{overpic}	
\end{figure}

\vspace{-0.5cm}

\begin{figure}[H]	
	\quad\quad\quad 
	\begin{overpic}[width=0.29\textwidth]{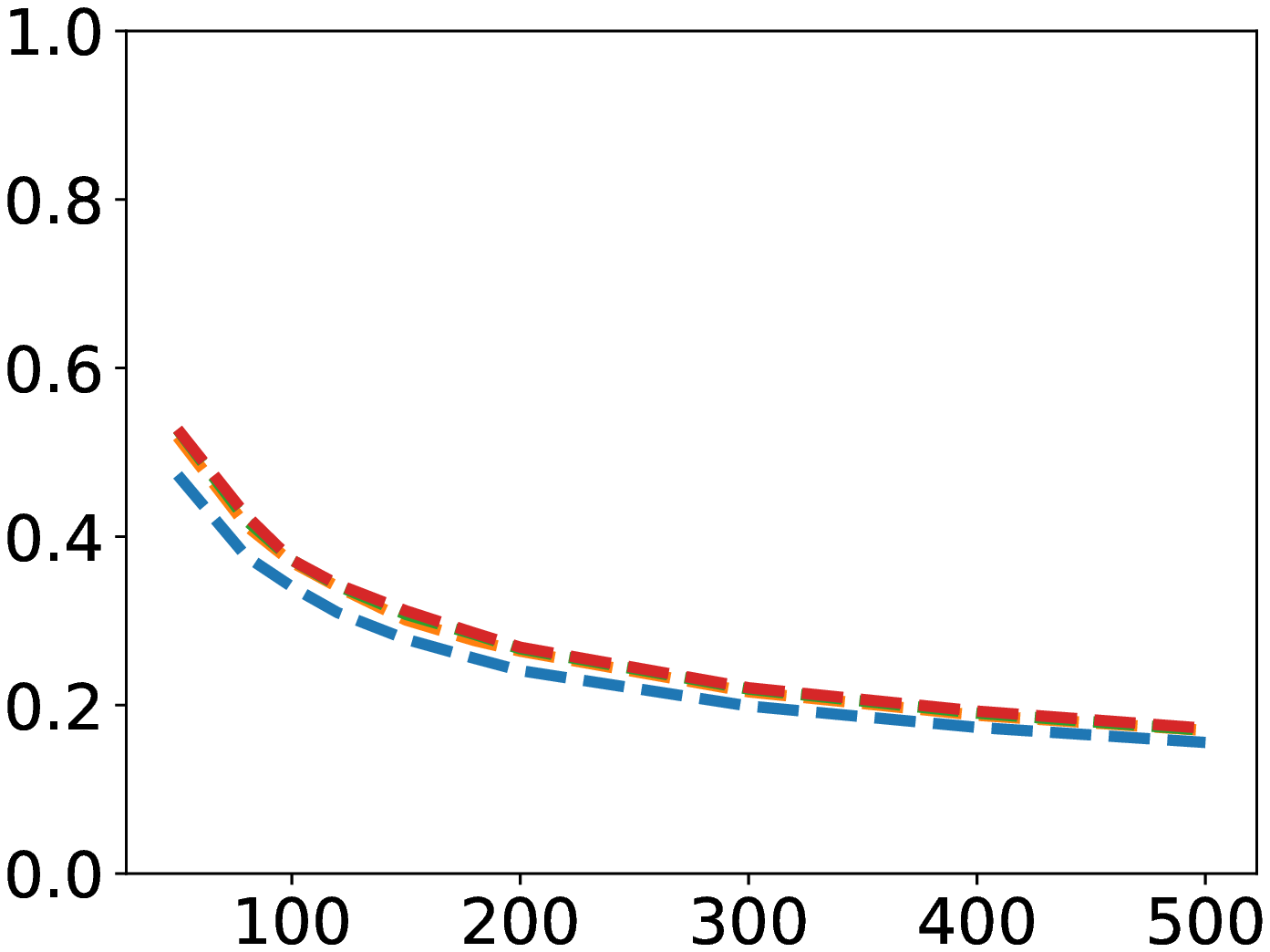} 
	\put(-20,-1){\rotatebox{90}{ {\small \ \ \ standardization \  \ \ }}}
	\end{overpic}
	~
	\DeclareGraphicsExtensions{.png}
	\begin{overpic}[width=0.29\textwidth]{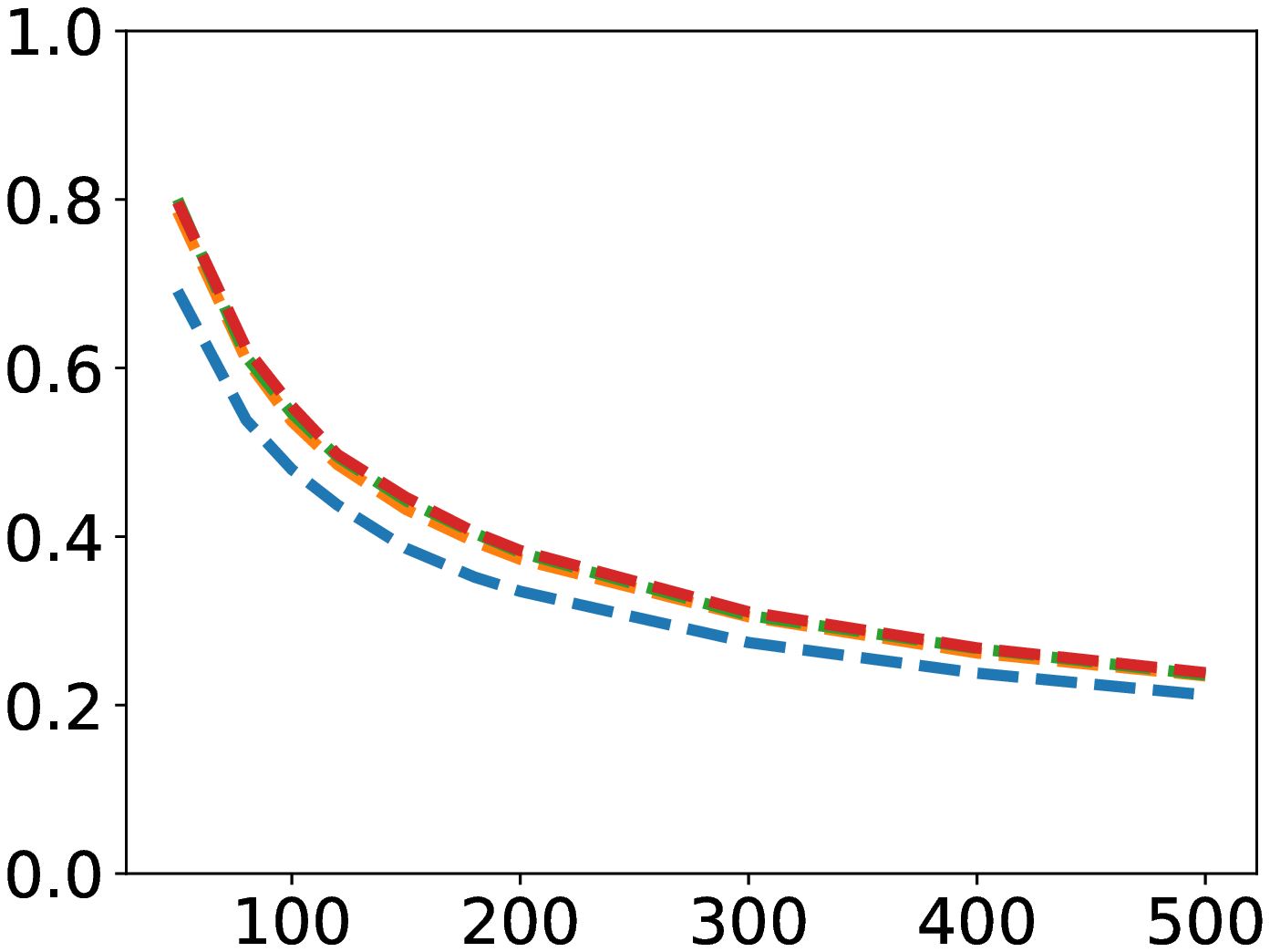} 
	\end{overpic}
	~	
	\begin{overpic}[width=0.29\textwidth]{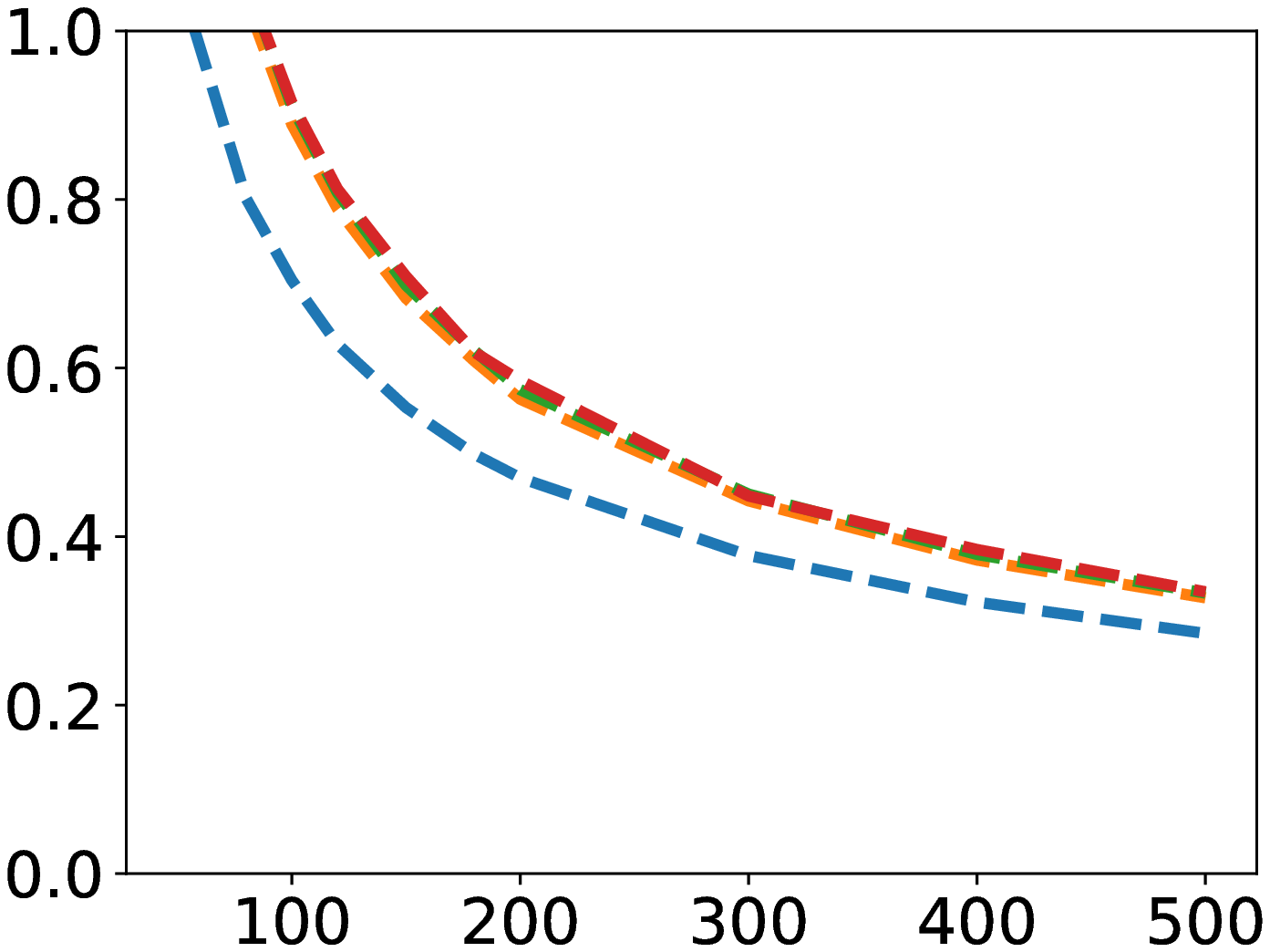} 
		 				
	\end{overpic}	
\end{figure}

\vspace{-0.5cm}

\begin{figure}[H]	
	\quad\quad\quad 
	\begin{overpic}[width=0.29\textwidth]{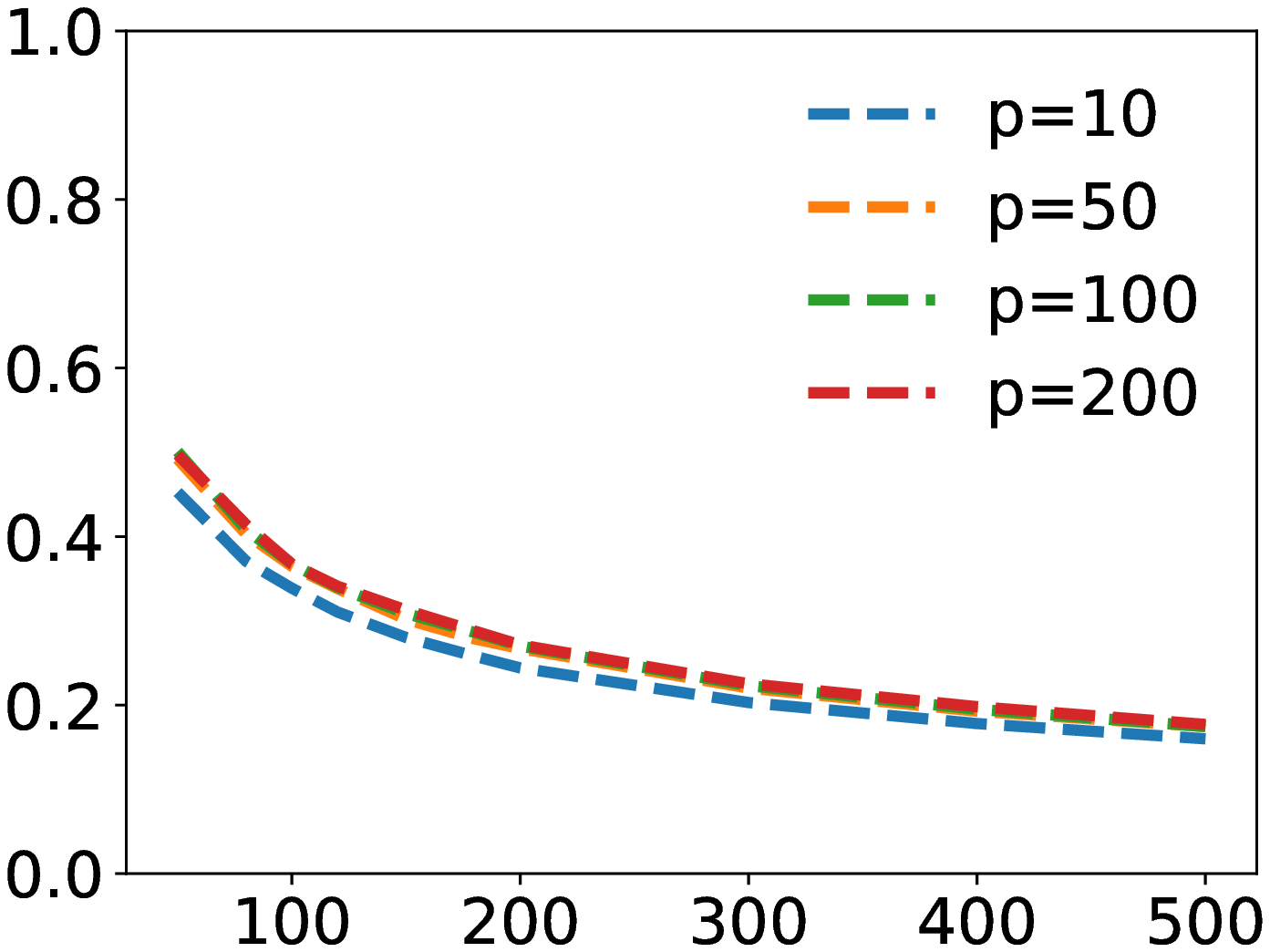} 
	\put(-21,1){\rotatebox{90}{\ $\sqrt{ \ \ }$}}
	\put(-20,-3){\rotatebox{90}{  { \ \ \ \ \ \ \small transformation \ \ } }}

	\end{overpic}
	~
	\DeclareGraphicsExtensions{.png}
	\begin{overpic}[width=0.29\textwidth]{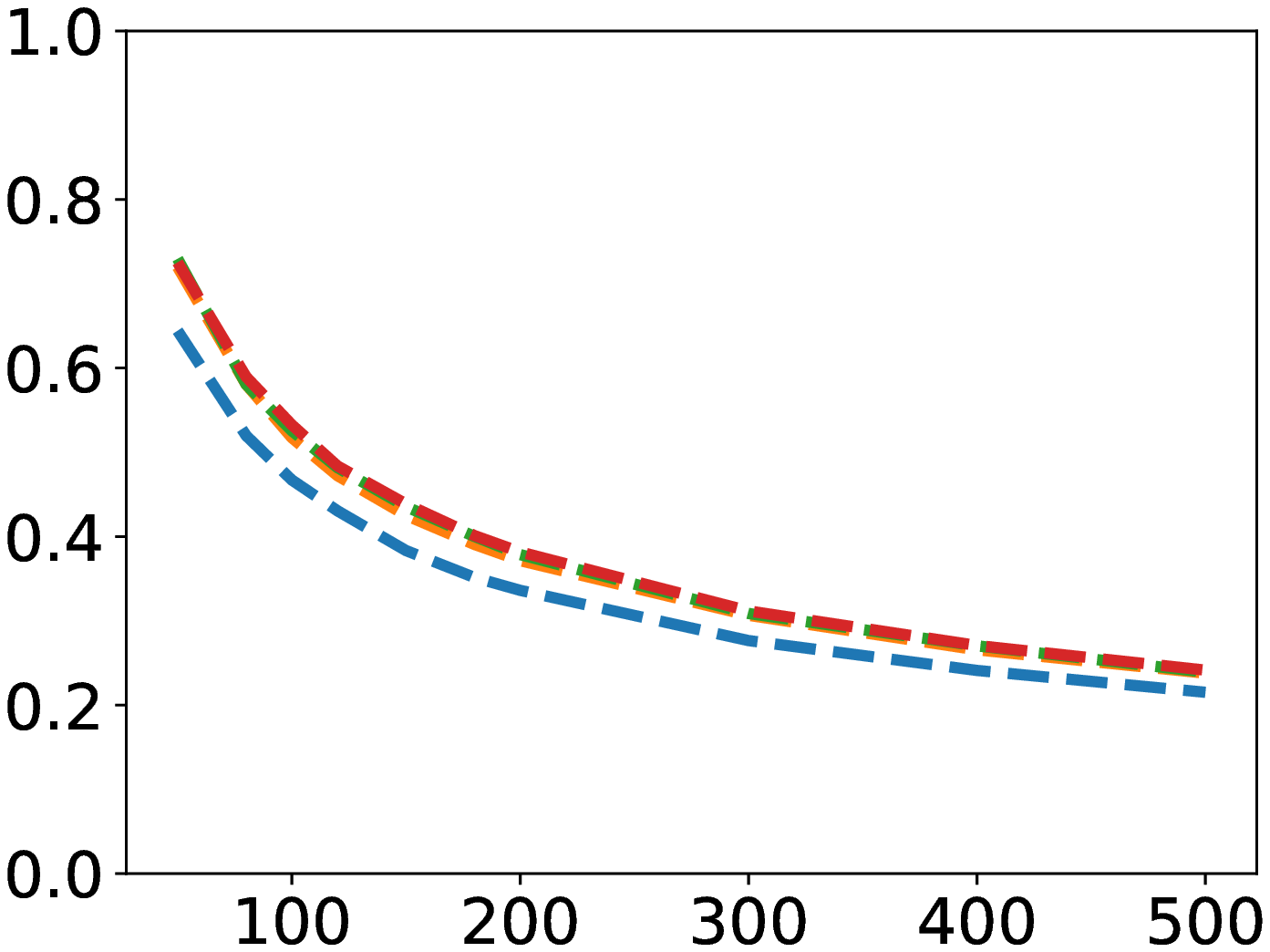} 
	\end{overpic}
	~	
	\begin{overpic}[width=0.29\textwidth]{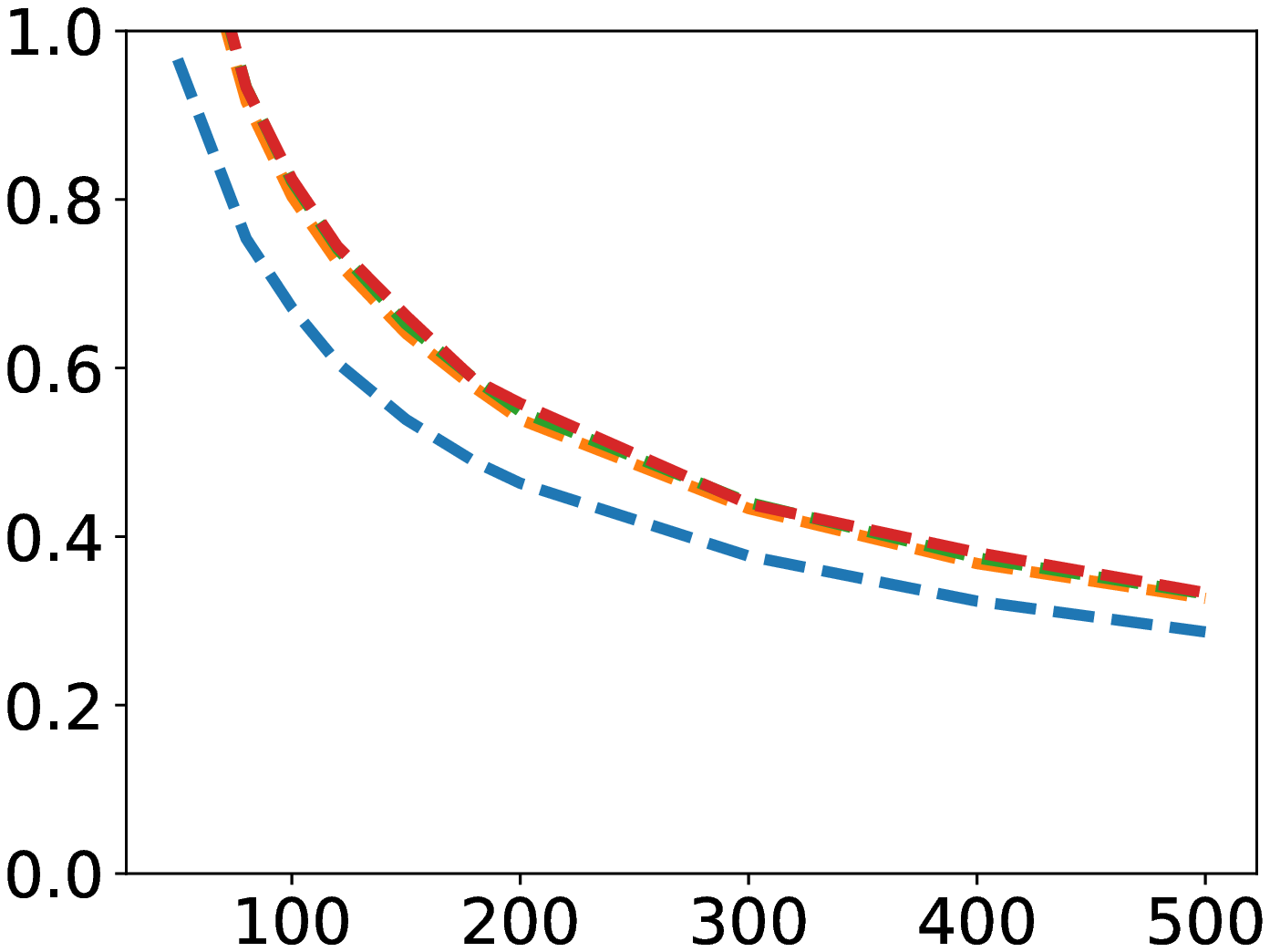} 
		 				
	\end{overpic}	
	\vspace{+0.2cm}	
	\caption{(Average width versus $n$ in simulation model (i) with an exponential decay profile). The plotting scheme is the same as described in the caption of Figure~\ref{SUPP:fig5} above, except that the three columns correspond to values of the eigenvalue decay parameter $\delta$.} 
	\label{SUPP:fig7}
\end{figure}

\newpage
\begin{figure}[H]	
\vspace{0.5cm}
	\quad\quad\quad 
	\begin{overpic}[width=0.29\textwidth]{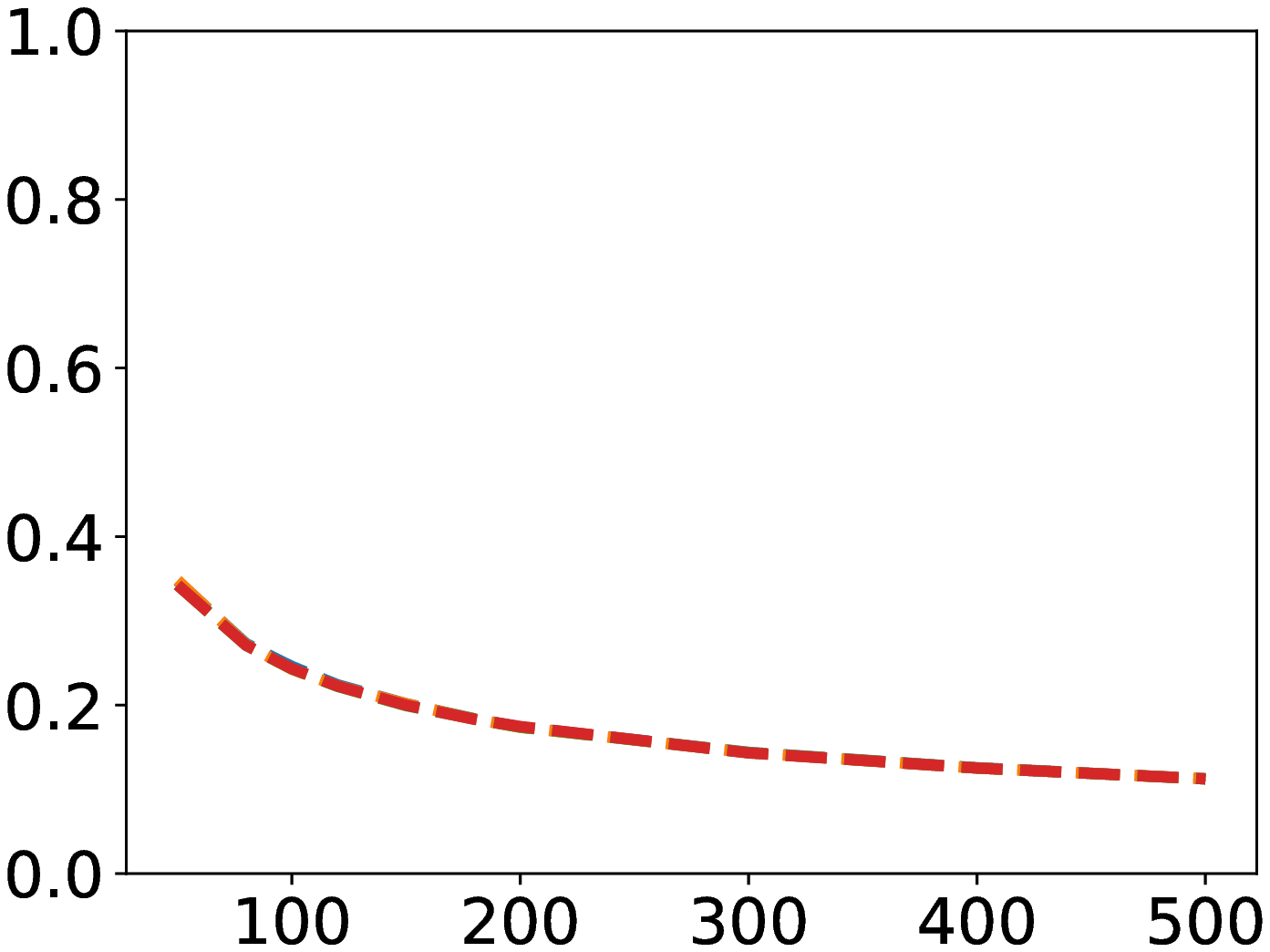} 
    \put(25,80){ \ul{\ \ \  \ $\delta=0.7$ \ \ \ \    }}
	\put(-20,-5){\rotatebox{90}{ {\small \ \ \ log transformation  \ \ }}}
\end{overpic}
	~
	\DeclareGraphicsExtensions{.png}
	\begin{overpic}[width=0.29\textwidth]{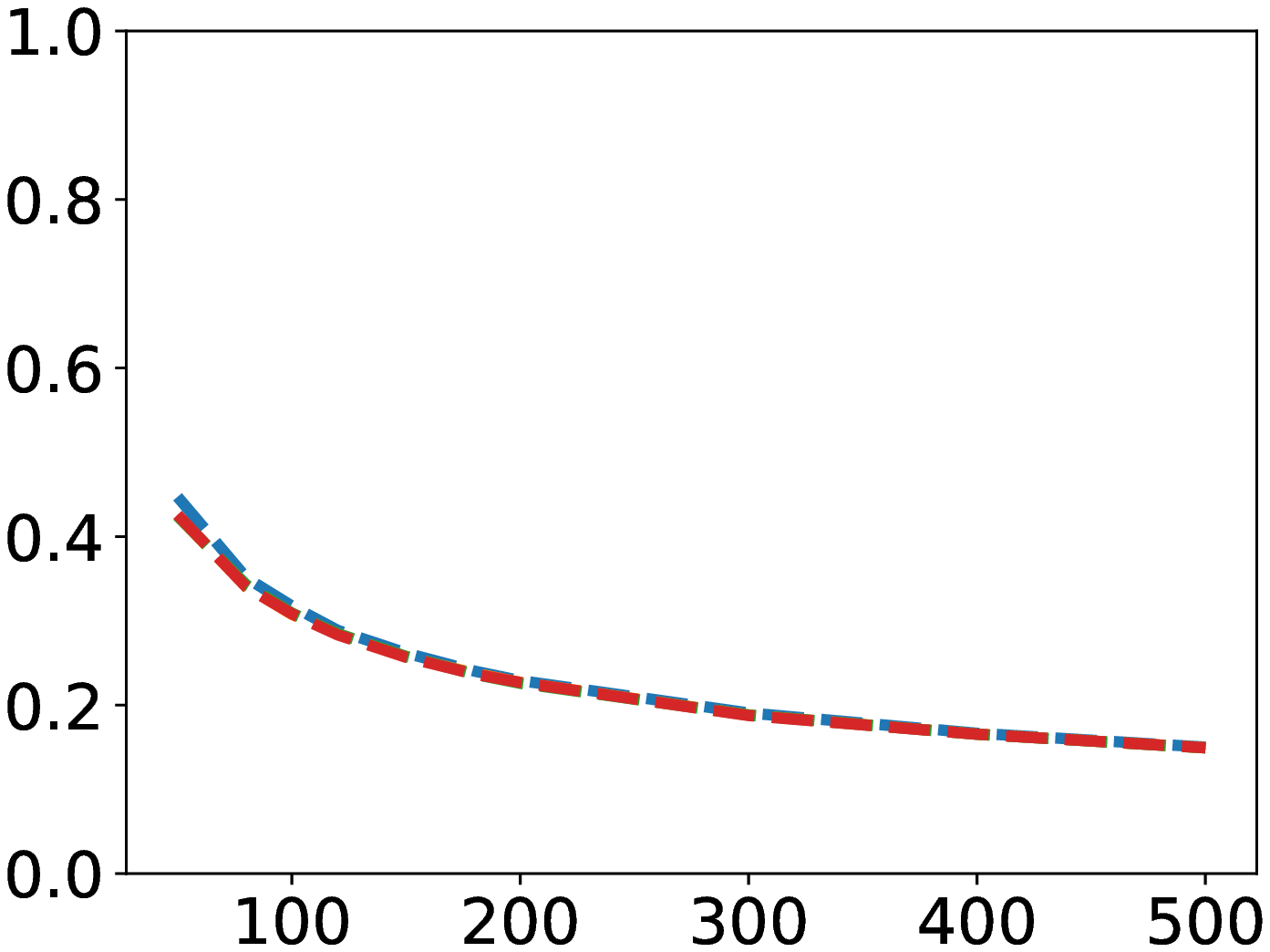} 
	\put(25,80){ \ul{\ \ \  \ $\delta=0.8$ \ \ \ \    }}
	\end{overpic}
	~	
	\begin{overpic}[width=0.29\textwidth]{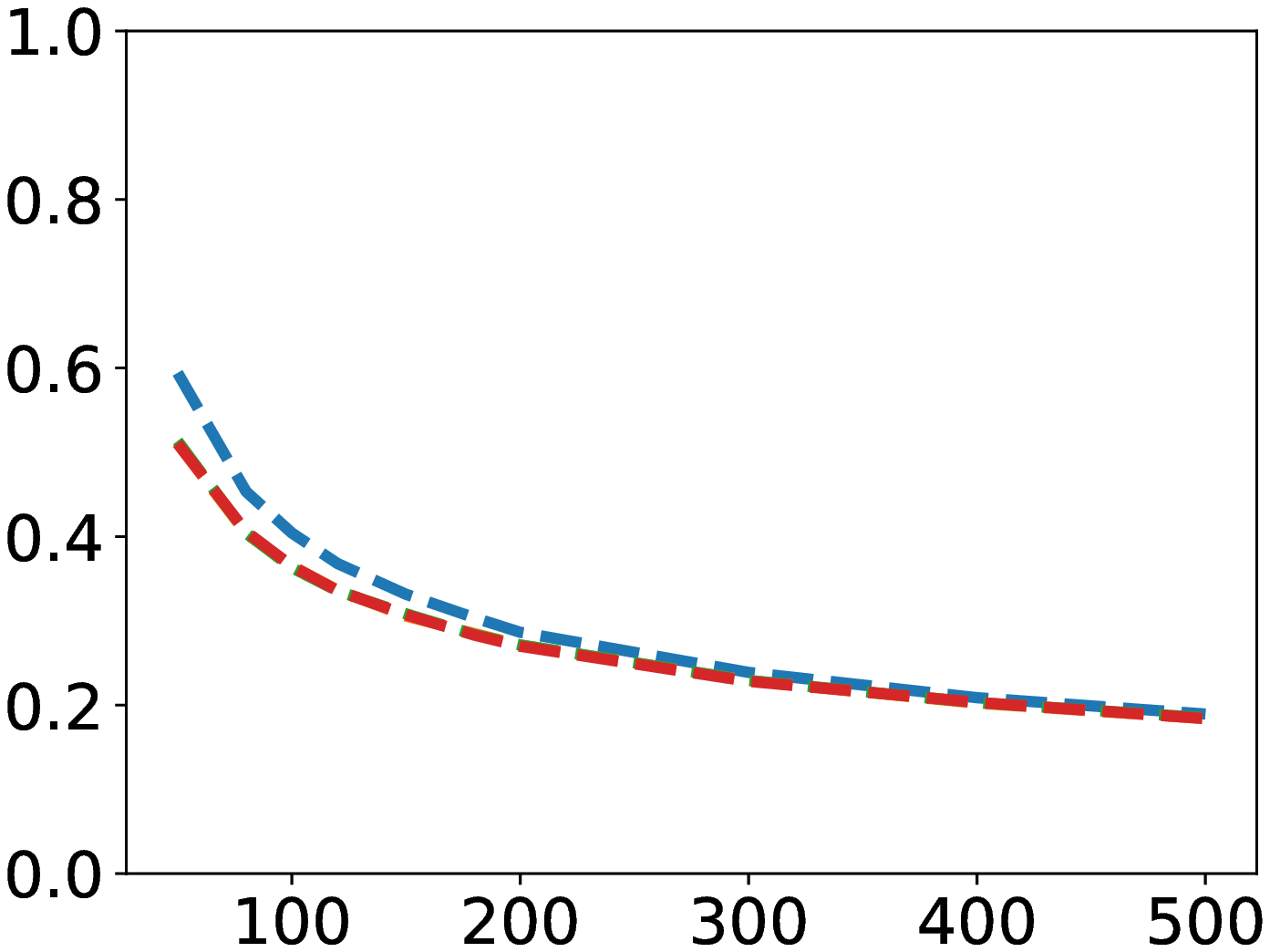} 
	\put(25,80){ \ul{\ \ \  \ $\delta=0.9$ \ \ \ \    }}
	\end{overpic}	
\end{figure}

\vspace{-0.5cm}

\begin{figure}[H]	
	\quad\quad\quad 
	\begin{overpic}[width=0.29\textwidth]{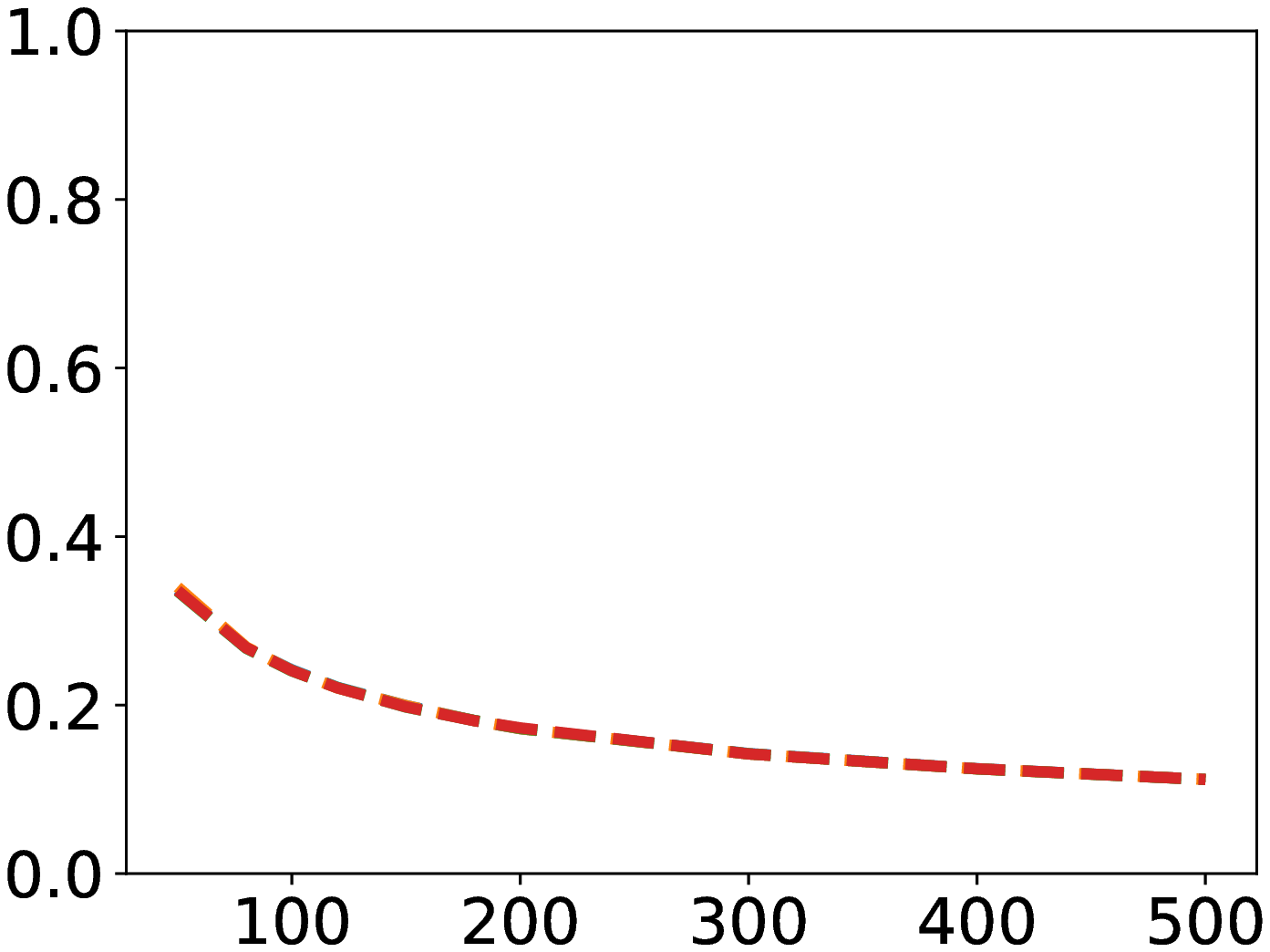} 
	\put(-20,-1){\rotatebox{90}{ {\small \ \ \ standardization \  \ \ }}}
	\end{overpic}
	~
	\DeclareGraphicsExtensions{.png}
	\begin{overpic}[width=0.29\textwidth]{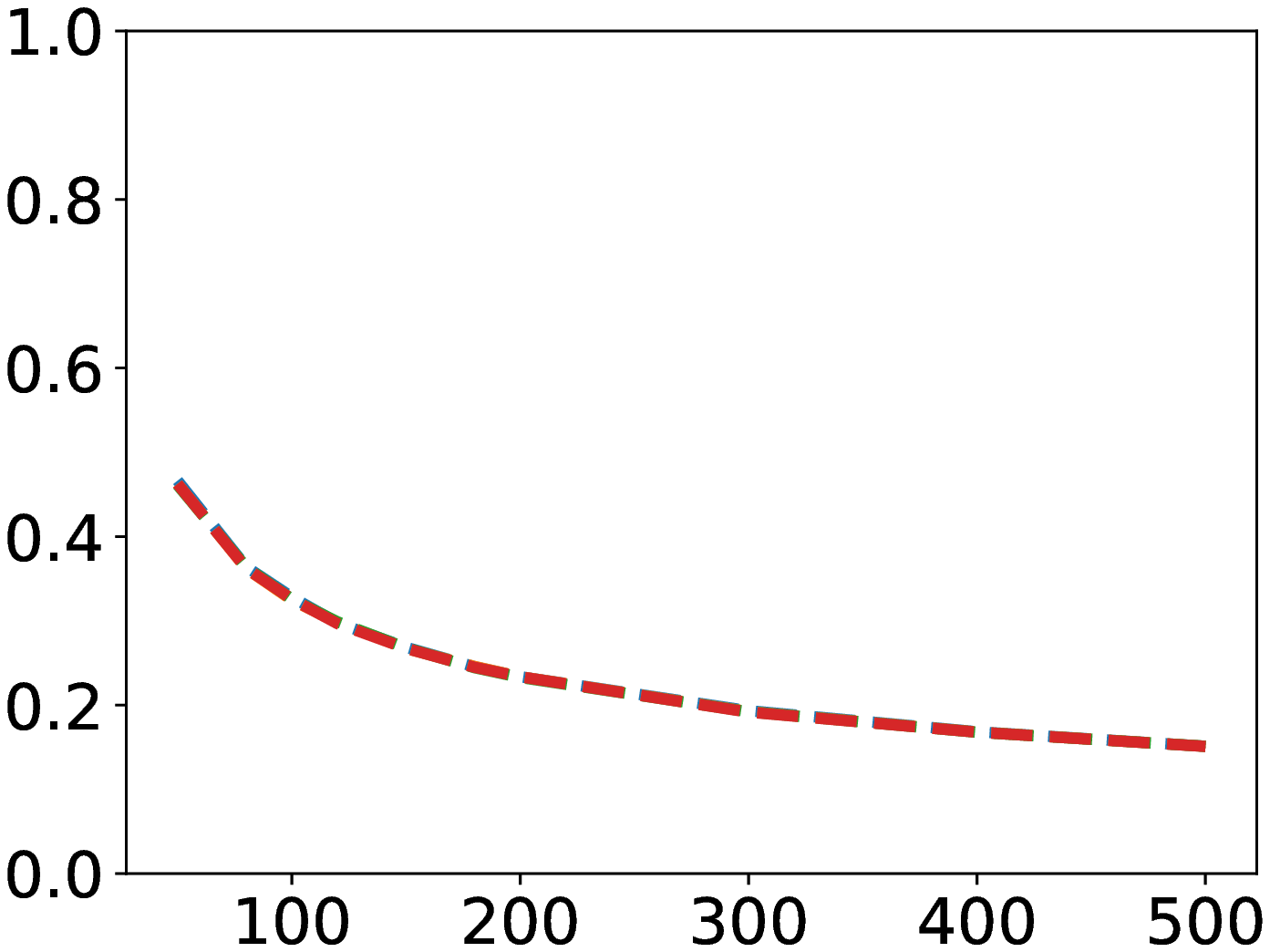} 
	\end{overpic}
	~	
	\begin{overpic}[width=0.29\textwidth]{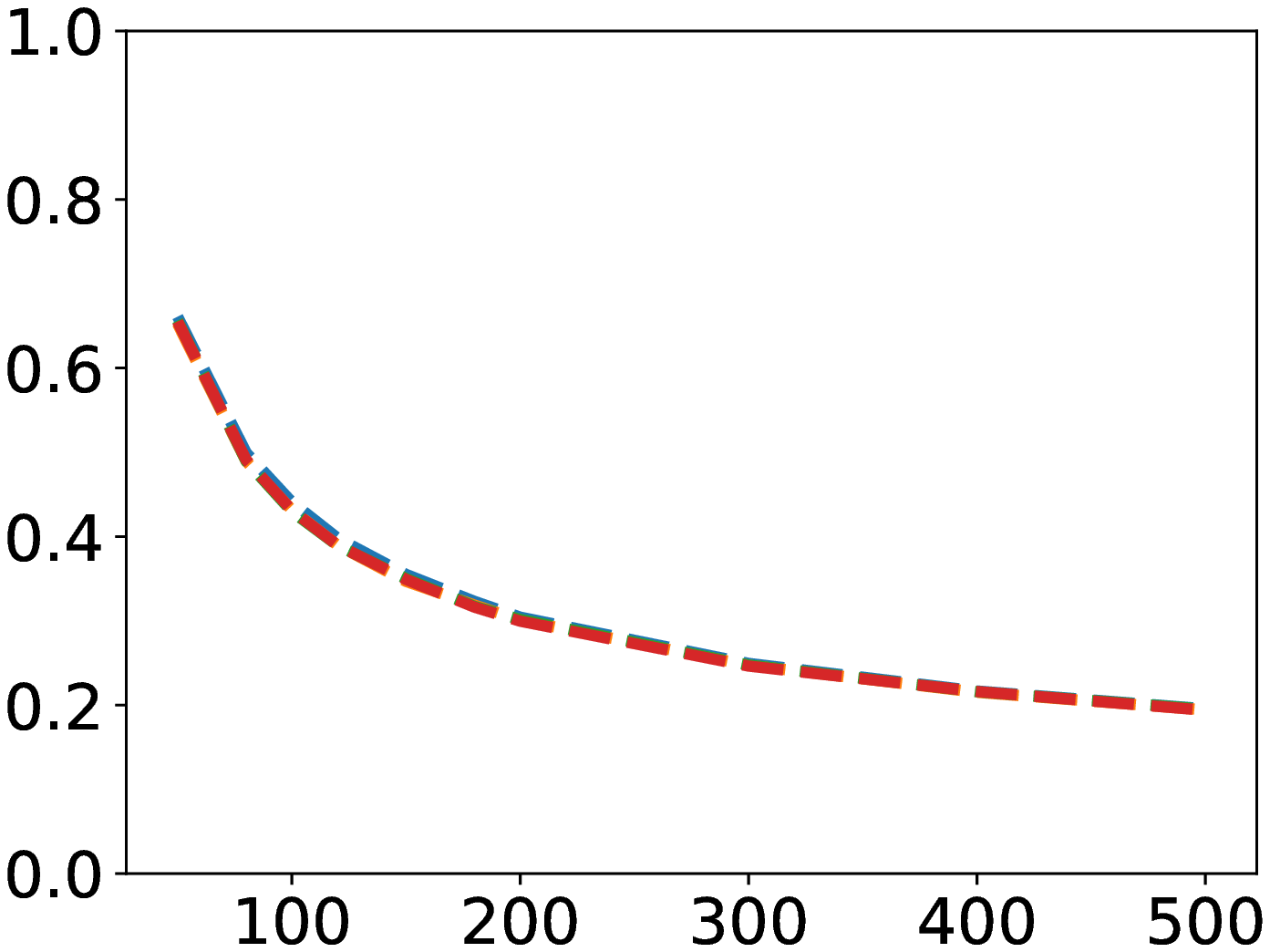} 
	\end{overpic}	
\end{figure}

\vspace{-0.5cm}

\begin{figure}[H]	
	\quad\quad\quad 
	\begin{overpic}[width=0.29\textwidth]{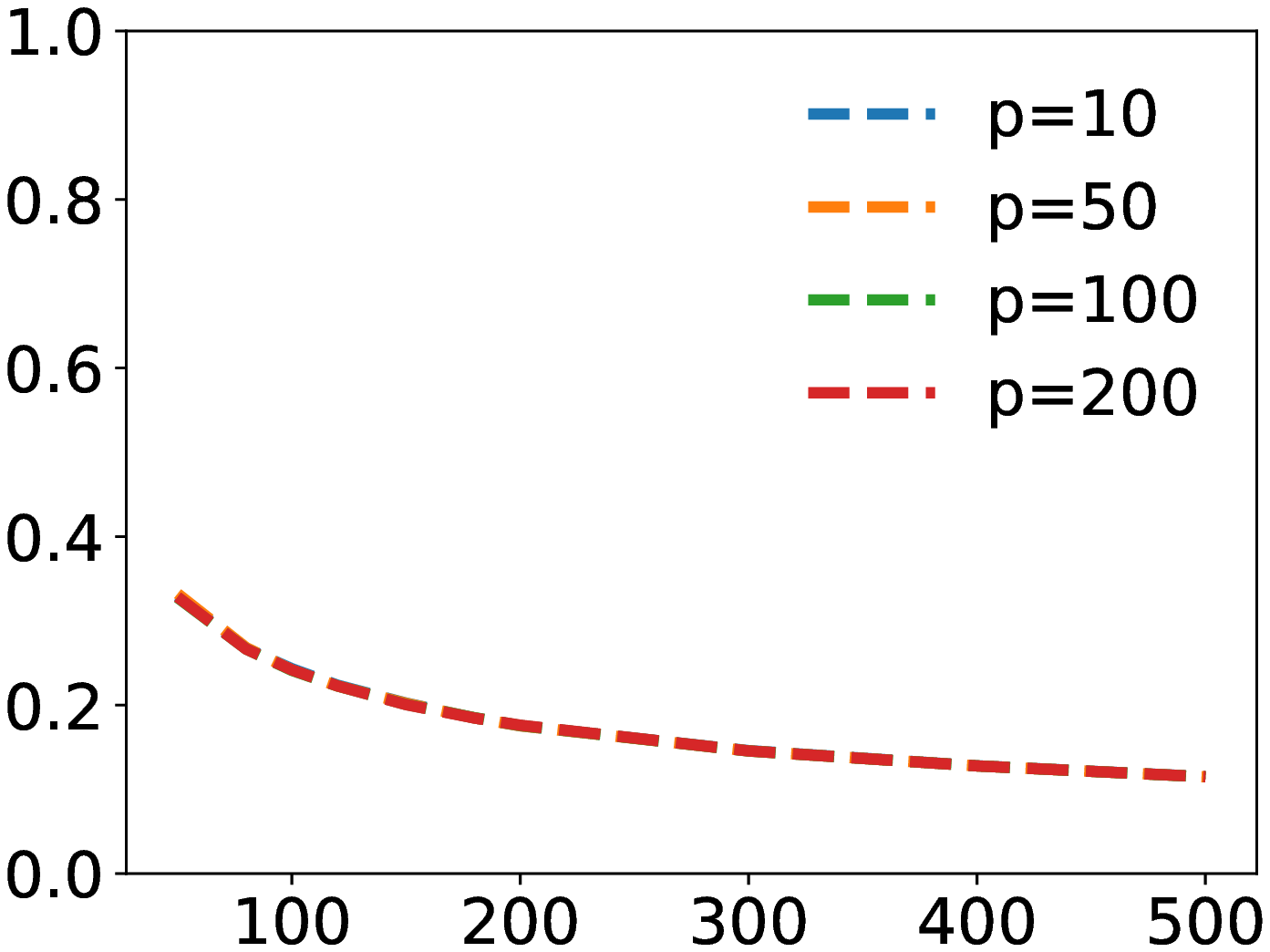} 
	\put(-21,1){\rotatebox{90}{\ $\sqrt{ \ \ }$}}
	\put(-20,-3){\rotatebox{90}{  { \ \ \ \ \ \ \small transformation \ \ } }}

	\end{overpic}
	~
	\DeclareGraphicsExtensions{.png}
	\begin{overpic}[width=0.29\textwidth]{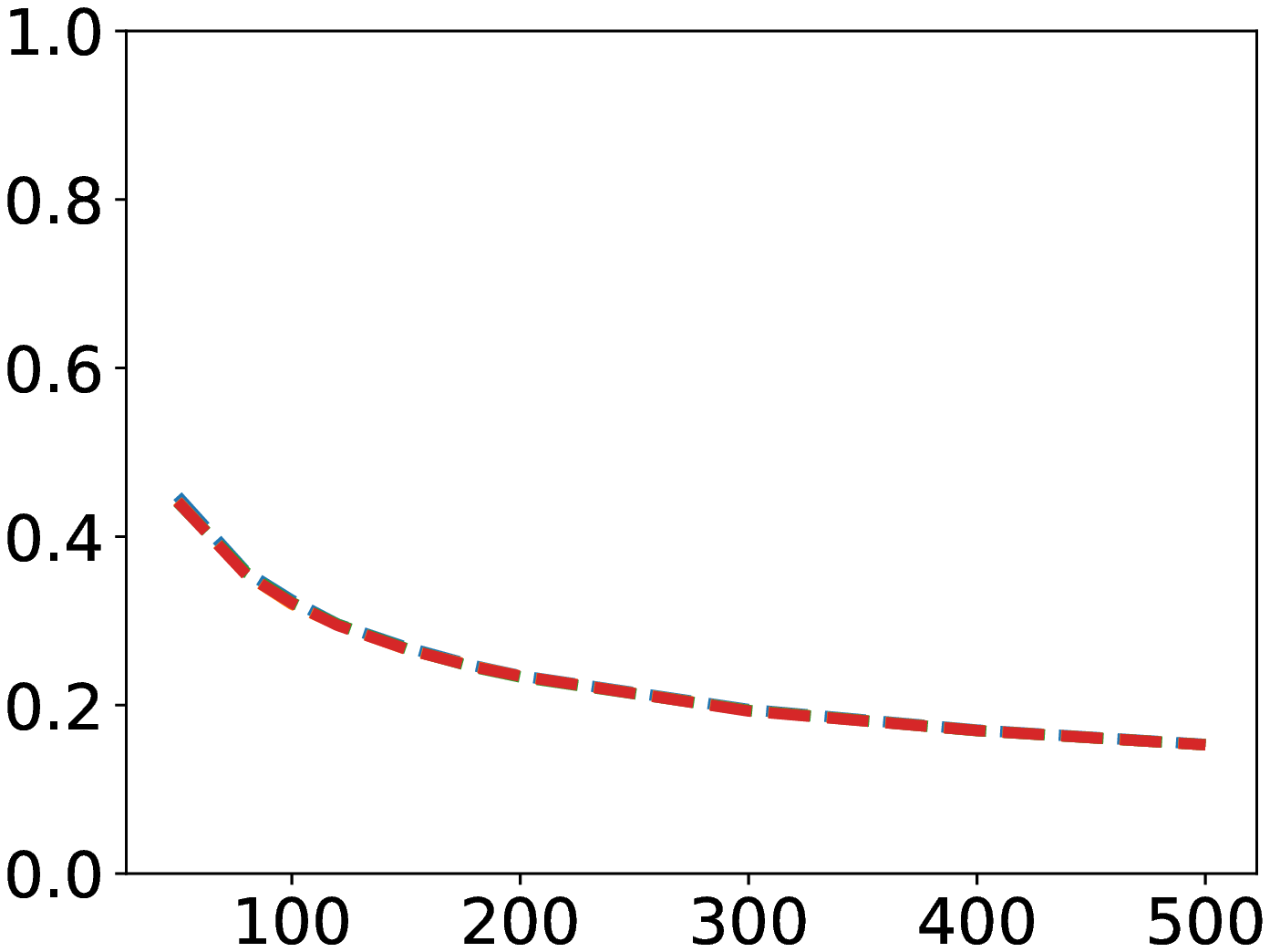} 
	\end{overpic}
	~	
	\begin{overpic}[width=0.29\textwidth]{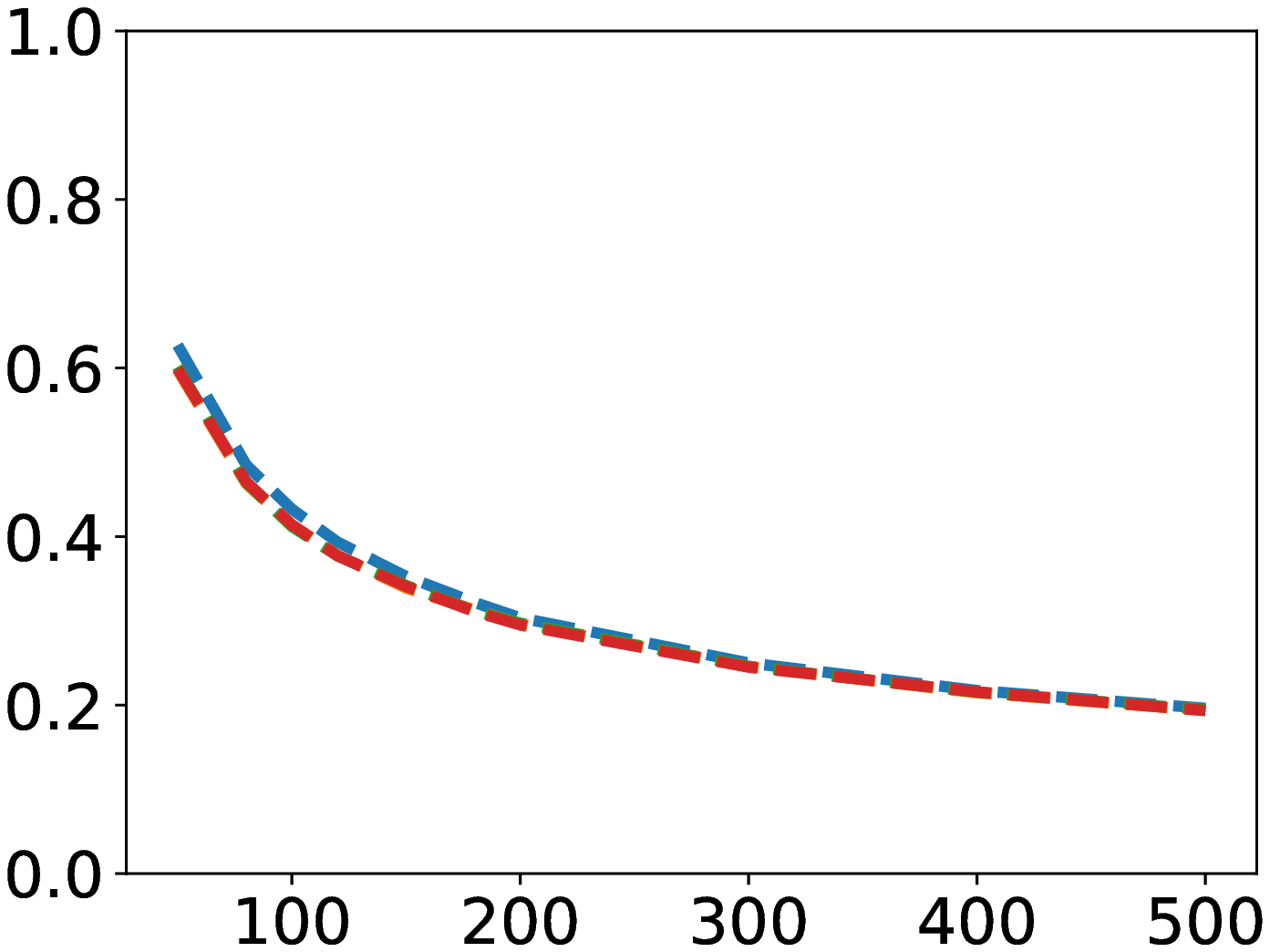} 
		 				
	\end{overpic}	
	\vspace{+0.2cm}	
	\caption{(Average width versus $n$ in simulation model (ii) with an exponential decay profile). The plotting scheme is the same as described in the caption of Figure~\ref{SUPP:fig7} above.}
	\label{SUPP:fig8}
\end{figure}

\newpage

\section{Illustration with stock market data}\label{sec:stocks}

Within the context of finance, PCA is often applied to stock market return data for the purposes of risk analysis and portfolio selection~\citep{Fabozzi,Ruppert:2015}. Here, we look at several high-dimensional datasets of stock market returns to illustrate how the bootstrap can be applied to do inference on parameters of interest in PCA.

Starting from one dataset of S\&P 500 returns during the period February 2013 to December 2017 \citep{kaggle17}, we isolated four distinct datasets in the following way. First, we ranked the 500 stocks based on their average monthly trading volume over the stated time period. Second, we selected four subsets of the 500 stocks, corresponding to the top 50, 150, 200, and 300 members of the ranked list. Third, for each stock, we extracted its biweekly log returns over the time period, resulting in $118$ log return values per stock. (The use of log returns rather than ordinary returns is a standard practice in finance~\citep{Ruppert:2015}.) Altogether, this produced four data matrices of size $n\times p$ with the same number of rows $n=118$, but differing numbers of columns $p = 50, 150, 200, 300$.
In addition to being high-dimensional, these datasets also conform with our interest in settings that are well suited to PCA, since the empirical effective rank satisfies ${\tt{r}}(\hat\Sigma)\leq 4$ for every dataset.

\subsection{Inference on population eigenvalues}
When PCA is used to analyze stock market returns, the leading eigenvectors and eigenvalues of the population covariance matrix $\Sigma$ have special interpretations. Namely, the eigenvector corresponding to $\lambda_1(\Sigma)$ is often viewed as representing an overall ``market portfolio'', while subsequent eigenvectors represent ``principal portfolios'', which produce returns that are uncorrelated with the overall market return~\citep{Laloux}. Also, the eigenvalues can be interpreted as the variances (or volatilities) of the returns  associated with the principal portfolios. For this reason, the population eigenvalues are important for risk assessment, and so it is of interest to quantify the uncertainty in these unknown parameters. 

For each of the four datasets described above, we applied the bootstrap method with square-root transformation from Section~\ref{sec:bootci} to construct simultaneous confidence intervals for the leading ten eigenvalues $\lambda_1(\Sigma),\dots,\lambda_{10}(\Sigma)$.
The bootstrap intervals are plotted in Figure~\ref{fig9}, based on a simultaneous coverage probability of 95\%, with a black dot representing the sample eigenvalue $\lambda_j(\hat\Sigma)$ in the $j$th interval for $j=1,\dots10$.
Upon close inspection, it can be seen that $\lambda_j(\hat\Sigma)$ tends to sit slightly above the midpoint of the $j$th interval. This is encouraging, because it means that the bootstrap intervals are able to counteract the well-known phenomenon that the leading sample eigenvalues tend to be biased upwards in high-dimensional settings~\citep[][Ch.11]{Yao:2015}. In addition, as a way to gain extra empirical support for the bootstrap intervals, we carried out the following exercise with estimates of $\lambda_1(\Sigma),\dots,\lambda_{10}(\Sigma)$ computed via the method of QuEST~\citep{LedoitW15}, which is designed for use in high-dimensional settings, and has been  adopted frequently in the literature.  Specifically, we verified that the QuEST estimate of $\lambda_j(\Sigma)$ was contained in the $j$th bootstrap interval for every $j=1,\dots,10$ and $p=50,150,200,300$. Hence, this makes it more plausible that the bootstrap intervals also contain the population eigenvalues.

To comment further on the numerical results in Figure~\ref{fig9}, first note that in every panel, the interval for $\lambda_1(\Sigma)$ is well separated from the intervals for $\lambda_2(\Sigma),\dots,\lambda_{10}(\Sigma)$, while there is substantial overlap among the latter intervals.  This type of situation occurs frequently when PCA is applied to stock market return data, and this is generally interpreted to mean that the overall behavior of the market has a much more dominant effect on returns than other types of economic factors~\cite[][]{Laloux,Ruppert:2015}. A second observation is that the bootstrap intervals can provide some additional insight into the relationship between  $\lambda_2(\Sigma)$ and $\lambda_3(\Sigma)$. On one hand, a user who only looks at the sample eigenvalues $\lambda_2(\hat\Sigma)$ and $\lambda_3(\hat \Sigma)$ might be tempted to conclude that there is a clear difference between the population eigenvalues $\lambda_2(\Sigma)$ and $\lambda_3(\Sigma)$. On the other hand, a user who looks at the overlap of the second and third intervals would have more information to see that the difference between $\lambda_2(\Sigma)$ and $\lambda_3(\Sigma)$ might actually be negligible. Lastly, one more aspect of Figure~\ref{fig9} to mention is that the relative positions of the ten intervals stay approximately the same for each of the four dimensions $p=50,150,200,300$. Given that the empirical effective rank satisfies ${\tt{r}}(\hat\Sigma)\leq 4\ll p$ for every dataset, this makes sense from the standpoint of our theoretical results, which indicate that the bootstrap should be relatively insensitive to the ambient dimension compared to the effective rank.

\begin{figure}[H]	
\vspace{0.5cm}
	\quad\quad 
	\begin{overpic}[width=0.40\textwidth]{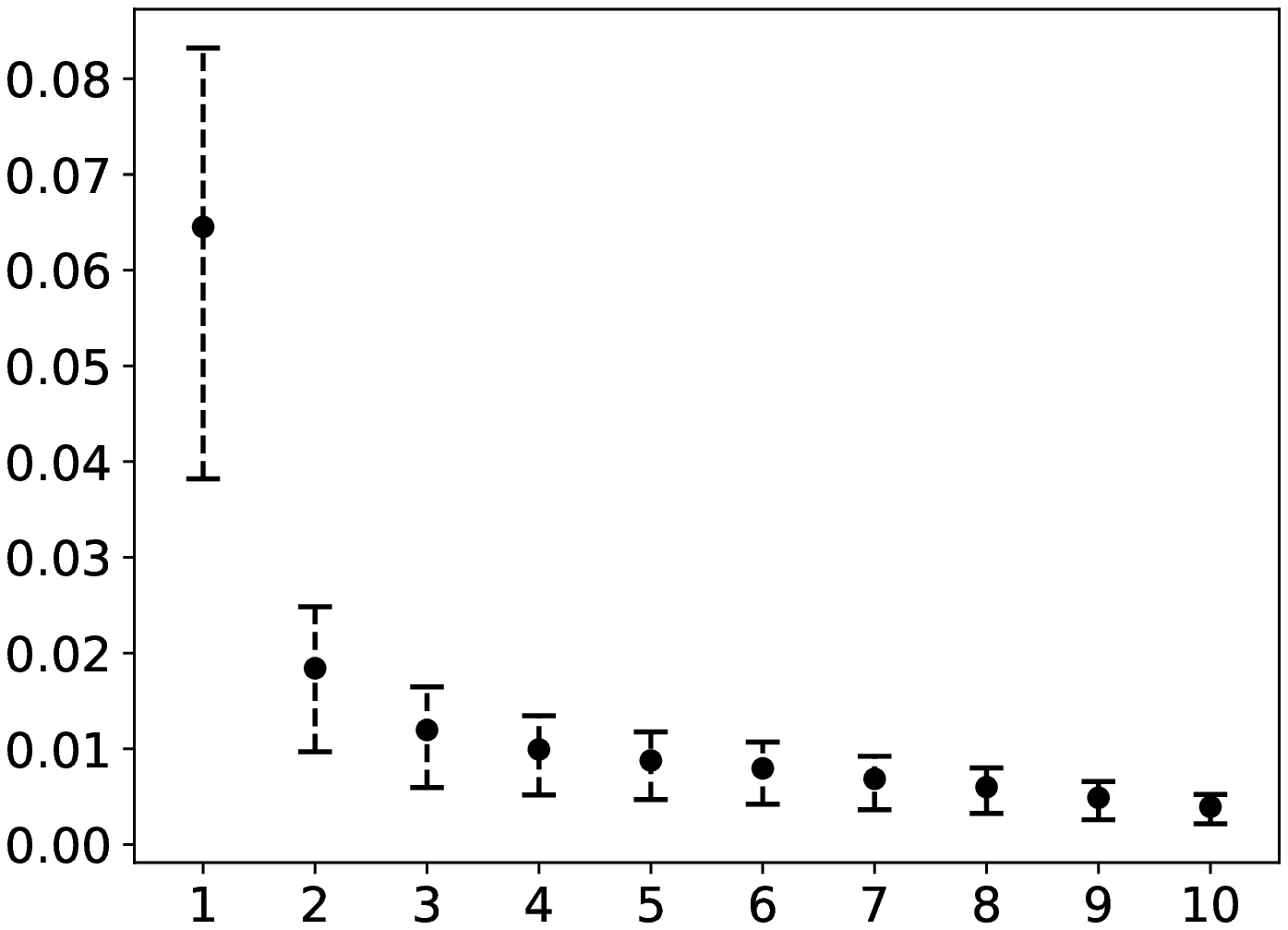} 

	\put(31,80){ \ul{\ \ \  \ $p = 50$ \ \ \ \    }}
\end{overpic}
	~
	\DeclareGraphicsExtensions{.png}
	\begin{overpic}[width=0.40\textwidth]{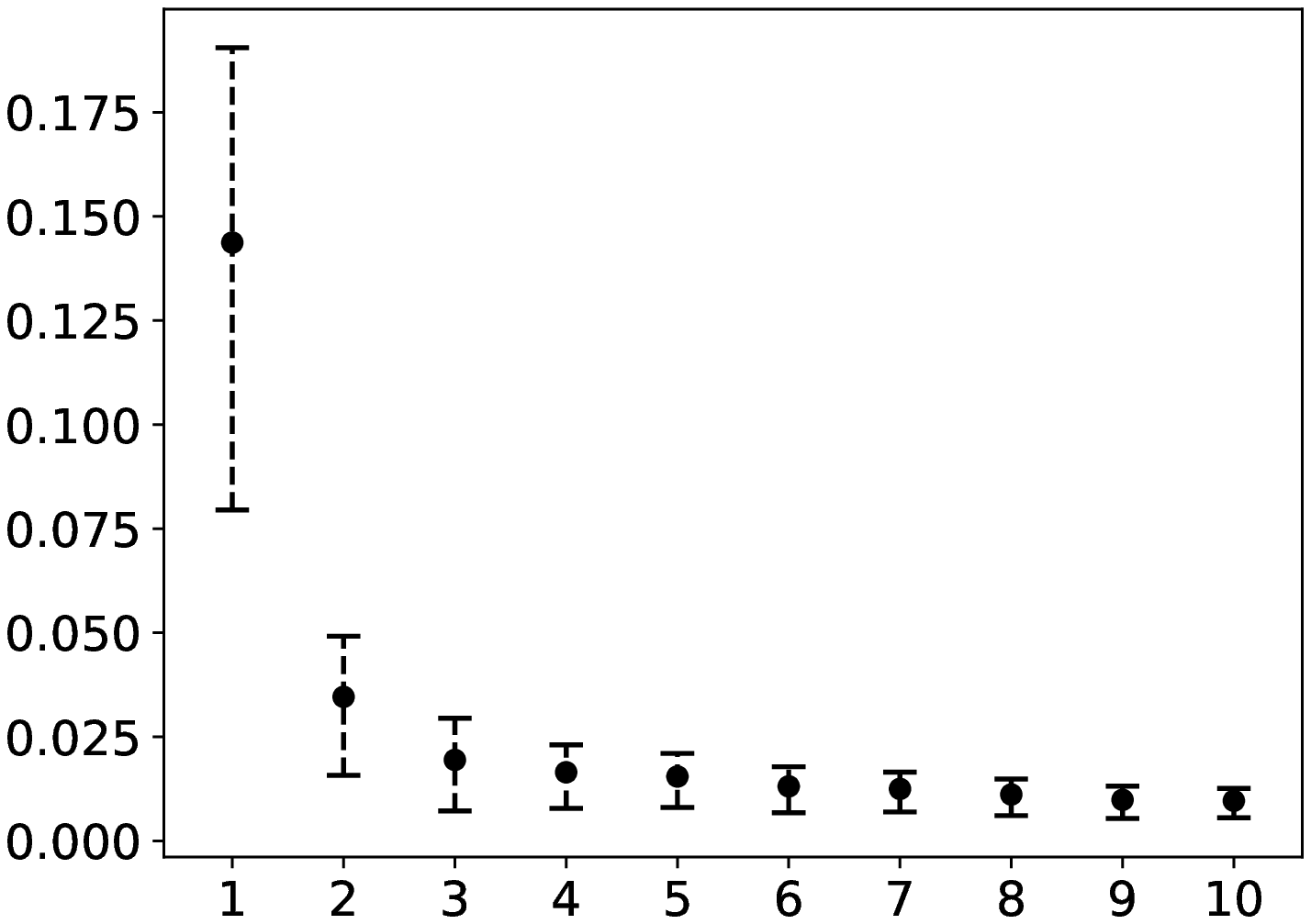} 
		\put(31,80){ \ul{\ \ \  \ $p = 150$ \ \ \ \    }}
	\end{overpic}
\end{figure}

\begin{figure}[H]	
	\quad\quad
	\begin{overpic}[width=0.40\textwidth]{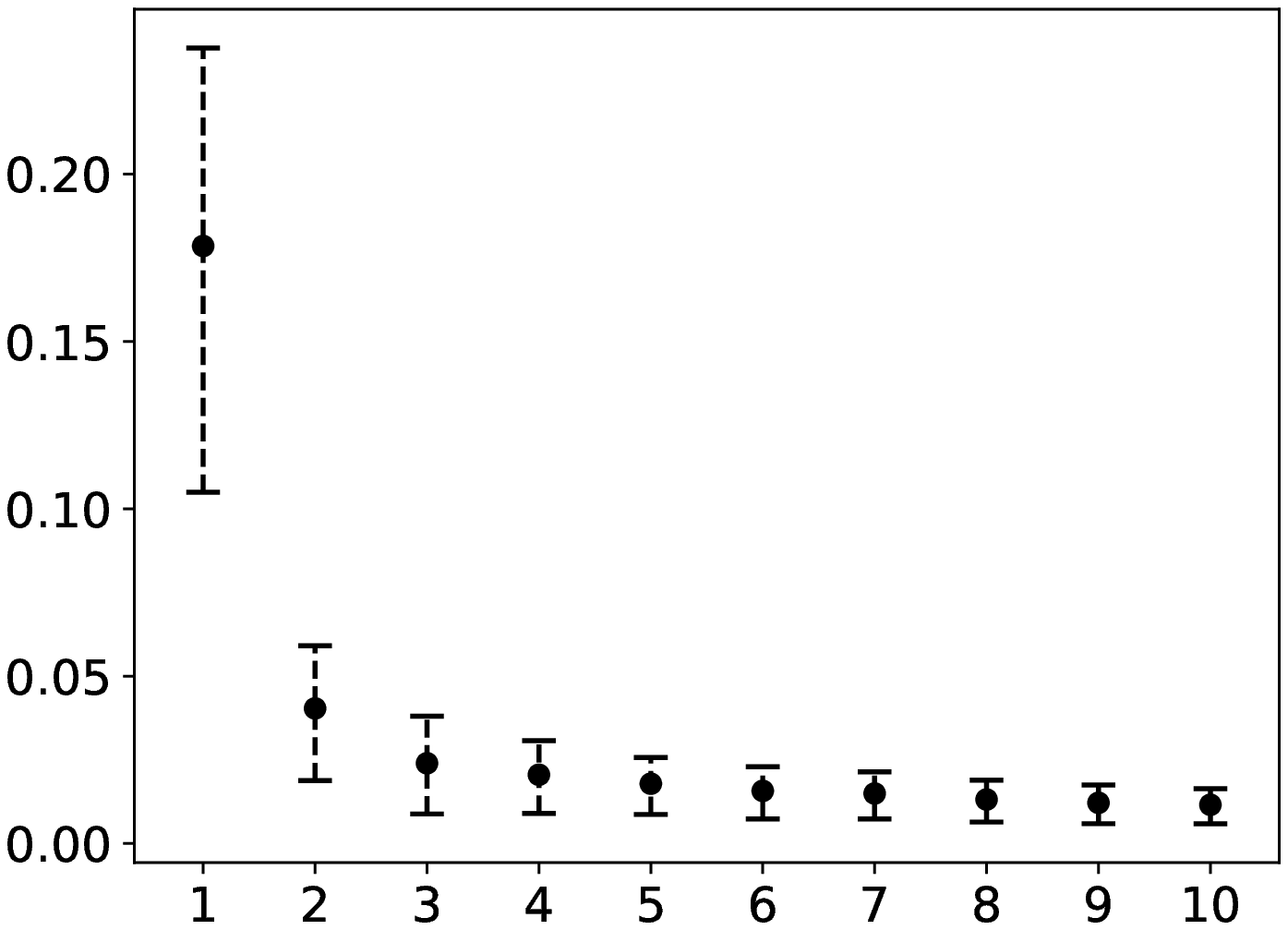} 
	\put(31,80){ \ul{\ \ \  \ $p = 200$ \ \ \ \    }}
	\end{overpic}
	~
	\DeclareGraphicsExtensions{.png}
	\begin{overpic}[width=0.40\textwidth]{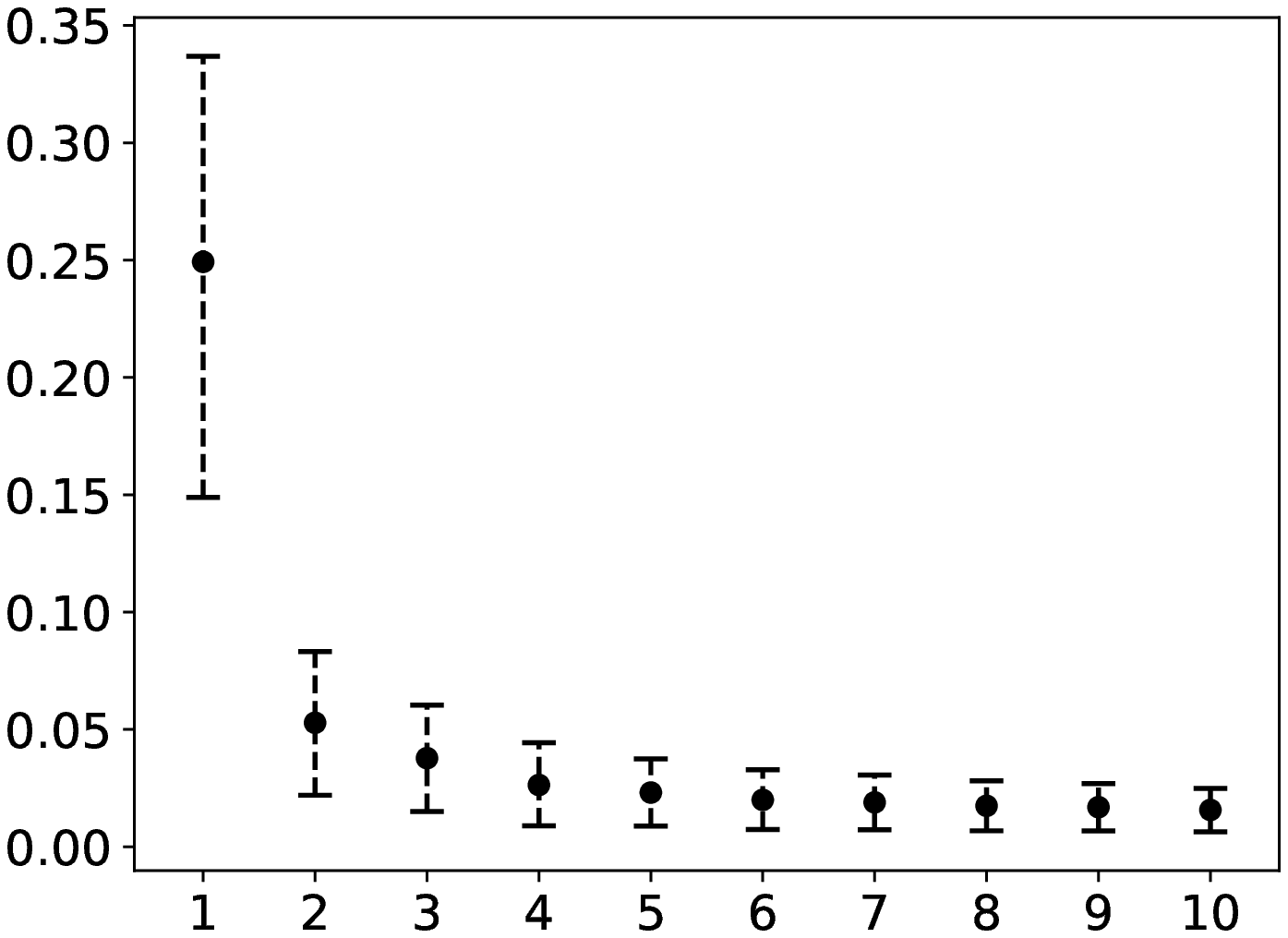} 
	\put(31,80){ \ul{\ \ \  \ $p = 300$ \ \ \ \    }}
	\end{overpic}
	\vspace{-0.2cm}	
	\caption{(Simultaneous confidence intervals for $\lambda_1(\Sigma),\dots,\lambda_{10}(\Sigma)$.) In each panel, the $y$-axis corresponds to the magnitude of eigenvalues, and the $x$-axis corresponds to the index $j=1,\dots,10$. The intervals are based on a simultaneous coverage probability of 95\%, and the black dots represent the sample eigenvalues $\lambda_j(\hat\Sigma)$ for each $j$.}
	\label{fig9}
\end{figure}

\subsection{Inference on proportions of explained variance}

The proportions of explained variance, denoted $\pi_j(\Sigma)=\sum_{i=1}^j \lambda_i(\Sigma) / \tr(\Sigma)$ for $j=1,\dots,p$, often play a decisive role in applications of PCA, since they form of the basis of standard decision rules for selecting an appropriate number of components.
To complement our previous example dealing with the population eigenvalues $\lambda_1(\Sigma),\dots,\lambda_{10}(\Sigma)$, this subsection looks instead at inference with simultaneous confidence intervals for the parameters $\pi_1(\Sigma),\dots,\pi_{10}(\Sigma)$. As before, the bootstrap method with square-root transformation from Section~\ref{sec:bootci} was used to construct the intervals based on a simultaneous coverage probability of 95\%. The results are given in Figure~\ref{fig10}, with black dots showing the locations of the empirical proportions $\pi_j(\hat\Sigma)$ within the $j$th interval. %

Due to the fact that the proportions $\pi_1(\Sigma),\dots,\pi_p(\Sigma)$ are unknown, one of the most widely used rules for selecting the number of components is to choose the smallest number $k$ for which $\pi_k(\hat\Sigma)$ exceeds a given threshold. Although this rule may be appropriate when dealing with low-dimensional data, it is known in the literature that this rule can be unreliable in high-dimensional settings, because it tends to select too few components~\citep{LedoitW15}. One way of avoiding this pitfall is to consider the following simple modification, based on simultaneous bootstrap confidence intervals for the proportions: If $\hat l_j$ denotes the lower endpoint of the $j$th interval, then the bootstrap-based rule selects the smallest number $k$ for which $\hat l_k$ exceeds the threshold. Consequently, if the intervals perform properly with a simultaneous coverage probability of 95\%, then the bootstrap-based rule will select a sufficient number of components with at least 95\% probability.

\begin{figure}[H]	
\vspace{0.5cm}
	\quad\quad 
	\begin{overpic}[width=0.40\textwidth]{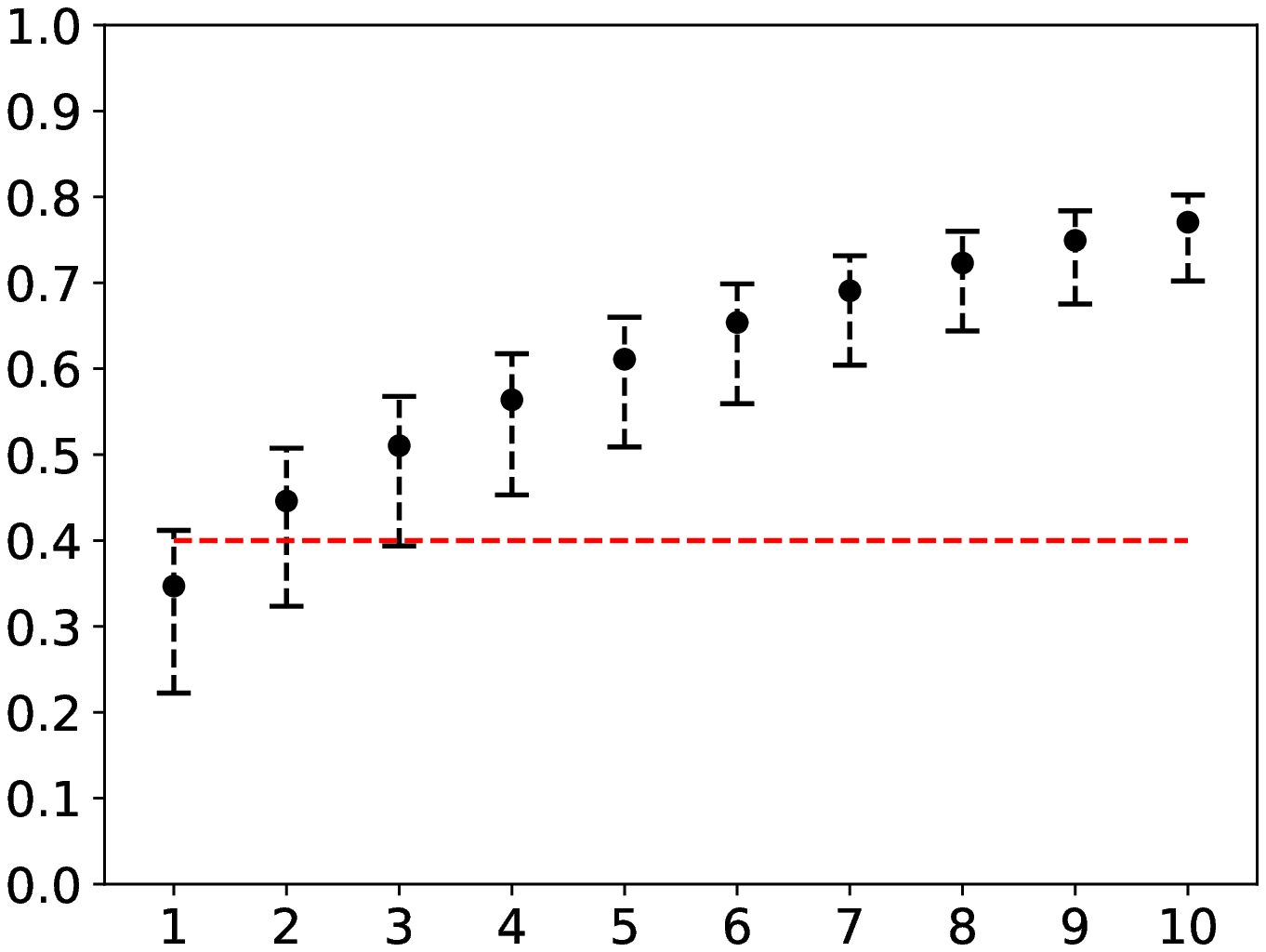} 

	\put(31,80){ \ul{\ \ \  \ $p = 50$ \ \ \ \    }}
\end{overpic}
	~
	\DeclareGraphicsExtensions{.png}
	\begin{overpic}[width=0.40\textwidth]{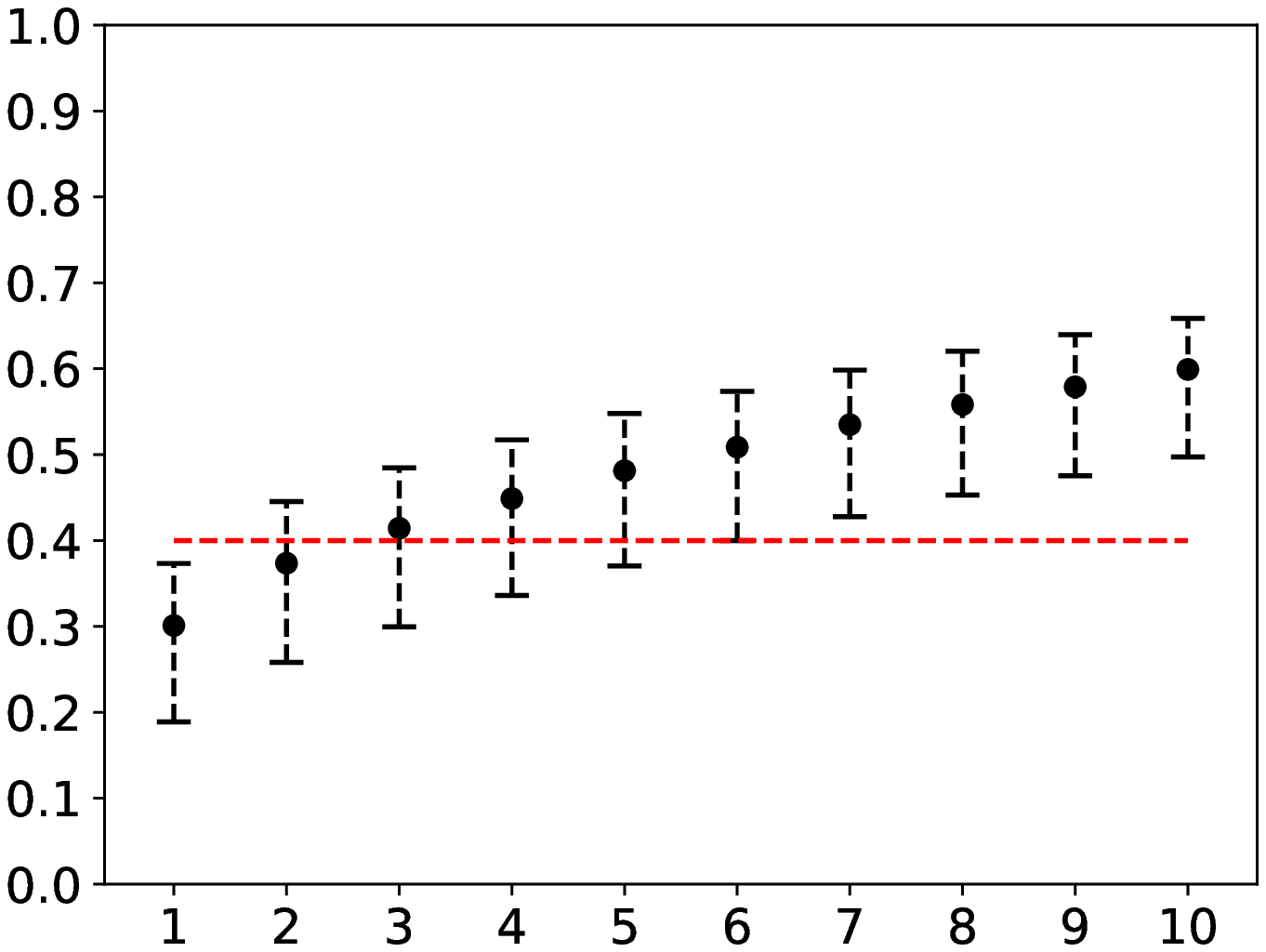} 
		\put(31,80){ \ul{\ \ \  \ $p = 150$ \ \ \ \    }}
	\end{overpic}
\end{figure}

\begin{figure}[H]	
	\quad\quad
	\begin{overpic}[width=0.40\textwidth]{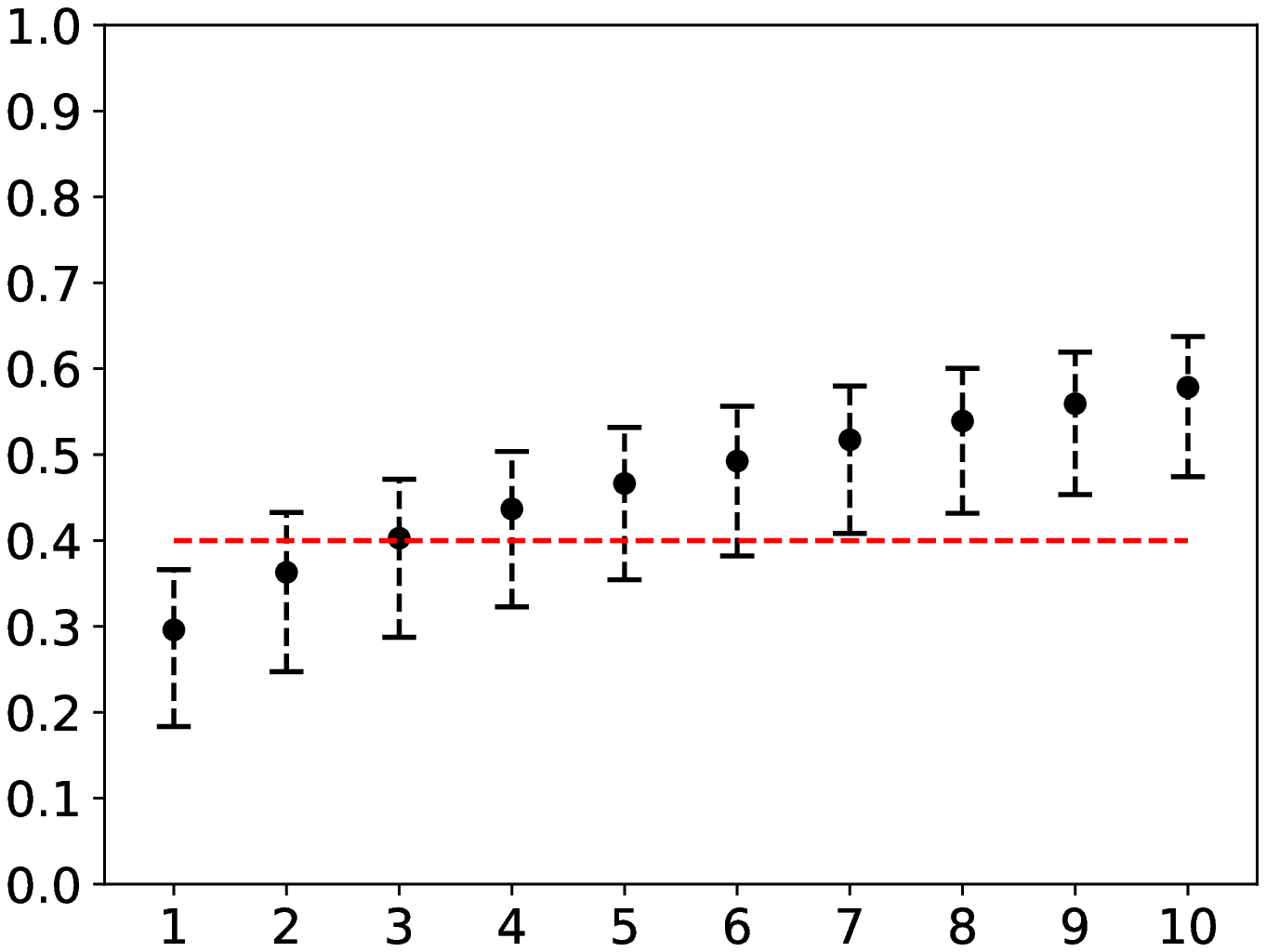} 
	\put(31,80){ \ul{\ \ \  \ $p = 200$ \ \ \ \    }}
	\end{overpic}
	~
	\DeclareGraphicsExtensions{.png}
	\begin{overpic}[width=0.40\textwidth]{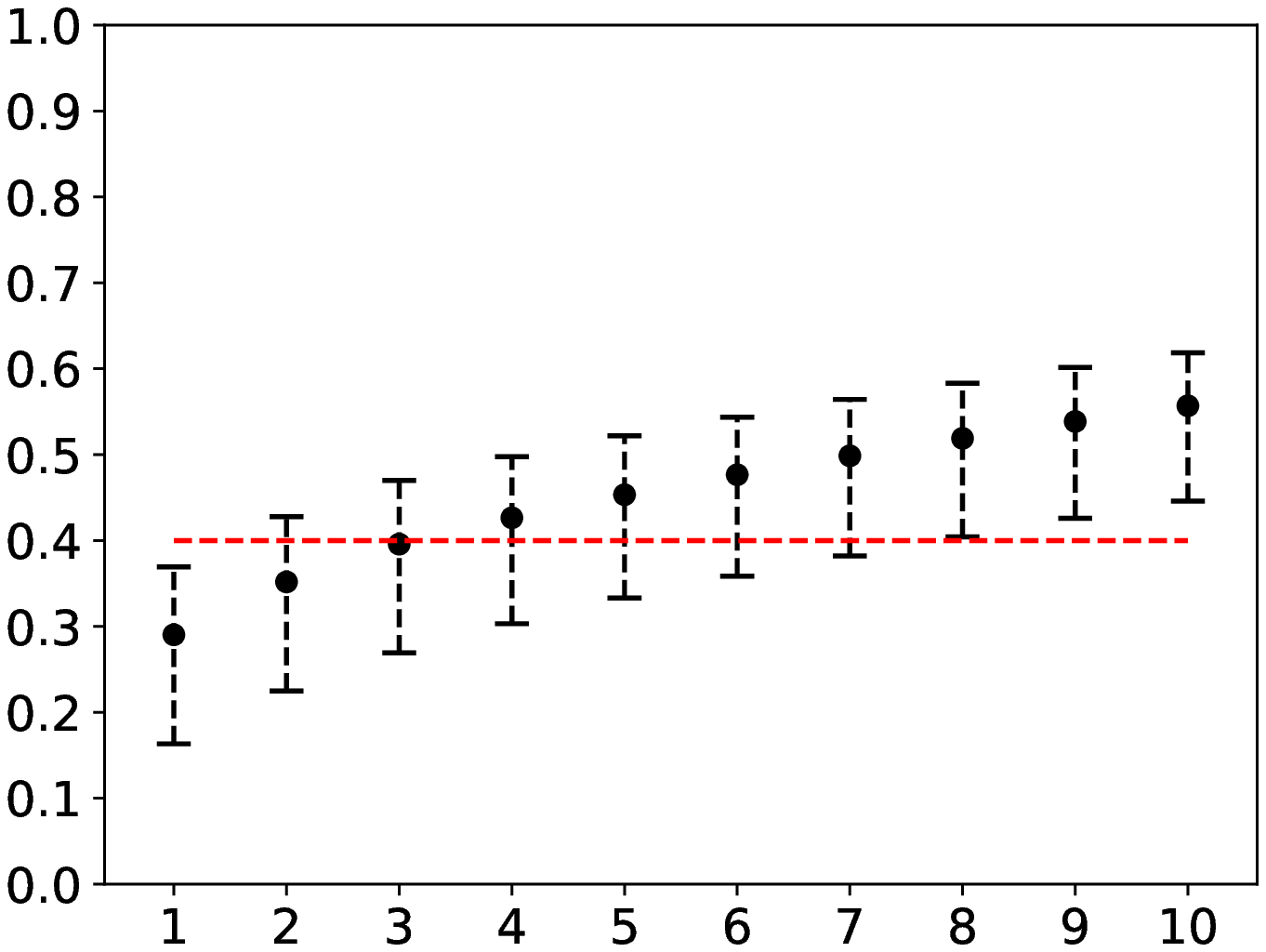} 
	\put(31,80){ \ul{\ \ \  \ $p = 300$ \ \ \ \    }}
	\end{overpic}
	\vspace{-0.2cm}	
	\caption{(Simultaneous confidence intervals for the proportions of explained variance $\pi_1(\Sigma),\dots,\pi_{10}(\Sigma)$.) In each panel, the $y$-axis corresponds to the magnitude of proportions, and the $x$-axis corresponds to the index $j=1,\dots,10$. The intervals are based on a simultaneous coverage probability of 95\%, and the black dots represent the empirical proportions $\pi_j(\hat\Sigma)$ for each $j$.}

	\label{fig10}
\end{figure}

To illustrate the difference between the two rules, a red line corresponding to a particular threshold of 0.4 has been drawn in each panel of Figure~\ref{fig10}. (The value of 0.4 has no special importance, and is used only for ease of presentation.) In each of the four cases $p=50,150,200,300$, the original rule selects the respective values $k=2, 3, 3, 4$. By contrast, the bootstrap-based rule selects the respective values $k=4, 6, 7, 8$, and hence, it clearly counteracts the problem of selecting too few components. Moreover, from a financial standpoint, it makes sense that extra components are needed as $p$ increases, because as a wider variety of stocks are included in the data, there is greater opportunity for the returns to be influenced by economic factors that are not captured by the previously leading components.

\section{Computational cost}\label{app:computation}
The cost to compute a single bootstrap sample of $\boldsymbol \lambda_k(\hat\Sigma^{\star})-\boldsymbol\lambda_k(\hat\Sigma)$ can be broken into two steps. The first step is to sample $n$ points with replacement from the original set of $n$ observations, which has a cost of $O(n\log(n))$. The second step consists of computing the largest $k$ eigenvalues of $\hat\Sigma^{\star}$. This can be done by computing the largest $k$ singular values of the $n\times p$ matrix of resampled observations, which has a cost of $\mathcal{O}(npk)$~\citep{Halko}. (Note that the largest $k$ eigenvalues of $\hat \Sigma$ only need to be computed once, before any resampling is done.) So, the cost to compute $B$ bootstrap samples of $\boldsymbol \lambda_k(\hat\Sigma^{\star})-\boldsymbol\lambda_k(\hat\Sigma)$ on a single processor is $\mathcal{O}(B(n\log(n)+npk))$. However, it is common to compute bootstrap samples in a parallel manner, across say $m$ processors, and in this case the cost per processor becomes $\mathcal{O}(\frac{B}{m}(n\log(n)+npk)$. In order to simplify this expression, it is natural to consider a scenario where $B/m=\mathcal{O}(1)$ and $\log(n)=\mathcal{O}(p)$, which leads to a cost per processor that is $\mathcal{O}(npk)$.

\subsection{Empirical computational cost}

Table~\ref{SUPP:num_exp_cost} displays the time, in seconds, to compute a single bootstrap sample of $\boldsymbol \lambda_k(\hat\Sigma^{\star})-\boldsymbol\lambda_k(\hat\Sigma)$ for different choices of $n$ and $p$. Each entry reflects an average over $1000$ trials with data generated under simulation model (i). The computations were done on a single Intel Xeon E5-2699v3 processor. Notably, even when $(n,p)=(500,200)$, the computing time for each bootstrap sample is on the order of just $10^{-3}$ seconds. Hence, it is possible to generate hundreds bootstrap samples within about 1 second.
\begin{table}[H]	
\vspace{-.2cm}
\centering
\begin{tabular}{c c c c c} 
& \multicolumn{4}{c}{$p$}\\ 
\cmidrule{2-5} 
$n$   & 10 & 50 & 100 & 200 \\
\cmidrule{1-5} 
10  & 1.9e-04 & 4.2e-04 & 8.7e-04 & 2.6e-03\\
50  & 2.0e-04 & 4.4e-04 & 9.0e-04 & 2.8e-03 \\
100 & 2.1e-04 & 4.6e-04 & 9.6e-04 & 2.9e-03 \\
200 & 2.2e-04 & 5.1e-04 & 1.1e-03 & 3.0e-03 \\
500 & 2.4e-04 & 6.5e-04 & 1.5e-03 & 4.0e-03 \\
\bottomrule
\end{tabular}
\vspace{+.1cm}
\caption{(Average time, in seconds, to compute one bootstrap sample from simulated data.) The rows correspond to the values of $n= 10, 50, 100, 200, 500$, and the columns correspond to the values of $p= 10, 50, 100, 200$.}
\label{SUPP:num_exp_cost}
\end{table}

Table~\ref{SUPP:real_data_time} displays analogous computing times for $\boldsymbol \lambda_k(\hat\Sigma^{\star})-\boldsymbol\lambda_k(\hat\Sigma)$, in seconds, based on the stock market data.  Here, the computations were done on a single 3.5 GHz Dual-Core Intel Core i7 processor. The results show that there is little difference in computing time compared to the setting of synthetic data in  Table~\ref{SUPP:num_exp_cost}.

\begin{table}[H]
\centering
\begin{tabular}{c  c c c c} 
 & \multicolumn{4}{c}{$p$}\\ 
\cmidrule{2-5} 
$n$ & 50 & 150 & 200 & 300 \\
\cmidrule{1-5} 
118 & 3.8e-04 & 1.4e-03 & 2.0e-03 & 4.0e-03 \\
\bottomrule
\end{tabular}
\vspace{+.1cm}
\caption{(Average time, in seconds, to compute one bootstrap sample from the stock market data.) The single row corresponds to the value of $n= 118$, and the columns correspond to the values of $p= 50, 150, 200, 300$.}
\label{SUPP:real_data_time}
\end{table}

\section{Additional discussion on Assumption 1(b)}\label{app:assumption}
Recall that Assumption 1(b) requires the condition $\min_{1\leq j\leq k}(\lambda_j(\Sigma)-\lambda_{j+1}(\Sigma))\gtrsim \lambda_1(\Sigma)$, which ensures that there are  gaps between the leading eigenvalues $\lambda_1(\Sigma),\dots,\lambda_{k+1}(\Sigma)$.
To inspect whether or not this type of condition holds in practice, the paper \citep{Hall:2009} proposes a diagnostic method that constructs a preliminary set of \emph{conservative} simultaneous confidence intervals for $\lambda_1(\Sigma),\dots,\lambda_{k+1}(\Sigma)$. In exchange for their conservatism, these intervals have the property that they are not sensitive to the existence of gaps. (See also~\citep{LopesEM19} for theoretical analysis related to this technique.) If none of the intervals overlap, then the user may conclude that the eigenvalues are adequately separated. However, if some of the intervals do overlap, then the paper~\citep{Hall:2009} recommends that an adjusted form of bootstrapping be used to approximate the distribution of the statistic $\boldsymbol \lambda_k(\hat\Sigma)-\boldsymbol\lambda_k(\Sigma)$.

\subsection{Sensitivity analysis}
To study the sensitivity of the bootstrap to Assumption 1(b), we now discuss some numerical experiments with varying gaps between the leading population eigenvalues.
The first three eigenvalues were specified as $(\lambda_1(\Sigma), \lambda_2(\Sigma), \lambda_3(\Sigma)) = (1+g, 1, 1-g)$ for a gap parameter $g\in\{0,0.1,0.2\}$, and the remaining eigenvalues were chosen to follow the polynomial decay profile $\lambda_j(\Sigma)=j^{-1}$ for $j\geq 4$. Next, we applied the bootstrap (as in Section~\ref{sec:experiments}) to construct simultaneous confidence intervals for $(\lambda_1(\Sigma),\dots,\lambda_5(\Sigma))$ based on a nominal level of 95\%.  Figure~\ref{SUPP:extra_sim_elliptical_sqrtexp_95} contains a grid of plots in which the columns correspond to increasing values of the gap parameter $g$, and the rows correspond to the three transformation rules considered in Section~\ref{sec:experiments}. These plots are based on data generated from simulation model (i), and an analogous set of plots based on simulation model (ii) are given in Figure~\ref{SUPP:extra_sim_iid_gaussian_95}.

As expected, Figures~\ref{SUPP:extra_sim_elliptical_sqrtexp_95} and~\ref{SUPP:extra_sim_iid_gaussian_95} show that the coverage accuracy of the  bootstrap confidence intervals improves as the gap parameter increases. In the case when $g=0$, the coverage accuracy is poor for all three transformation rules, even at large sample sizes. Next, when $g=0.1$, the coverage is mostly accurate for sample sizes $n\geq 400$, but at smaller sample sizes the coverage generally falls below the desired level. Lastly, when $g=0.2$, the coverage becomes more accurate when $n<400$, especially when the square-root transformation is used.

\newpage

\begin{figure}[H]	
\vspace{0.5cm}
	\quad\quad\quad 
	\begin{overpic}[width=0.29\textwidth]{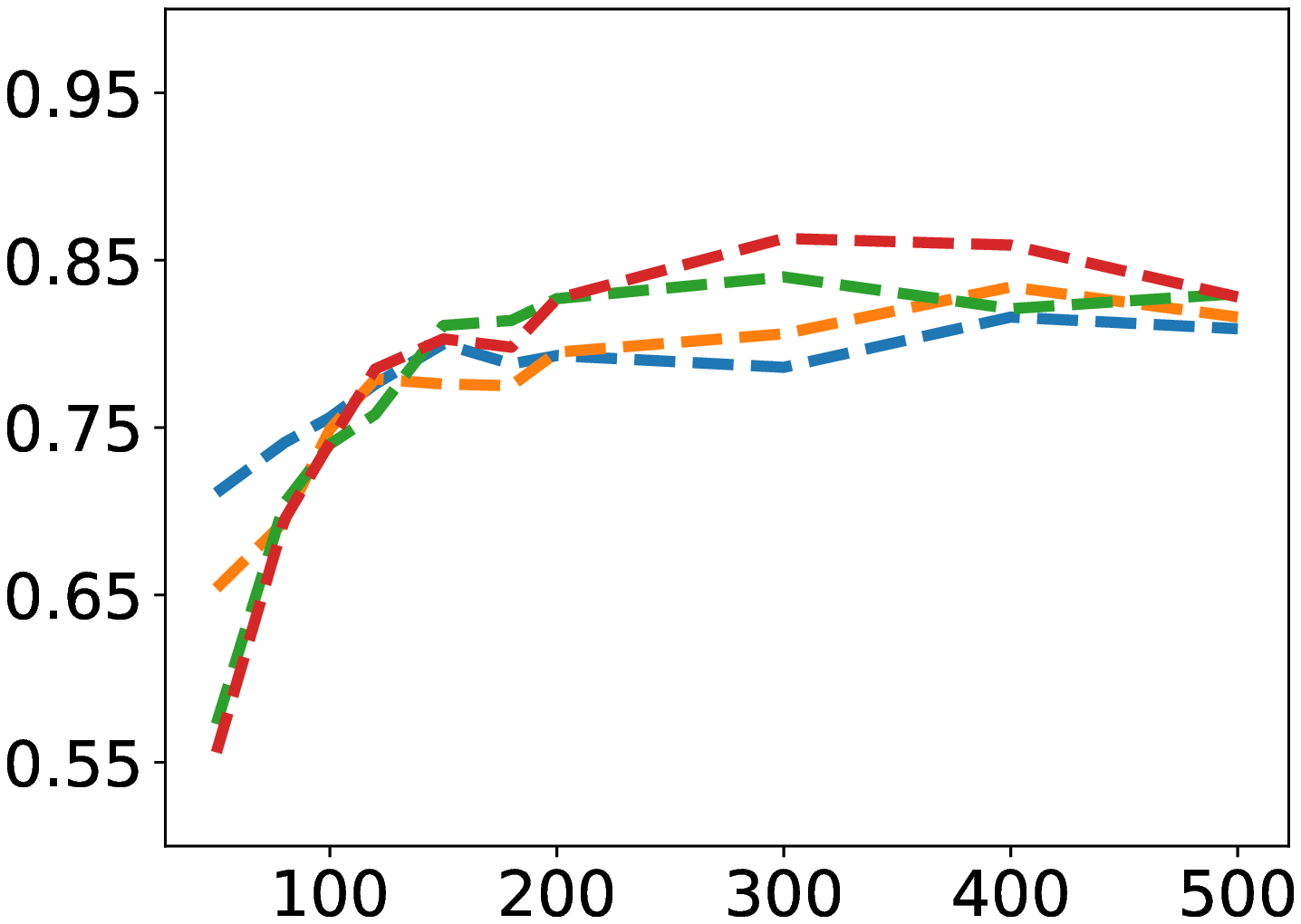} 
    \put(25,80){ \ul{\ \ \  \ $g = 0$ \ \ \ \    }}
	\put(-20,-5){\rotatebox{90}{ {\small \ \ \ log transformation  \ \ }}}
\end{overpic}
	~
	\DeclareGraphicsExtensions{.png}
	\begin{overpic}[width=0.29\textwidth]{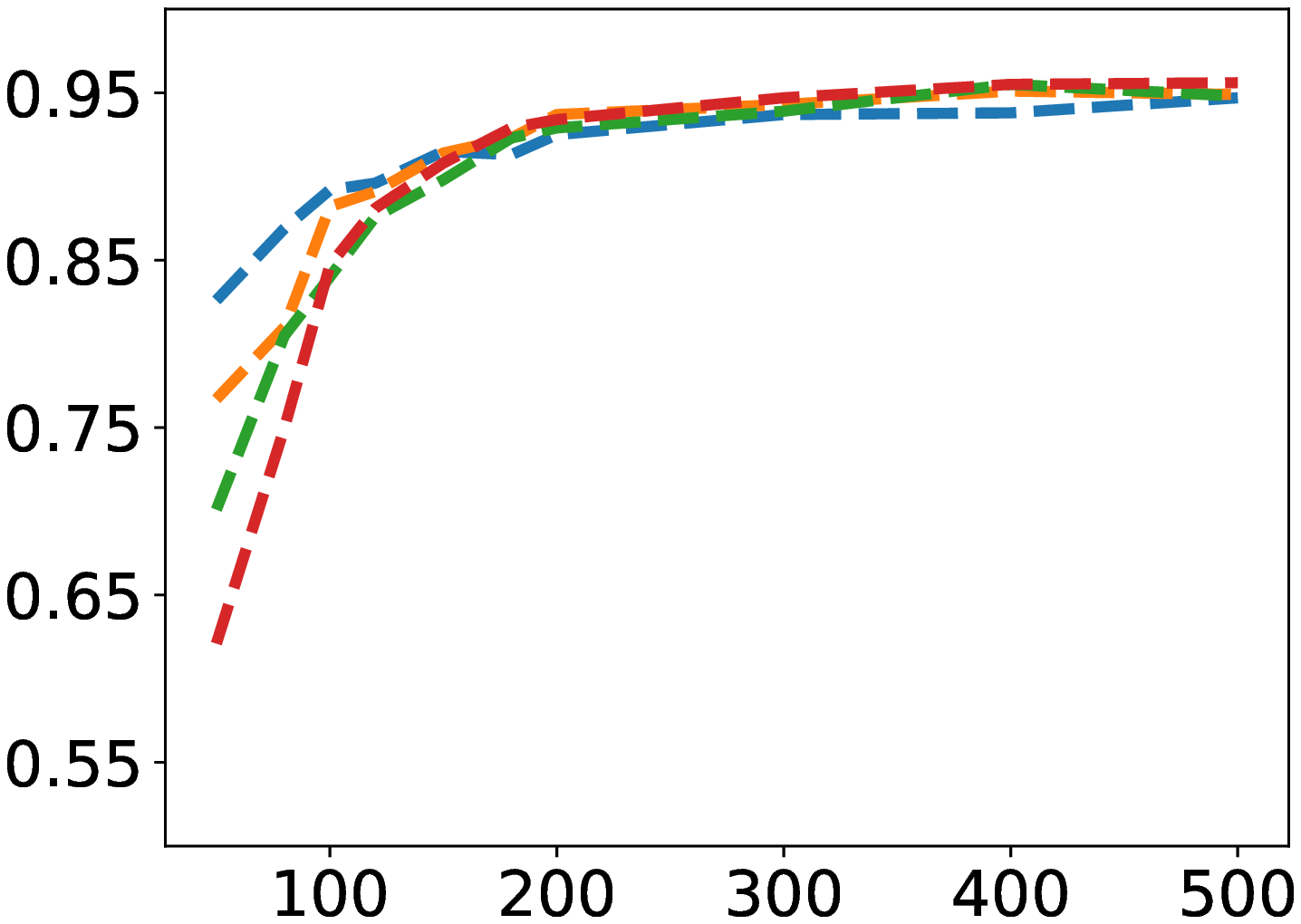} 
	\put(25,80){ \ul{\ \ \  \ $g = 0.1$ \ \ \ \    }}
	\end{overpic}
	~	
	\begin{overpic}[width=0.29\textwidth]{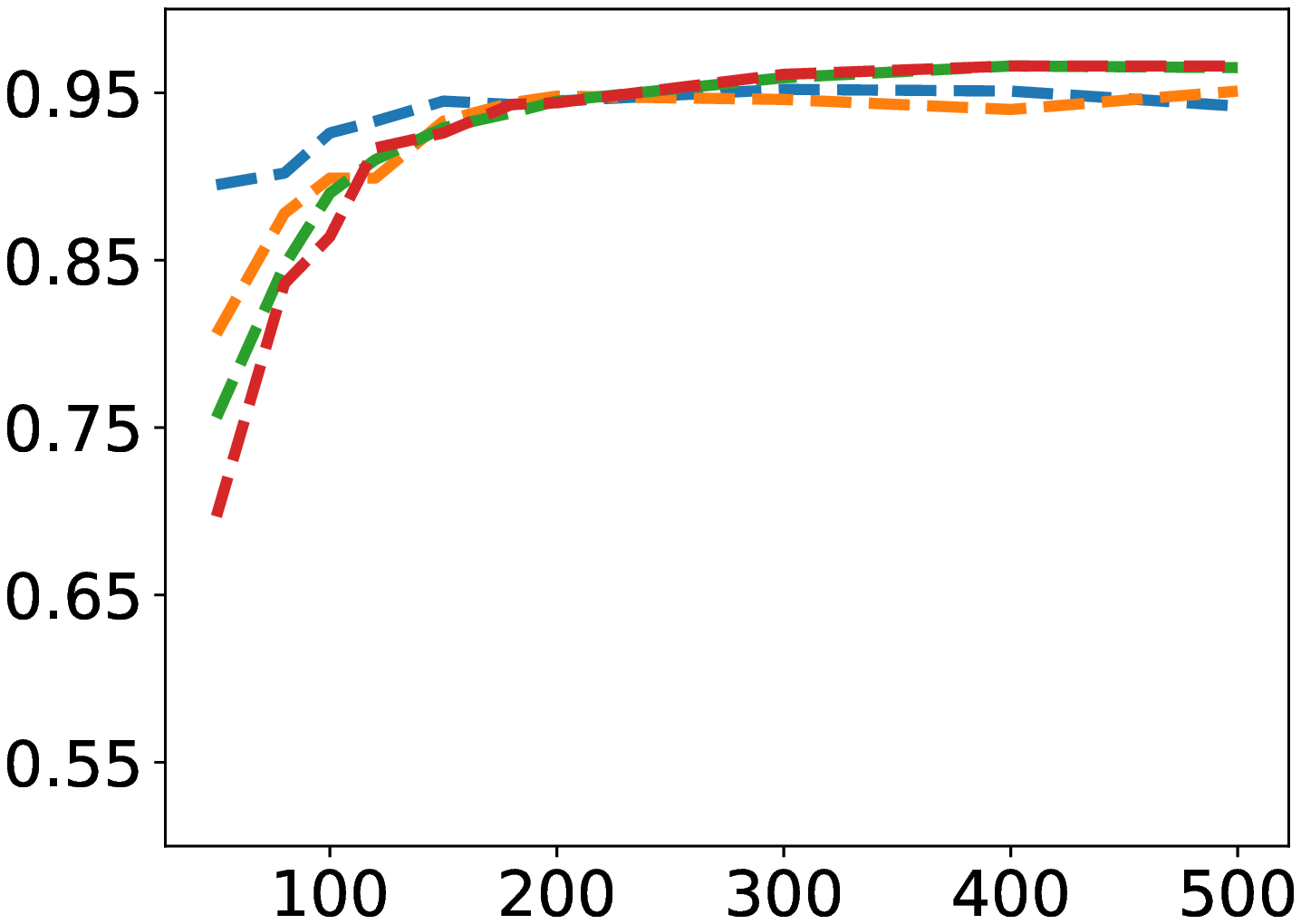} 
	\put(25,80){ \ul{\ \ \  \ $g = 0.2$ \ \ \ \    }}
	\end{overpic}	
\end{figure}

\vspace{-0.5cm}

\begin{figure}[H]	
	\quad\quad\quad 
	\begin{overpic}[width=0.29\textwidth]{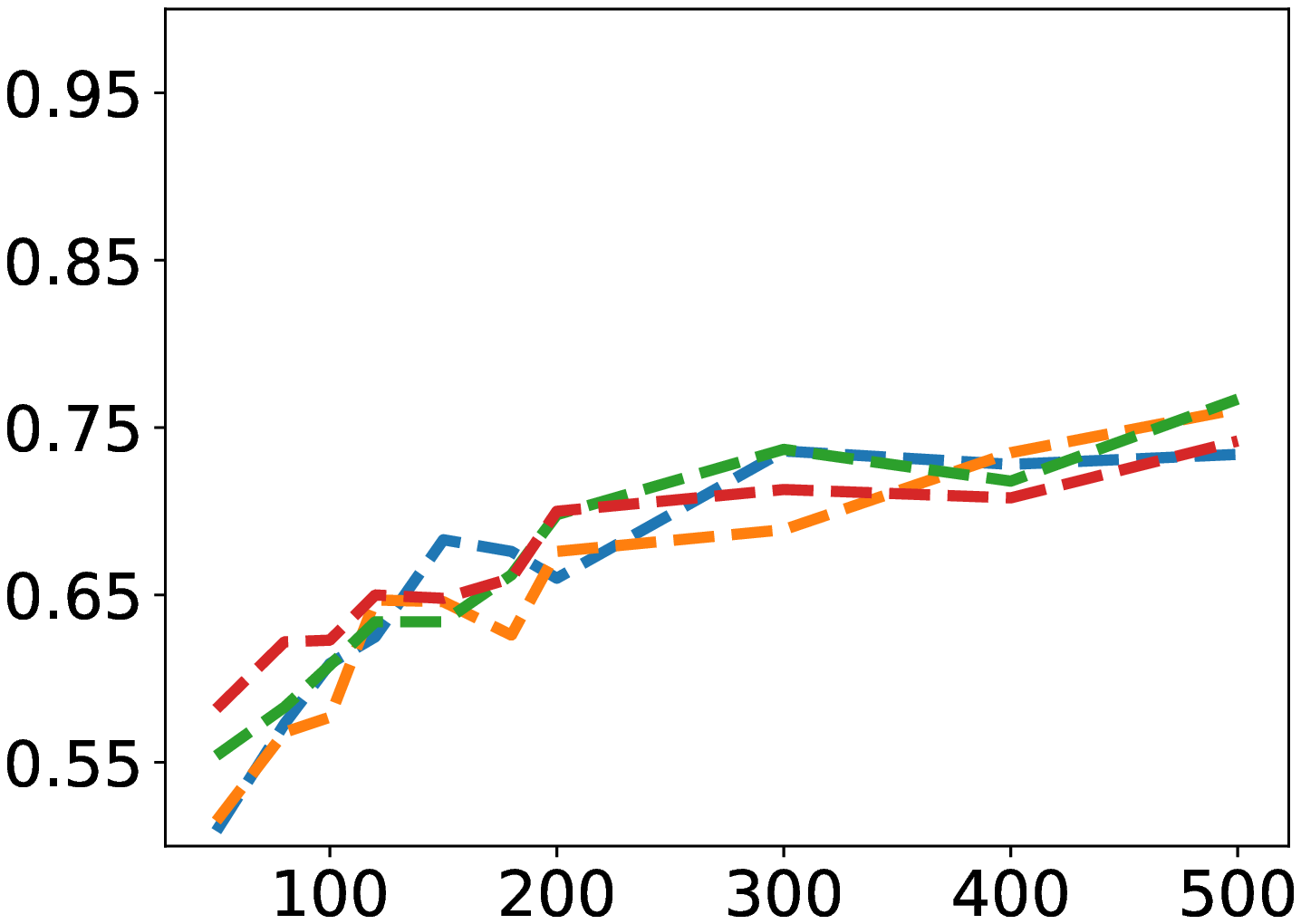} 
	\put(-20,-1){\rotatebox{90}{  {\small \ \ \ standardization \ \ \ } }}

	\end{overpic}
	~
	\DeclareGraphicsExtensions{.png}
	\begin{overpic}[width=0.29\textwidth]{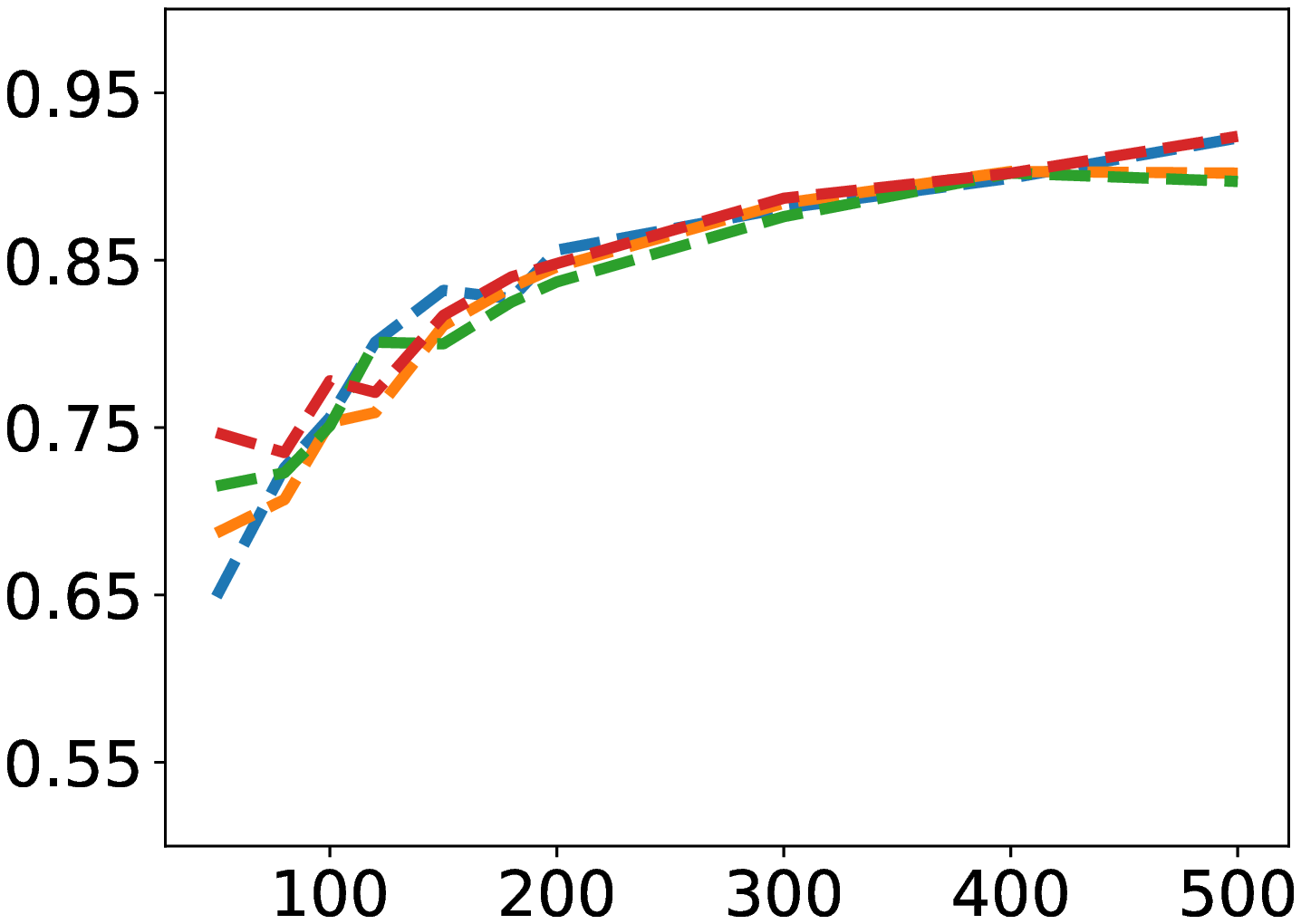} 
	\end{overpic}
	~	
	\begin{overpic}[width=0.29\textwidth]{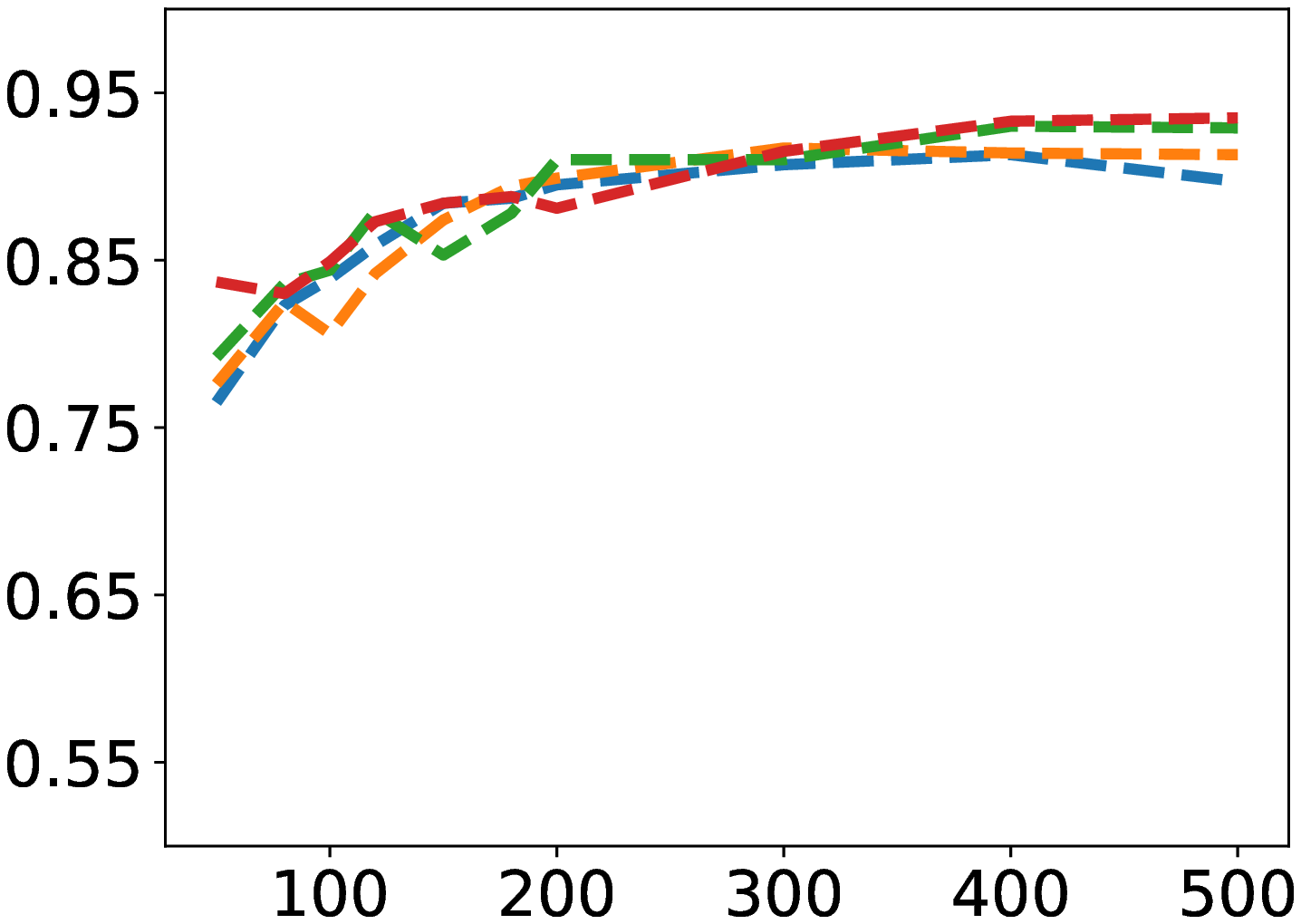} 
		 				
	\end{overpic}	
%
%
%
\end{figure}

\vspace{-0.5cm}

\begin{figure}[H]	
	\quad\quad\quad 
	\begin{overpic}[width=0.29\textwidth]{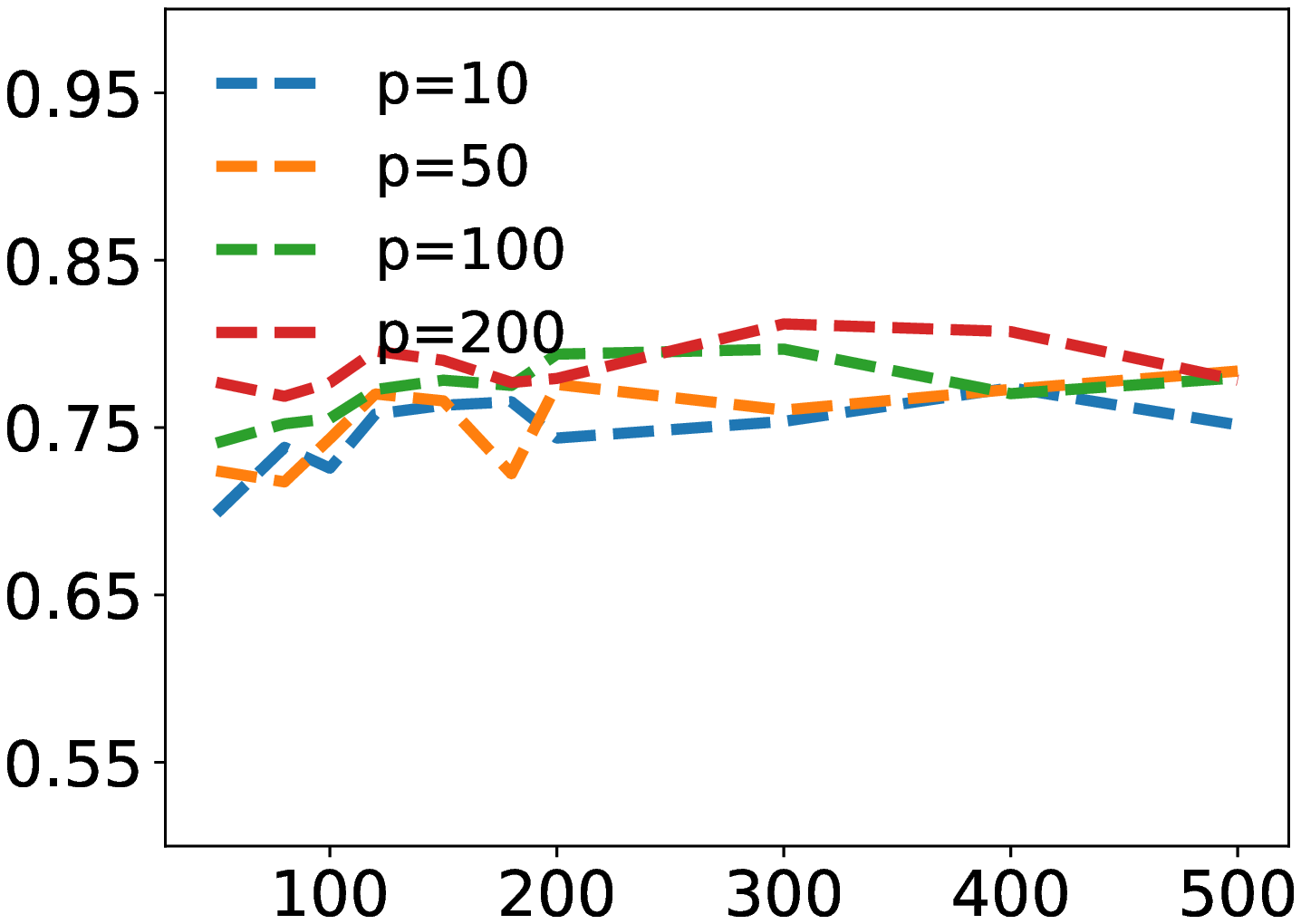} 
    \put(-21,1){\rotatebox{90}{\ $\sqrt{ \ \ }$}}
	\put(-20,-3){\rotatebox{90}{ { \ \ \ \ \ \ \small transformation \ \ }}}
	\end{overpic}
	~
	\DeclareGraphicsExtensions{.png}
	\begin{overpic}[width=0.29\textwidth]{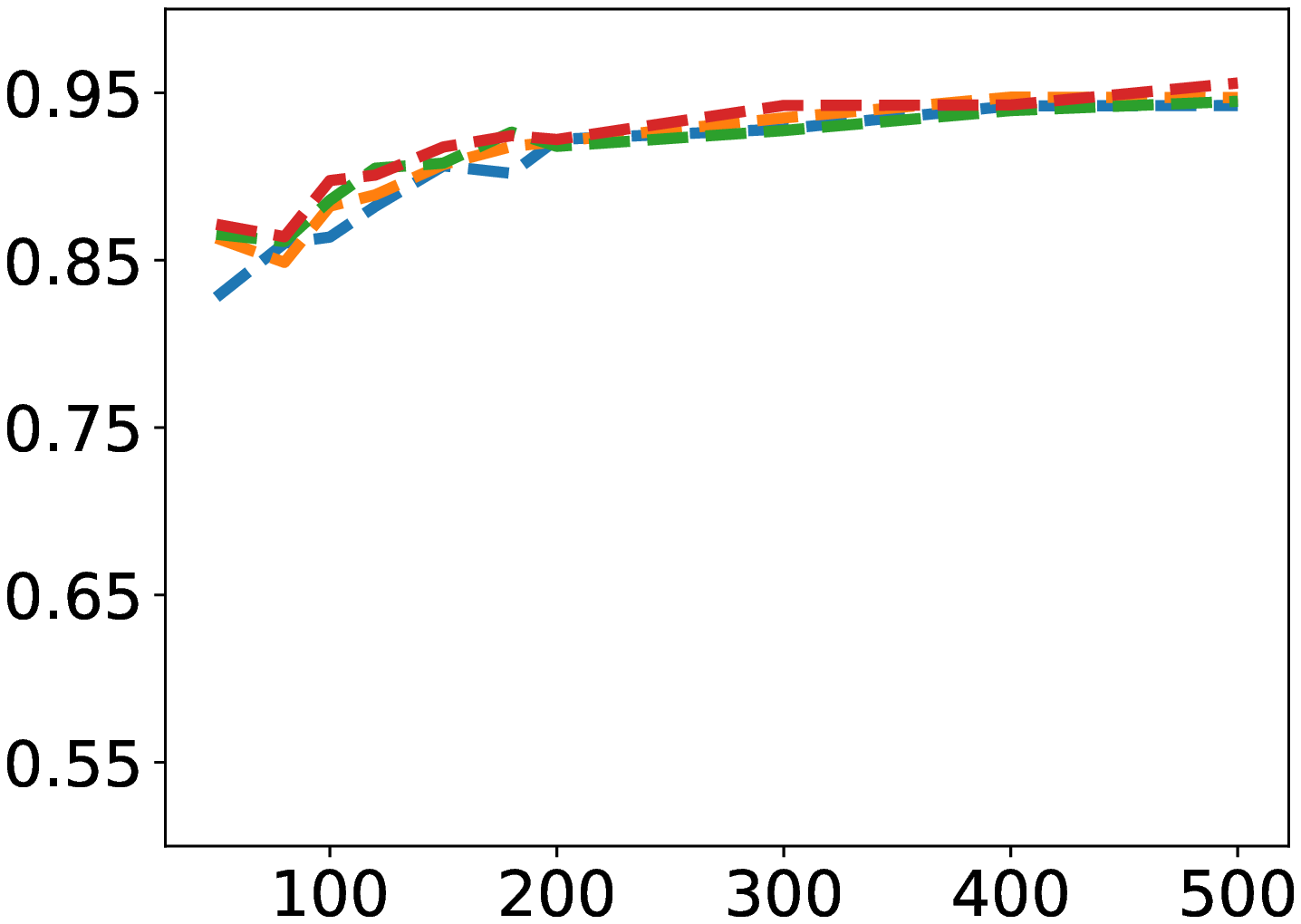} 
	\end{overpic}
	~	
	\begin{overpic}[width=0.29\textwidth]{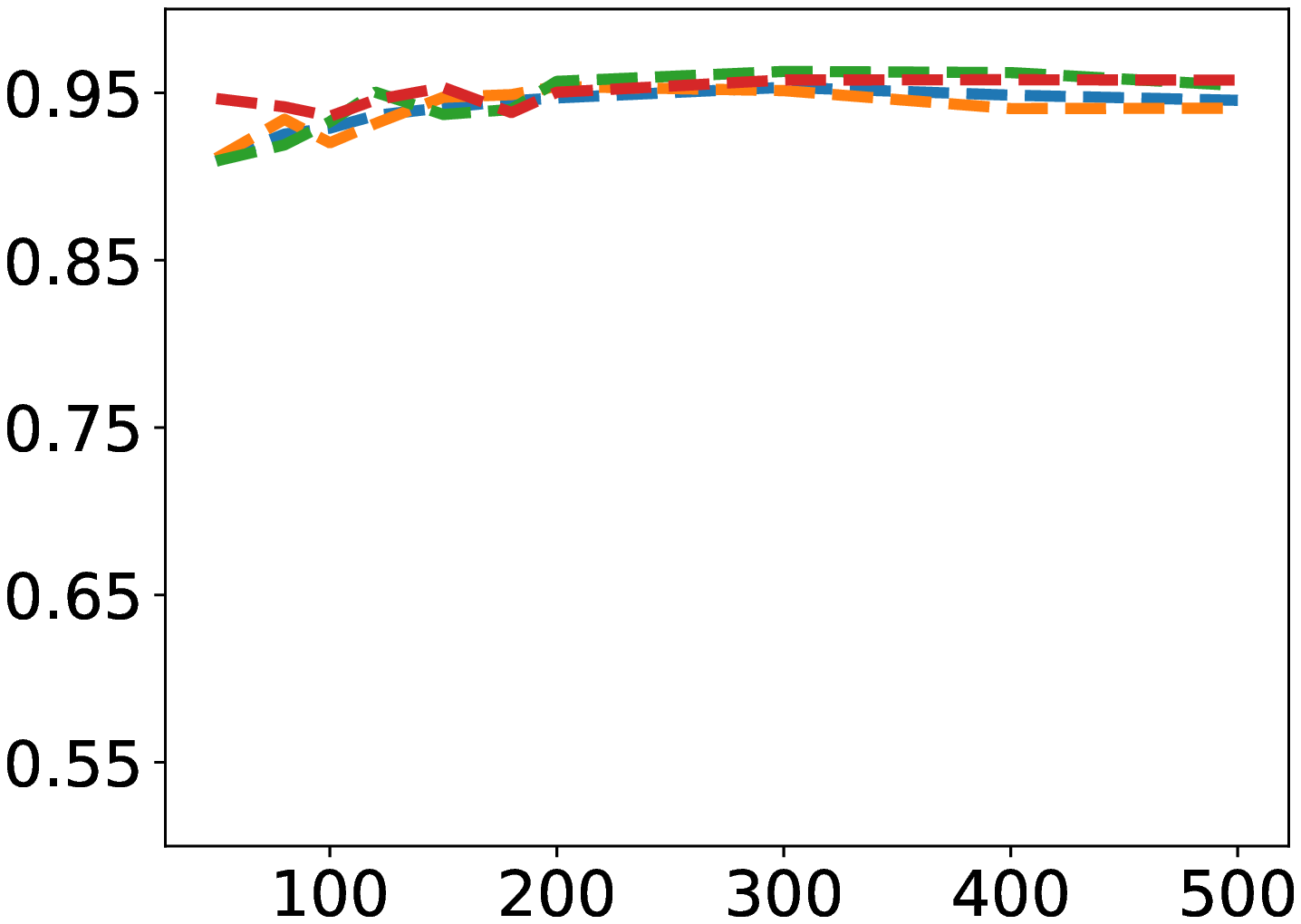} 
	\end{overpic}
	\vspace{+.2cm}
	\caption{(Simultaneous coverage probability versus $n$ in simulation model
(i) with the decay profile: $(\lambda_1)(\Sigma), \lambda_2(\Sigma), \lambda_3(\Sigma)) = (1+g, 1, 1 - g)$ and $\lambda_j(\Sigma) = j^{-1}$ for $j \geq 4$). In each panel, the $y$-axis measures $\P (\cap_{j=1}^5 \{ \lambda_j(\Sigma) \in \hat{\mathcal{I}}_j \})$ based on a nominal level of 95\%, and the $x$-axis measures $n$. The colored curves correspond to the different values of $p$, indicated in the legend.}
	\label{SUPP:extra_sim_elliptical_sqrtexp_95}
\end{figure}

\newpage

\begin{figure}[H]	
\vspace{0.5cm}
	\quad\quad\quad 
	\begin{overpic}[width=0.29\textwidth]{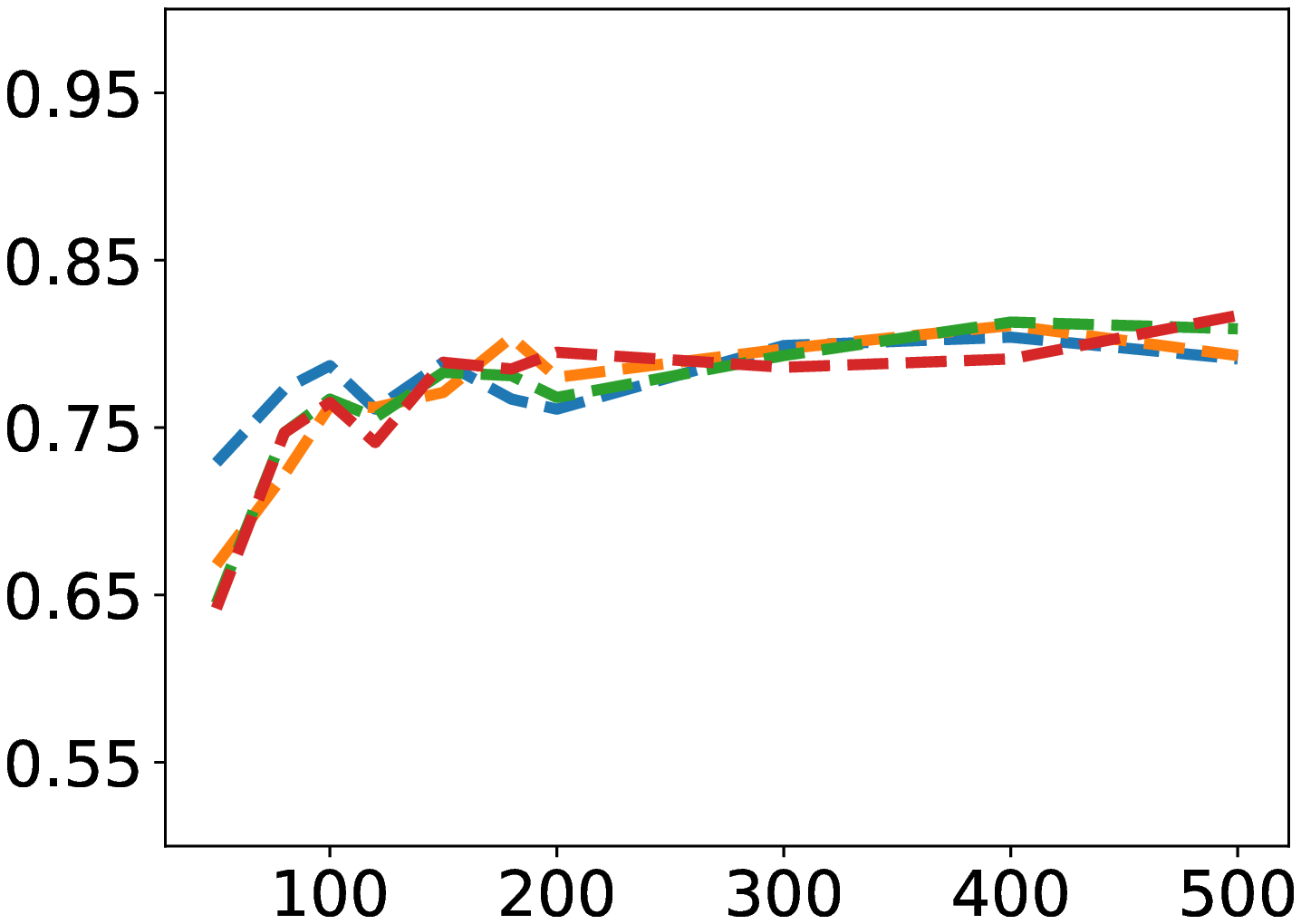} 
    \put(25,80){ \ul{\ \ \  \ $g = 0$ \ \ \ \    }}
	\put(-20,-5){\rotatebox{90}{ {\small \ \ \ log transformation  \ \ }}}
\end{overpic}
	~
	\DeclareGraphicsExtensions{.png}
	\begin{overpic}[width=0.29\textwidth]{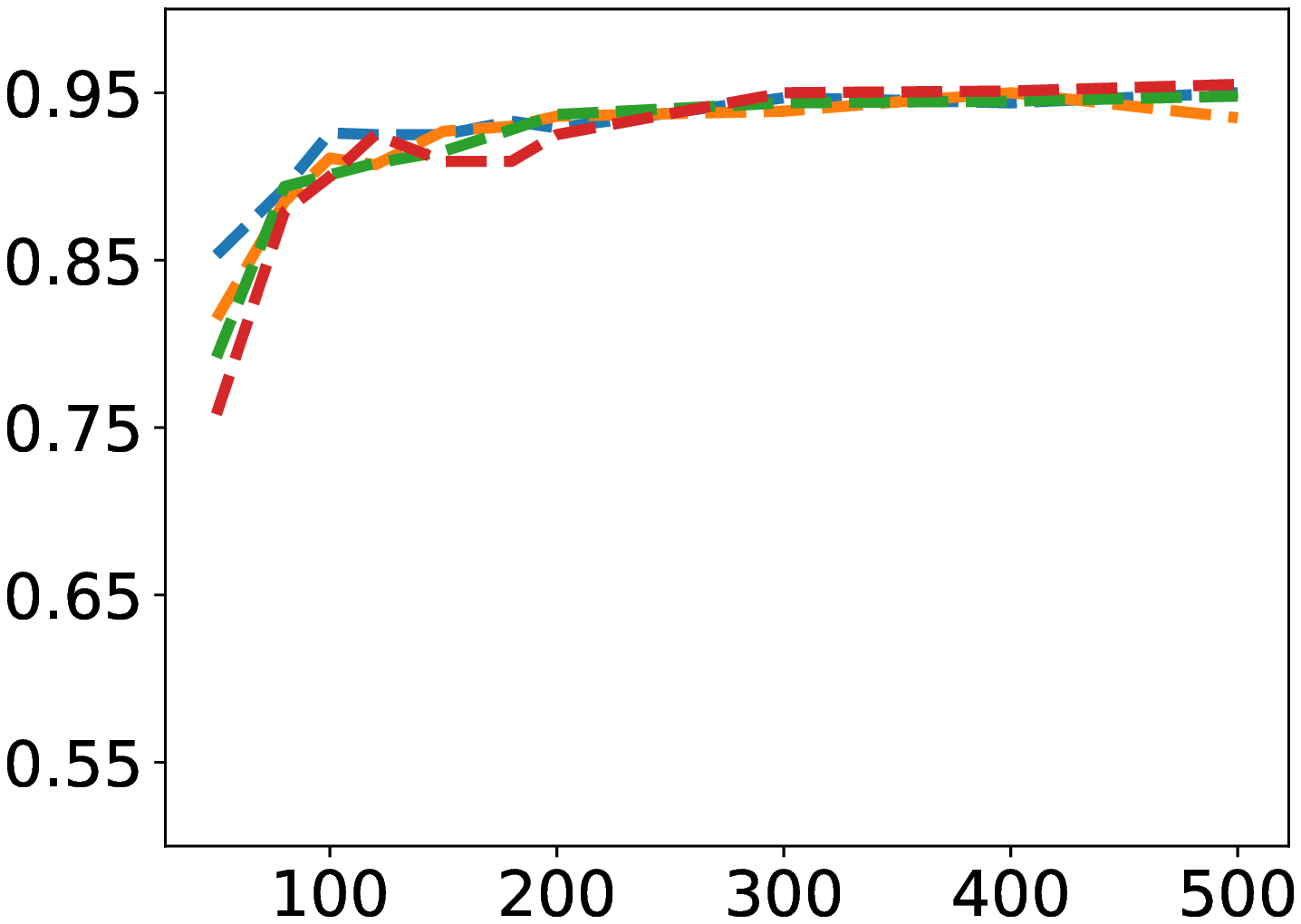} 
	\put(25,80){ \ul{\ \ \  \ $g = 0.1$ \ \ \ \    }}
	\end{overpic}
	~	
	\begin{overpic}[width=0.29\textwidth]{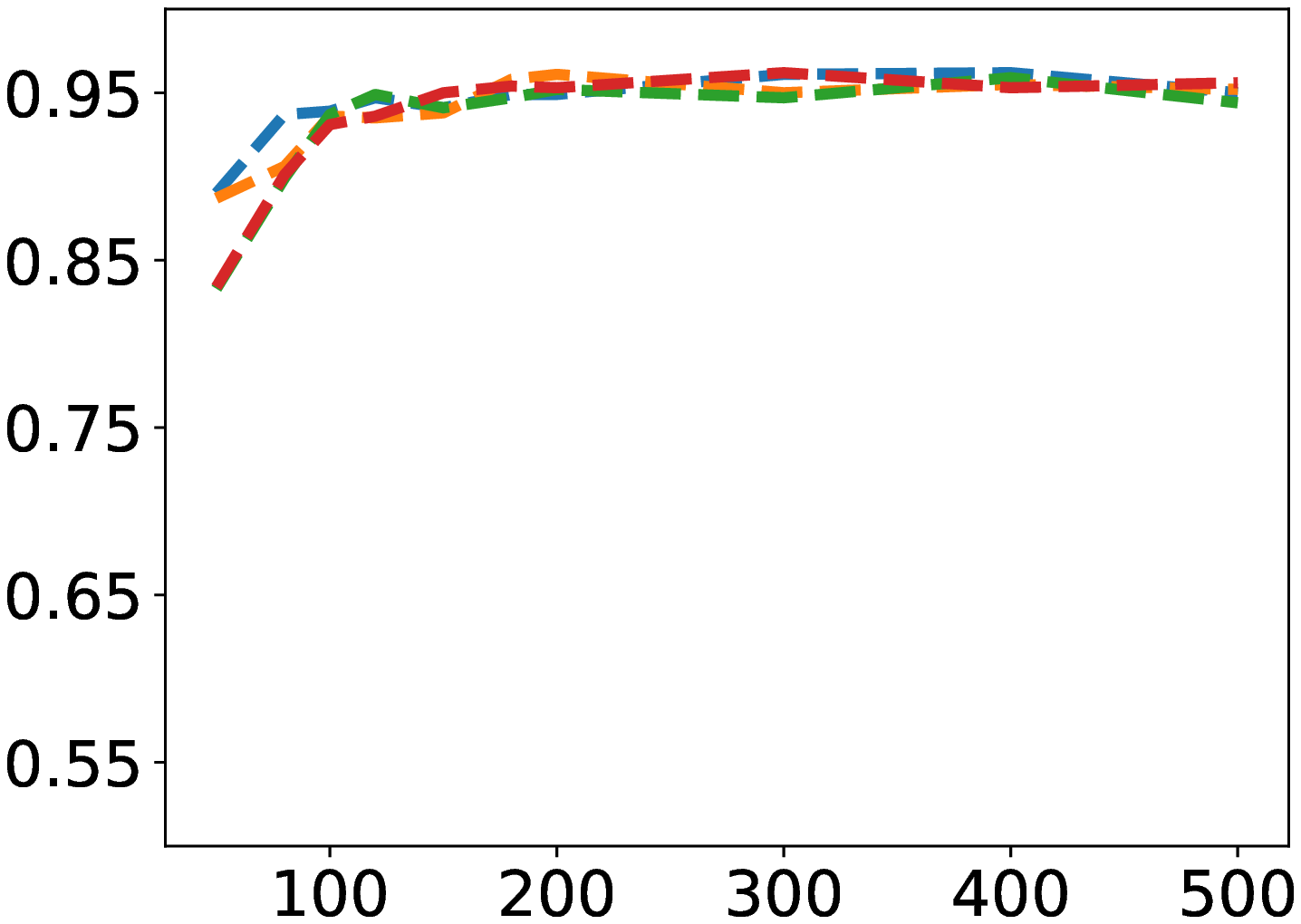} 
	\put(25,80){ \ul{\ \ \  \ $g = 0.2$ \ \ \ \    }}
	\end{overpic}	
\end{figure}

\vspace{-0.5cm}

\begin{figure}[H]	
	\quad\quad\quad 
	\begin{overpic}[width=0.29\textwidth]{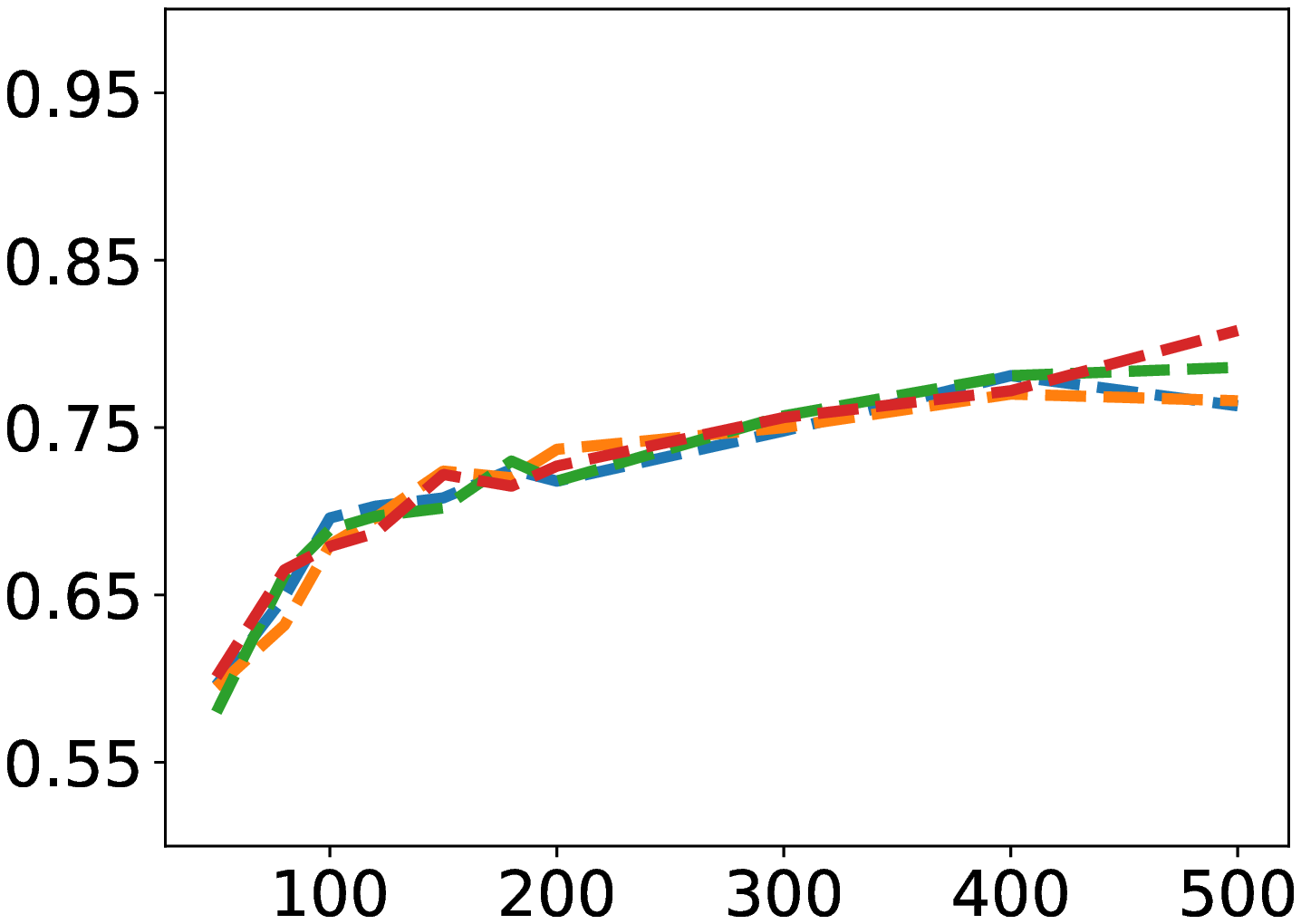} 
	\put(-20,-1){\rotatebox{90}{  {\small \ \ \ standardization \ \ \ } }}

	\end{overpic}
	~
	\DeclareGraphicsExtensions{.png}
	\begin{overpic}[width=0.29\textwidth]{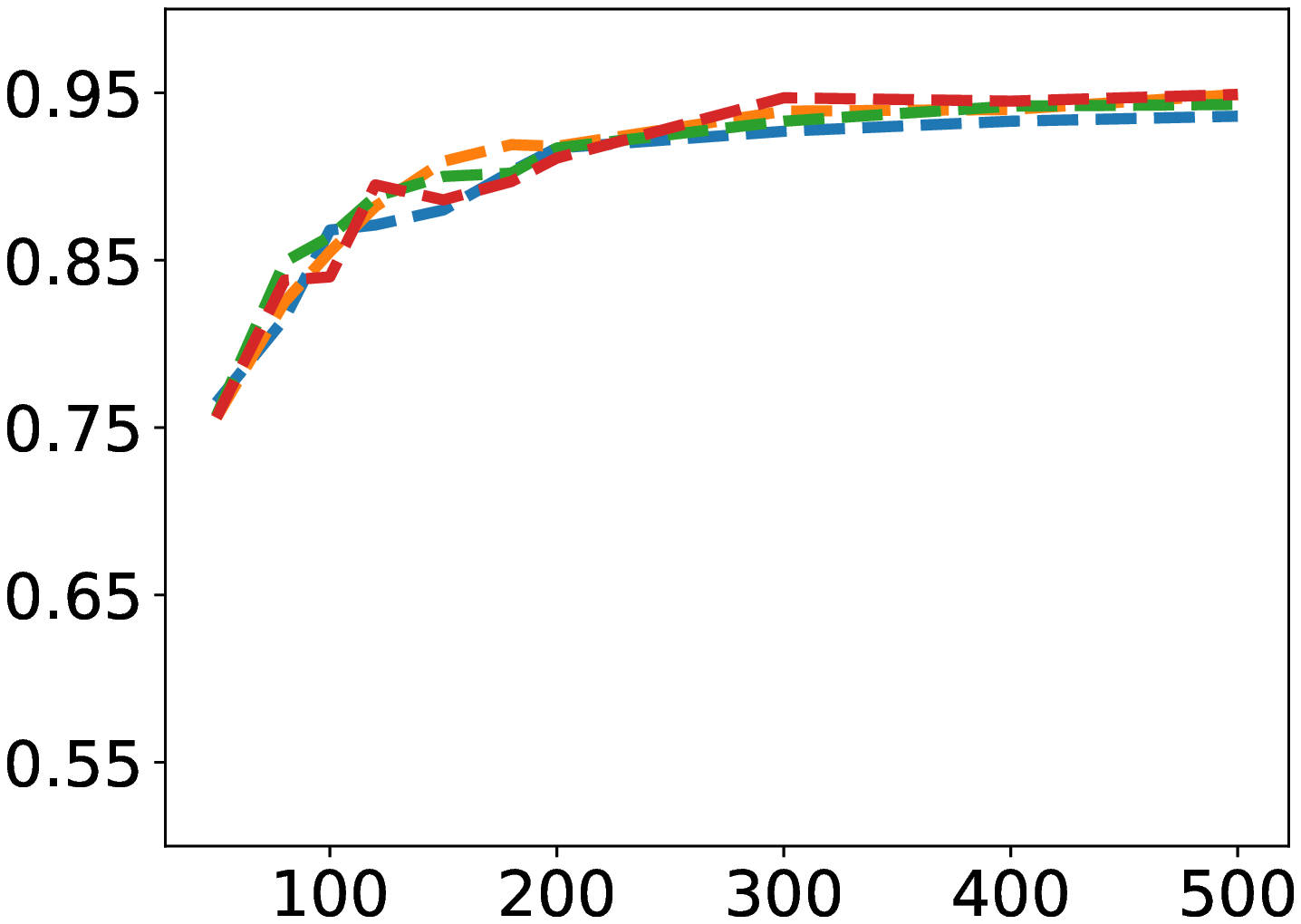} 
	\end{overpic}
	~	
	\begin{overpic}[width=0.29\textwidth]{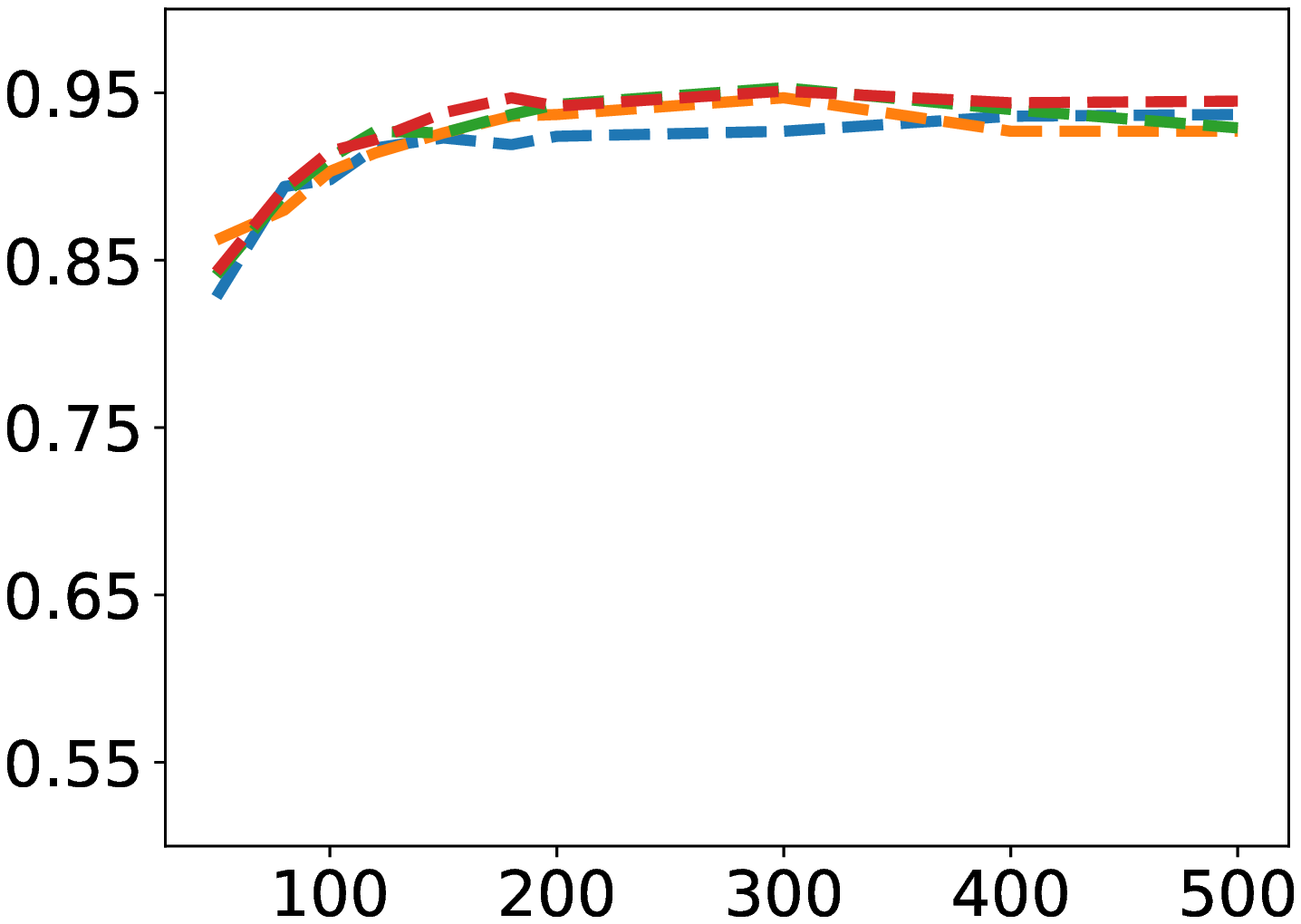} 
		 				
	\end{overpic}	
%
%
%
\end{figure}

\vspace{-0.5cm}

\begin{figure}[H]	
	\quad\quad\quad 
	\begin{overpic}[width=0.29\textwidth]{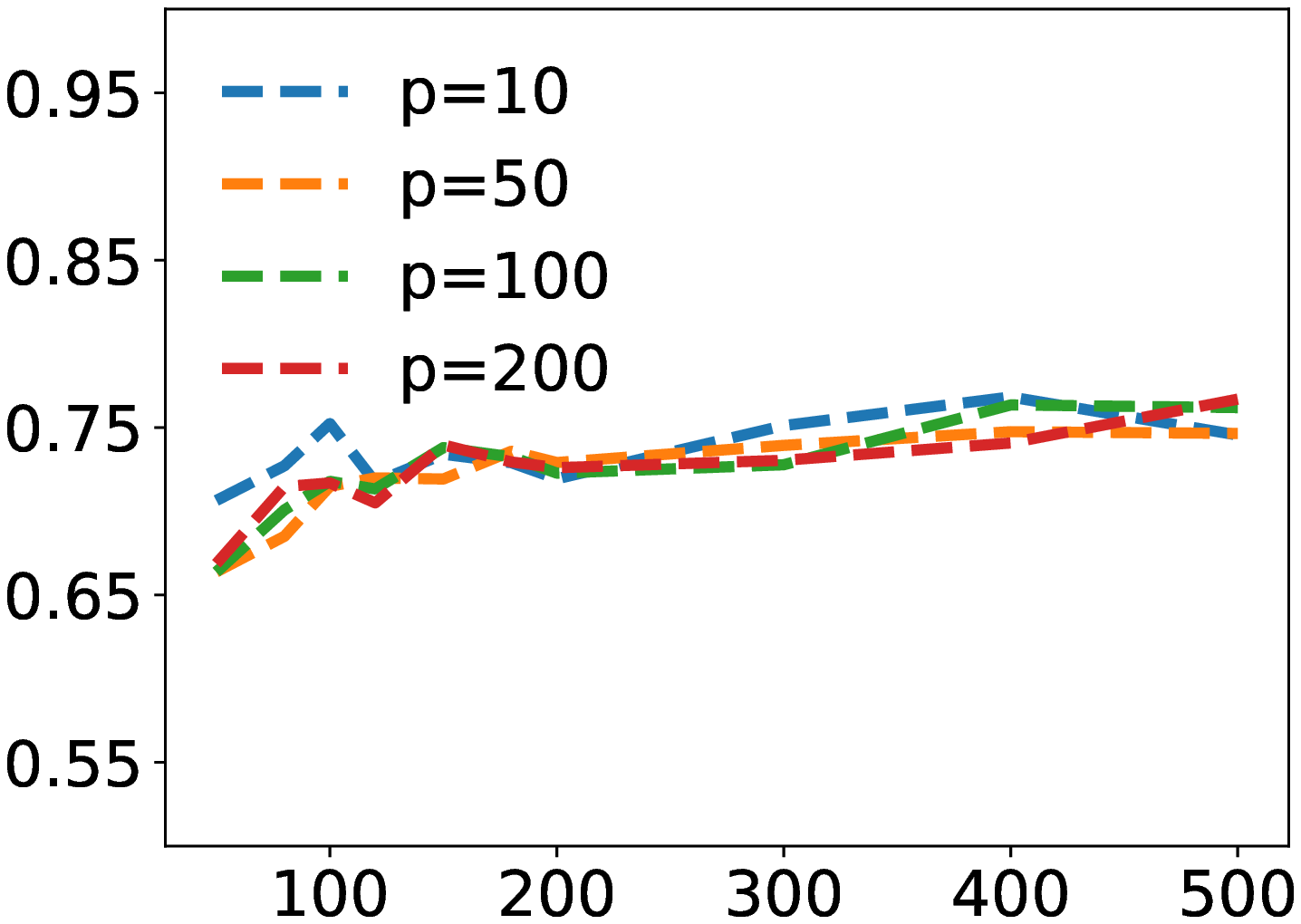} 
    \put(-21,1){\rotatebox{90}{\ $\sqrt{ \ \ }$}}
	\put(-20,-3){\rotatebox{90}{ { \ \ \ \ \ \ \small transformation \ \ }}}
	\end{overpic}
	~
	\DeclareGraphicsExtensions{.png}
	\begin{overpic}[width=0.29\textwidth]{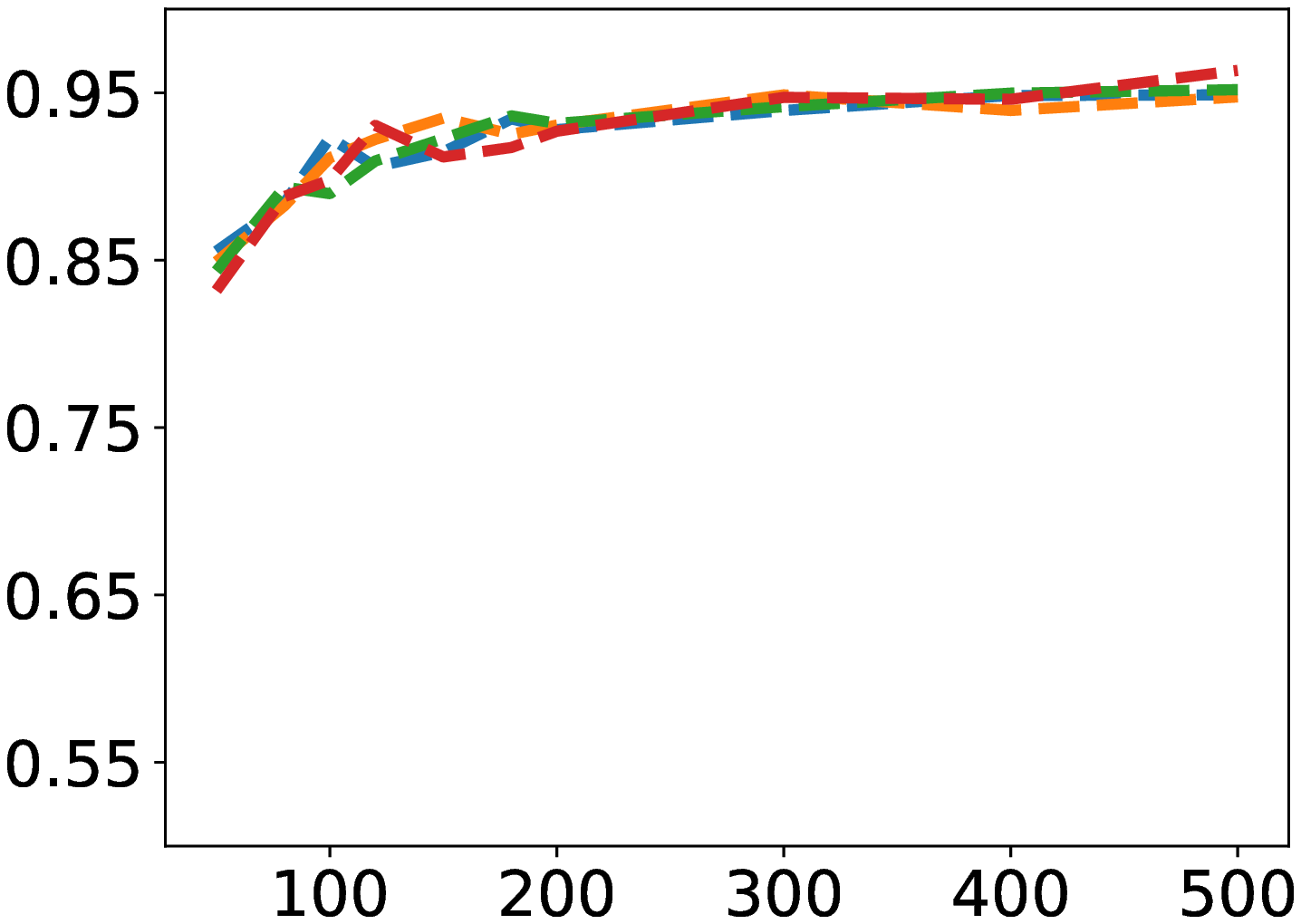} 
	\end{overpic}
	~	
	\begin{overpic}[width=0.29\textwidth]{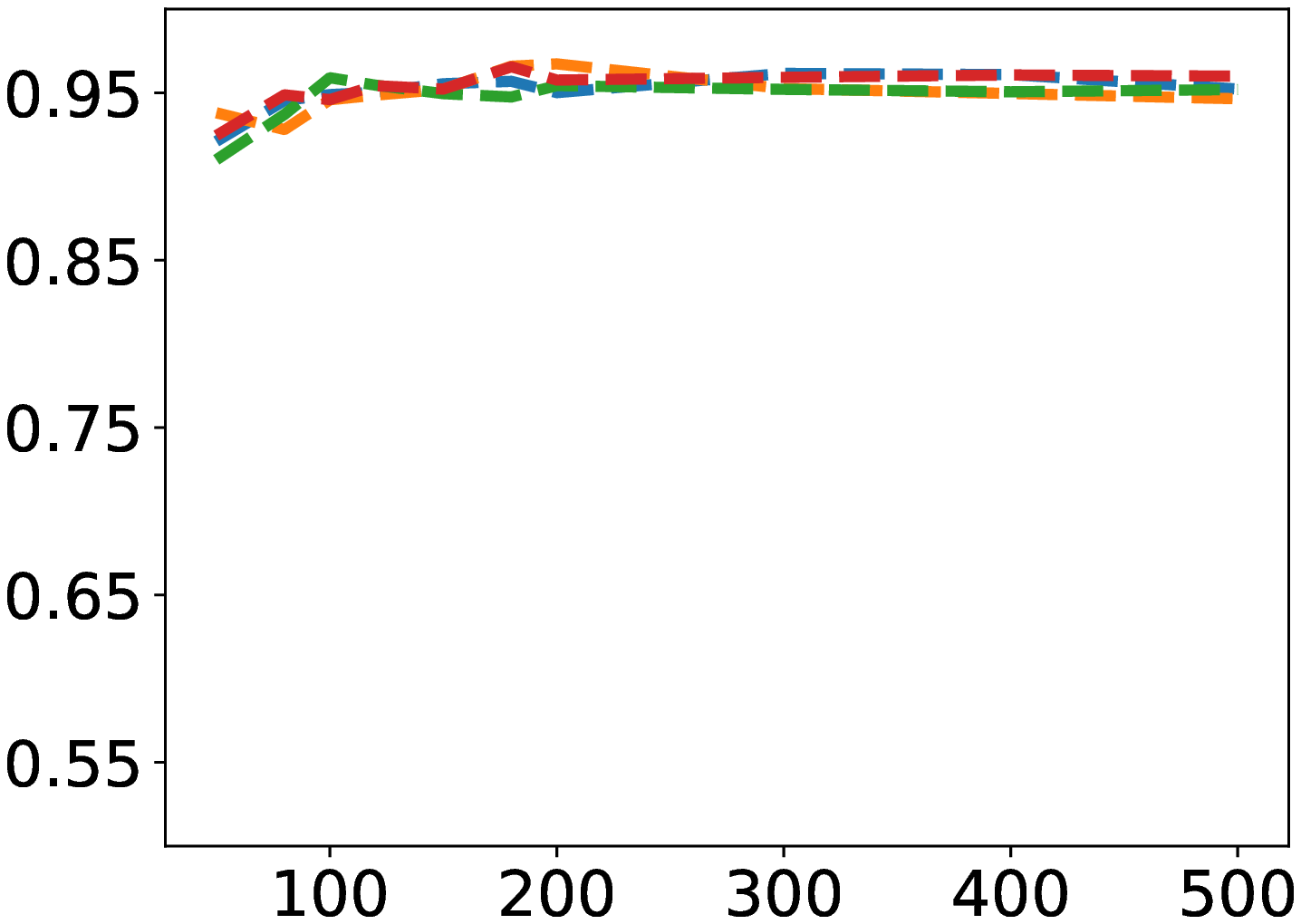} 
	\end{overpic}	
	\vspace{+.2cm}
	\caption{(Simultaneous coverage probability versus $n$ in simulation model
(ii) with the decay profile: $(\lambda_1)(\Sigma), \lambda_2(\Sigma), \lambda_3(\Sigma)) = (1+g, 1, 1 - g)$ and $\lambda_j(\Sigma) = j^{-1}$ for $j \geq 4$). The plotting scheme is the same as described in the caption of Figure~\ref{SUPP:extra_sim_elliptical_sqrtexp_95}.}
	\label{SUPP:extra_sim_iid_gaussian_95}
\end{figure}

\end{document}